Московский физико-технический институт
(государственный университет)
Кафедра Математических основ управления
Факультета управления и прикладной математики

На правах рукописи

УДК 519.8

**Гасников Александр Владимирович**

# Эффективные численные методы поиска равновесий в больших транспортных сетях

Диссертация на соискание ученой степени
доктора физико-математических наук
по специальности 05.13.18 – Математическое моделирование,
численные методы и комплексы программ

Научный консультант:
доктор физико-математических наук,
чл.-корр. РАН, профессор Шананин А.А.

МОСКВА – 2016

# Оглавление

# Введение

## Актуальность темы и степень ее разработанности

Настоящая диссертация посвящена разработке новых подходов к построению многостадийных моделей транспортных потоков и эффективных численных методов поиска равновесий в таких моделях. Начиная с 50-х годов XX века вопросам поиска равновесий в транспортных сетях стали уделять большое внимание в связи с ростом городов и необходимостью соответствующего транспортного планирования. В 1955 г. появилась первая модель равновесного распределения потоков по путям: BMW-модель, также называемая моделью Бэкмана [147]. В этой модели при заданных корреспонденциях (потоках из источников в стоки) решалась задача поиска равновесного распределения этих корреспонденций по путям, исходя из принципа Вардропа, т.е. исходя из предположения о том, что каждый пользователь транспортной сети рационален и выбирает кратчайший маршрут следования. Таким образом, поиск равновесия в такой модели сводился к поиску равновесия Нэша в популяционной игре [302] (популяций столько, сколько корреспонденций). Поскольку в модели предполагалось, что время прохождения ребра есть функция от величины потока только по этому ребру, то получившаяся игра была игрой загрузки, следовательно, потенциальной. Последнее означает, что поиск равновесия сводится к решению задачи оптимизации. Получившуюся задачу выпуклой оптимизации решали с помощью метода условного градиента [203]. Описанная модель и численный метод и по настоящее время используются в подавляющем большинстве продуктов транспортного моделирования для описания блока равновесного распределения потоков по путям [112, 286, 289, 308]. Однако в работе [275] было указано на ряд существенных недостатков модели Бэкмана, и предложена альтернативная модель, которую авторы назвали моделью Стабильной Динамики.

Отмеченные выше модели равновесного распределения потоков по путям могут быть использованы при решении различных задач долгосрочного планирования. Например, такой задачи. Имея заданный бюджет, нужно решить, на каких участках графа транспортной сети стоит увеличить полосность дороги / построить новые дороги. Заданы несколько сценариев, нужно отобрать лучший. Задачу можно решить, найдя равновесные распределения потоков, отвечающие каждому из сценариев, и сравнивая найденные решения, например, по критерию суммарного времени, потерянного в пути всеми пользователями сети в данном равновесии. При значительных изменениях графа транспортной сети необходимо в приведенную выше цепочку рассуждений включать дополнительный контур, связанный с тем, что изменения приведут не только к перераспределению потоков на путях, но и поменяют корреспонденции. Таким образом, корреспонденции также необходимо моделировать. В 60-е годы XX века появилось сразу несколько различных моделей для расчета матрицы корреспонденций, исходя из информации о численностях районов и числе рабочих мест в них. Наибольшую популярность приобрела энтропийная модель расчета матрицы корреспонденций [26]. В этой модели поиск матрицы корреспонденций сводился к решению задачи энтропийно-линейного программирования.

К сожалению, при этом в энтропийную модель явным образом входит информация о матрице затрат на кратчайших путях по всевозможным парам районов. Возникает "порочный круг": чтобы посчитать эту матрицу затрат, нужно сначала найти равновесное распределение потоков по путям, а чтобы найти последнее, необходимо знать матрицу корреспонденций, которая рассчитывается по матрице затрат. На практике отмеченную проблему решали методом простых итераций. Как-то "разумно" задавали начальную матрицу корреспонденций, по ней считали распределение потоков по путям, по этому распределению считали матрицу затрат, на основе которой пересчитывали матрицу корреспонденций, и процесс повторялся. Повторялся он до тех пор, пока не выполнялся критерий останова. К сожалению, до сих пор для описанной процедуры не известно никаких гарантий ее сходимости и, тем более, оценок скорости сходимости [286].

Описанный выше подход можно назвать двухстадийной моделью транспортных потоков, потому что модель состоит из прогонки двух разных блоков. В действительности, в реальных приложениях, число блоков 3-4 [286]. В частности, как правило, всегда включают блок расщепления потоков по типу передвижения (например, личный и общественный транспорт) – этот блок описывается моделью аналогичной модели Бэкмана. В математическом плане это уточнение не существенно. Все основные имеющиеся тут сложности хорошо демонстрирует уже двухстадийная модель. Отметим также, что часто в приложениях вместо модели Бэкмана используется её "энтропийно-регуляризованный" вариант, который отражает ограниченную рациональность пользователей транспортной сети [137, 308]. Равновесие в такой модели часто называют "стохастическим" равновесием.

Резюмируем написанное выше. До настоящего момента не существовало строгого научного обоснования используемого повсеместно на практике (и зашитого во все современные пакеты транспортного моделирования) способа формирования многостадийных моделей транспортных потоков. Не существовало также никаких гарантий сходимости численных методов, используемых для поиска равновесий в многостадийных моделях. В используемых сейчас повсеместно многостадийных моделях в качестве основных блоков фигурируют блоки с моделями типа Бэкмана, а не более современные блоки Стабильной Динамики. Таким образом, актуальной является задача обоснования современной многостадийной модели и разработка эффективных численных методов поиска (стохастического) равновесия в такой модели.

**Цели и задачи**

Многие законы природы могут быть записаны в форме “вариационных принципов”. В моделировании транспортных потоков это также имеет место. Однако, на текущий момент с помощью вариационных принципов описываются только отдельные блоки многостадийной транспортной модели, и эволюционный вывод вариационных принципов имеется только для блоков с моделью Бэкмана в основе. Одной из целей данной диссертационной работы является эволюционный вывод всех блоков многостадийной транспортной модели (прежде всего, речь идет о блоке расчета матрицы корреспонденций), и получение (с помощью эволюционного вывода) вариационного принципа для описания равновесия в многостадийной модели.

Целью также является разработка “алгебры” над блоками-моделями (каждый блок описывается своим вариационным принципом), которая позволит, как в конструкторе, собирать (формируя общий вариационный принцип) сколько угодно сложные модели из небольшого числа базисных элементов конструктора (блоков).

Описанный выше формализм приводит в итоге к решению задач выпуклой оптимизации в пространствах огромных размеров и имеющих довольно специальную иерархическую (многоуровневую) структуру функционала задачи. Чтобы подчеркнуть нетривиальность таких задач, отметим, что переменные, по которым необходимо оптимизировать, – это, в частности, компоненты вектора распределения потоков по путям. Для графа в виде двумерной квадратной решетки (Манхэттенская сеть) с числом вершин порядка нескольких десятков тысяч (это число соответствует транспортному графу Москвы) такой вектор с большим запасом нельзя загрузить в память любого современного суперкомпьютера, не говоря уже о том, чтобы как-то работать с такими векторами.

Важной целью диссертации является разработка (с теоретическими гарантиями) эффективных численных методов, способных за несколько часов на персональном компьютере с хорошей точностью (и с высокой вероятностью) найти равновесие в многостадийной модели транспортных потоков крупного мегаполиса.

В частности, целью является разработка “алгебры” над численными методами, используемыми для расчета отдельных блоков многостадийной модели, которая позволит, как в конструкторе, собирать итоговую эффективную численную процедуру (для поиска равновесия в многостадийной модели) с помощью правильного чередования / комбинации работы алгоритмов, используемых для отдельных блоков.

Более общей целью является выделение небольшого набора “оптимальных” базисных численных методов выпуклой оптимизации и операций над ними, чтобы с помощью всевозможных сочетаний можно было получать “оптимальные” методы для структурно сложных задач выпуклой оптимизации. То есть цель – научиться раскладывать (декомпозировать) сложно составленную задачу выпуклой оптимизации на простые блоки, чтобы численный метод для общей задачи можно было бы “собирать” из простых базисных блоков (методов). Желательно также, чтобы разработанный формализм позволял автоматизировать эту процедуру.

**Научная новизна, методология и методы исследования**

В основе предложенного в диссертации эволюционного формализма обоснования многостадийной транспортной модели лежит часто используемая в популяционной теории игр марковская logit-динамика, отражающая ограниченную рациональность агентов (водителей) [302]. Новым элементом является понимание этой динамики как модели стохастической химической кинетики с унарными реакциями и рассмотрение сразу нескольких разных типов таких унарных реакций, происходящих с разной (по порядку величины) интенсивностью, и отвечающих разномасштабным процессам, протекающим в транспортной системе. Например, для двухстадийной модели динамика, отвечающая формированию корреспонденций, идет, по терминологии А.Н. Тихонова, в медленном времени (годы), а отвечающая распределению потоков по путям – в быстром времени (дни). Тогда с некоторыми оговорками функционал в вариационном принципе с точностью до множителя и аддитивной константы можно, с одной стороны, понимать как функционал Санова (действие), отвечающий за концентрацию стационарной (инвариантной) меры введенной марковской динамики, а, с другой стороны, как функционал Ляпунова–Больцмана для кинетической динамики, полученной при (каноническом) скейлинге (по числу агентов) введённой марковской динамики.

При разработке эффективных численных процедур в качества базиса были выбраны два метода [93, 162, 259]: метод зеркального спуска – МЗС (А.С. Немировский, 1977) и быстрый градиентный метод – БГМ (Ю.Е. Нестеров, 1982). Практически все используемые в диссертации алгоритмы являются некоторыми вариантами (производными) этих двух методов.

В диссертации часто использовались следующие операции над алгоритмами: процедура рестартов, регуляризация функционала задачи, mini-batc'инг [162, 259]. По отдельности эти процедуры были давно и хорошо известны. Однако в работе были предложены различные сочетания отмеченных операций, позволившие получить часть результатов.

Важное место в диссертации занимает "игра" между стоимостью итерации численного метода и необходимым числом итераций. Здесь было рассмотрено два направления.

1. В большинстве приложений "стоимость" (время) получения от оракула (роль которого, как правило, играют нами же написанные подпрограммы вычисления градиента) градиента функционала заметно превышает время, затрачиваемое на то, чтобы сделать шаг итерации, исходя из выданного оракулом вектора. Желание сбалансировать это рассогласование (усложнить итерации, сохранив при этом старый порядок их сложности, и выиграть за счет этого в сокращении числа итераций) привело к возникновению композитной оптимизации [266], в которой (аддитивная) часть функционала задачи переносится без лианеризации (запроса градиента) в итерации. Здесь остается еще много степеней свободы, позволяющих играть на том, насколько "дорогим" окажется оракул и соответствующая этому оракулу "процедура проектирования", и на том, сколько (внешних) итераций потребуется методу для достижения заданной точности. В частности, если обращение к оракулу за градиентом и последующее проектирование требуют, в свою очередь, решения вспомогательных оптимизационных задач, то можно "сыграть" на том, насколько точно надо решать эти вспомогательные задачи, пытаясь найти компромисс между "стоимостью" итерации и числом итераций. Также можно сыграть и на том, как выделять эти вспомогательные задачи – иными словами, что понимать под оракулом и что – под итерацией метода. Общая идея подхода "разделяй и властвуй" к численным методам выпуклой оптимизации может принимать довольно неожиданные и при этом весьма эффективные формы (как например, в методах внутренней точки Нестерова–Немировского, 1989 [93, 259]). В диссертации рассмотрены разные оригинальные варианты описанной игры в связи с транспортно-сетевыми приложениями рассмотрены в диссертации. Заметим, что для обоснования упомянутых конструкций (подходов) в диссертации существенным образом используется (и развивается) концепция неточного оракула, выдающего градиент и значение функции [183]. Эта же концепция существенным образом используется (и развивается) в диссертации при разработке линейки универсальных методов [274], самонастраивающихся на гладкость задачи. Общая идея таких методов – искусственно вводить неточность, чтобы правильно подбирать класс гладкости задачи.

2. Достаточно часто имеет смысл (с точки зрения минимизации общего времени работы метода) вместо градиента в итерационную процедуру подставлять некоторый (легко вычислимый) аналог градиента. Как правило, это несмещенная (или мало смещенная) оценка градиента (стохастический градиент). Число итераций при этом возрастает, но зато итерации становятся более “дешевыми”. В диссертации описываются общие “рецепты” получения из детерминированных (полноградиентных) методов их стохастических (рандомизированных) вариантов, позволяющие достаточно просто изучать “наследуемые” при этом свойства исходных методов. К таким рандомизированным методам можно отнести, например, спуски по направлению, покомпонентные методы и методы нулевого порядка. Нетривиально уже то, что это оказалось возможным. Например, покомпонентный БГМ был предложен практически независимо от полноградиентного варианта БГМ [264]. В диссертации продемонстрировано, как можно достаточно просто получить покомпонентный БГМ с наследованием всех основных свойств из БГМ в специальной форме [135]. Полезно заметить, что в такой форме БГМ можно проинтерпретировать, как выпуклую комбинацию МЗС и обычного метода проекции градиента. Важное место в диссертации занимает изучение описанного в этом пункте формализма вместе с концепцией эффективного пересчета используемого варианта градиента. Не расчета, а именно пересчета, т.е. расчета с учетом результатов вычислений, проделанных на предыдущих итерациях. Поясним последнее примером. При безусловной минимизации квадратичной формы с помощью (неускоренного) покомпонентного спуска новая точка итерационного процесса отличается от старой только в одной компоненте, поэтому можно так организовать вычисления, чтобы итерация в среднем занимала $O(s)$ арифметических операций ($s$ – среднее число ненулевых элементов в столбце матрицы квадратичной формы), независимо от размеров матрицы квадратичной формы. Развивая идеи Ю.Е. Нестерова [273], в рамках описанного выше формализма в диссертации изучены способы наилучшего учета разреженности задачи.

Ранее уже отмечалась проблема огромной размерности пространства потоков по путям, в котором ставятся задачи выпуклой оптимизации, возникающие при поиске равновесий в транспортных сетях. С помощью метода условного градиента эта проблема решается благодаря эффективно вычислимому с помощью алгоритма Дейкстры поиска кратчайших путей в графе линейному минимизационному оракулу [173]. В диссертации выбран иной формализм [275], более удобный для перенесения на многостадийные модели, связанный с переходом к двойственной задаче, и ее решением прямо-двойственным мето-дом, позволяющим практически “бесплатно” восстанавливать решение прямой задачи. Отметим, что оба базисных метода (МЗС и БГМ) – прямо-двойственные. Этот формализм распространяется в диссертации на многостадийные модели. Именно с помощью переходов к двойственным задачам в части блоков удалось свести поиск равновесия в многостадийной модели к эффективно решаемой (с помощью описанных выше конструкций) задаче выпуклой оптимизации. Однако все это потребовало серьезного погружения в прямо-двойственность численных методов выпуклой оптимизации, особенно для задач выпуклой оптимизации на неограниченных областях. В частности, потребовалось изучение сочетания в одной задаче на одном методе основных (базисных) прямо-двойственных конструкций (способов восстановления решения сопряженной задачи), используемых ранее только по отдельности [261,266, 269, 281].

**Теоретическая и практическая значимость работы**

Настоящая диссертация мотивирована транспортными приложениями. Выбор задач был обусловлен общением со специалистами из НИиПИ Генплана г. Москвы, Департамента транспорта г. Москвы, компании А+С (основным дистрибьютором в России пакетов транспортного моделирования линейки PTV), ЦИТИ г. Москвы, Института экономики транспорта и транспортной политики НИУ ВШЭ.

За последние 10 лет резко возрос объем и качество доступных для моделирования транспортных данных (GPS-треки, данные сотовых операторов, данные видеокамер, всевозможные опросы населения, данные геоинформационных систем). При этом задачи поиска равновесий в транспортных сетях “вышли на передний план”, поскольку появилась возможность ставить такие задачи на больших (детализированных) транспортных сетях. В результате сложность таких задач резко возросла. С другой стороны, увеличилась и потребность в многократном решении таких задач с целью просмотра различных сценариев развития транспортной инфраструктуры. Как следствие, появилась необходимость в разработке нового (адаптированного под эти реалии) аппарата моделирования и соответствующих вычислительных процедур. Именно этому – поиску равновесий в больших (реальных) транспортных сетях, прежде всего, и посвящена диссертационная работа.

Большой акцент в диссертации сделан на теоретическое исследование оптимальных вычислительных процедур. Основные подходы здесь были заложены в Советском Союзе в работах Б.Т. Поляка, А.С. Немировского, Ю.Е. Нестерова и др. [91, 99]. В диссертации удалось “овыпуклить” многие классические результаты и посмотреть на многообразие этих результатов с единых позиций: с помощью базисного набора методов и введенных операций над ними удалось получить более простым способом как известные результаты, так и новые. В частности, исследовать вопросы: о равномерной ограниченности последовательностей, генерируемых численными методами (в том числе рандомизированными); о практической реализации всех рассматриваемых методов (с теоретическим обоснованием) в условиях отсутствия априорной информации о свойствах функционала задачи и свойствах решения (в частности, ограниченности нормы решения известным числом); о практически эффективных критериях останова рассматриваемых методов (с установленными теоретическими гарантиями скорости сходимости, что особенно нетривиально для прямо-двойственных методов на неограниченных областях).

Отметим также, что подавляющая часть изложения в диссертации ведется на современном уровне [93, 162, 259] – с точными константами в оценках числа итераций и с оценками вероятностей больших уклонений для стохастических (рандомизированных) методов.

Алгоритмы поиска равновесий в транспортных сетях, разработанные в диссертации, вошли в комплекс программ, созданных коллективом при участии автора диссертации. Комплекс был успешно принят в рамках отчета по гранту федеральной целевой программы «Исследования и разработки по приоритетным направлениям развития научно-технологического комплекса России на 2014 – 2020 годы», Соглашение № 14.604.21.0052 от 30.06.2014 г. с Минобрнаукой (идентификатор проекта RFMEFI60414X0052).

Также на основе разработок главы 3 диссертации, в том числе автором диссертации, был создан комплекс программ, проданный компании Huawei Russia (Договор Huawei с МФТИ №: YB2014120038 от 23 декабря 2014 года).

**Положения, выносимые на защиту (основные положения 1 – 4)**

1. Предложен новый эволюционный вывод энтропийной модели расчета матрицы корреспонденций с помощью марковской logit-динамики и транспортных потенциалов Канторовича–Гавурина. Таким образом, показано, что на модель расчета матрицы корреспонденций можно смотреть как на (энтропийно-регуляризованную) разновидность модели Бэкмана.
2. Получен оригинальный вывод модели Стабильной Динамики из модели Бэкмана с помощью вырождения функций затрат на прохождения ребер методом внутренних штрафов.
3. Используя то, что все блоки многостадийной модели транспортных потоков – есть вариации модели Бэкмана, являющейся, в свою очередь, популяционной игрой загрузки, предложен общий способ формирования вариационных принципов для поиска равновесий в многостадийных транспортных моделях. Полученные результаты распространены на общие иерархические популяционные игры загрузки. Исследован класс сетевых рынков (частным случаем которых является модель грузоперевозок РЖД) поиск конкурентных равновесий Вальраса в которых может быть осуществлен в рамках описанного выше (вариационного) формализма с заменой задачи выпуклой оптимизации на задачу поиска седловой точки с правильной выпукло-вогнутой структурой.

4. Впервые предложена вариация универсального быстрого градиентного метода Ю.Е. Нестерова (БГМ), самонастраивающегося на гладкость задачи для сильно выпуклых задач композитной оптимизации. Предложенную вариацию удалось распространить и на задачи стохастической оптимизации. Ранее считалось, что такая вариация либо невозможна, либо будет весьма сложной – полученные результаты опровергли эти опасения. Получены оценки скорости сходимости, показывающие, что предложенный метод является равномерно оптимальным (по числу итераций, с точностью до числового множителя) для общего класса задач выпуклой оптимизации. Описываемая линейка универсальных методов активно использовалась в диссертации при разработке комплекса программ для расчета различных блоков многостадийной транспортной модели.
5. Для поиска равновесия в модели Стабильной Динамики предложена специальная рандомизированная версия прямо-двойственного метода зеркального спуска (МЗС). Впервые получены оценки скорости сходимости в терминах вероятностей больших уклонений общего стохастического варианта МЗС для случая, когда оптимизация проводится на неограниченном множестве (как в модели Стабильной Динамики). Также впервые получены оценки вероятностной локализации итерационной последовательности. Используемая при этом оригинальная техника априорной формулировки гипотезы о характере хвостов распределений с апостериорной проверкой представляется полезной и для ряда других приложений.
6. Предложены два различных прямо-двойственных подхода для решения задачи энтропийно-линейного программирования, в частности, для задачи расчета матрицы корреспонденций. В основе обоих подходов лежит идея решения двойственной задачи с помощью БГМ. Оба подхода распространены на общие задачи минимизации сильно выпуклых функционалов простой (например, сепарабельной) структуры при аффинных ограничениях.

7. Исследованы общие способы приближенного восстановления решения сопряженной (двойственной) задачи по последовательности, генерируемой методом, решающим исходную задачу. В частности, было продемонстрировано сочетание различных способов восстановления на одной задаче – поиск равновесного распределения потоков по ребрам в Смешанной Модели, когда часть ребер – из модели Бэкмана, оставшаяся часть – из модели Стабильной Динамики.
8. Исследована роль неточностей неслучайной природы при расчете градиента и значения оптимизируемой функции, а также неточностей, возникающих при проектировании, на итоговые оценки скорости сходимости различных методов. Полученные здесь результаты используются для обоснования концепции суперпозиции (правильного чередования – чередования с правильными частотами) численных методов.
9. Разработана концепция суперпозиции численных методов выпуклой оптимизации. Рассмотрены конкретные примеры приложений (min max-задачи, min min-задачи), в частности, задачи поиска равновесия в многостадийных моделях транспортных потоков.
10. Исследованы покомпонентные методы, спуски по направлению и методы нулевого порядка, с точки зрения наследования “хороших” свойств своих полноградиентных аналогов. В частности, была исследована прямо-двойственность покомпонентных методов. Новые результаты удалось получить для спусков по направлению и их дискретных аналогов (методов нулевого порядка) в части сочетания структуры множества, на котором происходит оптимизация, и способа выбора случайного направления спуска. В зависимости от контекста, оказалось, что оптимально выбирать случайное направление либо среди координатных осей, либо равномерно на евклидовой сфере. Исследована возможность использования прямо-двойственных блочно-покомпонентных методов для поиска равновесия в модели Стабильной Динамики.
11. На примере задачи ранжирования web-страниц (Google problem) исследованы основные конструкции huge-scale-оптимизации. В частности, предложено учитывать разреженность задачи и использовать рандомизированные методы. Предложенные методы позволили не только эффективно находить вектор PageRank, но и решать целый ряд других задач huge-scale-оптимизации. Например, предложен эффективный численный метод решения равномерно разреженных (по строкам и столбцам) систем линейных уравнений огромных размеров. Доказана высокая эффективность метода, если решение состоит из разных по масштабу компонент (1-норма вектора решения и его 2-норма достаточно близки).

12. Алгоритм Григориадиса–Хачияна поиска равновесий в антагонистических матричных играх удалось проинтерпретировать как специальным образом рандомизированный МЗС, что позволило получить для него оценки вероятностей больших отклонений и распространить его на задачи с равномерно разреженной матрицей, показав, что при этом разреженность учитывается оптимальным образом. Также получена естественная онлайн-интерпретация этого алгоритма.

**Степень достоверности, публикации и апробация результатов**

Степень достоверности и обоснованности полученных результатов достаточно высоки. Они обеспечиваются строгостью и корректным использованием математических доказательств, подтверждением результатов работы экспертными оценками специалистов. Результаты, включенные в данную работу, представлены более чем в 40 работах, в числе которых одна монография, одно учебное пособие и 26 статей в журналах из перечня ВАК. В работах, написанных в соавторстве, соискателю принадлежат результаты, указанные в положениях, выносимых на защиту. По результатам диссертации было прочитано несколько оригинальных курсов лекций студентам МФТИ, НМУ, БФУ им. И. Канта. Большая часть этих лекций (более 50) доступна в сети Интернет: http://www.mathnet.ru/ и https://www.youtube.com/user/PreMoLab.

Результаты, включенные в диссертацию, неоднократно докладывались в 2010–2016 гг. на научных семинарах в ведущих университетах и институтах России: семинаре отдела Математического моделирования экономических систем ВЦ ФИЦ ИУ РАН (рук. чл.-корр. РАН И.Г. Поспелов), семинаре лаборатории Адаптивных и робастных систем им. Я.З. Цыпкина ИПУ РАН (рук. проф. Б.Т. Поляк), семинаре лаборатории Структурных методов анализа данных в предсказательном моделировании МФТИ–ИППИ (рук. проф. В.Г. Спокойный), семинаре лаборатории Больших случайных систем мехмата МГУ (рук. проф. В.А. Малышев), семинаре кафедры Исследования операций ВМиК МГУ (рук. проф. А.А. Васин), семинаре ЦЭМИ РАН (рук. проф. Л.А. Бекларян), семинаре ИАП РАН (рук. чл.-корр. РАН А.С. Холодов), семинаре ПОМИ РАН (рук. проф. А.М. Вершик), семинарах федерального профессора А.М. Райгородского в МФТИ и МГУ, семинарах Яндекса (Москва, Новосибирск), семинаре кафедры Высшей математики МФТИ (рук. проф. Е.С. Половинкин), семинаре кафедры Информатики МФТИ (рук. чл.-корр. РАН И.Б. Петров), семинаре НИиПИ Генплана г. Москвы, семинаре по Транспортному моделированию ИПМ РАН (рук. доц. В.П. Осипов, акад. Б.Н. Четверушкин), семинаре МАДИ–МТУСИ–МИАН (рук. проф. А.П. Буслаев, проф. М.В. Яшина, акад. В.В. Козлов), научных семинарах

ДФВУ г. Владивосток (рук. проф. Е.А. Нурминский), УФУ г. Екатеринбург (рук. проф. Н.Н. Субботина), БФУ г. Калининград (рук. проф. С.В. Мациевский), Иннополиса г. Иннополис (рук. доц. Я.А. Холодов). Результаты докладывались также на научном семинаре института К. Вейерштрасса в Берлине (рук. проф. В.Г. Спокойный).

Результаты диссертации докладывались на многочисленных международных конференциях и школах: Traffic and granular flows, Moscow 2011; OPTIMA (Optimization and applications) 2012–2015; ММРО (Математические методы распознавания образов) 2013, 2015; ИОИ (Интеллектуализация обработки информации) 2012, 2014; ИТиС (Информационные технологии и системы) 2013– 2015; ISMP (International Symposium on Mathematical Programming) 2012, 2015; International Conference Network Analysis (NET2015), Нижний Новгород, 2015; Летняя школа “Современная математика”, Дубна 2011–2014, 2016; Сириус-2016; Традиционная математическая школа “Управление, информация и оптимизация” Б.Т. Поляка 2012–2016; Summer School on Operation Research and Application, Нижний Новгород, 2015; Байкальские чтения Алексея Савватеева, 2014–2016; Skoltech Deep Machine Intelligence workshop, 2016. На конференции “Moscow International Conference on Operation Research (ORM)”, 2013, соискатель руководил работой транспортной секции.

По тематике диссертации соискателем организовано несколько научных семинаров: “Математическое моделирование транспортных потоков” (совместный семинар МФТИ–НМУ), “Стохастический анализ в задачах” (совместный семинар МФТИ–НМУ), Математический кружок (МФТИ), на которых в 2012–2016 гг. докладывались результаты работы.

**Благодарности**

Настоящая работа была инициирована в конце 2007 года проф. А.А. Шананиным, предложившим соискателю заняться математическим моделированием транспортных потоков и параллельно начать читать студентам МФТИ одноименный курс (такой курс читается с 2008 года). Полезным при написании работы оказался фундамент численных методов выпуклой оптимизации, заложенный во время обучения на базовой кафедре ФУПМ МФТИ в ВЦ РАН. Важным этапом для научного роста соискателя стало создание в 2011 г. на ФУПМ МФТИ лаборатории Структурных методов анализ данных в предсказательном моделировании (рук. проф. В.Г. Спокойный). Благодаря лаборатории и лично проф. В.Г. Спокойному и акад. А.П. Кулешову у соискателя появились новые возможности и новые научные контакты, которые и определили окончательный выбор темы диссертации (основные результаты этой работы были получены в МФТИ и ИППИ РАН в период апрель 2012 года – апрель 2016 года). Принципиально важную роль сыграло научное общение с проф. Ю.Е. Нестеровым. Хотелось бы поблагодарить проф. А.Б. Юдицкого, приглашавшего соискателя в 2015, 2016 гг. на зимние школы по оптимизации в Лез–Уше (Франция). Участие в этих школах существенно упростило знакомство с состоянием мировой науки в области численных методов выпуклой оптимизации. Для соискателя при работе над главой 1 важной оказалась поддержка со стороны ИПМ РАН, особенно, доц. В.П. Осипова и акад. Б.Н. Четверушкина и со стороны ВШЭ, особенно, проф. М.Я. Блинкина. Также хотелось бы поблагодарить своих коллег (соавторов) и учеников.

Написание диссертации в 2015–2016 гг. было финансово поддержано грантами РФФИ: 15-31-70001 мол_а_мос (Глава 1), 15-31-20571-мол_а_вед (Главы 2, 3), 14-01-00722-а (глава 4); и грантом Российского научного фонда, проект № 14-50-00150 (Главы 5, 6).

**Структура и объем работы**

Диссертация состоит из введения, шести глав, заключения, одного приложения и списка литературы. Полный объем диссертации составляет 487 страниц, список литературы содержит 331 наименование.

# Глава 1 Эволюционный вывод равновесных моделей распределения потоков в больших транспортных сетях и численные методы поиска равновесий в таких моделях

## 1.1 Трехстадийная модель равновесного распределения транспортных потоков

### 1.1.1 Введение

Одной из основных задач последнего времени, остро стоящих в Москве и ряде других крупных городов России (Санкт-Петербург, Пермь, Владивосток, Иркутск, Калининград и др.) является разработка транспортной модели города, позволяющей решать задачи долгосрочного планирования (развития) транспортной инфраструктуры города [287]. В частности, ожидается, что разработка такой модели поможет ответить на вопросы: какой из проектов дорожного строительства оптимален, где пропускная способность дороги недостаточна, как изменится транспортная ситуация, если построить в этом месте торговый центр (жилой район, стадион), как правильно определять маршруты и расписание движения общественного транспорта, какой эффект даст выделение полос для общественного транспорта и т.п. Уже имеется программное обеспечение, позволяющее частично решать указанные выше задачи. Однако имеется много вопросов к тому, какие модели и алгоритмы используются в большинстве программных продуктах. Например, не очевидным элементом почти всех этих продуктов является использование в качестве одного из блоков модели равновесного распределения транспортных потоков Бэкмана (1955) [49, 112, 147, 289, 310]. Эта во многом хорошая модель, тем не менее, имеет довольно много недостатков (см. [275]). Например, калибровка такой модели требует знания функций затрат на ребрах графа транспортной сети (эти функции связывают время в пути по ребру с величиной транспортного потока по этому ребру), причем сама модель оказывается довольно чувствительной к выбору этих функций, которые в модели Бэкмана, как правило, предполагаются выпуклыми, монотонно возрастающими. Достаточно сказать, что в случае наличия платных дорог, для расчета оптимальных плат за проезд требуется вычислять, например, производные этих неизвестных функций [301, 302]. Существование таких функций в модели Бэкмана является одним из основных предположений, и, одновременно, одним из самых слабых мест. Реальные данные показывают (см. Рисунок 1.1.1, полученный В.А. Данилкиным по данным ЦОДД в 2012 году), что предположение о классе функций затрат не выполняется.

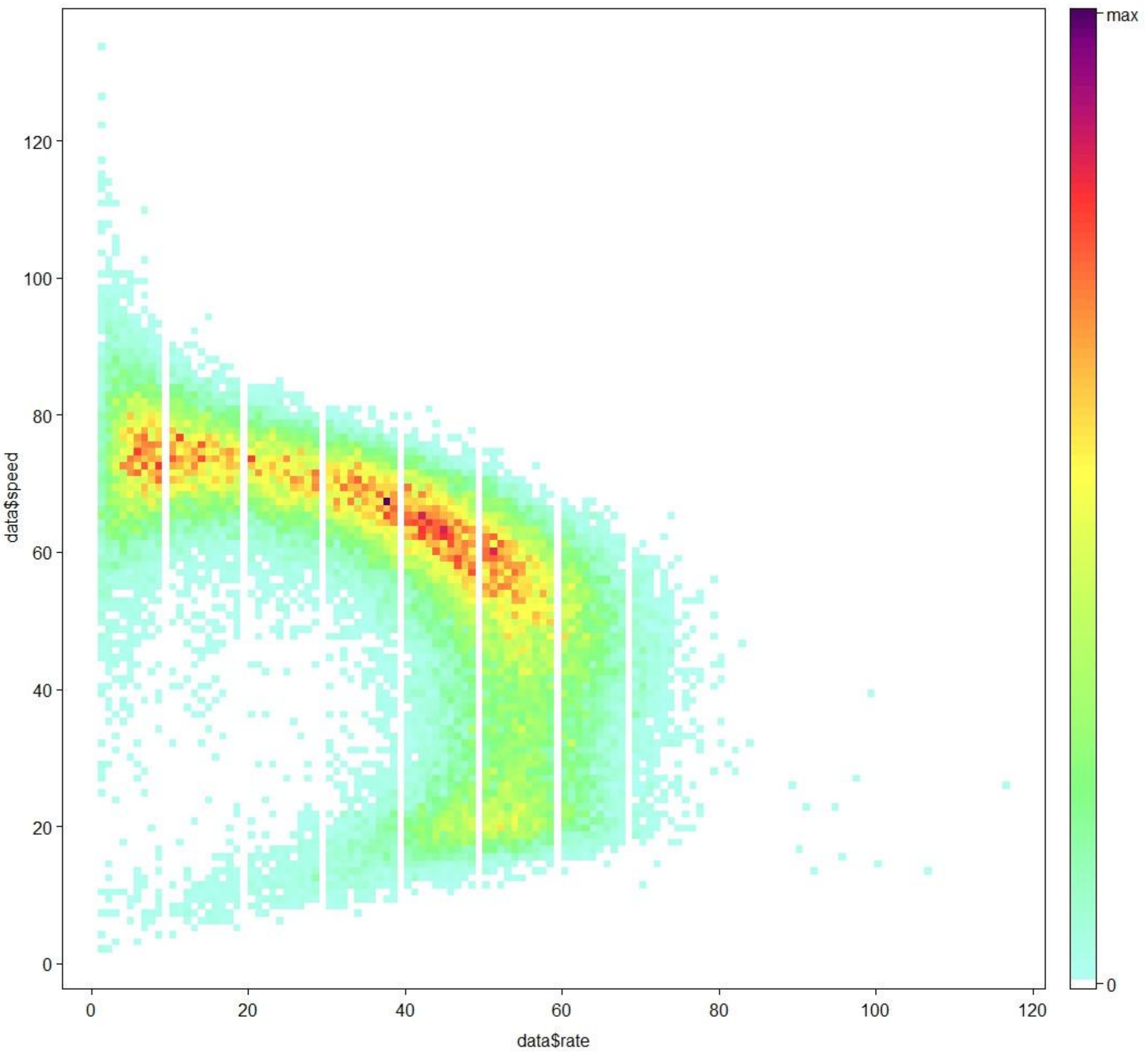

Рисунок 1.1.1. По оси абсцисс – поток по двум полосам (авт/мин), а по оси ординат – скорость (км/час)

Но даже если предположить, что такая зависимость все же существует,[1] то по-прежнему остается другая проблема: как калибровать модель Бэкмана, то есть откуда брать эти зависимости. Не получится ли переобучения у создаваемой нами модели? То есть, не получится ли так, что распоряжаясь большим произволом при калибровке по обучающей выборке (историческим данным) мы “переподгоним” модель: исторические данные за счет большого числа подкручиваемых параметров мы действительно можем хорошо научиться описывать, но использовать такую модель для планирования будет опасно, поскольку не будет контроля переобучения. Обычным средством борьбы с переобучением в этом месте является параметризация функций затрат, например, в классе BPR функций [49, 112, 147, 287, 289, 310]. К сожалению, какого бы то ни было научного обоснования, почему именно такая параметризация используется, нам не известно.

[1] Это допущение можно оправдать, например, тем, что, как правило, мы рассматриваем равновесные конфигурации с точки зрения пользователей сети в модели Бэкмана, которым соответствует только одна из веток – верхняя (рис. 1.1.1). В модели Бэкмана пользователи сети при принятии решений оценивают время в пути в зависимости от величины потока, которая интерпретируется, как число желающих проехать по этому ребру в единицу времени. Нижняя ветка, отвечающая приблизительно линейному росту скорости с ростом потока $V \approx f/\rho_{\max}$ , соответствует ситуации, когда есть узкое место, пропускная способность которого по каким-то причинам определяется не типичными характеристиками рассматриваемого ребра, а скажем, пробкой, пришедшей с впереди идущего ребра (по ходу движения). И в таких ситуациях величина потока $f$ интерпретируется не как число желающих, а как число могущих проехать. Такие ситуации просто исключаются в модели Бэкмана.

Другим, неочевидным элементом этих программных продуктов являются используемые вычислительные алгоритмы: контроль их робастности [151] к неточности (неполноте) данных, ошибкам округления (поскольку возникают задачи огромных размеров, то такие ошибки могут интенсивно накапливаться). Наконец, сама философия, использующаяся в таких продуктах при построении равновесной модели города [287], также вызывает много вопросов, о которых немного подробнее будет написано в подразделе 1.1.5. Несмотря на отмеченные выше проблемы, разработчики программного обеспечения часто находят вполне разумные инженерные компромиссы (нам известно о положительном опыте PTV VISSUM и TRANSNET), неплохо работающие практике.

Целью данного раздела является предложить математическую трехстадийную транспортную модель города, в которой один из блоков (модель равновесного распределения потоков) предлагается заменить с модели Бэкмана на модель стабильной динамики [272, 275]. К сожалению, целый ряд проблем, свойственных ранее известным многостадийным моделям будет присущ и модели, предлагаемой в данном разделе. Однако несколько важных недостатков, по-видимому, удалось устранить. Прежде всего, речь идет о возможности калибровки модели по реальным данным, контроле переобучения и существовании эффективного робастного вычислительного алгоритма с гарантированными (не улучшаемыми) оценками числа затраченных арифметических операций для достижения требуемой точности. Последнее обстоятельство представляется особенно важным в контексте того, как обычно используются такие модели. А именно, с помощью таких моделей просматривается множество различных сценариев. К сожалению, оптимизационные задачи вида: "где и какую дорогу стоит построить при заданных бюджетных ограничениях", а также многие другие задачи, решаются непосредственным перебором различных вариантов, где для расчета каждого варианта потребуется запускать модель, меняя каждый раз что-то на входе. Кстати сказать, предложенный для расчета модели алгоритм позволяет также учитывать масштаб изменения входных данных. Если эти изменения небольшие (точнее говоря, меняется небольшое количество входных параметров), то для выполнения перерасчета по модели потребуется значительно меньше времени, чем при первом запуске.

### 1.1.2 Структура раздела и предварительные сведения

Опишем вкратце структуру раздела. В подразделе 1.1.3 описывается эволюционный способ вывода популярного на практике статического способа расчета матрицы корреспонденций (энтропийной модели). Также приводится основная идея, базирующаяся на теореме Тихонова о разделении времен [7], получения трехстадийной модели из отдельных блоков (расчет матрицы корреспонденций + равновесное расщепление потоков + равновесное распределение потоков). Точнее говоря, использование теоремы Тихонова – лишь часть идеи, которая сводит решение задачи к поиску (единственного) притягивающего положения равновесия системы в медленном времени, отвечающей за формирование корреспонденций, при подстановке в неё зависимостей времен в пути от корреспонденций (такие зависимости получаются из системы в быстром времени). Другая её часть, заключается в том, что задачи поиска (единственного) притягивающего положения равновесия системы в медленном времени и поиска зависимостей времен в пути от корреспонденций с помощью некоторых вариационных принципов, о которых говорится в подразделах 1.1.3 и 1.1.4, 1.1.6–1.1.8, сводятся к задачам выпуклой оптимизации, которые можно объединить в одну общую задачу поиска седловой точки негладкой выпукло-вогнутой функции. В конце пункта указывается возможность обобщения трехстадийной модели до четырехстадийной (есть некоторые нюансы в трактовках – стоит иметь в виду, что в разных литературных источниках к таким моделям могут предъявляться немного разные требования), в которой учитываются различные типы пользователей и различные типы передвижений (моделирование идет на больших масштабах времени). В подразделе 1.1.4 описывается модель равновесного распределения потоков Бэкмана. В конце пункта приводится эволюционный способ интерпретации возникающего в этой модели равновесия Нэша–Вардропа. В подразделе 1.1.5 приводится краткий обзор многостадийных моделей, построенных на основе моделей описанных в подразделах 1.1.3 и 1.1.4. В подразделе 1.1.6 описывается модель равновесного распределения потоков, которую мы далее будем называть моделью стабильной динамики. Модель стабильной динамики требует намного меньше данных для своей калибровки, наследует практически все основные “хорошие” свойства модели Бэкмана, и не наследует ряд недостатков. В подразделе 1.1.7 модель стабильной динамики выводится с помощью предельного перехода из модели Бэкмана. В подразделе 1.1.8 строится обобщение модели стабильной динамики на случай, когда есть несколько способов передвижения (личный транспорт и общественный; отметим при этом, что в Москве и области более 70% пользователей сети используют общественный транспорт). Таким образом, в подразделе 1.1.8 в модель стабильной динамики органично встраивается модель равновесного расщепления потоков. В подразделе 1.1.9 энтропийная модель расчета матрицы корреспонденций из подраздела 1.1.3 объединяется с моделью из подраздела 1.1.8. В результате получается трехстадийная модель, в которой учитывается и формирование корреспонденций, и расщепление потоков, и равновесное распределение потоков по графу транспортной сети. Примечательно, что поиск равновесия в полученной трехстадийной модели (из задачи поиска седловой точки негладкой выпукло-вогнутой функции – см. подраздел 1.1.3) в итоге сводится к решению задачи негладкой выпуклой оптимизации, с ограниченной константой Липшица функционала, но с неконтролируемой начальной невязкой. Важно отметить, что помимо самого решения, нужно определять и часть двойственных переменных, имеющих содержательный физический смысл. В подразделе 1.1.10 модель подраздела 1.1.9 обобщается на случай поиска стохастического равновесия, что можно проинтерпретировать как ограниченную рациональность водителей или их неполную информированность. Это допущение делает модель подразделе 1.1.9 более приближенной к практике. С вычислительной точки зрения, сделанная модификация сводит задачу к задаче гладкой выпуклой оптимизации с немного более громоздким функционалом. Полученная задача во многом наследует все вычислительные минусы и плюсы негладкого случая. Обе задачи выпуклой оптимизации из подразделов 1.1.9, 1.1.10 требуют разработки адекватных прямо-двойственных субградиентных алгоритмов решения. Заметим, что использовать методы с оракулом (сообщающим, в зависимости от своего порядка, значения функций в выбранных точках, их градиенты, и т.д.) порядка выше первого [93] не представляется возможным в виду размеров задач. В заключительном подразделе 1.1.11 приводятся (с объяснениями) несколько практических рецептов по калибровке предложенных моделей.

На протяжении всего этого раздела (и всех последующих разделов этой диссертации) мы будем активно использовать элементы выпуклого анализа и методы численного решения задач выпуклой оптимизации, ориентируясь на априорное знакомство читателя с этими дисциплинами, например, в объеме книг [85, 93]. Одним из основных инструментов этого раздела будет теорема фон Неймана о минимаксе для выпукло-вогнутых функций [119]. Причем использоваться эта теорема будет не только для функций, заданных на произведении компактов, но и на неограниченных множествах. Надо лишь иметь гарантию, что максимумы и минимумы существуют (достигаются). Проблема сводится к существованию неподвижной точки у многозначного отображения. В этой области имеется большое количество результатов, с запасом покрывающих потребности данного раздела, для обоснования возможности перемены порядка взятия максимума и минимума. В частности, в этом разделе мы будем пользоваться вариантом минимаксной теоремы, называемой в зарубежной литературе “Sion’s minimax theorem” [313], в которой предполагается компактность лишь одного из множеств, отвечающих выпуклым или вогнутым переменным, и непрерывность функции. Тем не менее, далее мы будем использовать более привычное название “минимаксная теорема фон Неймана”.

Все необходимые обозначения и ссылки будут вводиться в разделе по мере необходимости.

### 1.1.3 Энтропийная модель расчета матрицы корреспонденций

Приведем, во многом следуя книге [49] (см. также раздел 1.2 этой главы 1), обоснование (ориентированное, скорее, на западный уклад жизни, поскольку предполагается, что человек за жизнь несколько раз меняет место жительство и место работы), пожалуй, одному из самых популярных способов расчета матрицы корреспонденций, имеющему более чем сорокалетнюю историю, – энтропийной модели [26, 105, 123].

Пусть в некотором городе имеется $n$ районов, $L_i > 0$ – число жителей $i$-го района, $W_j > 0$ – число работающих в $j$-м районе. При этом $N = \sum_{i=1}^{n} L_i = \sum_{j=1}^{n} W_j$, – общее число жителей города. В последующих пунктах, под $L_i \geq 0$ будет пониматься число жителей района, выезжающих в типичный день за рассматриваемый промежуток времени из $i$-го района, а под $W_j \geq 0$ – число жителей города, приезжающих на работу в $j$-й район в типичный день за рассматриваемый промежуток времени. Обычно, так введенные, $L_i$, $W_j$ рассчитываются через число жителей $i$-го района и число работающих в $j$-м районе с помощью более менее универсальных (в частности, не зависящих от $i$, $j$) коэффициентов пропорциональности. Эти величины являются входными параметрами модели, т.е. они не моделируются (во всяком случае, в рамках выбранного подхода). Для долгосрочных расчетов с разрабатываемой моделью требуется иметь прогноз изменения значений этих величин.

Обозначим через $d_{ij}(t) \geq 0$ – число жителей, живущих в $i$-м районе и работающих в $j$-м в момент времени $t$. Со временем жители могут только меняться квартирами, поэтому во все моменты времени $t \geq 0$

$$d_{ij}(t)\ge 0,\ \sum_{j=1}^{n} d_{ij}(t)\equiv L_i,\ \sum_{i=1}^{n} d_{ij}(t)\equiv W_j,\ i,j=1,...,n,\ 1\ll n^2\ll N.^2 \qquad \textbf{(A)}$$

Опишем основной стимул к обмену: работать далеко от дома плохо из-за транспортных издержек. Будем считать, что эффективной функцией затрат [49] будет $R(T)=\beta T/2$, где $T>0$ – время в пути от дома до работы, а $\beta>0$ – настраиваемый параметр модели (который также можно проинтерпретировать и даже оценить).

Теперь опишем саму динамику (детали см. в [36, 49]). Пусть в момент времени $t\ge 0$ $r$-й житель живет в $k$-м районе и работает в $m$-м, а $s$-й житель живет в $p$-м районе и работает в $q$-м. Тогда $\lambda_{k,m;\,p,q}(t)\Delta t+o(\Delta t)$ – есть вероятность того, что жители с номерами $r$ и $s$ ( $1\le r<s\le N$ ) “поменяются” квартирами в промежутке времени $(t,t+\Delta t)$. Вероятность обмена местами жительства зависит только от мест проживания и работы обменивающихся:

$$\lambda_{k,m;\,p,q}(t)\equiv\lambda_{k,m;\,p,q}=\lambda N^{-1}\exp\Big(\underbrace{R(T_{km})+R(T_{pq})}_{\text{суммарные затраты до обмена}}-\underbrace{\big(R(T_{pm})+R(T_{kq})\big)}_{\text{суммарные затраты после обмена}}\Big)>0,$$

где коэффициент $\lambda>0$ характеризует интенсивность обменов. Совершенно аналогичным образом можно было рассматривать случай “обмена местами работы”. Здесь стоит оговориться, что “обмены” не стоит понимать буквально – это лишь одна из возможных интерпретаций. Фактически используется так называемое “приближение среднего поля” [32],

---

[2] Для Москвы (и других крупных мегаполисов) часто выбирают $n\sim 10^2-10^3$. Следовательно, корреспонденций будет $n^2\sim 10^4-10^6$, и чтобы каждую из корреспонденций определить с точностью 10% (относительной точностью $\varepsilon\sim 10^{-1}$) потребуется (это оценка снизу) опросить не менее $n^2/\varepsilon^2\sim 10^6-10^8$ жителей города, что не представляется возможным – это обстоятельство является одной из причин (не единственной), почему стараются снизить размерность пространства параметров, считая, что матрица корреспонденций задается не $n^2$ параметрами, а только, например, $2n$ ($L$ и $W$), которые и надо определять (см. также подраздел 1.1.5, в котором параметров $2n+1$). На самом деле выбирать критерием качества “относительную точность” для каждой (в том числе очень маленькой) корреспонденции не очень разумно. Более естественно восстанавливать матрицу (вектор) корреспонденций в 1-норме. Другой причиной использования описанной ниже энтропийной модели является возможность прогнозирования с её помощью того, как будет меняться матрица корреспонденций при изменении инфраструктуры города, собственно эта одна из тех задач, которые необходимо уметь решать, для получения ответов на вопросы, приведенные в подразделе 1.1.1.

Кроме того, важно заметить, что крупный Мегаполис, как правило, представляется в такого рода моделях вместе со всеми своими окрестными территориями. Скажем, для Москвы – это Московская область. Отметим также, что для Москвы и области $N\sim 10^7$.

т.е. некое равноправие агентов (жителей) внутри фиксированной корреспонденции и их независимость.[3]

Согласно эргодической теореме для марковских цепей (в независимости от начальной конфигурации $\{d_{ij}(0)\}_{i=1,\, j=1}^{n,\, n}$) [16, 49] предельное распределение совпадает со стационарным (инвариантным), которое можно посчитать (получается проекция прямого произведение распределений Пуассона на гипергрань некоторого многогранника):

$$\exists\ c_n > 0:\ \forall\ \{d_{ij}\}_{i=1,\, j=1}^{n,\, n} \in (\mathrm{A}),\ t \ge c_n N \ln N$$

$$P\left(d_{ij}(t) = d_{ij}, i, j = 1,\dots,n\right) \simeq Z^{-1} \prod_{i,j=1}^{n} \exp\left(-2R\left(T_{ij}\right) d_{ij}\right) \cdot \left(d_{ij}!\right)^{-1} \stackrel{def}{=} p\left(\{d_{ij}\}_{i=1,\, j=1}^{n,\, n}\right),$$

где “статсумма” $Z$ находится из условия нормировки получившейся “пуассоновской” вероятностной меры. Отметим, что стационарное распределение $p\left(\{d_{ij}\}_{i=1,\, j=1}^{n,\, n}\right)$ удовлетворяет условию детального равновесия:

$$\left(d_{km}+1\right)\left(d_{pq}+1\right) p\left(\left\{d_{11},\dots,d_{km}+1,\dots,d_{pq}+1,\dots,d_{pm}-1,\dots,d_{kq}-1,\dots,d_{nn}\right\}\right) \lambda_{k,m;\,p,q} =$$

$$= d_{pm} d_{kq} p\left(\{d_{ij}\}_{i=1,\, j=1}^{n,\, n}\right) \lambda_{p,m;\,k,q}.$$

При $N \gg 1$ распределение $p\left(\{d_{ij}\}_{i=1,\, j=1}^{n,\, n}\right)$ экспоненциально сконцентрировано на множестве (**A**) в $\mathrm{O}\left(\sqrt{N}\right)$ окрестности наиболее вероятного значения $\{d_{ij}^*\}_{i=1,\, j=1}^{n,\, n}$, которое определяется как решение задачи энтропийно-линейного программирования [49] (подробнее о таких задачах и о том, как их решать см. [6, 26, 32, 37, 55, 105, 123, 207, 221, 230], см. International workshops on Bayesian inference and maximum entropy methods in science and engineering (AIP Conf. Proceedings, проводится каждый год с 1980 г.), а также в разде-лах 3.1, 3.2 главы 3 – при небольших $\beta$ успешно работает метод балансировки [198], при больших $\beta$ доминирует подход из [6, 138]):[4]

---

[3] Отметим также, что конечной цели (получение задачи (1.1.1)) можно добиться отличными способами. Скажем, используя формализм Л.И. Розоноэра “систем обмена и распределения ресурсов” со структурной функцией энтропией (кстати, это можно не постулировать, структурная функция появляется, при работе с условием интегрируемости дифференциальной формы возможных обменов ресурсами) – см., например, [18] и цитированную там литературу. При таком подходе стохастика не нужна, и вариационный принцип (максимизации энтропии при аффинных ограничениях) получается мало чувствительным к особенностям возможных превращений в системе. Другие способы получения (1.1.1) связаны с информационно-статистическими соображениями, например, принципом максимума правдоподобия [26, 105, 304].

$$\ln p\left(\left\{d_{ij}\right\}_{i=1,\, j=1}^{n,\, n}\right) \sim -\sum_{i,j=1}^{n} d_{ij}\ln\left(d_{ij}/e\right) - \beta\sum_{i,j=1}^{n} d_{ij}T_{ij} \to \max_{\left\{d_{ij}\right\}_{i=1,\, j=1}^{n,\, n}\in(\mathrm{A})}. \qquad (1.1.1)$$

Естественно принимать решение этой задачи $\left\{d_{ij}^{*}\right\}_{i=1,\, j=1}^{n,\, n}$ за равновесную конфигурацию [10]. Однако имеется проблема: $T_{ij}$ – неизвестны, и зависят от $\left\{d_{ij}\right\}$. Эту проблему мы постараемся решить в дальнейшем.

Обратим внимание, что предложенный выше вывод энтропийной модели расчета матрицы корреспонденций отличается от классического [26]. В монографии А.Дж. Вильсона [26] $\beta$ интерпретируется как множитель Лагранжа к ограничению на среднее “время в пути”: $\sum_{i,j=1}^{n} d_{ij}T_{ij} = C$.[5] Сделано это нами для того, чтобы контролировать знак параметра $\beta > 0$ и лучше понимать его физический смысл (нам это понадобится в дальнейшем).[6] Отметим также, что характерный временной масштаб формирования корреспонденций – годы (это не совсем так для корреспонденций типа дом–торговля, дом–отдых). В то время как характерное время установления равновесных значений $T_{ij}\left(d\right)$ – недели (см. ниже).

---

[4] При получении этой формулы использовалась асимптотическая формула Стирлинга [26, 49, 105], то есть предполагалось, что если $d_{ij} > 0$, то $d_{ij} \gg 1$ (если все $L_i > 0$, $W_j > 0$, то и все $d_{ij} > 0$ [49, 105]) и, как следствие, целочисленностью переменных $d_{ij}$ можно пренебречь, то есть решать не *NP*-полную задачу выпуклого целочисленного программирования, а обычную задачу выпуклой оптимизации (1.1.1). Сделанное предположение “если $d_{ij} > 0$, то $d_{ij} \gg 1$” во многом будет следовать из дальнейших рассуждений (см. также [49, 105]).

[5] При этом остальные ограничения имеют такой же вид, а функционал $F\left(d\right) = -\sum_{i,j=1}^{n} d_{ij}\ln\left(d_{ij}/e\right)$. Тогда, согласно экономической интерпретации двойственных множителей Л.В. Канторовича [49, 105]:

$$\beta(C) = \partial F\left(d\left(C\right)\right)/\partial C.$$

Из такой интерпретации иногда делают вывод о том, что $\beta$ можно понимать, как цену единицы времени в пути. Чем больше $C$, тем меньше $\beta$.

[6] Отметим, что также как и в [26] из принципа ле Шателье–Самуэльсона [18] следует, что с ростом $\beta$ среднее время в пути $\sum_{i,j=1}^{n} d_{ij}\left(\beta\right)T_{ij}\left(d\left(\beta\right)\right)$ будет убывать. В связи с этим обстоятельством, а также исходя из соображений размерности, вполне естественно понимать под $\beta$ величину, обратную к характерному (среднему) времени в пути [26] – физическая интерпретация. Собственно, такая интерпретация параметра $\beta$, как правило, и используется в многостадийных моделях (см., например, [287], а также подраздел 1.1.5).

Хочется сказать, что мы здесь находимся в условиях теоремы Тихонова о разделении времен [7], точнее говоря, что можно просто подставить в (1.1.1) зависимости $T_{ij}(d)$ и решать полученную задачу. Но теорема Тихонова должна применяться для системы ОДУ, представляющей собой, в данном случае, динамику квазисредних введенной выше стохастической динамики [32, 49]:[7]

$$\frac{d}{dt}c_{ij} = \sum_{k,p=1}^{n} \lambda \exp\left(\frac{\beta}{2}\left(\left[T_{ip}+T_{kj}\right]-\left[T_{ij}+T_{kp}\right]\right)\right)c_{ip}c_{kj} -$$

$$-\sum_{k,p=1}^{n} \lambda \exp\left(\frac{\beta}{2}\left(\left[T_{ij}+T_{kp}\right]-\left[T_{ip}+T_{kj}\right]\right)\right)c_{ij}c_{kp},\ c_{ij}=d_{ij}/N;$$

$$\varepsilon\frac{d}{dt}T_{ij} = \left[\text{сложный оператор, зависящий от } c\right],\ \varepsilon \sim 10^{-2}-10^{-3} \ll 1.$$

Для обоснования возможности применения здесь такого рода результата как теорема Тихонова требуется много усилий: во-первых, изначально введенные динамики стохастические (см. также подразделы 1.1.4, 1.1.7, 1.1.8 для $T_{ij}$), поэтому, в конечном итоге, нужно все обосновывать именно для них, во-вторых, нам интересна асимптотика по времени первой системы – в медленном времени (то есть нельзя ограничиться классическим случаем: ограниченного отрезка времени), в-третьих, сложный характер оператора, стоящего в правой части второй системы – в быстром времени, не позволяет явно его выписать. Тем не менее, процедура построения этого оператора, которая, в свою, очередь, предполагает переход к пределу (по $\mu \to 0+$ и по числу пользователей транспортной сети (аналог $N(\to\infty)$) – см. подраздел 1.1.7), может быть при необходимости получена из того, что будет далее приведено в подразделах 1.1.4, 1.1.7, 1.1.8. Мы не будем здесь подробно описывать, как можно бороться с указанными проблемами, заметим лишь, что в виду “хороших” свойств зависимостей $T_{ij}(c)$ (см. подразделы 1.1.7, 1.1.8), получающихся из приравнивания нулю левой части системы в быстром времени, функция[8]

$$H(c) = \sum_{i,j=1}^{n} c_{ij}\ln\left(c_{ij}\right) + \beta\Phi(c),\ \text{(следует сравнить с (1.1.1))}$$

---

[7] Точнее говоря, из теоремы Т. Куртца [195] следует (с некоторыми, довольно общими, оговорками относительно зависимостей $T_{ij}(d/N)$): если предположить, что делается такой предельный переход, что

$$\exists\ c_{ij}(0) = \lim_{N\to\infty} d_{ij}(0)/N,\ \text{то}\ \forall\ t\ge 0\ \exists\ c_{ij}(t) \overset{п.н.}{=} \lim_{N\to\infty} d_{ij}(t)/N$$

(в общем случае пределы существуют не равномерно по времени, но в нашем случае равномерно), где функции $c_{ij}(t)$ (не случайные) удовлетворяют выписанной системе ОДУ.

[8] Полученная функция будет выпуклой. Далее это будет проясняться.

где $\partial\Phi(c)/\partial c_{ij} = T_{ij}(c)$ (то, что такая функция $\Phi$ существует, показано ниже), будет функцией Ляпунова, и, одновременно, функцией Санова (действием), то есть функцией, которая с точностью до знака и аддитивной постоянной характеризует экспоненциальную концентрацию стационарной меры введенной марковской динамики. Это довольно общий факт, справедливый для широкого класса систем [32, 49]. Поиск равновесной конфигурации $\left\{c_{ij}^*\right\}_{i=1,\, j=1}^{n,\, n}$ представляет собой решение задачи

$$\min\left\{H(c): c \geq 0,\ \sum_{j=1}^{n} c_{ij} = l_i,\ \sum_{i=1}^{n} c_{ij} = \mathrm{w}_j\right\},$$

где $l = L/N$, $\mathrm{w} = W/N$, то есть (1). По $\left\{c_{ij}^*\right\}_{i=1,\, j=1}^{n,\, n}$ ($\left\{d_{ij}^*\right\}_{i=1,\, j=1}^{n,\, n} = N\left\{c_{ij}^*\right\}_{i=1,\, j=1}^{n,\, n}$) уже можно будет определить и равновесные $T_{ij} = T_{ij}\left(c^*\right)$.

В дальнейшем нам будет удобно привести задачу (1.1.1) к следующему виду (при помощи метода множителей Лагранжа [85], теоремы фон Неймана о минимаксе [119] и перенормировки $d := d/N$):

$$\min_{\lambda^L,\lambda^W} \max_{\sum_{i,j=1}^{n} d_{ij}=1,\, d_{ij}\geq 0} \left[-\sum_{i,j=1}^{n} d_{ij}\ln d_{ij} - \beta\sum_{i=1,\, j=1}^{n,\, n} d_{ij}T_{ij} + \sum_{i=1}^{n}\lambda_i^L\left(l_i - \sum_{j=1}^{n} d_{ij}\right) + \sum_{i=1}^{n}\lambda_j^W\left(\mathrm{w}_j - \sum_{i=1}^{n} d_{ij}\right)\right], \quad (1.1.2)$$

где, напомним, $l = L/N$, $\mathrm{w} = W/N$.[9] Вместо того, чтобы подставлять сюда зависимость $T_{ij}(d)$, мы используем следующий трюк. В подразделах 1.1.6–1.1.8 задача поиска зависимости $T_{ij}(d)$ будет сведена к задаче вида:

---

[9] Обратим внимание, что мы берем максимум в (1.1.2) при дополнительном ограничении $\sum_{i,j=1}^{n} d_{ij} = 1$, которого изначально не было. Однако легко показать, что это ограничения является следствием системы ограничений (А), точнее следствием одновременно двух подсистем в (А) – отвечающих $l$ и $w$. Как следствие, можно считать, что $\sum_{i=1}^{n}\lambda_i^L = 0$, $\sum_{j=1}^{n}\lambda_j^W = 0$. Отметим также, что задача (1.1.2) может быть упрощена, поскольку внутренний максимум явно находится, однако мы отложим соответствующие выкладки до подраздела 1.1.9. Здесь же мы рассмотрим случай, когда мы заносим в функционал с помощью множителей Лагранжа только часть ограничений (соответствующих $l$ или $w$). В таком случае, вычисляя внутренний максимум по $d$, получим, соответственно, задачи:

$$\min_{\lambda^W}\left[\sum_{i=1}^{n} l_i \ln\left[\sum_{j=1}^{n}\exp\left(-\lambda_j^W - \beta T_{ij}\right)\right] + \sum_{j=1}^{n}\lambda_j^W \mathrm{w}_j\right];\quad \min_{\lambda^L}\left[\sum_{j=1}^{n} w_j \ln\left[\sum_{i=1}^{n}\exp\left(-\lambda_i^L - \beta T_{ij}\right)\right] + \sum_{i=1}^{n}\lambda_i^L l_i\right],$$

где, соответственно

$$\min_{t\ge\bar t}\left\{\beta\left\langle \bar f, t-\bar t\right\rangle-\beta\sum_{i=1,\,j=1}^{n,\,n} d_{ij}T_{ij}(t)\right\}\overset{def}{=}-\Phi(d),$$

где $\bar t$ и $\bar f$ – известные векторы (входные параметры модели) временных затрат в пути на ребрах графа транспортной сети и максимальных пропускных способностей ребер графа транспортной сети, $T_{ij}(t)$ – длина кратчайшего пути из района $i$ в район $j$ (на графе транспортной сети имеется много вершин – намного больше числа районов, но мы считаем, что в каждом районе есть только одна вершина транспортного графа, являющаяся источником/стоком для пользователей сети, именно между такими представителями вершин районов $i$ и $j$ берется кратчайшие расстояние), если веса ребер графа транспортной сети задаются вектором $t$. Не стоит путать $T_{ij}(t)$ с искомой зависимостью $T_{ij}(d)$ – это разные функциональные зависимости. Решив указанную задачу $t^*(d)$ минимизации, считая $d$ параметрами, мы получили бы искомую зависимость $T_{ij}(d):=T_{ij}(t^*(d))$. Но явно это нельзя сделать в типичных ситуациях, да и не нужно, потому что в итоге все равно необходимо работать с потенциалом $\Phi(d)$. Поэтому предлагается (это предложение было сделано Ю.Е. Нестеровым в конце 2012 года) ввести в задачу (1.1.1) подзадачу и добавить слагаемое $\min\limits_{\lambda^L,\lambda^W}\max\limits_{\sum_{i,j=1}^{n} d_{ij}=1,\, d_{ij}\ge 0}\min\limits_{t\ge\bar t}\left[\,\cdot+\beta\left\langle \bar f, t-\bar t\right\rangle\right]$.

Решение такой задачи сразу даст все, что нужно. Это следует из формулы Демьянова–Данскина [85]. Подраздел 1.1.9 посвящен упрощению только что сформулированной конструкции. Более полное обоснование и дальнейшее развитие можно отследить по работам [9, 31, 35], а также в разделах 1.3, 1.4 этой главы 1.

Подобно тому, как мы рассматривали в этом пункте трудовые корреспонденции (в утренние и вечерние часы более 70% корреспонденций по Москве и области именно такие), можно рассматривать перемещения, например, из дома к местам учебы, отдыха, в магазины и т.п. (по-хорошему, еще надо было учитывать перемещения типа работа–магазин–детский_сад–дом) – рассмотрение всего этого вкупе приведет также к задаче (1.1.1). Только будет больше типов корреспонденций $d$: помимо пары районов, еще

---

$$d_{ij}=l_i\exp\left(-\lambda_j^W-\beta T_{ij}\right)\left[\sum_{k=1}^{n}\exp\left(-\lambda_k^W-\beta T_{ik}\right)\right]^{-1};\quad d_{ij}=w_j\exp\left(-\lambda_i^L-\beta T_{ij}\right)\left[\sum_{k=1}^{n}\exp\left(-\lambda_k^L-\beta T_{kj}\right)\right]^{-1}.$$

Отсюда можно усмотреть интерпретацию двойственных множителей как соответствующих “потенциалов притяжения/отталкивания районов” [26]. К этому мы еще вернемся в замечании п. 4. В принципе далее удобнее было бы работать с одной из этих задач, а не с (1.1.2), однако для сохранения симметрии, мы оставляем задачу в форме (1.1.2).

нужно будет учитывать тип корреспонденции [глава 2, 26]. Все это следует из того, что инвариантной мерой динамики с несколькими типами корреспонденций по-прежнему будет прямое произведение пуассоновских мер. Важно отметить при этом, что в таком контексте (равновесные) $T_{ij}$ будут определяться парами районов, а не типом передвижения, то есть с точки зрения последующего изложения это означает, что ничего по сути не поменяется. Другое дело, когда мы рассматриваем разного типа пользователей транспортной сети, например: имеющих личный автомобиль и не имеющих личный автомобиль. Первые могут им воспользоваться равно, как и общественным транспортом, а вторые нет. И на рассматриваемых масштабах времени пользователи могут менять свой тип. То есть время в пути может для разных типов пользователей быть различным [глава 2, 26]. Считая, подобно тому как мы делали раньше, что желание пользователей корреспонденции[10] $(i, j)$ сети сменить свой тип (вероятность в единицу времени) есть

$$\tilde{\lambda}\exp\Big( \underbrace{\tilde{R}\left(T_{ij}^{old}\right)}_{\substack{\text{суммарные затраты} \\ \text{до смены типа}}} - \underbrace{\tilde{R}\left(T_{ij}^{new}\right)}_{\substack{\text{суммарные затраты} \\ \text{после смены типа}}} \Big), \text{ где } \tilde{R}(T) = \beta T,$$

и учитывая в “обменах” тип пользователя (будет больше типов корреспонденций $d$, но “меняются местами работы” только пользователи одного типа), можно показать, что все это вкупе приведет также к задаче типа (1.1.1), и все последующие рассуждения распространяются на этот случай. Но с оговоркой, что в подразделе 1.1.8 расщепление потоков надо делать только для тех пользователей, для которых имеется возможность использовать личный транспорт. При этом часть пользователей (в зависимости от текущего $d$) распределяются только по сети общественного транспорта. Несложно понять, что это ничего принципиально не изменит. Фактически, в этом абзаце мы описали, как сделать из трехстадийной модели полноценную четырехстадийную модель, которыми обычно и пользуются на практике. Чтобы не делать изложение излишне громоздким, далее мы не учитываем нюансы, описанные в этом абзаце.

---

[10] Здесь, конечно, надо учитывать не только трудовые миграции, но и все остальные, поскольку для заметной части жителей Москвы и области решение о покупки автомобиля напрямую связано с желанием ездить на нем в основном только на дачу. В этом месте возникает необходимость моделирования (учета) перемещений пользователей сети не только в рамках установленного диапазона времени в течение типичного дня, но и в целом все возможные перемещения. Это обстоятельство вынуждает использовать здесь различного рода эвристики, дабы не отказываться от важного предположения: “в рамках установленного диапазона времени в течение типичного дня”.

### 1.1.4 Модель равновесного распределения потоков Бэкмана

Следуя книге [49], опишем наиболее популярную на протяжении более чем полувека модель равновесного распределения потоков Бэкмана [49, 112, 147, 289, 310].

Пусть транспортная сеть города представлена ориентированным графом $\Gamma=(V,E)$, где $V$ – узлы сети (вершины), $E\subset V\times V$ – дуги сети (рёбра графа). В современных моделях равновесного распределения потоков в крупном мегаполисе число узлов графа транспортной сети обычно выбирают порядка $|V|\sim 10^4-10^5$. Число ребер $|E|$ получается в три – четыре раза больше. Пусть $W\subseteq\{w=(i,j):i,j\in V\}$ – множество корреспонденций, т.е. возможных пар «исходный пункт» – «цель поездки» ($|W|$ по порядку величины может быть $n^2$); $p=\{v_1,v_2,...,v_m\}$ – путь из $v_1$ в $v_m$, если $(v_k,v_{k+1})\in E$, $k=1,...,m-1$, $m>1$; $P_w$ – множество путей, отвечающих корреспонденции $w\in W$, то есть если $w=(i,j)$, то $P_w$ – множество путей, начинающихся в вершине $i$ и заканчивающихся в $j$; $P=\bigcup_{w\in W}P_w$ – совокупность всех путей в сети $\Gamma$ (число “разумных” маршрутов $|P|$, которые потенциально могут использоваться, обычно растет с ростом числа узлов сети не быстрее чем $\mathrm{O}\left(|V|^3\right)$, однако теоретически может быть экспоненциально большим); $x_p$ [автомобилей/час] – величина потока по пути $p$, $x=\{x_p:\ p\in P\}$; $f_e$ [автомобилей/час] – величина потока по дуге $e$:

$$f_e(x)=\sum_{p\in P}\delta_{ep}x_p\text{, где }\delta_{ep}=\begin{cases}1, & e\in p\\ 0, & e\notin p\end{cases};$$

$\tau_e(f_e)$ – удельные затраты на проезд по дуге $e$. Как правило, предполагают, что это – (строго) возрастающие, гладкие функции от $f_e$ (в конце этого раздела нам потребуется еще и выпуклость). Точнее говоря, под $\tau_e(f_e)$ правильнее понимать представление пользователей транспортной сети об оценке собственных затрат (обычно временных в случае личного транспорта и комфортности пути (с учетом времени в пути) в случае общественного транспорта) при прохождении дуги $e$, если поток желающих оказаться на этой дуге будет $f_e$.

Зависимость $\tau_e(f_e)$ можно попробовать и вывести, например, из следующих соображений [49] (вариация на тему модели Бобкова–Буслаева–Танака). Рассматривается одна полоса длины $L$ и транспортный поток, характеризующийся максимальной скоростью $v_{\max}$ и следующей зависимостью безопасного расстояния (“комфортного” расстояния до впереди идущего транспортного средства) от скорости: $d(v)=l+\tilde{\tau}v+cv^2$ – динамический габарит, где $l$ – средняя длина автомобиля в “стоячей” пробке (эта длина немного больше средней “физической” длины автомобиля $\approx 6.5$ м), $\tilde{\tau}$ – время реакции водителей (эксперименты показывают, что для европейских водителей эта величина обычно равна одной секунде, для российских водителей она, как правило, не превышает полсекунды), $c$ – характеризует коэффициент трения шин о поверхность дороги, поскольку слагаемое $cv^2$ характеризует в $d(v)$ тормозной путь. Действительно, пока водитель среагирует на ситуацию он в среднем проедет (не изменяя своей скорости) путь $\tilde{\tau}v$. Когда уже реакция произошла, водитель начинает тормозить, и кинетическая энергия автомобиля $mv^2/2$ должна быть “погашена” работой силой трения на участке $h$ тормозного пути $\mu mgh$ ($\mu$ – коэффициент трения, $g$ – ускорение свободного падения). Отсюда можно найти $c = 1/(2\mu g)$. Имея зависимость $d(v)$, можно ввести зависимость $\rho(v) = 1/d(v)$, которая порождает зависимость потока от скорости $f(v) = v\rho(v)$. Используя то, что $\tau(f(v)) = L/v$ и $0 \le v \le v_{\max}$, можно явно выписать искомую зависимость $\tau(f)$. Несмотря на предложенный вывод, еще раз подчеркнем, что под $\tau_e(f_e)$ правильнее понимать представление пользователей транспортной сети об оценке собственных затрат при прохождении дуги $e$, если поток желающих оказаться на этой дуге будет $f_e$, поэтому функции $\tau_e(f_e)$, выбираемые на практике (см., например, подраздел 1.1.7), сильно отличаются от описанной в этом абзаце зависимости, приводящей, на самом деле, к наличию двух веток (см. в этой связи Рисунок 1.1.1 из подраздела 1.1.1).

Рассмотрим теперь $G_p(x)$ – затраты временные или финансовые на проезд по пути $p$. Естественно считать, что $G_p(x) = \sum_{e \in E} \tau_e(f_e(x))\delta_{ep}$. В приложениях часто требуется учитывать также затраты на прохождения вершин графа, которые могут зависеть от величин всех потоков через рассматриваемую вершину.

Пусть также известно, сколько перемещений в единицу времени $d_w$ осуществляется согласно корреспонденции $w \in W$. Тогда вектор $x$, характеризующий распределение потоков, должен лежать в допустимом множестве:

$$X = \left\{ x \ge 0 : \sum_{p \in P_w} x_p = d_w,\ w \in W \right\}.$$

Это множество может иметь и более сложный вид, если дополнительно учитывать, например, конечность пропускных способностей рёбер (ограничения сверху на $f_e$).

Рассмотрим игру, в которой каждому элементу $w \in W$ соответствует свой, достаточно большой ($d_w \gg 1$), набор однотипных “игроков”, осуществляющих передвижение согласно корреспонденции $w$. Чистыми стратегиями игрока служат пути, а выигрышем – величина $-G_p(x)$. Игрок “выбирает” путь следования $p \in P_w$, при этом, делая выбор, он пренебрегает тем, что от его выбора также “немного” зависят $|P_w|$ компонент вектора $x$ и, следовательно, сам выигрыш $-G_p(x)$. Можно показать, что отыскание равновесия Нэша–Вардропа $x^* \in X$ (макро описание равновесия) равносильно решению задачи нелинейной комплементарности (принцип Вардропа):

$$\textit{для любых } w \in W,\ p \in P_w \textit{ выполняется } x_p^* \cdot \left( G_p\left(x^*\right) - \min_{q \in P_w} G_q\left(x^*\right) \right) = 0.$$

Действительно допустим, что реализовалось какое-то другое равновесие $\tilde{x}^* \in X$, которое не удовлетворяет этому условию. Покажем, что тогда найдется водитель, которому выгодно поменять свой маршрут следования. Действительно, тогда

$$\text{существуют такие } \tilde{w} \in W,\ \tilde{p} \in P_{\tilde{w}}, \text{ что } \tilde{x}_{\tilde{p}}^* \cdot \left( G_{\tilde{p}}\left(\tilde{x}^*\right) - \min_{q \in P_{\tilde{w}}} G_q\left(\tilde{x}^*\right) \right) > 0.$$

Каждый водитель (множество таких водителей не пусто $\tilde{x}_{\tilde{p}}^* > 0$), принадлежащий корреспонденции $\tilde{w} \in W$, и использующий путь $\tilde{p} \in P_{\tilde{w}}$, действует не разумно, поскольку существует такой путь $\tilde{q} \in P_{\tilde{w}}$, $\tilde{q} \ne \tilde{p}$, что $G_{\tilde{q}}\left(\tilde{x}^*\right) = \min_{q \in P_{\tilde{w}}} G_q\left(\tilde{x}^*\right)$. Этот путь $\tilde{q}$ более выгоден, чем $\tilde{p}$. Аналогично показывается, что при $x^* \in X$ никому из водителей уже не выгодно отклоняться от своих стратегий. Но это по определению и называется равновесием Нэша, которое ввел в своей диссертации в конце 40-х годов XX века Джон Нэш, получивший именно за эту концепцию в 1994 г. нобелевскую премию по экономике [60]. Мы также добавляем фамилию Дж.Г. Вардропа, которой чуть позже Дж. Нэша привнес к этой концепции условие “конкурентного рынка”: игрок, принимающий решение, пренебрегает тем, что его решение сколько-нибудь значительно поменяет ситуацию на “рынке”. Когда игроков двое, трое (ситуации, рассматриваемые Дж. Нэшем), то, очевидно, что так делать нельзя. Но когда игроков (водителей) десятки и сотни тысяч… Вся эта конструкция неявно предполагает, что $x_p^* > 0 \Rightarrow x_p^* \gg 1$. Поэтому, не боясь сильно ошибиться, можно искать решение задачи нелинейной комплементарности, не предполагая целочисленности компонент вектора $x^* \in X$. Такая релаксация изначально целочисленной задачи заметно упрощает её с вычислительной точки зрения!

Хотя мы и смогли выписать условие равновесия в виде задачи нелинейной комплементарности, это не сильно продвинуло нас в понимании того, как его находить. Пытаться честно решить задачу в таком виде – вычислительно бесперспективная задача. С другой стороны современные вычислительные методы позволяют эффективно решать задачи выпуклой оптимизации. Постараемся свести нашу задачу к таковой.

Для этого, прежде всего, заметим, что рассматриваемая нами игра принадлежит к классу, так называемых, потенциальных игр. В нашем случае это означает, что существует такая функция

$$\Psi(x) = \sum_{e\in E} \int_0^{\sum_{p\in P}\delta_{ep}x_p} \tau_e(z)\,dz = \sum_{e\in E} \sigma_e\left(f_e(x)\right),$$

где $\sigma_e(f_e) = \int_0^{f_e(x)} \tau_e(z)\,dz$, что $\partial\Psi(x)/\partial x_p = G_p(x)$ для любого $p \in P$. Таким образом, мы имеем дело с потенциальной игрой. Оказывается, что $x^* \in X$ – равновесие Нэша–Вардропа тогда и только тогда, когда оно доставляет минимум $\Psi(x)$ на множестве $X$. Действительно, предположим, что $x^* \in X$ – точка минимума. Тогда, в частности, для любых $w \in W$, $p,q \in P_w$ ($x_p^* > 0$) и достаточно маленького $\delta x_p > 0$ выполняется:

$$-\frac{\partial\Psi(x^*)}{\partial x_p}\delta x_p + \frac{\partial\Psi(x^*)}{\partial x_q}\delta x_p \ge 0.$$

Иначе, заменив $x^*$ на

$$\breve{x}^* = x^* + \Big(\underbrace{\overbrace{0,\ldots,0,-\delta x_p}^{p},0,\ldots,0,\delta x_p}_{q},0,\ldots,0\Big) \in X,$$

мы пришли бы к вектору $\breve{x}^*$, доставляющему меньшее значение $\Psi(x)$ на множестве $X$:

$$\Psi(\breve{x}^*) \approx \Psi(x^*) - \frac{\partial\Psi(x^*)}{\partial x_p}\delta x_p + \frac{\partial\Psi(x^*)}{\partial x_q}\delta x_p < \Psi(x^*).$$

Вспоминая, что $\partial\Psi(x)/\partial x_p = G_p(x)$, и учитывая, что $q$ можно выбирать произвольно из множества $P_w$, получаем:

*для любых* $w \in W$, $p \in P_w$, *если* $x_p^* > 0$, *то выполняется* $\min_{q\in P_w} G_q(x^*) \ge G_p(x^*)$.

Но это и есть по-другому записанное условие нелинейной комплементарности. Строго говоря, мы показали сейчас только то, что точка минимума $\Psi(x)$ на множестве $X$ будет равновесием Нэша–Вардропа. Аналогично рассуждая, можно показать и обратное: равновесие Нэша–Вардропа доставляет минимум $\Psi(x)$ на множестве $X$. Этот минимум мож-но искать, например, с помощью метода условного градиента [40, 123] (Франк–Вульф). Подробнее об этом методе также написано в разделе 1.5 этой главы.

**Теорема 1.1.1 [**49, 112, 289, 310**].** *Вектор* $x^*$ *будет равновесием Нэша–Вардропа тогда и только тогда, когда*

$$x \in \operatorname{Arg}\min_x \left[ \Psi\left(f\left(x\right)\right) = \sum_{e\in E} \sigma_e\left(f_e\left(x\right)\right):\ f = \Theta x,\ x \in X \right]. \tag{1.1.3}$$

*Если преобразование* $G(\cdot)$ *строго монотонное, то равновесие* $x$ *единственно. Если* $\tau_e'(\cdot) > 0$*, то равновесный вектор распределения потоков по ребрам* $f$ *– единственный* (*это еще не гарантирует единственность вектора распределения потоков по путям* $x$ *[*49*]*).

Проинтерпретируем, следуя [36, 49, 137, 302] (см. также разделы 1.2, 1.4 этой главы 1), эволюционным образом полученное равновесие, попутно отвечая на вопрос (поскольку задача (1.1.4) ниже имеет единственное решение, детали см. в [36, 143]): какому из равновесий стоит отдать предпочтение, в случае неединственности? Опишем марковскую логит динамику (также говорят гиббсовскую динамику) в повторяющейся игре загрузки графа транспортной сети [302]. Пусть каждой корреспонденции отвечает $d_w M$ агентов ($M \gg 1$), $x := x/M$, $f := f/M$, $\tau_e\left(f_e\right) := \tau_e\left(f_e/M\right)$. Пусть имеется $tN$ шагов ($N \gg 1$). Пусть $k$-й агент, принадлежащий корреспонденции $w \in W$, независимо от остальных на шаге $m+1$ с вероятностью с вероятность $1-\lambda/N$ выбирает путь $p^{k,m}$, который использовал на шаге $m = 0,...,tN$, а с вероятностью $\lambda/N$ ($\lambda > 0$) решает “поменять” путь, и выбирает (возможно тот же самый) зашумленный кратчайший путь

$$p^{k,m+1} = \arg\max_{q\in P_w} \left\{ -G_q\left(x^m\right) + \xi_q^{k,m+1} \right\},$$

где независимые случайные величины $\xi_q^{k,m+1}$, имеют одинаковое двойное экспоненциальное распределение, также называемое распределением Гумбеля[11] [137, 302]:

$$P\left(\xi_q^{k,m+1} < \zeta\right) = \exp\left\{-e^{-\zeta/T-E}\right\},\ T > 0.$$

Отметим, что

---

[11] Распределение Гумбеля можно объяснить исходя из идемпотентного аналога центральной предельной теоремы (вместо суммы случайных величин – максимум) для независимых случайных величин с экспоненциальным и более быстро убывающим правым хвостом [84]. Распределение Гумбеля возникает в данном контексте, например, если при принятии решения водитель собирает информацию с большого числа разных (независимых) зашумленных источников, ориентируясь на худшие прогнозы по каждому из путей.

$$P\left(p^{k,m+1}=p \,\middle|\, \text{агент решил "поменять" путь}\right)=\frac{\exp\left(-G_p\left(x^m\right)/T\right)}{\sum_{q\in P_w}\exp\left(-G_q\left(x^m\right)/T\right)}.$$

Кроме того, если взять $E\approx 0.5772$ – константа Эйлера, то

$$M\left[\xi_q^{k,m+1}\right]=0,\ D\left[\xi_q^{k,m+1}\right]=T^2\pi^2/6.$$

Такая динамика отражает ограниченную рациональность агентов (см. подраздел 1.1.5), и часто используется в популяционной теории игр [302] и теории дискретного выбора [137]. Оказывается, эта марковская динамика в пределе $N\to\infty$ превращается в марковскую логит динамику в непрерывном времени (вырождающуюся при $T\to 0+$ в динамику наилучших ответов [302] – последующие рассуждения, в частности, формулы (1.1.4), (1.1.5), допускают переход к пределу $T\to 0+$). Марковская логит динамика в непрерывном времени допускает два предельных перехода (обоснование перестановочности этих пределов см. в [10, 195]): $t\to\infty$, $M\to\infty$ или $M\to\infty$, $t\to\infty$. При первом порядке переходов мы сначала ($t\to\infty$) согласно эргодической теореме для марковских процессов (в нашем случае марковский процесс – модель стохастической химической кинетики с унарными реакциями в условиях детального баланса [10, 36]) приходим к финальной (=стационарной) вероятностной мере, имеющей в основе мультиномиальное распределение. С ростом числа агентов ($M\to\infty$) эта мера

$$\sim\exp\left(-\frac{M}{T}\cdot\left(\Psi_T\left(x\right)+o\left(1\right)\right)\right)$$

экспоненциально концентрируется около наиболее вероятного состояния, поиск которого сводится к решению энтропийно регуляризованной задачи (1.1.3) (как численно решать эту задачу описано в работе [33], см. также [6])

$$\Psi_T\left(x\right)=\Psi\left(f\left(x\right)\right)+T\sum_{w\in W}\sum_{p\in P_w}x_p\ln x_p\to\min_{\substack{f(x)=\Theta x\\ x\in X}}. \qquad (1.1.4)$$

Функционал в этой задаче оптимизации с точностью до потенцирования и мультипликативных и аддитивных констант соответствует исследуемой стационарной мере – то есть это функционал Санова [10, 36]. При обратном порядке предельных переходов, мы сначала ($M\to\infty$) осуществляем, так называемый, канонический скейлинг [10, 195], приводящий к детерминированной кинетической динамике, описываемой СОДУ на $x$

$$\frac{dx_p}{dt}=d_w\frac{\exp\left(-G_p\left(x\right)/T\right)}{\sum_{l\in P_w}\exp\left(-G_l\left(x\right)/T\right)}-x_p,\ p\in P_w,\ w\in W, \qquad (1.1.5)$$

а затем ($t\to\infty$) ищем аттрактор получившейся СОДУ. Глобальным аттрактором оказывается неподвижная точка, которая определяется решением задачи (4). Более того, функцио-нал $\Psi_T(x)$ является функцией Ляпунова полученной кинетической динамики (5) (то есть является функционалом Больцмана). Последнее утверждение – достаточно общий факт (функционал Санова, является функционалом Больцмана), верный при намного более общих условиях (см. [10] и цитированную там литературу).

Хотелось бы подчеркнуть, что рассматриваемая выше “игра” – потенциальная (это общий факт для игр загрузок [302]; Розенталь 1973, Мондерер–Шэпли 1996), поэтому из общих результатов эволюционной теории игр [174, 302], следует, что любые разумные содержательно интерпретируемые (суб-)градиентные спуски приводят к равновесию Нэша–Вардропа (или его стохастическому варианту, даваемому решением задачи (1.1.4)). В частности, в [49, 51, 54] содержательно интерпретируется быстро сходящаяся динамика, связанная с методом зеркального спуска, которая порождается имитационной логит динамикой [302]. Хотелось бы также обратить внимание на эволюционную интерпретацию парадокса Браесса: когда неэффективное по Парето, единственное равновесие Нэша–Вардропа в специально сконструированной транспортной сети является, тем не менее, эволюционно устойчивым [49, 302].

Нетривиальным является следующее наблюдение. Если рассмотреть энтропийно-регуляризованный функционал $\Psi_T(x)$ и взять предел (см. подразделы 1.1.6, 1.1.7)

$$\tau_e^{\mu}(f_e)\xrightarrow[\mu\to 0+]{}\begin{cases}\bar{t}_e, & 0\le f_e<\bar{f}_e\\ [\bar{t}_e,\infty), & f_e=\bar{f}_e\end{cases},$$

$$d\tau_e^{\mu}(f_e)/df_e\xrightarrow[\mu\to 0+]{}0,\ 0\le f_e<\bar{f}_e,$$

то переход к двойственной задаче (для задачи минимизации этого функционала на множестве $X$) дает стохастический вариант модели стабильной динамики (см. подраздел 1.1.6), который используется в стохастическом варианте трехстадийной модели стабильной динамики (см. подраздел 1.1.10).

**Замечание 1.1.1 (“облачная модель” расчета матрицы корреспонденций [36], см. также раздел 1.2 этой главы 1).** В контексте написанного выше полезно отметить другой способ обоснования энтропийной модели расчета матрицы корреспонденций из подраздела 1.1.3. Предположим, что все вершины, отвечающие источникам корреспонденций, соединены ребрами с одной вспомогательной вершиной (облако № 1). Аналогично, все вершины, отвечающие стокам корреспонденций, соединены ребрами с другой вспомогательной вершиной (облако № 2). Припишем всем новым ребрам постоянные веса. И проинтерпретируем веса ребер, отвечающих источникам $\lambda_i^L$, например, как средние затраты на проживание (в единицу времени, скажем, в день) в этом источнике (районе), а веса ребер, отвечающих стокам $\lambda_j^W$, как уровень средней заработной платы со знаком минус (в единицу времени) в этом стоке (районе), если изучаем трудовые корреспонденции. Будем следить за системой в медленном времени, то есть будем считать, что равновесное распределение потоков по путям стационарно. Поскольку речь идет о равновесном распределении потоков, то нет необходимости говорить о затратах на путях или ребрах детализированного транспортного графа, достаточно говорить только о затратах (в единицу времени), отвечающих той или иной корреспонденции. Таким образом, у нас есть взвешенный транспортный граф с одним источником (облако 1) и одним стоком (облако 2). Все вершины этого графа, кроме двух вспомогательных (облаков), соответствуют районам в модели расчета матрицы корреспонденций из подраздела 1.3. Все ребра этого графа имеют стационарные (не меняющиеся и не зависящие от текущих корреспонденций) веса $\left\{T_{ij};\lambda_i^L;\lambda_j^W\right\}$. Если рассмотреть естественную в данном контексте логит динамику с $T=1/\beta$ (здесь полезно напомнить, что согласно подразделу 1.3 $\beta$ обратно пропорционально средним затратам, а $T$ имеет как раз физическую размерность затрат), описанную выше, то поиск равновесия рассматриваемой макросистемы приводит (в прошкалированных переменных) к задаче, сильно похожей на задачу (1.1.2) из подраздела 1.1.3

$$\max_{\sum_{i,j=1}^{n} d_{ij}=1,\, d_{ij}\ge 0}\left[-\sum_{i,j=1}^{n} d_{ij}\ln d_{ij}-\beta\sum_{i=1,\,j=1}^{n,n} d_{ij}T_{ij}-\beta\sum_{i=1}^{n}\left(\lambda_i^L\sum_{j=1}^{n} d_{ij}\right)-\beta\sum_{i=1}^{n}\left(\lambda_j^W\sum_{i=1}^{n} d_{ij}\right)\right].$$

Разница состоит в том, что здесь мы не оптимизируем по $2n$ двойственным множителям $\lambda^L$, $\lambda^W$. Более того, мы их и не интерпретируем здесь как двойственные множители, поскольку мы их ввели на этапе взвешивания ребер графа. Тем не менее, значения этих переменных, как правило, не откуда брать. Тем более, что приведенная выше (наивная) интерпретация вряд ли может всерьез рассматриваться, как способ определения этих параметров исходя из данных статистики. Более правильно понимать $\lambda_i^L$, $\lambda_j^W$ как потенциалы притяжения/отталкивания районов, включающиеся в себя плату за жилье и зарплату, но включающие также и многое другое, что сложно описать количественно. И здесь как раз помогает информация об источниках и стоках, содержащаяся в $2n$ уравнениях из формулы (**А**) подраздела 1.1.3. Таким образом, мы приходим ровно к той же самой задаче (1.1.2) с той лишь разницей, что мы получили дополнительную интерпретацию двойственных множителей в задаче (1.1.2). При этом двойственные множители в задаче (1.1.2) равны (с точностью до мультипликативного фактора $\beta$) введенным здесь потенциалам притяжения районов. Несложно распространить на изложенную здесь модель написанное в подразделе 1.1.3 по поводу того, как с помощью разделения времен можно учитывать обратную связь:

перераспределение потоков по путям изменяется (в быстром времени) при изменении корреспонденций, а также распространить на то, что написано в самом конце подраздела 1.1.3. Нам представляется такой способ рассуждения даже более привлекательным, чем способ, описанный в подразделе 1.1.3, и основанный на “обменах”. И связано это с тем, что для получения равновесия в многостадийной модели, мы можем рассмотреть всего одну (общую) логит динамику, в которой с малой интенсивностью (в медленном времени) происходят переходы, описанные в этом замечании (жители города меняют места жительства, работы), а с высокой интенсивностью (в быстром времени, изо дня в день) жители города перераспределяются по путям (в зависимости от текущих корреспонденций, подстраиваясь под корреспонденции) – это как раз и было описано непосредственно перед замечанием. Другая причина – большая вариативность модели построенной в этом замечании. Нам представляется очень плодотворной и перспективной идея перенесения имеющейся информации об исследуемой системе из обременительных законов сохранения динамики, описывающей эволюцию этой системы, в саму динамику путем введения дополнительных естественно интерпретируемых параметров. При таком подходе становится возможным, например, учитывать в моделях и рост транспортной сети. Другими словами, при таком подходе, например, можно естественным образом рассматривать также и ситуацию, когда число пользователей транспортной сетью меняется со временем (медленно).

В заключение отметим, что если штрафовать (назначать платы) за использование различных стратегий (маршрутов) по правилу

$$\bar{G}_p(x) = G_p(x) + \underbrace{\sum_{w\in W}\sum_{q\in P_w} x_q \frac{\partial G_q(x)}{\partial x_p}}_{штраф}, \quad p \in P,$$

то найдется такая функция $\bar{\Psi}(x) = \sum_{w\in W}\sum_{p\in P_w} x_p G_p(x)$, что $\partial\bar{\Psi}(x)/\partial x_p = \bar{G}_p(x)$. Поэтому из сказанного выше в этом пункте будет следовать, что возникающее в такой управляемой транспортной сети равновесие Нэша–Вардропа будет (единственным) глобально устойчивым и соответствовать социальному оптимуму в изначальной транспортной сети [301]. Чтобы в этом убедиться, достаточно (в виду линейной связи $f = \Theta x$) проверить строгую выпуклость функции $\bar{\Psi}(x) = \sum_{w\in W}\sum_{p\in P_w} x_p G_p(x) = \sum_{e\in E} f_e(x)\tau_e(f_e(x))$, для чего достаточно выпуклости функций $\tau_e(f_e)$, $e \in E$. Все это хорошо соответствует механизму Викри–Кларка–Гроуса [134] (VCG mechanism) – штраф (плата) за использование маршрута новым пользователем равен дополнительным потерям, которые понесут из-за этого все ос-тальные пользователи. Однако важно сделать две оговорки. Во-первых, все это хотя и можно попытаться практически осуществить (например, собирая транспортные налоги исходя из трековой информации, которую в перспективе можно будет иметь о каждом автомобиле), но довольно сложным оказывается механизм. Платы взимаются не за проезд по ребру, как хотелось бы, а именно за выбор (проезд) маршрута. Кроме того, плата является функцией состояния транспортной системы $x$, которое, в отличие от $f$, не наблюдаемо. Во-вторых, взимая платы за проезд, мы с одной стороны приводим систему в социальный оптимум, а с другой стороны для достижения этой цели вынуждены собирать с участников движения налоги. К сожалению, их размер может оказаться внушительным, и это уже необходимо учитывать с точки зрения расщепления участников движения по выбору типа передвижения. Относительно второй проблемы – готовых решений нам не известно. Это известная проблема в современном разделе теории игр: mechanism design [134]. А вот по первой проблеме (адаптивное) решение есть [301]:

$$\bar{\tau}_e\left(f_e\right)=\tau_e\left(f_e\right)+\underbrace{f_e\tau_e'\left(f_e\right)}_{штраф},\ e\in E.$$

Легко проверить, что это приведет к указанному выше пересчету $G_p\left(x\right)\to\bar{G}_p\left(x\right)$. Далее, если мы зафиксируем (знаем) социальный оптимум $f^{opt}$, то платы за проезд можно выбирать постоянными $f_e^{opt}\tau_e'\left(f_e^{opt}\right),\ e\in E$. Все сказанное выше об устойчивости останется в силе (без всяких дополнительных предположений о выпуклости функций $\tau_e\left(f_e\right),\ e\in E$) с одной лишь оговоркой, что транспортная система должна поддерживаться при заданных корреспонденциях $d_w$ и функциях затрат $\tau_e\left(f_e\right)$. В противном случае, возникающее равновесие уже может не соответствовать социальному оптимуму.

### 1.1.5 Краткий обзор подходов к построению многостадийных моделей транспортных потоков, с моделью типа Бэкмана в качестве модели равновесного распределения потоков[12]

Так называемая «4-стадийная» модель является наиболее употребительной методологией моделирования транспортных систем городов и агломераций (см., например, [287]). В рамках данной модели производится поиск равновесия спроса и предложения на поездки. При этом рассматриваются в единой совокупности модели генерации трафика, его распределения по типам передвижения и дальнейшее распределение по маршрутам.

---

[12] Этот и следующий подраздел были написаны совместно с Ю.В. Дорном. Также Ю.В. Дорном была оказана помощь в проверке формул, приведенных в подразделе 1.1.7.

Методология, как уже описывалось ранее, включает последовательное выполнение четырех этапов, последние три из которых закольцовываются для получения самосогласованных результатов. Исходными данными для модели являются: граф дорожной сети и сети общественного транспорта с заданными функциями издержек и других определяющих параметров, разделение города на транспортные зоны и параметры этих зон (например, количество рабочих мест или мест жительства). На выходе модели выдаются: оценка матрицы корреспонденций для каждого типа передвижения, загрузка элементов сети (например, конкретной дороги) и издержки, соответствующие данному уровню загрузки. Понятие «тип передвижения», используемое выше, является обобщением понятия тип транспорта и обозначает последовательность (или просто множество) используемых типов транспорта. Например, типом передвижения может считаться поездка на общественном транспорте (не важно, автобусе или троллейбусе или на них обоих поочередно с пересадками) или же использование схемы park-and-ride. Уровень детализации при этом определяется самим модельером. Обычно выбирается или модель с тремя типами передвижений (общественный и личный транспорт, пешие прогулки) или с детализацией до типа транспорта (автомобиль, наземный общественный транспорт, метро, пешие прогулки, и различные комбинации, перечисленные ранее). При этом стараются не учитывать типы передвижений, которые используются очень редко.

Структурно, расчет модели можно описать следующим образом:

1. На первом шаге из исходных данных о транспортных зонах получают вектора отправления и прибытия для транспортных зон, т.е. в наших обозначениях $\{L_i\}_{i=1}^n$ и $\{W_j\}_{j=1}^n$.
2. Рассчитывается первая оценка матрицы корреспонденций (обычно, с помощью гравитационной модели [49, 123]).
3. Рассчитывается расщепление корреспонденций по типам передвижений.
4. Для каждого типа передвижений рассчитывается распределение потоков по маршрутам.
5. Получается оценка матрицы издержек корреспонденций и вектор загрузки сети.
6. Проверяется критерий остановки.
7. Если критерий выполнен, решение получено.
8. Если критерий не выполнен, вернуться на шаг 2 с переоцененной матрицей издержек.

Первый этап модели мы разбирать не будем, так как он является достаточно обособленным от других в том смысле, что он не входит в итерационную часть метода и практи-чески «бесплатен» с точки зрения сложности операций. Отметим лишь, что вектора отправления и прибытия рассчитываются из параметров зон с помощью простейших регрессионных моделей.

Рассмотрим модели, лежащие в основе этого алгоритма более подробно. Данный алгоритм – итерационный, после «прогонки» очередной итерации мы возвращаемся на второй «шаг» схемы. На $m$-й итерации на 2-м шаге формируется $m$-я оценка матрицы корреспонденций. Для её построения необходимо знание о матрице издержек в сети (т.е. знание о стоимости проезда из каждой зоны в каждую). На всех шагах, кроме первого, данные значения получаются из предыдущей ($m$–1)-й итерации. На первой итерации используются издержки, соответствующие незагруженной сети или почти любая другая, разумная, оценка матрицы издержек. Если в сети для пользователей доступно несколько типов передвижений, то в качестве оценки издержек берутся или средние издержки (время в пути) для конкретной корреспонденции или же (если пользователи распределены «равновесно») значения издержек, например, для личного транспорта, так как в «равновесии» издержки зачастую (исключая некоторые особые случаи) должны быть равны для различных, используемых, типов передвижений.

Как правило, оценки матрицы корреспонденций строятся согласно гравитационной или энтропийной модели [49, 123].

Посмотрим на то, как замыкается модель при наличии различных слоев спроса и типов передвижений. После того, как был проведен первый этап, модельеру уже известно сколько людей выезжает и въезжает в каждую транспортную зону, а также известно какая доля этих людей совершает поездку того или иного типа, т.е. распределение поездок по слоям спроса. Для каждого слоя спроса формируется своя матрица корреспонденций, при этом используется одна и та же матрица издержек (так как используемая транспортная сеть для всех жителей одна и та же), однако параметры $\{\beta\}$ для каждого слоя спроса свои. После этого, данные матрицы корреспонденций суммируются. Полученная агрегированная матрица корреспонденций используется для проведения этапа расщепления корреспонденций (суммарных) по типам передвижений. Далее алгоритм работает только с ней (агрегированной матрицей) до начала следующей итерации, когда вся процедура повторяется. Более подробно вся процедура рассматривается ниже.

После построения очередной оценки матрицы корреспонденций, используя, опять же, матрицу издержек, происходит расщепление корреспонденций по типам передвижений. Для этого используются модели дискретного выбора, родственные уже описанной выше logit-choice модели. На выходе данного этапа (3-го шага в описанной схеме) получаются матрицы корреспонденций для каждого типа передвижений.

Отметим, что нет единой методологии и стандарта, какую из моделей дискретного выбора стоит использовать. Каждая из моделей имеет свои недостатки. Так, например, для модели logit-choice получаемое стохастическое равновесие даже в простейших постановках может сильно отличаться от равновесия Нэша. Также результаты данной модели очень чувствительны к тому, как задается набор альтернатив для выбора. Например, расщепление для дерева выбора (личный транспорт, автобус) и для дерева выбора (личный транспорт, зеленый автобус, красный автобус) будут разными (даже если выбор описывается для одной и той же транспортной системы и вся разница только в том, что во втором случае указан цвет автобуса). Наиболее используемыми на практике являются родственные logit-модели Nested Logit Model и Multinomial Logit Model, а также композитная модель Mixed Logit Model (см. [287]).

Важно сказать, что порядок шагов 2 и 3 зависит от того, какой слой спроса описывается, т.е., грубо говоря, от цели поездки и некоторых параметров жителей. Например, слой спроса может определяться как «люди предпенсионного возраста, совершающие поездки из дома на дачу». Иногда слои спроса определяются как цель поездки с учетом пункта отправления, например «поездка из дома на работу». Выше мы отмечали, что характерное время формирования корреспонденций – годы. Это справедливо для поездок дом–работа и обратно. Однако данный тезис кажется не совсем точным, например, для поездок за покупками. В этом случае можно предположить, что люди определяют место покупок уже после того, как будет определен тип передвижения, т.е. для данного слоя спроса, шаги 2 и 3 оказывается возможным поменять местами.

Наконец, последним этапом (шаг 4) является расчет равновесного распределения потоков для каждого типа передвижения. Для личного транспорта при этом используется модель Бэкмана. Для общественного транспорта не существует единого подхода к моделированию и, как правило, используются композитные модели, включающие элементы моделей дискретного выбора. Стоит лишь отметить, что при численном расчете моделей реального города именно этап расчета равновесного распределения потоков является самым затратным по времени. На выходе этого этапа получаются вектора загрузки сети (т.е. потоки по ребрам транспортного графа) и матрица издержек (шаг 5).

Шаги 6–8 в приведенной выше схеме отвечают за создание обратной связи в модели, которая позволит учесть взаимное влияние матрицы корреспонденций и матрицы издержек. Напомним, что мы использовали матрицу издержек для расчета матрицы корреспонденций на шаге 2 и матрицу корреспонденций для получения матрицы издержек на шаге 4. Логично требовать, чтобы матрица издержек, подаваемая на вход на шаге 2 и матрица издержек, получаемая на выходе на шаге 5, если не совпадали, то были «близки» по какому-либо разумному критерию.

Критерием остановки служит равенство (с требуемой точностью) средних издержек для матрицы издержек, полученной по модели расчета матрицы корреспонденций, и средних издержек, полученных эмпирическим путем. Мы к этому еще вернемся чуть позже. В статьях [304, 328] было показано, что данный критерий остановки в нашей постановке соответствует оценке матрицы корреспонденций методом максимального правдоподобия [319] для описанной далее модели.

Пусть имеются данные о реализации некоторой случайной матрицы корреспонденций $\{R_{ij}\}$, каждый элемент которой (независимо от всех остальных) распределен по закону Пуассона с математическим ожиданием $M\left(R_{ij}\right)=d_{ij}$, причем верно параметрическое предположение о том, что матрица корреспонденций представляется следующим образом:

$$d_{ij}=\tilde{A}_i\tilde{B}_j f\left(C_{ij}\right),$$

где $f\left(C_{ij}\right)$ – функция притяжения (часто выбирают $f\left(C_{ij}\right)=\exp\left(-\beta C_{ij}\right)$), а матрица издержек $\{C_{ij}\}$ считается известной.[13]

Зададимся целью найти оценку $\{d_{ij}\}$ при заданной (наблюдаемой) матрице $\{R_{ij}\}$ и известной матрице $\{C_{ij}\}$ методом максимума правдоподобия (на основе теоремы Фишера) [319]. Точнее говоря, оценивать требуется не саму матрицу $\{d_{ij}\}$, а неизвестные параметры $\tilde{A},\tilde{B},\beta$. Для возможности применять теорему Фишера [319], и таким образом гарантировать хорошие (асимптотические) свойства полученных оценок (в 1-норме), будем считать, что число оцениваемых параметров $p=2n+1$ и объем выборки $N=\sum_{i,j=1}^{n}R_{ij}$ удовлетворяют следующему соотношению: $p/N\ll 1$.

[13] Важно заметить, что мы допускаем, следуя В.Г. Спокойному, “model misspecification” [319, 101], т.е., что эти предположения не верны, и истинное распределение вероятностей $\{R_{ij}\}$ не лежит в этом семействе. Тогда полученное решение по методу максимума правдоподобия можно интерпретировать, как дающее асимптотически (по $N=\sum_{i,j=1}^{n}R_{ij}$) наиболее близкое (по расстоянию Кульбака–Лейблера) распределение в этом семействе к истинному распределению. Другими словами, так полученные параметры – являются “асимптотически наилучшими” оценками параметров проекции (по расстоянию Кульбака–Лейблера) истинного распределения вероятностей на выбранное нами параметрическое семейство. Детали см. в [319].

Из определения распределения Пуассона следует, что вероятность реализации корреспонденции $R_{ij}$ может быть посчитана как:

$$P\left(R_{ij} \mid d_{ij}\right) = \frac{\exp\left(-d_{ij}\right) \cdot d_{ij}^{R_{ij}}}{R_{ij}\,!}.$$

Функция правдоподобия (вероятность того, что “выпадет” матрица $\left\{R_{ij}\right\}$, если значения параметров $\tilde{A}, \tilde{B}, \beta$) будет иметь следующий вид:

$$\Lambda\left(\left\{R_{ij}\right\} \middle| \tilde{A}, \tilde{B}, \beta\right) = \prod_{i,j=1}^{n} \frac{\exp\left(-d_{ij}\right) \cdot d_{ij}^{R_{ij}}}{R_{ij}\,!} = \prod_{i,j=1}^{n} \frac{\exp\left(-\tilde{A}_i \tilde{B}_j f\left(C_{ij}\right)\right) \cdot \left(\tilde{A}_i \tilde{B}_j f\left(C_{ij}\right)\right)^{R_{ij}}}{R_{ij}\,!}.$$

Нам нужно найти точку $\left(\tilde{A}, \tilde{B}, \beta\right)$, в которой достигается максимум функции правдоподобия. Как известно, точка максимума неотрицательной функции (каковой по определению является вероятность) не изменится, если решать задачу максимизации не для исходной функции, а для ее логарифма. Перейдем к логарифму функции правдоподобия:

$$\ln \Lambda = \sum_{i,j=1}^{n} \left(-\tilde{A}_i \tilde{B}_j f\left(C_{ij}\right) + R_{ij} \cdot \left(\ln \tilde{A}_i + \ln \tilde{B}_j + \ln f\left(C_{ij}\right)\right) - \ln R_{ij}\,!\right).$$

Получаем следующую задачу оптимизации:

$$\sum_{i,j=1}^{n} \left(-\tilde{A}_i \tilde{B}_j f\left(C_{ij}\right) + R_{ij} \cdot \left(\ln \tilde{A}_i + \ln \tilde{B}_j + \ln f\left(C_{ij}\right)\right) - \ln R_{ij}\,!\right) \to \max_{\tilde{A} \geq 0,\, \tilde{B} \geq 0,\, \beta}.$$

Выпишем условия оптимальности:[14]

$$\frac{\partial \ln \Lambda}{\partial \tilde{A}_i} = \sum_{j=1}^{n} \left(-\tilde{B}_j \exp\left(-\beta C_{ij}\right) + \frac{R_{ij}}{\tilde{A}_i}\right) = 0,$$

$$\frac{\partial \ln \Lambda}{\partial \tilde{B}_j} = \sum_{i=1}^{n} \left(-\tilde{A}_i \exp\left(-\beta C_{ij}\right) + \frac{R_{ij}}{\tilde{B}_j}\right) = 0,$$

$$\frac{\partial \ln \Lambda}{\partial \beta} = \sum_{i,j=1}^{n} \left(\tilde{A}_i \tilde{B}_j C_{ij} \exp\left(-\beta C_{ij}\right) - R_{ij} C_{ij}\right) = 0.$$

Мы получили $2n+1$ уравнение максимума правдоподобия:

$$\sum_{j=1}^{n} \tilde{A}_i \tilde{B}_j \exp\left(-\beta C_{ij}\right) = \sum_{j=1}^{n} R_{ij},$$

---

[14] Легко понять, что максимум не может достигаться на границе. Если допустить, что, скажем, $L_i$ равно нулю в точке максимума, то поскольку $\partial \ln \Lambda / \partial L_i = \infty$ в этой точке, сдвинувшись немного перпендикулярно гиперплоскости $L_i = 0$ внутрь области определения, мы увеличили бы значения функционала – то есть пришли бы к противоречию с предположением о равенстве нулю $L_i$ в точке максимума.

$$\sum_{i=1}^{n} \tilde{A}_i \tilde{B}_j \exp\left(-\beta C_{ij}\right) = \sum_{i=1}^{n} R_{ij} ,$$

$$\sum_{i,j=1}^{n} \left( \tilde{A}_i \tilde{B}_j C_{ij} \exp\left(-\beta C_{ij}\right) - R_{ij} C_{ij} \right) = 0 .$$

Положим по определению $L_i = \sum_{j=1}^{n} R_{ij}$ и $W_j = \sum_{i=1}^{n} R_{ij}$, $\tilde{A}_i = L_i A_i$ и $\tilde{B}_j = W_j B_j$. Тогда

$$d_{ij} = A_i L_i B_j W_j f\left(C_{ij}\right),$$

где $A_i$, $B_j$ – структурные параметры (гравитационной) модели. Их рассчитывают с помощью простого итерационного алгоритма (метод балансировки = метод простых итераций [49, 55, 105, 123, 198]) по значениям $L_i$, $W_j$ следующим образом:[15]

$$A_i = \frac{1}{\sum_{j=1}^{n} B_j W_j f\left(C_{ij}\right)}, \; B_j = \frac{1}{\sum_{i=1}^{n} A_i L_i f\left(C_{ij}\right)}. \tag{1.1.6}$$

Последнее же уравнение системы уравнений максимума правдоподобия даст нам:

$$\sum_{i,j=1}^{n} \left( d_{ij} C_{ij} - R_{ij} C_{ij} \right) = 0 \text{ или } \sum_{i,j=1}^{n} R_{ij} C_{ij} = \sum_{i,j=1}^{n} d_{ij} C_{ij} . \tag{1.1.7}$$

**Замечание 1.1.2.** Обратим внимание, что к тем же самым соотношениям можно было прийти (при той же параметрической гипотезе), если вместо предположения: $R_{ij}$ – независимые случайные величины, распределенные по законам Пуассона с математическим ожиданием $M\left(R_{ij}\right) = d_{ij}$, считать, что $\left\{R_{ij}\right\} \gg 1$ имеют мультиномиальное распределение с параметрами $\left\{ d_{ij} \Big/ \sum_{i,j} d_{ij} \right\}$. Такой подход представляется более естественным, чем изложенный выше. Кроме того, поскольку константа строгой выпуклости энтропии в 1-норме равна 1 [93] (неравенство Пинскера, оценивающее снизу расстояние Кульбака–Лейблера с помощью квадрата 1-нормы), то на базе асимптотического представления логарифма функции правдоподобия в окрестности точки максимума [319], можно построить доверительный интервал для разности оценок и оцениваемых параметров в (наиболее естественной) 1-норме.

[15] При $f\left(C_{ij}\right) = \exp\left(-\beta C_{ij}\right)$ эту модель называют также энтропийной моделью [123] или моделью А.Дж. Вильсона [26, 49]. Энтропийная модель, использованная нами ранее (см. подраздел 1.1.3), для связи матрицы корреспонденций и матрицы издержек является частным случаем гравитационной модели при указанном выборе функции притяжения. Точнее, говоря, если бы мы считали, что $C_{ij} = T_{ij}$ не зависит от $\left\{d_{ij}\right\}$ и $R\left(C_{ij}\right) = f\left(C_{ij}\right)/2$, то решение задачи (1.1.2) в точности давало бы гравитационную модель. При этом $A_i$, $B_j$ выражались бы через множители Лагранжа (двойственные переменные), соответственно, $\lambda_i^L$, $\lambda_j^W$. Система (1.1.6) при этом получалась бы при подстановке решения задачи (1.1.2), зависящего от этих $2n$ неизвестных параметров, в ограничения (А) (которых тоже $2n$), точнее говоря, независимых параметров и уравнений было бы $2n$–1 [49].

Соответственно, если критерий остановки (1.1.7) не выполняется, то происходит перерасчет матрицы корреспонденций с учетом обновленной матрицы издержек.

Калибровку $\beta$ можно проводить в случае, когда имеется дополнительная информация о реальных издержках на дорогах. Для этого применяется следующий алгоритм.

Пусть $C_l^*$ – средние издержки на проезд в системе (известны, например, из опросов) для $l$-го слоя спроса [287], например, для трудовых миграций.

Алгоритм [219]:

1. Рассчитываем $\beta_l^0 = \dfrac{1}{C_l^*}$, $m = 0$;
2. Рассчитываем $\left\{d_{ij}\left(\beta_l^m\right)\right\}_{i,j=1}^{n,n}$ – матрицу корреспонденций при $\beta_l = \beta_l^m$;
3. Пересчет «4-стадийной модели»;
4. Рассчитываем $C_l^m$ – средние издержки на проезд, соответствующие матрице корреспонденций $\left\{d_{ij}\left(\beta_l^m\right)\right\}_{i,j=1}^{n,n}$;
5. При $m \geq 1$ проверяем условие: $\left|C_l^{m-1} - C_l^*\right| \leq \varepsilon$ (критерий остановки);
6. Рассчитываем $C_l^1$ по $\beta_l^1 = \dfrac{\beta_l^0 C_l^0}{C_l^*}$ (при $m = 0$) и полагаем
$$\beta_l^{m+1} = \frac{\left(C_l^* - C_l^{m-1}\right)\beta_l^m + \left(C_l^m - C_l^*\right)\beta_l^{m-1}}{C_l^m - C_l^{m-1}} \quad \text{(при } m \geq 1\text{)};$$
7. Переходим на шаг 2.

Данный алгоритм (и его вариации) дает оценку коэффициента $\beta$ при известных средних издержках в сети и является наиболее часто применяемым на практике [287].

Критерий остановки выбран именно таким по следующей причине. Если критерий остановки из подраздела 1.1.5) выполняется, то, подставляя $C_l^{m-1} \simeq C_l^*$ в формулу расчета $\beta_l^{m+1}$ на шаге 6, получаем:

$$\beta_l^{m+1} \simeq \frac{\left(C_l^* - C_l^*\right)\beta_l^m + \left(C_l^m - C_l^*\right)\beta_l^{m-1}}{C_l^m - C_l^*} = \beta_l^{m-1}.$$

Полученное значение параметра $\beta\left(C^*\right)$ «соответствует» средним издержкам $C^*$, наблюдаемым «в жизни». Сходимость алгоритма и монотонная зависимость $\beta(C)$ (следовательно, и взаимно-однозначное соответствие) между средними издержками $c$ и значением параметра $\beta$ были показаны в работе [196].

Стоит, однако, отметить, что при численной реализации алгоритма возникает множество проблем. В частности, при «плохом» выборе точки старта или при плохой калибровке параметров функций издержек возникают ситуации, в которых алгоритм требует вычисления экстремально больших чисел и превышает размер доступной памяти даже на самых современных машинах (кластерах).

Другой проблемой является отсутствие оценки скорости сходимости алгоритма. Зачастую, при моделировании критерий остановки устанавливается жестко. Например, прописывается, что алгоритм калибровки $\beta$ должен быть использован не более 20 раз (или другое, разумное, по мнению модельера, количество итераций). При этом отсутствует четкий критерий качества оценок, полученных таким образом.

Еще одной проблемой описанного подхода является высокая чувствительность к параметрам модели с одной стороны, и невозможность учесть реальные подтвержденные данные на отдельных элементах сети с другой. Другими словами, допустим нам известны величины потоков на той или иной дороге каждый день, к сожалению, данную информацию использовать в такой модели не представляется возможным.

### 1.1.6 Модель стабильной динамики

Приведем основные положения модели стабильной динамики, следуя [49, 275]. В рамках модели предполагается, что водители действуют оппортунистически, т.е. выполнен первый принцип Вардропа. Рассмотрим ориентированный граф $\Gamma(V,E)$. В модели каждому ребру $e \in E$ ставятся в соответствия параметры $\bar{f}_e$ и $\bar{t}_e$. Они имеют следующую трактовку: $\bar{f}_e$ – максимальная пропускная способность ребра $e$, $\bar{t}_e$ – минимальные временные издержки на прохождение ребра $e$. Таким образом, сама модель задается графом $\Gamma\left(V,E,\bar{f},\bar{t}\right)$, где $\bar{f} = \left\{\bar{f}_1,\dots,\bar{f}_{|E|}\right\}^T$, $\bar{t} = \left\{\bar{t}_1,\dots,\bar{t}_{|E|}\right\}^T$. Пусть $f$ – вектор распределения потоков по ребрам, инициируемый равновесным распределением потоков по маршрутам, а $t$ – вектор временных издержек, соответствующий распределению $f$. Тогда, если транспортная система находится в стабильном состоянии, всегда выполняются неравенства $f \le \bar{f}$ и $t \ge \bar{t}$. При этом считается, что, если поток по ребру $f_e$ меньше, чем максимальная пропу-скная способность ребра $\overline{f}_e$, то все автомобили в потоке двигаются с максимальной скоростью, а их временные издержки $t_e$ минимальны и равны $\overline{t}_e$. Если же поток по ребру $f_e$ становится равным пропускной способности ребра $\overline{f}_e$, то временные издержки водителей $t_e$ могут быть сколь угодно большими. Это удобно объяснить следующим образом. Допустим, на некоторое ребро $e$ стало поступать больше автомобилей, чем оно способно обслужить. Тогда на этом ребре начинает образовываться очередь (пробка). Временные издержки на прохождение ребра $t_e$ складываются из минимальных временных издержек $\overline{t}_e$ и времени, которое водитель вынужден отстоять в пробке. При этом, очевидно, если входящий поток автомобилей на ребро $e$ не снизится до максимально допустимого уровня (пропускной способности ребра), то очередь будет продолжать расти и система не будет находиться в стабильном состоянии. Если же в какой-то момент входящий на ребро $e$ поток снизится до уровня пропускной способности ребра, то в системе наступит равновесие. При этом пробка на ребре $e$ (если входящий поток $f_e$ будет равен $\overline{f}_e$) не будет рассасываться, т.е. временные издержки так и останутся на уровне $t_e$ ($t_e > \overline{t}_e$). Рассмотрим это на примере из статьи [275].

**Пример 1.1.1.** Рассмотрим (в рамках модели стабильной динамики) граф $\Gamma\left(V,E,\overline{t},\overline{f}\right)$ (см. Рисунок 1.1.2). Пункты 1 и 2 – потокообразующая пара. При этом выполнено: $d_{12}$ – поток из 1 в 2, $\overline{t}_{up} < \overline{t}_{down}$. Если выполнено $d_{12} < \overline{f}_{up}$, то все водители будут использовать ребро up, причем пробка образовываться не будет, так как пропускная способность ребра больше, чем количество желающих проехать (в единицу времени) водителей. В момент, когда $d_{12} = \overline{f}_{up}$ возможности ребра up будут использоваться на пределе. Если же в какой-то момент величина $f$ станет больше $\overline{f}_{up}$, то на ребре *up* начнет образовываться пробка. Пробка будет расти до тех пор, пока издержки от использования маршрута *up* с не сравняются с издержками от использования маршрута *down*. В этот момент оставшаяся часть начнет использовать маршрут *down*. Если же корреспонденция из 1 в 2 превысит суммарную пропускную способность ребер *up* и *down*, то пробки будут расти неограниченно (входящий поток на ребро будет больше, чем исходящий, соответственно количество автомобилей в очереди будет расти постоянно), т.е. стабильное распределение в системе никогда не установится. Более подробно рассмотрение данной задачи (и модели) стоит смотреть в работе [275]. ■

Рассмотрим другой модельный пример из статьи [275].

**Пример 1.1.2 (Парадокс Браесса).** Задан граф $\Gamma\left(V,E,\bar{t},\bar{f}\right)$ (см. Рисунок 1.1.3). При этом выполнено: $\bar{t}_{13}=1$ час, $\bar{t}_{12}=15$ минут, $\bar{t}_{23}=30$ минут; (1,3) и (2,3) – потокообразующие пары, $d_{13}=1500$ авт/час и $d_{23}=1500$ авт/час – соответствующие корреспонденции. Пусть пропускные способности всех ребер одинаковы и равны 2000 авт/час. Тогда равновесие будет такое: 500 авт/час из 1 будут направляться в 3 через 2, а 1000 авт/час будут ехать напрямую (для выезжающих из 2 никаких альтернатив нет). При этом все водители, выезжающие из 1, потратят 1 час, а водители, выезжающие из 2, потратят 45 минут. Таким образом, водителям, едущим из 2 в 3 выгодно, чтобы ребро 1–2 отсутствовало, в то время как для водителей, которые едут из 1 в 3, наличие ребра 1–2 безразлично. Т.е., другими словами, если бы мы имели власть запретить проезд по ребру 1–2, то часть водителей выиграла бы от такого решения, и никто бы не проиграл. Интересно заметить, что для модели Бэкмана это не выполняется. Действительно, в модели Бэкмана, как и для модели стабильной динамики, издержки для водителей, следующих из 1 в 3 равны издержкам на ребре 1–3. Однако они монотонно возрастают от потока на данном ребре. Если бы мы запретили проезд по ребру 1–2, то поток на 1–3 увеличился бы, следовательно, возросли бы и издержки на ребре 1–3. Другими словами, улучшение ситуации для водителей, следующих из 2 в 3, привело бы к ухудшению ситуации для водителей, следующих из 1 в 3. ■

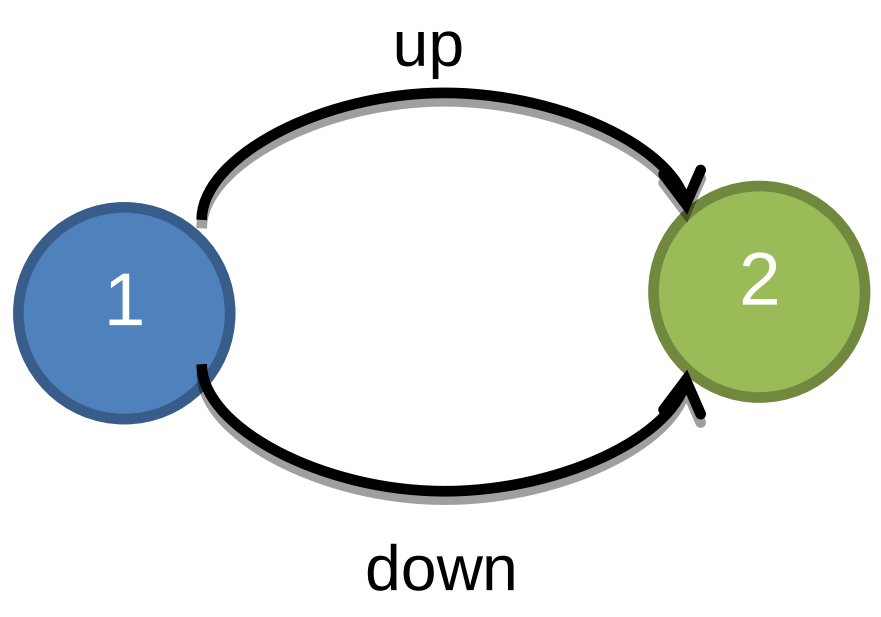

Рисунок 1.1.2

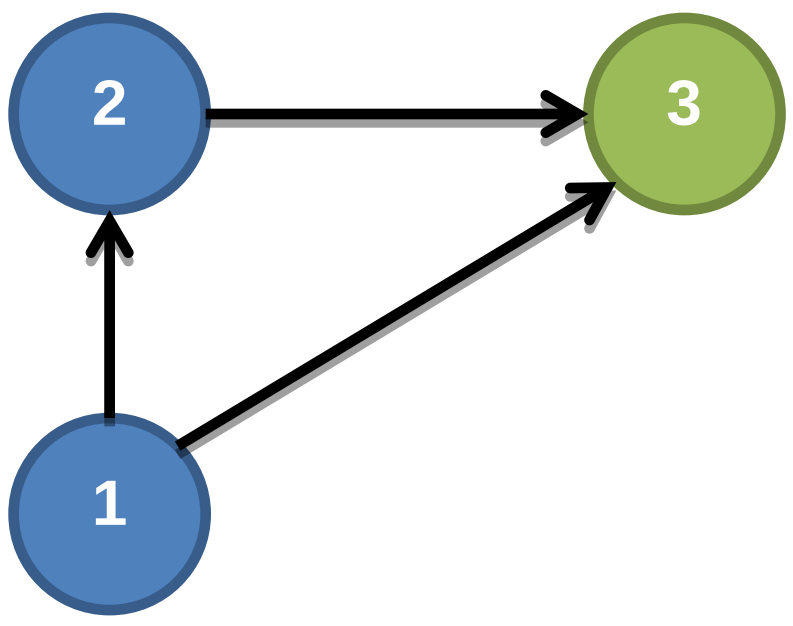

Рисунок 1.1.3

Введем ряд новых обозначений. Пусть $w=(i,j)$ – потокообразующая пара (источник – сток) графа $\Gamma(V,E,\bar{t},\bar{f})$, $P_w$ – множество соответствующих $w$ маршрутов, а $t$ – установившийся на графе вектор временных издержек. Тогда временные издержки, соответствующие самому "быстрому" (наикратчайшему) маршруту из $P_w$ равны: $T_w(t)=\min_{p\in P_w}\sum_{e\in E}\delta_{ep}t_e$. Функция $T_w(t)$ зависит от вектора временных издержек $t$. Важно заметить, что функция $T_w(t)$ и её супердифференциал эффективно вычисляются, например, алгоритмом Дейкстры [79, 132] или более быстрыми современными методами [325].

Пусть корреспонденция потокообразующей пары $w$ равна $d_w$. Тогда, если $f^*$ является равновесным распределением потоков для заданного графа $\Gamma(V,E,\bar{t},\bar{f})$, а $t^*$ – соответствующий равновесному распределению вектор временных издержек, то $f^*=\sum_{w\in W}f_w^*$, где $f_w^*$ – вектор распределения потока, порождаемого потокообразующей парой $w$ ($w\in OD$). При этом $f_w^*$ удовлетворяет

$$f_w^*\in d_w\partial T_w(t),$$

где $\partial T_w(t)$ – супердифференциал, который можно посчитать с помощью теоремы Милютина–Дубовицкого [85] используя структуру функции $T_w(t)$. Следовательно, по теореме Моро–Рокафеллара [85] $f^*\in\sum_{w\in W}d_w\partial T_w(t)$. Остается вопрос: как искать равновесный вектор временных издержек $t^*$.

**Теорема 1.1.2 [275].** *Распределение потоков $f^*$ и вектор временных издержек $t^*$ являются равновесными для графа $\Gamma(V,E,\bar{t},\bar{f})$, заданного множества потокообразующих пар $W$ и соответствующих им потоков $d_w$, тогда и только тогда, когда $t^*$ является решением следующий задачи вогнутой оптимизации (максимизация количества свободно движущихся автомобилей, то есть автомобилей не стоящих в пробках):*

$$\max_{t\geq\bar{t}}\left\{\sum_{w\in W}d_wT_w(t)-\left\langle\bar{f},t-\bar{t}\right\rangle\right\}, \qquad (1.1.8)$$

*где $f^*=\bar{f}-s^*$, где $s^*$ – (оптимальный) вектор двойственных множителей для ограничений $t\geq\bar{t}$ в задаче (1.1.8) или, другими словами, решение двойственной задачи.*

**Схема доказательства.** Резюмируем далее основные постулаты модели стабильной динамики:

1. $t_e \ge \bar{t}_e$, $f_e \le \bar{f}_e$, $e \in E$;
2. $(t_e - \bar{t}_e)(\bar{f}_e - f_e) = 0$, $e \in E$;
3. $f \in \partial \sum_{w \in W} d_w T_w(t)$ (принцип (Нэша–)Вардропа).

Наша цель, подобрать такую задачу выпуклой (вогнутой) оптимизации, решение которой давало бы пару $(t, f)$, удовлетворяющую этим постулатам. С задачей выпуклой оптимизации намного удобнее работать (разработаны эффективные способы решения таких задач [93]), чем с описанием 1 – 3. Если бы мы могли ограничиться только п. 3, то тогда этот пункт можно было бы переписать

$$\max_{t \ge 0} \left\{ \sum_{w \in W} d_w T_w(t) - \langle f, t - \bar{t} \rangle \right\}, \qquad (1.1.9)$$

но есть еще пп. 1, 2. Подсказка содержится в п. 2, который удобно понимать как условие дополняющей нежесткости [85] ограничения $t \ge \bar{t}$ (в задаче (1.1.9)) с множителями Лагранжа $s = \bar{f} - f \ge 0$. Причем, нет необходимости оговаривать дополнительно, что $f \le \bar{f}$, поскольку множители Лагранжа к ограничениям вида неравенство автоматически неотрицательны. Таким образом, мы приходим к (1.1.8). Все же необходимо оговориться, что получившаяся задача (1.1.8) хотя и является довольно простой задачей (не сложно заметить, что это просто задача линейного программирования (ЛП), записанная в более компактной форме, действительно: $\max_{t \ge 0;\, \{T_w\};\, C} \left\{ C - \langle \bar{f}, t - \bar{t} \rangle;\, T_w \le \sum_{e \in E} \delta_{ep} t_e,\ p \in P_w,\ w \in W;\, C = \sum_{w \in W} d_w T_w \right\}$), но требуется также найти двойственные множители (чтобы определить $f$). □

Задача (1.1.8) имеет (конечное) решение тогда и только тогда, когда существует хотя бы один вектор потоков $f \le \bar{f}$, соответствующий заданным корреспонденциям. Другими словами, существует способ так организовать движение, чтобы при заданных корреспонденциях не нарушались условия $f \le \bar{f}$. В противном случае, равновесного (стационарного) режима в системе не будет, и со временем весь граф превратится в одну большую пробку, которая начнет нарушать условия заданных корреспонденций, не позволяя новым пользователям приходить в сеть с той интенсивностью, с которой они этого хотят.

В отличие от модели Бэкмана, в которой при монотонно возрастающих функциях затрат на ребрах равновесие (по потокам на ребрах) единственно, в модели стабильной динамики это может быть не так. Хотя типичной ситуацией (общим положением) будет единственность равновесия, поскольку поиск стабильной конфигурации эквивалентен решению задачи ЛП, все же в определенных “вырожденных” ситуациях может возникать неединственность. Скажем, если в примере 1.1.1 (см. Рисунок 1.1.2) $\overline{f}_{up} = d_{12}$. Возникает вопрос: как выбрать единственное равновесие, то которое реализуется? Один из ответов имеется в [275]: следить за динамикой во времени $d_{12}(t)$, а именно за

$$\int_0^T \left(d_{12}(t) - \overline{f}_{up}\right) dt < \infty .$$

Другой ответ связан с прямым осуществлением для данного примера предельного перехода, описанного в следующем пункте.

По сути, с помощью этого и следующего пункта мы установим, что если

$$\tau_e^\mu(f_e) \xrightarrow[\mu \to 0+]{} \begin{cases} \overline{t}_e, & 0 \le f_e < \overline{f}_e \\ [\overline{t}_e, \infty), & f_e = \overline{f}_e \end{cases},$$

$$d\tau_e^\mu(f_e)/df_e \xrightarrow[\mu \to 0+]{} 0, \ 0 \le f_e < \overline{f}_e,$$

$x(\mu)$ – равновесное распределение потоков по путям в модели Бэкмана при функциях затрат на ребрах $\tau_e^\mu(f_e)$, то

$$\tau_e^\mu\left(f_e\left(x(\mu)\right)\right) \xrightarrow[\mu \to 0+]{} t_e,$$

$$f_e\left(x(\mu)\right) \xrightarrow[\mu \to 0+]{} f_e,$$

где пара $(t, f)$ – равновесие в модели стабильной динамики с тем же графом и матрицей корреспонденций, что и в модели Бэкмана, и с ребрами, характеризующимися набором $(\overline{t}, \overline{f})$ из определения функций $\tau_e^\mu(f_e)$.

Численно решать задачу (1.1.1) можно различными способами (см., например, [40]).

### 1.1.7 Связь модели стабильной динамики с моделью Бэкмана

Получим модель стабильной динамики в результате предельного перехода из модели Бэкмана. Для этого будем считать, что $\tau_e(f_e) = \overline{t}_e \cdot \left(1 - \mu \ln\left(1 - f_e / \overline{f}_e\right)\right)$.

Перепишем исходную задачу поиска равновесного распределения потоков (здесь используются обозначения: $\sigma_e^*(t_e)$ – сопряженная функция к $\sigma_e(f_e)$):

$$\min_{f,x}\left\{\sum_{e\in E}\sigma_e(f_e):\ f = \Theta x,\ x \in X\right\} = \min_{f,x}\left\{\sum_{e\in E}\max_{t_e \in \mathrm{dom}\,\sigma_e^*}\left[f_e t_e - \sigma_e^*(t_e)\right]:\ f = \Theta x,\ x \in X\right\} =$$

$$= \max_{t \in \mathrm{dom}\,\sigma^*}\left\{\min_{f,x}\left[\sum_{e\in E} f_e t_e :\ f = \Theta x,\ x\in X\right] - \sum_{e\in E}\sigma_e^*(t_e)\right\}.$$

Найдем $\sum_{e\in E}\sigma_e^*(t_e)$:

$$\sigma_e^*(t_e)=\sup_{f_e}\left(t_e\cdot f_e-\int_0^{f_e}\tau(z)\,dz\right)=\sup_{f_e}\left(t_e\cdot f_e-\int_0^{f_e}\bar{t}_e\cdot\left(1-\mu\ln\left(1-z/\bar{f}_e\right)\right)dz\right)=$$

$$=\sup_{f_e}\left((t_e-\bar{t}_e)\,f_e+\bar{t}_e\mu\int_0^{f_e}\ln\left(1-z/\bar{f}_e\right)dz\right).$$

Из принципа Ферма [85] найдем $f_e$, доставляющее максимум:

$$\frac{\partial}{\partial f_e}\left((t_e-\bar{t}_e)\,f_e+\bar{t}_e\mu\int_0^{f_e}\ln\left(1-z/\bar{f}_e\right)dz\right)=0\ \Rightarrow$$

$$\exp\left(-\frac{t_e-\bar{t}_e}{\bar{t}_e\mu}\right)=1-f_e/\bar{f}_e\ \Rightarrow$$

$$f_e=\bar{f}_e\cdot\left(1-\exp\left(-\frac{t_e-\bar{t}_e}{\bar{t}_e\mu}\right)\right).$$

Подставляем в $\sigma_e^*(t_e)$:

$$(t_e-\bar{t}_e)\,f_e+\bar{t}_e\mu\int_0^{f_e}\ln\left(1-z/\bar{f}_e\right)dz=(t_e-\bar{t}_e)\,f_e-\bar{t}_e\cdot\bar{f}_e\cdot\mu\int_1^{1-f_e/\bar{f}_e}\ln z\,dz=$$

$$=(t_e-\bar{t}_e)\,f_e-\bar{t}_e\cdot\bar{f}_e\cdot\mu\left(\left(1-\frac{f_e}{\bar{f}_e}\right)\ln\left(1-\frac{f_e}{\bar{f}_e}\right)-1+\frac{f_e}{\bar{f}_e}+1\right)=$$

$$=(t_e-\bar{t}_e)\,f_e-\left(-\left(\bar{f}_e-f_e\right)(t_e-\bar{t}_e)+f_e\cdot\bar{t}_e\cdot\mu\right)=$$

$$=(t_e-\bar{t}_e)\,\bar{f}_e-f_e\cdot\bar{t}_e\cdot\mu=(t_e-\bar{t}_e)\,\bar{f}_e-\bar{f}_e\cdot\bar{t}_e\cdot\mu\left(1-\exp\left(-\frac{t_e-\bar{t}_e}{\bar{t}_e\mu}\right)\right).$$

Возвращаясь к исходной задаче имеем:

$$\max_{t\in\operatorname{dom}\sigma^*}\left\{\min_{f,x}\left[\sum_{e\in E}f_et_e:\ f=\Theta x,\ x\in X\right]-\sum_{e\in E}\sigma_e^*(t_e)\right\}=$$

$$=\max_{t\in\operatorname{dom}\sigma^*}\left\{\sum_{w\in W}d_wT_w(t)-\left\langle\bar{f},t-\bar{t}\right\rangle+\mu\sum_{e\in E}\bar{f}_e\cdot\bar{t}_e\left(1-\exp\left(-\frac{t_e-\bar{t}_e}{\bar{t}_e\mu}\right)\right)\right\}\overset{\mu\to0+}{=}$$

$$\overset{\mu\to0+}{=}\max_{t\ge\bar{t}}\left\{\sum_{w\in W}d_wT_w(t)-\left\langle\bar{f},t-\bar{t}\right\rangle\right\}.$$

Возможность поменять местами порядок взятия максимума и минимума следует из теоремы фон Неймана о минимаксе [119].

Выбор именно функции $\tau_e(f_e)=\bar{t}_e\cdot\left(1-\mu\ln\left(1-f_e/\bar{f}_e\right)\right)$ не является определяющим. Вместо этой функции может быть взят любой другой гладкий внутренний барьер области $f_e < \overline{f}_e$, например, $\tau_e\left(f_e\right) = \overline{t}_e \cdot \left(1 + \mu \frac{\overline{f}_e}{\overline{f}_e - f_e}\right)$ – такую функцию вполне естественно использовать при поиске равновесного распределения пользователей сети общественного транспорта [275]. Найдем $\sigma_e^*\left(t_e\right)$:

$$\sigma_e^*\left(t_e\right) = \sup_{f_e}\left(t_e \cdot f_e - \int_0^{f_e} \overline{t}_e \cdot \left(1 + \mu \frac{z}{\overline{f}_e - z}\right) dz\right) = \sup_{f_e}\left(\left(t_e - \overline{t}_e\right) \cdot f_e + \overline{t}_e \cdot \overline{f}_e \cdot \mu \cdot \ln\left(1 - \frac{f_e}{\overline{f}_e}\right) + \overline{t}_e \cdot f_e \cdot \mu\right).$$

Вновь выписывая условия оптимальности первого порядка имеем:

$$\frac{\partial}{\partial f_e}\left(\left(t_e - \overline{t}_e\right) \cdot f_e + \overline{t}_e \cdot \overline{f}_e \cdot \mu \cdot \ln\left(1 - \frac{f_e}{\overline{f}_e}\right) + \overline{t}_e \cdot f_e \cdot \mu\right) = \left(t_e - \overline{t}_e\right) - \overline{t}_e \cdot \mu \cdot \frac{\overline{f}_e}{\overline{f}_e - f_e} + \overline{t}_e \cdot \mu = 0 \Rightarrow$$

$$f_e = \overline{f}_e \cdot \left(1 - \frac{\overline{t}_e \cdot \mu}{t_e - (1-\mu)\overline{t}_e}\right).$$

Подставляя в $\sigma_e^*\left(t_e\right)$, имеем:

$$\sigma_e^*\left(t_e\right) = \left(t_e - \overline{t}_e\right) \cdot \overline{f}_e \cdot \left(1 - \frac{\overline{t}_e \cdot \mu}{t_e - (1-\mu)\overline{t}_e}\right) + \overline{t}_e \cdot \overline{f}_e \cdot \mu \cdot \ln\left(\frac{\overline{t}_e \cdot \mu}{t_e - (1-\mu)\overline{t}_e}\right) + \overline{t}_e \cdot \mu \cdot \overline{f}_e \cdot \left(1 - \frac{\overline{t}_e \cdot \mu}{t_e - (1-\mu)\overline{t}_e}\right) =$$

$$= \left(t_e - \overline{t}_e\right) \cdot \overline{f}_e + \overline{t}_e \cdot \overline{f}_e \cdot \mu \cdot \ln\left(\frac{\overline{t}_e \cdot \mu}{t_e - (1-\mu)\overline{t}_e}\right).$$

В итоге получаем:

$$\max_{t \in \operatorname{dom} \sigma^*}\left\{\min_{f,x}\left[\sum_{e \in E} f_e t_e :\ f = \Theta x,\ x \in X\right] - \sum_{e \in E} \sigma_e^*\left(t_e\right)\right\} =$$

$$= \max_{t \in \operatorname{dom} \sigma^*}\left\{\sum_{w \in W} d_w T_w(t) - \left\langle \overline{f}, t - \overline{t}\right\rangle + \mu \sum_{e \in E} \overline{f}_e \cdot \overline{t}_e \cdot \ln\left(1 + \frac{t_e - \overline{t}_e}{\overline{t}_e \cdot \mu}\right)\right\} \overset{\mu \to 0+}{=} \max_{t \geq \overline{t}}\left\{\sum_{w \in W} d_w T_w(t) - \left\langle \overline{f}, t - \overline{t}\right\rangle\right\}.$$

Наконец, рассмотрим третий вариант выбора функции $\tau_e\left(f_e\right)$:

$$\tau_e\left(f_e\right) = \overline{t}_e \cdot \left(1 + \gamma \cdot \left(\frac{f_e}{\overline{f}_e}\right)^{\frac{1}{\mu}}\right).$$

Функции такого вида называются BPR-функциями и наиболее часто применяются при моделировании. Так, при использовании модели Бэкмана обычно полагают $\mu = 0.25$, а значение параметра $\gamma$ варьируется от 0.15 до 2 и определяется типом дороги [287]. Найдем $\sigma_e^*\left(t_e\right)$:

$$\sigma_e^*\left(t_e\right) = \sup_{f_e}\left(t_e \cdot f_e - \int_0^{f_e} \overline{t}_e \cdot \left(1 + \gamma \cdot \left(\frac{z}{\overline{f}_e}\right)^{\frac{1}{\mu}}\right) dz\right) = \sup_{f_e}\left(\left(t_e - \overline{t}_e\right) \cdot f_e - \overline{t}_e \cdot \gamma \cdot \int_0^{f_e}\left(\frac{z}{\overline{f}_e}\right)^{\frac{1}{\mu}} dz\right) =$$

$$= \sup_{f_e} \left( \left(t_e - \bar{t}_e\right) \cdot f_e - \bar{t}_e \cdot \frac{\mu}{1+\mu} \cdot \gamma \cdot \frac{f_e^{1+\frac{1}{\mu}}}{\bar{f}_e^{\frac{1}{\mu}}} \right).$$

Аналогично рассуждая:

$$\frac{\partial}{\partial f_e} \left( \left(t_e - \bar{t}_e\right) \cdot f_e - \bar{t}_e \cdot \frac{\mu}{1+\mu} \cdot \gamma \cdot \frac{f_e^{1+\frac{1}{\mu}}}{\bar{f}_e^{\frac{1}{\mu}}} \right) = t_e - \bar{t}_e - \bar{t}_e \cdot \gamma \cdot \frac{f_e^{\frac{1}{\mu}}}{\bar{f}_e^{\frac{1}{\mu}}} = 0 \;\Rightarrow\; f_e = \bar{f}_e \cdot \left( \frac{t_e - \bar{t}_e}{\bar{t}_e \cdot \gamma} \right)^{\mu}.$$

Тогда

$$\sigma_e^*\left(t_e\right) = \sup_{f_e} \left( \left(t_e - \bar{t}_e\right) \cdot f_e - \bar{t}_e \cdot \frac{\mu}{1+\mu} \cdot \gamma \cdot \frac{f_e^{1+\frac{1}{\mu}}}{\bar{f}_e^{\frac{1}{\mu}}} \right) = \bar{f}_e \cdot \left( \frac{t_e - \bar{t}_e}{\bar{t}_e \cdot \gamma} \right)^{\mu} \frac{\left(t_e - \bar{t}_e\right)}{1+\mu}.$$

В итоге получаем:

$$\max_{t \in \mathrm{dom}\, \sigma^*} \left\{ \sum_{w \in W} d_w T_w\left(t\right) - \sum_{e \in E} \bar{f}_e \cdot \left( \frac{t_e - \bar{t}_e}{\bar{t}_e \cdot \gamma} \right)^{\mu} \frac{\left(t_e - \bar{t}_e\right)}{1+\mu} \right\} \overset{\mu \to 0+}{=} \max_{t \ge \bar{t}} \left\{ \sum_{w \in W} d_w T_w\left(t\right) - \left\langle \bar{f}, t - \bar{t} \right\rangle \right\}.$$

**Замечание 1.1.3.** Естественно в контексте всего написанного выше теперь задаться вопросом: а можно ли получить модель стабильной динамики эволюционным образом, то есть подобно тому, как в конце подраздела 1.1.4 была эволюционно проинтерпретирована модель Бэкмана. Действительно, рассмотрим логит динамику (с $T \to 0+$) или, скажем, просто имитационную логит динамику [302]. Предположим, что гладкие, возрастающие, выпуклые функции затрат на всех ребрах $\tau_e\left(f_e\right)$ при $\mu \to 0+$ превращаются в “ступеньки” (многозначные функции)

$$\tau_e\left(f_e\right) = \begin{cases} \bar{t}_e, & 0 \le f_e < \bar{f}_e \\ \left[\bar{t}_e, \infty\right), & f_e = \bar{f}_e \end{cases}.$$

Согласно подразделу 1.1.4 равновесная конфигурация при таком переходе $\mu \to 0+$ должна находиться из решения задачи

$$\sum_{e \in E} \int_0^{f_e} \tau_e\left(z\right) dz \to \min_{f = \Theta x,\, x \in X}.$$

Считая, что в равновесии не может быть $\tau_e\left(f_e\right) = \infty$ (иначе, равновесие просто не достижимо, и со временем весь граф превратится в одну большую пробку, см. конец подраздела 1.1.6), можно не учитывать в интеграле вклад точек $\bar{f}_e$ (в случае попадания в промежуток интегрирования), то есть переписать задачу следующим образом

$$\min_{f = \Theta x,\, x \in X} \sum_{e \in E} \int_0^{f_e} \left( \bar{t}_e + \delta_{\bar{f}_e}\left(z\right) \right) dz \Leftrightarrow \min_{\substack{f = \Theta x,\, x \in X \\ f \le \bar{f}}} \sum_{e \in E} f_e \bar{t}_e,$$

где

$$\delta_{\bar{f}_e}(z) = \begin{cases} 0, & 0 \le z < \bar{f}_e \\ \infty, & z \ge \bar{f}_e \end{cases}.$$

Мы получили задачу линейного программирования. Интересно заметить, что двойственные множители $\eta \ge 0$ к ограничениям $f \le \bar{f}$ имеют размерность (“физический смысл”) времени (см. интерпретацию двойственных множителей из п. 3). Чему равно это время сразу может быть не ясно из решения только что выписанной задачи линейного программирования. Но если перейти к двойственной задаче, то мы получим уже рассматриваемую нами ранее задачу (1.1.9) из подраздела 1.1.6. Причем двойственные множители $\eta$ связаны с временами проезда ребер $t$ следующим образом: $\eta = t - \bar{t} \ge 0$, то есть, действительно, получают содержательную интерпретацию времени, потерянного на ребрах (дополнительно к времени проезда по свободной дороге) из-за наличия пробок. Еще одним подтверждением только что сказанному является то, что если $f_e = \bar{f}_e$, то на ребре $e$ пробка и время прохождения этого ребра $t_e \ge \bar{t}_e$ ($\eta \ge 0$), а если пробки нет, то $f_e < \bar{f}_e$ и $t_e = \bar{t}_e$ ($\eta = 0$). То есть имеют место условия дополняющей нежесткости.

Из данного замечания становится ясно, как оптимально назначать платы за проезд в модели стабильной динамики. Переведем предварительно время в деньги. Назначая платы за проезд по ребрам графа в соответствии с вектором $\eta$, получим, что равновесное распределение пользователей транспортной сети по этой сети будет описываться парой $(f, \bar{t})$, что, очевидно, лучше, чем до введения плат $(f, t)$. Действительно, задачу поиска равновесных потоков $f$:

$$\min_{\substack{f = \Theta x,\, x \in X \\ f \le \bar{f}}} \sum_{e \in E} f_e \bar{t}_e$$

можно переписать с помощью метода множителей Лагранжа, следующим образом

$$\min_{f = \Theta x,\, x \in X} \sum_{e \in E} \left( f_e \bar{t}_e + \eta_e \cdot \left( f_e - \bar{f}_e \right) \right).$$

Другими словами, у этих двух задач одинаковые решения $f$. Но $t = \bar{t} + \eta$, поэтому тот же самый вектор $f$ будет доставлять решение задаче

$$\min_{f = \Theta x,\, x \in X} \sum_{e \in E} f_e t_e .$$

Поскольку мы знаем, что $f \le \bar{f}$, то $f$ будет доставлять решение также и этой задаче

$$\min_{\substack{f = \Theta x,\, x \in X \\ f \le \bar{f}}} \sum_{e \in E} f_e t_e .$$

Таким образом, один и тот же вектор $f$ отвечает оптимальному (с точки зрения социального оптимума) распределению потоков в графе с ребрами, взвешенными согласно векторам $\bar{t}$ и $t$. Более того, не сложно понять, что все это остается верным не только для двух векторов $\bar{t}$ и $t$, но и для целого “отрезка” векторов: $\tilde{t} = \alpha\bar{t} + (1-\alpha)t$, $\alpha \in [0,1]$. К сожалению, если не взимать платы за проезд, то всегда реализуется сценарий $\alpha = 0$, то есть вектор затрат будет $t$ (объяснению причин были посвящены подразделы 1.1.6 и 1.1.7), но если мы взимаем платы согласно вектору $\eta_\alpha = \alpha \cdot (t - \bar{t})$, то вектор реальных затрат (не учитывающих платы) будет $t - \eta_\alpha = (1-\alpha)t + \alpha\bar{t}$. Следовательно, оптимально выбирать $\alpha = 1$, то есть платы согласно $\eta = t - \bar{t}$.

Собственно, всю эту процедуру (назначения плат) можно делать адаптивно: постепенно увеличивая плату за проезд на тех ребрах, на которых наблюдаются пробки. И делать это нужно до тех пор, пока пробки не перестанут появляться.

### 1.1.8 Учет расщепления потоков по способам передвижений

Распространим модель стабильной динамики на тот случай, когда (все) пользователи (игроки) транспортной сети имеют возможность выбирать между двумя альтернативными видами транспорта: личным и общественным [14]. Соответственно, теперь у нас имеется информация о $\Gamma_л\left(V, E^л, \bar{t}^л, \bar{f}^л\right)$ – для личного транспорта и $\Gamma_о\left(V, E^о, \bar{t}^о, \bar{f}^о\right)$ – для общественного транспорта. Аналогично рассматривается общий случай. Имеет место следующий результат, получаемый аналогично [275].

**Теорема 1.1.3**. *Распределение потоков $f^* = \left(f^л, f^о\right)$ и вектор временных издержек $t^* = \left(t^л, t^о\right)$ являются равновесными для графа $\Gamma\left(V, E, \bar{t}, \bar{f}\right)$, заданного множества потокообразующих пар $W$ и соответствующих им потоков $d_w$, тогда и только тогда, когда $t^*$ является решением следующий задачи выпуклой оптимизации (в дальнейшем для нас будет удобно писать эту задачу как задачу минимизации, а не максимизации, как раньше):*

$$\min_{\substack{t^л \geq \bar{t}^л \\ t^о \geq \bar{t}^о}} \left\{ \left\langle \bar{f}^л, t^л - \bar{t}^л \right\rangle + \left\langle \bar{f}^о, t^о - \bar{t}^о \right\rangle - \sum_{w \in W} d_w \min\left\{ T_w^л\left(t^л\right), T_w^о\left(t^о\right) \right\} \right\}, \tag{1.1.10}$$

*где $f^л = \bar{f}^л - s^л$, $f^о = \bar{f}^о - s^о$, а $s^л$, $s^о$ – (оптимальные) векторы двойственных множителей для ограничений $t^л \geq \bar{t}^л$, $t^о \geq \bar{t}^о$ в задаче (1.1.10) или, другими словами, решение двойственной задачи.*

Если мы хотим учитывать возможность пересаживания пользователей сети в пути с личного транспорта на общественный (и наоборот), то выражение $\min\left\{T_w^л\left(t^л\right), T_w^o\left(t^o\right)\right\}$ нужно будет немного изменить. Мы не будем здесь вдаваться в детали, скажем лишь, что с точки зрении всего дальнейшего это не принципиально. Более того, хотя на время работы алгоритма это и скажется (время увеличится), но не критическим образом, поскольку решение задачи о кратчайшем маршруте, которое возникает на каждом шаге субградиентного спуска, естественным образом (небольшим раздутием графа) обобщается на случай когда есть несколько весов ребер, отвечающих разным типам транспортных средств, и в вершинах графа есть затраты на пересадку (изменение транспортного средства), при условии наличия возможности её осуществления.

Отметим, что единого подхода для моделирования общественного транспорта пока не существует. Указанный нами метод соответствует подходу, когда предполагается, что пассажиры выбирают только “оптимальные стратегии”, т.е. соответствует концепции равновесия по Нэшу. Другой, альтернативный подход, который сейчас чаще используется на практике, следует концепции «стохастического равновесия». В нем предполагается, что пассажиры выбирают каждый из возможных маршрутов с некоторой вероятностью, зависящей от ожидаемых издержек в сети и ожидаемого времени ожидания соответствующего маршрута на остановочной станции. Об этом немного подробнее будет рассказано в подразделе 1.1.10.

Обратим внимание на то, что общественный транспорт все же более естественно описывать моделью Бэкмана, а не моделью стабильной динамики. Во всяком случае, к такому выводу склоняются авторы модели стабильной динамики [275]. В этой связи далее приводится смешанная модель, в которой личный транспорт описывается моделью стабильной динамики, а общественный – моделью Бэкмана (см. подраздел 1.1.7):

$$\min_{\substack{t^л \ge \bar{t}^л \\ t^o \ge \bar{t}^o\cdot(1-\mu)}} \left\{ \left\langle \bar{f}^л, t^л - \bar{t}^л \right\rangle + \left\langle \bar{f}^o, t^o - \bar{t}^o \right\rangle - \mu \sum_{e\in E} \bar{f}_e^o \cdot \bar{t}_e^o \cdot \ln\left(1 + \frac{t_e^o - \bar{t}_e^o}{\bar{t}_e^o \cdot \mu}\right) - \sum_{w\in W} d_w \min\left\{T_w^л\left(t^л\right), T_w^o\left(t^o\right)\right\} \right\},$$

где $f^л = \bar{f}^л - s^л$, $f_e^o = \bar{f}_e^o \cdot \left(1 - \frac{\bar{t}_e^o \cdot \mu}{t_e^o - (1-\mu)\bar{t}_e^o}\right)$, $s^л$ – (оптимальный) вектор двойственных множителей для ограничений $t^л \ge \bar{t}^л$.

Несложно распространить все проводимые в следующих пунктах рассуждения именно на такую версию модели стабильной динамики с учетом расщепления по типам передвижений. Однако мы не будем здесь этого делать. Укажем лишь, что для разработки эффективного способа решения возникающих задач оптимизации потребуется привлекать прямо-двойственные субградиентные спуски для задач композитной оптимизации [4]. Подробнее об этом будет написано в разделе 3.3 главы 3.

### 1.1.9 Трехстадийная модель стабильной динамики

Объединим теперь задачи (1.1.2) и (1.1.10) в одну задачу:[16]

$$\min_{\lambda^L,\lambda^W} \max_{\sum_{i,j=1}^n d_{ij}=1,\, d_{ij}\ge 0} \left[ -\sum_{i,j=1}^n d_{ij}\ln d_{ij} + \sum_{i=1}^n \lambda_i^L \left( l_i - \sum_{j=1}^n d_{ij} \right) + \sum_{j=1}^n \lambda_j^W \left( \mathrm{w}_j - \sum_{i=1}^n d_{ij} \right) + \right.$$

$$\left. +\beta \min_{\substack{t^л \ge \bar{t}^л \\ t^o \ge \bar{t}^o}} \left\{ \left\langle \bar{f}^л, t^л - \bar{t}^л \right\rangle + \left\langle \bar{f}^o, t^o - \bar{t}^o \right\rangle - \sum_{i,j=1}^n d_{ij} \min\left\{ T_{ij}^л\left(t^л\right), T_{ij}^o\left(t^o\right) \right\} \right\} \right] =$$

$$= \min_{\substack{t^л \ge \bar{t}^л \\ t^o \ge \bar{t}^o \\ \lambda^L,\lambda^W}} \left\{ \max_{\sum_{i,j=1}^n d_{ij}=1,\, d_{ij}\ge 0} \left[ -\sum_{i,j=1}^n d_{ij}\ln d_{ij} - \beta \sum_{i,j=1}^n d_{ij} \min\left\{ T_{ij}^л\left(t^л\right), T_{ij}^o\left(t^o\right) \right\} - \right. \right.$$

$$\left. \left. -\sum_{i,j=1}^n d_{ij} \cdot \left( \lambda_i^L + \lambda_j^W \right) \right] + \sum_{i=1}^n \lambda_i^L l_i + \sum_{j=1}^n \lambda_j^W \mathrm{w}_j + \beta \left\langle \bar{f}^л, t^л - \bar{t}^л \right\rangle + \beta \left\langle \bar{f}^o, t^o - \bar{t}^o \right\rangle \right\} =$$

$$= \min_{\substack{t^л \ge \bar{t}^л \\ t^o \ge \bar{t}^o \\ \lambda^L,\lambda^W}} \left\{ \ln \left( \sum_{i,j=1}^n \exp\left( -\beta \min\left\{ T_{ij}^л\left(t^л\right), T_{ij}^o\left(t^o\right) \right\} - \lambda_i^L - \lambda_j^W \right) \right) + \right.$$

$$\left. + \sum_{i=1}^n \lambda_i^L l_i + \sum_{j=1}^n \lambda_j^W \mathrm{w}_j + \beta \left\langle \bar{f}^л, t^л - \bar{t}^л \right\rangle + \beta \left\langle \bar{f}^o, t^o - \bar{t}^o \right\rangle \right\}, \qquad (1.1.11)$$

причем $d_{ij} = \tilde{Z}^{-1} \exp\left( -\beta \min\left\{ T_{ij}^л\left(t^л\right), T_{ij}^o\left(t^o\right) \right\} - \lambda_i^L - \lambda_j^W \right)$, где $\tilde{Z}^{-1}$ – ищется из условия нормировки, $f^* = \left( f^л, f^o \right)$, $f^л = \bar{f}^л - s^л$, $f^o = \bar{f}^o - s^o$, а $s^л$, $s^o$ – (оптимальные) векторы двойственных множителей для ограничений $t^л \ge \bar{t}^л$, $t^o \ge \bar{t}^o$ в задаче (1.1.11). Таким образом, все что осталось сделать, это решить задачу негладкой выпуклой оптимизации (1.1.11) прямо-двойственным методом. Заметим при этом, что на переменные $t^л$, $t^o$ – ограничения сверху возникают из вполне естественных соображений. Пусть $\bar{f}$ максимально возможный поток на ребре, поскольку нас интересует оценка сверху, то можно считать $t > \bar{t}$, стало быть, $f = \bar{f}$. Представим себе самую плохую ситуацию: все ребро стоит в пробке. Пусть длина ребра $L$, а средняя длина автомобиля $l$, число полос $r$. Тогда $t \le Lr / \left( l\bar{f} \right) + \bar{t}$ не может превышать нескольких часов. Задачу (1.1.11) можно ре-шать, например, прямо-двойственным универсальным композитным методом треугольника [53] (см. также раздел 2.2 главы 2) или методом из работ [41, 43] (см. также раздел 1.6 этой главы 1).

---

[16] Отметим, что все приводимые далее выкладки (а также выкладки следующего пункта) можно провести, отталкиваясь не от задачи (1.1.2), а от её упрощенного варианта, описанного в сноске 9.

Отметим, что если известна информация о потоках по ряду дуг $f_e^л = \tilde{f}_e^л$, $e \in \mathrm{A}^л$; $f_e^o = \tilde{f}_e^o$, $e \in \mathrm{A}^o$, то эту информацию можно "зашить" в модель (ЗС) подобно тому, как это делается в работе [275]: а именно брать минимум по множеству $t_e^л \geq \bar{t}_e^л$, $e \in E \setminus \mathrm{A}^л$, $t_e^л \geq 0$, $e \in \mathrm{A}^л$, $t_e^o \geq \bar{t}_e^o$, $e \in E \setminus \mathrm{A}^o$, $t_e^o \geq 0$, $e \in \mathrm{A}^o$, а слагаемые $+\beta\langle \bar{f}^л, t^л - \bar{t}^л\rangle + \beta\langle \bar{f}^o, t^o - \bar{t}^o\rangle$ стоит заменить на (параметры $\bar{t}_e^л$, $e \in \mathrm{A}^л$ и $t_e^o$, $e \in \mathrm{A}^o$ – неизвестны, но они и не нужны для расчетов, поскольку входят в виде аддитивных констант в функционал)

$$+\beta \sum_{e \in E \setminus \mathrm{A}^л} \bar{f}_e^л \cdot \left(t_e^л - \bar{t}_e^л\right) + \beta \sum_{e \in \mathrm{A}^л} \tilde{f}_e^л \cdot \left(t_e^л - \bar{t}_e^л\right) + \beta \sum_{e \in E \setminus \mathrm{A}^o} \bar{f}_e^o \cdot \left(t_e^o - \bar{t}_e^o\right) + \beta \sum_{e \in \mathrm{A}^o} \tilde{f}_e^o \cdot \left(t_e^o - \bar{t}_e^o\right).$$

### 1.1.10 Стохастический вариант трехстадийной модели стабильной динамики

Рассуждая аналогично тому, как мы делали выше, можно обобщить результаты подразделе 1.1.9 на случай, когда вместо модели стабильной динамики используется её стохастический вариант (с параметром $T > 0$) [35, 49, 262] (развитие подхода этого пункта будет описано в разделе 1.4 главы 1):

$$\min_{\substack{t^л \geq \bar{t}^л \\ t^o \geq \bar{t}^o \\ \lambda^L, \lambda^W}} \left\{ \ln\left( \sum_{i,j=1}^{n} \exp\left( \beta T \psi_{ij}\left( \frac{t^л}{T}, \frac{t^o}{T} \right) - \lambda_i^L - \lambda_j^W \right) \right) + \right.$$

$$\left. + \sum_{i=1}^{n} \lambda_i^L l_i + \sum_{j=1}^{n} \lambda_j^W \mathrm{w}_j + \beta\langle \bar{f}^л, t^л - \bar{t}^л \rangle + \beta\langle \bar{f}^o, t^o - \bar{t}^o \rangle \right\}, \qquad (1.1.12)$$

где

$$\psi_{ij}\left(t^л, t^o\right) = \ln\left( \sum_{p \in P_{(i,j)}^Л} \exp\left( -\sum_{e \in E} \delta_{ep} t_e^л \right) + \sum_{p \in P_{(i,j)}^O} \exp\left( -\sum_{e \in E} \delta_{ep} t_e^o \right) \right).$$

Причем здесь, также как и в предыдущем пункте,

$$d_{ij} = \breve{Z}^{-1} \exp\left( \beta T \psi_{ij}\left( t^л/T, t^o/T \right) - \lambda_i^L - \lambda_j^W \right),$$

где $\breve{Z}^{-1}$ – ищется из условия нормировки, $f^* = \left(f^л, f^o\right)$, $f^л = \bar{f}^л - s^л$, $f^o = \bar{f}^o - s^o$, а $s^л$, $s^o$ – (оптимальные) векторы двойственных множителей для ограничений $t^л \geq \bar{t}^л$, $t^o \geq \bar{t}^o$ в задаче (1.1.11). Таким образом, все что осталось сделать, это решить задачу гладкой выпуклой оптимизации (1.1.12) прямо-двойственным методом, например, [41, 43, 44]. Мы не будем здесь подробно описывать возможные способы решения. Обратим только вни-мание на то, что задача вычисления значений и градиента функции $\psi\left(t^{л},t^{o}\right)$ – вычислительно сложная, и, на первый взгляд, даже кажется бесперспективной из-за потенциально экспоненциально большого числа возможных маршрутов. Однако с помощью “аппарата характеристических функций на ориентированных графах” [33, 49, 262] и быстрого автоматического дифференцирования [65, 76] пересчитывать значения функции $\psi\left(t^{л},t^{o}\right)$ и ее градиента можно довольно эффективно – см. алгоритм из работ [33, 49, 262], который вырождается при $T \to 0+$ в алгоритм Беллмана–Форда [33, 49, 262]. Отметим, что при этом предельном переходе модель подраздела 1.1.10 перейдет в модель подраздела 1.1.9.

Отметим также, что все сказанное здесь переносится и на функции $\psi_{ij}\left(t^{л},t^{o}\right)$ более общего вида, соответствующие различным иерархическим способам выбора (типа Nested Logit) [79, 132, 325]. Например, “практически бесплатно” можно сделать такую замену:

$$T\psi_{ij}\left(\frac{t^{л}}{T},\frac{t^{o}}{T}\right)\to$$

$$\to \eta \ln\left(\exp\left(T\ln\left(\sum_{p\in P^{л}_{(i,j)}}\exp\left(-\sum_{e\in E}\delta_{ep}\frac{t^{л}_{e}}{T}\right)\right)\Big/\eta\right)+\exp\left(T\ln\left(\sum_{p\in P^{O}_{(i,j)}}\exp\left(-\sum_{e\in E}\delta_{ep}\frac{t^{o}_{e}}{T}\right)\right)\Big/\eta\right)\right).$$

### 1.1.11 Калибровка модели стабильной динамики

О практическом использовании модели стабильной динамики написано, например, в [172] и немного в [275]. В этом пункте мы сконцентрируем внимание на потенциальных плюсах описанной модели с точки зрения её калибровки. Далее мы опускаем нижней индекс $e$.

Способ оценки $\bar{t}$ довольно очевиден: $\bar{t}$ определяется длиной участка (ребра) и типом ребра (в сельской местности, в городе, шоссе). Вся эта информация обычно бывает доступной. С оценкой $\bar{f}$ немного посложнее. Пусть в конце дуги, состоящей из $r$ полос, стоит светофор, который пускает поток с этой дуги долю времени $\chi$. Пусть $q_{\max} \approx 1800\left[авт/час\right]$ – максимально возможное значение потока по одной полосе (это понятие не совсем корректное, но в первом приближении, им можно пользоваться [49], и оно довольно универсально). Тогда [49] $\bar{f} \approx \chi r q_{\max}$. Тут имеется важный нюанс. Во время зеленой фазы (в зависимости от того, что это за фаза, скажем, движение прямо и налево или движение только направо) “работают” не все полосы: $\bar{f} \approx \chi_1 r_1 q_{\max} + ... + \chi_l r_l q_{\max}$, где $l$ – число фаз, “пускающих” поток с рассматриваемой дуги, $\chi_k$ – доля времени отводимое фазе $k$, $r_k$ – эффективное число полос рассматриваемой дуги, задействованных на фазе $k$.

Раздобыть информацию по работе светофоров, как правило, не удается. Однако, существует определенные регламенты, согласно которым и устанавливаются фазы работы светофора. Поскольку все довольно типизировано, то часто бывает достаточно иметь информацию только о полосности дорог. Такая информация с 2012 года уже включается в коммерческие системы ГИС.

К сожалению, в реальных транспортных сетях в часы пик бывают пробки, которые приводят к тому, что пропускная способность ребра определяется пропускной способностью не данного ребра, а какого-то из впереди идущих (по ходу движения) ребер. Другими словами, пробка с ребра полностью “заполнила” это ребро и распространилась на ребра, входящие в это ребро. В таких ситуациях $\bar{f}$ нужно считать исходя из ограничений на пропускную способность на впереди идущих ребрах. При этом в выборе разбиения транспортного графа на ребра и в определении на этих ребрах значений $\bar{f}$, $\bar{t}$ стоит исходить из того, чтобы в типичной ситуации пробка если и переходила с ребра на ребро, то, желательно, чтобы это происходило без изменения пропускной способности входящего ребра. Этого не всегда можно добиться, поскольку часто приходится осуществлять разбиение без особого произвола, исходя из въездов/съездов, перекрестков. В таких случаях, формула $\bar{f} \approx \chi_1 r_1 q_{\max} + ... + \chi_l r_l q_{\max}$ является лишь оценкой сверху для “среднего” (типичного) значения, которое надо подставлять в модель (часть этого бремени придется переносить и на $\bar{t}$, увеличивая его).

Заметим, что также как и для обычных многостадийных моделей, для калибровки предложенной модели необходимо использовать один и тот же промежуток времени (скажем, обеденные часы) каждого типового дня (например, буднего) – в зависимости от целей. Кроме того, для Москвы в утренние часы пик значения потоков существенно нестационарные, и в течение час они, как правило, меняются сильно. Вместе с этим среднее время в пути оказывается больше часа. Таким образом, возникает вопрос: что понимается, например, под равновесным распределением потоков $f$? Ответ: “средние значения” нуждается в пояснении. Средние значения не только по дням, но и по исследуемому промежутку времени (в несколько часов) внутри каждого дня. Как следствие получаем, что выводы, сделанные по предложенной модели, нельзя, например, напрямую использовать для краткосрочного прогнозирования ситуации на дорогах или адаптивного управления светофорной сигнализацией. Таким образом, модель работает и выдает не реальные данные, хотя многие переменные модели и обозначают реальные физические параметры транспортного потока, а лишь некоторые средние (агрегированные) показатели, которые, тем не менее, многое могут сказать о ситуации на дорогах.

В данном разделе была предложена модель, описывающая равновесие в транспортной сети с точки зрения макро масштабов времени (месяцы). Для исследования поведения транспортного потока внутри одного дня требуются микро модели. Контекст, в котором такие модели используются, часто связан с необходимостью краткосрочного прогнозирования (на несколько часов вперед) и оптимального управления, например, светофорами. Как правило, в таких микро моделях численно исследуется начально-краевая задача для нелинейных УЧП (системы законов сохранения), в которой начальные условия имеются в наличии, а краевые условия (характеристики источников и стоков во времени, матрицы перемешивания в узлах графа транспортной сети) определяются исходя из исторической информации. Однако в последнее время в Москве все чаще можно слышать предложения о том, чтобы максимально информировать участников дорожного движения с целью лучшей маршрутизации. Такая полная информированность приводит к необходимости определять, скажем, матрицы перемешивание не из (и не только из) исторической информации, но и исходя из равновесных принципов, заложенных в модели Бэкмана и стабильной динамики. Причем нам представляется, что модель стабильной динамики подходит намного лучше для этих целей (мотивация имеется в работе [275], в которой на простом примере показывается, как возникает равновесная конфигурация модели стабильной динамики из внутредневной динамики водителей). Как отмечалось в работе [275], за исключением модельных примеров (Рисунок 1.1.2), очень сложной с вычислительной точки зрения представляется задача описания динамического режима функционирования модели, имеющей в своей основе модель стабильной динамики (или какую-то другую равновесную модель). Более того, пока, насколько нам известно, не было предложено ни одной более-менее обоснованной модели такого типа (иногда такие модели называют моделями динамического равновесия [287]). Здесь мы лишь укажем на одно, на наш взгляд, перспективное направление: объединить модель стабильной динамики с микромоделью CTM К. Даганзо [49] в современном ее варианте [63].

## 1.2 Эволюционные выводы энтропийной модели расчета матрицы корреспонденций

### 1.2.1. Введение

Настоящий раздел 1.2 развивает подраздел 1.1.3 раздела 1.1 этой главы 1. В этом разделе предлагается два способа объяснения популярной в урбанистике энтропийной модели расчета матрицы корреспонденций [26]. Обе предложенные модели в своей основе имеют марковский процесс, живущий в пространстве огромной размерности (говорят, что такой процесс порождает макросистему). Более точно, этот марковский процесс представляет собой ветвящийся процесс специального вида "модель стохастической химической кинетики" [86, 195]. Более того, в обоих случаях имеет место условие детального равновесия [12, 32, 86]. Из общих результатов [12, 32, 86] отсюда сразу можно заключить, что инвариантная (стационарная) мера такого процесса будет иметь вид мультиномиального распределения, сосредоточенного на аффинном многообразии, определяемым линейными законами сохранения введенной динамики. Исходя из явления концентрация меры [159] можно ожидать концентрацию этой меры около наиболее вероятного состояния с ростом размера макросистемы. Такое состояние, естественно, и принимать за равновесие изучаемой макросистемы, поскольку с большой вероятностью на больших временах мы найдем систему в малой окрестности такого состояния. Поиск равновесия сводится по теореме Санова [111] к задаче энтропийно-линейного программирования (1.2.3) [198]. Собственно, именно таким образом и планируется объяснить, почему для описания равновесной матрицы корреспонденций необходимо решить задачу ЭЛП.

В подразделе 1.2.2 на базе бинарных реакций обменного типа, популярных в различного рода физических и социально-экономических приложениях моделей стохастической химической кинетики [20, 28], будет приведен первый способ эволюционного вывода энтропийной модели расчета матрицы корреспонденций. Второй способ будет описан в подразделе 1.2.3. Он базируется на классических понятиях популяционной теории игр [302]. Например, таких как логит динамика и игра загрузки (именно к такой игре сводится поиск равновесного распределения потоков по путям). Второй способ также базируется на редукции задачи расчета матрицы корреспонденций к задаче поиска равновесного распределения потоков по путям. Нетривиальным и, по-видимому, новым здесь является эволюционно-экономическая интерпретация двойственных множителей, обобщающая известную стационарную конструкцию, восходящую к Л.В. Канторовичу.

Отметим, что результаты данного раздела являются обобщением результатов работ [47, 49]. Обобщение заключается в том, что в формулировках основных результатов (теоремы 1.2.1 и 1.2.2) фигурируют точные (неулучшаемые) оценки времени сходимости изу-чаемых макросистем на равновесие и плотности концентрации инвариантной меры в окрестности равновесий.

Мы намеренно опускаем вопросы численного решения, возникающих задач ЭЛП. Подробнее об этом можно посмотреть, например, в [37, 198] и других частях диссертации.

### 1.2.2 Вывод на основе обменов

Приведем, базируясь на работе [49], эволюционное обоснование одного из самых популярных способов расчета матрицы корреспонденций, имеющего более чем сорокалетнюю историю, – энтропийной модели [26].

Пусть в некотором городе имеется $n$ районов, $L_i > 0$ – число жителей $i$-го района, $W_j > 0$ – число работающих в $j$-м районе. При этом $N = \sum_{i=1}^{n} L_i = \sum_{j=1}^{n} W_j$, – общее число жителей города, $n^2 \ll N$. Далее под $L_i \geq 0$ будет пониматься число жителей района, выезжающих в типичный день за рассматриваемый промежуток времени из $i$-го района, а под $W_j \geq 0$ – число жителей города, приезжающих на работу в $j$-й район в типичный день за рассматриваемый промежуток времени. Обычно, так введенные, $L_i$, $W_j$ рассчитываются через число жителей $i$-го района и число работающих в $j$-м районе с помощью более менее универсальных (в частности, не зависящих от $i$, $j$) коэффициентов пропорциональности. Эти величины являются входными параметрами модели, т.е. они не моделируются (во всяком случае, в рамках выбранного подхода). Для долгосрочных расчетов с разрабатываемой моделью требуется иметь прогноз изменения значений этих величин.

Обозначим через $d_{ij}(t) \geq 0$ – число жителей, живущих в $i$-м районе и работающих в $j$-м в момент времени $t$. Со временем жители могут только меняться квартирами, поэтому во все моменты времени $t \geq 0$

$$d_{ij}(t) \geq 0, \ \sum_{j=1}^{n} d_{ij}(t) \equiv L_i, \ \sum_{i=1}^{n} d_{ij}(t) \equiv W_j, \ i, j = 1, ..., n.$$

Определим

$$\mathrm{A} = \left\{ d_{ij} \geq 0: \ \sum_{j=1}^{n} d_{ij} = L_i, \sum_{i=1}^{n} d_{ij} = W_j, i, j = 1, ..., n \right\}.$$

Опишем основной стимул к обмену: работать далеко от дома плохо из-за транспортных издержек. Будем считать, что эффективной функцией затрат [49] будет $R(T) = \beta T/2$, где $T > 0$ – время в пути от дома до работы (в общем случае под $T$ стоит понимать затра-ты, в которые может входить не только время), а $\beta > 0$ – настраиваемый параметр модели (который также можно проинтерпретировать и даже оценить, что и будет сделано ниже).

Теперь опишем саму динамику. Пусть в момент времени $t \ge 0$ $r$-й житель живет в $k$-м районе и работает в $m$-м, а $s$-й житель живет в $p$-м районе и работает в $q$-м. Тогда $\lambda_{k,m;\,p,q}(t)\Delta t + o(\Delta t)$ – есть вероятность того, что жители с номерами $r$ и $s$ ($1 \le r < s \le N$) "поменяются" квартирами в промежутке времени $(t, t+\Delta t)$. Вероятность обмена местами жительства зависит только от мест проживания и работы обменивающихся:

$$\lambda_{k,m;\,p,q}(t) \equiv \lambda_{k,m;\,p,q} = \lambda N^{-1} \exp\Big( \underbrace{R(T_{km}) + R(T_{pq})}_{\text{суммарные затраты до обмена}} - \underbrace{\big(R(T_{pm}) + R(T_{kq})\big)}_{\text{суммарные затраты после обмена}} \Big) > 0,$$

где коэффициент $0 < \lambda = \mathrm{O}(1)$ характеризует интенсивность обменов. Совершенно аналогичным образом можно было рассматривать случай "обмена местами работы". Здесь стоит оговориться, что "обмены" не стоит понимать буквально – это лишь одна из возможных интерпретаций. Фактически используется, так называемое, "приближение среднего поля" [32, 195], т.е. некое равноправие агентов (жителей) внутри фиксированной корреспонденции и их независимость.

Согласно эргодической теореме для марковских цепей (в независимости от начальной конфигурации $\{d_{ij}(0)\}_{i=1,\, j=1}^{n,\, n}$) [16, 20, 28, 49, 86, 246, 302] предельное распределение совпадает со стационарным (инвариантным), которое можно посчитать (получается проекция прямого произведение распределений Пуассона на A):

$$\lim_{t\to\infty} P\big(d_{ij}(t) = d_{ij},\, i,j = 1,\dots,n\big) = Z^{-1} \prod_{i,j=1}^{n} \exp\big(-2R(T_{ij})d_{ij}\big)\cdot\big(d_{ij}!\big)^{-1} \stackrel{def}{=} p\Big(\{d_{ij}\}_{i=1,\, j=1}^{n,\, n}\Big), \qquad (1.2.1)$$

где $\{d_{ij}\}_{i=1,\, j=1}^{n,\, n} \in \mathrm{A}$, а "статсумма" $Z$ находится из условия нормировки получившейся "пуассоновской" вероятностной меры. Отметим, что стационарное распределение $p\Big(\{d_{ij}\}_{i=1,\, j=1}^{n,\, n}\Big)$ удовлетворяет условию детального равновесия [32, 302]:

$$(d_{km}+1)(d_{pq}+1)\, p\Big(\big\{d_{11},\dots,d_{km}+1,\dots,d_{pq}+1,\dots,d_{pm}-1,\dots,d_{kq}-1,\dots,d_{nn}\big\}\Big)\lambda_{k,m;\,p,q} =$$

$$= d_{pm} d_{kq}\, p\Big(\{d_{ij}\}_{i=1,\, j=1}^{n,\, n}\Big)\lambda_{p,m;\,k,q}.$$

При $N \gg 1$ распределение $p\Big(\{d_{ij}\}_{i=1,\, j=1}^{n,\, n}\Big)$ экспоненциально сконцентрировано на множестве A в $\mathrm{O}\big(\sqrt{N}\big)$ окрестности наиболее вероятного значения $d^* = \{d_{ij}^*\}_{i=1,\, j=1}^{n,\, n}$, которое определяется, как решение задачи энтропийно-линейного программирования [47, 49, 86]:

$$\ln p\left(\left\{d_{ij}\right\}_{i=1,\, j=1}^{n,\, n}\right) \sim -\sum_{i,j=1}^{n} d_{ij} \ln\left(d_{ij}\right) - \beta \sum_{i,j=1}^{n} d_{ij} T_{ij} \to \max_{\left\{d_{ij}\right\}_{i=1,\, j=1}^{n,\, n} \in (\mathrm{A})} . \qquad (1.2.2)$$

Это следует из теоремы Санова о больших уклонениях для мультиномиального распределения [111] (на распределение (1.2.1) можно смотреть также как на проекцию мультиномиального распределения на А)

$$\frac{N!}{d_{11}! \cdot \ldots \cdot d_{ij}! \cdot \ldots \cdot d_{nn}!} p_{11}^{d_{11}} \cdot \ldots \cdot p_{ij}^{d_{ij}} \cdot \ldots \cdot p_{nn}^{d_{nn}} = \exp\left(-N \sum_{i,j=1}^{n} \nu_{ij} \ln\left(\nu_{ij} / p_{ij}\right) + \bar{R}\right),$$

где $\nu_{ij} = d_{ij}/N$, $\left|\bar{R}\right| \le \frac{n^2}{2}\left(\ln N + 1\right)$, и формулы Тейлора (с остаточным членом второго порядка в форме Лагранжа), примененной к энтропийному функционалу в точке $d^*$, заданному на А [26].

Сформулируем более точно полученный результат (схожий результат возникнет в другом контексте в разделе 4.1 главы 4, см. также приложение в конце диссертации).

**Теорема 1.2.1.** *Для любого* $d(0) \in \mathrm{A}$ *существует такая константа* $c_n\left(d(0)\right) > 0$, *что для всех* $\sigma \in (0, 0.5)$, $t \ge c_n\left(d(0)\right) N \ln N$ *имеет место неравенство*

$$P\left(\frac{\left\|d(t) - d^*\right\|_2}{N} \ge \frac{2\sqrt{2} + 4\sqrt{\ln\left(\sigma^{-1}\right)}}{\sqrt{N}}\right) \le \sigma,$$

*где* $d(t) = \left\{d_{ij}(t)\right\}_{i=1,\, j=1}^{n,\, n} \in \mathrm{A}$.

**Схема доказательства.** Установим сначала оценку для плотности концентрации стационарной меры

$$\lim_{t\to\infty} P\left(d_{ij}(t) = d_{ij}, i, j = 1, \ldots, n\right) = \frac{N!}{d_{11}! \cdot \ldots \cdot d_{ij}! \cdot \ldots \cdot d_{nn}!} p_{11}^{d_{11}} \cdot \ldots \cdot p_{ij}^{d_{ij}} \cdot \ldots \cdot p_{nn}^{d_{nn}} .$$

Из неравенства Хефдинга в гильбертовом пространстве [159] следует ( $\varepsilon \ge \sqrt{2/N}$ )

$$\lim_{t\to\infty} P\left(\left\|d(t) - d^*\right\|_2 \ge \varepsilon N\right) \le \exp\left(-\frac{1}{4N}\left(\varepsilon N - \sqrt{2N}\right)^2\right).$$

Беря в этом неравенстве

$$\varepsilon = \frac{\sqrt{2} + 2\sqrt{\ln\left(\sigma^{-1}\right)}}{\sqrt{N}},$$

получим

$$\lim_{t\to\infty} P\left(\frac{\left\|d(t) - d^*\right\|_2}{N} \ge \varepsilon\right) \le \sigma.$$

Однако если не переходить к пределу по времени, а лишь обеспечить при $t \ge T_n\left(\varepsilon, N; d(0)\right)$ выполнение условия

$$E\left\|\frac{d(t)-d^*}{N}\right\|_2 \le \frac{\varepsilon}{2},$$

то при $t \ge T_n\left(\varepsilon, N; d(0)\right)$

$$P\left(\left\|d(t)-d^*\right\|_2 \ge \varepsilon N\right) \le \exp\left(-\frac{1}{4N}\left(\varepsilon N-\left(\sqrt{2N}+\frac{\varepsilon N}{2}\right)\right)^2\right) = \exp\left(-\frac{1}{4N}\left(\frac{\varepsilon N}{2}-\sqrt{2N}\right)^2\right).$$

Беря в этом неравенстве

$$\varepsilon = \frac{2\sqrt{2}+4\sqrt{\ln\left(\sigma^{-1}\right)}}{\sqrt{N}},$$

получим при $t \ge T_n\left(\varepsilon, N; d(0)\right)$ аналогичное неравенство

$$P\left(\frac{\left\|d(t)-d^*\right\|_2}{N} \ge \varepsilon\right) \le \sigma.$$

Отметим, что с точностью до мультипликативной константы эта оценка не может быть улучшена. Это следует из неравенства Чебышёва

$$P\left(X \ge EX-\varepsilon\right) \ge 1-\frac{\operatorname{Var}(X)}{\varepsilon^2},$$

при

$$X = \frac{\left\|d(t)-d^*\right\|_2}{N}, \; d^* = N\cdot\left(n^{-2},\dots,n^{-2}\right)^T, \; N \gg n^2.$$

Осталось оценить зависимость $T_n\left(\varepsilon, N; d(0)\right)$. Для этого воспользуемся неравенством Чигера [246] (см. также раздел 4.1 главы 4). Поставим в соответствие нашему марковскому процессу его дискретный аналог (в дискретном времени) – случайное блуждание (со скачками, соответствующими парным реакциям, введенным выше) на целочисленных точках части гиперплоскости, задаваемой множеством А. Граф, на котором происходит блуждание, будем обозначать $G = \langle V_G, E_G\rangle$. Пусть $\pi(\cdot)$ стационарная мера этого блуждания (в нашем случае – мультиномиальная мера, экспоненциально сконцентрированая на множестве А в $\mathrm{O}\left(\sqrt{N}\right)$ окрестности наиболее вероятного значения $d^*$), а $P = \left\|p_{ij}\right\|_{i,j=1}^{|V_G|,|V_G|}$ – матрица переходных вероятностей. Тогда, вводя константу Чигера ( $\bar{S} = V_G \setminus S$ )

$$h(G)=\min_{S\subseteq V_G:\pi(S)\le 1/2} P\left(S\to\bar{S}\,\middle|\,S\right)=\min_{S\subseteq V_G:\pi(S)\le 1/2}\frac{\sum_{(i,j)\in E_G:i\in S,\,j\in\bar{S}}\pi(i)\,p_{ij}}{\sum_{i\in S}\pi(i)}$$

и время выхода блуждания (стартовавшего из $i\in V_G$ ) на стационарное распределение

$$T(i,\varepsilon)=\mathrm{O}\left(h(G)^{-2}\left(\ln\left(\pi(i)^{-1}\right)+\ln\left(\varepsilon^{-1}\right)\right)\right),$$

получим для любых $i=1,\ldots,|V_G|,\ t\ge T(i,\varepsilon)$

$$\left\|P^t(i,\cdot)-\pi(\cdot)\right\|_2\le\left\|P^t(i,\cdot)-\pi(\cdot)\right\|_1=\sum_{j=1}^{n}\left|P^t(i,j)-\pi(j)\right|\le\frac{\varepsilon}{2},$$

где $P^t(i,j)$ – условная вероятность того, что стартуя из состояния $i\in V_G$ через $t$ шагов, марковский процесс окажется в состоянии $j\in V_G$.

В нашем случае можно явно геометрически описать множество $S$ в виде целочисленных точек множества А, попавших внутрь эллипсоида (сферы) размера $\mathrm{O}\left(\sqrt{N}\right)$ с $\pi(S)\simeq 1/2$ и с центром в точке $d^*$. На этом множестве достигается решение изопериметрической задачи в определении константы Чигера. Значение константы Чигера при этом будет пропорционально отношению площади поверхности этого эллипсоида к его объему, т.е. $h(G)\sim N^{-1/2}$. Используя это наблюдение, можно получить, что

$$T_n\left(\varepsilon,N;d(0)\right)\sim N\cdot\left(\ln\left(\pi(i)^{-1}\right)+\ln N\right).$$

Отсюда видно, что время выхода зависит от точки старта. Если ограничиться точками старта $i\in V_G$, для которых равномерно по $N\to\infty$ и $i,j=1,\ldots,n$ имеет место неравенство $d_{ij}(0)/N\ge\varsigma_n>0$, то $\ln\left(\pi(i)^{-1}\right)\sim\ln N$.

Объединяя приведенные результаты, получаем утверждение теоремы 1.2.1. ●

Эта теорема уточняет результат работы [49], явно указывая скорость сходимости и плотность концентрации. Приведенная скорость сходимости характерна для более широкого класса моделей стохастической химической кинетики, приводящих к равновесию вида неподвижной точки [32]. Плотность концентрации оценивалась на базе конструкции работы [46].

Естественно в виду теоремы 1.2.1 принимать решение задачи (1.2.2) $\left\{d_{ij}^*\right\}_{i=1,\,j=1}^{n,\,n}$ за равновесную конфигурацию. Обратим внимание, что предложенный выше вывод известной энтропийной модели расчета матрицы корреспонденций отличается от классического [26]. В монографии А.Дж. Вильсона [26] $\beta$ интерпретируется как множитель Лагранжа к ограничению на среднее “время в пути”: $\sum_{i,j=1}^{n} d_{ij}T_{ij} = C$.

**Замечание 1.2.1.** При этом остальные ограничения имеют такой же вид, а функционал имеет вид $F(d) = -\sum_{i,j=1}^{n} d_{ij}\ln(d_{ij})$. Тогда, согласно экономической интерпретации двойственных множителей Л.В. Канторовича: $\beta(C) = \partial F(d(C))/\partial C$. Из такой интерпретации иногда делают вывод о том, что $\beta$ можно понимать, как цену единицы времени в пути. Чем больше $C$, тем меньше $\beta$.

Приведенный нами вывод позволяет контролировать знак параметра $\beta > 0$ и лучше понимать его физический смысл.

**Замечание 1.2.2.** Отметим, что также как и в [26] из принципа ле Шателье–Самуэльсона следует, что с ростом $\beta$ среднее время в пути $\sum_{i,j=1}^{n} d_{ij}(\beta)T_{ij}$ будет убывать. В связи с этим обстоятельством, а также исходя из соображений размерности, вполне естественно понимать под $\beta$ величину, обратную к характерному (среднему) времени в пути [26] – физическая интерпретация. Собственно, такая интерпретация параметра $\beta$, как правило, и используется в многостадийных моделях (см., например, [287]).

В дальнейшем, нам будет удобно привести задачу (1.2.2) к следующему виду (при помощи метода множителей Лагранжа [85], теоремы фон Неймана о минимаксе [313] и перенормировки $d := d/N$):

$$\max_{\lambda^L,\lambda^W} \min_{\sum_{i,j=1}^{n} d_{ij}=1,\, d_{ij}\ge 0} \left[ \sum_{i,j=1}^{n} d_{ij}\ln d_{ij} + \beta \sum_{i=1,\, j=1}^{n,n} d_{ij}T_{ij} + \sum_{i=1}^{n} \lambda_i^L \left( \sum_{j=1}^{n} d_{ij} - l_i \right) + \sum_{j=1}^{n} \lambda_j^W \left( \mathrm{w}_j - \sum_{i=1}^{n} d_{ij} \right) \right], \quad (1.2.3)$$

где $l = L/N$, $\mathrm{w} = W/N$.

**Замечание 1.2.3.** Используя принцип Ферма [85] не сложно проверить, что решение задачи (1.2.3) можно представить в виде

$$d_{ij} = \exp(-\lambda_i^L)\exp(\lambda_j^w)\exp(-\beta T_{ij}),$$

где

$$\exp(\lambda_i^L) = \frac{1}{l_i}\sum_{j=1}^{n}\exp(\lambda_j^w)\exp(-\beta T_{ij}), \quad \exp(-\lambda_j^w) = \frac{1}{w_j}\sum_{i=1}^{n}\exp(-\lambda_i^L)\exp(-\beta T_{ij}). \quad (1.2.4)$$

Отсюда можно усмотреть интерпретацию двойственных множителей как соответствующих “потенциалов притяжения/отталкивания районов” [26]. К этому мы еще вернемся в подразделе 1.2.3.

**Замечание 1.2.4.** Департамент транспорта г. Москвы несколько лет назад поставил следующую задачу. Сколько человек надо обзвонить и опросить на предмет того какой корреспонденции они принадлежат (где живут и где работают), чтобы восстановить матрицу корреспонденций с достаточной точностью и доверительным уровнем? Формализуем задачу. Предположим, что истинная (пронормированная) матрица корреспонденций $\left\{d_{ij}^*\right\}_{i=1,\, j=1}^{n,\, n}$ – ее и надо восстановить. В результате опросов населения получилась матрица (вектор) $r=\left\{r_{ij}\right\}_{i=1,\, j=1}^{n,\, n}$, где $r_{ij}$ – количество опрошенных респондентов, принадлежащих корреспонденции $\left(i,j\right)$, $\sum_{i,j=1}^{n} r_{ij}=N$. Задачу формализуем следующим образом: найти наименьшее $N$, чтобы

$$P_{d^*}\left(\left\|\frac{r_{ij}}{N}-d^*\right\|_2\geq\varepsilon\right)\leq\sigma.$$

Нижний индекс $d^*$ у вероятности означает, что случайный вектор $r$ имеет мультиномиальное распределение с параметром $d^*$, т.е. считаем, что опрос проводился в идеальных условиях. Из теоремы 1.2.1 не сложно усмотреть, что достаточно опросить $N=\left(4+8\ln\left(\sigma^{-1}\right)\right)\varepsilon^{-2}$ респондентов. Скажем, опрос ~50 000 респондентов, который и был произведен, действительно позволяет неплохо восстановить матрицу корреспонденций $d^*\approx\bar{d}\overset{def}{=}r_{ij}/N$. Однако мы привели здесь это замечание для других целей. В ряде работ (см., например [304, 328]) исходя из данных таких опросов также восстанавливают матрицу корреспонденций, но при другой параметрической гипотезе (меньшее число параметров): $d_{ij}=\exp\left(-\lambda_i^L\right)\exp\left(\lambda_j^w\right)\exp\left(-\beta T_{ij}\right)$. То есть дополнительно предполагают, что $n^2$ неизвестных параметров в действительности однозначно определяются $2n$ параметрами (иногда к ним добавляют еще один параметр $\beta$), которые у нас ранее (замечание 1.2.3) интерпретировались как множители Лагранжа. Встает вопрос: как оптимально оценить эти параметры? Ответ дает теорема Фишера (в современном не асимптотическом варианте изложение этой теоремы можно найти в [319]) об оптимальности оценок максимального правдоподобия (ОМП). Собственно, для выборки из мультиномиального распределения оценкой максимального правдоподобия как раз и будет выборочное среднее $\bar{d}$. Для описанной параметрической модели (с $2n$ параметрами) поиск такой оценки сводится разрешению системы (1.2.4) (замечание 1.2.3).

**Замечание 1.2.5.** Подобно тому, как мы рассматривали в этом пункте трудовые корреспонденции (в утренние и вечерние часы более 70% корреспонденций по Москве и области именно такие), можно рассматривать перемещения, например, из дома к местам учебы, отдыха, в магазины и т.п. (по хорошему, еще надо было учитывать перемещения типа работа–магазин–детский_сад–дом) – рассмотрение всего этого вкупе приведет также к задаче (1.2.2). Только будет больше типов корреспонденций $d$: помимо пары районов, еще нужно будет учитывать тип корреспонденции [26]. Все это следует из того, что инвариантной мерой динамики с несколькими типами корреспонденций по-прежнему будет прямое произведение пуассоновских мер. Другое дело, когда мы рассматриваем разного типа пользователей транспортной сети, например: имеющих личный автомобиль и не имеющих личный автомобиль. Первые могут им воспользоваться равно, как и общественным транспортом, а вторые нет. И на рассматриваемых масштабах времени пользователи могут менять свой тип. То есть время в пути может для разных типов пользователей быть различным [26]. Считая, подобно тому как мы делали раньше, что желание пользователей корреспонденции $(i, j)$ сети сменить свой тип (вероятность в единицу времени) есть

$$\tilde{\lambda}\exp\Big( \underbrace{\tilde{R}\left(T_{ij}^{old}\right)}_{\text{суммарные затраты до смены типа}} - \underbrace{\tilde{R}\left(T_{ij}^{new}\right)}_{\text{суммарные затраты после смены типа}} \Big), \text{ где } \tilde{R}(T) = \beta T,$$

и учитывая в “обменах” тип пользователя (будет больше типов корреспонденций $d$, но “меняются местами работы” только пользователи одного типа), можно показать, что все это вкупе приведет также к задаче типа (1.2.2).

### 1.2.3 Вывод на основе модели равновесного распределения потоков

Предварительно напомним, следуя [47, 302] эволюционный вывод модели равновесного распределения потоков [47, 49, 289, 310].

Задан ориентированный граф $\Gamma = (V, E)$ – транспортная сеть города ($V$ – узлы сети (вершины), $E \subset V \times V$ – дуги сети (ребра графа)). В графе имеются две выделенные вершины. Одна из вершин графа является источником, другая стоком. Из источника в сток ведет много путей, которые мы будем обозначать $p \in \mathrm{P}$, $|\mathrm{P}| = m$. Обозначим через

$x_p$ – величина потока по пути $p$, $x_p \in S_m(N) = \left\{ x = \{x\}_{p\in\mathrm{P}} \geq 0 : \sum_{p\in\mathrm{P}} x_p = N \right\}$;

$y_e$ – величина потока по ребру $e \in E$: $y_e = \sum_{p\in\mathrm{P}} \delta_{ep} x_p$ ($y = \Theta x$, $\Theta = \left\{\delta_{ep}\right\}_{e\in E, p\in\mathrm{P}}$), где

$$\delta_{ep} = \begin{cases} 1, & e \in p \\ 0, & e \notin p \end{cases};$$

$\tau_e(y_e)$ – удельные затраты на проезд по ребру $e$ (гладкие неубывающие функции);

$G_p(x)=\sum_{e\in E}\tau_e(y_e)\delta_{ep}$ – удельные затраты на проезд по пути $p$.

Рассмотрим, следуя [302], популяционную игру в которой имеется набор $N$ однотипных игроков (агентов). Множеством чистых стратегий каждого такого агента является $\mathrm{P}$, а выигрыш (потери со знаком минус) от использования стратегии $p\in\mathrm{P}$ определяются формулой $-G_p(x)$.

Опишем динамику поведения игроков. Пусть в момент времени $t\geq 0$ агент использует стратегию $p\in\mathrm{P}$, $\lambda_{p,q}(t)\Delta t+o(\Delta t)$ – вероятность того, что он поменяет свою стратегию на $q\in\mathrm{P}$ в промежутке времени $(t,t+\Delta t)$. Будем считать, что

$$\lambda_{p,q}(t)\equiv\lambda_{p,q}=\lambda P_q\left(\left\{G_p\left(x(t)\right)\right\}_{p\in\mathrm{P}}\right),$$

где (как и в подразделе 1.2.2) коэффициент $0<\lambda=\mathrm{O}(1)$ характеризует интенсивность “перескоков” агентов, а

$$P_q\left(\left\{G_p\left(x(t)\right)\right\}_{p\in\mathrm{P}}\right)=P\left(q=\arg\max_{p\in\mathrm{P}}\left\{-G_p\left(x(t)\right)+\xi_p\right\}\right).$$

Если $\xi_p\equiv 0$, то получаем динамику наилучших ответов [302], если $\xi_p$ – независимые одинаково распределенные случайные величины с распределением Гумбеля [137]: $P(\xi_p<\xi)=\exp\left\{-e^{-\xi/\omega-E}\right\}$, где $\omega\in(0,\omega_0]$ ($\omega_0=\mathrm{O}(1)$), $E\approx 0.5772$ – константа Эйлера, а $\mathrm{Var}\left[\xi_p\right]=\omega^2\pi^2/6$, то получаем логит динамику [302]

$$P_q\left(\left\{G_p\left(x(t)\right)\right\}_{p\in\mathrm{P}}\right)=\frac{\exp\left(-G_q\left(x(t)\right)/\omega\right)}{\sum_{p\in\mathrm{P}}\exp\left(-G_p\left(x(t)\right)/\omega\right)},$$

вырождающуюся в динамику наилучших ответов при $\omega\to 0+$. Далее мы будем считать, что задана логит динамика. Объясняется такая динамика совершенно естественно. Каждый агент имеет какую-то картину текущего состояния загрузки системы, но либо он не имеет возможности наблюдать ее точно, либо он старается как-то спрогнозировать возможные изменения загрузок (а соответственно и затраты на путях) в будущем (либо и то и другое).

Приведем соответствующий аналог теоремы 1.2.1.

**Теорема 1.2.2.** *Для любого* $x(0)\in S_m(N)$ *существует такая константа* $c_m(x(0))>0$, *что для всех* $\sigma\in(0,0.5)$, $t\geq c_m(x(0))N\ln N$ *имеет место неравенство*

$$P\left(\left\|\frac{x(t)}{N}-x^*\right\|_2 \ge \frac{2\sqrt{2}+4\sqrt{\ln\left(\sigma^{-1}\right)}}{\sqrt{N}}\right)\le\sigma,\ x(t)\in S_m(N),$$

*где (стохастическое равновесие Нэша–Вардропа)*

$$x^* = \arg\min_{x\in S_m(1)} \Psi\left(y(x)\right)+\omega\sum_{p\in\mathrm{P}} x_p\ln\left(x_p\right), \tag{1.2.5}$$

$$\Psi\left(y(x)\right)=\sum_{e\in E}\int_0^{y_e(x)}\tilde{\tau}_e(z)dz,\ \tilde{\tau}_e(z)=\tau_e(zN).$$

**Схема доказательства.** Стационарная мера описанного марковского процесса имеет (с точностью до нормирующего множителя) вид (теорема 11.5.12 [302])

$$\frac{N!}{x_1!\cdot\ldots\cdot x_m!}\exp\left(-\frac{\Psi\left(y(x)\right)}{\omega}\right),\ x\in S_m(N).$$

Для оценок плотности концентрации достаточно заметить, что $\Psi\left(y(x)\right)$ – выпуклая функция, как композиция линейной и выпуклой функции. Поэтому оценки плотности концентрации здесь не могут быть хуже оценок в теореме 1.2.1. Рассуждения для оценки времени выхода проводятся аналогично теореме 1.2.1. ●

**Следствие (Proposition 1 [177]).** *Если $\omega\to 0+$, то решение задачи (1.2.5) сходится к*

$$x^* = \arg\min_{x\in S_m(1):\ \Theta x=y^*}\sum_{p\in\mathrm{P}} x_p\ln\left(x_p\right), \tag{1.2.6}$$

*где* $y^* = \arg\min_{y=\Theta x,\ x\in S_m(1)}\Psi(y)$.

Это следствие решает проблему обоснования гипотезы Бар-Гира [143]. Напомним вкратце в чем состоит эта гипотеза. Известно, см., например, [289], что поиск равновесного распределения потоков по путям (равновесия Нэша–Вардропа) в модели с одним источником и стоком сводится к задаче оптимизации (1.2.5) с $\omega=0$. Хотя функционал этой задачи выпуклый, но он не всегда строго выпуклый, даже в случае, когда все функции $\tau_e\left(y_e\right)$ – строго возрастающие (тогда можно лишь говорить о единственности равновесного распределения потоков по ребрам $y^*$). Гипотеза Бар-Гира говорит, что “в жизни” с большой вероятностью реализуется то равновесие из множества равновесий, которое находится из решения задачи (1.2.6).

**Замечание 1.2.6.** Распределение потоков по путям $x$ называется равновесием (Нэша–Вардропа) в рассматриваемой популяционной игре $\left\langle \left\{x_p\right\}_{p\in\mathrm{P}}, \left\{G_p(x)\right\}_{p\in\mathrm{P}}\right\rangle$, если из $x_p>0$, $p\in\mathrm{P}$ следует $G_p(x)=\min_{q\in\mathrm{P}} G_q(x)$.

Или, что то же самое:

$$\text{для любых } p \in \mathrm{P} \text{ выполняется } x_p \cdot \left( G_p(x) - \min_{q \in \mathrm{P}} G_q(x) \right) = 0.$$

**Замечание 1.2.7.** Если при $\omega = 0$ рассмотреть предельный случай

$$\tau_e(y_e) := \tau_e^\mu(y_e) \xrightarrow[\mu \to 0+]{} T_e,$$

то поиск равновесия Нэша–Вардропа просто сводится к поиску социального оптимума, что приводит в данном случае к решению транспортной задачи линейного программирования [132]. Если делать предельный переход с учетом ограничений на пропускные способности ребер

$$\tau_e(y_e) := \tau_e^\mu(f_e) \xrightarrow[\mu \to 0+]{} \begin{cases} T_e, & 0 \le y_e < \bar{y}_e \\ [T_e, \infty), & y_e = \bar{y}_e \end{cases},$$

то получится более сложная задача, которая описывает равновесие в (стохастической, если $\omega > 0$) модели стабильной динамики [47, 275].

Из рассмотренных в замечании 1.2.7 случаев выпала ситуация $\omega > 0$ $\tau_e(y_e) := \tau_e^\mu(y_e) \xrightarrow[\mu \to 0+]{} T_e$. Ее мы сейчас отдельно и рассмотрим на примере другого эволюционного способа обоснования энтропийной модели расчета матрицы корреспонденций. По ходу обсуждения этого способа с коллегами (прежде всего, Ю.Е. Нестеровым и Ю.В. Дорном) у предложенного подхода появилось название: “облачная модель”.

Предположим, что все вершины, отвечающие источникам корреспонденций, соединены ребрами с одной вспомогательной вершиной (облако № 1). Аналогично, все вершины, отвечающие стокам корреспонденций, соединены ребрами с другой вспомогательной вершиной (облако № 2). Припишем всем новым ребрам постоянные веса. И проинтерпретируем веса ребер, отвечающих источникам $\lambda_i^L$, например, как средние затраты на проживание (в единицу времени, скажем, в день) в этом источнике (районе), а веса ребер, отвечающих стокам $\lambda_j^W$, как уровень средней заработной платы (в единицу времени) в этом стоке (районе), если изучаем трудовые корреспонденции. Будем следить за системой в медленном времени, то есть будем считать, что равновесное распределение потоков по путям стационарно. Поскольку речь идет о равновесном распределении потоков, то нет необходимости говорить о затратах на путях или ребрах детализированного транспортного графа, достаточно говорить только затратах (в единицу времени), отвечающих той или иной корреспонденции. Таким образом, у нас есть взвешенный транспортный граф с одним источником (облако 1) и одним стоком (облако 2). Все вершины этого графа, кроме двух вспомогательных (облаков), соответствуют районам в модели расчета матрицы кор-респонденций. Все ребра этого графа имеют стационарные (не меняющиеся и не зависящие от текущих корреспонденций) веса $\left\{T_{ij};\lambda_i^L;-\lambda_j^W\right\}$. Если рассмотреть естественную в данном контексте логит динамику ($x \equiv d$), описанную выше, с $\omega = 1/\beta$ (здесь полезно напомнить, что согласно замечанию 1.2.2 $\beta$ обратно пропорционально средним затратам, а $\omega$ имеет как раз физическую размерность затрат), то поиск равновесия рассматриваемой макросистемы согласно теореме 1.2.1 приводит (в прошкалированных переменных) к задаче, сильно похожей на задачу (1.2.3) из подраздела 1.2.2

$$\min_{\sum_{i,j=1}^{n} d_{ij}=1,\, d_{ij}\ge 0}\left[\sum_{i,j=1}^{n} d_{ij}\ln d_{ij} + \beta\sum_{i=1,\, j=1}^{n,n} d_{ij}T_{ij} + \beta\sum_{i=1}^{n}\left(\lambda_i^L\sum_{j=1}^{n} d_{ij}\right) - \beta\sum_{i=1}^{n}\left(\lambda_j^W\sum_{i=1}^{n} d_{ij}\right)\right].$$

Разница состоит в том, что здесь мы не оптимизируем по $2n$ двойственным множителям $\lambda^L$, $\lambda^W$ (множителям Лагранжа). Более того, мы их и не интерпретируем здесь как двойственные множители, поскольку мы их ввели на этапе взвешивания ребер графа. Тем не менее, значения этих переменных, как правило, не откуда брать. Тем более что приведенная выше (наивная) интерпретация вряд ли может всерьез рассматриваться, как способ определения этих параметров исходя из данных статистики. Более правильно понимать $\lambda_i^L$, $\lambda_j^W$ как потенциалы притяжения/отталкивания районов (см. также замечание 3), включающиеся в себя плату за жилье и зарплату, но включающие также и многое другое, что сложно описать количественно. И здесь как раз помогает информация об источниках и стоках, содержащаяся в $2n$ уравнениях задающих множество А. Таким образом, мы приходим ровно к той же самой задаче (1.2.3) с той лишь разницей, что мы получили дополнительную интерпретацию двойственных множителей в задаче (1.2.3). При этом двойственные множители в задаче (1.2.3) равны (с точностью до мультипликативного фактора $\beta$) введенным здесь потенциалам притяжения районов.

Нам представляется очень плодотворной и перспективной идея перенесения имеющейся информации об исследуемой системе из обременительных законов сохранения динамики, описывающей эволюцию этой системы, в саму динамику путем введения дополнительных естественно интерпретируемых параметров. При таком подходе становится возможным, например, учитывать в моделях и рост транспортной сети. Другими словами, при таком подходе, например, можно естественным образом рассматривать также и ситуацию, когда число пользователей транспортной сетью меняется со временем (медленно).

## 1.3 Об эффективной вычислимости конкурентных равновесий в транспортно-экономических моделях

### 1.3.1 Введение

В данном разделе мы развиваем конструкцию подраздела 1.1.9 раздела 1.1 этой главы 1. А именно в этом разделе планируется сосредоточить внимание на транспортно-экономических моделях, объединяющих в себе, в частности, модели из недавних работ [23, 47]. Этот раздел мотивирован обоснованием существующих и созданием новых моделей транспортного планирования, включающих модели роста транспортной инфраструктуры городов, формирования матрицы корреспонденций и равновесного распределения потоков.

Имеется ориентированный транспортный граф, каждое ребро которого характеризуется неубывающей функцией затрат $\tau_e\left(f_e\right)$ на прохождение этого ребра, в зависимости от потока по этому ребру. Можно еще ввести затраты на прохождения вершин графа $E$, но это ничего не добавляет с точки зрения последующих математических выкладок [23]. Часть вершин графа является источниками, часть стоками (эти множества вершин могут пересекаться). В источниках $O$ и стоках $D$ имеются (соответственно) пункты производства и пункты потребления. Для большей наглядности в первой половине этого раздела мы будем считать, что производится и потребляется лишь один продукт. Несложно все, что далее будет написано, обобщить на многопродуктовый рынок.

Задача разбивается на две подзадачи разного уровня [180]. На нижнем уровне, соответствующем быстрому времени, при заданных корреспонденциях $\left\{d_{ij}\right\}$ (сколько товара перевозится в единицу времени из источника $i$ в сток $j$) идет равновесное формирование способов транспортировки товаров [47]. В результате формируются функции затрат $T_{ij}\left(\left\{d_{ij}\right\}\right)$. Исходя из этих затрат на верхнем уровне, соответствующему медленному времени, решается задача поиска конкурентного равновесия [8, 213] между производителями и потребителями с учетом затрат на транспортировку. В данном случае (см., например, [47]) мы будем иметь дело с адиабатическим исключением быстрых переменных (принцип подчинения Г. Хакена) в случае стохастических динамик. Обоснование имеется в [28].

Различные частные случаи такого рода постановок задач встречались в литературе. Так, например, в классической монографии [26] рассматривается большое количество моделей верхнего уровня, связанных с расчетом матрицы корреспонденций. В монографиях [49, 112, 289], напротив, внимание сосредоточено на моделях нижнего уровня, в которых с помощью принципа Нэша–Вардропа рассчитывается $T_{ij}\left(\left\{d_{ij}\right\}\right)$. В статье [47] эти модели объединяются для создания единой многоуровневой (в транспортной науке чаще используется термин “многостадийной”) модели, включающей в себя и расчет матрицы корреспонденций, и равновесное распределение потоков по путям. В препринте [23] введена терминология, которой мы будем придерживаться и в данном разделе, связанная с пунктами производства и потребления, и в отличие от [26, 47] внешняя задача в [23] больше привязана непосредственно к экономике. Но во всех этих случаях можно было обойтись (с некоторыми оговорками в случае [47]) и без понятия конкурентного равновесия, поскольку получающиеся в итоге (популяционные) игры были потенциальными[17], причем имелась и эволюционная интерпретация [302]. Поиск равновесия сводился к решению задачи выпуклой оптимизации, а цены определялись из решения двойственной задачи. В препринте [279] для задачи верхнего уровня была предложена оригинальная конструкция, сводящая поиск конкурентного равновесия к поиску седловой точки (причем, получившаяся игра не была потенциальной в обычном смысле). Тем не менее, в [279] не рассматривалась транспортировка, т.е. не было задачи нижнего уровня.

Целью настоящего раздела является предложить такое описание задачи верхнего уровня, включающее в себя описанные выше примеры, которое сводит в итоге поиск конкурентного транспортно-экономического равновесия к поиску седловой точки в выпукло-вогнутой игре. Отметим здесь, что в общем случае поиск конкурентного равновесия сводится к решению задачи дополнительности или (при другой записи) вариационному неравенству [8, 213]. При этом известно, что в общем случае это вычислительно трудные задачи. Однако в ряде случаев экономическая специфика задачи позволяет гарантировать, что полученное вариационное неравенство монотонное. Тогда задача становится существенно привлекательнее в вычислительном плане. В данном разделе мы рассматриваем класс задач, в которых вариационные неравенства, возникающие при поиске конкурентных равновесий, переписываются в виде седловых задач. Монотонность автоматически обеспечивается правильной выпукло-вогнутой структурой седловой задачи.

Опишем вкратце структуру раздела. В подразделе 1.3.2 описывается решение “транспортной” задачи нижнего уровня (ищется равновесное распределение потоков по путям). В подразделе 1.3.3 описывается конструкция равновесного формирования корреспонденций при заданных функциях транспортных затрат. Отметим, что в этих двух пунктах мы фактически работаем только с одним экономическим агентом “Перевозчик” (если смот-реть с точки зрения популяционной теории игр, то с агентами одного типа “Перевозчиками”). В подразделе 1.3.3 этот агент(-ы) могут производить товар, неся затраты, и его потреблять, получая выгоду. В подразделе 1.3.4 модель из подраздела 1.3.3 переносится на случай, когда помимо экономического агента “Перевозчик(-и)” в источниках и стоках транспортного графа располагаются независимые от “Перевозчика” новые экономические агенты “Производители” и “Потребители”, решающие свои задачи. В заключительном подразделе 1.3.5 модели верхнего и нижнего уровня объединяются в одну общую модель, конкурентное равновесие в которой сводится к поиску седловой точки в выпукло-вогнутой игре.

---

[17] Вектор-функция затрат, характеризующая затраты при выборе различных стратегий как функция от распределения игроков по стратегиям, является градиентом некоторой скалярной функции.

### 1.3.2 Равновесное распределение потоков по путям

Обозначим множество пар $w=(i,j)$ источник-сток $OD$, $x_p$ – поток по пути $p$; $P_w$ – множество путей, отвечающих корреспонденции $w$, $P=\bigcup_{w\in OD} P_w$ – множество всех путей; $f_e(x)=\sum_{p\in P}\delta_{ep}x_p$ – поток на ребре $e$ (здесь и далее $x=\{x_p\}$, $f=\Theta x$), где $\delta_{ep}=\begin{cases}1, & e\in p\\ 0, & e\notin p\end{cases}$; $\tau_e(f_e)$ – затраты на ребре $e$ ($\tau_e'(f_e)\geq 0$); $G_p(x)=\sum_{e\in E}\tau_e(f_e(x))\delta_{ep}$ – затраты на пути $p$; $X=\left\{x\geq 0:\ \sum_{p\in P_w}x_p=d_w,\ w\in OD\right\}$ – множество допустимых распределений потоков по путям, где $d_w$ – корреспонденция, отвечающая паре $w$.

**Определение 1.3.1.** *Распределение потоков по путям* $x=\{x_p\}\in X$ *называется равновесием (Нэша–Вардропа) в популяционной игре* $\left\langle \{x_p\}\in X, \{G_p(x)\}\right\rangle$, *если из* $x_p>0$ *(*$p\in P_w$*) следует* $G_p(x)=\min_{q\in P_w}G_q(x)$. *Или, что то же самое:*

$$\text{для любых } w\in OD,\ p\in P_w \text{ выполняется } x_p\cdot\left(G_p(x)-\min_{q\in P_w}G_q(x)\right)=0.$$

**Теорема 1.3.1 (см. [47, 49, 112, 289]).** *Популяционная игра* $\left\langle \{x_p\}\in X, \{G_p(x)\}\right\rangle$ *является потенциальной. Равновесие* $x^*$ *в этой игре всегда существует, и находится из решения задачи выпуклой оптимизации*

$$x^*\in \operatorname{Arg}\min_{x\in X}\Psi(f(x)), \tag{1.3.1}$$

*где*

$$\Psi(f(x))=\sum_{e\in E}\int_0^{f_e(x)}\tau_e(z)\,dz=\sum_{e\in E}\sigma_e(f_e(x)).$$

Мы оставляем в стороне вопрос единственности равновесия (детали см., например, в [47, 49]). Отметим лишь, что при естественных условиях равновесное распределение потоков по ребрам $f^*$ единственно. В частности, для этого достаточно, чтобы $\tau_e'(f_e)>0$ для всех $e\in E$. Если дополнительно $f^*=\Theta x$ однозначно разрешимо относительно $x$ (в реальных транспортных сетях часто случается, что число допустимых для перевозки путей меньше числа ребер, это как раз и приводит к однозначной разрешимости), то отсюда будет следовать, что равновесное распределение потоков по путям $x^*$ единственно.

Удобно считать [23, 47], что возрастающие функции затрат $\tau_e(f_e):=\tau_e^\mu(f_e)$ зависят от параметра $\mu>0$, причем

$$\tau_e^\mu(f_e)\xrightarrow[\mu\to0+]{}\begin{cases}\bar t_e, & 0\le f_e<\bar f_e\\ [\bar t_e,\infty), & f_e=\bar f_e\end{cases},$$

$$d\tau_e^\mu(f_e)/df_e\xrightarrow[\mu\to0+]{}0,\ 0\le f_e<\bar f_e.$$

В таком пределе задачу выпуклой оптимизации можно переписать как задачу ЛП [275]:

$$\min_{\substack{f=\Theta x,\,x\in X\\ f\le\bar f}}\sum_{e\in E}f_e\bar t_e.$$

Такого типа транспортные задачи достаточно хорошо изучены [72, 132].

Для дальнейшего будет важно переписать задачу $\min_{x\in X}\Psi(f(x))$ через двойственную [47]:

$$\min_{x\in X}\Psi(f(x))=\min_{f,x}\left\{\sum_{e\in E}\sigma_e(f_e):\ f=\Theta x,\ x\in X\right\}=$$

$$=\min_{f,x}\left\{\sum_{e\in E}\max_{t_e\in\operatorname{dom}\sigma_e^*}\left[f_et_e-\sigma_e^*(t_e)\right]:\ f=\Theta x,\ x\in X\right\}=$$

$$=\max_{t\in\operatorname{dom}\sigma^*}\left\{\min_{f,x}\left[\sum_{e\in E}f_et_e:\ f=\Theta x,\ x\in X\right]-\sum_{e\in E}\sigma_e^*(t_e)\right\}=$$

$$=\max_{t\ge\bar t}\left\{\sum_{w\in OD}d_wT_w(t)-\left\langle\bar f,t-\bar t\right\rangle-\mu\sum_{e\in E}h\left(t_e-\bar t_e,\bar t_e,\bar f_e,\mu\right)\right\}=$$

$$\overset{\mu\to0+}{=}\max_{t\ge\bar t}\left\{\sum_{w\in OD}d_wT_w(t)-\left\langle\bar f,t-\bar t\right\rangle\right\},\qquad(1.3.2)$$

где $\sigma_e^*(t_e)$ – сопряженная функция к $\sigma_e(f_e)$, $T_w(t)=\min_{p\in P_w}\sum_{e\in E}\delta_{ep}t_e$ – длина кратчайшего пути из $i$ в $j$ ($w=(i,j)$) на графе, взвешенном согласно вектору $t$, $h(t_e-\bar{t}_e,\bar{t}_e,\bar{f}_e,\mu)$ – сильно выпуклая функция по первому аргументу. При этом

$$\tau_e^{\mu}\left(f_e\left(x(\mu)\right)\right)\xrightarrow[\mu\to 0+]{}t_e,$$

где $x(\mu)$ – равновесное распределение потоков, рассчитывающееся по формуле (1.3.1), а $t=\{t_e\}$ – решение задачи (1.3.2), при естественных условиях единственное [47]. Описанный предельный переход позволяет переходить к задачам, в которых вместо функции затрат на ребрах $\tau_e(f_e)$ заданны ограничения на пропускные способности $0\le f_e\le \bar{f}_e$ и затраты $\bar{t}_e$ на прохождения ребер, когда на ребрах нет “пробок” ( $f_e<\bar{f}_e$ ) [47, 275].

Основным для дальнейшего выводом из этого всего является способ (основанный на применении теоремы Демьянова–Данскина [61, 62], как правило, в градиентном варианте в виду единственности $t$) потенциального описания набора $T(d):=\{T_w(t(d))\}$:

$$T(d)=\nabla_d \min_{x\in X(d)}\Psi(f(x))=\nabla_d \max_{t\ge\bar{t}}\left\{\sum_{w\in OD}d_wT_w(t)-\left\langle \bar{f},t-\bar{t}\right\rangle-\mu\sum_{e\in E}h\left(t_e-\bar{t}_e,\bar{t}_e,\bar{f}_e,\mu\right)\right\}. \quad (1.3.3)$$

В [47, 49, 302] приведены эволюционные динамики, приводящие к описанным здесь равновесиям. Отметим, однако, что если рассматривать Logit динамику [47, 302] (ограниченно рациональных агентов с параметром $\tilde{\gamma}>0$ [137]), то задачу (1.3.1) необходимо будет переписать в виде (говорят, что вместо равновесия Нэша–Вардропа ищется стохастическое равновесие [47, 310]):

$$\min_{x\in X}\left\{\Psi(f(x))+\tilde{\gamma}\sum_{w\in OD}\sum_{p\in P_w}x_p\ln\left(x_p/d_w\right)\right\}. \quad (1.3.4)$$

Это замечание понадобится нам в дальнейшем.

В заключение этого раздела отметим, что теорема 1.3.1 может быть распространена и на случай, когда затраты на ребрах $\tau_e\left(f_e;\{f_{\tilde{e}}\}\right)$ удовлетворяют условию потенциальности (частный случай – это когда $\tau_e\left(f_e;\{f_{\tilde{e}}\}\right)\equiv\tau_e(f_e)$) [302]:

$$\frac{\partial\tau_e\left(f_e;\{f_{\tilde{e}}\}\right)}{\partial e'}=\frac{\partial\tau_{e'}\left(f_{e'};\{f_{\tilde{e}}\}\right)}{\partial e}.$$

Тогда

$$\Psi(f(x))=\sum_{e\in E}\int_0^{f_e(x)}\tau_e\left(z;\left\{\{f_{\tilde{e}}\}_{\tilde{e}\ne e}\cup\{f_e=z\}\right\}\right)dz.$$

Такого рода обобщение нужно, например, когда пропускные способности узлов зависят от потоков, пересекающих узлы. В случае транспортных потоков такими узлами являются, в частности, перекрестки. Тогда, путем раздутия исходного графа, мы с одной стороны “развязываем узел”, сводя затраты прохождения узла по разным путям к затратам прохождения фиктивных (введенных нами) ребер, а с другой стороны приобретаем более общую зависимость $\tau_e\left(f_e;\left\{f_{\tilde{e}}\right\}\right)$.

### 1.3.3 Равновесный расчет матрицы корреспонденций

В подразделе 1.3.2 матрица корреспонденций $\left\{d_w\right\}_{w\in OD}$ была задана по постановке задачи. В данном пункте мы откажемся от этого условия, вводя в источники $O$ производство, а в стоки $D$ – потребление. Агенты “появляются” в тех пунктах производства, произведя товар в которых его можно с выгодой для себя реализовать в каком-нибудь из пунктов потребления. Это означает, что затраты на производство и затраты на транспортировку полностью окупаются последующей выручкой от реализации продукции в пункте потребления. Агенты, которых мы здесь считаем маленькими, будут “приходить” в систему до тех пор, пока существует цепочка (пункт производства–маршрут–пункт потребления), обеспечивающая им положительную прибыль. Важно отметить, что по ходу “наплыва” агентов транспортная сеть становится все более и более загруженной, что может сказываться на затратах на перевозку. В результате прибыль ранее пришедших агентов падает, что побуждает их перераспределяться, т.е. искать более выгодные цепочки. Возникает ряд вопросов. Например, сходится ли такая динамика (точнее семейство динамик, отражающих рациональность агентов) к равновесию? Если сходится, то единственно ли равновесие? Если равновесие единственно, то как его можно эффективно найти (описать)? Попытка ответить на эти вопросы (но, прежде всего, на последний вопрос) для достаточно широкого и важного в приложениях класса задач предпринята в последующей части раздела.

Рассмотрим сначала для большей наглядности отдельно потенциальный случай. А именно тот случай, когда в источнике $i\in O$ производственная функция имеет вид $\sigma_i\left(f_i\right)$, где $f_i=\sum_{k:(i,k)=e\in E} f_e=\sum_{j:(i,j)\in OD} d_{ij}$, аналогично для стоков $j\in D$ определим функции полезности со знаком минус $\sigma_j\left(f_j\right)$, $f_j=\sum_{k:(k,j)=e\in E} f_e=\sum_{i:(i,j)\in OD} d_{ij}$. Все эти функции считаем выпуклыми. Мы обозначаем эти функции одинаковыми буквами, однако, это не должно вызвать в дальнейшем путаницы в виду характерных нижних индексов. Редуцируем рассматриваемую задачу к задаче подраздела 1.3.2. Рассмотрим новый граф с множеством вершин $O \bigcup D$, соединенных теми же ребрами, что и в изначальном графе, и с одним дополнительным фиктивным источником и одним дополнительным фиктивным стоком. Этот фиктивный источник соединим со всеми источниками $O$, аналогично фиктивный сток соединим со всеми стоками $D$. Если существует путь из источника $i$ в сток $j$ в исходном графе, то в новом графе прочертим соответствующее ребро с функцией затрат $T_{ij}(d)$. Проведем дополнительное (фиктивное) ребро, соединяющее фиктивный источник с фиктивным стоком, затраты на прохождения которого тождественный ноль. Получим в итоге ориентированный граф путей из источника в сток. Легко понять, что мы оказываемся “почти” в условиях предыдущего пункта (причем с более частным графом – с одним источником и стоком) с точностью до обозначений:

$$\{x_p\} \to \{d_{ij}\},\ \{\tau_e(f_e)\} \to \{\sigma_i'(f_i), T_{ij}(d), \sigma_j'(f_j)\}.$$

“Почти”, потому что, во-первых, затраты $T_{ij}(d)$ зависят от всего набора $\{d_{ij}\}$, а не только от $d_{ij}$, а во-вторых, не ясно что в данном случае играет роль матрицы корреспонденций (в нашем случае это матрица $1 \times 1$, т.е. просто число). Начнем с ответа на второй вопрос. Мы считаем, что в источниках имеется потенциальная возможность производить неограниченное количество продукта, просто в какой-то момент, перестает быть выгодным что-то производить и перевозить. Для этого, собственно, и было введено нулевое ребро, поток по которому обозначим $d_0$. То есть, другими словами, мы должны считать, что $\sum_{(i,j)\in W} d_{ij} + d_0 = \bar{d}$. Если $\bar{d}$ – достаточно большое, то равновесная конфигурация не зависит от того, чему именно равно $\bar{d}$, поскольку не требуется определять $d_0$. С первой проблемой можно разобраться, немного обобщив теорему 1.3.1. Предположим, что

$$\exists\ \Phi(d) \text{ - выпуклая}:\ T(d) = \nabla\Phi(d). \qquad (1.3.5)$$

Тогда имеет место

**Теорема 1.3.2.** *Популяционная игра*

$$\left\langle \{d_{ij}, d_0 \ge 0\}, \{G_{ij}(d) = \sigma_i'(f_i) + T_{ij}(d) + \sigma_j'(f_j), G_0(d) \equiv 0\} \right\rangle,$$

*является потенциальной. Равновесие* $d^*$ *в этой игре всегда существует (если* $\sigma(\cdot)$ *– сильно выпуклые функции, то равновесие гарантировано единственно), и находится из решения задачи выпуклой оптимизации*

$$d^* \in \arg\min_{d \ge 0} \tilde{\Psi}(d),$$

$$\tilde{\Psi}\left(d=\{d_{ij}\}\right)=\sum_{i\in O}\sigma_i\left(\sum_{j:(i,j)\in OD}d_{ij}\right)+\sum_{j\in D}\sigma_j\left(\sum_{i:(i,j)\in OD}d_{ij}\right)+\Phi(d). \quad (1.3.6)$$

**Доказательство.** Выпишем условие нелинейной комплиментарности (то есть, по сути, определения равновесия Нэша в популяционной игре, заданной в условии). Для этого занумеруем все индексы $ij$ и $0$ одним индексом $k$:

для любых $k$ выполняется $d_k^*\cdot\left(G_k\left(d^*\right)-\min_{k'}G_{k'}\left(d^*\right)\right)=0.$

Действительно допустим, что реализовалось какое-то другое равновесие $\tilde{d}$, которое не удовлетворяет этому условию. Покажем, что тогда найдется агент, которому выгодно поменять свой выбор. Действительно, тогда существуют такой $\tilde{k}$, что $\tilde{d}_{\tilde{k}}\cdot\left(G_{\tilde{k}}\left(\tilde{d}\right)-\min_{k'}G_{k'}\left(\tilde{d}\right)\right)>0$.

Каждый агент (множество таких агентов не пусто $\tilde{d}_{\tilde{k}}>0$) использующий стратегию $\tilde{k}$ действует не разумно, поскольку существует такая стратегия $\bar{k}$, $\bar{k}\neq\tilde{k}$, что $G_{\bar{k}}\left(\tilde{d}\right)=\min_{k'}G_{k'}\left(\tilde{d}\right)$. Этот стратегия $\bar{k}$ более выгодна, чем $\tilde{k}$. Аналогично показывается, что при распределении $d^*$ никому агентов уже не выгодно отклоняться от своих стратегий.

Покажем, что рассматриваемая нами игра принадлежит к классу, так называемых, потенциальных игр. В нашем случае это означает, что существует такая функция $\tilde{\Psi}\left(d=\{d_{ij}\}\right)$, что $\partial\tilde{\Psi}(d)/\partial d_k=G_k(d)$ для любого $k$. В нашем случае легко проверить, что этому условию удовлетворяет (по условию) функция

$$\tilde{\Psi}\left(d=\{d_{ij}\}\right)=\sum_{i\in O}\sigma_i\left(\sum_{j:(i,j)\in OD}d_{ij}\right)+\sum_{j\in D}\sigma_j\left(\sum_{i:(i,j)\in OD}d_{ij}\right)+\Phi(d).$$

Таким образом, мы имеем дело с потенциальной игрой. Оказывается, что $d^*$ – равновесие Нэша–Вардропа тогда и только тогда, когда оно доставляет минимум $\tilde{\Psi}(d)$ на множестве $d\geq 0$. Действительно, предположим, что $d^*\geq 0$ – точка минимума. Тогда, в частности, для любых $p,q$ ($d_p^*>0$) и достаточно маленького $\delta d_p>0$ выполняется:

$$-\frac{\partial\tilde{\Psi}\left(d^*\right)}{\partial d_p}\delta d_p+\frac{\partial\tilde{\Psi}\left(d^*\right)}{\partial d_q}\delta d_p\geq 0.$$

Иначе, заменив $d^*$ на

$$\breve{d}^* = d^* + \Big( \underbrace{\overbrace{0,...,0,-\delta d_p}^{p},0,...,0,\delta d_p}_{q},0,...,0 \Big) \geq 0,$$

мы пришли бы к вектору $\breve{d}^*$, доставляющему меньшее значение $\tilde{\Psi}(d)$ на множестве $d \geq 0$:

$$\tilde{\Psi}\left(\breve{d}^*\right) \approx \tilde{\Psi}\left(d^*\right) - \frac{\partial \tilde{\Psi}\left(d^*\right)}{\partial d_p}\delta d_p + \frac{\partial \tilde{\Psi}\left(d^*\right)}{\partial d_q}\delta d_p < \tilde{\Psi}\left(d^*\right).$$

Вспоминая, что $\partial \tilde{\Psi}(d)/\partial d_p = G_p(d)$, и учитывая, что $q$ можно выбирать произвольно, получаем:

для любого $p$, если $d_p^* > 0$, то выполняется $\min\limits_q G_q\left(d^*\right) \geq G_p\left(d^*\right)$.

Но это и есть по-другому записанное условие нелинейной комплементарности. Строго говоря, мы показали сейчас только то, что точка минимума $\tilde{\Psi}(d)$ на множестве $d \geq 0$ будет равновесием Нэша–Вардропа. Аналогично рассуждая, можно показать и обратное: равновесие Нэша–Вардропа доставляет минимум $\tilde{\Psi}(d)$ на множестве $d \geq 0$. □

Именно такая конструкция и была рассмотрена в препринте [23] (для многопродуктового рынка). Если искать стохастическое равновесие, то функционал в теореме 1.3.2 необходимо энтропийно регуляризовать[18]. Такие конструкции рассматривались (в нерегуляризованном случае), например, в работах [26, 47, 49]. Уже в этих работах можно углядеть необходимость искусственного введения потенциалов (двойственных множителей) в сами функции $\sigma$. А именно, в этих работах предполагается, что все эти функции $\sigma$ – линейные с неизвестными наклонами. Тем не менее, считается, что при этом известно, чему должны равняться в равновесии следующие суммы (см. также раздел 2 этой главы):

$$\sum_{j:(i,j)\in OD} d_{ij} = L_i, \quad \sum_{i:(i,j)\in OD} d_{ij} = W_j \quad \left(\sum_{i\in O} L_i = \sum_{j\in D} W_j = N\right). \qquad (1.3.7)$$

То есть имеются скрытые от нас (модельера) потенциалы [72] (параметры) $\left\{\lambda_i^L, \lambda_j^W\right\}$, которые могут быть рассчитаны исходя из дополнительной информации. Применительно к модели расчета матрицы корреспонденций [26, 47, 49] выписанные дополнительные условия (1.3.7) однозначно определяют все неизвестные потенциалы. Однако при этом вместо задачи выпуклой оптимизации мы получаем минимаксную (седловую) задачу выпуклую по $\left\{d_{ij}\right\} \geq 0$ и вогнутую, точнее даже линейную, по потенциалам $\left\{\lambda_i^L, \lambda_j^W\right\}$:

[18] К сожалению, строгое обоснование (теорема 11.5.12 [302]) имеется только в случае известного (фиксированного) значения $\bar{d}$ (при этом можно считать $d_0 = 0$).

$$\min_{\substack{\{d_{ij}\}\geq 0 \\ \sum_{(i,j)\in W} d_{ij}=N}} \max_{\{\lambda_i^L,\lambda_j^W\}} \left[ \sum_{i\in O} \lambda_i^L \cdot \left( \sum_{j:(i,j)\in OD} d_{ij} - L_i \right) + \sum_{j\in D} \lambda_j^W \cdot \left( W_j - \sum_{i:(i,j)\in OD} d_{ij} \right) + \right.$$

$$\left. + \Phi(d) + \gamma \sum_{(i,j)\in OD} d_{ij} \ln\left(d_{ij}/N\right) \right]. \qquad (1.3.8)$$

Эта задача всегда имеет решение.

### 1.3.4 Сетевая модель алгоритмического рыночного поведения

В данном пункте мы предложим сетевой вариант модели поиска конкурентного равновесия из препринта [279]. Однако в контексте изложенного в конце прошлого пункта, нам будет удобнее стартовать с двухстадийной модели транспортных потоков [47], приводящей к равновесию, рассчитываемому по формуле (1.3.8).

Предположим теперь, что имеется $m$ видов товаров и, дополнительно, имеется $q$ типов материалов (количества которых можно использовать в единицу времени ограничены вектором $b$), использующихся в производстве. В источниках располагаются производители товаров, а в стоках потребители. Мы считаем, что любой производитель товара, одновременно, является и потребителем, т.е. $O \subseteq D$. Обозначим через $y$ – вектор цен (руб.) на материалы; $\lambda_i^L$ – вектор цен (руб.), по которым производитель продает товары перевозчику в пункте производства *i*, $\lambda_j^W$ – вектор цен (руб.), по которым потребитель покупает товары у перевозчика в пункте потребления *j*. Опишем каждого экономического агента:

***i*-й Производитель**

$U_i \subset \mathbb{R}_+^m$ – максимальное технологическое множество (замкнутое, выпуклое);

$\alpha_i \in [0,1]$– уровень участия;

$\chi_i(\alpha_i U_i) = \alpha_i \chi_i(U_i)$ – постоянные технологические производственные затраты (руб.) при уровне участия $\alpha_i$ (в единицу времени);

$L_i \in \alpha_i U_i$ , $[L_i]_k$ – количество произведенного *k*-го продукта (в единицу времени);

$A_i$, $[A_i]_{kl}$ – количество затраченного *l*-го продукта при производстве единицы *k*-го продукта;

$c_i$, $[c_i]_k$ – затраты (руб.) на производство единицы *k*-го продукта;

$R_i$, $[R_i]_{kl}$ – количество затраченного *k*-го материала для приготовления единицы *l*-го продукта.

Описанный “Производитель” решает задачу:

$$\max_{\substack{L_i\in\alpha_i U_i\\ \alpha_i\in[0,1]}}\left\{\left\langle\lambda_i^L,L_i\right\rangle-\chi_i\left(\alpha_i U_i\right)-\left\langle\lambda_i^W,A_i L_i\right\rangle-\left\langle c_i,L_i\right\rangle-\left\langle y,R_i L_i\right\rangle\right\}=$$

$$=\max_{L_i\in U_i}\left\{\left(\left\langle\lambda_i^L-c^i-A_i^T\lambda_i^W-R_i^T y,L_i\right\rangle-\chi_i\left(U_i\right)\right)_+\right\}.$$

***j*-й Потребитель**

Предположим, что каждый продукт имеет *s* различных свойств (своеобразных полезностей). Это может быть, например, содержание витаминов, белков, жиров, углеводов.

$Q_j$, $\left[Q_j\right]_{kl}$ – вклад единицы *l*-го продукта в удовлетворение *k*-го свойства;

$\sigma_j$, $\left[\sigma_j\right]_k$ – минимально допустимый уровень удовлетворения *k*-го свойства (в единицу времени);

$\beta_j\in[0,1]$– уровень участия;

$V_j=\left\{W_j\in\mathbb{R}_+^m:\ Q_jW_j\ge\sigma_j\right\}$ – допустимое множество при полном участии;

$W_j\in\beta_jV_j$ , $\left[W_j\right]_k$ – количество потребленного *k*-го продукта (в единицу времени);

$\tau_j$ – постоянный доход (руб.) при полном участии (в единицу времени).

Описанный “Потребитель” решает задачу:

$$\max_{\substack{W_j\in\beta_jV_j\\ \beta_j\in[0,1]}}\left\{\beta_j\tau_j-\left\langle\lambda_j^W,W_j\right\rangle\right\}=\max_{W_j\in V_j}\left\{\left(\tau_j-\left\langle\lambda_j^W,W_j\right\rangle\right)_+\right\}.$$

**Перевозчик**

Этот агент решает задачу типа (1.3.6), (1.3.8):

$$\min_{\{d_{ij}\}\ge0}\left[\sum_{i\in O}\left\langle\lambda_i^L,\sum_{j:(i,j)\in OD}d_{ij}\right\rangle-\sum_{j\in D}\left\langle\lambda_j^W,\sum_{i:(i,j)\in OD}d_{ij}\right\rangle+\Phi(d)+\gamma\sum_{(i,j)\in OD}\left(\sum_{k=1}^m\left[d_{ij}\right]_k\right)\ln\left(\sum_{k=1}^m\left[d_{ij}\right]_k\right)\right],$$

в которой корреспонденции формируются “Производителями” и “Потребителями”. Мы считаем, что все товары одинаковы с точки зрения “Перевозчика”, т.е. $\Phi(d):=\Phi\left(\left\{\sum_{k=1}^m\left[d_{ij}\right]_k\right\}\right)$ (можно рассматривать и другие варианты). Здесь и далее нам будет удобно писать энтропийную регуляризацию в виде $\gamma\sum_{(i,j)\in OD}\left(\sum_{k=1}^m\left[d_{ij}\right]_k\right)\ln\left(\sum_{k=1}^m\left[d_{ij}\right]_k\right)$, т.е. опускать $N=\sum_{(i,j)\in OD}\sum_{k=1}^m\left[d_{ij}\right]_k$. Точнее полагать $N=1$ с той же потоковой (физической) размерностью, что и $d$, чтобы под логарифмом была безразмерная величина.[19] При естественных балансовых условиях $\sum_{i\in O} L_i = \sum_{i\in O} A_i L_i + \sum_{i\in D} W_j$ это никак не повлияет на решение задачи.

Проблема здесь в том, что все эти три типа задач завязаны друг на друга посредством векторов цен. Выпишем, как это принято при поиске конкурентных равновесий [8, 213], все имеющиеся **законы Вальраса** (балансовые ограничения + условия дополняющей нежесткости), которые накладывают совершенно естественные ограничения на эти векторы цен:

$$\sum_{j:(i,j)\in OD} d_{ij} \le L_i, \ \left\langle \lambda_i^L, \sum_{j:(i,j)\in OD} d_{ij} - L_i \right\rangle = 0, \ \lambda_i^L \ge 0;$$

$$\sum_{k:(k,i)\in OD} d_{ki} \ge W_i + A_i L_i, \ \left\langle \lambda_i^W, W_i + A_i L_i - \sum_{k:(k,i)\in OD} d_{ki} \right\rangle = 0, \ \lambda_i^W \ge 0, \ i \in O;$$

$$\sum_{i:(i,j)\in OD} d_{ij} \ge W_j, \ \left\langle \lambda_j^W, W_j - \sum_{i:(i,j)\in OD} d_{ij} \right\rangle = 0, \ \lambda_j^W \ge 0, \ j \in D \setminus O;$$

$$\sum_{i\in O} R_i L_i \le b, \ \left\langle y, b - \sum_{i\in O} R_i L_i \right\rangle = 0, \ y \ge 0.$$

**Определение 1.3.2.** *Набор* $\left\langle \{d_{ij}\}, \{L_i\}, \{W_j\}; y, \{\lambda_i^L\}, \{\lambda_j^W\} \right\rangle$ *называется конкурентным равновесием (Вальраса–Нестерова–Шихмана) если он доставляет решения задачам всех агентов и удовлетворяет законам Вальраса.*

Для того чтобы установить корректность этого определения, подобно [279], введем **условие продуктивности**:

*существуют такие* $\bar{L}_i \in U_i$, $\bar{W}_j \in V_j$, *что* $\sum_{i\in O} \bar{L}_i > \sum_{i\in O} A_i \bar{L}_i + \sum_{i\in D} \bar{W}_j$ *и* $\sum_{i\in O} R_i \bar{L}_i < b$.

**Теорема 1.3.3.** *В условиях продуктивности конкурентное равновесие существует и находится из решения правильной выпукло-вогнутой седловой задачи:*

$$\min_{\{d_{ij}\}\ge 0} \max_{\substack{\{\lambda_i^L, \lambda_j^W\}\ge 0 \\ y \ge 0}} \min_{\substack{\{L_i \in \alpha_i U_i,\ \alpha_i \in [0,1]\} \\ \{W_j \in \beta_j V_j,\ \beta_j \in [0,1]\}}} \left[ \sum_{i\in O} \left( \left\langle \lambda_i^L, \sum_{j:(i,j)\in OD} d_{ij} - L_i \right\rangle + \left\langle \lambda_i^W, \sum_{k:(k,i)\in OD} A_i L_i \right\rangle + \chi_i(\alpha_i U_i) \right) + \right.$$

[19] В случае микроскопического обоснования такого рода вариационных принципов (см. подраздел 1.3.5, а также [302]) исходя из рассмотрения соответствующей марковской динамики нащупывания равновесия, мы должны полагать $N \gg 1$, чтобы сделать соответствующий (канонический) скейлинг и перейти к детерминированной постановке.

$$+\sum_{j\in D}\left(\left\langle \lambda_j^W, W_j - \sum_{i:(i,j)\in OD} d_{ij}\right\rangle - \beta_j \tau_j\right) + \left\langle y, b - \sum_{i\in O} R_i L_i\right\rangle +$$

$$\left.+\Phi(d) + \gamma \sum_{(i,j)\in OD}\left(\sum_{k=1}^{m}\left[d_{ij}\right]_k\right)\ln\left(\sum_{k=1}^{m}\left[d_{ij}\right]_k\right)\right] =$$

$$= \max_{\substack{\{\lambda_i^L, \lambda_j^W\}\geq 0\\ y\geq 0}}\left[\langle y, b\rangle - \sum_{i\in O}\max_{L_i\in U_i}\left\{\left(\left\langle \lambda_i^L - c^i - A_i^T\lambda_i^W - R_i^T y, L_i\right\rangle - \chi_i(U_i)\right)_+\right\} -\right.$$

$$-\sum_{j\in D}\max_{W_j\in V_j}\left\{\left(\tau_j - \left\langle \lambda_j^W, W_j\right\rangle\right)_+\right\} + \min_{\{d_{ij}\}\geq 0}\left\{\sum_{i\in O}\left\langle \lambda_i^L, \sum_{j:(i,j)\in OD} d_{ij}\right\rangle - \sum_{j\in D}\left\langle \lambda_j^W, \sum_{i:(i,j)\in OD} d_{ij}\right\rangle +\right.$$

$$\left.\left.+\Phi(d) + \gamma \sum_{(i,j)\in OD}\left(\sum_{k=1}^{m}\left[d_{ij}\right]_k\right)\ln\left(\sum_{k=1}^{m}\left[d_{ij}\right]_k\right)\right\}\right]. \qquad (1.3.9)$$

Порядок взятия минимума и максимумов можно менять согласно “Sion’s minimax theorem” [47, 61, 62, 313].

### 1.3.5 Общее конкурентное равновесие

Для того чтобы объединить модели подразделов 1.3.2, 1.3.3, 1.3.4 в одну модель рассмотрим формулы (1.3.3), (1.3.5), (1.3.8), (1.3.9). Легко понять, что формула (1.3.3) как раз и задает тот самый потенциал, существование которого (формула (1.3.5)) требуется для справедливости теоремы 1.3.2 (неявно это предполагается и в теореме 1.3.3), фактически, сводящей поиск конкурентного равновесия к задаче (1.3.8), а в общем случае (1.3.9).

**Определение 1.3.3.** *Набор* $\left\langle \{x_p\}, \{d_{ij}\}, \{L_i\}, \{W_j\}; y, \{\lambda_i^L\}, \{\lambda_j^W\}\right\rangle$ *называется полным (общим) конкурентным равновесием (Вальраса–Нэша–Вардропа–Нестерова–Шихмана) если* $\left\langle \{d_{ij}\}, \{L_i\}, \{W_j\}; y, \{\lambda_i^L\}, \{\lambda_j^W\}\right\rangle$ *– конкурентное равновесие, а* $\{x_p\}$ *является равновесием (Нэша–Вардропа) при заданном конкурентным равновесием наборе* $\{d_{ij}\}$.

**Теорема 1.3.4.** *В условиях продуктивности полное конкурентное равновесие существует и находится из решения правильной выпукло-вогнутой седловой задачи:*

$$\max_{\substack{\{\lambda_i^L, \lambda_j^W\}\geq 0\\ y\geq 0}}\left[\langle y, b\rangle - \sum_{i\in O}\max_{L_i\in U_i}\left\{\left(\left\langle \lambda_i^L - c^i - A_i^T\lambda_i^W - R_i^T y, L_i\right\rangle - \chi_i(U_i)\right)_+\right\} -\right.$$

$$-\sum_{j\in D}\max_{W_j\in V_j}\left\{\left(\tau_j - \left\langle \lambda_j^W, W_j\right\rangle\right)_+\right\} + \min_{\{d_{ij}\}\geq 0}\left\{\sum_{i\in O}\left\langle \lambda_i^L, \sum_{j:(i,j)\in OD} d_{ij}\right\rangle - \sum_{j\in D}\left\langle \lambda_j^W, \sum_{i:(i,j)\in OD} d_{ij}\right\rangle +\right.$$

$$+\max_{t\geq \bar{t}}\left\{\sum_{(i,j)\in OD}\left(\sum_{k=1}^{m}\left[d_{ij}\right]_k\right)T_{ij}(t) - \left\langle \bar{f}, t-\bar{t}\right\rangle - \mu\sum_{e\in E} h\left(t_e - \bar{t}_e, \bar{f}_e, \mu\right)\right\} +$$

$$+\gamma \sum_{(i,j)\in OD} \left( \sum_{k=1}^{m} \left[ d_{ij} \right]_k \right) \ln \left( \sum_{k=1}^{m} \left[ d_{ij} \right]_k \right) \Bigg\} \Bigg] . \qquad (1.3.10)$$

Таким образом, поиск общего конкурентного равновесия также сводится к седловой задаче (если мы вынесем все маскимумы и минимумы за квадратную скобку, то получим минимаксную = седловую задачу), имеющей правильную структуру с точки зрения того, что минимум берется по переменным, по которым выражение в квадратных скобках выпукло, а максимум по переменным, по которым выражение вогнуто. Порядок взятия всех максимумов и минимума можно менять согласно “Sion’s minimax theorem”. В частности, это дает возможность явно выполнить минимизацию по $\{d_{ij}\} \geq 0$, “заплатив” за этого некоторым усложнением получившегося в итоге функционала, который также сохранит правильные выпукло-вогнутые свойства [47].

Мы не будем здесь приводить, что получается после подстановки формулы (1.3.3) в формулы (1.3.6) и (1.3.8). Все выкладки аналогичны, и даже проще. Тем не менее, ссылаясь далее на задачи (1.3.6) и (1.3.8) мы будем считать, что такая подстановка была сделана.

Такого рода задачи можно эффективно численно решать (причем содержательно интерпретируемым способом), если транспортный граф задачи нижнего уровня (поиска равновесного распределения потоков) не сверх большой [260, 269, 271, 280]. Если же этот граф имеет, скажем, порядка $10^5$ ребер, как транспортный граф Москвы и области [49], то требуется разработка новых эффективных методов, учитывающих разреженность задачи и использующих рандомизацию. Мы не будем здесь на этом останавливаться, поскольку планируется посвятить численным методам решения таких задач больших размеров отдельную публикацию. Впрочем, некоторые возможные подходы отчасти освещены в [23, 47]. К сожалению, численный метод, предложенный в [23], не совсем корректен.

Сделаем несколько замечаний в связи с полученным результатом.

Во-первых, в рассмотренных в разделе задачах с помощью штрафных механизмов (типа платных дорог) можно добиваться, чтобы возникающие равновесия соответствовали социальному оптимуму. Для этого можно использовать VCG-механизм [134, 301], см. также раздел 1.1 этой главы 1.

Во-вторых, используя аппарат [23, 47, 279] несложно вычленить из выписанных задач (1.3.6), (1.3.8), (1.3.10) всевозможные цены, тарифы, длины очередей (пробок) – если делаем предельный переход $\mu \to 0+$ и т.п., понимаемые в смысле Л.В. Канторовича, как двойственные множители.

В-третьих, рассматривая два разномасштабных по времени марковских процесса нащупывания равновесной конфигурации можно прийти к решению задач (1.3.6), (1.3.8) и, с некоторыми оговорками, (1.3.10). Например, если в быстром времени динамика перераспределения потоков по путям задается имитационной Logit динамикой [302], а в медленном времени процесс перераспределения корреспонденций (исходя из быстро подстраивающихся затрат $\{T_{ij}(d)\}$) задается просто Logit динамикой [302], то выражение в квадратных скобках (1.3.8)) будет играть роль действия в теореме типа Санова, т.е. описывать экспоненциальную концентрацию инвариантной меры марковского процесса с оговоркой, что речь идет о переменных $d$ и $x$ [302]. Аналогичное можно сказать и про $\tilde{\Psi}$ в теореме 1.3.2 после подстановки (1.3.3). Кроме того, эти же самые функции будут играть роль функций Ляпунова соответствующих прошкалированным (каноническим скейлингом) марковских динамик, приводящих к СОДУ Тихоновского типа [21, 47, 302]. Это также следует из общих результатов работы [10, 32] (см. также приложение в конце диссертации). Отметим, что относительно нащупывания цен (потенциалов) в задачах (1.3.8) и, особенно, (1.3.10) работают механизмы похожие на те, которые описаны в классической работе [72]. Другими словами, при фиксированных ценах (потенциалах) динамика соответствует классическим популяционным динамикам нащупывания равновесий [302]. Но из-за того, что потенциалы не известны и, в свою очередь, должны как-то параллельно подбираться предполагается, что в медленном времени экономические агенты переоценивают эти потенциалы исходя из обратной связи (пример имеется в [72]) на то, что они ожидают видеть и то, что они реально видят.

В-четвертых, упомянутая выше эволюционная динамика при правильно дискретизации дает разумный численный способ поиска конкурентного равновесия. В частности, упоминаемая имитационная Logit динамика при правильной дискретизации даст метод зеркального спуска / двойственных усреднений, представляющий собой метод проекции градиента с усреднением [269] и без [280], где проекция понимается в смысле "расстояния" Кульбака–Лейблера. Зеркальный спуск можно получить также из дискретизации Logit динамики, если ориентироваться не только на предыдущую итерацию, а на среднее арифметическое всех предыдущих итераций [47, 54]. В работе [54] поясняется некоторая привилегированность этих двух Logit динамик (см. также [47, 51, 137, 302] и подраздел 1.1.4 раздела 1.1 этой главы 1). Отметим при этом, что Logit динамики может быть проинтерпретирована также как имитационная Logit динамика для потенциальной игры с энтропийно регуляризованным потенциалом [54] (см. также подраздел 1.1.4 раздела 1.1 этой главы 1).

В-пятых, везде выше мы исходили из того, что есть разные масштабы времени. Из-за этого задачи подразделов 1.3.2 – 1.3.4 удалось завязать с помощью формул (1.3.3), (1.3.5). Однако к аналогичным выводам можно было прийти, если вместо введения разных масштабов времени ввести иерархию в принятии решений [137]. Скажем, сначала пользователь транспортной сети выбирает тип транспорта (личный/общественный), а потом маршрут [47]. Здесь особенно актуальным становятся такие модели дискретного выбора как Nested Logit [137]. А именно, если использовать энтропийную регуляризацию только в одной из этих двух задач разного уровня (иерархии), описанных в подразделах 1.3.2, 1.3.3, то получается обычная (Multinomial) Logit модель выбора (например, в (1.3.10) мы регуляризовали только задачу верхнего уровня), но если энтропийно регуляризовать обе задачи, то получится двухуровневая Nested Logit модель выбора [137]. Это означает, что соответствующая Nested Logit динамика в популяционной иерархической игре приводит к равновесию, которое описывается решением задач типа (1.3.8), (1.3.10) с дополнительной энтропийной регуляризацией задачи нижнего уровня. Несложно показать, что хорошие выпукло-вогнутые свойства задач (1.3.8), (1.3.10) при этом сохраняются. Да и в вычислительном плане задача не становится принципиально сложнее, особенно если учесть конструкцию "The shortest path problem", описанную в пятой главе монографии [250], см. также [47]. Подробнее об этом будет написано в следующем разделе.

Резюмируем полученные в разделе результаты. На конкретных семействах примеров (но, тем не менее, достаточно богатых в смысле встречаемости в приложениях), была продемонстрирована некоторая "алгебра" над различными конструкциями равновесия. Было продемонстрировано, как их можно сочетать друг с другом, чтобы получать все более и более содержательные задачи. Ключевым местом стал переход, связанный с формулой (1.3.3), который можно понимать как произведение (суперпозицию) транспортно-экономических моделей, и конструкция задач (1.3.8), (1.3.9), которую можно понимать как "сумму" моделей. Представляется, что в этом направлении, может возникнуть довольно интересное движение, связанное с вычленением той "минимальной алгебры операций" над моделями, с помощью которой можно было бы описывать большое семейство равновесных конфигураций, встречающихся в различных приложениях.

## 1.4 О связи моделей дискретного выбора с разномасштабными по времени популяционными играми загрузок

### 1.4.1 Введение

В работах [9, 31, 47, 49], а также в подразделе 1.1.10 раздела 1.1 этой главы 1 было анонсировано, что в последующем цикле публикаций будет приведен общий способ вариационного описания равновесий (стохастических равновесий) в популярных моделях распределения транспортных потоков. Также отмечалось, что планируется предложить эффективные численные методы поиска таких равновесий. В данном разделе предпринята попытка погрузить известные нам подходы к многостадийному моделированию потоков на иерархических сетях (реальных транспортных сетях или сетях принятия решение – не важно) в одну общую схему, сводящую поиск равновесия к решению многоуровневой задачи выпуклой оптимизации. В основе схемы получения вариационного принципа для описания равновесия лежит популяционная игра загрузки с соответствующими логит динамиками (отвечающими моделям дискретного выбора [137]) пользователей на каждом уровне иерархии [302] (см. подраздел 1.4.2). Для решения описанной задачи выпуклой оптимизации в разделе изучается двойственная задача, представляющая самостоятельный интерес (см. подраздел 1.4.3). Основным инструментом изучения двойственной задачи является аппарат характеристических функций на графе [33, 49, 262] и ускоренные прямо-двойственные методы в композитном варианте [266, 277].

Отметим, что общность результатов раздела достигается за счет введения большого числа параметров, которые можно вырождать, стремя их к нулю или бесконечности. Игра на выборе этих параметров позволяет, например, получать различные многостадийные модели транспортных потоков [9, 23, 31, 47, 287]. Приводимые далее результаты можно обобщать и на потоки товаров, в случае, когда имеется более одного наименования товара [23]. Однако в данном разделе мы не будем касаться этого обобщения. Мы также не планируем приводить конкретные примеры получения многостадийных транспортных моделей согласно изложенной в разделе общей схеме.

### 1.4.2 Постановка задачи

Рассмотрим транспортную сеть, заданную ориентированным графом $\Gamma^1 = \left\langle V^1, E^1 \right\rangle$. Часть его вершин $O^1 \subseteq V^1$ является источниками, часть стоками $D^1 \subseteq V^1$. Множество пар источник-сток, обозначим $OD^1 \subseteq O^1 \otimes D^1$. Пусть каждой паре $w^1 \in OD^1$ соответствует своя корреспонденция: $d_{w^1}^1 := d_{w^1}^1 \cdot M$ ($M \gg 1$) пользователей, которые хотят в единицу времени перемещаться из источника в сток, соответствующих заданной корреспонденции $w^1$.

Пусть ребра $\Gamma^1$ разделены на два типа $E^1 = \tilde{E}^1 \coprod \bar{E}^1$. Ребра типа $\tilde{E}^1$ характеризуются неубывающими функциями затрат $\tau^1_{e^1}\left(f^1_{e^1}\right) := \tau^1_{e^1}\left(f^1_{e^1}/M\right)$. Затраты $\tau^1_{e^1}\left(f^1_{e^1}\right)$ несут те пользователи, которые используют в своем пути ребро $e^1 \in \tilde{E}^1$, в предположении, что поток пользователей по этому ребру равен $f^1_{e^1}$. Пары вершин, задающие ребра типа $\bar{E}^1$, являются, в свою очередь, парами источник-сток $OD^2$ (с корреспонденциями $d^2_{w^2} = f^1_{e^1}$, $w^2 = e^1 \in \bar{E}_1$) в транспортной сети следующего уровня, $\Gamma^2 = \left\langle V^2, E^2 \right\rangle$, ребра которой, в свою очередь, разделены на два типа $E^2 = \tilde{E}^2 \coprod \bar{E}^2$. Ребра типа $\tilde{E}^2$ характеризуются неубывающими функциями затрат $\tau^2_{e^2}\left(f^2_{e^2}\right) := \tau^2_{e^2}\left(f^2_{e^2}/M\right)$. Затраты $\tau^2_{e^2}\left(f^2_{e^2}\right)$ несут те пользователи, которые используют в своем пути ребро $e^2 \in \tilde{E}^2$, в предположении, что поток пользователей по этому ребру равен $f^2_{e^2}$. Пары вершин, задающие ребра типа $\bar{E}^2$, являются, в свою очередь, парами источник-сток $OD^3$ (с корреспонденциями $d^3_{w^3} = f^2_{e^2}$, $w^3 = e^2 \in \bar{E}^2$) в транспортной сети более высокого уровня, $\Gamma^3 = \left\langle V^3, E^3 \right\rangle$, … и т.д. Будем считать, что всего имеется $m$ уровней: $\tilde{E}^m = E^m$. Обычно в приложениях число $m$ небольшое [9, 31, 47, 49, 137]: 2 – 10.

Каждый пользователь в графе $\Gamma^1$ выбирает путь $p^1_{w^1} \in P^1_{w^1}$ (последовательный набор проходимых пользователем ребер), соответствующий его корреспонденции $w^1 \in OD^1$ ($P^1_{w^1}$ – множество всех путей, отвечающих в $\Gamma^1$ корреспонденции $w^1$). Задав $p^1_{w^1}$ можно однозначно восстановить ребра типа $\bar{E}^1$, входящие в этот путь. На каждом из этих ребер $w^2 \in \bar{E}^1$ пользователь может выбирать свой путь $p^2_{w^2} \in P^2_{w^2}$ ($P^2_{w^2}$ – множество всех путей, отвечающих в $\Gamma^2$ корреспонденции $w^2$), … и т.д. Пусть каждый пользователь сделал свой выбор. Обозначим через $x^1_{p^1}$ величину потока пользователей по пути $p^1 \in P^1 = \coprod_{w^1 \in OD^1} P^1_{w^1}$, $x^2_{p^2}$ – величина потока пользователей по пути $p^2 \in P^2 = \coprod_{w^2 \in OD^2} P^2_{w^2}$, … и т.д. Заметим, что

$$x^k_{p^k_{w^k}} \geq 0,\ p^k_{w^k} \in P^k_{w^k},\ \sum_{p^k_{w^k} \in P^k_{w^k}} x^k_{p^k_{w^k}} = d^k_{w^k},\ w^k \in OD^k,\ k = 1,\dots,m,$$

что для компактности мы будем далее записывать

$$\left\{x^k_{p^k_{w^k}}\right\}_{p^k_{w^k} \in P^k_{w^k}} \in S_{\left|P^k_{w^k}\right|}\left(d_{w^k}\right).$$

Отметим, что здесь и везде в дальнейшем

$$w^{k+1}\left(= e^k\right) \in OD^{k+1}\left(= \bar{E}^k\right),\ d^{k+1}_{w^{k+1}} = f^k_{e^k},\ k = 1,\dots,m-1.$$

Введем для графа $\Gamma^k$ и множества путей $P^k$ матрицу (Кирхгофа)

$$\Theta^k = \left\| \delta_{e^k p^k} \right\|_{e^k \in E^k, p^k \in P^k}, \ \delta_{e^k p^k} = \begin{cases} 1, & e^k \in p^k \\ 0, & e^k \notin p^k \end{cases}, \ k = 1,\dots,m .$$

Тогда вектор потоков на ребрах $f^k$ на графе $\Gamma^k$ однозначно определяется вектором потоков на путях $x^k = \left\{ x^k_{p^k} \right\}_{p^k \in P^k}$:

$$f^k = \Theta^k x^k, \ k = 1,\dots,m .$$

Обозначим через

$$x = \left\{ x^k \right\}_{k=1}^m, \ f = \left\{ f^k \right\}_{k=1}^m, \ \Theta = \operatorname{diag}\left\{ \Theta^k \right\}_{k=1}^m,$$

$$d^k = \left\{ d^k_{w^k} \right\}_{w^k \in OD^k}, \ X^k\left(d^k\right) = \coprod_{w^k \in OD^k} S_{\left|P^k_{w^k}\right|}\left(d^k_{w^k}\right), \ X = \coprod_{k=1}^m X^k\left(d^k\right),$$

а через

$$\breve{p}^k_{w^k} = \left( p^k_{w^k}, \left\{ p^{k+1}_{w^{k+1}} \right\}_{w^{k+1} \in p^k_{w^k} \cap \bar{E}^k}, \dots \right), \ k = 1,\dots,m$$

полное описание возможного пути (в графе $\Gamma^k$ и графах следующих уровней), соответствующего корреспонденции $w^k \in OD^k$. Множество всех таких путей будем обозначать $\breve{P}^k_{w^k}$. Введем также множество путей $\breve{P}^k = \coprod_{w^k \in OD^k} \breve{P}^k_{w^k}$ и соответствующий вектор распределения потоков по этим путям $x_{\breve{P}^k}$. Определим функции затрат пользователей на пути $\breve{p}^k_{w^k}$ по индукции:

$$G^m_{\breve{p}^m_{w^m}}\left(x_{\breve{P}^m}\right) = \sum_{e^m \in \tilde{E}^m} \delta_{e^m p^m} \tau^m_{e^m}\left(f^m_{e^m}\right),$$

$$G^k_{\breve{p}^k_{w^k}}\left(x_{\breve{P}^k}\right) = \sum_{e^k \in \tilde{E}^k} \delta_{e^k \breve{p}^k_{w^k}} \tau^k_{e^k}\left(f^k_{e^k}\right) + \sum_{w^{k+1} \in \bar{E}^k} G^{k+1}_{\breve{p}^{k+1}_{w^{k+1}}}\left(x_{\breve{P}^{k+1}}\right), \ k = 1,\dots,m-1 .$$

Опишем *марковскую логит динамику* (также говорят гиббсовскую динамику) в повторяющейся игре загрузки графа транспортной сети [33, 302] (см. также подраздел 1.1.4 раздела 1.1 этой главы 1). Пусть имеется $TN$ шагов ($N \gg 1$). Каждый пользователь транспортной сети, использовавший на шаге $t$ путь $\breve{p}^1_{w^1}$, независимо от остальных, на шаге $t+1$ (все введенные новые параметры положительны)

- с вероятностью

$$\frac{\lambda^1}{N} \frac{\exp\left(-G^{t,1}_p / \gamma^1\right)}{\sum_{\tilde{p} \in \breve{P}^1_{w^1}} \exp\left(-G^{t,1}_{\tilde{p}} / \gamma^1\right)} .$$

пытается изменить свой путь $\breve{p}^1_{w^1}$ на $p \in \breve{P}^1_{w^1}$, где $G^{t,1}_p = G^1_p\left(x^t_{\breve{P}^1}\right)$ – затраты на пути $p$ на шаге $t$ ($G^{0,1}_p \equiv 0$);

- равновероятно выбирает $w^2 \in p^1_{w^1} \cap \bar{E}^1$ и затем с вероятностью

$$\frac{\lambda^2\left|p^1_{w^1} \cap \bar{E}^1\right|}{N} \frac{\exp\left(-G^{t,2}_p/\gamma^2\right)}{\sum\limits_{\tilde{p}\in\breve{P}^2_{w^2}} \exp\left(-G^{t,2}_{\tilde{p}}/\gamma^2\right)}$$

пытается изменить в своем пути $\breve{p}^1_{w^1} = \left(p^1_{w^1}, \left\{\breve{p}^2_{w^2}\right\}_{w^2\in p^1_{w^1}\cap\bar{E}^1}\right)$ участок пути $\breve{p}^2_{w^2}$, выбирая путь $p \in \breve{P}^2_{w^2}$, где $G^{t,2}_p = G^2_p\left(x^t_{\breve{P}^2}\right)$ – затраты на пути $p$ на шаге $t$ ($G^{0,2}_p \equiv 0$);

- ... и т.д.;
- с вероятностью

$$1 - \sum_{k=1}^{m} \lambda^k \Big/ N$$

решает не менять тот путь, который использовал на шаге $t$.

Такая динамика отражает ограниченную рациональность агентов (см. замечание 1.4.5 подраздела 1.4.3), и часто используется в теории дискретного выбора [137] и популяционной теории игр [302]. В основном нас будет интересовать поведение такой системы в предположении

$$\lambda^2/\lambda^1 \to \infty,\ \lambda^3/\lambda^2 \to \infty,\ \ldots,\ \lambda^m/\lambda^{m-1} \to \infty,\ N/\lambda^m \to \infty. \qquad (1.4.1)$$

Эта марковская динамика в пределе $N \to \infty$ превращается в марковскую динамику в непрерывном времени [195]. Далее мы, как правило, будем считать, что такой предельный переход был осуществлен.

В пределе $M \to \infty$ эта динамика (концентраций) описывается зацепляющейся системой обыкновенных дифференциальных уравнений (*СОДУ*)

$$\frac{dx_{\breve{p}^1_{w^1}}}{dt} = \lambda^1 \cdot \left( d^1_{w^1} \frac{\exp\left(-G^1_{\breve{p}^1_{w^1}}\left(x_{\breve{P}^1}\right)\Big/\gamma^1\right)}{\sum\limits_{\tilde{p}\in\breve{P}^1_{w^1}} \exp\left(-G^1_{\tilde{p}}\left(x_{\breve{P}^1}\right)\Big/\gamma^1\right)} - x_{\breve{p}^1_{w^1}} \right),\ \breve{p}^1_{w^1} \in \breve{P}^1_{w^1},\ w^1 \in OD^1,$$

$$\frac{dx_{\breve{p}^2_{w^2}}}{dt} = \lambda^2 \cdot \left( d^2_{w^2} \frac{\exp\left(-G^2_{\breve{p}^2_{w^2}}\left(x_{\breve{P}^2}\right)\Big/\gamma^2\right)}{\sum\limits_{\tilde{p}\in\breve{P}^2_{w^2}} \exp\left(-G^2_{\tilde{p}}\left(x_{\breve{P}^2}\right)\Big/\gamma^2\right)} - x_{\breve{p}^2_{w^2}} \right),\ \breve{p}^2_{w^2} \in \breve{P}^2_{w^2},\ w^2 \in OD^2 = \bar{E}^1,$$

……………………………………………………………………

Применяя по индукции (в виду условия (1.4.1)) теорему Тихонова к этой СОДУ [107, 116], можно получить описание аттрактора СОДУ – глобально устойчивой (при $T \to \infty$) неподвижной точки. Для того чтобы это сделать, введем обозначение

$$\sigma_{e^k}^k\left(f_{e^k}^k\right) = \int_0^{f_{e^k}^k} \tau_{e^k}^k(z)\,dz,\ k = 1,\dots,m.$$

Рассмотрим задачу

$$\Psi(x,f) := \Psi^1(x) = \sum_{e^1 \in \tilde{E}^1} \sigma_{e^1}^1\left(f_{e^1}^1\right) + \Psi^2(x) + \gamma^1 \sum_{w^1 \in OD^1} \sum_{p^1 \in P_{w^1}^1} x_{p^1}^1 \ln\left(x_{p^1}^1 / d_{w^1}^1\right) \to \min_{f=\Theta x,\, x \in X}, \quad (1.4.2)$$

$$\Psi^2(x) = \sum_{e^2 \in \tilde{E}^2} \sigma_{e^2}^2\left(f_{e^2}^2\right) + \Psi^3(x) + \gamma^2 \sum_{w^2 \in \bar{E}^1} \sum_{p^2 \in P_{w^2}^2} x_{p^2}^2 \ln\left(x_{p^2}^2 / d_{w^2}^2\right),\ d_{w^2}^2 = f_{w^2}^1,$$

……………………………………………………………………………………..

$$\Psi^k(x) = \sum_{e^k \in \tilde{E}^k} \sigma_{e^k}^k\left(f_{e^k}^k\right) + \Psi^{k+1}(x) + \gamma^k \sum_{w^k \in \bar{E}^{k-1}} \sum_{p^k \in P_{w^k}^k} x_{p^k}^k \ln\left(x_{p^k}^k / d_{w^k}^k\right),\ d_{w^{k+1}}^{k+1} = f_{w^{k+1}}^k,$$

……………………………………………………………………………………..

$$\Psi^m(x) = \sum_{e^m \in E^m} \sigma_{e^m}^m\left(f_{e^m}^m\right) + \gamma^m \sum_{w^m \in \bar{E}^{m-1}} \sum_{p^m \in P_{w^m}^m} x_{p^m}^m \ln\left(x_{p^m}^m / d_{w^m}^m\right),\ d_{w^m}^m = f_{w^m}^{m-1}.$$

Эта задача эквивалентна следующей цепочке зацепляющихся задач выпуклой (многоуровневой [257]) оптимизации

$$\Phi^1\left(d^1\right) = \min_{\substack{f^1=\Theta^1 x^1,\, x^1 \in X^1\left(d^1\right)\\ d_{e^1}^2 = f_{e^1}^1,\, e^1 \in \bar{E}^1}} \left\{ \sum_{e^1 \in \tilde{E}^1} \sigma_{e^1}^1\left(f_{e^1}^1\right) + \Phi^2\left(d^2\right) + \gamma^1 \sum_{w^1 \in OD^1} \sum_{p^1 \in P_{w^1}^1} x_{p^1}^1 \ln\left(x_{p^1}^1 / d_{w^1}^1\right) \right\}, \quad (1.4.3)$$

$$\Phi^2\left(d^2\right) = \min_{\substack{f^2=\Theta^2 x^2,\, x^2 \in X^2\left(d^2\right)\\ d_{e^2}^3 = f_{e^2}^2,\, e^2 \in \bar{E}^2}} \left\{ \sum_{e^2 \in \tilde{E}^2} \sigma_{e^2}^2\left(f_{e^2}^2\right) + \Phi^3\left(d^3\right) + \gamma^2 \sum_{w^2 \in \bar{E}^1} \sum_{p^2 \in P_{w^2}^2} x_{p^2}^2 \ln\left(x_{p^2}^2 / d_{w^2}^2\right) \right\},$$

……………………………………………………………………………………..

$$\Phi^k\left(d^k\right) = \min_{\substack{f^k=\Theta^k x^k,\, x^k \in X^k\left(d^k\right)\\ d_{e^k}^{k+1} = f_{e^k}^k,\, e^k \in \bar{E}^k}} \left\{ \sum_{e^k \in \tilde{E}^k} \sigma_{e^k}^k\left(f_{e^k}^k\right) + \Phi^{k+1}\left(d^{k+1}\right) + \gamma^k \sum_{w^k \in \bar{E}^{k-1}} \sum_{p^k \in P_{w^k}^k} x_{p^k}^k \ln\left(x_{p^k}^k / d_{w^k}^k\right) \right\},$$

……………………………………………………………………………………..

$$\Phi^m\left(d^m\right) = \min_{f^m=\Theta^m x^m,\, x^m \in X^m\left(d^m\right)} \left\{ \sum_{e^m \in E^m} \sigma_{e^m}^m\left(f_{e^m}^m\right) + \gamma^m \sum_{w^m \in \bar{E}^{m-1}} \sum_{p^m \in P_{w^m}^m} x_{p^m}^m \ln\left(x_{p^m}^m / d_{w^m}^m\right) \right\}.$$

То, что эти задачи выпуклые, сразу может быть не очевидно. Чтобы это понять, заметим, что ограничения $x^k \in X^k\left(d^k\right)$ с помощью метода множителей Лагранжа можно убрать, добавив в функционал слагаемые

$$\sum_{w^k \in \bar{E}^{k-1}} \max_{\lambda^k_{w^k}} \left\langle \lambda^k_{w^k}, \sum_{p^k \in P^k_{w^k}} x^k_{p^k} - d^k_{w^k} \right\rangle, \ k = 1, ..., m.$$

Каждое такое слагаемое – есть выпуклая функция по совокупности параметров $x^k$, $d^k$ (см., например, формулу (3.1.8) стр. 96 [68]). Следовательно (см., например, теорему 3.1.2 стр. 92 и формулу (3.1.9) стр. 96 [68]), $\Phi^m\left(d^m\right)$ – выпуклая функция, но тогда и $\Phi^k\left(d^k\right)$ – выпуклая функция, поскольку (по индукции) $\Phi^{k+1}\left(d^{k+1}\right)$ выпуклая функция ($k = 1, ..., m-1$).

**Теорема 1.4.1.** 1. *Задачи (1.4.2) и (1.4.3) являются эквивалентными задачами выпуклой оптимизации, имеющими единственное решение.*

2. *Введенная марковская логит динамика при $N \to \infty$ – эргодическая. Ее финальное распределение (возникающее в пределе $T \to \infty$) совпадает со стационарным. В предположении (1.4.1) стационарное распределение экспоненциально сконцентрировано в окрестности решения задачи (1.4.3) (в пределе $M \to \infty$ стационарное распределение полностью сосредотачивается на решении задачи (3)).*

3. *Введенная марковская логит динамика при пределах $N \to \infty$, $M \to \infty$ описывается СОДУ. В предположении (1.4.1) любая допустимая траектория СОДУ (соответствующая вектору корреспонденций $d^1$) сходится при $T \to \infty$ к решению задачи (1.4.3).*

**Замечание 1.4.1.** Утверждения 1, 2 теоремы 1.4.1 (кроме единственности решения) остаются верными и в предположении, что по части параметров $\gamma^k$ сделаны предельные переходы (от стохастических равновесий к равновесиям Нэша) $\gamma^k \to 0+$ (важно, что эти переходы осуществляются после предельных переходов, указанных в соответствующих пунктах теоремы 1.4.1). К этому же результату (с точки зрения того, к какой задаче оптимизации в итоге сводится поиск равновесий) приводит рассмотрение на соответствующих уровнях вместо логит динамик имитационных логит динамик [51, 302].

**Замечание 1.4.2.** Утверждения теоремы 1.4.1 и замечания 1.4.1 остаются верными, если на части ребер (любого уровня) сделать предельные переходы (важно, что эти переходы осуществляются после предельных переходов, указанных в соответствующих пунктах теоремы) вида (предел стабильной динамики [23, 47], см. также подраздел 1.1.7 раздела 1.1 этой главы 1)

$$\tau_e\left(f_e\right) := \tau_e^{\mu}\left(f_e\right) \xrightarrow[\mu \to 0+]{} \begin{cases} \bar{t}_e, & f_e < \bar{f}_e \\ \left[\bar{t}_e, \infty\right), & f_e = \bar{f}_e \end{cases},$$

$$d\tau_e^{\mu}\left(f_e\right)/df_e \xrightarrow[\mu \to 0+]{} 0, \ 0 \le f_e < \bar{f}_e,$$

с дополнительной оговоркой, что существует такой $x \in X$, что условие $f = \Theta x$ совместно с $\left\{ f_e < \bar{f}_e \right\}_e$. При этом

$$\sigma_e\left(f_e\right) = \lim_{\mu \to 0+} \int_0^{f_e} \tau_e^\mu\left(z\right) dz = \begin{cases} f_e \bar{t}_e, & f_e \le \bar{f}_e \\ \infty, & f_e > \bar{f}_e \end{cases}.$$

Величину $t_e = \lim_{\mu \to 0+} \tau_e^\mu\left(f_e^\mu\right) \ge \bar{t}_e$ можно понимать как затраты на проезд по ребру $e$ (см. также подраздел 1.4.3), а $\lim_{\mu \to 0+} \tau_e^\mu\left(f_e^\mu\right) - \bar{t}_e$ – как дополнительные затраты, приобретенные из-за наличия “пробки” на ребре $e$ [47, 272], возникшей из-за функционирования ребра на пределе пропускной способности $\bar{f}_e$. Эти дополнительные затраты в точности совпадают с множителем Лагранжа к ограничению $f_e \le \bar{f}_e$ [47, 272]. Их также можно понимать как оптимальные платы за проезд (для обычных ребер эти платы равны $f_e \, d\tau_e^\mu\left(f_e\right) / df_e$ [31, 134, 301]), взимаемые согласно механизму Викри–Кларка–Гроуса [134].

Если для некоторых $1 \le p \le q \le m$ имеют место равенства $\gamma^p = \ldots = \gamma^q$, то можно свернуть $\Gamma^p$, …, $\Gamma^q$ в один граф $\coprod_{k=p}^{q} \Gamma^k$. Это следует из свойств энтропии (см. свойство 3 § 4 главы 2 [126]).

Далее мы отдельно рассмотрим специальный случай $\gamma^1 = \ldots = \gamma^m = \gamma$. В этом случае мы имеем граф

$$\Gamma = \coprod_{k=1}^{m} \Gamma^k = \left\langle V, E = \coprod_{k=1}^{m} \tilde{E}^k \right\rangle,$$

который имеет всего один уровень, а задача (1.4.1) может быть переписана следующим образом

$$\Psi\left(x\right) = \sum_{e \in E} \sigma_e\left(f_e\right) + \gamma \sum_{w^1 \in OD^1} \sum_{p \in \tilde{P}^1_{w^1}} x_p \ln\left(x_p / d^1_{w^1}\right) \to \min_{f = \Theta x, \, x \in X}. \qquad (1.4.4)$$

**Теорема 1.4.2.** *При* $\gamma^1 = \ldots = \gamma^m = \gamma$

1. *задачи (1.4.2) и (1.4.4) являются эквивалентными задачами выпуклой оптимизации, имеющими единственное решение;*
2. *введенная марковская логит динамика при* $N \to \infty$ *– эргодическая. Ее финальное распределение (возникающее в пределе* $T \to \infty$*) совпадает со стационарным, которое представимо в виде (представление Санова)*

$$\sim \exp\left(-\frac{M}{\gamma} \cdot \left(\Psi\left(x\right) + o\left(1\right)\right)\right), \; M \gg 1.$$

*Как следствие, получаем, что стационарное распределение экспоненциально сконцентрировано в окрестности решения задачи (1.4.4) (в пределе $M \to \infty$ стационарное распределение полностью сосредотачивается на решении задачи (1.4.4));*

3. *введенная марковская логит динамика при пределах $N \to \infty$, $M \to \infty$ описывается СОДУ. Функция $\Psi(x)$ является функцией Ляпунова этой СОДУ (принцип Больцмана). То есть убывает на траекториях СОДУ. Как следствие, любая допустимая траектория СОДУ (соответствующая вектору корреспонденций $d^1$) сходится при $T \to \infty$ к решению задачи (1.4.4).*

**Замечание 1.4.3.** К теореме 1.4.2 можно сделать замечания аналогичные замечаниям 1.4.1, 1.4.2 к теореме 1.4.1.

### 1.4.3 Двойственная задача

Рассмотрим граф

$$\Gamma = \coprod_{k=1}^{m} \Gamma^k = \left\langle V, E = \coprod_{k=1}^{m} \tilde{E}^k \right\rangle.$$

Обозначим через $t_e = \tau_e(f_e)$ (здесь специально упрощаем обозначения, поскольку в виду предыдущего раздела контекст должен восстанавливаться однозначным образом). Запишем в пространстве $t = \{t_e\}_{e \in E}$ *двойственную задачу* к (1.4.3) [33, 47, 49] (далее мы используем обозначение $\operatorname{dom} \sigma^*$ – область определения сопряженной к $\sigma$ функции)

$$\min_{f,x} \left\{ \Psi(x, f) :\ f = \Theta x,\ x \in X \right\} =$$

$$= -\min_{t \in \operatorname{dom} \sigma^*} \left\{ \gamma^1 \psi^1 \left( t/\gamma^1 \right) + \sum_{e \in E} \sigma_e^* (t_e) \right\}, \qquad (1.4.5)$$

где

$$\sigma_e^*(t_e) = \max_{f_e} \left\{ f_e t_e - \int_0^{f_e} \tau_e(z) dz \right\},$$

$$\frac{d\sigma_e^*(t_e)}{dt_e} = \frac{d}{dt_e} \max_{f_e} \left\{ f_e t_e - \int_0^{f_e} \tau_e(z) dz \right\} = f_e : \ t_e = \tau_e(f_e),\ e \in E;$$

$$g_{p^m}^m(t) = \sum_{e^m \in \tilde{E}^m} \delta_{e^m p^m} t_{e^m} = \sum_{e^m \in E^m} \delta_{e^m p^m} t_{e^m},$$

$$g_{p^k}^k(t) = \sum_{e^k \in \tilde{E}^k} \delta_{e^k p^k} t_{e^k} - \sum_{e^k \in \bar{E}^k} \delta_{e^k p^k} \gamma^{k+1} \psi_{e^k}^{k+1} \left( t/\gamma^{k+1} \right),\ k = 1, \dots, m-1,$$

$$\psi_{w^k}^k(t)=\ln\left(\sum_{p^k\in P_{w^k}^k}\exp\left(-g_{p^k}^k(t)\right)\right),\ k=1,\dots,m,$$

$$\psi^1(t)=\sum_{w^1\in OD^1}d_{w^1}^1\psi_{w^1}^1(t).$$

**Теорема 1.4.3.** *Имеет место сильная двойственность (1.4.5). Решение задачи выпуклой оптимизации (1.4.5)* $t\geq 0$ *существует и единственно. По этому решению однозначно можно восстановить решение исходной задачи (1.4.3) (если какой-то из* $\gamma^k\to 0+$*, то однозначность восстановления* $x$ *может потеряться)*

$$f=\Theta x=-\nabla\psi^1\left(t/\gamma^1\right),$$

$$x_{p^k}^k=d_{w^k}^k\frac{\exp\left(-g_{p^k}^k(t)/\gamma^k\right)}{\sum_{\tilde{p}^k\in P_{w^k}^k}\exp\left(-g_{\tilde{p}^k}^k(t)/\gamma^k\right)},\ p^k\in P_{w^k}^k,\ w^k\in OD^k,\ k=1,\dots,m. \tag{1.4.6}$$

*Верен и обратный результат. Пусть* $f=\Theta x$ *– решение задачи (1.4.3), тогда* $t=\left\{\tau_e\left(f_e\right)\right\}_{e\in E}$ *– единственное решение задачи (1.4.3) (если какой-то из* $\gamma^k\to 0+$*, то решение* $x$ *может быть не единственно, однако это никак не сказывается на возможности однозначного восстановления* $t$*).*

**Замечание 1.4.4.** К теореме 1.4.3 можно сделать замечания аналогичные замечаниям 1.4.1, 1.4.2 к теореме 1.4.1. При этом оговорки, возникающие при $\gamma^k\to 0+$, частично уже были сделаны в формулировке самой теоремы. Дополним их следующим наблюдением. Слагаемое $\gamma^1\psi^1\left(t/\gamma^1\right)$ в двойственной задаче (1.4.5) имеет равномерно ограниченную константу Липшица градиента в 2-норме:

$$L_2\leq\frac{1}{\min\limits_{k=1,\dots,m}\gamma^k}\sum_{w^1\in OD^1}d_{w^1}^1\max_{\breve{p}_{w^1}^1\in\breve{P}_{w^1}^1}\left|\breve{p}_{w^1}^1\right|^2,$$

где $\left|\breve{p}_{w^1}^1\right|$ – число ребер в пути $\breve{p}_{w^1}^1$. Эта гладкость теряется при $\gamma^k\to 0+$:

$$-\lim_{\gamma^k\to 0+}\gamma^k\psi_{w^k}^k\left(t/\gamma^k\right)=\min_{p^k\in P_{w^k}^k}g_{p^k}^k(t)$$

– длина кратчайшего пути в графе $\Gamma^k$, отвечающего корреспонденции $w^k\in OD^k$, ребра $e^k\in\bar{E}^k$, которого взвешены величинами $\gamma^{k+1}\psi_{e^k}^{k+1}\left(t/\gamma^{k+1}\right)$, которые можно понимать как “средние” затраты на $e^k\in\bar{E}^k$ (см. замечание 5). Заметим также, что в пределе (стабильной динамики) $\mu\to 0+$ (см. замечание 2) получаем:

$$\sigma_e^*(t_e)=\lim_{\mu\to 0+}\max_{f_e}\left\{f_et_e-\int_0^{f_e}\tau_e^\mu(z)\,dz\right\}=\begin{cases}\bar{f}_e\cdot\left(t_e-\bar{t}_e\right), & t_e\geq\bar{t}_e\\ \infty, & t_e<\bar{t}_e\end{cases}.$$

При этом $\bar{f}_e - f_e$ в точности совпадает с множителем Лагранжа к ограничению $t_e \ge \bar{t}_e$ [47, 272].

**Замечание 1.4.5.** Формулу (1.4.6) можно получить и из других соображений. Предположим, что каждый пользователь $l$ транспортной сети, использующий корреспонденцию $w^k \in OD^k$ на уровне $k$ (ребро $e^{k-1}\left(=w^k\right) \in \bar{E}^{k-1}$ на уровне $k-1$), выбирает маршрут следования $p^k \in P_{w^k}^k$ на уровне $k$, если

$$p^k = \arg\max_{q^k \in P_{w^k}^k} \left\{ -g_{q^k}^k(t) + \xi_{q^k}^{k,l} \right\},$$

где независимые случайные величины $\xi_{q^k}^{k,l}$, имеют одинаковое двойное экспоненциальное распределение, также называемое распределением Гумбеля [33, 137, 302] (см. также подраздел 1.1.4 раздела 1.1 этой главы 1):

$$P\left(\xi_{q^k}^{k,l} < \zeta\right) = \exp\left\{-e^{-\zeta/\gamma^k - E}\right\}.$$

Отметим также, что если взять $E \approx 0.5772$ – константа Эйлера, то

$$M\left[\xi_{q^k}^{k,l}\right] = 0, \; D\left[\xi_{q^k}^{k,l}\right] = \left(\gamma^k\right)^2 \pi^2 / 6.$$

Распределение Гиббса (логит распределение) (1.4.6) получается в пределе, когда число агентов на каждой корреспонденции $w^k \in OD^k$, $k = 1,\ldots,m$ стремится к бесконечности (случайность исчезает и описание переходит на средние величины). Полезно также в этой связи иметь в виду, что [137, 266]

$$\gamma^k \psi_{w^k}^k\left(t/\gamma^k\right) = M_{\left\{\xi_{p^k}^k\right\}_{p^k \in P_{w^k}^k}} \left[\max_{p^k \in P_{w^k}^k} \left\{-g_{p^k}^k(t) + \xi_{p^k}^k\right\}\right].$$

Таким образом, если каждый пользователь сориентирован на вектор затрат $t$ на ребрах $E$ (одинаковый для всех пользователей) и на каждом уровне (принятия решения) пытается выбрать кратчайший путь исходя из зашумленной информации и исходя из усреднения деталей более высоких уровней (такое усреднение можно обосновывать, если, например, как в подразделе 1.4.2, ввести разный масштаб времени (частот принятия решений) на разных уровнях, а можно просто постулировать, что пользователь так действует, как это принято в моделях типа Nested Logit [137]), то такое поведение пользователей (в пределе, когда их число стремится к бесконечности) приводит к описанию распределения пользователей по путям/ребрам (1.4.6). Равновесная конфигурация характеризуется тем, что вектор $t$ породил согласно формуле (1.4.6) такой вектор $f$, что имеет место соотношение $t = \left\{\tau_e\left(f_e\right)\right\}_{e \in E}$. Поиск такого $t$ (неподвижной точки) приводит к задаче (1.4.5).

**Замечание 1.4.6.** Сопоставить формуле (1.4.4), теореме 1.4.2 и замечанию 1.4.3 (отвечающих случаю $\gamma^1 = ... = \gamma^m = \gamma$) вариант двойственной задачи (1.4.5) чрезвычайно просто (мы здесь опускаем соответствующие выкладки). Собственно, понять формулу (1.4.4) как раз проще не из свойств энтропии (как это было описано в подразделе 1.4.2), а с помощью обратного перехода от двойственного описания (1.4.5). Теорема 1.4.3 и замечание 1.4.5 в случае $\gamma^1 = ... = \gamma^m = \gamma$ наглядно демонстрируют отсутствие какой бы то ни было иерархии, и возможность работать на одном графе с естественной интерпретацией функций затрат на путях $g_p(t)$ (без всяких “средних” оговорок).

Перейдем к конспективному обсуждению численных аспектов решения задачи (1.4.5). Как правило, выгоднее решать именно задачу (1.4.5), а не (1.4.3) [33]. На эту задачу удобно смотреть, как на гладкую (с Липшицевым градиентом) задачу композитной оптимизации [266, 277] с евклидовой прокс-структурой (задаваемой 2-нормой). При этом, даже если по ряду параметров $\gamma^k$ требуется сделать предельный переход $\gamma^k \to 0+$, то, как правило, лучше считать, что численно мы все равно решаем задачу со всеми $\gamma^k > 0$ [33]. Этого можно добиться обратным процессом: энтропийной регуляризацией прямой задачи = сглаживанием двойственной. Некоторые детали того, как именно и в каких случаях полезно сглаживать задачу (1.4.5) описаны в работе [33] (см. также [271]).

Композитный быстрый градиентный метод (и различные его вариации с адаптивным подбором константы Липшица градиента, универсальный метод и др. [266, 277, 301], в частности УМТ из раздела 2.2 главы 2) обладает прямо-двойственной структурой [33, 42, 44, 47, 261, 269]. Это означает, что генерируемые этим методом последовательности $\{t^i\}$ и $\{\tilde{t}^i\}$ обладают следующим свойством

$$\gamma^1\psi^1\left(\tilde{t}^N/\gamma^1\right) + \sum_{e\in E}\sigma_e^*\left(\tilde{t}_e^N\right) -$$

$$- \min_{t\in \operatorname{dom}\sigma^*}\left\{\frac{1}{A_N}\left[\sum_{i=0}^{N} a_i\cdot\left(\gamma^1\psi^1\left(t^i/\gamma^1\right) + \left\langle\nabla\psi^1\left(t^i/\gamma^1\right), t - t^i\right\rangle\right)\right] + \sum_{e\in E}\sigma_e^*\left(t_e\right)\right\} \le \frac{CL_2R_2^2}{A_N}, \qquad (1.4.7)$$

где константа $C \le 10$ зависит от метода,

$$a_N \sim N,\ A_N = \sum_{i=0}^{N} a_i,\ A_N \sim N^2,$$

$$R_2^2 = \max\left\{\tilde{R}_2^2, \hat{R}_2^2\right\},\ \tilde{R}_2^2 = \frac{1}{2}\left\|\bar{t} - t^*\right\|_2^2,\ \hat{R}_2^2 = \frac{1}{2}\sum_{e\in E}\left(\tau_e\left(\bar{f}_e^N\right) - t_e^*\right)^2,$$

$\bar{f}^N$ определяется в теореме 1.4.4, метод стартует с $t^0 = \bar{t}$, $t^*$ – решение задачи (1.4.5).

**Теорема 1.4.4.** *Пусть задача (1.4.5) решается прямо-двойственным методом, генерирующим последовательности* $\{t^i\}$ *и* $\{\tilde{t}^i\}$, *с оценкой скорости сходимости (1.4.7), тогда*

$$0 \le \left\{ \gamma^1 \psi^1\left(\tilde{t}^N / \gamma^1\right) + \sum_{e \in E} \sigma_e^*\left(\tilde{t}_e^N\right) \right\} + \Psi\left(\bar{x}^N, \bar{f}^N\right) \le \frac{C L_2 R_2^2}{A_N},$$

*где*

$$f^i = \Theta x^i = -\nabla \psi^1\left(t^i / \gamma^1\right), \; x^i = \left\{ x_{p^k}^{k,i} \right\}_{p^k \in P_{w^k}^k, w^k \in OD^k}^{k=1,\dots,m},$$

$$x_{p^k}^{k,i} = d_{w^k}^k \frac{\exp\left(-g_{p^k}^k\left(t^i\right)/\gamma^k\right)}{\sum_{\tilde{p}^k \in P_{w^k}^k} \exp\left(-g_{\tilde{p}^k}^k\left(t^i\right)/\gamma^k\right)}, \; p^k \in P_{w^k}^k, \; w^k \in OD^k, \; k = 1,\dots,m,$$

$$\bar{f}^N = \frac{1}{A_N} \sum_{i=0}^N a_i f^i, \; \bar{x}^N = \frac{1}{A_N} \sum_{i=0}^N a_i x^i.$$

**Замечание 1.4.7.** В общем случае описанный выше подход, представляется наиболее предпочтительным. Однако для различных специальных случаев приведенные оценки, по-видимому, можно немного улучшить [9, 42]. Подробнее об этом будет написано в разделе 1.6 этой главы 1.

Приведенная теорема 1.4.4 (доказательство более общего утверждения приведено по ссылке http://arxiv.org/pdf/1606.08988.pdf ) оценивает число необходимых итераций. Но на каждой итерации необходимо считать $\nabla \psi^1\left(t/\gamma^1\right)$, а для ряда методов и $\psi^1\left(t/\gamma^1\right)$ (например, для всех адаптивных методов, настраивающихся на параметры гладкости задачи [266, 274, 277]). Подобно [49, 262, 302] можно показать (с помощью сглаженного варианта метода Форда–Беллмана), что для этого достаточно сделать $O\left(\left|O^1\right|\left|E\right| \max_{\breve{p}^1 \in \breve{P}^1} \left|\breve{p}^1\right|\right)$ арифметических операций. Однако необходимо обговорить один нюанс. Для возможности использовать сглаженный вариант метода Форда–Беллмана [49, 262, 302] необходимо предположить, что любые движения по ребрам графа с учетом их ориентации являются допустимыми, т.е. множество путей, соединяющих заданные две вершины (источник и сток), – это множество всевозможных способов добраться из источника в сток по ребрам имеющегося графа с учетом их ориентации. Сделанная оговорка не сильно обременительная, поскольку нужного свойства всегда можно добиться раздутием исходного графа в несколько раз за счет введения дополнительных вершин и ребер. Детали изложены в магистерском дипломе Е.И. Ершова (МФТИ, 2014) и кандидатской диссертации, подготовленной А.А. Лагуновской (МФТИ, 2016).

В целом, хотелось бы отметить, что прием, связанный с искусственным раздутием исходного графа путем добавления новых вершин, ребер, источников, стоков является весьма полезным для ряда приложений [9, 23, 31, 47]. В частности, достаточно популярным является введение фиктивных (с нулевыми затратами) путей-ребер, которые дают возможность ничего не делать пользователям (не перемещаться [31], не торговать [21] и т.п.), что, в свою очередь, позволяет рассматривать ситуации с нефиксированными корреспонденциями $d^1$ [21, 31]. Также популярным приемом является перенесение затрат на преодоление вершин (узлов) графа (перекрестков [31], сортировочных станций [23]) в затраты на прохождение дополнительных ребер, появившиеся при "распутывании" узлов. Но, пожалуй, наиболее важным для большинства приложений является введение фиктивного общего источника и общего стока, соединенных дополнительными ребрами с уже имеющимися вершинами графа [9, 23, 31].

## 1.5 Численные методы поиска равновесного распределения потоков в модели Бэкмана и модели стабильной динамики

### 1.5.1 Введение

В недавних работах [21, 47] (см. также подразделы 1.1.4, 1.1.6, 1.1.8 раздела 1.1 этой главы 1), посвященных сведению поиска равновесного распределения транспортных потоков на сетях к решению задач выпуклой оптимизации, было поставлено несколько таких задач со специальной “сетевой” структурой. Это означает, что, скажем, расчет градиента (стохастического градиента) функционала сводится к поиску кратчайших путей в графе транспортной сети. Эта специфика задач с одной стороны говорит о том, что в реальных приложениях размерность задачи может быть колоссально большой. Это связано с тем, что число реально используемых путей даже в планарном графе, как правило, пропорционально кубу числа вершин, а число вершин в реальных приложениях обычно не меньше тысячи, отметим, что в худших случаях число путей может расти экспоненциально с ростом числа вершин. С другой стороны, такого рода задачи имеют хорошую геометрическую интерпретацию, что позволяет эффективно снижать их размерность.

В частности, в подразделе 1.5.2 этого раздела 1.5 мы описываем метод Франк–Вульфа [142, 203, 243] поиска равновесного распределения потоков в модели Бэкмана, который на каждой итерации требует решения задачи минимизации линейной функции от потоков на ребрах в сети (число ребер порядка нескольких тысяч) на прямом произведении симплексов (симплексов столько, сколько корреспонденций). Для реальных транспортных сетей (тысяча вершин) получается задача минимизации линейной функции в пространстве размерности миллиард, поскольку она зависит от распределения потоков по путям, число которых порядка миллиарда. Ясно, что если смотреть на эту задачу формально с точки зрения оптимизации, то все сводится к полному перебору миллиарда вершин всех симплексов (причем проработка одной вершины – это расчет соответствующего скалярного произведения, то есть порядка нескольких тысяч умножений). К счастью, транспортная специфика задачи позволяет с помощью алгоритма Дейкстры и более современных подходов [145] (в том числе учитывающих “планарность” сети: A*, ALT, SHARC, Reach based routing, Highway hierarchies, Contraction hierarchies и т.п.) решать описанную задачу делая не более десятка миллионов операций (типа умножения двух чисел с плавающей запятой), что намного быстрее. Такого рода конструкции возникают не только в связи с сетевой спецификой задачи [173], но именно для ситуаций, когда в задаче имеется сетевая структура, возможность такой редукции наиболее естественна и типична.

В подразделе 1.5.3 мы предлагаем другой способ поиска равновесия в модели Бэкмана и аналогичных моделях (в модели стабильной динамики, в промежуточных моде-лях). Для этого мы переходим (следуя Ю.Е. Нестерову) к двойственной задаче, в которой целевой функционал оказывается зависящим только от потоков по ребрам, а не от распределения потоков по путям. Таким образом, задача сводится к поиску равновесного распределения потоков по ребрам. При этом в ходе вычислений потоков на ребрах, мы попутно (без дополнительных затрат) вычисляем порождающие их потоки по путям. Также в виду транспортно-сетевой специфики появляется возможность содержательной интерпретации [47] (подобно интерпретации Л.В. Канторовичем цен в экономике [71]), возникающих двойственных множителей, которые в ряде приложений представляют независимый самостоятельный интерес (например, в задаче о тарифной политике грузоперевозок РЖД [21] двойственные множители – тарифы, которые и надо рассчитывать). Также сетевая структура задачи дает возможность не рассчитывать градиент целевой функции на каждой итерации заново, а пересчитывать его, используя градиент, полученный на предыдущей итерации. Грубо говоря, найдя кратчайшие пути, и посчитав на их основе градиент, мы сделаем шаг по антиградиенту, немного изменив веса ребер. Ясно, что большая часть кратчайших путей при этом останется прежними, и можно специально организовать их пересчет, чтобы ускорить вычисления. Похожая философия используется в покомпонентных спусках и в современных подходах к задачам huge-scale оптимизации [264, 273]. Однако сетевая структура задачи требует переосмысления этой техники, рассчитанной изначально в основном на свойства разреженности матриц, возникающих в условии задачи.

Применимость современных вариантов рандомизированных покомпонентных спусков к поиску равновесия в модели стабильной динамики, записанной в новой специальной форме, предложенной недавно Ю.В. Дорном и Ю.Е. Нестеровым, изучается в подразделе 1.5.4.

В заключительном подразделе 1.5.5 приводится краткий сравнительный анализ описанных в разделе методов.

Настоящий раздел представляет собой одну из первых попыток сочетать современные эффективные численные методы выпуклой оптимизации с сетевой структурой задачи, на примере задач, пришедших из поиска равновесного распределения потоков в транспортных сетях и сетях грузовых перевозок РЖД [23, 47].

### 1.5.2 Метод Франк–Вульфа поиска равновесия в модели Бэкмана

Для удобства чтения напомним здесь наиболее популярную на протяжении более чем полувека модель равновесного распределения потоков Бэкмана [47, 49, 112, 147, 289, 310] ранее описанную нами в подразделе 1.1.4 раздела 1.1 этой главы 1. Таким образом, первая половина этого раздела во многом будет повторять текст подразделе 1.1.4 раздела 1.1 и текст раздела 4 работы [47].

Пусть транспортная сеть города представлена ориентированным графом $\Gamma=(V,E)$, где $V$ – узлы сети (вершины), $E\subset V\times V$ – дуги сети (рёбра графа), $O\subseteq V$ – источники корреспонденций ($S=|O|$), $D\subseteq V$ – стоки. В современных моделях равновесного распределения потоков в крупном мегаполисе число узлов графа транспортной сети обычно выбирают порядка $n=|V|\sim 10^3-10^4$. Число ребер $|E|$ получается в три четыре раза больше. Пусть $W\subseteq\{w=(i,j):i\in O, j\in D\}$ – множество корреспонденций, т.е. возможных пар «источник» – «сток»; $p=\{v_1,v_2,...,v_m\}$ – путь из $v_1$ в $v_m$, если $(v_k,v_{k+1})\in E$, $k=1,...,m-1$, $m>1$; $P_w$ – множество путей, отвечающих корреспонденции $w\in W$, то есть если $w=(i,j)$, то $P_w$ – множество путей, начинающихся в вершине $i$ и заканчивающихся в $j$; $P=\bigcup_{w\in W}P_w$ – совокупность всех путей в сети $\Gamma$ (число “разумных” маршрутов $|P|$, которые потенциально могут использоваться, обычно растет с ростом числа узлов сети не быстрее чем $\mathrm{O}(n^3)$ [112, 289, 310]); $x_p$ [автомобилей/час] – величина потока по пути $p$, $x=\{x_p:\ p\in P\}$; $f_e$ [автомобилей/час] – величина потока по дуге $e$:

$$f_e(x)=\sum_{p\in P}\delta_{ep}x_p\text{, где }\delta_{ep}=\begin{cases}1, & e\in p\\ 0, & e\notin p\end{cases};$$

$\tau_e(f_e)$ – удельные затраты на проезд по дуге $e$. Как правило, предполагают, что это – (строго) возрастающие, гладкие функции от $f_e$. Точнее говоря, под $\tau_e(f_e)$ правильнее понимать представление пользователей транспортной сети об оценке собственных затрат (обычно временных в случае личного транспорта и комфортности пути (с учетом времени в пути) в случае общественного транспорта) при прохождении дуги $e$, если поток желающих проехать по этой дуге будет $f_e$.

Рассмотрим теперь $G_p(x)$ – затраты временные или финансовые на проезд по пути $p$. Естественно считать, что $G_p(x)=\sum_{e\in E}\tau_e(f_e(x))\delta_{ep}$.

Пусть также известно, сколько перемещений в единицу времени $d_w$ осуществляется согласно корреспонденции $w\in W$. Тогда вектор $x$, характеризующий распределение потоков, должен лежать в допустимом множестве:

$$X=\left\{x\geq 0:\ \sum_{p\in P_w}x_p=d_w,\ w\in W\right\}.$$

Рассмотрим игру, в которой каждой корреспонденции $w \in W$ соответствует свой, достаточно большой ($d_w \gg 1$), набор однотипных "игроков", осуществляющих передвижение согласно корреспонденции $w$. Чистыми стратегиями игрока служат пути, а выигрышем – величина $-G_p(x)$. Игрок "выбирает" путь следования $p \in P_w$, при этом, делая выбор, он пренебрегает тем, что от его выбора также "немного" зависят $|P_w|$ компонент вектора $x$ и, следовательно, сам выигрыш $-G_p(x)$. Можно показать (см., например, [47]), что отыскание равновесия Нэша–Вардропа $x^* \in X$ (макро описание равновесия) равносильно решению задачи нелинейной комплементарности (принцип Вардропа):

$$\textit{для любых } w \in W,\ p \in P_w \textit{ выполняется } x_p^* \cdot \left( G_p\left(x^*\right) - \min_{q \in P_w} G_q\left(x^*\right)\right) = 0.$$

Действительно допустим, что реализовалось какое-то другое равновесие $\tilde{x}^* \in X$, которое не удовлетворяет этому условию. Покажем, что тогда найдется водитель, которому выгодно поменять свой маршрут следования. Действительно, тогда

$$\text{существуют такие } \tilde{w} \in W,\ \tilde{p} \in P_{\tilde{w}}, \text{ что } \tilde{x}_{\tilde{p}}^* \cdot \left( G_{\tilde{p}}\left(\tilde{x}^*\right) - \min_{q \in P_{\tilde{w}}} G_q\left(\tilde{x}^*\right)\right) > 0.$$

Каждый водитель (множество таких водителей не пусто, так как $\tilde{x}_{\tilde{p}}^* > 0$), принадлежащий корреспонденции $\tilde{w} \in W$, и использующий путь $\tilde{p} \in P_{\tilde{w}}$, действует не разумно, поскольку существует такой путь $\tilde{q} \in P_{\tilde{w}}$, $\tilde{q} \neq \tilde{p}$, что $G_{\tilde{q}}\left(\tilde{x}^*\right) = \min_{q \in P_{\tilde{w}}} G_q\left(\tilde{x}^*\right)$. Этот путь $\tilde{q}$ более выгоден, чем $\tilde{p}$. Аналогично показывается, что при $x^* \in X$ никому из водителей уже не выгодно отклоняться от своих стратегий.

Условие равновесия может быть переписано следующим образом [49, 112, 147, 289, 310]

$$\text{для всех } x \in X \text{ выполняется } \left\langle G\left(x^*\right), x - x^*\right\rangle \geq 0.$$

Рассматриваемая нами игра принадлежит к классу, так называемых, потенциальных игр [254, 302], поскольку $\partial G_p(x)/\partial x_q = \partial G_q(x)/\partial x_p$. Существует такая функция

$$\Psi\left(f(x)\right) = \sum_{e \in E} \int_0^{f_e(x)} \tau_e(z)\,dz = \sum_{e \in E} \sigma_e\left(f_e(x)\right),$$

где $\sigma_e\left(f_e\right) = \int_0^{f_e(x)} \tau_e(z)\,dz$, что $\partial \Psi(x)/\partial x_p = G_p(x)$ для любого $p \in P$. Таким образом, $x^* \in X$ – равновесие Нэша–Вардропа в этой игре тогда и только тогда, когда оно доставляет минимум $\Psi\left(f(x)\right)$ на множестве $X$.

**Теорема 1.5.1 [47, 49, 112, 289, 310].** *Вектор $x^*$ будет равновесием Нэша–Вардропа тогда и только тогда, когда*

$$x \in \operatorname{Arg}\min_x \left[ \Psi\big(f(x)\big) = \sum_{e\in E} \sigma_e\big(f_e(x)\big):\ f = \Theta x,\ x \in X \right].$$

*Если преобразование $G(\cdot)$ строго монотонное, то равновесие $x$ единственно. Если $\tau_e'(\cdot) > 0$, то равновесный вектор распределения потоков по ребрам $f$ – единственный* (*это еще не гарантирует единственность вектора распределения потоков по путям $x$ [49]*).

Итак, будем решать задачу ( $\Psi_*$ – оптимальное значение функционала)

$$\Psi(f) = \sum_{e\in E} \sigma_e(f_e) \to \min_{\substack{f=\Theta x \\ x\in X}}$$

методом условного градиента [165, 214, 220, 263] (Франк–Вульфа). Этот метод нам также будет встречаться несколько раз в 4 главе.

**Начальная итерация**

*Положим $\tilde{t}_e^0 = \partial\Psi(0)/\partial f_e = \tau_e(0)$ и рассмотрим задачу*

$$\sum_{e\in E} \tilde{t}_e^0 f_e \to \min_{\substack{f=\Theta x \\ x\in X}} \quad .$$

*Эту задачу можно переписать, как*

$$\min_{x\in X} \sum_{e\in E} \tilde{t}_e^0 \sum_{p\in P} \delta_{ep} x_p = \sum_{w\in W} d_w \min_{p\in P_w} \left\{ \sum_{e\in E} \delta_{ep} \tilde{t}_e^0 \right\} = \sum_{w\in W} d_w T_w\big(\tilde{t}^0\big),$$

*где $T_w\big(\tilde{t}^0\big)$ – длина кратчайшего пути из $i$ в $j$ (где $w = (i, j)$) на графе, ребра которого взвешены вектором $\tilde{t}^0 = \left\{\tilde{t}_e^0\right\}_{e\in E}$. Таким образом, выписанную задачу можно решить с учетом того, что $n = |V| \sim |E|$, за $\tilde{O}(Sn)$ (здесь и далее $\tilde{O}(\ ) = O(\ )$ с точностью до логарифмического фактора) и быстрее современными вариациями алгоритма Дейкстры [127, 145, 150, 200]. Обозначим решение этой задачи через $f^0$.*

Можно интерпретировать ситуацию таким образом, что в начальный момент водители посчитали, что все дороги абсолютно свободны и выбрали согласно этому предположению кратчайшие пути, соответствующие их целям, и поехали по этим маршрутам (путям). На практике более равномерное распределение водителей по путям в начальный момент может оказаться более предпочтительным.

Поняв, что в действительности из-за наличия других водителей время в пути не соответствует первоначальной оценке, выраженной весами ребер $\tilde{t}_e^0$, доля $\gamma^k$ водителей (обнаруживших это и готовых что-то менять) на следующем $(k+1)$-м шаге изменят свой выбор исходя из кратчайших путей, посчитанных по распределению водителей на предыдущем $k$-м шаге. Таким образом, возникает процедура “нащупывания” равновесия. Если выбирать специальным образом $\gamma^k$ (в частности, необходимо $\gamma^k \xrightarrow[k\to\infty]{} 0$, чтобы избежать колебания вокруг равновесия (minority game [169]), и $\sum_{k=0}^{\infty}\gamma^k = \infty$, чтобы до равновесия дойти), то система, действительно, сойдется в равновесие. Опишем теперь более формально сказанное.

**Итерации** $k = 0,1,2,...$

*Пусть* $f^k$ *– вектор потоков на ребрах, полученный на предыдущей итерации с номером* $k$. *Положим* $\tilde{t}_e^k = \partial\Psi\left(f^k\right)/\partial f_e = \tau_e\left(f^k\right)$ *и рассмотрим задачу*

$$\sum_{e\in E}\tilde{t}_e^k y_e \to \min_{\substack{y=\Theta x\\ x\in X}} .$$

*Так же, как и раньше задача сводится к поиску кратчайших путей на графе, ребра которого взвешены вектором* $\tilde{t}^k = \left\{\tilde{t}_e^k\right\}_{e\in E}$.

*Обозначим решение задачи через* $y^k$. *Положим*

$$f^{k+1} = \left(1-\gamma^k\right)f^k + \gamma^k y^k,\ \gamma^k = \frac{2}{k+1}.$$

Заметим, что возникающую здесь задачу поиска кратчайших путей на графе можно попробовать решать быстрее, чем за $\tilde{O}(Sn)$. Связано это с тем, что мы уже решали на предыдущей итерации аналогичную задачу для этого же графа с близкими весами ребер [127, 145, 150, 200] (веса ребер графа с ростом $k$ меняются все слабее от шага к шагу, поскольку $\gamma^k \xrightarrow[k\to\infty]{} 0$). Тем не менее, далее в разделе мы будем считать, что одна итерация этого метода занимает $\tilde{O}(Sn)$.

Заметим также, что решая задачи поиска кратчайших путей мы находим (одновременно, т.е. без дополнительных затрат) не только вектор распределения потоков по ребрам $y$, но и разреженный вектор распределения потоков по путям $x$.

Строго говоря, нужно найти вектор $y^k$, а не кратчайшие пути. Чтобы получить вектор $y^k$ за $\tilde{O}(Sn)$ стоит для каждого из $S$ источников построить (например, алгоритмом Дейкстры) соответствующее дерево кратчайших путей (исходя из принципа динамического программирования “часть кратчайшего пути сама будет кратчайшим путем” несложно понять, что получится именно дерево, с корнем в рассматриваемом источнике). Это мож-но сделать для одного источника за $\tilde{O}(n)$. Однако, главное, правильно взвешивать ребра (их не больше $n$) такого дерева, чтобы за один проход этого дерева можно было восстановить вклад (по всем ребрам) соответствующего источника в общий вектор $y^k$. Ребро должно иметь вес равный сумме всех проходящих через него корреспонденций с заданным источником (корнем дерева). Имея значения соответствующих корреспонденций (их также не больше $n$) за один обратный проход (то есть с листьев к корню) такого дерева можно осуществить необходимое взвешивание (с затратами не более $O(n)$). Делается это по правилу: вес ребра равен сумме корреспонденции (возможно, равной нулю), в соответствующую вершину, в которую ребро входит и сумме весов всех ребер (если таковые имеются), выходящих из упомянутой вершины.

**Теорема 1.5.2 [165, 214, 220, 263].** *Имеет место следующая оценка*

$$\Psi\left(f^N\right)-\Psi_*\le\Psi\left(f^N\right)-\Psi_N\le\frac{2L_pR_p^2}{N+1},\; f^N\in\Delta=\left\{f=\Theta x:\; x\in X\right\},$$

*где*

$$\Psi_N=\max_{k=0,\dots,N}\left\{\Psi\left(f^k\right)+\left\langle\nabla\Psi\left(f^k\right),y^k-f^k\right\rangle\right\},$$

$$R_p^2=\max_{k=0,\dots,N}\left\|y^k-f^k\right\|_p^2\le\max_{f,\breve{f}\in\Delta}\left\|\breve{f}-f\right\|_p^2,\; L_p=\max_{\|h\|_p\le 1}\max_{f\in\operatorname{conv}\left(f^0,f^1,\dots,f^N\right)}\left\langle h,\operatorname{diag}\left\{\tau_e'\left(f_e\right)\right\}h\right\rangle,\; 1\le p\le\infty.$$

**Замечание 1.5.1.** Из доказательства этой теоремы [165, 214, 220, 263] можно усмотреть немного более тонкий способ оценки $L_p$, в котором вместо $f\in\operatorname{conv}\left(f^0,f^1,\dots,f^N\right)$ можно брать $f\in\operatorname{conv}\left(f^0,f^1\right)\cup\operatorname{conv}\left(f^1,f^2\right)\cup\dots\cup\operatorname{conv}\left(f^{N-1},f^N\right)$. Однако для небольшого упрощения выкладок мы будем использовать приведенный в формулировке теоремы огрубленный вариант.

К сожалению, в приложениях (см., пример о расщеплении потоков на личный и общественный транспорт в подразделе 1.5.3) функции $\tau_e\left(f_e\right)$ могут иметь вертикальные асимптоты, что не позволяет равномерно по $N$ ограничить $L_p$ (даже если более тонко оценивать $L_p$, см. замечание 1.5.1). Такие случаи мы просто исключаем из рассмотрения для метода, описанного в этом разделе. Другими словами, мы считаем, что функции $\tau_e\left(f_e\right)$ заданы на положительной полуоси. К таким функциям относятся, например, BPR-функции (см. подраздел 1.5.3).

Обратим внимание на то, что сам метод никак не зависит от выбора параметра $p$, от того какие получаются $R_p^2$ и $L_p$, в то время как оценка на число итераций, которые необходимо сделать для достижения заданной по функции (функционалу) точности, от этого выбора зависит. Как следствие, от этого выбора зависит и критерий останова (значение $\Psi_*$ нам априорно не известно.

Будем считать $p = 2$ (сопоставимые оценки, получаются и при выборе $p = \infty$):

$$L_2\left(f^0, f^1, ..., f^N\right) = \max_{f \in \mathrm{conv}\left(f^0, f^1, ..., f^N\right)} \max_{e \in E} \tau'_e\left(f_e\right) = \max_{e \in E} \tau'_e\left(\max_{k=0,...,N} f_e^k\right),\; R_2^2 = \max_{f, \breve{f} \in \Delta} \left\| \breve{f} - f \right\|_2^2.$$

Величину $R_2^2$ мы можем оценить априорно (при это, к сожалению, получается довольно грубая оценка), т.е. можно считать её нам известной. Труднее обстоит дело с $L_2$. Далее предлагается оригинальный способ запуска метода Франк–Вульфа, критерий останова которого не требует априорного знания $L_2$.

Задаемся точностью $\varepsilon > 0$. Оцениваем $R_2^2$. Полагаем $L_2 = 1$ (для определенности). Запускаем метод Франк–Вульфа с $N\left(L_2\right) = 2L_2R_2^2/\varepsilon$. На каждом шаге проверяем условие (это делается за $\mathrm{O}\left(n\right)$)

$$L_2\left(f^0, f^1, ..., f^k\right) = \max_{e \in E} \tau'_e\left(\max_{l=0,...,k} f_e^l\right) \le L_2.$$

Если на всех шагах условие выполняется, то сделав $N\left(L_2\right)$ шагов, гарантированно получим решение с нужной точностью. Если же на каком-то шаге $k < N\left(L_2\right)$ условие нарушилось, т.е. $L_2\left(f^0, f^1, ..., f^k\right) > L_2$, то полагаем $L_2 := L_2\left(f^0, f^1, ..., f^k\right)$, пересчитываем $N\left(L_2\right)$ и переходим к следующему шагу. Таким образом, по ходу итерационного процесса мы корректируем критерий останова, оценивая необходимое число шагов по получаемой последовательности $\left\{f^k\right\}$. Специфика данной постановки, которая позволила так рассуждать, заключается в наличии явного представления

$$L_2\left(f^0, f^1, ..., f^k\right) = \max_{e \in E} \tau'_e\underbrace{\left(\max_{l=0,...,k} f_e^l\right)}_{\breve{f}_e},$$

и независимости используемого метода от выбора $L_2$ (шаг метода Франк–Вульфа $\gamma^k = 2\left(k+1\right)^{-1}$ вообще ни от каких параметров не зависит).

На практике, однако, приведенный способ работает не очень хорошо из-за использования завышенных оценок для $L_2$ и $R_2^2$. Более эффективным оказался другой способ, который использует неравенство (см. теорему 2) $\Psi\left(f^N\right) - \Psi_N \le \varepsilon$. В этом способе важно, что $\Psi_k$, $k = 0, ..., N$ – автоматически рассчитываются на каждой итерации без дополнительных затрат, а $\Psi\left(f^k\right)$ может быть рассчитано на каждой итерации по известному $f^k$ за $\tilde{O}(n)$. Однако нет необходимости проверять этот критерий на каждой итерации, можно это делать, например, с периодом $\alpha\tilde{\varepsilon}^{-1}$, где $\tilde{\varepsilon}$ – относительная точность по функции (скажем, $\tilde{\varepsilon}=0.01$ – означает, что $\varepsilon = 0.01\Psi\left(f^0\right)$), $\alpha \approx 1$ подбирается эвристически, исходя из задачи. Следуя [220], можно еще немного упростить рассуждения за счет небольшого увеличения числа итераций. А именно, можно использовать оценки (при этом следует полагать $\gamma^k = 2(k+2)^{-1}$)

$$\Psi\left(f^k\right)-\Psi_* \le \left\langle \nabla\Psi\left(f^k\right), f^k - y^k \right\rangle,\ \min_{k=1,\dots,N} \left\langle \nabla\Psi\left(f^k\right), f^k - y^k \right\rangle \le \frac{7L_2R_2^2}{N+2}.$$

Таким образом, в данном разделе был описан способ поиска равновесного распределения потоков по ребрам $f$, который за время

$$\tilde{O}\left(SnL_2R_2^2/\varepsilon\right)$$

находит такой $f^{N(\varepsilon)}$, что

$$\Psi\left(f^{N(\varepsilon)}\right)-\Psi_* \le \varepsilon .$$

Численные эксперименты с описанным методом были произведены А.С. Аникиным в 2015 году (использовались данные https://github.com/bstabler/TransportationNetworks ). Эксперименты показали высокую эффективность метода. Заметим, что именно этот метод используется в большинстве коммерческих продуктах транспортного моделирования для поиска равновесного распределения потоков по путям.

### 1.5.3 Рандомизированный метод двойственных усреднений поиска равновесия в модели стабильной динамики (Нестерова–де Пальма)

В ряде постановок задач вместо функций затрат на ребрах $\tau_e\left(f_e\right)$ заданы ограничения на пропускные способности $f_e \le \bar{f}_e$ и затраты на прохождения свободного (не загруженного $f_e < \bar{f}_e$) ребра $\bar{t}_e$. В модели стабильной динамики это сделано для всех ребер [47, 275], а в модели грузоперевозок РЖД – только для части [23]. Согласно работе [47], такую новую модель можно получить предельным переходом из модели Бэкмана, с помощью введения внутренних штрафов в саму модель. А именно, будем считать, что (как и в модели Бэкмана) у всех ребер есть свои функции затрат $\tau_e^\mu\left(f_e\right)$, но для части ребер $e \in E'$ (какой именно части, зависит от задачи) осуществляется предельный переход

$$\tau_e^\mu\left(f_e\right) \xrightarrow[\mu\to 0+]{} \begin{cases} \bar{t}_e, & 0 \le f_e < \bar{f}_e \\ \left[\bar{t}_e, \infty\right), & f_e = \bar{f}_e \end{cases},$$

$$d\tau_e^{\mu}\left(f_e\right)/df_e \xrightarrow[\mu\to 0+]{} 0,\ 0\le f_e<\bar{f}_e.$$

Обозначив через $x(\mu)$ – равновесное распределение потоков по путям в модели Бэкмана при функциях затрат на ребрах $\tau_e^{\mu}\left(f_e\right)$, получим, что при $e\in E'$

$$\tau_e^{\mu}\left(f_e\left(x(\mu)\right)\right)\xrightarrow[\mu\to 0+]{} t_e,\ \ f_e\left(x(\mu)\right)\xrightarrow[\mu\to 0+]{} f_e,$$

где пара $(t, f)$ – равновесие в модели стабильной динамики и ее вариациях [23, 47, 275] с тем же графом и матрицей корреспонденций, что и в модели Бэкмана, и с ребрами $e\in E'$, характеризующимися набором $\left(\bar{t}, \bar{f}\right)$ из определения функций $\tau_e^{\mu}\left(f_e\right)$. Заметим, что если $t_e>\bar{t}_e$, то $t_e-\bar{t}_e$ можно интерпретировать, например, как время, потерянное в пробке на этом ребре [47, 275].

Согласно подразделу 1.5.2 равновесная конфигурация при таком переходе $\mu\to 0+$ должна находиться из решения задачи

$$\Psi(f)=\sum_{e\in E\setminus E'}\int_0^{f_e}\tau_e^{\mu}(z)\,dz+\lim_{\mu\to 0+}\sum_{e\in E'}\int_0^{f_e}\tau_e^{\mu}(z)\,dz\to\min_{f=\Theta x,\,x\in X}.$$

Считая, что в равновесии не может быть $\lim\limits_{\mu\to 0+}\tau_e^{\mu}\left(f_e\right)=\infty$ (иначе, равновесие просто не достижимо, и со временем весь граф превратится в одну большую пробку), можно не учитывать в интеграле вклад точек $\bar{f}_e$ (в случае попадания в промежуток интегрирования), то есть переписать задачу следующим образом

$$\min_{f=\Theta x,\,x\in X}\left\{\sum_{e\in E\setminus E'}\int_0^{f_e}\tau_e^{\mu}(z)\,dz+\sum_{e\in E'}\int_0^{f_e}\left(\bar{t}_e+\delta_{\bar{f}_e}(z)\right)dz\right\}\Leftrightarrow\min_{\substack{f=\Theta x,\,x\in X\\ f_e\le\bar{f}_e,\,e\in E'}}\left\{\sum_{e\in E\setminus E'}\int_0^{f_e}\tau_e^{\mu}(z)\,dz+\sum_{e\in E'}f_e\bar{t}_e\right\},$$

где

$$\delta_{\bar{f}_e}(z)=\begin{cases}0, & 0\le z<\bar{f}_e\\ \infty, & z\ge\bar{f}_e\end{cases},\ e\in E'.$$

**Теорема 1.5.3 [47, 275].** *Двойственная задача к выписанной выше задаче может быть приведена к следующему виду:*

$$\Upsilon(t)=-\sum_{w\in W}d_wT_w(t)+\left\langle\bar{f},t-\bar{t}\right\rangle-\mu\sum_{e\in E\setminus E'}h_e^{\mu}\left(t_e\right)\to\min_{\substack{t_e\ge\bar{t}_e,\,e\in E'\\ t_e\in\operatorname{dom}h_e^{\mu}(t_e),\,e\in E\setminus E'}},\qquad(1.5.1)$$

*где $T_w(t)$ – длина кратчайшего пути из $i$ в $j$ ($w=(i,j)\in W$) на графе, ребра которого взвешены вектором $t=\left\{t_e\right\}_{e\in E}$, а функции $h_e^{\mu}\left(t_e\right)$ – гладкие и вогнутые.*

*При этом решение изначальной задачи $f$ можно получить из формул:*

$f_e=\bar{f}_e-s_e,\ e\in E'$, *где* $s_e\ge 0$ *– множитель Лагранжа к ограничению* $t_e\ge\bar{t}_e$;

$$\tau_e^{\mu}\left(f_e\right)=t_e,\ e\in E\setminus E'.$$

Приведем пример модели (типа стабильной динамики) расщепления пользователей на личный и общественный транспорт [47] (см. также подраздел 1.1.8 раздела 1.1 этой главы 1), в которой каждое ребро $e\in E$ изначального графа продублировано для личного ("*л*") и общественного ("*о*") транспорта, при этом для общественного транспорта [275]

$$\tau_e\left(f_e^o\right)=\bar{t}_e^o\cdot\left(1+\mu\frac{\bar{f}_e^o}{\bar{f}_e^o-f_e^o}\right),$$

а для личного транспорта был осуществлен предельный переход $\mu\to 0+$ в аналогичных формулах

$$\tau_e\left(f_e^л\right)=\bar{t}_e^л\cdot\left(1+\mu\frac{\bar{f}_e^л}{\bar{f}_e^л-f_e^л}\right).$$

Поиск равновесного расщепления на [личный] и [общественный] транспорт приводит к следующей задаче [47]:

$$-\sum_{w\in W}d_w\min\left\{T_w^л\left(t^л\right),T_w^o\left(t^o\right)\right\}+\left\langle\bar{f}^л,t^л-\bar{t}^л\right\rangle+\left\langle\bar{f}^o,t^o-\bar{t}^o\right\rangle-\mu\sum_{e\in E}\bar{f}_e^o\cdot\bar{t}_e^o\cdot\ln\left(1+\frac{t_e^o-\bar{t}_e^o}{\bar{t}_e^o\mu}\right)\to\min_{\substack{t^л\ge\bar{t}^л\\ t^o\ge\bar{t}^o\cdot(1-\mu)}},$$

при этом $f^л=\bar{f}^л-s^л$, где $s^л$ – вектор множителей Лагранжа для ограничений $t^л\ge\bar{t}^л$,

$$f_e^o=\bar{f}_e^o\cdot\left(1-\frac{\bar{t}_e^o\cdot\mu}{t_e^o-(1-\mu)\bar{t}_e^o}\right).$$

Для упрощения рассуждений далее будем считать, что $E'=E$, т.е. на всех ребрах перешли к пределу $\mu\to 0+$. Если это не так, то надо будет далее использовать не метод зеркального спуска [165] (см. формулу (1.5.2) ниже), а его композитный вариант [266] с сепарабельным композитом $-\mu\sum_{e\in E\setminus E'}h_e^{\mu}\left(t_e\right)$. Возникающая на каждой итерации задача (поиска градиентного отображения, см. формулу (1.5.2)), в виду сепарабельности ограничений, распадается в $|E|=O(n)$ одномерных задач выпуклой оптимизации, которые можно решить за $\tilde{O}(n)$ с машинной точностью. Интересно заметить, что в случае, когда

$$\tau_e^{\mu}\left(f_e\right)=\bar{t}_e\cdot\left(1+\gamma\cdot\left(\frac{f_e}{\bar{f}_e}\right)^{\frac{1}{\mu}}\right)$$

– BPR-функции с $\mu=0.25$ (наиболее часто встречающаяся на практике ситуация [254]), то имеются явные формулы, поэтому итерацию можно сделать быстрее – за $O(n)$.

Численно решать задачу негладкой выпуклой оптимизации (1.5.1) с $E'=E$ (общий случай изложен в разделе 3.3 главы 3) будем с помощью специальным образом рандоми-зированных вариантов метода зеркального спуска [54, 165, 225, 259, 269] с евклидовой прокс-структурой (об этом методе будет много написано в главах 3, 4, 6). Выбор такой прокс-структуры, прежде всего, связан с наличием ограничения $t \ge \bar{t}$ [259].

Пусть известно такое число $R$, что

$$R^2 \ge \frac{1}{2}\left\|t_* - t^0\right\|_2^2,$$

где $t_*$ – решение задачи (1). Выберем $N$ – число шагов алгоритма (далее, см. формулу (1.5.3), мы опишем, как можно выбирать $N$ исходя из желаемой точности $\varepsilon$). Положим начальное приближение $t^0 = \bar{t}$. Пусть на шаге с номером $k$ получен вектор $t^k$ и стохастический субградиент $\nabla\Upsilon\left(t^k, \xi_k\right)$, зависящий от случайного вектора $\xi_k$, и удовлетворяющий условию несмещенности $E_{\xi_k}\left[\nabla\Upsilon\left(t^k, \xi_k\right)\right] = \nabla\Upsilon\left(t^k\right)$. Здесь и далее под $\nabla\Upsilon(t)$ мы имеем в виду какой-то измеримый селектор [309] многозначного отображения $\partial\Upsilon(t)$ (субдифференциала функции $\Upsilon(t)$). Следующая точка вычисляется из соотношения

$$t^{k+1} = \arg\min_{t \ge \bar{t}}\Big\{ \underbrace{\Upsilon\left(t^k\right)}_{\text{можно не писать}} + \gamma_k \underbrace{\left\langle \nabla\Upsilon\left(t^k, \xi_k\right), t - t^k \right\rangle}_{\nabla\Upsilon\left(t^k, \xi_k\right) \text{ - стохастический субградиент}} + \frac{1}{2}\left\|t - t^k\right\|_2^2 \Big\},\ k = 0, \ldots, N. \qquad (1.5.2)$$

Здесь (используя метод двойственных усреднений [269] с $\beta_k \simeq R^{-1}\sqrt{\sum_{l=0}^{k}\left\|\nabla\Upsilon\left(t^l, \xi_l\right)\right\|_2^2}$, $\lambda_k \equiv 1$, можно избавиться от потенциально неизвестной оценки параметра $M$ в шаге $\gamma_k$)

$$\gamma_k = \frac{R}{M_k}\sqrt{\frac{2}{N+1}},$$

где

$$\max\left\{\left\|\nabla\Upsilon\left(t^k, \xi_k\right)\right\|_2, \left\|\nabla\Upsilon\left(t^k\right)\right\|_2\right\} \le M_k \le M.$$

Отметим, что такой выбор $\gamma_k$ обусловлен решением задачи

$$2R^2/\tilde{\gamma} + (N+1)M^2\tilde{\gamma} \to \min_{\tilde{\gamma} \ge 0}, \text{ где } M = \max_{k=0,\ldots,N} M_k,$$

возникающей при доказательстве теоремы 1.5.4, см. ниже.

В виду сепарабельности ограничений $t \ge \bar{t}$ и сепарабельности выражения, стоящего в фигурных скобках в (1.5.2), задача (1.5.2) на каждом шаге итерационного процесса декомпозируется на $|E| = O(n)$ одномерных подзадач, каждая из которых представляет собой задачу минимизации параболы на полуоси. Следовательно, каждая такая подзадача решается по явным формулам, т.е. на каждой итерации за $O(n)$ можно решить задачу (1.5.2) в предположении, что мы нашли стохастический субградиент $\nabla\Upsilon\left(t^k,\xi_k\right)$. Чтобы оценить насколько много потребуется итераций $N = N(\varepsilon)$ для достижения заданной точности $\Upsilon\left(\bar{t}^{N(\varepsilon)}\right) - \Upsilon_* \le \varepsilon$, где $\Upsilon_* = \min_{t\ge\bar{t}} \Upsilon(t)$, сформулируем теорему о сходимости метода (1.5.2).

**Теорема 1.5.4.** *Пусть*

$$\bar{t}^N = \frac{1}{S_N}\sum_{k=0}^{N}\gamma_k t^k, \quad S_N = \sum_{k=0}^{N}\gamma_k, \quad \bar{f}^N = \bar{f} - s^N,$$

*где* $s^N$ *– есть множитель Лагранжа к ограничению* $t \ge \bar{t}$ *в задаче*

$$\frac{1}{S_N}\left\{\sum_{k=0}^{N}\gamma_k\left\{\Upsilon\left(t^k\right) + \left\langle\nabla\Upsilon\left(t^k,\xi_k\right), t - t^k\right\rangle\right\} + \frac{1}{2}\left\|t - t^0\right\|_2^2\right\} \to \min_{t\ge\bar{t}}.$$

*Тогда с вероятностью* $\ge 1-\sigma$

$$0 \le \Upsilon\left(\bar{t}^N\right) - \Upsilon_* \le \frac{16\sqrt{2}MR}{\sqrt{N}}\ln\left(\frac{4N}{\sigma}\right)$$

*с* $R^2 = \frac{1}{2}\left\|t_* - t^0\right\|_2^2$ *и с вероятностью* $\ge 1-\sigma$

$$0 \le \Upsilon\left(\bar{t}^N\right) + \Psi\left(\bar{f}^N\right) \le \frac{16\sqrt{2}MR}{\sqrt{N}}\ln\left(\frac{4N}{\sigma}\right) \tag{1.5.3}$$

*с*

$$R^2 \ge \max\left\{\frac{1}{2}\left\|t_* - t^0\right\|_2^2, \frac{1}{2}\left\|\tilde{t}^N - t^0\right\|_2^2\right\}, \quad \tilde{t}^N = \arg\max_{t\ge 0}\left\{\sum_{w\in W} d_w T_w(t) - \left\langle\bar{f}^N, t - \bar{t}\right\rangle\right\}.$$

**Доказательство.** Согласно [165] , п. 3.4 [269]

$$0 \le \Upsilon\left(\bar{t}^N\right) - \Upsilon_* \le$$

$$\le \frac{1}{S_N}\left\{\frac{1}{2}\left\|t_* - t^0\right\|_2^2 - \frac{1}{2}\left\|t_* - t^{N+1}\right\|_2^2 - \sum_{k=0}^{N}\gamma_k\left\langle\nabla\Upsilon\left(t^k,\xi_k\right) - \nabla\Upsilon\left(t^k\right), t^k - t_*\right\rangle + \frac{1}{2}\sum_{k=0}^{N}\gamma_k^2 M_k^2\right\} \le$$

$$\le \frac{R^2}{S_N} - \frac{1}{2S_N}\left\|t_* - t^{N+1}\right\|_2^2 + \frac{1}{S_N}\sum_{k=0}^{N}\gamma_k\left\langle\nabla\Upsilon\left(t^k,\xi_k\right) - \nabla\Upsilon\left(t^k\right), t_* - t^k\right\rangle + \frac{R^2}{S_N} \le$$

$$\le \frac{\sqrt{2}MR}{\sqrt{N+1}} + \frac{1}{S_N}\sum_{k=0}^{N}\gamma_k\left\langle\nabla\Upsilon\left(t^k,\xi_k\right) - \nabla\Upsilon\left(t^k\right), t_* - t^k\right\rangle. \tag{1.5.4}$$

При формировании последнего слагаемого предпоследнего неравенства было учтено определение $\gamma_k$.

Считая, что $\frac{1}{2}\left\|t_*-t^k\right\|_2^2 \le R_\sigma^2$ для всех $k=0,...,N$ с вероятностью $\ge 1-\sigma/2$, получим из неравенства Азума–Хефдинга [159, 238] (см. также раздел 6.1 главы 6) для ограниченной мартингал–разности

$$\left|\left\langle \nabla\Upsilon\left(t^k,\xi_k\right)-\nabla\Upsilon\left(t^k\right),t_*-t^k\right\rangle\right| \le \left\|\nabla\Upsilon\left(t^k,\xi_k\right)-\nabla\Upsilon\left(t^k\right)\right\|_2\left\|t_*-t^k\right\|_2 \le 2\sqrt{2}M_kR_\sigma$$

следующее неравенство

$$P\left(\sum_{k=0}^N\gamma_k\left\langle \nabla\Upsilon\left(t^k,\xi_k\right)-\nabla\Upsilon\left(t^k\right),t_*-t^k\right\rangle \ge 2\sqrt{2}R_\sigma\Lambda\sqrt{\sum_{k=0}^N\gamma_k^2M_k^2}\right) \le \exp\left(-\Lambda^2/2\right)+\sigma/2.$$

Следовательно,

$$P\left(\sum_{k=0}^N\gamma_k\left\langle \nabla\Upsilon\left(t^k,\xi_k\right)-\nabla\Upsilon\left(t^k\right),t_*-t^k\right\rangle \ge 4\sqrt{2}R_\sigma R\sqrt{\ln\left(2/\sigma\right)}\right) \le \sigma,$$

$$P\left(\frac{1}{S_N}\sum_{k=0}^N\gamma_k\left\langle \nabla\Upsilon\left(t^k,\xi_k\right)-\nabla\Upsilon\left(t^k\right),t_*-t^k\right\rangle \ge \frac{4MR_\sigma\sqrt{\ln\left(2/\sigma\right)}}{\sqrt{N+1}}\right) \le \sigma. \qquad (1.5.5)$$

Из неравенства (1.5.4) имеем

$$\frac{1}{2}\left\|t_*-t^k\right\|_2^2 \le 2R^2+\sum_{l=0}^{k-1}\gamma_l\left\langle \nabla\Upsilon\left(t^l,\xi_l\right)-\nabla\Upsilon\left(t^l\right),t_*-t^l\right\rangle.$$

Неравенство (1.5.5) представим в виде $k=1,...,N$

$$P\left(\sum_{l=0}^{k-1}\gamma_l\left\langle \nabla\Upsilon\left(t^l,\xi_l\right)-\nabla\Upsilon\left(t^l\right),t_*-t^l\right\rangle \ge 4\sqrt{2}R_\sigma R\sqrt{\ln\left(2/\sigma\right)}\right) \le \sigma.$$

Отсюда по неравенству Буля (вероятность суммы событий, соответствующих $k=1,...,N$, не больше суммы вероятностей событий) имеем с вероятностью $\ge 1-\sigma$ для всех $k=0,...,N$

$$\frac{1}{2}\left\|t_*-t^k\right\|_2^2 \le 2R^2+4\sqrt{2}R_\sigma R\sqrt{\ln\left(2N/\sigma\right)}.$$

Положим $R_\sigma^2=2R^2+4\sqrt{2}R_\sigma R\sqrt{\ln\left(4N/\sigma\right)}$. Отсюда получаем следующую оценку

$$R_\sigma \le 4R\left(1+\sqrt{2\ln\left(4N/\sigma\right)}\right).$$

Следовательно, из неравенств (1.5.4), (1.5.5) имеем с вероятностью $\ge 1-\sigma$ неравенство

$$0 \le \Upsilon\left(\bar{t}^N\right)-\Upsilon_* \le \frac{\sqrt{2}MR}{\sqrt{N+1}}\left(1+8\sqrt{2}\left(1+\sqrt{2\ln\left(4N/\sigma\right)}\right)\sqrt{\ln\left(2/\sigma\right)}\right) \le \frac{16\sqrt{2}MR}{\sqrt{N}}\ln\left(\frac{4N}{\sigma}\right).$$

Положим теперь (см. также конец п. 3 статьи [42])

$$R^2 := \max\left\{\frac{1}{2}\left\|t_*-t^0\right\|_2^2, \frac{1}{2}\left\|\tilde{t}^N-t^0\right\|_2^2\right\},$$

где

$$\tilde{t}^N = \arg\max_{t\geq 0}\left\{\sum_{w\in W} d_w T_w(t) - \left\langle \bar{f}^N, t-\bar{t}\right\rangle\right\}.$$

Тогда

$$\Upsilon\left(\bar{t}^N\right) \leq \frac{1}{S_N}\sum_{k=0}^{N}\gamma_k \Upsilon\left(t^k\right) \leq \frac{1}{2S_N}\sum_{k=0}^{N}\gamma_k^2 M_k^2 +$$

$$+\frac{1}{S_N}\min_{t\geq\bar{t}}\left\{\sum_{k=0}^{N}\gamma_k\left\{\Upsilon\left(t^k\right)+\left\langle\nabla\Upsilon\left(t^k,\xi_k\right),t-t^k\right\rangle\right\}+\frac{1}{2}\left\|t-t^0\right\|_2^2\right\} \leq \frac{1}{2S_N}\sum_{k=0}^{N}\gamma_k^2 M_k^2 +$$

$$+\min_{t}\left\{\frac{1}{S_N}\left\{\sum_{k=0}^{N}\gamma_k\left\{\Upsilon\left(t^k\right)+\left\langle\nabla\Upsilon\left(t^k,\xi_k\right),t-t^k\right\rangle\right\}+\frac{1}{2}\left\|t-t^0\right\|_2^2\right\}+\left\langle s^N,\bar{t}-t\right\rangle\right\} \leq$$

$$\leq \frac{1}{2S_N}\sum_{k=0}^{N}\gamma_k^2 M_k^2 + \min_{t\geq 0}\left\{\frac{1}{S_N}\left\{\sum_{k=0}^{N}\gamma_k\left\langle\nabla\Upsilon\left(t^k,\xi_k\right)-\nabla\Upsilon\left(t^k\right),t-t^k\right\rangle+\frac{1}{2}\left\|t-t^0\right\|_2^2\right\}+\right.$$

$$\left.+\Upsilon(t)+\left\langle s^N,\bar{t}-t\right\rangle\right\} \leq$$

$$\leq \frac{1}{2S_N}\sum_{k=0}^{N}\gamma_k^2 M_k^2 + \frac{1}{S_N}\sum_{k=0}^{N}\gamma_k\left\langle\nabla\Upsilon\left(t^k,\xi_k\right)-\nabla\Upsilon\left(t^k\right),\tilde{t}^N-t^k\right\rangle+\frac{1}{2S_N}\left\|\tilde{t}^N-t^0\right\|_2^2-\Psi\left(\bar{f}^N\right).$$

К сожалению, в данном случае неравенство Азума–Хефдинга уже не будет выполняться, поскольку $\tilde{t}^N$ зависит от $\xi_k$ – нарушается условие, что последовательность $\left\{\gamma_k\left\langle\nabla\Upsilon\left(t^k,\xi_k\right)-\nabla\Upsilon\left(t^k\right),\tilde{t}^N-t^k\right\rangle\right\}$ будет последовательностью мартингал-разностей. Тем не менее, аналогично рассуждая, можно показать, что существует такое

$$R^2 \geq \max\left\{\frac{1}{2}\left\|t_*-t^0\right\|_2^2, \frac{1}{2}\left\|\tilde{t}^N-t^0\right\|_2^2\right\},$$

что с вероятностью $\geq 1-\sigma$

$$0 \leq \Upsilon\left(\bar{t}^N\right)+\Psi\left(\bar{f}^N\right) \leq \frac{16\sqrt{2}MR}{\sqrt{N}}\ln\left(\frac{4N}{\sigma}\right). \ \square$$

**Замечание 1.5.2.** Выписанная в теореме 1.5.4 оценка (1.5.3) была получена с помощью теоремы 3 работы [269] и замечания 4 работы [44]. Более аккуратный способ рассуждений позволяет немного улучшить оценку (1.5.3). Для практических приложений метода (1.5.2) удобнее записать его не в терминах неизвестного $R$, а сразу в терминах другого неизвестного

$$\bar{R}^2 \geq \max\left\{\max_{k=0,\ldots,N}\left\|t_*-t^k\right\|_2^2, \frac{1}{2}\left\|\tilde{t}^N-t^0\right\|_2^2\right\}.$$

Речь идет о выборе шагов

$$\gamma_k = \frac{\bar{R}}{M_k}\sqrt{\frac{2}{N+1}},$$

и об оценке [54]

$$0 \le \Upsilon\left(\bar{t}^N\right) + \Psi\left(\bar{f}^N\right) \le \frac{\sqrt{2}M\bar{R}}{\sqrt{N}}\left(1+\sqrt{8\ln\left(2/\sigma\right)}\right)$$

с вероятностью $\ge 1-\sigma$. С практической точки зрения от этого ничего не меняется, и даже становится лучше (оценка улучшается). Этот параметр как был априорно неизвестным, так и остался таковым, поменялась немного только его интерпретация. Тем не менее, оценка (1.5.3) представляет определенный теоретический интерес, поскольку является хорошей демонстрацией тех теоретических трудностей, которые возникают при описанном выше способе восстановления решения прямой задачи. Существует, однако, и другой способ восстановления решения прямой задачи, описанный в следующем замечании (см. также разделы 2, 3 главы 3).

**Замечание 1.5.3.** Пусть $\nabla\Upsilon\left(t^k, \xi_k\right) \equiv \nabla\Upsilon\left(t^k\right)$, т.е. вместо стохастического субградиента в алгоритме (1.5.2) используется обычный субградиент. Тогда вектор равновесного распределения потоков по ребрам $f$ можно считать по-другому (при таком способе подсчета возможно нарушение неравенств $\bar{f}^N \le \bar{f}$ и $0 \le \Upsilon\left(\bar{t}^N\right) + \Psi\left(\bar{f}^N\right)$):

$$f^k \in \partial \sum_{w \in W} d_w T_w\left(t^k\right),\ \bar{f}^N = \frac{1}{S_N}\sum_{k=0}^{N}\gamma_k f^k .$$

Решая возникающие здесь задачи поиска кратчайших путей, мы находим (одновременно, т.е. без дополнительных затрат) не только вектор распределения потоков по ребрам $f^k$, но и разреженный вектор распределения потоков по путям $x^k$ (при этом $f^k = \Theta x^k$). При этом

$$\Upsilon\left(\bar{t}^N\right) + \Psi\left(\bar{f}^N\right) + 2\sqrt{2}R\left\|\left(\bar{f}^N - \bar{f}\right)_+\right\|_2 \le \frac{2\sqrt{2}MR}{\sqrt{N}},$$

где $R^2 = \frac{1}{2}\left\|t_* - t^0\right\|_2^2$. Отсюда можно получить, что

$$0 \le \Upsilon\left(\bar{t}^N\right) - \Upsilon_* \le \frac{\sqrt{2}MR}{\sqrt{N}},\ \left|\Psi\left(\bar{f}^N\right) - \Psi_*\right| \le \frac{2\sqrt{2}MR}{\sqrt{N}},\ \left\|\left(\bar{f}^N - \bar{f}\right)_+\right\|_2 \le \frac{2M}{\sqrt{N}}.$$

Из теоремы 1.5.4 и слабой двойственности ($-\Upsilon_* \le \Psi_*$) имеем

$$\Upsilon\left(\bar{t}^N\right) - \Upsilon_* + \Psi\left(\bar{f}^N\right) - \Psi_* \le \Upsilon\left(\bar{t}^N\right) + \Psi\left(\bar{f}^N\right) \le \varepsilon.$$

Это обосновывает следующее следствие теоремы 1.5.4.

**Следствие 1.5.1.** *В условиях теоремы 1.5.4*

$$0 \le \Upsilon\left(\bar{t}^N\right) - \Upsilon_* \le \varepsilon, \; 0 \le \Psi\left(\bar{f}^N\right) - \Psi_* \le \varepsilon.$$

Осталось описать, как можно случайно выбирать быстро вычислимый, равномерно ограниченный по норме, стохастический субградиент функции $\Upsilon(t)$ со свойством несмещенности. Подробнее о различных способах рандомизации (именно об этом сейчас идет речь) будет написано в разделе 1.1 главы 1, главе 4 и разделе 5.1 главы 5.

Прежде всего, опишем, как формируется субградиент функции $\Upsilon(t)$. Заметим, что субградиент выпуклой негладкой функции $\partial\Upsilon(t)$ – есть выпуклое множество, превращающееся в точках гладкости $\Upsilon(t)$ в один вектор – обычный градиент $\nabla\Upsilon(t)$. Из определения $\Upsilon(t)$ (см. формулу (1.5.1)) имеем

$$\partial\Upsilon(t) = -\sum_{w\in W} d_w \partial T_w(t) + \bar{f},$$

где $\partial T_w(t)$ – супердифференциал негладкой вогнутой (как минимум выпуклых, в нашем случае аффинных) функции $T_w(t) = \min_{p\in P_w} \sum_{e\in E} \delta_{ep} t_e$ (следует сравнить с задачами, возникающими на каждом шаге метода Франк–Вульфа из подраздела 1.5.2). Супердифференциал $\partial T_w(t)$ представляет собой выпуклую комбинацию векторов с компонентами $\delta_{ep}$, отвечающих кратчайшим путям (если их несколько) для заданной корреспонденции $w$ на графе, ребра которого взвешены вектором $t = \{t_e\}_{e\in E}$. Каждый такой вектор (с числом компонент равным числу ребер) можно описать следующем образом: если ребро входит в кратчайший путь, то в компоненте вектора, отвечающего этому ребру, стоит 1, иначе 0.

Теперь опишем два варианта (первый вариант был сообщен нам Ю.Е. Нестеровым в 2013 г.) выбора несмещенного стохастического субградиента (мы вводим случайность, говорят также рандомизацию, чтобы за счет этого сократить стоимость вычисления)

**Вариант 1**

$$\nabla\Upsilon(t,\xi) = -d\nabla T_\xi(t) + \bar{f},$$

*где с.в.* $\xi = w$ *с вероятностью* $d_w/d$, $w \in W$, $\nabla T_\xi(t)$ *– произвольный элемент супердифференциала* $\partial T_\xi(t)$.

**Вариант 2**

$$\nabla\Upsilon(t,\xi) = -d \sum_{j:(\xi,j)\in W} \frac{d_{\xi j}}{d_{\xi\cdot}} \nabla T_{\xi j}(t) + \bar{f},$$

*где с.в.* $\xi = i$ *с вероятностью* $d_{i\cdot}/d$ *,* $d_{i\cdot} = \sum_{j:(i,j)\in W} d_{ij}$ *,* $i \in V$ *,* $\nabla T_{\xi}(t)$ *– произвольный элемент супердифференциала* $\partial T_{\xi}(t)$*.*

Заметим, что в варианте 2 $\nabla \Upsilon(t,\xi)$ может быть вычислен алгоритмом Дейкстры за $\tilde{O}(n)$ в силу особенность алгоритма Дейкстры поиска кратчайших путей, заключающейся в том, что за время $\tilde{O}(n)$ он находит кратчайшие пути из заданной вершины во все остальные [145]. Однако подсчет суммы в определении $\nabla \Upsilon(t,\xi)$, если это делать напрямую, может занять время $O(n^{3/2})$, поскольку число слагаемых равно $O(n)$, а число ненулевых компонент в векторах $\nabla T_{\xi j}(t)$ (число ребер в соответствующем кратчайшем пути) может быть порядка $O(\sqrt{n})$ (сеть типа двумерной решетки). Однако можно по-другому организовать вычисление компонент вектора $\nabla \Upsilon(t,\xi)$ (описываемая далее в этом абзаце конструкция была предложена Ю.В. Максимовым). Алгоритм Дейкстры строит (ориентированное) дерево кратчайших путей с корнем в $\xi$ за $\tilde{O}(n)$. Для каждой вершины этого дерева $j \neq \xi$, припишем ребру, ведущему в эту вершину вес $d_{\xi j}$. Получим таким образом взвешенное дерево. Далее припишем листьям дерева веса, равные весам ребер, ведущим в эти листья. А далее по индукции: припишем вершине дерева сумму весов всех потомков. Припишем теперь ребрам дерева новые веса: вес ребра равен весу вершины, в которую это ребро входит. Далее нужно просто пробежаться по ребрам этого дерева, считывая веса и нормируя их на $-d/d_{\xi\cdot}$. С точностью до $\bar{f}$ получим таким образом $\nabla \Upsilon(t,\xi)$. Все это можно сделать за $\tilde{O}(n)$. Аналогичные конструкции позволяют вычислять $\nabla \Upsilon(t)$ и $\Upsilon(\bar{t}^N)$ за $\tilde{O}(Sn)$.

Заметим также, что в варианте 2 оценка константы $M$ (см. теорему 4) получается заметно лучше, чем в варианте 1.

Таким образом, вариант 2 кажется более предпочтительным. В чем может быть минус использования варианта 2 – в возможности использовать более быстрые алгоритмы поиска кратчайших путей, которые находят кратчайший путь ровно между двумя вершинами. Хотя нижняя оценка затрат на поиск кратчайшего пути здесь также $\tilde{O}(n)$ для реальных (почти планарных) транспортных (и не только) сетей эта оценка может быть существенно редуцирована. В частности, в определенных ситуациях до $\tilde{O}(1)$. Как правило, это, в свою очередь, требует затрат порядка $\tilde{O}(n)$ на подготовку специальной структуры данных [145] (будем называть такой процесс препроцессингом). Тем не менее, учитывая, что веса ребер меняются от шага к шагу не сильно, такой препроцессинг не обязательно осуществлять на каждой итерации. К тому же совершенно не обязательно находить всегда кратчайшие пути, считая таким образом точный субдифференциал. Желая решить задачу с точностью $\varepsilon$ по зазору двойственности $\Upsilon\left(\bar{t}^N\right)+\Psi\left(\bar{f}^N\right)\le\varepsilon$ достаточно вычислять $O(\varepsilon)$-субдифференциал на каждой итерации [99].

Выше мы пояснили, что нельзя однозначно, исходя из теоретических оценок, отдать предпочтение одному из описанных вариантов выбора $\nabla\Upsilon(t,\xi)$. Более того, в контексте проводимого выше анализа методов может быть не очевидна и сама необходимость в рандомизации.

Действительно, на первый взгляд кажется, что если считать полностью $\nabla\Upsilon(t)$ (такой подход с методом двойственных усреднений [269] вместо зеркального спуска предлагался Ю.Е. Нестеровым в 2012 г.), то оценка (1.5.3) примет вид (заметим, что с точностью до множителя $\sqrt{2}$ эта оценка – не улучшаема [99, 165, 259])

$$\Upsilon\left(\bar{t}^N\right)+\Psi\left(\bar{f}^N\right)\le\frac{\sqrt{2}MR}{\sqrt{N}},$$

Поскольку вычисление $\nabla\Upsilon(t)$ требует $\tilde{O}(Sn)$ операций (в определенных ситуациях возможно и быстрее), то, кажется, что при таком подходе мы просто теряем фактор $n$ в оценке сложности метода. Однако в действительности это не совсем так. Во-первых, константа $M$ здесь может быть заметно меньше своего аналога в варианте 2. Во-вторых, поскольку на каждой итерации мы должны пересчитывать все кратчайшие пути, то в этой постановке также как и в подразделе 1.5.2 на передний план выходит пересчет кратчайших путей вместо расчета, что с учетом допустимости использования приближенно вычисленных кратчайших путей может редуцировать оценку сложности итерации $\tilde{O}(Sn)$ [127, 145, 150, 200].

Далее мы, тем не менее, ограничимся рассмотрением рандомизированных методов, считая стоимость одной итерации равной $\tilde{O}(n)$.

Свойство несмещенности стохастических субградиентов следует из построения (вариант 1 и 2). Сложнее обстоит дело с определением параметра метода $R$ (или $R_\sigma$, см. замечание 1.5.2; далее мы ограничимся рассмотрением случая, когда в качестве параметра выбран $R$, с помощью замечания 1.5.2 можно провести аналогичные рассуждения и для $R_\sigma$), который явно входит в итерационный процесс (1.5.2) (следует сравнить с методом

Франк–Вульфа, для которого параметры метода не входили в сам метод, только в один из вариантов критерия останова).

Критерий останова можно задавать явной формулой для числа итераций (см. теорему 1.5.4 и замечание 1.5.2), в которую входит неизвестный параметр $R$ (или $\bar{R}$), но лучше его задавать немного по-другому (см. ниже). Проблема выбора критерия останова и априорная неизвестность ряда параметров, необходимых методу для работы, и входящих в критерий останова, еще не раз будет обсуждаться по ходу диссертации в самых разных контекстах. Обзорно об этом написано в разделе 2.1 главы 2.

Прежде всего, отметим, что если известна оценка (сверху) на $M$, то ее можно использовать при выборе шага метода

$$\gamma_k \equiv \frac{R}{M}\sqrt{\frac{2}{N+1}} \text{ или } \gamma_k \equiv \frac{\bar{R}}{M}\sqrt{\frac{2}{N+1}},$$

при этом теорема 1.5.4 и замечание 1.5.2 останутся верными. Далее сконцентрируемся на замечании 1.5.2. В формуле для $\gamma_k$ стоит неизвестное $\bar{R}$, которое исчезает при подстановке сюда зависимости $N\left(\varepsilon;\bar{R},M,\sigma\right)$, определяемой в замечании 1.5.2 из условия

$$\frac{\sqrt{2}M\bar{R}}{\sqrt{N\left(\varepsilon;\bar{R},M,\sigma\right)}}\left(1+\sqrt{8\ln\left(2/\sigma\right)}\right)=\varepsilon .$$

Таким образом, шаг метода не зависит от неизвестного $\bar{R}$. В качестве, критерия остановки метода используется условие (проверяемое за $\tilde{O}\left(Sn\right)$)

$$\Upsilon\left(\bar{t}^N\right)+\Psi\left(\bar{f}^N\right)\le\varepsilon .$$

Если априорно оценка константы $M$ не известна, то процедура усложняется. Сначала предполагается, что $M = M_0$. Далее предлагается следующая процедура адаптивного подбора неизвестного параметра $\bar{R}$ (отчасти являющаяся оригинальной, см., например, [44]). Задаем какое-то начальное значение, скажем, $\bar{R}=\left\|t^0\right\|_2$, запускаем итерационный процесс (1.5.2) с этими параметрами. В какой-то момент мы либо обнаружим, что $M > M_0$, либо будет сделано предписанное (теоремой 1.5.4) для выбранных значений параметров число итераций. Предположим, что имеет место вторая альтернатива. Далее проверяем условие

$$\Upsilon\left(\bar{t}^N\right)+\Psi\left(\bar{f}^N\right)\le\varepsilon .$$

Если оно выполняется, то мы нашли решение с требуемой точностью. Если не выполняется, то запускаем процесс заново, полагая $\bar{R}:=\sqrt{2}\bar{R}$ (такой выбор константы также оптимален, см. [37]). Это дополнительно может привести к не более чем логарифмиче-скому (от отношения истинного значения $\bar{R}$ к $\left\|t^0\right\|_2$) числу перезапусков. Здесь в рассуждениях мы пренебрегли вероятностными оговорками, поскольку $\sigma$ можно считать очень малым (см. формулу (1.5.3)).

Если мы вышли из описанного цикла из-за того, что на каком-то шаге получили $M > M_0$, то полагаем $M := \sqrt{2}M$ (такой выбор константы также оптимален, см. [37]) и запускаем итерационный процесс заново с новым значением $M$. Это дополнительно может привести к не более чем логарифмическому (от отношения истинного значения $M$ к $M_0$) числу перезапусков.

В итоге, ожидаемая оценка времени работы метода (1.5.2) – есть $\tilde{O}\left(nM^2R^2/\varepsilon^2\right)$.

### 1.5.4 Рандомизированный покомпонентный спуск для модели стабильной динамики в форме Ю.Е. Нестерова

Предположим, что число источников намного меньше числа вершин

$$S = |O| \ll |V| = n.$$

Описанные в подразделах 1.5.2, 1.5.3 методы не сильно учитывают такую разреженность. В мае 2014 года Ю.В. Дорном и Ю.Е. Нестеровым была предложена следующая эквивалентная переформулировка (см. Теорему 1.5.5 в Дополнении к этому разделу 1.5) двойственной задачи (1.5.1) для модели стабильной динамики [187, 275] (т.е. при $E' = E$)

$$\min_{t \ge \bar{t}} \left\{ -\sum_{w \in W} d_w T_w(t) + \left\langle \bar{f}, t - \bar{t} \right\rangle \right\} =$$

$$= \min_{T} \left\{ -\sum_{\substack{s \in O,\, k \in D \\ (s,k) \in W}} d_{sk} \cdot \left(T_{sk} - T_{ss}\right) + \sum_{(i,j) \in E} \bar{f}_{ij} \max_{s \in O} \left(T_{sj} - T_{si} - \bar{t}_{ij}, 0\right) \right\}, \qquad (1.5.6)$$

которая, по-видимому, позволяет в большей степени учесть свойство $S \ll n$. По решению задачи (1.5.6) можно восстанавливать решение двойственной задачи (1.5.1):

$$t_{ij} = \max\left\{ \max_{s \in O}\left(T_{sj} - T_{si}\right), \bar{t}_{ij} \right\},\ (i,j) \in E,$$

но, к сожалению, нельзя восстанавливать вектор равновесного распределения потоков по путям $x$. Учитывая что

$$\min_{f \in \Delta,\, f \le \bar{f}} \sum_{e \in E} f_e \bar{t}_e = \min_{\substack{f = \Theta x,\, x \in X \\ f \le \bar{f}}} \sum_{e \in E} f_e \bar{t}_e = \min_{f = \Theta x,\, x \in X} \max_{\tau \ge 0} \left\{ \sum_{e \in E} f_e \cdot \left(\bar{t}_e + \tau_e\right) - \left\langle \bar{f}, \tau \right\rangle \right\} =$$

$$= \max_{\tau \ge 0} \left\{ \sum_{w \in W} d_w T_w\left(\bar{t} + \tau\right) - \left\langle \bar{f}, \tau \right\rangle \right\} \overset{t = \bar{t} + \tau}{=} -\min_{t \ge \bar{t}} \left\{ -\sum_{w \in W} d_w T_w(t) + \left\langle \bar{f}, t - \bar{t} \right\rangle \right\},$$

можно получить по формуле Демьянова–Данскина–Рубинова [47, 259] и решение прямой задачи

$$f \in \partial_{\bar t}\left(\min_{f\in\Delta,\, f\le \bar f}\sum_{e\in E} f_e \bar t_e\right) = \partial_{\bar t}\left(-\min_T\left\{-\sum_{\substack{s\in O,\,k\in D\\ (s,k)\in W}} d_{sk}\cdot\left(T_{sk}-T_{ss}\right)+\sum_{(i,j)\in E}\bar f_{ij}\max_{s\in O}\left(T_{sj}-T_{si}-\bar t_{ij},0\right)\right\}\right).$$

Следуя [187, 271], с помощью техники двойственного сглаживания, запишем функционал, равномерно аппроксимирующий целевой функционал задачи (1.5.6) (здесь $\eta_{ij} = \varepsilon/\left(4n\bar f_{ij}\ln(S+1)\right)$, $\varepsilon$ – точность, с которой хотим решить задачу (1.5.1)) в виде

$$-\sum_{s\in O,k\in D} d_{sk}\cdot\left(T_{sk}-T_{ss}\right)+\sum_{(i,j)\in E}\bar f_{ij}\eta_{ij}\ln\left(\frac{1}{S+1}\left[\sum_{s\in O}\exp\left(\frac{T_{sj}-T_{si}-\bar t_{ij}}{\eta_{ij}}\right)+1\right]\right)\to\min_T, \qquad (1.5.7)$$

$$f_{ij} = \bar f_{ij}\frac{\sum_{s\in O}\exp\left(\left(T_{sj}-T_{si}-\bar t_{ij}\right)/\eta_{ij}\right)}{\sum_{s\in O}\exp\left(\left(T_{sj}-T_{si}-\bar t_{ij}\right)/\eta_{ij}\right)+1} = \frac{\bar f_{ij}}{1+\left(\sum_{s\in O}\exp\left(\left(T_{sj}-T_{si}-\bar t_{ij}\right)/\eta_{ij}\right)\right)^{-1}},\ (i,j)\in E.$$

Решая задачу (1.5.7) с точностью $\varepsilon/4$ по функции прямо-двойственным методом (недавно было установлено [261], что при правильном взгляде любой разумный численный метод является прямо-двойственным, точнее имеет соответствующую модификацию) можно восстановить $\tilde t^N$ и $\tilde f^N$ так, чтобы зазор двойственности был меньше $\varepsilon$

$$\Upsilon\left(\tilde t^N\right)+\Psi\left(\tilde f^N\right)\le\varepsilon.$$

Таким образом, достаточно научиться эффективно решать задачу (1.5.7). Это можно сделать, например, с помощью ускоренного покомпонентного спуска Ю.Е. Нестерова, или более современных вариантов ускоренных покомпонентных спусков APPROX, ALPHA [199, 264, 295], ACRCD* (к сожалению, на данный момент, нам удалось показать только прямо-двойственность метода ACRCD* и его вариаций – см. раздел 5.1 главы 5). Разобьем все компоненты $\left\{T_{sk}\right\}_{s\in O,k\in V}$ на блоки $\left\{T_{(k)}\right\}_{k\in V}$, где $T_{(k)}=\left\{T_{sk}\right\}_{s\in O}$. Чтобы получить оценку скорости сходимости, нужно оценить $L_{ij,k}$ – константу липшица в 2-норме градиента функции (по переменным блока $T_{(k)}$)

$$\bar f_{ij}\eta_{ij}\ln\left(\frac{1}{S+1}\left[\sum_{s\in O}\exp\left(\frac{T_{sj}-T_{si}-\bar t_{ij}}{\eta_{ij}}\right)+1\right]\right)=$$

$$=\bar f_{ij}\max_{\substack{u_0+\sum_{s\in O}u_s=1\\ u_0,u_s\ge 0,\ s\in O}}\left\{\sum_{s\in O}\left(T_{sj}-T_{si}-\bar t_{ij}\right)u_s-\eta_{ij}u_0\ln u_0-\eta_{ij}\sum_{s\in O}u_s\ln u_s\right\}.$$

Из выписанного представления и теоремы 1 работы [271] имеем, что

$$L_{ij,k} = \bar{f}_{ij} / \eta_{ij} ,\ k = i, j ;\ L_{ij,k} = 0,\ k \neq i, j .$$

Введя

$$R_T^2 = \frac{1}{2} \left\| T^0 - T_* \right\|_2^2 ,$$

$C_{ij} \le 10$ – число “соседей” в транспортном графе у вершин $i$ и $j$,

получим, что алгоритм 2 APPROX (можно распространить все последующее и на ACRCD*) из работы [199] с $n$ блоками (размер каждого блока $S$), с евклидовой прокс-структурой в композитном варианте, с композитом

$$-\sum_{\substack{s \in O,\, k \in D \\ (s,k) \in W}} d_{sk} \cdot \left( T_{sk} - T_{ss} \right),$$

имеет, согласно теореме 3 [199], следующую оценку (в среднем) общей сложности решения задачи (1.5.7) с точностью по функции $\varepsilon/4$ ($\bar{L}$, $\bar{C}$ – специальным образом “взвешенные средние” констант $\bar{f}_{ij}/\eta_{ij}$, $C_{ij}$, причем $\bar{L} \le \max_{(i,j)\in E} \bar{f}_{ij}/\eta_{ij}$, $\bar{C} \le \max_{(i,j)\in E} C_{ij}$)

$$\underbrace{\mathrm{O}\left( n\sqrt{\frac{\bar{C}\bar{L}R_T^2}{\varepsilon}} \right)}_{\text{число итераций}} \cdot \underbrace{\mathrm{O}\left( \bar{C}S \right)}_{\substack{\text{стоимость} \\ \text{итерации}}} = \mathrm{O}\left( \bar{C}Sn\sqrt{\frac{\bar{C}\bar{L}R_T^2}{\varepsilon}} \right).$$

Если бы мы использовали не APPROX, а быстрый градиентный метод в композитном варианте [266, 271], то оценка была бы хуже

$$\mathrm{O}\left( Sn\sqrt{\frac{n\bar{L}R_T^2}{\varepsilon}} \right).$$

Заметим, что алгоритм можно распараллелить на $\mathrm{O}\left( n/\bar{C} \right)$ процессорах, каждый процессор при этом должен сделать

$$\mathrm{O}\left( \bar{C}^2 S\sqrt{\frac{\bar{C}\bar{L}R_T^2}{\varepsilon}} \right)$$

арифметических операций. При $S \ll n$ получается довольно оптимистичная оценка.

### 1.5.5 Заключение

В разделе были описаны различные способы поиска равновесного распределения потоков в моделях Бэкмана (подраздел 1.5.2), в модели стабильной динамики Нестерова–де Пальма (раздел 4) и промежуточных моделях (подраздел 1.5.3). Промежуточные модели получаются из модели Бэкмана при предельном переходе $\mu \to 0+$ по части ребер (см. подраздел 1.5.2). Модель стабильной динамики получается, когда предельный переход осуществлен на всех ребрах (см. подраздел 31.5.). Подход подраздела 1.5.2 не применим к моделям подразделов 1.5.3, 1.5.4, поскольку (см. подраздел 1.5.2)

$$L_2 = \max_{e \in E} \tau_e'^{\mu}\left(\hat{f}_e\right) \to \infty \text{ при } \mu \to 0+.$$

Подход раздела 1.4 не применим к моделям Бэкмана, поскольку, по сути, ограничивается только решением специальной задачи линейного программирования. Подход подраздела 1.5.3 применим ко всем рассмотренным в разделе моделям.

Введя относительную точность по функции $\tilde{\varepsilon}$ (см. подраздел 1.5.2), с которой мы хотим искать равновесие (в реальных приложениях относительной точности $\tilde{\varepsilon} \sim 0.01 - 0.05$ оказывается более чем достаточно), запишем оценки общего времени работы методы из подразделов 1.5.2 – 1.5.4:

Таблица 1.5.1

| метод / модель | **Бэкман** | **Нестеров–де Пальма** | **Промежуточная** |
|---|---|---|---|
| **Франк–Вульф** | $K_2(S,n)Sn/\tilde{\varepsilon}$ | $\varnothing$ | $\varnothing$ |
| **Зеркальный спуск** | $K_3(S,n)n/\tilde{\varepsilon}^2$ | $K_3(S,n)n/\tilde{\varepsilon}^2$ | $K_3(S,n)n/\tilde{\varepsilon}^2$ |
| **Покомпонентный спуск** | $\varnothing$ | $K_4(S,n)Sn/\tilde{\varepsilon}$ | $\varnothing$ |

К сожалению, константы $K$ этих методов могут существенно отличаться и зависеть, в частности, от $n$. В худших случаях может быть $K = O(n)$.

Заметим, что ожидаемую (среднюю) сложность $\tilde{O}\left(\max\{n,|W|\}^3\right)$ поиска равновесия в модели стабильной динамики дает симплекс метод [158, 303, 318] (для прямой задачи – задачи линейного программирования). Из Таблицы 1.5.1 можно заключить, что симплекс метод будет конкурентоспособным лишь при небольшом числе корреспонденций $|W| \le n$. Отметим при этом, что общую сложность $\tilde{O}\left(n^3\right)$ имеет заметно более простая транспортная задача линейного программирования с $n$ пунктами производства и потребления, причем эта оценка является не улучшаемой [132, 293].

Интересно сравнить приведенные в Таблице 1.5.1 оценки с оценками сложности решения задачи поиска стохастического равновесия в модели Бэкмана из работы [33].

В заключение сделаем существенную для практики оговорку. На полученные в данном разделе оценки следует смотреть исключительно с точки зрения качественного понимания сложности того или иного метода, но ни коим образом не с точки зрения отбора лучшего метода. В данном разделе отсутствуют сравнительный анализ констант методов из Таблицы 1.5.1. Такой анализ требует проработки некоторых технических деталей, свя-занных с дополнительным погружением в специфику постановки. При выбранном в данном разделе уровне грубости получения оценок можно считать, что все три метода (из подразделов 1.5.2 – 1.5.4) конкурентоспособны. Помочь отобрать наилучший метод здесь могут численные эксперименты. Такие эксперименты сейчас активно проводит группа А.Ю. Горнова из Иркутска. Результаты работы по этому направлению должны войти в кандидатскую диссертацию А.С. Аникина (ИПМ РАН, 2016).

### 1.5.6 Дополнение

Приведем доказательство использованного нами в подразделе 1.5.4 результата (1.5.6) (Дорна–Нестерова, 2014).

**Теорема 5.** *Задача (1.5.1) эквивалентна задаче*

$$\min_{T}\left\{-\sum_{\substack{s\in O,\,k\in D\\(s,k)\in W}} d_{sk}\cdot\left(T_{sk}-T_{ss}\right)+\sum_{(i,j)\in E}\bar{f}_{ij}\max_{s\in O}\left(T_{sj}-T_{si}-\bar{t}_{ij},0\right)\right\}. \tag{1.5.8}$$

*Решение задачи (1.5.8) связано с решением задачи (1.5.1) следующим образом*

$$t_{ij}=\max\left\{\max_{s\in O}\left(T_{sj}-T_{si}\right),\bar{t}_{ij}\right\},\ (i,j)\in E. \tag{1.5.9}$$

**Доказательство.** Двойственная задача для задачи (1.5.1) имеет вид

$$\min_{t\geq\bar{t}}\left\{-\sum_{(s,k)\in W} d_{sk}\cdot T_{sk}(t)+\sum_{(i,j)\in E}\bar{f}_{ij}\cdot\left(t_{ij}-\bar{t}_{ij}\right)\right\}. \tag{1.5.10}$$

Введем набор переменных $T_{sk}$ – время проезда по кратчайшему пути из вершины $s$ в вершину $k$. Тогда для любых трех вершин $s$, $i$, $j$, таких, что существует ребро $(i,j)\in E$ выполнено соотношение

$$T_{sj}\leq T_{si}+t_{ij}.$$

Следовательно, задачу (1.5.10) можно переписать в виде

$$\min_{t,T}\left\{-\sum_{(s,k)\in W} d_{sk}\cdot T_{sk}+\sum_{(i,j)\in E}\bar{f}_{ij}\cdot\left(t_{ij}-\bar{t}_{ij}\right)\,\middle|\, t\geq\bar{t};\ T_{sj}\leq T_{si}+t_{ij},\, s\in V,\,(i,j)\in E\right\}. \tag{1.5.11}$$

Это задача ЛП. В её целевом функционале все компоненты вектора $t$ неотрицательны, при этом ограничение на каждую компоненту имеет вид

$$t_{ij}\geq\max\left\{\max_{s}\left[T_{sj}-T_{si}\right],\bar{t}_{ij}\right\}. \tag{1.5.12}$$

Следовательно, задача (1.5.11) может быть явно решена относительно вектора $t$ – в (1.5.12) имеет место равенство, т.е. имеет место (1.5.9). Исключая вектор $t$, приходим к задаче (1.5.8), что завершает доказательство. □

## 1.6 Поиск равновесий в многостадийных транспортных моделях

### 1.6.1 Введение

Поиск (стохастических) равновесий в многостадийных моделях транспортных потоков приводит к решению следующей седловой задачи с правильной (выпукло-вогнутой) структурой [9, 31, 35, 47] (см. также подраздел 1.1.10 раздела 1.1 этой главы 1, а также разделы 1.3, 1.4 этой главы 1):

$$\min_{\substack{\sum_{j=1}^{n} x_{ij}=L_i,\ \sum_{i=1}^{n} x_{ij}=W_j \\ x_{ij}\ge 0,\ i,j=1,\dots,n}} \max_{y\in Q} \left\{ \sum_{i,j=1}^{n} x_{ij}\ln x_{ij} + \sum_{i,j=1}^{n} c_{ij}(y)x_{ij} + g(y) \right\} =$$

$$= \max_{y\in Q} \max_{\lambda,\mu\in\mathbb{R}^n} \left\{ \langle \lambda, L\rangle + \langle \mu, W\rangle - \sum_{i,j=1}^{n} \exp\left(-c_{ij}(y) - 1 + \lambda_i + \mu_j\right) + g(y) \right\}, \qquad (1.6.1)$$

где $g(y)$ и $c_{ij}(y)\ge 0$ – вогнутые гладкие функции (если ищутся не стохастические равновесия, то $c_{ij}(y)$ могут быть негладкими), $Q$ – множество простой структуры, например,

$$Q = \{y:\ y \ge \bar{y}\}.$$

Легко понять, что система балансовых ограничений в (1.6.1) либо несовместна $\sum_{i=1}^{n} L_i \ne \sum_{j=1}^{n} W_j$ , либо вырождена (имеет не полный ранг). В последнем случае это приводит к тому, что двойственные переменные $(\lambda,\mu)$ определены с точностью до произвольной постоянной $C$:

$$(\lambda + Ce, \mu - Ce),\ e = \underbrace{(1,\dots,1)}_{n}.$$

Задачу (1.6.1) также можно переписать следующим образом (не ограничивая общности, считаем $\sum_{i=1}^{n} L_i = \sum_{j=1}^{n} W_j = 1$)

$$\min_{\substack{\sum_{j=1}^{n} x_{ij}=L_i,\ \sum_{i=1}^{n} x_{ij}=W_j \\ x_{ij}\ge 0,\ i,j=1,\dots,n;\ \sum_{i,j=1}^{n,n} x_{ij}=1}} \max_{y\in Q} \left\{ \sum_{i,j=1}^{n} x_{ij}\ln x_{ij} + \sum_{i,j=1}^{n} c_{ij}(y)x_{ij} + g(y) \right\} =$$

$$= \max_{y\in Q} \max_{\lambda,\mu\in\mathbb{R}^n} \left\{ \langle \lambda, L\rangle + \langle \mu, W\rangle - \ln\left( \sum_{i,j=1}^{n} \exp\left(-c_{ij}(y) + \lambda_i + \mu_j\right) \right) + g(y) \right\} = -\min_{y\in Q} f(y), \quad (1.6.2)$$

где выпуклая функция $f(y)$ определяется как

$$f(y)=\max_{\substack{\sum_{j=1}^{n}x_{ij}=L_i,\ \sum_{i=1}^{n}x_{ij}=W_j\\ x_{ij}\ge 0,\ i,j=1,\dots,n;\ \sum_{i,j=1}^{n,n}x_{ij}=1}}\left\{-\sum_{i,j=1}^{n}x_{ij}\ln x_{ij}-\sum_{i,j=1}^{n}c_{ij}(y)x_{ij}-g(y)\right\}=$$

$$=\min_{\lambda,\mu}\left\{\ln\left(\sum_{i,j=1}^{n}\exp\left(-c_{ij}(y)+\lambda_i+\mu_j\right)\right)-\langle\lambda,L\rangle-\langle\mu,W\rangle-g(y)\right\}. \qquad (1.6.3)$$

Поскольку мы добавили в ограничения условие $\sum_{i,j=1}^{n,n}x_{ij}=1$, являющееся следствием балансовых уравнений, то это привело к тому, что двойственные переменные $(\lambda,\mu)$ определены с точностью до двух произвольных постоянных $C_\lambda$, $C_\mu$: $(\lambda+C_\lambda e,\mu+C_\mu e)$.

В данном разделе мы покажем, как можно решать задачу (1.6.2).

Заметим также, что расчет градиента $\nabla f(y)$ (в ряде транспортных приложений вогнутые функции $c_{ij}(y)$ – негладкие, тогда вместо градиентов стоит понимать суперградиенты $c_{ij}(y)$ и субградиент $f(y)$) осуществляется по следующей формуле (Демьянова–Данскина–Рубинова, см., например, [37, 47])

$$\nabla f(y)=-\frac{\sum_{i,j=1}^{n}\exp\left(-c_{ij}(y)+\lambda_i^*+\mu_j^*\right)\nabla c_{ij}(y)}{\sum_{i,j=1}^{n}\exp\left(-c_{ij}(y)+\lambda_i^*+\mu_j^*\right)}-\nabla g(y)=-\sum_{i,j=1}^{n}x_{ij}\left(\lambda^*,\mu^*\right)\nabla c_{ij}(y)-\nabla g(y), \qquad (1.6.4)$$

где $(\lambda^*,\mu^*)$ – решение задачи (1.6.3), не важно какое именно, градиент $\nabla f(y)$ от выбора $C_\lambda$, $C_\mu$ (см. выше) не зависит. В данном разделе мы ограничимся изучением только полноградиентных методов для задачи (1.6.2), т.е. не будем рассматривать, например, рандомизацию при вычислении градиента по формуле (1.6.4). Планируется отдельно исследовать вопрос о возможности ускорения вычислений за счет введения рандомизации для внешней задачи. На данный момент нам представляется (см. формулу (1.6.13) в подразделе 1.6.3), что это может принести дивиденды только в случае, когда вспомогательная задача расчета $\nabla c_{ij}(y)$ достаточно сложная, в свою очередь. Тут требуется много оговорок, в частности, в большинстве приложений умение рассчитывать $\nabla c_{ij}(y)$ для конкретной пары $(i,j)$ без дополнительных затрат позволяет заодно рассчитать и все $\nabla c_{ij}(y)$, $j=1,\dots,n$. Также отдельно планируется исследовать вопрос о том какие подходы и насколько хорошо допускают распараллеливание. Вопрос о целесообразности рандомизации оказывается завязанным и на вопрос о возможности распараллеливания.

Структура раздела следующая. В подразделе 1.6.2 мы рассматриваем популярную в последнее время (в связи с большим числом приложений) задачу вычисления барицентра Вассерштейна различных вероятностных мер. Эта задача оказывается тесно связанной с задачей (1.6.1). Мы разбираем в разделе этот пример, потому что он хорошо проясняет возможные альтернативы предлагаемому нами основному подходу решения задач (1.6.1), (1.6.2), изложенному в подразделе 1.3. В основе оригинального подхода подразделе 1.6.3 лежит сочетание метода балансировки для решения внутренней задачи оптимизации по двойственным множителям и универсального метода с неточным оракулом для внешней задачи (1.6.2). В заключительном подразделе 1.6.4 делаются замечания относительно возможности ускорения подхода, описанного в подразделе 1.6.3.

### 1.6.2 Поиск барицентра Вассерштейна[20]

К похожей на (1.6.1) задаче приводит поиск барицентра Монжа–Канторовича (в западной литературе чаще говорят барицентра Вассерштейна [131, 176])[21]. Изложим вкратце постановку задачи [149, 176, 178]. Вводится энтропийно сглаженное транспортное расстояние (см. рис. 1 в [175]), с матрицей $\left\|c_{ij}\right\|_{i,j=1}^{n.n}$, сформированной из квадратов попарных расстояний $c_{ij}=l_{ij}^2$ от носителя меры $i$ до носителя меры $j$ ($L,W\in S_n(1)$):

$$\Delta(L,W)=\min_{\substack{\sum_{j=1}^{n}x_{ij}=L_i,\ \sum_{i=1}^{n}x_{ij}=W_j\\ x_{ij}\geq 0,\ i,j=1,\dots,n}}\left\{\gamma\sum_{i,j=1}^{n}x_{ij}\ln x_{ij}+\sum_{i,j=1}^{n}c_{ij}x_{ij}\right\}=$$

$$=\max_{\lambda,\mu}\left\{\langle\lambda,L\rangle+\langle\mu,W\rangle-\gamma\sum_{i,j=1}^{n}\exp\left(\frac{-c_{ij}+\lambda_i+\mu_j}{\gamma}-1\right)\right\}=$$

$$=\max_{\lambda}\left\{\langle\lambda,L\rangle-\underbrace{\gamma\sum_{j=1}^{n}W_j\ln\left(\frac{1}{W_j}\sum_{i=1}^{n}\exp\left(\frac{-c_{ij}+\lambda_i}{\gamma}\right)\right)}_{H_W^*(\lambda)}\right\}. \quad (1.6.5)$$

---

[20] Задача этого пункта была поставлена В.Г. Спокойным. Интересные результаты численных экспериментов по сравнительному анализу алгоритмов из этого пункта с известными алгоритмами имеются у А.Л. Сувориковой, П.Е. Двуреченского и А.В. Чернова. Эти результаты должны войти в кандидатские диссертации А.Л. Сувориковой и А.В. Чернова.

[21] Строго говоря, мы будем искать барицентр вероятностных мер не согласно настоящему (негладкому) расстоянию Вассерштейна (на наш взгляд исторически более правильно это расстояние называть расстоянием Монжа–Канторовича–Добрушина), как это можно было подумать из названия, а согласно энтропийно-сглаженному расстоянию Вассерштейна [176].

Определим при $L \in S_n(1)$ функцию $H_W(L) = \Delta(L,W)$. Эта гладкая на $L \in S_n(1)$ функция с градиентом (см. утверждение 3 [178]):

$$\nabla H_W(L) = \lambda^*,$$

где $\lambda^*$ единственное решение (1.6.5), удовлетворяющее условию[22] $\langle \lambda^*, e \rangle = 0$. Отсюда следует, что

$$H_W^*(\lambda) = \max_{L \in S_n(1)} \left\{ \langle \lambda, L \rangle - H_W(L) \right\} = \gamma \sum_{j=1}^n W_j \ln\left( \frac{1}{W_j} \sum_{i=1}^n \exp\left( \frac{-c_{ij} + \lambda_i}{\gamma} \right) \right).$$

Теперь можно перейти к изложению основной конструкции. Задача поиска барицентра Вассерштейна[23] записывается следующим образом:

$$\sum_{k=1}^m H_{W_k}(L) \to \min_{L \in S_n(1)}. \tag{1.6.6}$$

К сожалению, в такой формулировке мы не можем оценить константу Липшица градиента функционала (1.6.6), явно входящую в большинство современных быстрых (ускоренных) численных методов. Однако оказывается (см. подразделе 1.6.3), что существуют быстрые методы, которым для работы не требуется такая информация (константа Липшица градиента).

Перепишем задачу (1.6.6), следуя п. 3 работы [178], следующим образом

$$-\sum_{k=1}^m H_{W_k}(L_k) \to \max_{\substack{L_1 = L_m \mid \lambda^1 \\ \dots\dots\dots \\ L_{m-1} = L_m \mid \lambda^{m-1} \\ L_1, \dots, L_m \in S_n(1)}},$$

$$\sum_{k=1}^{m-1} \max_{L_k \in S_n(1)} \left\{ \langle \lambda^k, L_k \rangle - H_{W_k}(L_k) \right\} + \max_{L_m \in S_n(1)} \left\{ \left\langle -\sum_{k=1}^{m-1} \lambda^k, L_m \right\rangle - H_{W_m}(L_m) \right\} \to \min_{\lambda^1, \dots, \lambda^{m-1} \in \mathbb{R}^n},$$

$$\sum_{k=1}^{m-1} H_{W_k}^*(\lambda^k) + H_{W_m}^*\left( -\sum_{k=1}^{m-1} \lambda^k \right) \to \min_{\lambda^1, \dots, \lambda^{m-1} \in \mathbb{R}^n}, \tag{1.6.7}$$

---

[22] Решая задачу (1.6.6), каким-нибудь прокс-методом с KL прокс-структурой [91], легко понять, что от того, как именно выбирать $\lambda^*$, задаваемое с точностью до сдвига всех компонент на одно и то же произвольное число, метод зависеть не будет. Единственное для чего имеет смысл стремиться к выполнению этого нормирующего условия, так это для лучшей практической обработки (меньшее накопление ошибок округления из-за конечной длинны мантиссы) экспоненциального взвешивания компонент градиента, возникающего на каждом шаге итерационного процесса при выборе KL прокс-структуры.

[23] К сожалению, пока не так много известно о статистической обоснованности использования расстояния Вассерштейна. Другими словами, хотелось бы иметь связь барицентра Вассерштейна с оценками максимального правдоподобия, ну или хотя бы с состоятельными оценками для соответствующих схем экспериментов. Пока установлена только связь с состоятельными оценками [153, 156].

$$L_* = \nabla H^*_{W_k}\left(\lambda^k_*\right) \text{ для любого } k = 1,...,m-1,$$

где $L_*$ – единственное решение задачи (1.6.6), $\left\{\lambda^k_*\right\}_{k=1}^{m-1}$ – единственное решение задачи (1.6.7). Важное свойство функционала задачи (1.6.7) – равномерная ограниченность константы Липшица градиента (следует из [271]). Задача безусловной минимизации (1.6.7) может быть эффективно решена различными способами (в зависимости от того насколько велики $n$ и $m$). В частности, для больших $n$ и $m$ неплохо с задачей справляются различные модификации метода сопряженных градиентов и быстрых градиентных методов [178]. Структура задачи (1.6.7) позволяет эффективно использовать покомпонентные методы (см., например, [199, 295]), которые к тому же хорошо параллелятся для данной задачи. Задача (1.6.7) хорошо также решается с помощью распределенных вычислений [160].

В приложениях к поиску разладки требуется много раз перерешивать задачу (1.6.7), которую для симметричности перепишем следующим образом

$$\sum_{k=1}^{m} H^*_{W_k}\left(\lambda^k\right) \to \min_{\substack{\lambda^1,...,\lambda^m \in \mathbb{R}^n \\ \sum_{k=1}^{m} \lambda^k = 0}},$$

немного смещая окно, т.е. заменяя каждый раз несколько первых слагаемых в сумме

$$H^*_{W_1}\left(\lambda^1\right),..,H^*_{W_r}\left(\lambda^r\right)$$

на столько же новых (которые, как ожидается, близки к $H^*_{W_m}\left(\lambda^m\right)$). В таком случае предлагается в итерационном процессе стартовать при сдвиге окошка с того, на чем остановились на прошлом положении окошка, экстраполируя $\lambda^m_*$ на вновь пришедшие слагаемые. Ясно, что для новой задачи набор

$$\left( \lambda^{r+1}_*,...,\lambda^{m-1}_*,\lambda^m_*,\underbrace{\lambda^m_*,...,\lambda^m_*}_{r} \right),$$

с которого стартуем, уже не будет оптимальным, однако, мы вправе надеяться на его близость к оптимальному набору, что существенно сокращает число последующих итераций. Интересной, особенно в данном контексте, представляется возможность использования (и интерпретации) распределенных вычислений [160].

Может показаться, что подход, сводящий поиск решения задачи (1.6.6) к задаче (1.6.7), не доминируем, поскольку в отличие от задачи (1.6.6), в задаче (1.6.7) мы можем явно выписать функционал и по простым формула рассчитать градиент, который к тому же имеет равномерно ограниченную константу Липшица. С одной стороны, это, действительно, преимущество, но получено оно дорогой ценной – ценной раздутия пространства, в котором происходит оптимизация почти в $m$ раз. И это раздутие скажется не только на сложности одной итерации. Для задачи (1.6.6) осуществление одной итерации будет еще более дорогим в виду необходимости на каждой итерации решать $m$ отдельных подзадач расчета $\nabla H_{W_k}(L)$. Скажется это, прежде всего, на числе необходимых итераций. В следующем разделе будет отмечено, что расчет $\nabla H_{W_k}(L)$ с помощью метода балансировки не намного сложнее расчета $\nabla H^*_{W_k}(\lambda^k)$. При этом задача (1.6.6) решается в пространстве намного меньшей размерности, и мы вправе ожидать, что необходимое число итераций может быть намного меньше, чем для задачи (1.6.7). Кроме того, задача (1.6.6) решается на компакте, т.е. размер решения (если быть точным, то расстояние от точки старта до решения), входящий в оценку необходимого числа итераций, заведомо ограничен размером симплекса. Задача (1.6.7) – задача безусловной оптимизации, причем без свойства сильной выпуклости функционала. Размер ее решения может быть большим, и входит он в оценки необходимого числа итераций также как и для задачи (1.6.6), к сожалению, степенным образом (для быстрых (ускоренных) методов можно ожидать линейной зависимости необходимого числа итераций от этого размера). Наконец, для постановок задач об обнаружении разладки (см. выше) также ожидается, что использовать близость решений прямых задач (1.6.6) при смещении окошка удастся намного лучше, чем близость в решении двойственных задач (1.6.7). В итоге, выгода от подхода, связанного с переходом к задаче (1.6.7), уже не столь очевидна, и требует отдельного и более аккуратного исследования, с решающей ролью численных экспериментов.

В ряде задач требуется искать параметрический барицентр Вассерштейна. В таком случае в одном из вариантов постановки предполагают наличие параметрической зависимости $L(\theta) \in S_n(1)$, $\theta \in \Theta$, где размерность вектора параметров $\dim\theta \ll n$. К сожалению, в этом случае нельзя гарантировать с помощью стандартных приемов (стр. 86 [161]) выпуклости задачи

$$\sum_{k=1}^{m} H_{W_k}\left(L(\theta)\right) \to \min_{\theta \in \Theta}, \qquad (1.6.8)$$

за исключением случая, когда $L(\theta) = A\theta + b$, а $\Theta$ – выпуклое множество. В этом случае конструкция (1.6.7) видоизменяется следующим образом

$$-\sum_{k=1}^{m} H_{W_k}(L_k) \to \max_{\substack{L_1 = L_m \mid \lambda^1 \\ \dots\dots\dots \\ L_{m-1} = L_m \mid \lambda^{m-1} \\ L_m = A\theta + b \mid \tilde{\lambda} \\ L_1,\dots,L_m \in S_n(1),\ \theta \in \Theta}},$$

$$\sum_{k=1}^{m-1} H_{W_k}^{*}\left(\lambda^k\right) + H_{W_m}^{*}\left(-\sum_{k=1}^{m-1}\lambda^k - \tilde{\lambda}\right) + \left\langle \tilde{\lambda}, A\theta + b \right\rangle \to \min_{\substack{\lambda^1,\ldots,\lambda^{m-1},\tilde{\lambda}\in\mathbb{R}^n \\ \theta\in\Theta}}, \tag{1.6.9}$$

$$L_* = \nabla H_{W_k}^{*}\left(\lambda_*^k\right) \text{ для любого } k = 1,\ldots,m-1.$$

Как следствие, нет никаких гарантий, что изложенная выше конструкция, связанная с переходом к двойственной задаче (1.6.7), и восстановлению решения прямой задачи (1.6.6) по явным формулам через двойственные множители, в общем случае будет работать хотя бы для поиска локальных решений.[24] Другими словами, необходимо искать глобальный оптимум задачи (1.6.8), исходя из работы с прямой задачей (1.6.8). Один из вариантов того, как это можно делать, будет описан в следующем пункте.[25]

Однако при другом варианте постановки (более предпочтительном) можно задавать зависимость $L(\theta)$ с помощью аффинных равенств и выпуклых неравенств

$$\sum_{k=1}^{m} H_{W_k}\left(L\right) \to \min_{\substack{A\theta + BL = c \\ g(\theta, L)\le 0 \\ L\in S_n(1);\, \theta\in\Theta}}.$$

Многие параметрические зависимости можно загнать в такое представление [259]. В частности, отметим полиэдральные представления Фурье–Моцкина [259], возникающие, например, в робастной оптимизации

$$\sum_{k=1}^{m} H_{W_k}\left(L\right) \to \min_{L\in\left\{L\in S_n(1):\, \exists\,\theta:\, A\theta + BL\le c\right\}}.$$

Можно переписать задачи (1.6.7), (1.6.9) и на эти случаи, причем сделать это корректно в том смысле, что правомочность подхода полностью сохранится. При этом принципиально ничего из сказанного выше не поменяется.

В действительности, в приложениях наиболее интересен случай, когда ищется барицентр именно расстояний Вассерштейна,[26] а не энтропийно-сглаженных расстояний [149,

---

[24] Впрочем, есть результаты (см. формулу (8) п. 3 § 2 главы 8 [99]) о локальной сходимости обычного градиентного спуска для задачи (1.6.9) при некоторых дополнительных предположениях.

[25] При этом правомочность подхода п. 3 для постановки задачи (1.6.8) имеет место при дополнительном предположении, что метод стартует из выпуклой окрестности точки минимума с небольшим запасом, допускающим возможность по ходу итерационного процесса оказаться дальше от решения, чем в начальный момент.

[26] Численные методы поиска “честного” барицентра Вассерштейна вероятностных мер в основном строятся на том, что когда меры заданы на прямой, существуют эффективные способы решения задачи поиска барицентра [131]. Далее проектируют (считают преобразования Радона) меры на случайные прямые и решают одномерные задачи. По их решениям восстанавливают решение исходной задачи [157, 296]. В отличие от других подходов, здесь существенно используется структура матрицы затрат $c_{ij} = l_{ij}^2$ (в дискретном случае). Интересно было бы исследовать вопрос о применимости этого подхода к постановкам

175, 176, 178]. Другими словами, интересно изучать предельное поведение $\gamma \to 0+$ (см. п. 3.1 [149], утверждение 1 [178], п. 3 и конец п. 4 [175]). К сожалению, методы из подразделов 1.6.2, 1.6.3 оказываются весьма чувствительными к этому предельному переходу. Для метода из этого раздела константа Липшица градиента в задаче (1.6.7) будет расти как $\gamma^{-1}$, соответственно, число итераций будет увеличиваться (при использовании быстрых (ускоренных) методов) как $\gamma^{-1/2}$. Еще более плохое поведение (см. [37]) можно ожидать от метода балансировки, использующегося в подходе из подраздела 1.6.3. Вопрос о том, как действовать при малых $\gamma > 0$ (в частности, в вырожденном случае $\gamma = 0$) изучался П.Е. Двуреченским, А.Л. Суворикова и А.В. Черновым. Имеется гипотеза, что в этом случае поможет философия искусственного сглаживания[27] [271], в которой искусственно введенная энтропийная регуляризация уже задается с четко заданным коэффициентом регуляризации $\gamma > 0$, зависящим от итоговой точности по функции, с которой требуется решить задачу. Другой способ – использовать менее чувствительные (чем метод балансировки) способы решения двойственной задачи, например, при небольших значениях $n$ ожидается, что лучше сработает $r$-метод Н.З. Шора и некоторые его обобщения [113, 114]. В данном разделе мы фиксируем $\gamma > 0$ и далее уже не будем возвращаться к подобного рода вопросам.

В заключение этого раздела отметим, что поиск барицентра Вассерштейна в случае $m=1$ может быть осуществлен явно: $L = W$. Обоснование этого частного результата представляется довольно полезным для понимания основной конструкции этого раздела.

---

задач о разладках, в которых требуется много раз пересчитывать барицентр (см. выше). Также интересно было бы сравнить описанные подходы с остальными. Этому планируется посвятить отдельную публикацию.

[27] Это сглаживание правильно называть двойственным сглаживанием, поскольку для того чтобы добиться гладкости в прямой негладкой задаче, которая имеет Лежандровское (седловое) представление [259, 271], в это представление, которое также можно понимать как двойственное, вводят аддитивным образом с небольшим коэффициентом сильно выпуклый (вогнутый) функционал. Этот функционал и обеспечивает гладкость (а еще точнее Липшицевость градиента) в прямой задаче. В нашем случае мы исходим из задачи о перемещении масс (Монж–Канторович), являющейся задачей ЛП. Однако для удобства вычисления расстояний Вассерштейна мы перешли к двойственной задаче. Мы хотим сделать гладкой двойственную задачу, потому что именно с ней в дальнейшем и идет работа. Для этого двойственное сглаживание (в нашем случае энтропийное) применяется к двойственной задаче для двойственной задачи к транспортной задаче, т.е. применяется просто к транспортной задаче.

### 1.6.3 Универсальный метод с неточным оракулом

Из подраздела 1.6.2 следует, что внутренняя задача максимизации по $(\lambda,\mu)$ может быть явно решена по $\mu$ при фиксированном $\lambda$, и наоборот (это верно для задач (1.6.1) и (1.6.2), и приводит к одним и тем же формулам). Собственно, таким образом, получается метод балансировки расчета матрицы корреспонденций по энтропийной модели, см., например, [37] (тесно связанный с методом Синхорна [149, 175]), как метод простой итерации для явно выписываемых условий экстремума (принципа Ферма): $\lambda = \Lambda(\mu)$, $\mu = \mathrm{M}(\lambda)$.

Метод балансировки имеет вид ($[\lambda]_0 = [\mu]_0 = 0$) – см., например, подраздел 3.1.7 раздела 3.1 главы 3:

$$[\lambda_i]_{k+1} = -\ln\left(\frac{1}{L_i}\sum_{j=1}^{n}\exp\left(-c_{ij}(y) - 1 + [\mu_j]_k\right)\right),$$

$$[\mu_j]_{k+1} = -\ln\left(\frac{1}{W_j}\sum_{i=1}^{n}\exp\left(-c_{ij}(y) - 1 + [\lambda_i]_k\right)\right)$$

или

$$[\mu_j]_{k+1} = -\ln\left(\frac{1}{W_j}\sum_{i=1}^{n}\exp\left(-c_{ij} - 1 + [\lambda_i]_{k+1}\right)\right).$$

В этих формулах “–1” в экспоненте для метода (1.6.2) (в отличие от метода (1.6.1)) можно не писать, поскольку двойственные множители определяются неоднозначным образом с бо́льшим произволом для задачи (1.6.2) (см. выше), достаточным для справедливости этого замечания.

Оператор $(\lambda,\mu) \to (\Lambda(\mu), \mathrm{M}(\lambda))$ является сжимающим в метрике Биркгофа–Гильберта $\rho$ [204]. Это означает, после $N \sim \ln(\sigma^{-1})$ итераций метода балансировки можно получить такие $(\lambda_N, \mu_N)$, что ($\{(\lambda_*(y), \mu_*(y))\}$ – двумерное аффинное множество решений, см. подраздел 1.6.1)

$$\rho\left((\lambda_N, \mu_N); \{(\lambda_*(y), \mu_*(y))\}\right) \le \sigma. \qquad (1.6.10)$$

Причем на практике наблюдается очень быстрая сходимость, т.е. коэффициент пропорциональности небольшой [37]. Таким образом, мы можем приближенно решить внутреннюю задачу.[28]

[28] Заметим, что в пп. 3.1, 3.2 работы [149] предлагается за счет раздутия прямого пространства с помощью обобщения описанного метода балансировки Брэгмана (метода проекций Брэгмана) решать задачу поиска

Далее предлагается воспользоваться прямо-двойственным (эта важно, поскольку нужно восстанавливать двойственные переменные) универсальным методом [274] (можно использовать и УМТ, описанный в разделе 2.2 главы 2) для решения внешней задачи оптимизации по $y$ (имеется видео/презентация Ю.Е.Нестерова с описанием этого метода). К сожалению, в формулировке (1.6.1) (в отличие от формулировки (1.6.2)) кроме того что внешняя задача гладкая (при условии гладкости $c_{ij}(y)$ [33, 40]), больше ничего о ней сказать нельзя (константа Липшица градиента не ограничена). Также не понятна гладкость задачи (1.6.6). Поэтому и по ряду других причин, о которых будет сказано далее, было отдано предпочтение универсальному методу, оптимально адаптивно настраивающемуся на гладкость функционала $f(y)$ на текущем участке пребывания итерационного процесса.[29] Однако нам потребуется использовать этот метод в варианте (описанном П.Е. Двуреченским) с неточным оракулом, выдающим градиент [41].

Напомним (см. подраздел 1.6.1), что мы решаем задачу (1.6.2), представимую в виде (здесь в max представлении $x = x$, $\bar{Q} = \left\{ x_{ij} \geq 0,\ i, j = 1, ..., n : \sum_{j=1}^{n} x_{ij} = L_i, \sum_{i=1}^{n} x_{ij} = W_j \right\}$, а в min представлении $x = (\lambda, \mu)$, $\bar{Q} = \mathbb{R}^{2n}$ – см. формулу (1.6.3)):

$$f(y) = \max_{x \in \bar{Q}} \Psi(x, y) = \min_{x \in \bar{Q}} \Phi(x, y) \to \min_{y \in Q}.$$

Далее везде будем предполагать, что $y \in Q$.

---

барицентра напрямую, т.е. отпадает необходимость в решении внешней задачи. Плата за это достаточно большая – увеличение размера прямого пространства в $m$ раз, но метод при этом хорошо параллелится.

[29] Бытует мнение, что любой универсальный метод должен чем-то платить за свою универсальность, и в этой связи возникает много вопросов, в частности: насколько дорога эта плата? В принципе, в статье [274] (см. также раздел 2.2 главы 2) довольно подробно проясняется этот момент. Тем не менее, мы повторим здесь соображения из [274]. Действительно, плата за универсальность есть. Универсальный метод из работы [274] может сделать где-то в 4 раза больше обращений к оракулу (что можно понимать как увеличения числа итераций в 4 раза) для задачи с, более менее, одинаковой константой Липшица градиента во всей области (где довелось пройти итерационному процессу), по сравнению с обычным быстрым градиентным методом [271]. Тем не менее, замечательная особенность универсального метода не только в том, что он настраивается на гладкость задачи и применимым к любым задачам, но и в том, что (в отличие от подавляющего большинства методов) этот метод локально настраивается на гладкость функционала. И для сильно неоднородных функционалов типично, что универсальный метод делает заметно меньше итераций, чем, скажем, быстрый градиентный метод (плата за это уже учтена в отмеченном выше потенциально возможном увеличении числа итераций в 4 раза в худшем случае). Примеры, поясняющие сказанное, имеются в работе [274].

**Определение 1.6.1 (см. главу 4 [181]).** $(\delta, L)$*-оракул выдает (на запрос, в котором указывается только одна точка* $y$ *) такие* $(F(y), G(y))$*, что и для любых* $y, y' \in Q$

$$0 \le f(y') - F(y) - \langle G(y), y' - y \rangle \le \frac{L}{2}\|y' - y\|^2 + \delta .$$

Из определения 1.6.1 разу следует, что для любого $x \in Q$

$$F(x) \le f(x) \le F(x) + \delta$$

и для любых $x, y \in Q$

$$f(y) \ge f(y) - \langle G(x), y - x \rangle - \delta .$$

Из последнего свойства получаем, что определение $(\delta, L)$-оракула можно понимать как обобщение на гладкие задачи классического понятия негладкой выпуклой оптимизации: $\delta$-субградиента (см. п. 5 § 1 главы 5 [99]). В приводимом далее утверждении в первой его части следует сохранить обозначения для задачи (1.6.2), (1.6.3) и следует обозначить $x = \lambda$, $y = L$ для задачи (1.6.5), (1.6.6); а во второй части утверждения следует обозначить $x = (\lambda, \mu)$, $y = y$ для задачи (1.6.2), (1.6.3). Таким образом, на задачу (1.6.2), (1.6.3) можно посмотреть с двух разных ракурсов, однако второй ракурс менее привлекателен в виду необходимости рассмотрения ограниченых множеств $\bar{Q}$, что в интересующих нас приложениях место не имеет.

**Утверждение 1.6.1.** *Если* $\psi(y) = \max_{x \in \bar{Q}} \Psi(x, y)$*, где* $\Psi(x, y)$ *– выпуклая по* $y$ *и вогнутая по* $x$ *функция, и найден такой* $\tilde{x} \in \bar{Q}$*, что*

$$\psi(y) - \Psi(\tilde{x}, y) \le \delta ,$$

*то субградиент* $\partial_y \Psi(\tilde{x}, y)$ *– есть* $\delta$*-субградиент функции* $\psi(y)$ *в точке* $y$.

*Если* $\varphi(y) = \min_{x \in \bar{Q}} \Phi(x, y)$*, где* $\Phi(x, y)$ *– выпуклая по совокупности переменных функция, и найден такой* $\tilde{x} \in \bar{Q}$*, что*

$$\max_{z \in \bar{Q}} \langle \Phi_x(\tilde{x}, y), \tilde{x} - z \rangle \le \delta ,$$

*то*

$$\Phi(\tilde{x}, y) - \varphi(y) \le \delta$$

*и субградиент* $\Phi_y(\tilde{x}, y) = \partial_y \Phi(\tilde{x}, y)$ *– есть* $\delta$*-субградиент функции* $\varphi(y)$ *в точке* $y$.

**Доказательство.** Ограничимся доказательством только второй части этого утверждения. Доказательство первой части см. на стр. 124 (лемма 13) книги [99].

Из выпуклости $\Phi(x,y)$ по совокупности переменных имеем

$$\Phi(x',y') \ge \Phi(x,y) + \left\langle \Phi_x(x,y), x'-x \right\rangle + \left\langle \Phi_y(x,y), y'-y \right\rangle . \tag{1.6.11}$$

Определим зависимость $x(y)$ из соотношения

$$\varphi(y) = \min_{x \in Q} \Phi(x,y) = \Phi(x(y),y).$$

Заметим, что $\Phi(\tilde{x},y) \ge \varphi(y)$. Положим в (11) $x' = x(y')$, $x = \tilde{x}$. Тогда

$$\varphi(y') = \Phi(x',y') \ge \Phi(\tilde{x},y) + \left\langle \Phi_x(\tilde{x},y), x'-\tilde{x} \right\rangle + \left\langle \Phi_y(\tilde{x},y), y'-y \right\rangle \ge$$

$$\ge \varphi(y) + \left\langle \Phi_x(\tilde{x},y), x(y')-\tilde{x} \right\rangle + \left\langle \Phi_y(\tilde{x},y), y'-y \right\rangle \ge \varphi(y) + \left\langle \Phi_y(\tilde{x},y), y'-y \right\rangle - \delta .$$

В последней формуле мы использовали, что

$$\left\langle \Phi_x(\tilde{x},y), \tilde{x} - x(y') \right\rangle \le \delta .$$

В свою очередь, из выпуклости $\Phi(x,y)$ по $x$ (для всех допустимых $y$), имеем

$$\Phi(\tilde{x},y) - \Phi(x(y'),y) \le \left\langle \Phi_x(\tilde{x},y), \tilde{x} - x(y') \right\rangle .$$

Беря в этой формуле $y' = y$ и воспользовавшись определением $x(y)$, получаем, что

$$\Phi(\tilde{x},y) - \varphi(y) \le \left\langle \Phi_x(\tilde{x},y), \tilde{x} - x(y) \right\rangle . \bullet$$

Однако не хочется довольствоваться возможностью находить только $\delta$-субградиент (из утверждения 1.6.1 эта возможность очевидна), поскольку в определенных ситуациях явно можно рассчитывать на некоторую гладкость итоговой (внешней) задачи (1.6.2). Понятие $(\delta, L)$-оракула в некотором смысле налагает наиболее слабые условия на возможные неточности в вычислении функции и градиента, при которых можно рассчитывать, что скорость сходимости метода, учитывающего гладкость (Липшицевость градиента функционала) задачи, сильно не пострадает (см. теорему 1 ниже).

На первый взгляд может показаться, что применимость описанной концепции $(\delta, L)$-оракула к задаче (1.6.1) следует из следующего результата (см. п. 4.2.2 [181]).

**Утверждение 1.6.2.** *Пусть подзадача энтропийно-линейного программирования (ЭЛП) в (1.6.2) решена (по функции) с точностью* $\delta$, *т.е. найден такой* $\tilde{x}(c)$, *удовлетворяющей балансовым ограничениям, что*

$$\sum_{i,j=1}^{n} \tilde{x}_{ij}(c) \ln \tilde{x}_{ij}(c) + \sum_{i,j=1}^{n} c_{ij} \tilde{x}_{ij}(c) - \min_{\substack{\sum_{j=1}^{n} x_{ij} = L_i,\ \sum_{i=1}^{n} x_{ij} = W_j \\ i,j=1,\dots,n}} \left\{ \sum_{i,j=1}^{n} x_{ij} \ln x_{ij} + \sum_{i,j=1}^{n} c_{ij} x_{ij} \right\} \le \delta .$$

*Тогда для функции*

$$\bar{f}(c) = -\min_{\substack{\sum_{j=1}^{n} x_{ij}=L_i,\ \sum_{i=1}^{n} x_{ij}=W_j \\ i,j=1,\dots,n}} \left\{ \sum_{i,j=1}^{n} x_{ij}\ln x_{ij} + \sum_{i,j=1}^{n} c_{ij}x_{ij} \right\}$$

*набор*

$$-\left( \sum_{i,j=1}^{n} \tilde{x}_{ij}(c)\ln \tilde{x}_{ij}(c) + \sum_{i,j=1}^{n} c_{ij}\tilde{x}_{ij}(c), \left\{\tilde{x}_{ij}(c)\right\}_{i,j=1}^{n,n} \right)$$

*является* $\left(\delta, 2\cdot \max_{i,j=1,\dots,n} c_{ij}\right)$*-оракулом.*

К сожалению, большинство методов (в том числе метод балансировки) не удовлетворяют одному пункту утверждения 1.6.2, а именно, они выдают вектор $\tilde{x}$, который лишь приближенно удовлетворяет балансовым ограничениям (в утверждении требование точного удовлетворения балансовых ограничений является существенным, и не может быть как-то равнозначно релаксировано). Связанно это с тем, что для задачи ЭЛП, когда ограничений намного меньше числа прямых переменных, обычно решается двойственная задача, по которой восстанавливается решение прямой задачи [37, 138]. Как следствие приобретается невязка и в ограничениях. Собственно, в представлении градиента функционала по формуле (1.6.4) имеются два способа. Первый через двойственные множители $(\lambda, \mu)$, второй через решение прямой задачи $x$. Функционал прямой задачи сильно выпуклый по $x$, поскольку энтропия 1-сильно выпуклая функция в 1-норме [271]. Поэтому сходимость в решении прямой задачи по функции обеспечивает сходимость и по аргументу, что и означает возможность определения с хорошей точностью градиента по формуле (1.6.4) через $x$. Другая ситуация возникает, если смотреть на двойственную задачу к задаче ЭЛП (в приводимом далее утверждении следует обозначить $x = (\lambda, \mu)$, $y = y$ для задачи (1.6.2), (1.6.3)).

**Утверждение 1.6.3.** *Пусть* $\varphi(y) = \min_{x\in \bar{Q}} \Phi(x,y)$, *где* $\Phi(x,y)$ *– такая достаточно гладкая, выпуклая по совокупности переменных функция, что*[30]

$$\left\|\nabla\Phi(x',y') - \nabla\Phi(x,y)\right\|_2 \le L\left\|(x',y') - (x,y)\right\|_2 .$$

*Пусть для произвольного* $y \in Q$ *(считаем, что множество* $Q$ *содержит внутри себя шар радиуса более* $\sqrt{2\delta/L}$ *) можно найти такой* $\tilde{x} = \tilde{x}(y) \in \bar{Q}$, *что*

[30] Это утверждение имеет достаточно простую геометрическую интерпретацию. Проекция надграфика выпуклой функции будет выпуклым множеством, то есть, в свою очередь, надграфиком некоторой выпуклой функции. Кривизна границы у полученного при проектировании множества будет не больше, чем была у исходного множества. Это следует из того, что проектирование – сжимающий оператор.

$$\max_{z\in Q}\left\langle \Phi_x\left(\tilde{x},y\right),\tilde{x}-z\right\rangle \le \delta .$$

*Тогда*

$$\Phi\left(\tilde{x},y\right)-\varphi\left(y\right)\le \delta ,$$

$$\left\|\nabla\varphi\left(y'\right)-\nabla\varphi\left(y\right)\right\|_2 \le L\left\|y'-y\right\|_2 ,$$

*и* $\left(\Phi\left(\tilde{x},y\right)-2\delta,\Phi_y\left(\tilde{x},y\right)\right)$ *будет* $\left(6\delta,2L\right)$*-оракулом для* $\varphi\left(y\right)$ *на выпуклом множестве, полученном из множества* $Q$ *отступанием от границы* $\partial Q$ *во внутрь* $Q$ *на расстояние* $\sqrt{2\delta/L}$ *(по условию это множество не пусто).*

**Доказательство.** По условию задачи имеем при всех допустимых значениях аргументов $\Phi$:

$$\lambda_{\max}\left(\begin{Vmatrix}\Phi_{xx} & \Phi_{xy}\\ \Phi_{yx} & \Phi_{yy}\end{Vmatrix}\right)=\sup_{\|h\|_2\le 1}\left\langle h,\begin{Vmatrix}\Phi_{xx} & \Phi_{xy}\\ \Phi_{yx} & \Phi_{yy}\end{Vmatrix}h\right\rangle \le L . \qquad (1.6.12)$$

Заметим, что также по условию при всех допустимых значениях аргументов $\Phi$:

$$\begin{Vmatrix}\Phi_{xx} & \Phi_{xy}\\ \Phi_{yx} & \Phi_{yy}\end{Vmatrix}\succ 0,\ \Phi_{xx}\succ 0,\ \Phi_{yy}\succ 0,\ \Phi_{yx}=\Phi_{xy}^T,\ \Phi_{xx}=\Phi_{xx}^T,\ \Phi_{yy}=\Phi_{yy}^T .$$

Для упрощения последующих рассуждений (в частности, чтобы не работать с псевдообратными матрицами) будем, считать, что матрица $\Phi_{xx}\succ 0$ положительно определена (исходя из условий, гарантировать можно лишь неотрицательную определенность). Также будем считать (в интересующих нас приложениях к задачам (2), (1.6.6) это имеет место), что зависимость $x\left(y\right)$, определяемая из соотношения

$$\varphi\left(y\right)=\min_{x\in Q}\Phi\left(x,y\right)=\Phi\left(x\left(y\right),y\right)$$

однозначным образом, и удовлетворяет соотношению

$$\Phi_x\left(x\left(y\right),y\right)\underset{y}{\equiv}0 ,$$

из которого имеем

$$\Phi_{xx}\left(x\left(y\right),y\right)\left\|\frac{\partial x\left(y\right)}{\partial y}\right\|+\Phi_{xy}\left(x\left(y\right),y\right)\underset{y}{\equiv}0 ,$$

т.е.

$$\left\|\partial x/\partial y\right\|=\left\|\partial x_i/\partial y_j\right\|=-\Phi_{xx}^{-1}\Phi_{xy} .$$

Поскольку $\varphi\left(y\right)=\Phi\left(x\left(y\right),y\right)$, то

$$\varphi_{yy}=\left\|\partial x/\partial y\right\|^T\Phi_{xx}\left\|\partial x/\partial y\right\|+\left\|\partial x/\partial y\right\|^T\Phi_{xy}+\Phi_{yx}\left\|\partial x/\partial y\right\|+\Phi_{yy}=\Phi_{yy}-\Phi_{yx}\Phi_{xx}^{-1}\Phi_{xy} .$$

С учетом этой формулы и из формулы дополнения по Шуру [329], получаем

$$\begin{Vmatrix}\Phi_{xx} & \Phi_{xy}\\ \Phi_{yx} & \Phi_{yy}\end{Vmatrix} = \begin{Vmatrix}E_x & 0\\ \Phi_{yx}\Phi_{xx}^{-1} & E_y\end{Vmatrix}\begin{Vmatrix}\Phi_{xx} & 0\\ 0 & \varphi_{yy}\end{Vmatrix}\begin{Vmatrix}E_x & \Phi_{xx}^{-1}\Phi_{xy}\\ 0 & E_y\end{Vmatrix},$$

где $E_x$, $E_y$ – единичные матрицы соответствующих размеров. Поскольку

$$\begin{Vmatrix}E_x & 0\\ \Phi_{yx}\Phi_{xx}^{-1} & E_y\end{Vmatrix} = \begin{Vmatrix}E_x & \Phi_{xx}^{-1}\Phi_{xy}\\ 0 & E_y\end{Vmatrix}^T,$$

и эти матрицы полного ранга, то из (1.6.12) имеем, что

$$\sup_{\|h\|_2\le 1}\langle h,\varphi_{yy}\, h\rangle = \lambda_{\max}(\varphi_{yy}) \le \max\left\{\lambda_{\max}(\Phi_{xx}),\lambda_{\max}(\varphi_{yy})\right\} = \lambda_{\max}\left(\begin{Vmatrix}\Phi_{xx} & \Phi_{xy}\\ \Phi_{yx} & \Phi_{yy}\end{Vmatrix}\right) \le L.$$

Таким образом, установлено, что

$$\varphi(y) \le \varphi(x) + \langle \nabla\varphi(x), y-x\rangle + \frac{L}{2}\|y-x\|_2^2.$$

Согласно утверждению 1.6.1

$$\varphi(y) \ge \varphi(x) + \langle \Phi_y(\tilde{x},y), y-x\rangle - \delta.$$

Далее проведем рассуждения аналогично стр. 107 (и не много отлично от стр. 115) диссертации [181]. Вычитая из первого неравенства второе, получим

$$\langle \Phi_y(\tilde{x},y) - \nabla\varphi(x), y-x\rangle \le \frac{L}{2}\|y-x\|_2^2 + \delta.$$

Положим ($t>0$)

$$y - x = t\frac{\Phi_y(\tilde{x},y)-\nabla\varphi(x)}{\|\Phi_y(\tilde{x},y)-\nabla\varphi(x)\|_2},$$

получим

$$\|\Phi_y(\tilde{x},y)-\nabla\varphi(x)\|_2 \le \frac{Lt}{2} + \frac{\delta}{t}.$$

Минимизируя правую часть неравенства по $t>0$, получим (при $t=\sqrt{2\delta/L}$), что

$$\|\Phi_y(\tilde{x},y)-\nabla\varphi(x)\|_2 \le \sqrt{2\delta L}.$$

Отсюда и из утверждения 1 имеем, что

$$\varphi(y) \le \varphi(x) + \langle \nabla\varphi(x), y-x\rangle + \frac{L}{2}\|y-x\|_2^2 \le$$

$$\le \varphi(x) + \langle \Phi_y(\tilde{x},y), y-x\rangle + \sqrt{2\delta L}\,\|y-x\|_2 + \frac{L}{2}\|y-x\|_2^2 \le$$

$$\le \Phi(\tilde{x},y) - 2\delta + \langle \Phi_y(\tilde{x},y), y-x\rangle + \sqrt{2\delta L}\,\|y-x\|_2 + \frac{L}{2}\|y-x\|_2^2 + 2\delta \le$$

$$\le \Phi(\tilde{x}, y) - 2\delta + \langle \Phi_y(\tilde{x}, y), y - x \rangle + \frac{2L}{2}\|y - x\|_2^2 + 6\delta .$$

С учетом того, что (см. утверждение 1.6.1)

$$\varphi(y) \ge \varphi(x) + \langle \Phi_y(\tilde{x}, y), y - x \rangle - \delta \ge \Phi_y(\tilde{x}, y) - 2\delta + \langle \Phi_y(\tilde{x}, y), y - x \rangle,$$

из определения 1.6.1 получаем доказываемое утверждение. ●

Это утверждение, позволяет установить гладкость задачи (1.6.2), (1.6.3) (но не (1.6.5), (1.6.6)). Таким образом, для (1.6.5), (1.6.6) необходимость использования универсального метода для внешней задачи является отражением надежды сходиться быстрее, чем в негладком случае, в то время как для (1.6.2), (1.6.3) использование универсального метода для внешней задачи является скорее отражением желания настраиваться на правильную константу Липшица градиента. Можно, конечно, пытаться использовать приведенные выше формулы, однако из способа рассуждений (см., например, доказательство утверждения 1.6.3) видно, что полученная таким образом константа Липшица может оказаться завышенной.

К сожалению, практическое применение утверждения 1.6.3 натыкается на следующие сложности:

1) необходимости отступать от границы множества $Q$ во внутрь на $\sqrt{2\delta/L}$,

2) необходимости рассмотрения ситуации (см. доказательство утверждения 1.6.3)

$$\left\| \partial x_i / \partial y_j \right\| = -\Phi_{xx}^{-1} \Phi_{xy},$$

3) необходимости предположения о компактности множества $\bar{Q}$, иначе невозможно будет добиться выполнения условия

$$\max_{z \in \bar{Q}} \langle \Phi_x(\tilde{x}, y), \tilde{x} - z \rangle \le \delta .$$

Сложность 1, как правило, на практике преодолима за счет возможности доопределения функционала задачи с сохранением всех свойств на $\sqrt{2\delta/L}$ окрестность множества $Q$ (заметим, что доопределение часто не требуется, поскольку функционал и так задан "с запасом"). Например, для рассматриваемых нами транспортных приложений с $Q = \{y:\ y \ge \bar{y}\}$ это возможно [9, 31, 35, 47]. Сложность 2 часто вообще не возникает (разве что оговорка о существовании $\Phi_{xx}^{-1}$, впрочем, приведенные выше рассуждения можно провести, сохранив все результаты в идентичном виде, так, что эта оговорка будет не нужна), поскольку $\bar{Q}$ совпадает со всем (двойственным) пространством. А вот сложность 3, действительно, портит дело. К сожалению, простых теоретически обоснованных способов борьбы с этой сложностью мы пока не знаем. Тем не менее, полезно заметить, что в действительности нужно гарантировать выполнение (см. доказательство утверждения 1.6.1)

$$\left\langle \Phi_x\left(\tilde{x}(y), y\right), \tilde{x}(y) - x(y')\right\rangle \le \delta,$$

где точки $y$ и $y'$ близки, поскольку возникают на соседних итерациях внешнего метода. С учетом ожидаемой “близости” $\tilde{x} = \tilde{x}(y)$ и $x(y)$, мы можем заменить в этом критерии настоящее множество $\bar{Q}$, которое, как правило, совпадает со всем пространством, на шар конечного радиуса. Более детальные исследования (для задачи (1.6.2), (1.6.3)) и практические эксперименты показывают, что для выполнения приведенного выше условия достаточно обеспечить для внутреннего итерационного процесса $\{x_k\} \to x(y)$ условия

$$\left\|\Phi_x\left(x_k, y\right)\right\|_2 \left\|x_k\right\|_2 \le \delta/2,\ \left\|\Phi_x\left(x_k, y\right)\right\|_2 \le \delta.$$

Соответствующее $x_k = (\lambda_k, \mu_k)$ порождает нужное $\tilde{x}(y) = x_k$. С учетом специфики рассматриваемой нами задачи (1.6.2), (1.6.3), имеем следующий критерий (возвращаемся к обозначениям (1.6.2), (1.6.3))

$$\left\|Ax(\lambda_k, \mu_k) - b\right\|_2 \left\|(\lambda_k, \mu_k)\right\|_2 \le \delta/2, \left\|Ax(\lambda_k, \mu_k) - b\right\|_2 \le \delta,$$

где $x(\lambda_k, \mu_k)$ определяется в формуле (1.6.4), а введённая линейная система балансовых уравнений $Ax = b$ – есть общая запись аффинных (транспортных) ограничений:

$$\sum_{j=1}^{n} x_{ij} = L_i, \sum_{i=1}^{n} x_{ij} = W_j,\ i, j = 1,\ldots,n.$$

В связи со сказанным выше заметим, что (это следует из оценки (1.6.10)) метод балансировки обеспечивает сходимость и по аргументу, что для других методов (без введения регуляризации) решения двойственной задачи, вообще говоря, нельзя гарантировать. Это свойство наряду с линейной скоростью сходимости метода (со скоростью геометрической прогрессии) позволяет надеяться, что выбранный критерий является достаточно точным (точнее не слишком грубым).

Принципиально важно для гладкого случая ($c_{ij}(y)$ – функции с Липшицевым градиентом), как это будет следовать из дальнейших оценок (см. теорему 1), не просто уметь решать двойственную задачу, т.е. находить $(\lambda, \mu)$ так, чтобы была сходимость по аргументу, а делать это так, чтобы сложность решения задачи зависела от точности ее решения логарифмическим образом. Выше мы отмечали, что это имеет место для метода балансировки. Также это имеет место и для быстрых методов, примененных к регуляризованной двойственной задачи. При фиксации параметра регуляризации, исходя из итоговой желаемой точности, быстрые градиентные методы (для сильно выпуклых функций) ре-шают регуляризованную двойственную задачу так, что зависимость сложности от точности ее решения – логарифмическая.

Хочется, чтобы при решении внешней задачи в (1.6.2), т.е. задачи

$$\min_{y\in Q} f(y),$$

можно было не задумываться ни о какой гладкости. Если она есть, то метод бы это хорошо учитывал, не требуя знания констант Липшица градиента (это намного более существенно для возможности применять описанный подход к поиску барицентра Вассерштейна вероятностных мер, см. подраздел 1.6.2), если ее нет, то метод также бы работал оптимальным (для негладкого случая) образом. Именно таким свойством и обладает универсальный метод [274], работающий и в концепции неточного оракула [41] (см. определение 1.6.1).

Заметим [274], что можно погрузить задачу с гёльдеровым градиентом ($\nu \in [0,1]$)

$$\|\nabla f(y') - \nabla f(y)\|_* \le L_\nu \|y' - y\|^\nu$$

(в том числе и негладкую задачу с ограниченной нормой разности субградиентов при $\nu = 0$) в класс гладких задач с оракулом, характеризующимся точностью $\delta$ и (см. также разделы 2.1, 2.2 главы 2)

$$L = L_\nu \left[ \frac{L_\nu (1-\nu)}{2\delta(1+\nu)} \right]^{\frac{1-\nu}{1+\nu}}.$$

Это позволяет даже в случае, когда можно рассчитывать только на $\delta$-субградиент[31] (с ограниченной нормой субградиента (разности субградиентов), причем какой именно константой ограниченной, методу знать не обязательно), все равно работать в концепции $(\delta, L)$-оракула.

Итак, у нас есть внешняя задача (1.6.2)

$$\min_{y\in Q} f(y),$$

для которой обращение к $(\delta, L)$-оракулу за значением функции и градиента стоит $\sim \ln(\delta^{-1})$. Насколько быстро мы можем решить такую задачу, т.е. при каком $N(\varepsilon)$ можно гарантировать, что

$$f\left(y_{N(\varepsilon)}\right) - \min_{y\in Q} f(y) \le \varepsilon \,?$$

Ответ можно получить из следующего результата.

---

[31] На $\delta$-субградиент всегда можно рассчитывать согласно утверждению 1.6.1. Причем, как уже отмечалось раньше, для получения $\delta$-субградиента не нужна сходимость по аргументу для вспомогательной задачи.

**Теорема 1.6.1 (см. [41, 42, 259, 274]).** *Существует однопараметрическое семейство универсальных градиентных методов (параметр $p \in [0,1]$), не получающих на вход кроме $p$ больше никаких параметров (в частности, не использующих значения $L_\nu$ и $R$ – "расстояние" от точки старта до решения, априорно не известное), которое приводит к следующей оценке на требуемое число итераций*

$$N_p(\varepsilon) = \mathrm{O}\left(\inf_{\nu \in [0,1]} \left(\frac{L_\nu R^{1+\nu}}{\varepsilon}\right)^{\frac{2}{1+2p\nu+\nu}}\right),$$

*если* $\delta \le \mathrm{O}\left(\varepsilon / N_p(\varepsilon)^p\right)$.

Из теоремы 1.6.1 можно заключить, что если мы рассчитываем на некоторую гладкость $f(y)$, то стоит выбирать значение параметра $p = 1$, при этом общее трудозатраты машинного времени будут

$$\mathrm{O}\left(N_1(\varepsilon)\left(T\ln\left(\varepsilon^{-1}\right) + \tilde{T}\right)\right), \qquad (1.6.13)$$

где $\tilde{T}$ – время вычисления (суб-)градиента функционала (в основном это вычисления $\left\{\nabla c_{ij}(y)\right\}_{i,j=1}^{n,n}$ [33, 40]), $T$ – время решения вспомогательной задачи методом балансировки с относительной точностью 1%. Численные эксперименты показывают, что на одном современном ноутбуке при $n \sim 10^2$ время $T \approx 1$ сек. [41], что сопоставимо с временем $\tilde{T}$ для таких $n$ [40].

Выгода, от описанной выше конструкции, по сравнению с обычным способом решения исходной задачи минимизации (1.6.2), (1.6.3) сразу по совокупности всех переменных (см., например, [35]) заключается в гарантированном не увеличении константы Липшица градиента в оценке необходимого числа итераций (см. утверждение 1.6.3) и ожидаемое уменьшение в этой же оценке "расстояния от точки старта до (неизвестного априори) решения. Выгода здесь вполне может достигать одного порядка и более. При этом можно ожидать лишь незначительного увеличения стоимости одной итерации. Причем стоит иметь в виду, что при оптимизации сразу по всем переменным требуется рассчитывать градиента функционала по большему числу переменных, чем в описанном выше подходе, что также играет нам на пользу. В конечном итоге, сокращение числа итераций заметно превалирует над небольшим увеличением стоимости одной итерации.

Что касается задачи (1.6.5), (1.6.6), то описанный выше подход представляется естественным и не имеющих альтернатив в рассматриваемом классе. Альтернативные методы, с которыми можно сравнивать, мы упоминали по ходу этого раздела 1.6 и раздела 1.4 этой главы, но все они были предложены из принципиально других подходов.

Резюмируем ключевой результат этого раздела (и всего раздела) следующим образом.

*Для решения задачи (1.6.2) (или (1.6.6) или (1.6.8)) предлагается использовать универсальный метод из работы [274] (а точнее его модификацию из [41]). Если рассчитываем на гладкость*[32] $f(y)$, *то полагаем в методе* $p=1$. *Если на гладкость рассчитывать не приходится*[33], *то полагаем* $p=0$. *В обоих случаях кроме априорной подсказки относительно параметра* $p$, *методу больше ничего от нас знать не надо!*

### 1.6.4 Заключительные замечания

В приложениях часто возникают задачи, имеющие следующий вид (см., например, [44, 154])

$$f(x)=\Phi\big(x,y(x)\big)\to\min_x \tag{1.6.14}$$

при этом $y(x)$ и $\nabla y(x)$ могут быть получены из решения отдельной подзадачи лишь приближенно. Довольно типично, что существует способ, который выдает $\varepsilon$-приближенное решение за время, зависящее от $\varepsilon$, логарифмическим образом $\sim\ln\left(\varepsilon^{-1}\right)$. В данном разделе намечен общий способ решения таких задач. Его наиболее важной отличительной чертой является адаптивность (самонастраиваемость), т.е. методу на вход не надо подавать никаких констант Липшица (более того, метод будет работать и в негладком случае). Метод сам настраивается локально на оптимальную гладкость функции. Это свойство метода делает его привилегированным, поскольку в реальных приложениях, чтобы что-то знать о свойствах $f(x)$ нужно что-то знать о свойствах зависимости $y(x)$, а это часто не доступно по постановке задачи, или, при попытке оценить, приводит к сильно завышенным оценкам.

В задаче (1.6.14) важно уметь эффективно пересчитывать значения $y(x)$, а не рассчитывать их каждый раз заново (на каждой итерации внешнего цикла). Поясним сказан-

---

[32] В этом случае как раз существенна логарифмическая сложность приближенного вычисления градиента и значения функции $f(y)$ от точности, обеспеченная методом балансировки.

[33] В этом случае точность решения вспомогательной задачи расчета $\delta$-субградиента можно завязать на желаемую точность решения задачи (1.6.2) $\varepsilon$ (или (1.6.6) или (1.6.8)) по формуле $\delta=O(\varepsilon)$ (см. теорему 1.6.1 при $\nu=0$), с константой порядка 1.

ное. Предположим, что мы уже как-то посчитали, скажем, $y(x)$, решив, например, с какой-то точностью соответствующую задачу оптимизации. Тогда для вычисления $y(x+\Delta x)$ (на следующей итерации внешнего цикла) у нас будет хорошее начальное приближение $y(x)$. А, как известно (см., например, [44]), расстояние от точки старта до решения (не в случае сходимости метода со скоростью геометрической прогрессии или быстрей) существенным образом определяет время работы алгоритма. Тем не менее, известные нам приложения (см., [44, 154]) пока как раз всецело соответствуют сходимости процедуры поиска $y(x)$ со скоростью геометрической прогрессии. Связано это с тем, что если расчет $y(x)$ с точностью $\varepsilon$ осуществляется за $O\left(\ln\left(\varepsilon^{-1}\right)\right)$ операций, то для внешней задачи можно выбирать самый быстрый метод (а, стало быть, и самый требовательный к точности), и с точностью до того, что стоит под логарифмом, общая трудоемкость будет рассчитываться по формуле, аналогичной формуле (1.6.13). Как правило, такое сочетание оказывается недоминируемым. Здесь мы ограничимся ссылкой на пример 4 и последующий текст из работы [44] (см. также раздел 2.1 главы 2) и общим тезисом, который пока неплохо подтверждался на практике:

*если есть возможность в задаче оптимизации (в седловой задаче) явно прооптимизировать по части переменных (или как-то эффективно это сделать с хорошей точностью), то, как правило, это и надо сделать, и строить итерационный метод исходя из этого.*

В реальных транспортных приложениях [31, 33, 40, 47] достаточно сложным является расчет $c_{ij}(y)$ и их градиентов (особенно при поиске стохастических равновесий). Тем не менее, эти задачи имеют вполне четкую привязку к решению некоторых задач на графах типа поиска кратчайших путей (см. раздел 1.5 этой главы 1, а также раздел 3.3 главы 3). Также как и в предыдущем абзаце для внутренней задачи, для внешней задачи, можно не рассчитывать $c_{ij}(y)$ и их градиенты каждый раз заново, а пересчитывать; также можно допускать неточность в их вычислении (и ненулевую вероятность ошибки) надеясь на ускорение (благо метод работает в концепции неточного оракула, природа которой не принципиальна, см. [181]). Также здесь оказываются полезными идеи БАД (быстрого автоматического дифференцирования [65, 76]), которые позволяют практически за тоже время, что занимает вычисление самих функций, вычислять их градиенты.

# Глава 2 Градиентные методы с неточным оракулом для задач выпуклой оптимизации

## 2.1 Стохастические градиентные методы с неточным оракулом

### 2.1.1 Введение

В 1960-е годы численные методы выпуклой оптимизации переживали свою первую большую революцию. В работах того времени четко и последовательно развивалась линия градиентных методов. Основополагающим здесь можно признать вклад Бориса Теодоровича Поляка [99, 100], с работ которого во многом и началось активное и повсеместное использование градиентных методов в Советском Союзе. Следующая революция началась в конце 1970-х годов после фундаментальных работ А.С. Немировского, Д.Б. Юдина, Л.Г. Хачияна, N. Karmarkar-а и др. [91, 120]. В монографии [91] была предложена классификация задач выпуклой (и не только) оптимизации по степени гладкости и выпуклости. Были получены нижние оценки для соответствующих классов задач оптимизации с оракулом, выдающим по запросу градиент или стохастический градиент, его компоненту или просто значение функции в точке. Стало понятно, чего в принципе можно достичь. Стали строиться оптимальные методы, см., например, [94, 104, 122]. При этом на задачи стали смотреть более пристально с точки зрения теории сложности. Появилась битовая сложность. Была показана полиномиальная разрешимость задач линейного программирования в битовой сложности [120]. Началась разработка полиномиальных методов внутренней точки для задач выпуклой оптимизации на базе метода Ньютона, которая впоследствии привела к созданию общей теории [93, 163, 259, 276] и соответствующего пакета CVX (*Open software* http://cvxr.com/cvx/), способного решать широкий спектр задач выпуклой оптимизации в пространствах размерности до $n \sim 10^4 - 10^5$. Однако, вызовы нового тысячелетия заставляют снова вернуться к градиентным методам. Задачи, которые стали возникать в последние десять лет, отличаются огромными размерностями $n \sim 10^6 - 10^9$. Такие задачи (классифицируемые как задачи large-scale и huge-scale оптимизации) приходят из анализа данных, поиска равновесий в различных сетевых моделях (связанных с компьютерными и транспортными сетями), биоинформатики и многих других областей. Для таких размерностей шаг (итерация) метода Ньютона,[1] становится слишком дорогим, поэтому приходится

---

[1] Заметим, что рассматривать методы более высокого порядка, чем метод Ньютона практически в любой ситуации не имеет смысла, поскольку зависимость числа итераций от желаемой точности решения задачи

снова возвращаться к более медленным (в смысле скорости сходимости), но более дешевым (в смысле стоимости одной итерации) градиентным методам (см. [99]). Но для указанных размерностей даже градиентные методы могут испытывать проблемы. В этой связи оказалась очень полезной концепция “заглядывания в черный ящик”, т.е. использование структуры задачи с целью ускорения вычислений [92], и использование вместо градиента его легко вычислимой (стохастической) аппроксимации [163]. Как следствие, принято стало считать, что правильный способ эффективно решать ту или иную задачу – это отказаться от общих методов, оптимальных на больших классах, и погружаться в специфику конкретной задачи в надежде ускориться и получить оценки лучше, чем нижние границы [91]. Можно сказать, что началась новая революция. Поток работ на эту тему в основных профилирующих журналах (например, Math. Program.) резко возрос (см., например, обзор [163]). Тем не менее, параллельно стали появляться работы, показывающие, что многие эффективные методы решения современных задач выпуклой оптимизации в пространствах огромных размеров получаются сочетанием небольшого количества приемов и идей. Цель настоящей раздела состоит в том, чтобы собрать воедино набор основных таких идей и показать их связь с некоторыми концепциями 1960-х годов, многие из которых восходят к Б.Т. Поляку. Мы сосредоточимся на оценках числа итераций, требующихся различным методам для решения задачи выпуклой оптимизации с заданной точностью по функции. Эта информация не в полной мере характеризует эффективность метода, но она необходима для последующего его полного исследования. Мы также ограничимся рассмотрением методов проекции градиента [21, 83], в которые, например, не входят очень популярные в последнее время методы условного градиента [163, 220, 263]. В качестве основного инструментария для получения эффективных методов используется аппарат оценивающих последовательностей, восходящий к работам Ю.Е. Нестерова [92, 93 ,94]. Здесь имеются и альтернативные подходы, например, [135, 283, 323]. Из-за ограничений на объем и большого количества технических деталей мы ограничимся в этом разделе лишь изложением общей картины. В частности, в этом разделе не приводится псевдокод соответствующих алгоритмов, но, как правило, указываются источники, в которых его можно найти. Более подробно об упоминаемых далее в этом разделе алгоритмах будет написано в следующих разделах и главах диссертации. Мы также не претендуем здесь на полный обзор современного состояния исследований, посвященных градиент-

---

невозможно сделать лучше, чем у метода Ньютона (в окрестности его квадратичной сходимости), сколько бы старших (равномерно ограниченных) производных функционала не использовалось в методе [91].

ным методам. Более того, при ссылках на литературу мы далеко не всегда ссылались на первоисточники, иногда предпочитая ссылаться на удачно написанный более доступный и более современный обзор или монографию.

### 2.1.2 Стохастическая оптимизация

Рассматривается задача выпуклой стохастической оптимизации [67, 96, 99]:

$$f(x) = E_{\xi}\left[f(x,\xi)\right] \to \min_{x \in Q} \tag{2.1.1}$$

где $f(x)$ – выпуклая по $x \in \mathbb{R}^n$ ($n \gg 1$) функция. Будем называть $\nabla f(x,\xi)$ стохастическим субградиентом функции $f(x,\xi)$ в точке $x$ по первой переменной [309]. Будем считать, что[2] п.н. $\left\|\nabla f(x,\xi)\right\|_2 \le M$, $\nabla = \nabla_x$ и $E_{\xi}$ – перестановочны.[3] Предположим, что $Q$ – выпуклое замкнутое ограниченное множество. Обозначим через $R$ – диаметр множества $Q$: $R = \max_{x,y \in Q}\left\|x - y\right\|_2$. В действительности, достаточно считать, что $R$ – расстояние от точки старта до решения (ближайшего, если решение не единственно) задачи (2.1.1) (см. замечание 2.1.1). При этом множество $Q$ может быть не ограничено [269].[4] Мы будем считать, что множество $Q$ простой структуры, т.е. на него можно эффективно проектиро-ваться. В работах [101, 225, 227, 269, 294] рассматривались различные варианты методов проекции градиента с усреднением и длинными шагами[5] применительно к решению задачи (2.1.1). Общая оценка скорости сходимости этих методов есть ($\sigma > 0$ – малый доверительный уровень, $N$ – число итераций метода, на каждой итерации мы можем один раз обратиться к оракулу за субградиентом)[6]

$$P_{x_N}\left( f\left(x_N\right) - \min_{x\in Q} f\left(x\right) \geq CMR\sqrt{\frac{1+\ln\left(\sigma^{-1}\right)}{N}} \right) =$$

$$= P_{x_N}\left( E_\xi\left[ f\left(x_N,\xi\right)\right] - \min_{x\in Q} E_\xi\left[ f\left(x,\xi\right)\right] \geq CMR\sqrt{\frac{1+\ln\left(\sigma^{-1}\right)}{N}} \right) \leq \sigma ,$$

где $C$ – константа (здесь и далее константы в основном будут в диапазоне $\sim 10^0 - 10^2$), а случайный вектор $x_N$ – то, что выдает алгоритм (например, метод зеркального спуска [225] или метод двойственных усреднений [269] – сравнительный анализ и описание “физики” этих методов в детерминированном случае проводится в работе [135]) после $N$ итераций. Мы будем называть $x_N$ – $\left(\varepsilon,\sigma\right)$-решением задачи (2.1.1), если

$$P_{x_N}\left( f\left(x_N\right) - \min_{x\in Q} f\left(x\right) \geq \varepsilon \right) \leq \sigma .$$

---

[2] В действительности [268], здесь и практически в любом другом контексте, где возникает такого типа условия, достаточно требовать, что п.н. выполнено неравенство $\left\|\nabla f(y,\xi) - \nabla f(x,\xi)\right\|_2 \le M$. Это позволяет в ряде случаев понизить оценку константы $M$ и как следствие (см. (2.1.2)), ускорить метод. Отметим, что под $\nabla f(x,\xi)$ (аналогично под $\nabla f(y,\xi)$) понимается любой элемент соответствующего стохастического субградиента [309].

[3] Для задач онлайн оптимизации условие перестановочности необходимо записывать в более общем (мартингальном) виде [54].

[4] Впрочем, в случае неограниченного множества $Q$ даже когда выпуклая функция $f(x)$ имеет ограниченную вариацию и равномерно ограниченную норму субградиента на $Q$ (рассматривается детерминированная постановка), мы не имеем никаких гарантий на скорость сходимости используемого метода (под любой метод можно подобрать такую функцию из описываемого класса, что сходимость будет сколь угодно медленной), поскольку не можем никак априорно ограничить расстояние $R$ [91]. Интересный нюанс для выпуклой (но не обязательно сильно выпуклой и гладкой) функции $f(x)$ имеет место, если $f(x)$ задана на ограниченном множестве $Q$ (рассматривается детерминированная постановка). В этом случае размер $R$ множества $Q$ может не входить в оценку необходимого числа итераций (входит вариация функции на этом множестве). Например, это имеет место для метода центра тяжести [91, 99, 163] и для задач, решаемых с относительной точностью [92].

Таким образом, для достижения точности по функции $\varepsilon$ и доверительного уровня $\sigma$ методу потребуется (здесь и далее мы будем использовать $\mathrm{O}\left(\cdot\right)$, однако все эти формулы могут быть переписаны с точными константами, что важно, поскольку во многих ситуациях такие оценки используются для формирования критерия останова метода)

$$\mathrm{O}\left( M^2R^2 \ln\left(\sigma^{-1}\right) \big/ \varepsilon^2 \right) \tag{2.1.2}$$

итераций. На каждой итерации вычисляется стохастический субградиент и осуществляется проектирование.

Отметим, что если использовать метод Монте-Карло, заключающийся в замене исходной задачи (2.1.1) следующей задачей

$$\frac{1}{N}\sum_{k=1}^{N} f\left(x,\xi_k\right) \to \min_{x\in Q} , \tag{2.1.3}$$

---

[5] Б.Т. Поляком было показано [101], что такое сочетание позволяет получать эффективные методы для данного класса задач.

[6] Эта оценка является неулучшаемой с точностью до мультипликативной константы $C$ (при $N \leq n$ оценка является неулучшаемой и в детерминированном случае $f\left(x,\xi\right) \equiv f\left(x\right)$), см. [91].

где с.в. $\xi_k$ – i.i.d., и распределены также как и $\xi$, то для того, чтобы гарантировать, что абсолютно точное решение этой новой задачи является $(\varepsilon,\sigma)$-решением исходной задачи потребуется взять $N$ порядка [309]

$$\mathrm{O}\left(M^2R^2\left(n\ln\left(MR/\varepsilon\right)+\ln\left(\sigma^{-1}\right)\right)/\varepsilon^2\right).$$

Эта наблюдение хорошо поясняет, что подход, связанный с усреднением случайности за счет самого метода более предпочтителен, чем замена задачи (2.1.1) ее стохастической аппроксимацией (2.1.3).[7] Более предпочтителен не только тем, что допускает адаптивность постановки и легко переносится на онлайн модификации исходной задачи, но, прежде всего, лучшей приспособленностью к большим размерностям.

Здесь важно подчеркнуть фундаментальную идею[8], которую можно усмотреть, например, в [99] и в цикле работ Б.Т. Поляка с Я.З. Цыпкиным [59, 105], о том, что для получения (агрегирования) хороших оценок неизвестных параметров (особенно когда размерность пространства параметров велика) имеет смысл рассматривать задачу поиска оптимальных значений параметра, как задачу стохастической оптимизации и рассматривать выборку как источник стохастических градиентов. Например, истинное значение неизвестного вектора параметров в предположении верности исходной параметрической гипотезы может быть записано как решение задачи стохастической оптимизации [70, 319] (метод наибольшего правдоподобия Фишера)

$$\theta^* = \arg\max_{\theta\in Q} E_{\xi}\left[L(\theta,\xi)\right],$$

[7] Особенно ярко это проявляется в случае, бесконечномерных пространств, возникающих в статистической теории обучения (СТО = SLT, Statistical Learning Theory) [320]. Попытка обучиться за счет минимизации эмпирического риска (а именно так можно расшифровать формулу (2.1.3) в СТО) может не дать состоятельной оценки/решающего правила, в то время как соответствующий стохастический зеркальный спуск дает состоятельную оценку. Отметим, что в работе [320] приводится достаточно интересный общий результат: в задачах обучения (в частности, в задачах СТО, математической статистики и онлайн обучения) способ получения оптимальных (с точностью до логарифмических факторов) оценок/решающих правил (или, другими словами, способ наискорейшего обучения) базируется на применении соответствующего метода зеркального спуска. Правда, найти “соответствующий метод”, в свою очередь, представляет собой непростую задачу.

[8] Распространяемую и на непараметрическую статистику. Отметим, что начиная с 1980-х годов XX века в этом направлении был цикл работ А.С. Немировского, Б.Т. Поляка и А.Б. Цыбакова, оказавших заметное влияние и на текущие исследования в этой области.

где $L(\theta,\xi)$ – логарифм функции правдоподобия. Однако решать эту задачу обычными методами мы не можем, потому что математическое ожидание берется по с.в. $\xi$, распределение которой задается неизвестным параметром $\theta^*$. Обойти эту сложность можно, если решать ту же самую задачу

$$E_\xi\left[-L(\theta,\xi)\right]\to\min_{\theta\in Q}$$

методами стохастической оптимизации, получая на каждом шаге новую реализацию (элемент выборки) $\xi_k$ и рассчитывая значения стохастического градиента $\partial L(\theta,\xi_k)/\partial\theta$. То, что выдает алгоритм и будет оценкой вектора неизвестных параметров $\theta^*$. Как правило, дополнительно известно, что $L(\theta,\xi)$ – гладкая и $\mu$-сильно вогнутая (равномерно по $\xi$) функция от $\theta$. Последнее обстоятельство позволяет получить лучшую оценку скорости сходимости по функции [218, 228, 236, 298] (в [236, 298] используются специальная модификация метода проекции градиента с усреднением и выбором шагов $h_k=2\left(\mu\cdot(k+1)\right)^{-1}$ и $h_k=(\mu k)^{-1}$, где $k$ – номер итерации, о подходе [218] и близком к нему подходе [228] будет немного написано в подразделе 2.1.3)

$$\mathrm{O}\left(M^2\ln\left(\ln(N)/\sigma\right)/(\mu N)\right), \tag{2.1.4}$$

т.е. ($x=\theta$, $f=-L$ , $\bar{C}$ – некоторая константа)

$$P_{x_N}\left(f(x_N)-\min_{x\in Q}f(x)\ge\bar{C}M^2\frac{\ln\left(\ln(N)/\sigma\right)}{\mu N}\right)\le\sigma.$$

Из неравенства Рао–Крамера [70] ($Q=\mathbb{R}^n$) будет следовать, что оценка (2.1.4) – не улучшаемая (с точностью до слагаемого $\ln\left(\ln(N)\right)$). Правда, тут возникают некоторые тонкости, когда мы говорим о неулучшаемости оценок с учетом вероятностей больших отклонений. Строго говоря, классические результаты типа Рао–Крамера, Ван-Трисса и т.п. (см., например, [70]) позволяют лишь говорить о неулучшаемости в смысле сходимости полных математических ожиданий (без вероятностей больших отклонений), и именно в таком смысле можно получить (с помощью методов [163, 218, 228, 298]) неулучшаемую (с точностью до мультипликативной константы) оценку:

$$E_{\xi,x_N}\left[f(x_N,\xi)\right]-\min_{x\in Q}E_\xi\left[f(x,\xi)\right]\le\frac{\breve{C}M^2}{\mu N},$$

где $\breve{C}$ – некоторая константа.

Можно обобщить рассмотренную постановку задачи (2.1.1) на случай, когда $\left\|\nabla f(x,\xi)\right\|_2$ имеет субгауссовский хвост (определение см., например, в [225]). Оценки (2.1.2), (2.1.4) при этом сохранят прежний вид (см. arXiv:1601.07592). Если же $\left\|\nabla f(x,\xi)\right\|_2^2$ имеет степенной хвост [17], т.е.

$$P\left(\frac{\left\|\nabla f(x,\xi)\right\|_2^2}{M^2} \ge t\right) = \mathrm{O}\left(\frac{1}{t^\alpha}\right),$$

где $\alpha > 2$, то[9]

$$P_{x_N}\left(f(x_N) - \min_{x\in Q} f(x) \ge \frac{C_\alpha MR}{\sqrt{N}}\left(\sqrt{\ln\left(\sigma^{-1}\right)} + \frac{(N/\sigma)^{1/\alpha}}{N}\right)\right) \le \sigma .$$

Если дополнительно $f(x) = E_\xi\left[f(x,\xi)\right]$ – $\mu$-сильно выпуклая функция, то (при $\alpha > 2$)

$$P_{x_N}\left(f(x_N) - \min_{x\in Q} f(x) \ge \bar{C}_\alpha \frac{M^2}{\mu N}\left(\ln\left(\ln\left(N\ln\left(\sigma^{-1}\right)\right)/\sigma\right) + \frac{\left(\ln\left(N/\sigma\right)/\sigma\right)^{2/\alpha}}{N^{2(1-1/\alpha)}}\right)\right) \le \sigma .$$

Если ничего не известно о $\left\|\nabla f(x,\xi)\right\|_2^2$, кроме неравенства $E_\xi\left[\left\|\nabla f(x,\xi)\right\|_2^2\right] \le M^2$, то по неравенству Маркова

$$P_{x_N}\left(f(x_N) - \min_{x\in Q} f(x) \ge \frac{\hat{C}MR}{\sigma\sqrt{N}}\right) \le \sigma,$$

$$P_{x_N}\left(f(x_N) - \min_{x\in Q} f(x) \ge \frac{\breve{C}M^2}{\sigma\mu N}\right) \le \sigma,$$

второе неравенство подразумевает $\mu$-сильную выпуклость $f(x)$.

Можно задать вопрос: насколько вообще уместно рассматривать постановки, в которых возникают тяжелые хвосты. Ведь, если мы можем эффективно вычислять значения функции $f(x) = E_\xi\left[f(x,\xi)\right]$ в задаче (2.1.1), то ни о каких тяжелых хвостах можно не заботиться. Поскольку, выбрав число шагов так, чтобы метод находил $\varepsilon$-решение с вероятностью $\ge 1/2$, запустив $\log_2\left(\sigma^{-1}\right)$ реализаций такого метода и выбрав реализацию с ми-нимальным значением функции в конечной точке алгоритма, мы за дополнительную $\log_2\left(\sigma^{-1}\right)$ плату (мультипликативную) получим с вероятностью $1-\sigma$ среди выданных ответов хотя бы одно $\varepsilon$-решение [218, 227]. Однако предположение о возможности эффективно вычислять значения функции (при условии трудной вычислимости ее градиента), как правило, не встречается на практике. В некотором смысле типичным тут является пример 2.1.1 (см. ниже) вычисления вектора PageRank (при $n \sim 10^9$). Собственно, искусственность ситуации, в которой значение функции легко вычислимо, а градиент нет, неплохо соответствует философии быстрого автоматического дифференцирования (БАД) [65, 76]. Согласно теории БАД, если мы можем посчитать значение функции, то мы можем не более чем в 4 раза дороже посчитать и ее градиент.[10] Как следствие, если мы можем эффективно вычислить значение $f\left(x\right)$, то, как правило, мы и $\nabla f\left(x\right)$ можем эффективно вычислить. Тогда и на исходную задачу (2.1.1) можно смотреть уже не как на задачу стохастической оптимизации, а как на обычную задачу выпуклой оптимизации, что может существенно ускорить ее решение (см. подраздел 2.1.3 ниже). Впрочем, во многих интересных приложениях отмеченный прием (амплификация), как правило, весьма успешно работает [3, 45], поскольку время работы метода, как правило, оказывается заметно большим, чем расчет значения функции.

---

[9] Приводимые ниже неравенства стоит понимать так, что $x_N$ выдается методом [225, 269], а в сильно выпуклом случае, методом [218, 228, 298]. При этом для оценок вероятностей больших уклонений в случае тяжелых хвостов требуется некоторые оговорки и уточнения. К сожалению, мы не смогли найти соответствующий выписанным оценкам (в случае тяжелых хвостов) источник литературы.

Отметим также (следуя А.С. Немировскому), что с помощью концепции неточного оракула (см. подраздел 2.1.3 ниже) мы можем редуцировать задачу с тяжелыми хвостами $\left\|\nabla f\left(x,\xi\right)\right\|_2^2$ и компактным множеством $Q$ к ситуации, когда п.н. $\left\|\nabla f\left(x,\xi\right)\right\|_2 \le M\left(\varepsilon\right)$. Для этого нужно “обрезать” стохастический градиент

$$\nabla f\left(x,\xi\right) := \begin{cases} \nabla f\left(x,\xi\right), & \left\|\nabla f\left(x,\xi\right)\right\|_2 \le M\left(\varepsilon\right) \\ M\left(\varepsilon\right)\dfrac{\nabla f\left(x,\xi\right)}{\left\|\nabla f\left(x,\xi\right)\right\|_2}, & \left\|\nabla f\left(x,\xi\right)\right\|_2 > M\left(\varepsilon\right) \end{cases}.$$

Константа $M\left(\varepsilon\right)$ подбирается оптимальным образом, исходя из желаемой точности $\varepsilon$. Чем больше $M\left(\varepsilon\right)$, тем меньше смещение (bias) обрезанного стохастического градиента, как следствие, тем точнее можно восстановить решение исходной задачи, но при этом возрастает необходимое число итераций (см. (2.1.2), (2.1.4), в которые входит константа $M = M(\varepsilon)$). Оптимальный выбор этой константы (с точностью до логарифмического фактора) дает приведенные выше оценки.

---

[10] Это легко понять в случае $f\left(x\right)=\left\langle c,x\right\rangle$. В случае, когда $f\left(x\right)$ – многочлен, это также несложно понять (Баур–Штрассен). В общем случае рассуждения аналогичны.

Все сказанное выше[11] обобщается и на другие прокс-структуры [91] (не обязательно евклидовы, когда выбирается прокс-функция $d(x) = \|x\|_2^2/2$), согласно которым осуществляется (как правило, по явным формулам[12]) "проектирование" на $Q$ (причины появления прокс-функции также объясняются в следующем разделе). Например, для множества $Q = S_n(1)$ ($S_n(R) = \left\{ x \in \mathbb{R}^n : x_i \geq 0,\ i = 1,...,n,\ \sum_{i=1}^{n} x_i = R \right\}$ – единичный симплекс в $n$-мерном пространстве) часто рассматривается (см. пример 2.1.1 вычисления вектора PageRank ниже) KL-прокс-структура: $d(x) = \ln n + \sum_{i=1}^{n} x_i \ln x_i$. Эта прокс-функция $d(x) \geq 0$ сильно выпукла в 1-норме с константой сильной выпуклости $\alpha = 1$ на $S_n(1)$ – в силу неравенства Пинскера [225, 269]. Она "наилучшим" образом подходит для симплекса (с некоторыми оговорками [54, 166]). Выгода от ее использования в том, что норма стохастического субградиента всегда оценивается в сопряженном пространстве к пространству, в котором прокс-функция 1-сильно выпукла. В рассматриваемом случае получается $\|\nabla f(x,\xi)\|_\infty \leq M$, что в типичных ситуациях дает оценку константы $M$ в $\sim \sqrt{n}$ раз лучше, чем в 2-норме, а плата за это – увеличение оценки размера области (в этой ситуации в оценке числа итераций нужно использовать $R^2 = \max_{x \in Q} d(x)/\alpha$ ) в $\sim \ln n$ раз. Детали имеются, например, в статье [225]. Интересным также представляется выбор прокс-функции для прямого произведения симплексов [51]. Здесь мы отметим (следуя А.С. Немировскому), что в общем случае оптимальный выбор прокс-структуры (с точностью до умножения на степень логарифма размерности пространства) связан с симметризацией множества $Q$. Выпуклое центрально симметричное множество $B = (Q - Q)/2$ порождает по теореме Колмогорова норму, в которой $B$ является единичным шаром. Далее ищется оптимальная прокс-функция, согласованная с этой нормой. Говоря более формально, ищется такая сильно выпуклая в этой норме функция $d(x) \geq 0$ с константой сильной выпуклости $\alpha \geq 1$, чтобы число $R^2 = \max_{x \in Q} d(x)/\alpha$ было минимально возможным. Если $Q = B_2^n(1)$ – единичный евклидов шар, то значение $R^2 \leq 1$, т.е. не зависит от размерности пространства $n$, но если $Q = B_\infty^n(1)$ – единичный шар в $l_\infty^n$ норме, то $R^2 = \Omega(n)$ (т.е. существует такое число $\chi$, что при достаточно больших значениях $n$ имеет место неравенство $R^2 \geq \chi n$, причем можно добиться того, что $R^2 = O(n)$). Как будет видно из замечания 2.1.2 (на примере когда $Q = B_\infty^n(1)$), выбор $l_\infty^n$ нормы не всегда приводит к оптимальным во всех смыслах оценкам (аналогичные примеры нам встретятся и в следующих двух пунктах).

---

[11] В сильно выпуклом случае (если в прямом пространстве выбрана $q$-норма ($l_q^n$) и прокс-функция $d(x) \geq 0$, $d(x_0) = 0$) в оценку (4) дополнительно входит фактор $\omega = \sup_{x \in Q} 2V(x,x_0) \big/ \left(\alpha \|x - x_0\|_q^2\right) \geq 1$, где $V(x,x_0)$ определяется через $d(x)$ в замечании 2.1.1 (при $1 \leq q \leq 2$ удается найти такую прокс-функцию, что $\omega = O(\ln n)$, см. замечание 2), где $\alpha$ – константа сильной выпуклости $d(x)$ на $Q$ в $q$-норме [228]. Отметим, что при этом константы в отношение $M^2/\mu$ в оценке (2.1.4), считаются относительно $q$-нормы.

[12] Впрочем, в подавляющем большинстве случаев даже если нет возможности явно решить задачу проектирования, ее можно эффективно решить приближенно [259] (посредством перехода к двойственной задаче малой размерности). Как правило, при таком способе рассуждений необходимо использовать концепцию неточного оракула (см. подраздел 2.1.3), поскольку рассчитать градиент двойственного функционала можно лишь приближенно. Однако все эти выкладки обычно не изменяют по порядку сложность одной итерации метода, основной составляющей которой является расчет (пересчет) градиента или его стохастического аналога. Некоторые тонкости и оговорки тут возникают в случае разреженных постановок задач [3, 45].

**Замечание 2.1.1.** Стоит обратить внимание на то, что если выбрана евклидова прокс-структура, то $R^2$ – квадрат евклидова диаметра $Q$. При переходе к другой прокс-структуре в оценках числа итераций в качестве $R^2$ фигурирует прокс-диаметр $Q$ ($\operatorname{diam}(Q) = \max_{x \in Q} d(x)$), поделенный на константу сильной выпуклости $\alpha = \alpha(Q)$ прокс-функции, заданной на $Q$, относительно выбранной нормы в прямом пространстве. Скажем, в случае выбора KL-прокс-структуры, 1-нормы в прямом пространстве и $Q = S_n(r)$, имеем

$$R^2 = \operatorname{diam}(S_n(r))/\alpha(S_n(r)) = r \cdot \operatorname{diam}(S_n(1))/(\alpha(S_n(1))/r) = r^2 \cdot (\ln n)/1 = r^2 \ln n.$$

Для евклидовой прокс-структуры размер $Q = S_n(r)$ равнялся бы $2r^2$. Отсюда можно сделать вывод (верный и в общем случае), что выбор прокс-структуры имеет целью оптимально учесть структуру множества с точки зрения того как в итоговую оценку числа итераций будет входить размерность пространства, в котором происходит оптимизация. При гомотетичном увеличении/уменьшении множества оценки числа итераций будут меняться одинаково, независимо от выбранной прокс-структуры. Отметим также, что в формуле (2) для прокс-структуры, отличной от евклидовой точнее писать не $R^2\ln\left(\sigma^{-1}\right)$, где $R^2 = r^2\ln n$ (приводим для KL-прокс-структуры), а $r^2\left(\ln n + \ln\left(\sigma^{-1}\right)\right) = r^2\ln\left(n/\sigma\right)$. В действительности, в оценки скоростей сходимости (в среднем, но не в оценки вероятностей больших уклонений, см. замечание 2.1.4) всех упомянутых в данном разделе методов (кроме обычного (прямого) градиентного метода и метода Франк–Вульфа) входит не прокс-диаметр множества $Q$, на котором происходит оптимизация (если $Q = \mathbb{R}^n$ прокс-диаметр будет бесконечным), а брэгмановское “расстояние” $V\left(x_*, x_0\right)$ от решения $x_*$ до точки старта $x_0$ (часто выбирают $x_0 = \arg\min_{x\in Q} d\left(x\right)$, $d\left(x_0\right) = 0$, $\nabla d\left(x_0\right) = 0$ [269]), где $V\left(x,y\right) = d\left(x\right) - d\left(y\right) - \left\langle \nabla d\left(y\right), x - y\right\rangle$.

**Замечание 2.1.2.** Пусть $Q = B_q^n\left(1\right)$ – единичный шар в $q$-норме или, в более общем случае, $Q$ содержится в $B_q^n\left(1\right)$. Относительно оптимального выбора нормы и прокс-структуры можно заметить следующее (см., например, [91, 259, 128]): если $q \ge 2$, то в качестве нормы оптимально выбирать $\|\ \|_2$ (2-норму) и евклидову прокс-структуру. Определим $q'$ из $1/q + 1/q' = 1$. Пусть $1 \le q \le 2$, тогда $q' \ge 2$. Если при этом $q' = o\left(\log n\right)$, то оптимально выбирать $\|\ \| = \|\ \|_q$, а прокс-структуру задавать прокс-функцией $d\left(x\right) = \frac{1}{2\left(q-1\right)}\|x\|_q^2$. Во всех этих случаях $R^2 = \mathrm{O}\left(1\right)$. Для $q' \ge \Omega\left(\log n\right)$, выберем $a = 2\log n/\left(2\log n - 1\right)$, $\|\ \| = \|\ \|_a$, а прокс-структуру будем задавать прокс-функцией $d\left(x\right) = \frac{1}{2\left(a-1\right)}\|x\|_a^2$. В этом случае $R^2 = \mathrm{O}\left(\log n\right)$. Не сложно проверить, что для единичного симплекса, вложимого в единичный шар в 1-норме, выбор соответствующих прокс-структур из замечаний 2.1.1, 2.1.2 приводит к одинаковым оценкам числа итераций в категориях $\mathrm{O}(\ )$. В частности, для случая когда $Q = B_\infty^n\left(1\right)$, выбор 2-нормы и евклидовой прокс-структуры приводит к оценке (далее в замечании речь идет только об оценке (2.1.2)) 1) $\mathrm{O}\left(M_2^2 n\ln\left(\sigma^{-1}\right)/\varepsilon^2\right)$ вместо 2) $\mathrm{O}\left(M_\infty^2 n\ln\left(\sigma^{-1}\right)/\varepsilon^2\right)$ (здесь $E_\xi\left[\left\|\nabla f\left(x,\xi\right)\right\|_1^2\right] \le M_\infty^2$), получаемой при выборе $l_\infty^n$ нормы в прямом пространстве. Аналогично вышенаписанному можно отметить, что в типичных ситуациях оценка 2 может быть в $\sim n$ раз хуже оценки 1.

Тем не менее, оценка 2 $O\left(M_\infty^2 n \ln\left(\sigma^{-1}\right)/\varepsilon^2\right)$ не улучшаема в общем случае.[13] Потому что в общем случае нет гарантий, что $M_2^2 \ll M_\infty^2$, а если такие гарантии есть, то это уже сужает класс функций, для которого получена нижняя оценка с константой $M_\infty^2$.

Подчеркнем, что приведенные здесь оценки (2.1.2), (2.1.4) (в детерминированном случае при дополнительном условии, что требуемое число итераций для достижения точности $\varepsilon$ удовлетворяет неравенству $N(\varepsilon) \le n$ [91]) без дополнительных предположений являются неулучшаемыми (с точностью до мультипликативных констант) для класса задач стохастической оптимизации (2.1.1) и негладких детерминированных задач. Причем дополнительная гладкость функционала задачи (2.1.1) в стохастической постановке в общем случае не приводит к улучшению приведенных оценок (2.1.2), (2.1.4). Если делать дополнительные предположения о малости случайного шума (low noise conditions), то приведенные оценки можно улучшать (см. подраздел 2.1.3). Один пример того, как можно устанавливать неулучшаемость оценок был рассмотрен выше, следуя [99] (на основе неравенства Рао–Крамера), в общем случае следует смотреть монографию [91] и [128]. Отметим, что в работе [128] показывается, что для задач стохастической оптимизации (2.1.1) при оптимизации на шарах в $q$-норме оценки типа (2.1.2), даваемые методами зеркального спуска с выбором прокс-структуры согласно замечанию 2.1.2, соответствуют с точностью до логарифмического фактора нижним оценкам.

Следует, однако, различать задачи стохастической оптимизации и задачи, в которые мы сами искусственно привносим случайность (используя рандомизацию) с целью уменьшения числа арифметических операций на одну итерацию метода [3, 5, 45, 227]. К последнему можно отнести случай, когда (негладкий) выпуклый функционал в задаче является детерминированным, но представляет собой трудно вычислимый интеграл (сумму), зависящую от (оптимизируемых) параметров, который может быть компактно представлен в виде математического ожидания по некоторой простой вероятностной мере. Тогда выгоднее вычислять на каждой итерации метода стохастический градиент, существенно экономя на вычислениях на каждом шаге и лишь немного теряя на логарифмическом увеличении числа шагов ($\sim \ln\left(\sigma^{-1}\right)$). Подробнее об этом подходе будет сказано ниже в примере 2.1.3. Ярким примером на эту тему является Google problem (PageRank), которая подробно изучается в главе 4 (см. также приложение в конце диссертации). По-видимому, одними из первых на эту задачу посмотрели в указанном выше контексте А.В. Назин и Б.Т. Поляк в работе [88], см. также [3, 46, 264, 273].

---

[13] Общий результат здесь такой [91, 128]. Пусть необходимо найти минимум выпуклой функции $f(x)$ на множестве $Q = B_q^n(R)$. Оракул выдает несмещенные стохастические субградиенты со свойством $E_\xi\left[\left\|\nabla f(x,\xi)\right\|_{q'}^2\right] \le M_q^2$ $(1/q + 1/q' = 1)$. Тогда для того, чтобы найти такую точку $x^N$, что $E\left[f\left(x^N\right)\right] - \min_{x\in Q} f(x) \le \varepsilon$, необходимо обратиться к оракулу не менее $N = c_q M_q^2 R^2 / \varepsilon^{\max(2,q)}$ раз, при $N \ll n$, где $c_q = O(\ln n)$ (эта оценка $c_q$ становится точной при $q \to 1+0$); и не менее $N = c_q M_q^2 R^2 n^{1-2/\max(2,q)} / \varepsilon^2$ раз, при $N \gg n$. В детерминированном случае (когда оракул выдает субградиент $\left\|\nabla f(x)\right\|_{q'} \le M_q$) последняя оценка примет вид $N = cn\ln\left(M_q R/\varepsilon\right)$, при $N \gg n$.

**Пример 2.1.1 (PageRank).** Задача поиска вектора PageRank $p$ из уравнения $P^T p = p$ ($P$ – стохастическая матрица по строкам матрица), сводится [46, 264] к негладкой задаче выпуклой оптимизации (седловой задаче)

$$\max_{u \in S_n(1)} \left\langle u, P^T p - p \right\rangle \to \min_{p \in S_n(1)}.$$

Перепишем эту задачу в общем виде

$$\min_{x \in S_n(1)} \max_{y \in S_n(1)} \left\langle y, Ax \right\rangle,$$

где матрица $A$ большого размера $n \times n$ (вообще говоря, неразреженная) с элементами, ограниченными по модулю числом $M = 1$. Ключевое наблюдение для решения этой задачи состоит в том [87, 225], что:

$$Ax = E_{i[x]}\left[ A^{\langle i[x] \rangle} \right],$$

где $A^{\langle i \rangle}$ – $i$-й столбец матрицы $A$, вектор $x \in S_n(1)$, а с.в. $i[x]$ имеет категориальное распределение с вектором параметров $x$. Важным следствием является тот факт, что левая часть равенства, $Ax$, вычисляется за $O\left(n^2\right)$ арифметических операций, а выражение, стоящее в правой части под математическим ожиданием, $A^{\langle i[x] \rangle}$ – всего лишь за $O(n)$ арифметических операций. Используя это наблюдение (и аналогичное для умножения матрицы $A$ на вектор-строку слева), можно показать, что (рандомизированный) метод зеркального спуска [225] (с KL-прокс-структурой) и стохастическим градиентом по $x$ равным $A^{\langle i[x] \rangle}$ (аналогично по $y$) после

$$\mathrm{O}\left(\frac{nM^2\ln(n/\sigma)}{\varepsilon^2}\right)=\mathrm{O}\left(\frac{n\ln(n/\sigma)}{\varepsilon^2}\right)$$

элементарных арифметических операций выдает такие $x\in S_n(1)$ и $y\in S_m(1)$, что

$$\max_{\tilde{y}\in S_m(1)}\tilde{y}^T Ax-\min_{\tilde{x}\in S_n(1)}y^T A\tilde{x}\le\varepsilon$$

с вероятностью $\ge 1-\sigma$.

Аналогичные рассуждения [5, 227] позволяют получить с такими же затратами $\mathrm{O}\left(n\ln(n/\sigma)\varepsilon^{-2}\right)$ такой вектор $x\in S_n(1)$, что $\|Ax\|_\infty\le\varepsilon$. Кроме того, если дополнительно известно, что матрица $P$ – разрежена, то можно организовать поиск $(\varepsilon,\sigma)$-решения еще эффективнее – рандомизировать при проектировании на симплекс [5, 46, 54]. Тогда вместо фактора $n$ в оценках общего числа операций $\mathrm{O}\left(n+s\ln n\ln(n/\sigma)\varepsilon^{-2}\right)$ будет фигурировать $s$ – “среднее” число элементов матрицы $P$ (по строкам и столбцам) отличных от нуля (к сожалению, численные эксперименты Антона Аникина показали, что это “эффективное среднее” число на практике часто близко к максимальному по строкам и столбцам, т.е. от этого подхода можно получить гарантированную выгоду, только если имеет место равномерная разреженность матрицы по строкам и столбцам [5]).

Отметим, что в определенных ситуациях (например, при условии $n\gg\varepsilon^{-2}$ – типичном для задач huge-scale оптимизации) такому рандомизированному методу потребуется использовать относительно небольшое количество элементов матрицы $A$ за все время работы, в то время как для класса детерминированных алгоритмов потребуется считать как минимум половину элементов матрицы $A$ [120] для $\varepsilon=0.1$.

Хочется также отметить, что на задаче из примера 2.1.1 можно продемонстрировать большую часть современного инструментария, необходимого для решения задач huge-scale оптимизации. Так, в случае разреженной матрицы $A$ для решения поставленной негладкой задачи выпуклой оптимизации (и многих других) хорошо подходит метод Б.Т. Поляка [99, 273], работающий по нижним оценкам (2.1.2) (функционал негладкий) и при этом учитывающий разреженность $A$ при пересчете градиента [273]. Другой подход [264] (задача поиска вектора PageRank сводится к минимизации другого функционала), также нашедший широкое применение [45, 138, 199, 295], связан с заменой градиентного спуска на покомпонентный спуск. Такая замена увеличивает в среднем число итераций всегда не больше (а, как правило, намного меньше) чем в $n$ раз, но зато (благодаря разреженности) происходит экономия при пересчете одной компоненты градиента, как правило (но не всегда – особенности возникают в разреженных задачах), в $n$ раз по сравнению с расчетом полного градиента. В результате получается выгода, которая при определенных условиях может сократить объем вычислений в $\sim\sqrt{n}$ раз (см., например, [45]). Поясним это следующим примером [5], который можно понимать как вариацию неускоренного варианта покомпонентного метода с выбором максимальной компоненты [264].

**Пример 2.1.2 (разреженный PageRank по Ю.Е. Нестерову).** Задача поиска вектора PageRank также может быть сведена к следующей задаче выпуклой оптимизации (далее для определенности будем полагать $\gamma = 1$, в действительности, по этому параметру требуется прогонка)

$$f(x) = \frac{1}{2}\|Ax\|_2^2 + \frac{\gamma}{2}\sum_{k=1}^{n}(-x_k)_+^2 \to \min_{\langle x,e\rangle=1}$$

где как и в примере 1 $A = P^T - I$, $I$ – единичная матрица, $e = (1,\ldots,1)^T$,

$$(y)_+ = \begin{cases} y, & y \ge 0 \\ 0, & y < 0 \end{cases}.$$

При этом мы считаем, что в каждом столбце и каждой строке матрицы $P$ не более $s \ll \sqrt{n}$ элементов отлично от нуля ($P$ – разрежена). Эту задачу можно решать обычным градиентным методом[14], но не в евклидовой норме, а в 1-норме (см., например, [135]):

$$x_{k+1} = x_k + \operatorname*{arg\,min}_{h:\langle h,e\rangle=0}\left\{ f(x_k) + \langle \nabla f(x_k), h\rangle + \frac{L}{2}\|h\|_1^2 \right\},$$

где $L = \max_{i=1,\ldots,n}\left\|A^{\langle i\rangle}\right\|_2^2 + \gamma \le 3$ ($A^{\langle i\rangle}$ – $i$-й столбец матрицы $A$). Для достижения точности $\varepsilon^2$ по функции потребуется сделать $O\left(LR^2/\varepsilon^2\right) = O\left(1/\varepsilon^2\right)$ итераций [92]. Не сложно проверить, что пересчет градиента на каждой итерации заключается в умножении $A^T A h$, что может быть сделано за $O\left(s^2 \ln n\right)$. Связано это с тем, что вектор $h$ всегда имеет только две компоненты

$$\frac{1}{8}\left(\max_{i=1,\ldots,n}\partial f(x_k)/\partial x^i - \min_{i=1,\ldots,n}\partial f(x_k)/\partial x^i\right) \text{ и } -\frac{1}{8}\left(\max_{i=1,\ldots,n}\partial f(x_k)/\partial x^i - \min_{i=1,\ldots,n}\partial f(x_k)/\partial x^i\right)$$

---

[14] Выписанная далее оценка скорости сходимости (на число итераций) – неулучшаема с точностью до мултипликативного фактора. Речь идет не об оптимальности метода на классе гладких задач на симплексе, а о том, что конкретно для этого метода такая оценка если и может быть улучшена, то лишь на мультипликативный фактор. Это замечание касается практически всех известных сейчас градиентных методов. Показывается это приблизительно также (даже еще проще), как и в случае оптимальности оценок на классах [91]: строится конкретные примеры (семейства) функций.

отличные от нуля (такая разреженность получилась благодаря выбору 1-нормы), причем эти компоненты определяются, соответственно, как

$$\arg\min_{i=1,\dots,n} \partial f\left(x_k\right)/\partial x^i \text{ и } \arg\max_{i=1,\dots,n} \partial f\left(x_k\right)/\partial x^i ,$$

что пересчитывается (при использовании специального двоичного дерева (кучи) для поддержания максимальной и минимальной компоненты градиента [273]) за $O\left(s^2 \ln n\right)$ (Ю.В. Максимовым было замечено, что логарифмический фактор можно ослабить, если использовать, например, фибоначчиевы или бродалевы кучи [5]). Таким образом, общая трудоемкость предложенного метода будет $O\left(n + s^2 \ln n/\varepsilon^2\right)$, что заметно лучше многих известных методов [46]. Стоит также отметить, что функционал, выбранный в этом примере, обеспечивает намного лучшую оценку $\|Ax\|_2 \le \varepsilon$ по сравнению с функционалом из примера 2.1.1, который (в варианте [227]) обеспечивает $\|Ax\|_\infty \le \varepsilon$. Наилучшая (в разреженном случае без, условий на спектральную щель матрицы $P$ [46]) из известных нам на данный момент оценок $O\left(s \ln n \ln\left(n/\sigma\right)/\varepsilon^2\right)$ [46, 54] для $\|Ax\|_\infty$ может быть улучшена приведенной в этом примере оценкой, поскольку, как уже отмечалось ранее, $\|Ax\|_2$ может быть (и так часто бывает) в $\sim \sqrt{n}$ раз больше $\|Ax\|_\infty$, а $s \ll \sqrt{n}$.

Заметим, что в решении могут быть маленькие отрицательные компоненты. Также численные эксперименты показали [5], что для достижения выписанных оценок требуется препроцессинг (в нашем случае он заключается в представлении матрицы по строкам в виде списка смежности: в каждой строке отличный от нуля элемент хранит ссылку на следующий отличный от нуля элемент, аналогичное представление матрицы делается и по столбцам). Заметим, что препроцессинг помогает ускорять решение задач не только в связи с более полным учетом разреженности постановки, но и, например, в связи с более эффективной организацией рандомизации [45, 227, 264].

Пример 2.1.2 также характерным образом демонстрирует, как используется разреженность (см. также [2, 3, 273]). Обратим внимание на то, что число элементов в матрице $P$, отличных от нуля, даже при наложенном условии разреженности (по строкам и столбцам), все равно может быть достаточно большим $sn$. Удивляет то, что в оценке общей трудоемкости это число не присутствует. Это в перспективе (при правильной организации работы с памятью) позволяет решать задачи огромных размеров. Более того, даже в случае небольшого числа не разреженных ограничений вида $\langle a_i, x\rangle = b_i$, $i = 1,..,m = O(1)$, можно “раздуть” пространство (не более чем в два раза), в котором происходит оптимизация (во многих методах, которые учитывают разреженность такое раздутие не приведет к серьезным затратам), и переписать эту систему в виде $Ax = b$, где матрица будет иметь размеры $\mathrm{O}(n)\times\mathrm{O}(n)$, но число отличных от нуля элементов в каждой строке и столбце будет $\mathrm{O}(1)$. Таким образом, допускается небольшое число “плотных” ограничений.

Заметим, что если применить метод условного градиента [220] (Франк–Вульфа) к задаче из примера 2.1.2, то общая трудоемкость (для точности $\varepsilon^2$, как и в примере 2.1.2) будет [2, 5]

$$\mathrm{O}\left(n+\frac{s^2\ln\left(2+n/s^2\right)}{\varepsilon^2}\right).$$

В связи со сказанным выше, заметим, что задача может быть не разрежена, но свойство разреженности появляется в решении при использовании метода Франк–Вульфа, что также может заметно сокращать объем вычислений в постановках аналогичных примеру 2.1.2, но с матрицами $A$, у которой число столбцов на много порядков больше числа строк (см., например, п. 3.3 [163], [47]).

Приведем еще один пример, подсказывающий, как следует решать задачу (2.1.3), полученную из (2.1.1) с применением идеи метода Монте-Карло [67].

**Пример 2.1.3 (рандомизация суммы).** Пусть необходимо решить задачу выпуклой оптимизации (или ее композитный вариант, см., например, замечание 2.1.6)

$$f(x)=\frac{1}{N}\sum_{k=1}^{N}f_k(x)\to\min_{x\in Q}, \tag{2.1.5}$$

где $f_k(x)$ – негладкие выпуклые функции с ограниченной числом $M$ нормой субградиента, $Q$ – выпуклое замкнутое множество простой структуры (можем эффективно на него проектироваться, согласно заданной прокс-функции) прокс-диаметра $R$. Введем новую функцию

$$f(x,\xi)=\begin{cases} f_1(x), \text{с вероятностью } 1/N \\ \dots\dots\dots\dots\dots\dots\dots\dots \\ f_N(x), \text{с вероятностью } 1/N \end{cases}.$$

Ее стохастический субградиент легко вычислить. Для этого разыгрывается за $\mathrm{O}(\ln N)$ с.в. $\xi$, принимающая значения $1,\dots,N$ с равными вероятностями (см., например, [46]). Затем считается субградиент $f_\xi(x)$ (и выполняется прокс-проектирование на $Q$). Как уже отмечалось ранее, можно найти $(\varepsilon,\sigma)$-решение так понимаемой задачи (2.1.5) за

$$\mathrm{O}\left(\frac{M^2R^2\ln\left(\sigma^{-1}\right)}{\varepsilon^2}\right)$$

итераций, со стоимостью одной итерации равной $\mathrm{O}(\ln N)$ + затраты на вычисления субградиента $f_\xi(x)$ + затраты на вычисление проекции. Если решать задачу без рандомизации, то число итераций будет $\mathrm{O}\left(M^2R^2/\varepsilon^2\right)$, строго говоря, здесь $M$ должно быть немного меньше за счет того, что

$$\max_{x\in Q}\left\|\frac{1}{N}\sum_{k=1}^{N}\nabla f_k(x)\right\|_* \le \max_{\substack{k=1,\dots,N\\ x\in Q}}\left\|\nabla f_k(x)\right\|_*,$$

но мы считаем, что обе части неравенства одного порядка. Зато шаг итерации будет теперь почти в $N$ раз дороже. И если $N \gg 1$ это может оказаться существенным.

Приведенную постановку можно распространить на случай, когда взвешивание функций не равномерное (тогда первое разыгрывание с.в. $\xi$, имеющей категориальное распределение, или приготовление процедуры рандомизации займет $\mathrm{O}(N)$, а все последующие $\mathrm{O}(\ln N)$) и $f_k(x) := E_{\xi_k}\left[f_k(x,\xi_k)\right]$ с равномерно ограниченными (по $k$, $x$ и $\xi$) нормами стохастических субградиентов. При этом все приведенные оценки числа итераций сохранятся. Причем требование равномерной ограниченности норм стохастических субградиентов можно существенно ослабить за небольшую плату (см. выше).

Если на решение задачи (2.1.3) теперь посмотреть в контексте описанной рандомизации с $f_k(x) = f(x,\xi_k)$ (здесь $\xi_k$ – не случайная величина, а полученная в методе Монте-Карло $k$-я по порядку реализация с.в. $\xi$), то “все встанет на свои места” в смысле одинаковости (с точностью до логарифмического фактора) двух подходов к решению задачи (2.1.1), описанных в начале пункта.

Описанная рандомизация при вычислении субградиента суммы функций, по-видимому, была одной из первых, которые предлагались в стохастической оптимизации [67]. Однако она популярна и по сей день, например, в связи с приложениями к поиску равновесий в транспортных сетях [33, 40, 42, 43, 47] и анализу данных (см., например, работы P. Richtarik-а, S. Shalev-Shwartz-а, T. Zhang-а и др.). В частности, в [33, 45, 163, 223, 232, 234, 239, 242, 248] в предположении, что все функции в (2.1.5) гладкие с константой

Липшица градиента $L$, предложен специальный рандомизированный метод (на базе описанного выше способа рандомизации суммы), в котором число вычислений градиентов слагаемых[15]

$$O\left(\left(N+\min\left\{LR^2/\varepsilon,\sqrt{NLR^2/\varepsilon}\right\}\right)\left(\ln\left(\Delta f/\varepsilon\right)+\ln\left(\sigma^{-1}\right)\right)\right),$$

где $\Delta f$ разность значения функции в стартовой точке и в минимуме. Эта оценка с точностью до логарифмического множителя соответствует нижней оценке в классе детерминированных алгоритмов [129, 239]. Если дополнительно имеется еще и $\mu$-сильная выпуклость $f(x)$, то оценку можно переписать следующим образом

$$O\left(\left(N+\min\left\{L/\mu,\sqrt{NL/\mu}\right\}\right)\left(\ln\left(\Delta f/\varepsilon\right)+\ln\left(\sigma^{-1}\right)\right)\right).$$

Отметим, что вторая оценка переходит в первую при следующей квадратичной регуляризации. К выпуклому функционалу прибавляется регуляризирующее слагаемое $\mu\|x\|_2^2/2$. В результате функционал становится сильно выпуклым и справедлива вторая оценка на число вычислений градиента. Такая регуляризация изменяет исходную целевую функцию на число не больше $\mu R^2/2$ и чтобы итоговая погрешность по исходной функции была порядка $\varepsilon$ нужно выбирать $\mu\simeq\varepsilon/R^2$, и решать регуляризованную задачу с точностью $\varepsilon/2$. При подстановке этого значения во вторую оценку числа вычислений градиента последняя переходит в первую оценку. Подробнее об этой конструкции написано в следующем разделе.

Отметим также, что сначала (см., например, [138, 232]) получается результат о сходимости средних[16]

$$E\left(f\left(x_N\right)-\min_{x\in Q}f\left(x\right)\right)\le\varepsilon,$$

где

---

[15] Строго говоря, имеющиеся сейчас рассуждения для второго аргумента минимума [3, 307] позволяют получить только при дополнительных предположениях о структуре задачи оценку, аналогичную приведенной ниже, и то только в категориях общего числа арифметических операций. В разделе 5.1 главы 5 мы вернемся к обсуждению этой задачей. Впрочем, в недавнем цикле работ Z. Allen-Zhu удалось не только обосновать, но даже немного улучшить приведенную оценку

http://arxiv.org/pdf/1603.05953.pdf, http://arxiv.org/pdf/1603.05643.pdf.

[16] Описанная далее конструкция не зависит от того, изначально имела место сильная выпуклость или мы ее искусственно ввели должной регуляризацией.

$$N = N(\varepsilon) = \mathrm{O}\left(\left(N + \min\left\{L/\mu, \sqrt{NL/\mu}\right\}\right)\ln\left(\Delta f/\varepsilon\right)\right),$$

Потом из неравенства Маркова получают оценку больших уклонений

$$P\left(f\left(x_{N(\varepsilon)}\right) - \min_{x\in Q} f(x) \geq \sigma\right) \leq \varepsilon/\sigma,$$

которую переписывают в виде

$$P\left(f\left(x_{N(\varepsilon\sigma)}\right) - \min_{x\in Q} f(x) \geq \sigma\right) \leq \varepsilon,$$

где

$$N(\varepsilon\sigma) = \mathrm{O}\left(\left(N + \min\left\{L/\mu, \sqrt{NL/\mu}\right\}\right)\left(\ln\left(\Delta f/\varepsilon\right) + \ln\left(\sigma^{-1}\right)\right)\right).$$

Мы привели здесь это наблюдение, потому что оно оказывается полезным и во многих других контекстах, в которых рандомизированный метод сходится со скоростью геометрической прогрессии.

При наличии дополнительной структуры у задачи (2.1.5) приведенные оценки можно было получить (и даже немного улучшить, например, учитывая разреженность) исходя из рандомизированных покомпонентных методов (например, ALPHA или APPROX (см. работы P. Richtarik-а http://www.maths.ed.ac.uk/~richtarik/ ) или ACRCD* из замечания 8 [45], см. также раздел 5.1 главы 5) для “двойственной” к (2.1.5) задаче [33, 45, 138, 307, 331].[17] Заметим также, что в работе [45] (см. также раздел 5.1 главы 5) показывается, как можно просто получить часть выписанных оценок с помощью метода, работающего по оценкам (2.1.8) (см. ниже).

В книге [99] Б.Т. Поляк отмечает, что если рандомизация осуществляется каким-то специальным образом, например, таким, что[18]

$$E\left[\left\|\nabla f(x,\xi)\right\|_*^2\right] \leq C_n \left\|\nabla f(x)\right\|_*^2 + \Delta, \tag{2.1.6}$$

где $\Delta \geq 0$ некоторая малая погрешность, и в точке минимума $\nabla f(x) = 0$, то приведенные выше оценки (2.1.2), (2.1.4) можно существенно улучшить. Примеры будут приведены ниже в подразделе 2.1.4 (см. (2.1.22)). В частности, в сильно выпуклом случае можно по-лучить геометрическую скорость сходимости. Важно отметить, что при рандомизации, возникающей в покомпонентных спусках, спусках по направлению и безградиентных методах в гладком случае условие (2.1.6) выполняется [99, 264, 270]. Мы вернемся к этому кругу вопросов в подразделе 2.1.4. Описанная же выше конструкция (с довольно грубым неравенством Маркова) используется в данном контексте [264, 270] для (точной!) оценки больших уклонений. Причем за счет регуляризации функционала, о которой было сказано выше, все это переносится и просто на гладкий случай без предположения сильной выпуклости.

[17] Строго говоря, построение двойственной задачи предполагает возможность явного выделения в функционале в виде отдельного слагаемого сильно выпуклого композита – желательно сепарабельного.

[18] Если рассматривать приложения методов стохастической оптимизации к СТО [320], а в правой части неравенства вместо $\left\|\nabla f(x)\right\|_*^2$ писать $\left\|\nabla f(x)\right\|_*^{1/2}$, то выписанное неравенство будет соответствовать условиям малого шума Цыбакова–Массара, Бернштейна [144].

### 2.1.3 Стохастические градиентные методы с неточным оракулом

В этом пункте мы опишем, что можно получить, если дополнительно известно, что $f(x)$ – гладкая по $x$ функция, с константой Липшица градиента $L$ и(или) сильно выпуклая с константой $\mu \ge 0$, но вычисление стохастического градиента на каждом шаге происходит с неконтролируемой неточностью $\delta$, вообще говоря, не случайной природы.[19]

**Замечание 2.1.3.** И гладкости и сильной выпуклости можно добиться искусственно (см., например, http://arxiv.org/pdf/1603.05643.pdf ). Как уже отмечалось в подразделе 2.1.2, сильная выпуклость всегда легко получается регуляризацией функционала в исходной задаче. Как правило, это не дает ничего нового с точки зрения выписанных оценок (и даже может ухудшать эти оценки на логарифмический фактор), но в ряде специальных случаев (см. ниже) это может давать определенные преимущества. Кроме того, такая регуляризация иногда просто необходима для корректности постановки. Это связано с тем, что в общем случае даже для гладких детерминированных выпуклых задач мы можем гарантировать сходимость итерационного метода лишь по функции, но не по аргументу. Для сходимости по аргументу нужна сильная выпуклость функционала, которую и обеспечивают должной регуляризацией (см., например, конец подраздела 2.1.2) – при этом сходимость по аргументу имеет место к решению регуляризованной задачи. Идея регуляризации используется в популярном методе двойственного сглаживания [271] (регуляризация двойственной задачи с целью улучшения гладких свойств прямой). В отличие от прямой регуляризации, эта техника хорошо работает только для вполне конкретных задач, имеющих определенную (седловую – Лежандрову) структуру (модель), когда исходная задача имеет явное двойственное представление (см. пример 2.1.4), введя в которое регуляризацию, можно явно (эффективно) пересчитать во что превратится исходная прямая задача. Другой пример сглаживания будет приведен в подразделе 2.1.4.

---

[19] Особое внимание таким постановкам стали уделять после выхода книг [59, 99]. В них обстоятельно изучается “влияние помех”, в том числе не случайной природы, на методы выпуклой оптимизации.

Сформулируем более точно предположения об оракуле, выдающем стохастический градиент, следуя [181, 183].[20]

**Предположение 2.1.1.** $(\delta, L, \mu)$*-оракул выдает (на запрос, в котором указывается только одна точка* $x$*) такую пару* $(F(x,\xi), G(x,\xi))$ *(с.в.* $\xi$ *независимо разыгрывается из одного и того же распределения, фигурирующего в постановке (2.1.1)), что для всех* $x \in Q$ *ограничена дисперсия*

$$E_\xi\left[\left\|G(x,\xi) - E_\xi\left[G(x,\xi)\right]\right\|_*^2\right] \le D,$$

*и для любых* $x, y \in Q$

$$\frac{\mu}{2}\|y-x\|^2 \le E_\xi\left[f(y,\xi)\right] - E_\xi\left[F(x,\xi)\right] - \left\langle E_\xi\left[G(x,\xi)\right], y-x\right\rangle \le \frac{L}{2}\|y-x\|^2 + \delta.$$

Из недавних результатов [39, 140, 179, 184, 185, 193, 205, 206, 264, 270] можно получить общий метод (мы приводим огрубленный вариант оценки времени работы этого метода для большей наглядности), с такими оценками скорости сходимости[21]

$$\min\left\{O\left(\frac{LR^2}{N^{p+1}} + \sqrt{\frac{DR^2}{N}} + N^p\delta\right), O\left(LR^2\exp\left(-\Upsilon N\cdot\left(\frac{\mu}{L}\right)^{\frac{1}{p+1}}\right) + \frac{D}{\mu N} + \left(\frac{L}{\mu}\right)^{\frac{p}{p+1}}\delta\right)\right\}, \quad (2.1.7)$$

[20] В работе [183] собрано много различных мотиваций такому предположению (определению), обобщающему классическую концепцию $\delta$-субградиента [99]. В определенном смысле это предположение 2.1.1 наиболее общее и, одновременно, наиболее точно отражающее спектр всевозможных приложений [40, 42, 43, 45, 138]. В частности, об одном из таких приложений написано в разделе 1.6 главы 1, о другом приложении написано в разделе 2.3 этой главы 2. Здесь также можно отметить пример 2.1.4, приведенный ниже.

[21] Оценки характеризуют достигнутую в среднем точность (по оптимизируемому функционалу) после $N$ итераций. При этом в случае когда минимум достигается на втором аргументе (выгодно использовать факт наличия $\mu$-сильной выпуклости) под $N$ правильнее понимать не число итераций, а число обращений к $(\delta, L, \mu)$-оракулу [193]. Отметим, что при $p = 0$ оценку можно сделать непрерывной по параметру $\mu \ge 0$ (см. [181, 139]). Также заметим, что в метод (например, в размер шагов) не входит требуемое число итераций (или желаемая точность – одно через другое выражается). Таким образом, можно говорить об адаптивности метода. Отметим, что за это не приходится дополнительно платить логарифмическую плату [269, 266]. Отметим также, что при $p = 1$ метод наихудшим (а при $p = 0$ наилучшим) способом (среди всех разумных вариаций градиентного метода) накапливает неточность в вычислении градиента. Это переносится и на негладкие задачи (см. далее).

где $\Upsilon \geq 1$ – некоторая константа (можно считать $\Upsilon = O(\ln n)$), а параметр $p \in [0,1]$ подбирается "оптимально" перед запуском метода исходя из масштаба шума $\delta$. Для лучшего понимания оценки (2.1.7) полезно ее переписать в еще более огрубленном виде[22]

$$\min\left\{ O\left( \frac{LR^2}{N^{p+1}} + \sqrt{\frac{DR^2}{N}} + N^p\delta \right), O\left( LR^2 \exp\left( -\Upsilon N \cdot \left(\frac{\mu}{L}\right)^{\frac{1}{p+1}} \right) + \frac{D}{\mu N} + N^p\delta \right) \right\}.$$

Этот общий метод – есть в некотором смысле "выпуклая комбинация" двойственного градиентного метода (DGM) и быстрого градиентного метода[23] (FGM) [185, 193], оценки скорости сходимости для которых имеют соответственно вид:

$$\text{(DGM)} \qquad \min\left\{ O\left( \frac{LR^2}{N} + \sqrt{\frac{DR^2}{N}} + \delta \right), O\left( LR^2 \exp\left( -\Upsilon_1 N \frac{\mu}{L} \right) + \frac{D}{\mu N} + \delta \right) \right\},$$

$$\text{(FGM)} \qquad \min\left\{ O\left( \frac{LR^2}{N^2} + \sqrt{\frac{DR^2}{N}} + N\delta \right), O\left( LR^2 \exp\left( -\Upsilon_2 N \sqrt{\frac{\mu}{L}} \right) + \frac{D}{\mu N} + \sqrt{\frac{L}{\mu}}\delta \right) \right\}.$$

Комбинируя эти два метода можно непрерывно настраиваться (оптимально подбирая метод, регулируя $p \in [0,1]$) на шум (известного масштаба). В этой связи также полезно отметить (аналогичный факт имеет место и для покомпонентного варианта FGM [45]), что FGM есть специальная выпуклая комбинация прямого градиентного метода (PGM), оценки скорости сходимости которого совпадают с оценками DGM, и метода зеркального спуска / двойственных усреднений [6, 135] (по-видимому, здесь вместо зеркального спуска можно использовать и метод из работы [280]). Нельзя в этой связи не обратить внимание на то, что комбинация двух методов привела к новому методу, работающему лучше, чем каждый из методов по отдельности. Отметим здесь также недавнюю работу [280], в которой предлагается общий способ получения ускоренных (быстрых) методов.

---

[22] Насколько нам известно, для всех методов, которые используют только градиент и значение функции (или их стохастические аналоги) накопление шума методом со скоростью $N^p\delta$ с $p \in [0,1]$ – является общим местом.

[23] Отметим, что при $D = 0$ не улучшаемые оценки, которые дает метод FGM [163], были установлены Б.Т. Поляком [99] для ряда других многошаговых методов (метод тяжелого шарика, сопряженных градиентов). Отличие в том, что тогда оценки были установлены локально. Все приведенные в данном разделе оценки – глобальные, т.е. не требуют оговорок о близости точки старта к решению, для гарантии нужной скорости сходимости. Заметим также, что техника установления локальной сходимости основана, как правило, на первом методе Ляпунова [98, 99], в то время как глобальной – на втором [64, 99]. При этом функцию Ляпунова можно искать по непрерывному аналогу итерационного процесса – системе дифференциальных уравнений [64]. Скажем, для обычного градиентного метода это будет система [100] (Коши, 1847): $dx/dt = -\nabla f(x)$. Скорости сходимости у итерационного процесса и его непрерывного аналога могут отличаться. Скажем, непрерывный аналог метода Ньютона сходится за конечное время. Другой пример – метод зеркального спуска [91]. Недавно появились работы, посвященная и непрерывному аналогу FGM, см. список литературы в [326]. См. также https://arxiv.org/pdf/1611.02635v1.pdf.

Приведенный выше результат (формула (2.1.7)) был получен в цикле работ [181, 183, 184, 185] при $p \in (0,1)$ в не сильно выпуклом случае – без оценок вероятностей больших уклонений, а в сильно выпуклом случае рассматривались только ситуации, когда $p = 0,1$, причем рассматривался только детерминированный случай ($D = 0$). Формула (2.1.7) в более точном виде (с вероятностями больших уклонений и с явно выписанными числовыми константами) была получена в цикле совместных работ [39, 193]. Основную тяжесть в получении этой оценки (связанную с проработкой большого количества технических деталей в доказательстве оценки) взял на себя соавтор П.Е. Двуреченский, поэтому было решено не включать в диссертацию доказательство формулы (2.1.7). Впрочем, впоследствии автору удалось заметно упростить доказательство формулы (2.1.7), и даже существенно обобщить эту формулу (изложению наиболее интересной части полученного обобщения посвящен следующий раздел, содержащий доказательства всех необходимых утверждений).

Вся последующая часть подраздела 2.1.3 будет посвящена “обзорному” обсуждению отмеченного результата (формулы (2.1.7)) и его окрестностей.

Прежде всего, заметим, что дисперсию у первого аргумента минимума в (2.1.7) можно уменьшать в $m$ раз, запрашивая на одном шаге реализацию стохастического градиента не один раз, а $m$ раз, и заменяя стохастический градиент средним арифметическим [163, 181, 227] (в случае тяжелых хвостов у стохастических градиентов лучше пользоваться более робастными оценками, например, медианного типа [91]).[24] Это имеет смысл делать, если слагаемое, отвечающее стохастичности, доминирует. Важно, что мы при этом не увеличиваем число итераций, и слагаемое $N^p\delta$ остается прежним. Отметим, что число вызовов оракула при этом увеличивается, но тем не менее, в некоторых ситуациях такой подход может оказаться оправданным. Такая игра используется[25] в способе получения

[24] Этот прием в западной литературе часто называют “mini-batch” [163]. Подробнее об этом приеме написано в следующем разделе.

[25] Вместе с идеей рестартов [2, 39, 45, 193, 205, 206, 228], распространяющей (ускоряющей) практически любой итерационный метод (желательно с явной оценкой необходимого числа итераций $N(\varepsilon)$ для

второго аргумента оценки (2.1.7). В этой связи оценку (2.1.7) правильнее переписать следующим образом (здесь $N(\varepsilon)$ – число обращений к $(\delta, L, \mu)$-оракулу, необходимых для достижения в среднем по функции точности $\varepsilon$, индекс 1 соответствует просто выпуклому, а индекс 2 сильно выпуклому случаю):

$$N_1(\varepsilon) = \max\left\{\mathrm{O}\left(\left(\frac{LR^2}{\varepsilon}\right)^{\frac{1}{p+1}}\right), \mathrm{O}\left(\frac{DR^2}{\varepsilon^2}\right)\right\}, N_2(\varepsilon) = \max\left\{\mathrm{O}\left(\left(\frac{L}{\mu}\right)^{\frac{1}{p+1}} \ln\left(\frac{LR^2}{\varepsilon}\right)\right), \mathrm{O}\left(\frac{D}{\mu\varepsilon}\right)\right\} \quad (2.1.8)$$

при (условия на допустимый уровень шума, при котором оценки (2.1.8) имеют такой же вид, с точностью до $\mathrm{O}(1)$, как если бы шума не было)

$$\delta_1(\varepsilon) \le \mathrm{O}\left(\varepsilon \cdot \left(\frac{\varepsilon}{LR^2}\right)^{\frac{p}{p+1}}\right), \; \delta_2(\varepsilon) \le \mathrm{O}\left(\varepsilon \cdot \left(\frac{\mu}{L}\right)^{\frac{p}{p+1}}\right). \quad (2.1.9)$$

Как уже отмечалось, выписанные оценки (2.1.7) ((2.1.8), (2.1.9)) характеризуют скорость сходимости в среднем. Они с одной стороны не улучшаемы[26] с точностью до мультипликативной константы (см. подраздел 2.1.2 и [91, 128]), а с другой стороны достигаются. Все это (неулучшаемость оценок) справедливо и при $\delta = 0$ и(или) $D = 0$. При этом в случае $D = 0$, $\mu = 0$ необходимо считать, что требуемое число итераций для достижения точности $\varepsilon$ удовлетворяет неравенству $N(\varepsilon) \le n$ [91], в противном случае оценки улучшаемы – метод центров тяжести [91, 163], с оценкой числа итераций типа $\mathrm{O}(n \ln(B/\varepsilon))$, где $|f(x)| \le B$. В терминах больших отклонений возникают оценки, аналогичные тем, которые были приведены в подразделе 2.1.2, см. [193].

---

достижения заданной точности $\varepsilon$) на случай сильно выпуклого функционала. Нетривиально здесь то, что при довольно общих условиях при таком распространении сохраняется (и работает уже в условиях сильной выпуклости) свойство оптимальности исходного метода. Впрочем, к рестартам стоит относиться очень аккуратно (причины обсуждаются в следующем разделе, а также в разделе 3.2 главы 3 и в разделе 5.1 главы 5).

[26] В нижнюю оценку во втором выражении под знаком минимума при экспоненте вместо $LR^2$ входит $\mu R^2$, а константа $\Upsilon = 1$ в (2.1.7). Впрочем, получить вместо фактора $LR^2$ фактор $\mu R^2$ можно аккуратно проанализировав оценки, даваемые с помощью техники рестартов (см., например, [135, 193] и следующий раздел).

Отмеченные результаты переносятся и на прокс-структуры отличные от евклидовой [193]. При этом рассмотрение какой-либо другой $q$-нормы ($l_q$-нормы) в прямом пространстве ($q \geq 1$), отличной от евклидовой, в сильно выпуклом случае (когда минимум достигается на втором выражении в (2.1.7)), как правило, не имеет смысла. Связано это с тем, что квадрат евклидовой асферичности $q$-нормы, который может возникать в оценках числа обусловленности прокс-функции в $q$-норме (это число, в свою очередь, оценивает увеличение числа итераций метода при переходе от евклидовой норме к $q$-норме), больше либо равен 1. Равенство достигается на евклидовой норме. Скажем, для 1-нормы эта асферичность оценивается снизу размерностью пространства [92, 228]. Другими словами, действительно, можно выбирать в сильно выпуклом случае $q$-норму (отличную от евклидовой) и получать оценки на число итераций вида (см. (2.1.7), (2.1.8) и подраздел 2.1.2)

$$\mathrm{O}\left(\left(\frac{L}{\mu}\omega\right)^{\frac{1}{p+1}}\ln\left(\frac{LR^2}{\varepsilon}\right)\right),\ \omega = \sup_{x\in Q}\frac{2V(x,x_0)}{\alpha\|x-x_0\|_q^2},$$

где $L$ и $\mu$ считаются относительно $q$-нормы, а $R^2$ – брэгмановское “расстояние” от точки старта до решения (см. замечание 2.1.1) . Однако смысла, как правило, в этом нет, поскольку $\omega \geq 1$, а число обусловленности $\chi = L/\mu$ не меньше чем в случае выбора $2$-нормы. Например [228], для функции $\|x\|_2^2 = x_1^2 + \ldots + x_n^2$ в евклидовой норме число обусловленности $\chi = 1$, а в 1-норме $\chi = n$. Тем не менее, выгода от использования не евклидовой прокс-структуры в сильно выпуклом случае может быть, если рассматривать задачи композитной оптимизации, в которых сильная выпуклость приходит от композитного слагаемого (см. замечание 2.1.6). Так в приложениях, описанных в работах [33, 138], в качестве композитного слагаемого возникает сильно выпуклая в 1-норме энтропийная функция. Отметим, что энтропию при этом нельзя использовать в качестве прокс-функции. Нужно брать (и это можно сделать, см. замечание 2.1.2) другую прокс-функцию, соответствующую 1-норме, которая обеспечивает (по-видимому, оптимально возможную) оценку $\omega = \mathrm{O}(\ln n)$. Подробнее об этом написано в разделе 2.3 этой главы.

Заметим также, что обычный метод FGM в не стохастическом сильно выпуклом случае для задач безусловной оптимизации, в действительности, дает оценку (следует сравнить с (2.1.8)) [93]:

$$\mathrm{O}\left(\sqrt{\frac{L}{\mu}}\ln\left(\frac{f(x_0)-f(x_*)}{\varepsilon}\right)\right).$$

Поскольку $\nabla f(x_*)=0$, то $f(x_0)-f(x_*)\le LR^2/2$. Если рассматривается задача условной оптимизации (на выпуклом множестве $Q\subset\mathbb{R}^n$), то, вообще говоря, $\nabla f(x_*)\ne 0$, следовательно, нельзя утверждать, что $f(x_0)-f(x_*)\le LR^2/2$. В [181, 193] предлагается обобщение классического FGM для класса гладких сильно выпуклых задач, которое фактически позволяет вместо $f(x_0)-f(x_*)$ писать нижнюю оценку $\mu R^2/2\le f(x_0)-f(x_*)$ в том числе для задач условной оптимизации. Заметим, что это же наблюдение справедливо для описываемых в данной работе методов (мы не стали писать $\mu R^2$ вместо $LR^2$ в (2.1.7) и далее для сохранения непрерывности выписанных оценок по $\mu$, т.е. чтобы делать меньше оговорок о переключениях с сильно выпуклого случая на выпуклый при малых значениях $\mu$).

**Замечание 2.1.4.** Отметим, что пока нам не известно (для произвольной прокс-структуры, отличной от евклидовой) строгое обоснование оценок (2.1.7) ((2.1.8), (2.1.9)) с вероятностями больших отклонений для случая не ограниченного множества $Q$. В известном нам способе получения оценок вероятностей больших уклонений (см., например, [181, 193]), к сожалению, явно используется предположение об ограниченности множества $Q$. С другой стороны, для используемых в статье неускоренных методов (кроме Франк–Вульфа и кроме PGM в варианте [135], для PGM в варианте [181] все хорошо) оценки на скорость сходимости обычно получаются в следующем виде [6, 40, 42, 92, 181, 193, 269]:

$$\sum_{k=0}^{N}\lambda_k\cdot\left(f(x_k)-f(x_*)\right)\le V(x_*,x_0)-V(x_*,x_{N+1})+\sum_{k=0}^{N}\lambda_k\left\langle G(x_k,\xi_k)-\nabla f(x_k),x_*-x_k\right\rangle+$$
$$+\tilde{\Delta}_N\left(\{\lambda_k\}_{k=0}^{N},\left\{\left\|G(x_k,\xi_k)-\nabla f(x_k)\right\|_*^2\right\}_{k=0}^{N},\delta\right),\ \{\lambda_k\}\ge 0$$

или, в случае ускоренных методов (к которым относится FGM и его производные), в похожем, но немного более громоздком (с большим числом параметров и оценивающих последовательностей). Опуская в правой части $V(x_*,x_{N+1})$, далее оптимально подбирают параметры метода $\{\lambda_k\}$, получают оценку скорости сходимости метода по функции в среднем. Если считать, что $\|x_*-x_k\|=\mathrm{O}(R)$, то отсюда также получают оценки скорости сходимости и с вероятностями больших уклонений (используется обобщение неравенства

Азума–Хефдинга для последовательности мартингал-разностей [181, 309], см. также раздел 1.5 главы 1 и раздел 6.1 главы 6). В детерминированном случае соотношение $\|x_* - x_k\| = O(R)$ имеет место (всегда в евклидовом случае, и в зависимости от метода в общем случае) ввиду сходимости метода и того, что параметры оптимально подбираются так, что слагаемые $V(x_*, x_0)$ и $\tilde{\Delta}_N$ одного порядка (отличаются обычно не более чем в 10 раз):

$$\frac{1}{2}\|x_k - x_*\|^2 \le V(x_*, x_k) \le V(x_*, x_0) + \tilde{\Delta}_{k-1} \le V(x_*, x_0) + \tilde{\Delta}_N .$$

В случае стохастического оракула, к сожалению, такие рассуждения уже не проходят. Можно, однако, из таких соображений оценить $E_{x_k}\left[V(x_*, x_k)\right]$. Дальше угадывается хвост распределения случайной величины $\|x_* - x_k\|$ исходя из выписанного выше соотношения, которое стоит понимать как равенство, т.е. хвост распределения ищется как неподвижная точка (а точнее ее оценка). Задавшись определенным доверительным уровнем $\sigma \ge 0$ можно оценить "эффективный" $R$: с вероятностью $\ge 1-\sigma$ имеют место неравенства $V(x_*, x_k) \le R$, $k = 0,...,N$. В частности, для субгауссовских стохастических градиентов $R = O\left(V(x_*, x_0)\ln^2(N/\sigma)\right)$, а для равномерно ограниченных – $R = O\left(V(x_*, x_0)\ln(N/\sigma)\right)$. Детали можно посмотреть в доказательстве теоремы 4 работы [40] и в работе [45] (см. также раздел 1.5 главы 1, раздел 2.2 и раздел 3.1 главы 3). Примечательно, что все эти рассуждения в случае не ограниченного множества $Q$ не требуют равномерной ограниченности констант Липшица (функции, градиента) на всем $Q$ [6]. Похожим образом можно получать оценки вероятностей больших уклонений в сильно выпуклом случае в онлайн контексте (см. конец этого пункта). К сожалению, не все методы обладают такими же свойствами. Например, PGM [135] (в случае детерминированного оракула и не евклидовой прокс-структуры [45]) гарантирует лишь, что $\|x_* - x_k\| = O(R)$, $k = 0,...,N$, если [135] (прокс-диаметр здесь не нужен):

$$R = \max\left\{\|x - x_*\|:\ x \in Q,\ f(x) \le f(x_0)\right\}.$$

Хотя PGM и является релаксационным методом ( $f(x_{k+1}) \le f(x_k)$ ), возможно, что $R = \infty$. Требование $R < \infty$ (коэрцитивности) не является сильно обременительным. Его можно обеспечивать за счет регуляризации задачи [37].

Полезно также иметь в виду, что за счет допускаемой неточности оракула, можно погрузить задачу с гельдеровым градиентом, т.е. удовлетворяющим неравенству $\left\|\nabla f(x)-\nabla f(y)\right\|_{*}\le L_{\nu}\left\|x-y\right\|^{\nu}$, при некотором $\nu\in[0,1]$ (в том числе и негладкую задачу с ограниченной нормой разности субградиентов при $\nu=0$) в класс гладких задач с неточным оракулом, характеризующимся точностью $\delta$ и [183]

$$L=L_{\nu}\left[\frac{L_{\nu}(1-\nu)}{2\delta(1+\nu)}\right]^{\frac{1-\nu}{1+\nu}}. \tag{2.1.10}$$

Заметим, в этой связи, что если в предположении 2.1.1 считать

$$E_{\xi}\left[f(y,\xi)\right]-E_{\xi}\left[F(x,\xi)\right]-\left\langle E_{\xi}\left[G(x,\xi)\right],y-x\right\rangle\le\frac{L}{2}\left\|y-x\right\|^{2}+M\left\|y-x\right\|+\delta,$$

то вместо $D$ в (2.1.7) стоит писать $M^{2}+D$ (см. [237], а также другие работы G. Lan-a http://www.ise.ufl.edu/glan/publications/ ).

Таким образом, например, можно получить оценки (2.1.2), (2.1.4) из оценки (2.1.7). В частности, метод двойственных усреднений и зеркальный спуск (см. подраздел 2.1.2) можно получить из PGM в варианте [181] с неточным оракулом и $L=M^{2}/(2\delta)$. Но наряду с введенной нами искусственной неточностью оракула, можно допустить, что имеется также реальная неточность оракула. Несложно привести оценки (на базе формулы (2.1.7) и ((2.1.8), (2.1.9))) сочетающие наличие в задаче искусственной и реальной неточности [42].

Как мы предполагали выше, множество $Q$ должно быть достаточно простой структуры, чтобы на него можно было эффективно проектироваться. Однако в приложениях часто возникают задачи условной минимизации [99], в которых, например, есть ограничения вида $g(x)\le 0$, где $g(x)$ – выпуклые функции [83]. “Зашивать” эти ограничения в $Q$, как правило, не представляется возможным в виду вышесказанного требования о легкости проектирования. Тем не менее, на основе описанного выше можно строить (за дополнительную логарифмическую плату) двухуровневые методы (наверное, лучше говорить “методы уровней”, чтобы не было путаницы с многоуровневой оптимизацией, см. пример 2.1.5) условной оптимизации (см. [93], а также работы G. Lan-a). При этом на каждом шаге такого метода потребуется проектироваться на пересечение множества $Q$ с некоторым полиэдром, вообще говоря, зависящим от номера шага. Последнее обстоятельство в общем случае сужает класс задач, к которому применимы такие многоуровневые методы до класса задач, к которым применимы методы внутренней точки [93]. В частности, возникает довольно обременительное условие на размер пространства, в котором проходит опти-мизация: $n \sim 10^4 - 10^5$. Все это не удивительно, поскольку имеются нижние оценки [91] (рассматриваются аффинные ограничения в виде равенств, аналогично могут быть рассмотрены и неравенства), показывающие, что в общем случае для нахождения такого $x \in \mathbb{R}^n$, что $\|Ax - b\|_2 \le \varepsilon$ потребуется не меньше, чем $\Omega\left(\sqrt{L_x} R_x / \varepsilon\right)$ операций типа умножения $Ax$ ($L_x = \sigma_{\max}(A) = \lambda_{\max}\left(A^T A\right)$ с – максимальное собственное значение матрицы $A^T A$, $R_x = \|x^*\|_2 = \left\|\left(A^T A\right)^{-1} A^T b\right\|_2$ ). Аналогичное можно сказать и про седловые задачи: для отыскания такой пары $(x, y)$, что (левая часть этого неравенства всегда неотрицательная)

$$\max_{\tilde{y} \in S_n(1)} \tilde{y}^T A x - \min_{\tilde{x} \in S_n(1)} y^T A \tilde{x} \le \varepsilon$$

потребуется не меньше, чем $\Omega(\Lambda/\varepsilon)$ ($\Lambda$ – максимальный по модулю элемент матрицы $A$) операций типа умножения $Ax$ и $y^T A$. Заметим, что обе выписанные нижние оценки справедливы при условии, что число итераций (операций типа умножения $Ax$) $k \le n$. Как следствие, в общем случае даже для гладкой детерминированной сильно выпуклой постановки при наличии всего лишь аффинных ограничений $Ax = b$ нельзя надеяться на быстрое решение. Тем не менее, некоторые дополнительные предположения в ряде случаев позволяют ускорить решение таких задач (см., например, [37, 138]).

**Замечание 2.1.5 (Ю.Е. Нестеров, см. также [6]).** Задача поиска такого $x^*$, что $Ax^* = b$ сводится к задаче выпуклой гладкой оптимизации

$$f(x) = \|Ax - b\|_2^2 \to \min_x .$$

Нижняя оценка для скорости решения такой задачи [91] (см. также формулу (2.1.7) с $\delta = D = 0$, $p = 1$) имеет вид: $f(x_k) \ge \Omega\left(L_x R_x^2 / k^2\right)$. Откуда следует, что только при $k \ge \Omega\left(\sqrt{L_x} R_x / \varepsilon\right)$ можно гарантировать выполнение неравенства $f(x_k) \le \varepsilon^2$, т.е. $\|Ax_k - b\|_2 \le \varepsilon$. Заметим, что эта нижняя оценка для специальных матриц может быть улучшена. Причем речь идет не о недавних результатах D. Spielman-a [315] (премия Неванлины 2010 года), а о более простой ситуации. Вернемся к задаче поиска вектора PageRank (примеры 2.1.1, 2.1.2), которую мы перепишем как ($n$ – размерность вектора $x$)

$$Ax = \begin{pmatrix} \left(P^T - I\right) \\ n^{-1} \ldots\ldots n^{-1} \end{pmatrix} x = \begin{pmatrix} 0 \\ n^{-1} \end{pmatrix} = b \ , \ I \text{ – единичная матрица.}$$

По теореме Фробениуса–Перрона [95] решение такой системы с неразложимой матрицей $P$ единственно и положительно $x>0$. Сведем решение этой системы уравнений к вырожденной задаче выпуклой оптимизации

$$\frac{1}{2}\|x\|_2^2 \to \min_{Ax=b}.$$

Построим двойственную к ней задачу [259]

$$\min_{Ax=b}\frac{1}{2}\|x\|_2^2 = \min_x \max_\lambda \left\{\frac{1}{2}\|x\|_2^2 + \langle b-Ax,\lambda\rangle\right\} =$$

$$= \max_\lambda \min_x \left\{\frac{1}{2}\|x\|_2^2 + \langle b-Ax,\lambda\rangle\right\} = \max_\lambda \left\{\langle b,\lambda\rangle - \frac{1}{2}\left\|A^T\lambda\right\|_2^2\right\}.$$

Поскольку система $Ax=b$ совместна, то по теореме Фредгольма не существует такого $\lambda$, что $A^T\lambda=0$ и $\langle b,\lambda\rangle>0$, следовательно, двойственная задача имеет конечное решение (т.е. существует ограниченное решение двойственно задачи $\lambda^*$). Зная решение $\lambda^*$ двойственной задачи

$$\langle b,\lambda\rangle - \frac{1}{2}\left\|A^T\lambda\right\|_2^2 \to \max_\lambda$$

можно восстановить решение прямой задачи (из условия оптимальности по $x$): $x(\lambda)=A^T\lambda$. Однако важно здесь то, что FGM [93] для этой двойственной задачи дает возможность попутно получать следующую оценку на норму этого градиента [6]:

$$\|Ax_k-b\|_2 = O\left(\frac{L_yR_y}{k^2}\right),$$

где $x_k$ есть известная выпуклая комбинация

$$\{x(\lambda_i)\}_{i=1}^k,\ L_y=\sigma_{\max}\left(A^T\right)=\sigma_{\max}(A),\ R_y=\left\|\lambda^*\right\|_2,$$

где можно считать, что $\lambda^*$ – решение двойственной задачи с наименьшей евклидовой нормой. Кажется, что это противоречит нижней оценке $\|Ax_k-b\|_2 \ge \Omega\left(\sqrt{L_x}R_x/k\right)$. Однако, важно напомнить [91], что эта нижняя оценка установлена для всех $k\le n$, и она будет улучшена, в результате описанной процедуры только если дополнительно предположить, что матрица $A$ удовлетворяет следующему условию $L_yR_y \ll n\sqrt{L_x}R_x$, что сужает класс, на котором была получена нижняя оценка $\Omega\left(\sqrt{L_x}R_x/k\right)$. В типичных ситуациях можно ожидать, что $R_y \gg R_x$ ($R_x\le\sqrt{2}$). Это обстоятельство мешает выполнению требуемого условия.

**Пример 2.1.4.** Если имеется дополнительная информация о структуре седловой задачи, то можно её использовать для ускорения [92, 265]. Более того, многие современные постановки задач (негладкой) выпуклой оптимизации (в частности, связанные с compressed sensing и $l_1$-оптимизацией) в пространствах огромных размеров специально стараются представить седловым образом с целью получения эффективного решения (см. работы А.С. Немировского, А.Б. Юдицкого, например, [227]). Далее будет разобран один простой пример (немного обобщающий результаты [181, 183], см. также [6]), демонстрирующий возможности градиентных методов с неточным оракулом в седловом контексте. Рассматривается седловая задача ( $x \in \mathbb{R}^n$ , $y \in \mathbb{R}^m$ )

$$f(x) = \max_{\|y\|_2 \le R_y} \{G(y) + \langle By, x \rangle\} \to \min_{\|x\|_2 \le R_x},$$

где функция $G(y)$ – сильно вогнутая с константой $\kappa$ относительно 2-нормы и константой Липшица градиента $L_G$ (также в 2-норме). Тогда функция $f(x)$ будет гладкой, с константой Липшица градиента в 2-норме $L_f = \sigma_{\max}(B)/\kappa$. Казалось бы, что мы можем решить задачу минимизации функции $f(x)$ за $\mathrm{O}\left(\sqrt{\sigma_{\max}(B) R_x^2/(\kappa\varepsilon)}\right)$ итераций, где $\varepsilon$ - желаемая точность по функции. Но это возможно, только если мы можем абсолютно точно находить $\nabla f(x) = By^*(x)$, где $y^*(x)$ – решение вспомогательной задачи максимизации по $y$ при заданном $x$. В действительности, мы можем решать эту задачу (при различных $x$) лишь приближенно. Если мы решаем вспомогательную задачу быстрым градиентным методом [92] с точностью $\delta/2$ (на это потребуется $\mathrm{O}\left(\sqrt{L_G/\mu}\ln\left(L_G R_y^2/\delta\right)\right)$ итераций), то пара $\left(G\left(y_{\delta/2}(x)\right) + \left\langle By_{\delta/2}(x), x\right\rangle, By_{\delta/2}(x)\right)$, где $y_{\delta/2}(x)$ – $\delta/2$-решение вспомогательной задачи, будет $(\delta, 2L_f, 0)$-оракулом [181, 183]. Выбирая $\delta = \mathrm{O}\left(\varepsilon\sqrt{\varepsilon/\left(L_f R_x^2\right)}\right)$ (см. формулу (9) при $p = 1$), получим после

$$\mathrm{O}\left(\sqrt{\frac{L_G \sigma_{\max}(B) R_x^2}{\kappa^2 \varepsilon}} \ln\left(\frac{L_f L_G R_x^2 R_y^2}{\varepsilon}\right)\right)$$

итераций (на итерациях производится умножение матрицы $B$ на вектор/строчку и вычисление градиента $G(y)$) $\varepsilon$-решение задачи минимизации $f(x)$. Отметим, что если не использовать сильную вогнутость функции $G(y)$, то для получения пары $(x_N, y_N)$, удовлетворяющей неравенству

$$\max_{\|y\|_2 \le R_y} \left\{ G(y) + \langle By, x_N \rangle \right\} - \min_{\|x\|_2 \le R_x} \left\{ G(y_N) + \langle By_N, x \rangle \right\} \le \varepsilon,$$

потребуется $\Omega\Big(\max\big\{L_G R_y^2, \sigma_{\max}(B) R_x R_y\big\}/\varepsilon\Big)$ итераций (см., например, [91, 163, 259]).

Интересно отдельно разобрать ситуацию, когда вместо множества $\|y\|_2 \le R_y$ фигурирует симплекс $S_m(R_y)$, $G(y) = -\sum_{k=1}^{m} y_k \ln(y_k/R_y)$ – сильно вогнутая в 1-норме с константой $\kappa = 1$ функция и $R_x = \infty$ (энтропийно-линейное программирование [37]). В этом случае мы не можем обеспечить даже равномерной ограниченности градиента функции $G(y)$. Тем не менее, также можно рассчитывать [37] на зависимость $O\Big(\varepsilon^{-1/2}\ln\big(\varepsilon^{-1}\big)\Big)$ числа итераций от точности $\varepsilon$ для критерия:

$$\max_{y \in S_m(R_y)} \min_x \left\{ G(y) + \langle By, x \rangle \right\} - \min_x \left\{ G(y_N) + \langle By_N, x \rangle \right\} \le \varepsilon .$$

При этом вместо энтропии в качестве функции $G(y)$ можно брать любую сильно вогнутую в 1-норме функцию, для которой решение задачи максимизации (вычисление $f(x)$ с точностью $\varepsilon$) может быть осуществлено за $O\Big(\ln\big(\varepsilon^{-1}\big)\Big)$. В примере с энтропией, для $f(x)$ есть просто явная формула. Точнее, важно то, что есть явная формула[27] для оптимального решения $y^*(x)$.[28] К сожалению, имеется проблема вхождения в оценку необходимого числа итераций неизвестного размера решения $x_*$ задачи минимизации $f(x)$. Эта про-блема решаема [37]. В частности, в случае, когда $G(y)$ имеет ограниченную вариацию на множестве $S_m(R_y)$ (для энтропии эта вариация равна $R_y \ln m$), можно предложить метод, с оценкой числа итераций $O\left(\varepsilon^{-1}\ln\left(\varepsilon^{-1}\right)\right)$. В эту оценку уже никак не входит неизвестный размер решения $x_*$, который может оказаться большим [37]. Далее мы еще вернемся к вопросу о том, как действовать, в случае, когда тот или иной параметр задачи (в данном случае размер решения) априорно не известен.

---

[27] Сложность формулы оценивается числом ненулевых элементов в матрице $B$. При этом считаем, что градиент $G(y)$ рассчитывается быстрее, чем занимает умножение матрицы $B$ на столбец / строку.

[28] Отметим, что если $G(y)$ – сепарабельная вогнутая функция (но не обязательно сильно вогнутая) и вместо ограничения $\|y\|_2 \le R_y$ задано сепарабельное ограничение (например, $\|y\|_\infty \le R_y$), то $\varepsilon$-приближенный поиск $y^*(x)$ можно осуществить за $O\Big(\ln\big(\varepsilon^{-1}\big)\Big)$ умножений матрицы $B$ на столбец, решая соответствующие одномерные задачи. Немного более громоздкие рассуждения [138] позволяют и при наличии ограничения $\|y\|_2 \le R_y$ осуществить $\varepsilon$-приближенный поиск $y^*(x)$ также за $O\Big(\ln\big(\varepsilon^{-1}\big)\Big)$ умножений матрицы $B$ на столбец. К сожалению, отсутствие сильной вогнутости не позволяет использовать в том же виде концепцию $(\delta, L, \mu)$-оракула для внешней задачи, однако можно при этом использовать концепцию $\delta$-субградиента для внешней задачи [99]. Это приводит лишь к оценкам $O\big(\varepsilon^{-2}\big)$, которые уже не будут оптимальными (улучшаемы до $O\big(\varepsilon^{-1}\big)$).

Отметим, что сильной вогнутости можно добиться и искусственно [271], глава 3 [181], [6, 37] (см. также разделы 3.1, 3.2 главы 3). Подход отмеченных работ приводит к оптимальным для такого класса задач оценкам (с точностью до логарифмического фактора[29] $\ln\left(\varepsilon^{-1}\right)$), и позволяет, на самом деле, контролировать точность решения одновременно по $x$ и по $y$ без использования прямо-двойственности в классическом варианте (см. ниже), что может быть полезным в определенных ситуациях [138]. Здесь под оптимальными методами мы имеем в виду методы с проксимальным оракулом. Однако в ряде задач оптимизации огромных размеров оказывается эффективнее использовать линейный минимизационный оракул [173], пришедший из метода Франк–Вульфа (см., например, п. 3.3 [163] и раздел 1.5 главы 1). Грубо говоря, суть подхода в том, что сначала вычисляется не $f(x)$ согласно модели, описанной в примере 2.1.4, а в седловом представлении задачи меняется порядок взятия максимума и минимума, и вычисляется с помощью линейного минимизационного оракула сначала минимум по $x$. Причем это не обязательно делать точно (см. подраздел 5.1.5 раздел 5.1 главы 5 [99]). Получающаяся задача максимизации по $y$ уже не будет гладкой, поэтому с учетом сильной вогнутости $G(y)$ здесь можно рассчитывать только на зависимость числа итераций от желаемой точности $O\left(\varepsilon^{-1}\right)$. Получается вроде как хуже, чем раньше. Но тут надо учитывать, как входят размерности $n$ и $m$, которые могут быть огромными в приложениях, см. п. 3.3 [163], [33, 40, 42, 173]. Удивительным образом, в сложность внутренней задачи при таком подходе (минимизации по $x$) при определенной структуре (как правило, связанной с ограничениями симплексного типа и матрицей $B$, имеющей комбинаторую [173] или сетевую [40] природу, как в разделе 1.5 главы 1) может не входить размерность вектора $x$ (т.е. $n$), что позволяет решать задачи колоссальных размеров по $n$.

[29] Ниже мы обсудим, как можно избавиться от этого логарифмического фактора для задач с явной формулой для $y^*(x)$, например, для задач энтропийно-линейного программирования [37, 138].

Пример 2.1.4 был приведен, прежде всего, потому что он поясняет одно интересное и достаточно современное направление в численных методах выпуклой оптимизации (см., например, [6, 42, 43]). Грубо говоря, это направление можно охарактеризовать, как попытку ввести оптимальную “алгебру” над алгоритмами выпуклой оптимизации. А именно, если требуется оптимизировать функционал (искать седловую точку), который обладает разными свойствами (гладкости, сильной выпуклости, быстроты вычислимости частных производных и т.п.) по разным группам переменных (такие задачи часто в последнее время возникают в разных приложениях, в частности, в транспортных и экономических [4, 9, 21, 29, 35, 40, 42, 43, 47, 275, 279]) и(или) сам представляет собой некоторую суперпозицию других функционалов (с разными свойствами; наиболее популярен случай суммы двух функционалов [33, 45, 92, 138, 259]), то хотелось бы получить такую декомпозицию исходной задачи, чтобы правильное сочетание (правильное чередование с правильными частотами) оптимальных методов для получившихся отдельных подзадач позволило бы получить оптимальный метод для исходной задачи. В ряде интересных случаев такое оказывается возможным (с оговоркой, что оптимальность понимается с точностью до логарифмического фактора). По-видимому, новым в этом направлении является наблюдение, отмеченное в примере 2.1.4 (см. также [6, 42, 43], раздел 1.6 главы 1, раздел 2.3 этой главы 2 и раздел 5.1 главы 5), что при определенных условиях идея оптимального сочетания различных методов для решения одной сложной по структуре задачи оптимизации, может быть реализована на основе концепции неточного оракула.

Другой способ борьбы с дополнительными ограничениями типа равенств или неравенств в задачах выпуклой оптимизации базируется на прямо-двойственной структуре [269] всех обсуждаемых методов (поскольку они строят модель функции [92]).[30] Это означает, что ограничения вносятся во вспомогательную задачу оптимизации, возникающую на каждом шаге метода и отвечающую за проектирование. В результате на каждом шаге решается более сложная задача. Тем не менее, если такие вспомогательные задачи можно эффективно решать (что в общем случае также наталкивается на сложности, описанные ранее) с помощью метода множителей Лагранжа (найдя и сами множители) или когда у исходной задачи есть модель (см. пример 2.1.4 и [4, 6, 33, 37, 40, 41, 42, 138, 261, 263, 271]), то тогда описанные методы позволяют не только эффективно решать исходную за-дачу оптимизации с ограничениями, но и находить попутно (по явно выписываемым формулам) решение двойственной задачи.

---

[30] Существенно подробнее собранные далее сюжеты изложены в главе 3 данной диссертации.

Основная идея работы [269] состоит в том (здесь мы ограничимся рассмотрением детерминированного случая с точным оракулом, выдающим градиент; в стохастическом случае с неточным оракулом см., например, [138] и лемму 7.7 [181]), что метод генерирует в прямом пространстве на итерациях такую последовательность $\{x_k\}$,[31] что зазор двойственности (duality gap) $\Delta(\lambda, x; N)$ удовлетворяет условию

$$\Delta(\lambda, x; N) = \max_{u \in Q} \left\{ \frac{1}{S_N} \sum_{k=0}^{N} \lambda_k \left\langle \nabla f(x_k), x_k - u \right\rangle \right\} \le \varepsilon,$$

где $S_N = \sum_{k=0}^{N} \lambda_k$, $\lambda_k \ge 0$, поэтому

$$f\left( \frac{1}{S_N} \sum_{k=0}^{N} \lambda_k x_k \right) - \min_{x \in Q} f(x) \le \varepsilon.$$

Это следует из выкладки

$$f\left( \frac{1}{S_N} \sum_{k=0}^{N} \lambda_k x_k \right) - f(u) \le \frac{1}{S_N} \sum_{k=0}^{N} \lambda_k \cdot \left( f(x_k) - f(u) \right) \le \frac{1}{S_N} \sum_{k=0}^{N} \lambda_k \left\langle \nabla f(x_k), x_k - u \right\rangle.$$

Аналогичную точность (для двойственной задачи) дает следующая аппроксимация решения двойственной задачи

$$y = \frac{1}{S_N} \sum_{k=0}^{N} \lambda_k y_k.$$

Это сразу следует из того, что зазор двойственности оценивает сверху разность между получившимися значениями целевой функции в прямой задаче и двойственной [269], которую мы будем называть двойственным зазором. Эта разность всегда неотрицательна, и на точных решениях прямой и двойственной задачи (и только на них) равна нулю. Заметим, что контроль онлайн-зазора двойственности

---

[31] В двойственном пространстве при этом генерируется последовательность соответствующих множителей Лагранжа $\{y_k\}$ [269] или, в случае наличия модели у исходной прямой задачи (см. пример 2.1.4), последовательность $\{y_k\}$ генерируется по явным или расчетным формулам $\{y_k = y(x_k)\}$ согласно этой модели [4, 6, 33, 138, 261, 263]. Такой подход также позволяет убрать логарифмический фактор в задачах энтропийно-линейного программирования [37, 138] и аналогичных задачах, см. п. 5.2 [261] и [6, 45, 138, 171].

$$\Delta(\lambda, x; N) = \max_{u \in Q} \left\{ \frac{1}{S_N} \sum_{k=0}^{N} \lambda_k \left\langle \nabla f_k(x_k), x_k - u \right\rangle \right\}$$

позволяет в случае, когда удается выбрать $\lambda_k \equiv 1$, получать оценки регрета (псевдо регрета в стохастическом случае) в задачах онлайн оптимизации (см. [166, 320] и конец этого пункта). К сожалению, ограничение $\lambda_k \equiv 1$ существенно сужает класс методов. Скажем, для рассматриваемых в этом пункте быстрых градиентных методов $\lambda_k \sim k^p$. Кроме того, даже если в онлайн постановке допустить взвешивание с различными весами, все равно требуется, чтобы способ получения оценки на зазор двойственности допускал бы обобщение на онлайн-постановки. Быстрый градиентный метод, например, этого не допускает, что не сложно усмотреть из оценок работы [135].

Описанная выше конструкция, основанная на оценке зазора двойственности, работает в случае ограниченного множества $Q$. В случае неограниченного $Q$ (это типичная ситуация, когда необходимо решать двойственную задачу, по решению которой требуется восстанавливать решение прямой задачи) можно искусственно компактифицировать $Q$ [40, 41, 138, 261]. Однако, в большинстве случаев такая компактификация не позволяет очевидным образом оценивать настоящий (не обрезанный) двойственный зазор в исходной задаче, что часто представляется важным ввиду наличия простых явных формул для этого настоящего зазора, и возможности использования контроля зазора двойственности у качестве критерия останова метода. Несмотря на отмеченную теоретическую проблему, на практике проблема оказывается решаемой [33, 40, 138].

Более общий способ оценки разности между получившимися значениями целевых функций в прямой задаче и двойственной базируется на контроле сертификата точности (accuracy certificate) [261, 281], в который наряду с градиентами функционала входят градиенты нарушенных ограничений или в общем случае вектора нормалей к гиперплоскостям, отделяющим $x_k$ от множества $Q$ (в ряде постановок “градиенты” стоит заменить на “субградиенты”). Вектора двойственных множителей формируются из соответствующих (сертификату точности) взвешенных сумм векторов нормалей отделяющих гиперплоскостей [34, 261, 281]. Собственно, такая интерпретация двойственных множителей следует из способа обоснования принципа множителей Лагранжа на основе следствия теоремы Хана–Банаха (теоремы об отделимости) [85]. Причем в работе [34, 281] за счет слейтеровской релаксации ограничений (допущения возможности нарушения ограничений на $\varepsilon$ [37, 259]) получаются оценки скорости сходимости, не зависящие от размера двойственного решения, который может быть большим.

Во многих (транспортно-)экономических приложениях при поиске равновесных конфигураций (см. главу 1) требуется решать пару задач (прямую и двойственную), см., например, [33, 37, 40, 42, 43, 47, 244, 245, 245, 257, 257]. Причем интересны решения обеих задач (решения этих задач имеют содержательную интерпретацию и используются при принятии решений / управлении). Если у этой пары задач, на которую можно смотреть, как на одну седловую задачу, есть определенная структура (проявляющаяся, например, в сильной выпуклости функционала по части переменных, наличии эффективно вычислимого линейного минимизационного оракула и т.п.), то описанный выше формализм позволяет развить идею примера 2.1.4 таким образом, чтобы одновременно (без дополнительных затрат) получать решения обеих задач. Даже в случае огромного размера одной из этих задач, можно надеяться (при эффективном линейном минимизационном оракуле), что эта размерность не войдет в сложность поиска решения прямой и двойственной задачи [244, 245, 246, 252].

В действительности, выбранный в данном разделе класс проекционных методов с построением модели функции далеко не единственный возможный способ строить прямо-двойственные методы. Скажем, уже упоминавшиеся методы условного градиента также являются прямо-двойственными [220, 263]. Еще более удивительным может показаться, что прямо-двойственная интерпретация есть, например, у метода эллипсоидов [261] (и мы это будем существенно использовать в разделе 2.3 данной главы 2). Более того, в ряде ситуаций мы можем за линейное время (с геометрической скоростью сходимости) находить одновременно решение прямой и двойственной задачи. Причем речь идет не только о конструкциях типа [331], базирующихся на принципе (см. также замечание 2.1.6): сопряженная функция к выпуклой функции с липшицевым градиентом – сильно выпуклая, и обратно, сопряженная к сильно выпуклой функции – выпуклая функция с липшицевым градиентом; но и о более общем контексте [6, 45, 138, 261, 331].

Возвращаясь к сказанному выше в связи с оценками (2.1.7) – (2.1.10) интересно заметить, что если множество $Q \subset \mathbb{R}^n$ есть шар $B_q^n(R)$ радиуса $R$ в $q$-норме,[32] то нижние оценки (для случая $D = \delta = 0$) на точность (по функции), которую можно получить после $N \le n$ итераций, имеют вид [212] (считаем, что[33] $\left\|\nabla f(x) - \nabla f(y)\right\|_{q'} \le L_\nu \left\|x - y\right\|_q^\nu$, $1/q + 1/q' = 1$, $\nu \in (0,1]$):

$$\Omega\left(\frac{1}{\min\{q, \ln n\}^\nu} \frac{L_\nu R^{1+\nu}}{N^{\nu + (\nu+1)/q}}\right) \ (2 \le q \le \infty),$$

$$\Omega\left(\frac{1}{\ln^\nu (N+1)} \frac{L_\nu R^{1+\nu}}{N^{\nu + (\nu+1)/2}}\right) \ (1 \le q < 2).$$

---

[32] Если о выпуклом замкнутом множестве $Q$ известно только то, что оно содержит $B_q^n(R)$, то все сказанное далее также остается в силе.

Приведенный результат хорошо соответствует тому, что написано в замечании 2.1.2 (см. также [90]).

Для $q = \infty$ и $\nu = 1$ (гладкий случай) приведенная оценка с точностью до логарифмического фактора будет иметь вид $\Omega\left(LR^2/N\right)$. Эта оценка достигается, например, на методе условного градиента Франк–Вульфа [214, 220].[34] Исходя из только что написанного и

---

[33] В случае достаточной гладкости функции $f(x)$ можно выписать следующее представление для константы Липшица градиента (верхний индекс $q$ соответствует выбору нормы в прямом пространстве)

$$L^q = \max_{x \in Q, \|h\|_q \le 1} \left\langle h, \nabla^2 f(x) h \right\rangle.$$

В частности, $L^1 \le L^2 \le nL^1$, $L^2 \le L^\infty \le nL^2$. Эти формулы вместе со сказанным ранее относительно того, как может меняться $M$ – обычная константа Липшица $f(x)$, при изменении нормы в прямом пространстве, поясняют почему в “устойчивые сочетания” эти константы входят таким образом: $M^2R^2/\varepsilon^2$, $LR^2/\varepsilon$. Если ввести “физические размерности”, скажем, считать, что $f(x)$ это рубли (*руб*), а $x$ это килограммы (*кг*), то $\varepsilon$ [ *руб* ], $R$ [ *кг* ], $M$ [ *руб/кг* ], $L$ [ *руб/кг*$^2$ ]. Поскольку число итераций $N$ должно быть безразмерной величиной, то возникновение агрегатов $M^2R^2/\varepsilon^2$, $LR^2/\varepsilon$ вполне закономерно. Аналогичные рассуждения можно провести и для оценок в сильно выпуклом случае. Все это приводит к довольно интересным следствиям [6]. Например, что шаг метода в негладком случае $h \sim \varepsilon/M^2$, в гладком случае $h$ определяется из соотношения вида ($W(\ )$, $\tilde{W}(\ )$ – какие-то функции)

$$W\left(h\frac{M^2}{\varepsilon}, hL\right) = 1.$$

В стохастическом случае (вместо градиента получаем стохастический градиент с дисперсией $\sigma^2$) из

$$\tilde{W}\left(h\frac{M^2}{\varepsilon}, hL, h\frac{\sigma}{R}\right) = 1.$$

[34] Отметим также, что этот метод допускает обобщение на случай неточного оракула, и неулучшаемость оценок может быть проинтерпретирована с точки зрения сохранения свойства разреженности решения [220]. Это неудивительно, поскольку аналогичный метод (с линейным минимизационным оракулом, см.

тезиса о неулучшаемости оценок (2.1.7) ((2.1.8), (2.1.9)) (при $D=\delta=\mu=0$, $p=1$) может возникнуть ощущение противоречия. Это ощущение дополнительно усиливается примером 2.1.2 из подраздела 2.1.2. Действительно, исходя из этого примера, может сложиться ощущение, что проблема выбора прокс-структуры в задаче не очень актуальна, поскольку можно исходить просто из самой нормы. И это действительно так, если мы ограничиваемся не ускоренными градиентными методами (PGM, метод Франк–Вульфа), которые сходятся как $O\left(L^q R_q^2/N\right)$ (здесь $R_q=R$ – диаметр множества $Q$, посчитанный в $q$-норме, в нашем случае $q=\infty$). Если же мы хотим ускориться, и достичь оптимальной оценки $O\left(LR^2/N^2\right)$, то уже необходимо существенно использовать прокс-функцию $d(x)\ge 0$ со свойством сильной выпуклости относительно выбранной нормы и с константой сильной выпуклости $\alpha\ge 1$ [92]. Скажем (в связи с примером 2.1.2), квадрат 1-нормы – не есть сильно выпуклая функция относительно 1-нормы, т.е. $d(x)=\|x\|_1^2$ – не может быть прокс-функцией при выборе 1-нормы для симплекса. В классе удовлетворяющих условию 1-сильной выпуклости прокс-функций (относительно выбранной нормы) подбирается такая, которая минимизирует $R^2:=\max_{x\in Q} d(x)$. Именно это $R^2$ входит в оценку FGM $O\left(LR^2/N^2\right)$. И как уже отмечалось (см. подраздел 2.1.2) для $Q=B_\infty^n(R)$ имеет место следующая оценка на прокс-диаметр $R^2:=\max_{x\in Q} d(x)=R^2\Omega(n)$. Отсюда, с учетом того, что $N\le n$, получаем, что оценка $O\left(LR^2/N\right)$ и оценка $O\left(LR^2 n/N^2\right)$, приводят, в общем-то, к одному результату, но в случае использования FGM требуется дополнительно искать оптимальную прокс-структуру. Только в таком случае будет совпадение результатов. Более того, также как и в замечании 2, здесь хорошо видно, что при $q\ge 2$ можно ограничиться рассмотрением только евклидовой прокс-структуры для FGM и евклидовой нормы для метода Франк–Вульфа. В частности, для $Q=B_\infty^n(R)$ действуя так, мы получим для FGM оценку $O\left(L^2R^2 n/N^2\right)$ вместо ранее полученной оценки $O\left(L^\infty R^2 n/N^2\right)$, соответствующей при $N\le n$ неулучшаемой оценке $O\left(L^\infty R^2/N\right)$ (здесь мы проставили верхние индексы у $L$,

---

пример 2.1.4) с аналогичными оценками скорости можно получить (см. п. 5.5.1 [259], [278]) из композитного варианта FGM в концепции неточного оракула (суть метода в том, что в композитном варианте FGM на каждой итерации решается задача, в которой коэффициент при прокс-слагаемым равен нулю, т.е. оно просто отсутствует).

поскольку они различаются). Дальше можно написать все тоже по поводу неулучшаемости оценок, что и в конце замечания 2.1.2.

Отметим также, что параметры $R$ и $\mu$ могут быть не известны априорно или процедуры их оценивания приводят к слишком (соответственно) завышенным и заниженным результатам. Это может быть проблемой, поскольку в ряде случаев знание этих и других параметров требуется методу для расчета величин шагов и условий остановки. Из этой ситуации можно выйти за логарифмическое (по этим параметрам) число рестартов метода. Стартуя, скажем, с $R = 1$ и делая число шагов, вычисленное из оценки скорости сходимости при выбранном $R$, мы проверяем выполняется ли для вектора, выдаваемого алгоритмом, условие $\varepsilon$-близости по функции (при условии, что мы можем сделать такую проверку). Если условие $\varepsilon$-близости не выполняется, то полагаем $R := 2R$ и т.д. Все эти перезапуски увеличат общее число обращений к оракулу лишь в $O(1)$ раз [92, 227, 228].[35] Аналогичное можно сказать про[36] $L$, $M$, $\mu$ и $D$. Однако если убрать стохастичность, тогда $L$, $M$ можно не только эффективнее адаптивно подбирать (аналогично правилу Армихо [22] в независимости от того можем ли мы сделать проверку условия $\varepsilon$-близости значения функции в текущей точке к оптимальному) по ходу самих итераций (увеличив в среднем число обращений к оракулу не более чем в 4 раза), но и в некотором смысле оптимально самонастраиваться (используя формулу (2.1.10)) на гладкость функционала на текущем участке пребывания метода [274]. Это означает, что в детерминированном случае без учета сильной выпуклости функционала существует универсальный метод, работающий по оценкам (2.1.8) с $L$, рассчитанной по формуле (2.1.10) (в которой $\delta$ берется из (2.1.9)), и оптимальным в смысле скорости сходимости выбором параметра $\nu \in [0,1]$. Причем выбор $\nu$ осуществляется не нами заранее, исходя из знания всех констант и минимизации выписанных оценок, а самим алгоритмом (здесь выбрано $p=1$):

$$N(\varepsilon)=\inf_{\nu\in[0,1]}\left(\frac{2^{\frac{3+5\nu}{2}}L_\nu R^{1+\nu}}{\varepsilon}\right)^{\frac{2}{1+3\nu}}.$$

---

[35] В свою очередь можно поиграть и на этом $O(1)$, стараясь его минимизировать. Для этого шаг, который мы для простоты положили равным 2, подбирают оптимально исходя из того, с каким показателем степени входит неизвестный (прогоняемый) параметр в оценку числа итераций [2, 37]. Подробнее об этом будет написано в разделе 3.1 главы 3 и в разделе 4.3 главы 4.

[36] Впрочем, в детерминированных постановках мы можем явно наблюдать за последовательностью выдаваемых оракулом субградиентов и отслеживать условие на норму субградиентов. Как только наше предположение нарушилось (при этом мы не успели сделать предписанное текущему $M$ число шагов), мы увеличиваем $M$ в два раза и перезапускаем весь процесс с новым значением $M$. Число таких перезапусков будет не более чем логарифмическим от истинного значения $M$. Все эти рассуждения с небольшими оговорками (типа равномерной п.н. ограниченности стохастического субградиента) переносятся и на стохастические постановки, в которых наблюдается стохастический субградиент. Для определенного класса задач, в которые неизвестные параметры входят только в критерий останова, но не в сам метод (к таким задачам, например, относится задача поиска равновесия в модели Бэкмнана методом Франк–Вульфа и неизвестной константе $L$) можно обходиться и без перезапусков [40]. Заметим также, что у ряда популярных методов (например, метода зеркального спуска) есть варианты, в которые входит не оценка супремума нормы субградиента (или градиента), а норма субградиента на текущей итерации, которая известна [4, 40, 259, 269].

Это соответствует (с точностью до логарифмического фактора) нижним оценкам [261], выписанным выше для случая $q\in[1,2]$. Отметим, что здесь при определении $R$ используется соответствующая прокс-функция, см. замечание 2.1.2. К сожалению, пока не очень понятно можно ли что-то похожее сделать с параметром $\mu$ и с введенным нами в начале этого пункта параметром метода $p\in[0,1]$. Обзор других работ на тему самонастройки алгоритмов в гладком детерминированном случае имеется в [41, 48, 283], а в стохастическом случае в [53] и следующем разделе.

Приведенную оценку можно обобщить (см. следующий раздел 2.2), если дополнительно известно, что функция $f(x)$ – $\mu$-сильно выпукла. Можно дополнительно к искусственно введенной игре на неточности оракула допустить, что имеет место и настоящая неточность. В этом случае также можно выписать соответствующие оценки [42]. Не играя на выборе $\nu\in[0,1]$, можно распространить все, что описано выше в этом абзаце на стохастические постановки. Аналогичное можно сделать для стохастических безградиентных и покомпонентных методов с неточным оракулом (см. подраздел 2.1.4 и [45]). Соответствующие обобщающие формулы собраны в работе [41, 48, 53], мы не будем их здесь приводить (см. следующий раздел). Такие обобщения востребованы, например, в связи с приложениями к поиску равновесий в многостадийных моделях равновесного распределения транспортных потоков [9, 29, 33, 35, 42, 43, 47] (см. также раздел 1.6 главы 1). В основе этих приложений лежит конструкция, изложенная в примере 2.1.4, с универсальным методом [274] вместо FGM для решения внешней задачи.

Выше мы сделали обременительное предположение о возможности выполнять проверку условия $\varepsilon$-близости по функции. Такое заведомо возможно только при известном значении функционала в точке оптимума. Как правило, такой информации у нас априорно нет. Один из способов выхода из этой ситуации для задач стохастической оптимизации описан в п. 7.7 работы [181]. Другой способ – контролировать зазор двойственности (со стохастическими градиентами). Для применимости этого способа требуется, чтобы числовая последовательность $\{\lambda_k/S_N\}_{k=0}^N$ не зависела от неизвестных параметров. Во многих задачах, приходящих из транспортных и экономических приложений, нужно одновременно находить решения прямой и двойственной задачи, которые можно явно выписать. В таких случаях имеется эффективный способ проверки условия $\varepsilon$-близости по функции. Нужно проверить условие $\varepsilon$-малости разницы между полученным (приближенным) значением функционалов прямой и двойственной задачи (т.е. двойственного зазора) [33, 35, 40, 42].

Отметим также, что в детерминированном $\mu$-сильно выпуклом случае, когда в точке минимума $x_*$ выполняется условие[37] $\nabla f(x_*) = 0$, критерий $\varepsilon$-близости по функции может быть переписан в терминах малости рассчитываемого на итерациях градиента:

$$f(x_k) - f_* \le \frac{1}{2\mu}\|\nabla f(x_k)\|_*^2 \le \varepsilon.$$

В постановках с сильно выпуклой/вогнутой двойственной задачей (этого можно добиться искусственно, вводя регуляризацию в двойственную задачу, см. главу 3 [181], [6, 37], а также разделы 3.1, 3.2 главы 3) также можно оценивать точность решения прямой задачи по точности решения двойственной задачи, применяя к двойственной задаче неравенство[38]

$$\frac{1}{2L}\|\nabla f(x_k)\|_*^2 \le f(x_k) - f_*.$$

В частности, это обстоятельство используется в критерии остановки двойственного метода из [37], см. также пример 2.1.4 ниже, [6, 138] и раздел 3.1 главы 3.

---

[37] От этого условия можно избавиться, используя в приведенных далее формулах вместо градиента градиентное отображение [6, 93].

[38] В замечании 2.1.5 (см. также [6]) был приведен пример, когда $\|\nabla f(x_k)\|_2 = O(k^{-2})$. Из выписанного неравенства мы можем гарантировать лишь $\|\nabla f(x_k)\|_2^2 = O(k^{-2})$. Ситуацию можно улучшить, если регуляризовать функционал (см. конец подраздела 2.1.2 и [6, 37]), сделав его сильно выпуклым, и применить FGM [92, 93, 138, 267] к регуляризованной задаче, тогда $\|\nabla f(x_k)\|_2 = O\left((\ln k)^2/k^2\right)$ (если использовать, например, FGM с оценкой числа итераций $O\left(\sqrt{L/\mu}\left\lceil \ln\left(\mu R^2/\varepsilon\right)\right\rceil\right)$, $\mu \sim \varepsilon/R^2$ ). В негладком случае ситуация проще, см. [269].

Все сказанное выше, по-видимому, переносится в полной мере на задачи композитной оптимизации [92, 183, 259, 266] и некоторые их обобщения, см., например, работы

А.С. Немировского (http://www2.isye.gatech.edu/~nemirovs/ ),

А.Б. Юдицкого (http://ljk.imag.fr/membres/Anatoli.Iouditski/ ).

**Замечание 2.1.6.** Композитные задачи имеют вид: $f(x)+\lambda h(x)\to\min\limits_{x\in Q}$, где $\lambda>0$, $h(x)$ – выпуклая функция простой структуры, скажем $h(x)=\|x\|_1$. Хочется, чтобы сложность решения этой задачи всецело определялась только гладкостью выпуклого функционала $f(x)$, а сильная выпуклость – обоими слагаемыми. Если не лианеризовывать функцию $h(x)$ при подсчете на каждой итерации градиентного отображения [92], а просто оставлять это слагаемое как есть, то, конечно, сложность решения вспомогательной задачи на каждой итерации увеличится (впрочем, в виду простой структуры функции $h(x)$, ожидается, что не намного), зато в оценку необходимого числа итераций уже не будут входить никакие константы, характеризующие гладкость $h(x)$, только константы, характеризующие сильную выпуклость (если имеется). На такие задачи также можно смотреть следующим образом (принцип множителей Лагранжа): $f(x)\to\min\limits_{x\in Q,\, h(x)\le C(\lambda)}$. Поскольку функция $h(x)$ простой структуры, то проектироваться на множество Лебега этой функции несложно. Отсюда можно усмотреть независимость числа итераций от $h(x)$. Другой способ “борьбы” с композитным членом $h(x)$ (А.С. Немировский) заключается в переписывании задачи в “раздутом” (на одно измерение) пространстве: $f(x)+y\to\min\limits_{x\in Q,\, h(x)\le y}$. Норма в раздутом пространстве задается как $\|(x,y)\|=\|x\|+\alpha|y|$. Функционал имеет такой вид, что в независимости от гладкости $f(x)$ в оценки супремума нормы субградиента / константы Липшица градиента не будет входить что-либо, связанное с $y$. В гладком случае все ясно сразу из определения, а в случае негладкой $f(x)$ это связано с тем, что, в действительности, в оценку необходимого числа итераций входит не супремум нормы субградиента, а супремум нормы разностей субградиентов [268] (см. также начало подраздела 2.1.2). За счет возможности выбирать сколь угодно маленьким $\alpha>0$, можно считать независящим от $y$ и прокс-расстояние (от точки старта до решения раздутой задачи), входящее в оценку необходимого числа итераций. Таким образом, можно сделать оценку числа итераций независящей от $y$ и $h(x)$. В связи с написанным выше полезно заметить, что гладкость (липшицевость градиента) и сильная выпуклость функционала являются взаимодвойственными друг к другу для задач безусловной оптимизации (константа Липшица градиента переходит в константу сильной выпуклости и наоборот, отсюда, кстати сказать, можно усмотреть, что оценка скорости сходимости для таких задач должны зависеть от отношения $L/\mu$, другой способ понять это – соображения "физической" размерности), что активно используется в приложениях, см., например, [45, 138, 307, 331]. Однако для задач условной оптимизации остается только один переход: двойственная (сопряженная) задача к сильно выпуклой – гладкая (см., например, [271] и замечание 2.1.3), обратное не верно даже в случае сильно выпуклой функции $h(x)$. Собственно, мы уже сталкивались с "неравноправностью" гладкости и сильной выпуклости. При рассмотрении универсального метода, мы отмечали, что на гладкость можно настраиваться адаптивно, чего нельзя сказать про сильную выпуклость. Замечание 2.1.6 немного проясняет (с учетом Лагранжева формализма) соотношения между этими свойствами задачи. Впрочем, до окончательного понимания, к сожалению, сейчас еще довольно далеко. Не ясно даже принципиальны ли эти различия или их можно в перспективе устранить. По-видимому, принципиальны, но строгого обоснования мы здесь не имеем.

Также нам видится, что сказанное выше переносится на седловые задачи и монотонные вариационные неравенства [91, 92]. Причем, речь идет не о том, что было описано в примере 2.1.4, а о том, как скажется неточность оракула на оптимальные методы для седловых задач и монотонных вариационных неравенств [91, 92]. Ответ, более менее, известен: неточность оракула не будет накапливаться на оптимальных методах (в отличие от задач обычной выпуклой оптимизации). Отметим, что концепцию неточного оракула еще необходимо должным образом определить[39] – предположение 2.1.1 нуждается в корректировке для данного класса задач. Отсутствие накопления неточностей связано с тем, что для таких задач оценка (2.1.7) будет оптимальна (с некоторыми оговорками) при $p=0$.

---

[39] Например, для монотонного вариационного неравенства: найти такой $x\in Q$, что для всех $y\in Q$ выполняется $\langle g(y), y-x\rangle \ge 0$, достаточные условия на $(\delta, L)$-оракул будут иметь вид: для любых $x,y\in Q$

$$\langle g(y)-g(x), y-x\rangle \ge -\delta,\ \|g(y)-g(x)\|_* \le L\|x-y\|+\delta .$$

Вероятно, в ряде ситуаций эти условия можно ослабить (доказательства в этом случае нам не известны)

$$-\delta \le \langle g(y)-g(x), y-x\rangle \le L\|y-x\|^2+\delta .$$

Другие $p \in (0,1]$ рассматривать не стоит (такие оценки просто не достижимы). Впрочем, пока про это имеются лишь частичные результаты, напрмимер в работах: А.С. Немировского, А.Б. Юдицкого, G. Lan-а.

В оценку числа итераций для достижения заданной точности решения в описанных методах не входит явно размерность пространства $n$. Это наталкивает на мысль (подобные мысли, по-видимому, впервые были высказаны и реализованы для класса обычных градиентных методов в кандидатской диссертации Б.Т. Поляка [100], см. также [83, 99]) о возможности использовать эти методы, например, в гильбертовых пространствах [22, 64]. Оказывается это, действительно, можно делать при определенных условиях (см., например, [41] в контексте использованных в данном разделе обозначений). В частности, концепция неточного оракула позволяет привнести сюда элемент новизны, существенно мотивированный практическими нуждами – принципиальной невозможностью (в типичных случаях нет явных формул) решать с абсолютной (очень хорошей) точностью вспомогательную задачу на каждом шаге градиентного спуска. Например, решение такой вспомогательной задачи для класса задач оптимального управления со свободным правым концом приводит к двум начальным задачам Коши для систем обыкновенных дифференциальных уравнений (важно, чтобы СОДУ для фазовых переменных и сопряженных решались, скажем, методом Эйлера, на одной и той же сетке), которые необходимо решить для вычисления градиента функционала [22]. Однако, в действительности, почти все практически интересные задачи (за редким исключением, к коим можно отнести класс Ляпуновских задач [85]) в бесконечномерных пространствах не являются выпуклыми, поэтому здесь имеет смысл говорить лишь о поиске локальных минимумов (локальной теории) [55]. Если ограничиться неускоренными методами (например, PGM), то можно показать, что при весьма общих условиях эти методы могут быть использованы в гильбертовом пространстве в концепции неточного оракула и для невыпуклых (но гладких) функционалов, причем с аналогичными оценками скорости сходимости (отличие от выпуклого случая будет в том, что метод сходится лишь к стационарной точке (локальному экстремуму), в бассейне притяжения которой окажется точка старта). Заметим, что задачи оптимального управления можно численно решать, построив соответствующую (аппроксимирующую) задачу оптимального управления с дискретным временем, что приводит к конечномерным задачам, для решения которых можно использовать конечномерный вариант PGM в невыпуклом случае (с точным оракулом). Этот путь, как правило, и предлагается в большинстве пособий (см., например, [22]). Однако при таком подходе мы должны уметь (по возможности точно) решать сложную задачу оценки качества аппроксимации исходной задачи оптимального управления ее дискретным по времени вариантом. Более теоретически обоснованный способ рассуждений, по сути, приводящий к необходимости решать все те же конечномерные задачи, заключается в рассмотрении исходной задачи оптимального управления и ее решения бесконечномерным вариантом PGM в невыпуклом случае (с неточным оракулом). Неточность оракула существенна. Поскольку на каждой итерации этого градиентного метода необходимо решать две задачи Коши для СОДУ, что в общем случае можно сделать лишь приближенно, но с лучшим контролем точности, чем при подходе с дискретизацией задачи оптимального управления. Отметим, что во многих “физических” приложениях схема Эйлера имеет хорошие теоретические свойства сходимости (устойчивости). Связано это с тем, что на оптимальном режиме, как правило, наблюдается некоторая стабилизация поведения системы управления, что приводит к устойчивости якобиана прямой и обратной системы дифференциальных уравнений. Исследования в этом направлению проводит аспирант МФТИ Ахмед Мухамед.

Покажем в заключение этого пункта, как приведенные результаты переносятся на задачи стохастической онлайн оптимизации. Для этого напомним вкратце, в чем состоит постановка задачи (см., например, [50, 52, 54, 166, 216, 250, 259, 299, 306, 320], подробнее об этом написано в главе 6). Требуется подобрать последовательность[40] $\left\{x^k\right\}\in Q$ так, чтобы минимизировать псевдо регрет:

$$\frac{1}{N}\sum_{k=1}^{N}E_{\xi^k}\left[f_k\left(x^k,\xi^k\right)\right]-\min_{x\in Q}\frac{1}{N}\sum_{k=1}^{N}E_{\xi^k}\left[f_k\left(x,\xi^k\right)\right]$$

на основе доступной информации $\left\{\nabla\tilde{f}_1\left(x^1,\xi^1\right);\dots;\nabla\tilde{f}_{k-1}\left(x^{k-1},\xi^{k-1}\right)\right\}$ при расчете $x^k$, где[41]

$$\left\|\nabla\tilde{f}_k\left(x^k,\xi^k\right)-\nabla f_k\left(x^k,\xi^k\right)\right\|_*\le\delta,$$

$$E_{\xi_k}\left[\nabla f_k\left(x,\xi^k\right)\right]=\nabla f_k\left(x\right).$$

Здесь с.в. $\left\{\xi^k\right\}$ могут считаться независимыми одинаково распределенными. Онлайновость постановки задачи допускает, что на каждом шаге $k$ функция $f_k$ может подбирать-ся из рассматриваемого класса функций враждебно по отношению к используемому нами методу генерации последовательности $\{x^k\}$. В частности, $f_k$ может зависеть от

$$\left\{x^1,\xi^1,f_1(\cdot);...;x^{k-1},\xi^{k-1},f_{k-1}(\cdot);x^k\right\}.$$

---

[40] Если $Q=\left\{x:g\left(x\right)\le 0\right\}$ и на это множество сложно проектироваться, то можно обобщить, сохранив оценку (2.1.11), конструкцию прямо-двойственного метода из работы [281] (см. также [3]) на онлайн контекст с таким множеством $Q$.

[41] Это условие можно заменить, считая, что вместо субградиента мы получаем $\delta$-субградиент [50, 99].

В стохастической онлайн оптимизации с неточным оракулом можно получить следующие оценки псевдо регрета (см., например, [216, 227], случай неточного оракула в похожем контексте ранее уже частично прорабатывался в [227])

$$\min\left\{\mathrm{O}\left(\sqrt{\frac{M^2R^2}{N}}+R\delta\right),\mathrm{O}\left(\frac{M^2\ln N}{\mu N}+R\delta\right)\right\}, \qquad (2.1.11)$$

где $\|\nabla f_k(x,\xi)\|_* \le M$ – равномерно по $x$, $k$ и п.н. по $\xi$. Эти оценки достигаются (фактически на тех же методах, что и в подразделе 2.1.2 с небольшой оговоркой в сильно выпуклом случае [216, 218]) и неулучшаемы (в том числе для детерминированных постановок с $\delta = 0$ и с линейными функциями $f_k(x)$). Как видно из этих оценок, наличие гладкости не позволяет получить более высокую скорость сходимости. То есть никакого аналога формулы (2.1.7) здесь нет. Все что ранее говорилось про прокс-структуру[42] и большие отклонения, насколько нам известно, полностью и практически без изменений переносится и на задачи онлайн оптимизации.

### 2.1.4 Стохастические безградиентные и покомпонентные методы с неточным оракулом

Рассматривается задача стохастической выпуклой оптимизации (2.1.1)

$$f(x) = E_\xi\left[f(x,\xi)\right] \to \min_{x\in Q}.$$

Предположения те же, что и в первом абзаце подраздела 2.1.2. В частности, п.н. $\|\nabla f(x,\xi)\|_2 \le M$. Здесь важно, что функция $f(x)$ задана не только на множестве $Q$, но и в его $\tau_0$-окрестности (см. ниже), и все предположения делаются не для $x \in Q$, а для $x$ из $\tau_0$-окрестности множества $Q$ (аналогичная оговорка потребуется далее, при перенесении результатов подраздела 2.1.3). Однако теперь оракул не может выдавать стохастический субградиент функции. На каждой итерации мы можем запрашивать у оракула только зна-чения реализации функции $f(x,\xi)$ в нескольких точках $x$. Принципиальная разница есть только между запросом значения (реализации) функции в одной и запросом в двух точках [41, 190]. Здесь мы ограничимся рассмотрением случая двух точек – случай одной точки, по-видимому, представляет интерес только в онлайн контексте (см. [50] и цитированную там литературу, а также раздел 6.2 главы 6). Впрочем, есть достаточно большой и популярный класс одноточечных не онлайн постановок, которого мы здесь не касаемся (см., например, [104, 314]).

[42] За исключением сильно выпуклого случая, для которого нам известны только оценки в евклидовой прокс-структуре. Кроме того, в сильно выпуклом случае в оценках вероятностей больших уклонений $\ln(\ln(N)) \to \ln N$ – доказательство этого утверждения мы не смогли найти (впрочем, см. замечание 2.1.4).

**Предположение 2.1.2.** *$\delta$-оракул выдает (на запрос, в котором указывается только одна точка $x$) $f(x,\xi)+\delta(x,\xi)$, где с.в. $\xi$ независимо разыгрывается из одного и того же распределения, фигурирующего в постановке (2.1.1), случайная величина $\delta(x,\xi)=\tilde{\delta}(x)+\bar{\delta}(\xi)$, где $\bar{\delta}(\xi)$ – независимая от $x$ случайная величина с неизвестным распределением (случайность которой может быть обусловлена не только зависимостью от $\xi$), ограниченная по модулю $\delta/2$ (число $\delta$ – допустимый уровень шума), $\tilde{\delta}(x)/(R\delta)$ – неизвестная 1-липшицева функция.*

Далее в изложении мы будем во многом следовать [50, 52, 58, 59, 99, 104, 130, 166, 190, 270, 216, 314]. По полученным от оракула зашумленным значениям $f(x,\xi)+\delta(x,\xi)$ мы можем определить стохастический субградиент (важно, что можно обратиться с запросом к оракулу в двух разных точках при одной и той же реализации $\xi$):

$$g_{\tau,\delta}(x,s,\xi)=\frac{n}{\tau}\Big(f(x+\tau s,\xi)+\delta(x+\tau s,\xi)-\big(f(x,\xi)+\delta(x,\xi)\big)\Big)s, \qquad (2.1.12)$$

где $s$ – случайный вектор (независимый от $\xi$), равномерно распределенный на $S_2^n(1)$ – единичной сфере в 2-норме в пространстве $\mathbb{R}^n$.[43] Из этого представления можно усмотреть, что липшицева составляющая шума $\delta_2$ из предположения 2.1.2 и уровень шума $\delta_1$ из предположения 2.1.1 связаны соотношением $\delta_2 \sim \delta_1/n$ (см. формулу (2.1.24)). В действительности, для обоснования этой связи требуются значительно более громоздкие рассуждения.

Приведем одну из возможных мотивировок введенной в предположении 2 концепции $\delta$-оракула. Предположим, что оракул может считать абсолютно точно значение липшицевой функции, но вынужден нам выдавать лишь конечное (предписанное) число первых бит. Таким образом, в последнем полученном бите есть некоторая неточность (причем мы не знаем по какому правилу оракул формирует этот последний выдаваемый значащий бит). Однако мы всегда можем прибавить (по mod 1) к этому биту случайно приготовленный (независимый) бит. В результате, не ограничивая общности, можно считать, что оракул последний бит выбирает просто случайно в независимости от отброшенного остатка. То что в итоге выдает оракул соответствует концепции $\delta$-оракула.

---

[43] С помощью леммы Пуанкаре [244] такой вектор можно сгенерировать за $O(n)$, разыгрывая $n$ независимых одинаково распределенных стандартных нормальных случайных величин и нормируя их на корень из суммы их квадратов.

Перейдем к получению оценок (опущенные здесь детали можно посмотреть в главе 5 и разделе 6.2 главы 6). В отличие от подразделов 2.1.2, 2.1.3 везде далее в этом пункте мы будем считать, что имеет место обратное неравенство на требуемое число итераций $N \ge n$ [190]. Прежде всего, заметим, что[44]

$$E_{s,\xi,\delta}\left[g_{\tau,\delta}(x,s,\xi)\right] = \nabla f_\tau(x) + \nabla_x E_{\tilde{s},\xi,\delta}\left[\delta(x+\tau\tilde{s},\xi)\right], \tag{2.1.13}$$

где $\tilde{s}$ – случайный вектор, равномерно распределенный на $B_2^n(1)$ – единичном шаре в 2-норме, а $f_\tau(x) = E_{\tilde{s},\xi}\left[f(x+\tau\tilde{s},\xi)\right]$ – сглаженная[45] версия функции $f(x) = E_\xi\left[f(x,\xi)\right]$. Причем,

$$0 \le f_\tau(x) - f(x) \le M\tau, \tag{2.1.14}$$

$$\left\|g_{\tau,\delta}(x,s,\xi)\right\|_2 \le n\left(M + \frac{2\delta}{\tau}\right). \tag{2.1.15}$$

Основная идея [99] заключается в подмене задачи (2.1.1) следующей задачей

$$f_\tau(x) = E_{\tilde{s},\xi}\left[f(x+\tau\tilde{s},\xi)\right] \to \min_{x\in Q}, \tag{2.1.16}$$

$\varepsilon/2$-решение которой при $\tau = \varepsilon/(2M)$, будет $\varepsilon$-решением исходной задачи (2.1.1).

Считая $\delta = O(\varepsilon)$ и[46] (приведенное условие выполняется, если мы имеем доступ к $\delta$-оракулу из предположения 2.1.2)

---

[44] Взятие математического ожидания по $\delta$ подчеркивает, что $\delta(x,\xi)$ может быть случайной величиной не только потому, что может зависеть от $\xi$, но и потому, что может иметь собственную случайность.

[45] Все свойства функции $f(x)$ при переходе к $f_\tau(x)$ могут только улучшиться. В частности, $f_\tau(x)$ также выпуклая функция (можно перенести и на сильную выпуклость с не меньшей константой), с константой Липшица и константой Липшица градиента (если таковая существует у $f(x)$) не большей чем у $f(x)$.

[46] Если это условие не выполняется, то все что написано далее останется верным, правда, при более ограничительных условиях на допустимый уровень шума (это касается и всего последующего изложения).

$$\left\|\nabla_x E_{\tilde{s},\xi,\delta}\left[\delta\left(x+\tau\tilde{s},\xi\right)\right]\right\|_2 = O\left(\varepsilon/R\right),$$

можно получить для среднего числа итераций (используя те же алгоритмы для задачи (2.1.16), что и в подразделе 2.1.2, со стохастическим градиентом (2.1.12)), соответствующие аналоги оценок (2.1.2), (2.1.4):

$$O\left(n^2M^2R^2/\varepsilon^2\right),\ O\left(n^2M^2/\left(\mu\varepsilon\right)\right).^{47} \qquad (2.1.17)$$

Если дополнительно известно, что $f\left(x,\xi\right)$ – равномерно гладкая по $x$ функция (это условие можно ослабить [52]) и п.н. по $\xi$

$$\left\|\nabla f\left(x,\xi\right)-\nabla f\left(y,\xi\right)\right\|_2 \le L\left\|x-y\right\|_2, \qquad (2.1.18)$$

то вместо (2.1.14) будем иметь

$$0 \le f_\tau\left(x\right)-f\left(x\right) \le \frac{L\tau^2}{2}. \qquad (2.1.19)$$

Из формулы (2.1.19) следует, что можно ослабить требование к неточности: допускать неточность оракула масштаба[48] $\delta \sim M\tau \sim M\sqrt{\varepsilon/L}$.

При сделанных дополнительных предположениях о гладкости (2.1.18) за счет ужесточения требований к масштабу допускаемой неточности $\delta$ (как именно требуется это сделать можно усмотреть из формулы (2.1.21); ниже мы вернемся к этому вопросу) можно улучшить скорость сходимости (фактор $n^2$ перейдет в $n$):

$$O\left(nM^2R^2/\varepsilon^2\right),\ O\left(nM^2/\left(\mu\varepsilon\right)\right). \qquad (2.1.20)$$

Оценки (2.1.20) в общем случае не улучшаемы (даже при $\delta = 0$) для гладких стохастических и негладких задач [190]. Фактически это означает, что мы можем выбрать настолько малое $\tau$ (насколько малым мы можем его выбрать определяется $\varepsilon$ и $\delta$), что конечная разность в (2.1.12) “превращается” (с нужной точностью) в производную по направлению.

---

Так, если не налагать это ограничение, то потребуется считать $\delta = O\left(\varepsilon^2/\left(\sqrt{n}MR\right)\right)$ или $\delta = O\left(\varepsilon^{3/2}/\left(\sqrt{nL}R\right)\right)$ – в случае, если $f\left(x\right)$ имеет $L$-липшицев градиент. Это можно получить с помощью замечания 2.1.8.

[47] Аналогично (2.1.2), (2.1.4) можно переписать оценку (2.1.17) не в среднем, как сейчас, а с учетом вероятностей больших уклонений. Это замечание касается и последующих вариаций формулы (2.1.17). Нам не известно являются ли оценки (2.1.17) оптимальными при заданном уровне шума $\delta = O\left(\varepsilon\right)$.

[48] Здесь мы дополнительно считаем, что $\nabla_x E_{\tilde{s},\xi,\delta}\left[\delta\left(x+\tau\tilde{s},\xi\right)\right]=0$. В частности, это условие выполняется, если неточность $\delta\left(x,\xi\right)$ имеет независимое от $x$ распределение.

Для объяснения отмеченного перехода полезно заметить, что [50, 52, 190, 270] (см. также (2.1.22))

$$E_{s,\xi}\left[\left\|g_{\tau,\delta}\left(x,s,\xi\right)\right\|_2^2\right] \le 4nM^2 + L^2\tau^2n^2 + \frac{8\delta^2n^2}{\tau^2}. \qquad (2.1.21)$$

**Замечание 2.1.7 (техника двойного сглаживания негладких задач Б.Т. Поляка [99], см. также [190, 308]).** За счет подмены изначально негладкого функционала в задаче (2.1.1) на

$$f(x) := f_\gamma(x) = E_{\tilde{s}_1,\xi}\left[f\left(x+\gamma\tilde{s}_1,\xi\right)\right], \gamma \le \varepsilon/(2M),$$

где $\tilde{s}_1$ – случайный вектор (независимый от $\xi$), равномерно распределенный на $B_2^n(1)$, получим новую задачу ($\varepsilon/2$-аппроксимирующую исходную), для которой при достаточно малом $\tau$ будет иметь место оценка (2.1.21) с достаточно большим $L \ge 2nM^2/\varepsilon$ (см. раздел 5.2 главы 5 и раздел 6.2 главы 6). Далее решая с помощью уже описанной техники с точностью $\varepsilon := \varepsilon/2$ задачу стохастической оптимизации (2.1.1) с

$$\xi := (\tilde{s}_1,\xi),\ f(x,\xi) := f\left(x+\gamma\tilde{s}_1,\xi\right),$$

$$g_{\tau,\delta}(x;s_2,\xi) := \frac{n}{\tau}\Big(f\left(x+\gamma\tilde{s}_1+\tau s_2,\xi\right)+\delta\left(x+\gamma\tilde{s}_1+\tau s_2,\xi\right)-\Big(f\left(x+\gamma\tilde{s}_1,\xi\right)+\delta\left(x+\gamma\tilde{s}_1,\xi\right)\Big)\Big)s_2,$$

где $s_2$ – случайный вектор (независимый от $\xi$ и $\tilde{s}_1$), равномерно распределенный на $S_2^n(1)$, получим те же оценки (2.1.20), только при существенно более жестких условиях на уровень шума $\delta$. К сожалению, получить конструктивное описание этих условий на данный момент не удалось.

Оценки (2.1.17) и (2.1.20) переносятся и на задачи стохастической онлайн оптимизации (см., например, [45, 50, 52, 190, 216]) с возникновением дополнительного фактора $\ln N$ в сильно выпуклом случае (см. формулу (2.1.11)). При этом даже в гладком случае не обязательно требовать дополнительно стохастичность исходной постановки для оптимальности оценок (2.1.20).

Далее рассматривается не стохастический вариант постановки задачи (2.1.1) (не обобщаемый на онлайн постановки) с $Q = \mathbb{R}^n$ (обобщения на произвольные выпуклые множества $Q \subset \mathbb{R}^n$ представляются интересными, но на данный момент нам не известны такие обобщения[49] – в последующих рассуждениях существенным образом используется то, что в точке минимума $\nabla f(x_*)=0$). Так что теперь $R$ – расстояние от точки старта до решения в 2-норме. В этом варианте, выписанная оценка (2.1.21) может быть уточнена

$$E_s\left[\left\|g_{\tau,\delta}(x,s)\right\|_2^2\right]\le 4n\left\|\nabla f(x)\right\|_2^2+L^2\tau^2n^2+\frac{8\delta^2n^2}{\tau^2}. \tag{2.1.22}$$

Последняя оценка следует из явления концентрации равномерной меры на $S_2^n(1)$ с выделенными полюсами вокруг экватора (см. [244] – в случае покомпонентных методов эта оценка особенно просто получается [45], $s$ в приводимой формуле, и только в ней, соответствует покомпонентной рандомизации):

$$E_s\left[\left\langle\nabla f(x),s\right\rangle^2\right]=\frac{1}{n}\left\|\nabla f(x)\right\|_2^2.$$

Считая для простоты формулировок, что

$$\nabla_x E_{\tilde{s},\delta}\left[\delta\left(x+\tau\tilde{s}\right)\right]=0,$$

можно распространить метод [193], дающий оценки (2.1.8), на текущий контекст, и получить следующие оценки (с $\tau\sim\sqrt{\delta/L}$) числа итераций для достижения точности $\varepsilon$ для случая выпуклой и сильно выпуклой целевой функции соответственно

$$N_1(\varepsilon)=n\cdot\mathrm{O}\left(\left(\frac{LR^2}{\varepsilon}\right)^{\frac{1}{p+1}}\right),N_2(\varepsilon)=n\cdot\mathrm{O}\left(\left(\frac{L}{\mu}\right)^{\frac{1}{p+1}}\ln\left(\frac{LR^2}{\varepsilon}\right)\right) \tag{2.1.23}$$

при (условия на допустимый уровень шума, при котором оценки (2.1.23) имеют такой же вид, с точностью до $\mathrm{O}(1)$, как если бы шума не было[50])

$$\delta_1(\varepsilon)\le\frac{1}{n}\mathrm{O}\left(\varepsilon\cdot\left(\frac{\varepsilon}{LR^2}\right)^{\frac{p}{p+1}}\right),\ \delta_2(\varepsilon)\le\frac{1}{n}\mathrm{O}\left(\varepsilon\cdot\left(\frac{\mu}{L}\right)^{\frac{p}{p+1}}\right). \tag{2.1.24}$$

---

[49] По-видимому, такие обобщения возможны. Также, по-видимому, возможно перенесение концепции универсальных методов (см. подраздел 2.1.3) на безградиентные методы и спуски по направлению (покомпонентные спуски).

[50] В отсутствии шума, оракул нам фактически может выдавать производную по направлению $s$

$$g(x,s)=n\left\langle\nabla f(x),s\right\rangle s,$$

точнее $\left\langle\nabla f(x),s\right\rangle$, $s$ мы генерируем сами. Если, в свою очередь, считать, что $\left\langle\nabla f(x),s\right\rangle$ оракул выдает с аддитивным шумом (для простоты считаем, независящим от $s$) масштаба $\delta:=\sqrt{L\delta}$ ($\delta$ в правой части определяется исходя из формулы (2.1.24)), то формула (2.1.23) останется верной [41].

По-видимому (строгим доказательством мы не располагаем на данный момент), и в стохастическом случае имеет место аналог формул (2.1.23), (2.1.24) с заменой в формуле (2.1.23)

$$N_1(\varepsilon) = n \cdot \max\left\{ \mathrm{O}\left(\left(\frac{LR^2}{\varepsilon}\right)^{\frac{1}{p+1}}\right), \mathrm{O}\left(\frac{DR^2}{\varepsilon^2}\right)\right\},\ N_2(\varepsilon) = n \cdot \max\left\{ \mathrm{O}\left(\left(\frac{L}{\mu}\right)^{\frac{1}{p+1}} \ln\left(\frac{LR^2}{\varepsilon}\right)\right), \mathrm{O}\left(\frac{D}{\mu\varepsilon}\right)\right\}.$$

Можно продолжать переносить все написанное в подразделе 2.1.3 на рассматриваемую ситуацию (частично это уже сделано в [41, 45]). Однако мы остановимся лишь на наиболее интересном (на наш взгляд) месте. А именно, на согласовании прокс-структуры с рандомизацией, порождающей сглаживание.

Основным результатом (ввиду замечания 2.1.7) в негладком и(или) стохастическом случае с точным оракулом здесь является следующее наблюдение [50, 52]: независимо от выбора прокс-структуры рандомизацию всегда стоит выбирать согласно (2.1.12) (если ставить цель – минимизировать число итераций), т.е. с помощью разыгрывания случайного вектора $s$ равномерно распределенного на $S_2^n(1)$. В случае неточного оракула, сформулированное утверждение требует оговорок [52] (см. также раздел 5.2 главы 5). Ограничимся далее обобщением оценок (2.1.20) на случай использования общих прокс-структур.

Приведем соответствующее обобщение формулы (2.1.21) (здесь и далее нижний индекс “2” у констант Липшица подчеркивает, что они считаются согласно евклидовой норме из-за сделанного нами выбора способа рандомизации)

$$E_{s,\xi}\left[\left\|g_{\tau,\delta}(x,s,\xi)\right\|_{q'}^2\right] = \mathrm{O}\left(\left(4M_2^2 n + L_2^2\tau^2 n^2 + \frac{8\delta^2 n^2}{\tau^2}\right) E_s\left[\|s\|_{q'}^2\right]\right),$$

где в прямом пространстве выбрана $q$-норма и $1/q + 1/q' = 1$. Согласно замечанию 2.1.2 можно считать, что $2 \le q' \le \infty$ – выбирать другие нормы, как правило, бывает не выгодно. Для такого диапазона $E_s\left[\|s\|_{q'}^2\right] = \tilde{\mathrm{O}}\left(n^{2/q'-1}\right)$, в частности $E_s\left[\|s\|_{\Omega(\log n)}^2\right] = \tilde{\mathrm{O}}\left(n^{-1}\right)$ ($\tilde{\mathrm{O}}(\ )$ с точностью до логарифмического фактора от $n$ совпадает с $\mathrm{O}(\ )$, аналогично с $\tilde{\Omega}(\ )$). Исходя из такого обобщения, можно привести следующую Таблицу 2.1.1, распространяющую оценку (2.1.20) на произвольные прокс-структуры ($R^2$ – “расстояние” Брэгмана, согласованное с $q$-нормой, см. замечание 2.1.1, а также следующий раздел 2.2).

Выпишем условия, из которых можно получить требования на шум (нам представляется, что здесь это может хорошо прояснить суть дела):

- $\min\left\{M_2\tau, L_2\tau^2/2\right\} = \mathrm{O}(\varepsilon)$ – условие достаточной точности аппроксимации исходной функции ее сглаженной версией;
- $L_2^2\tau^2n^2 + \dfrac{8\delta^2n^2}{\tau^2} = \mathrm{O}\left(M_2^2n\right)$ – условие “правильной” ограниченности квадрата нормы аппроксимации стохастического градиента.

Таблица 2.1.1

| $f(x)$ – выпуклая | $f(x)$ – $\mu_q$-сильно выпуклая в $q$-норме |
|---|---|
| $\mathrm{O}\left(\dfrac{nM_2^2R^2}{\varepsilon^2}\right)\tilde{\mathrm{O}}\left(n^{2/q'-1}\right)$ | $\tilde{\mathrm{O}}\left(\dfrac{nM_2^2}{\mu_q\varepsilon}\right)\tilde{\mathrm{O}}\left(n^{2/q'-1}\right)$ |

Выписанные условия позволяют для всех полей Таблицы 2.1.1 (с оценками) написать соответствующие условия на допустимый уровень шума, и параллельно подобрать оптимальный размер параметра сглаживания $\tau$.

Аналогичную Таблицу 2.1.2 можно записать (на базе конструкций работ [45, 52]) и для оценки (2.1.23) ($p=1$)

Таблица 2.1.2

| $f(x)$ – гладкая выпуклая | $f(x)$ – гладкая $\mu_q$-сильно выпуклая в $q$-норме |
|---|---|
| $\tilde{\mathrm{O}}\left(n^{1/q'+1/2}\sqrt{\dfrac{L_2R^2}{\varepsilon}}\right)$ | $\tilde{\mathrm{O}}\left(n^{1/q'+1/2}\sqrt{\dfrac{L_2}{\mu_q}}\right)$ |

**Замечание 2.1.8 (см. [50]).** Если не делать никаких предположений о шуме $\delta(x,\xi)$ в предположении 2.1.2, кроме $|\delta(x,\xi)| \le \delta$, то для получения требований на уровень шума $\delta$, потребуется еще воспользоваться следующим утверждением.

Пусть последовательность независимых случайных векторов $\{s_k\}_{k=0}^N$, равномерно распределенных на $S_2^n(1)$, и $\left\{x_k\left(\{s_l\}_{l=0}^{k-1}\right)\right\}_{k=1}^N$ обладают свойством $E\left[\|x_k - x_*\|_2^2\right] \le R^2$. Тогда

$$E\left[\frac{1}{N}\sum_{k=1}^{N}\left|\langle s_k, x_k - x_*\rangle\right|\right] \le \frac{2R}{\sqrt{n}}.$$

В виду замечаний 2.1.2, 2.1.4 при использовании этого утверждения можно считать, что $R^2 = \mathrm{O}\left(V\left(x_*, x_0\right)\right)$. Причем константа в $\mathrm{O}(\ )$ может быть сделана $\sim 1$.

Приведенная в первом столбце Таблицы 2.1.1 оценка при определенных условиях может быть лучше нижней оценки [190]. [51] Здесь ситуация аналогична той, о которой написано в конце замечания 2.1.2 (будем использовать те же обозначения) и в подразделе 2.1.3. Например, если $Q = B_1^n(1)$, то в нижней оценке [190] стоит $\tilde{\Omega}\left(nM_1^2/\varepsilon^2\right)$ ($E_\xi\left[\left\|\nabla f\left(x,\xi\right)\right\|_\infty^2\right] \le M_1^2$), а в таблице будет стоять $\tilde{\mathrm{O}}\left(M_2^2/\varepsilon^2\right)$. Осталось заметить, что $M_2^2 \le nM_1^2$, причем в определенных ситуациях может быть $M_2^2 \ll nM_1^2$. Подробнее об этом написано в разделе 5.2 главы 5.

Основные конкурирующие рандомизации в гладком случае – это рандомизация на евклидовой сфере и покомпонентная рандомизация [45, 50, 52, 190], которая используется в основном только с евклидовой прокс-структурой [45] (см. также раздел 5.1 главы 5). Исследования последних нескольких лет показали (см., например, [45, 138, 264], а также ра-

[51] Аналогичное можно сказать и для Таблицы 2.1.2. А именно, рассмотрим задачу минимизации гладкого выпуклого функционала $f(x)$ с константой Липшица градиента в 2-норме равной $L_2$ на множестве $Q = B_1^n(R)$. Тогда нижняя оценка (при $N \le n$) будет иметь вид (А.С. Немировский, 2015)

$$f\left(x^N\right) - f_* \ge \frac{\tilde{C}_1 L_2 R^2}{N^3}.$$

С другой стороны, если использовать обычный FGM с KL-прокс-структурой для этой же задачи, то верхняя оценка будут иметь вид

$$f\left(x^N\right) - f_* \le \frac{\tilde{C}_2 L_1 R^2}{N^2},$$

где константа Липшица градиента в 1-норме $L_1$: $L_2/n \le L_1 \le L_2$, а $\tilde{C}_1$, $\tilde{C}_2 > 0$ – некоторые небольшие универсальные константы. Тем не менее, отсюда нельзя сделать вывод, что нижняя оценка достигается. Достигается ли эта нижняя оценка, и, если достигается, на каком методе? – насколько нам известно, это пока открытый вопрос, поставленный А.С. Немировским в 2015 году. Однако если оценивать не число итераций, а общее число арифметических операций и если ограничиться рассмотрением класса функций, для которых стоимость расчета производной по направлению (или значения функции в точке) в $\sim n$ раз меньше стоимости расчета полного градиента [67] (в виду БАД [65, 76] это предположение довольно обременительное, впрочем, если функция задана моделью черного ящика, выдающего только значение функции, а градиент восстанавливается при $n+1$ таком обращении, то сделанное предположение кажется вполне естественным), то при $N \le n$ выписанная выше нижняя оценка (в варианте для общего числа арифметических операций, необходимых для достижения заданной точности) будет соответствовать оценке из первого столбца Таблицы 2.

боты P. Richtarik-а и T. Zhang-а (http://www.stat.rutgers.edu/home/tzhang/ ) и раздел 5.1 главы 5), что для довольно большого класса гладких задач выпуклой оптимизации в пространствах огромной размерности,[52] возникающих в самых разных приложениях (см. [40, 138, 264] и работы P. Richtarik-а), покомпонентные методы являются наиболее эффективным способом решения (с точки зрения общего числа арифметических операций для достижения заданной точности по функции). Покомпонентные методы, безусловно, заслуживают отдельного подробного обзора. Поэтому здесь мы ограничимся только ссылкой на такие обзоры [45, 327] и раздел 5.1 главы 5.

Приведем далее несколько примеров, демонстрирующих важность изучения безградиентных методов (часто эти методы называют прямыми методами [99] или методами нулевого порядка [270]).

**Пример 2.1.5 (двухуровневая оптимизация [180]).** Требуется решить задачу, возникающую, например, при поиске равновесия по Штакельбергу [21]

$$\psi(x,u) \to \max_{u \in U},$$

$$f(x,u(x)) \to \min_{x}.$$

Из первой задачи находится зависимость $u(x)$, которая входит во вторую (внешнюю) задачу. Проблема здесь в том, что явная зависимость $u(x)$ в общем случае может быть недоступна. Как следствие, могут быть проблемы с расчетом $\nabla u(x)$. Поэтому предлагается приближенно решать первую задачу и использовать безградиентный метод с неточным оракулом для второй. Насколько точно надо решать первую задачу, и какой именно безградиентный метод (с точки зрения чувствительности к неточности) выбирать для второй – определяется сложностью решения первой задачи и свойствами второй.

Рассмотренная двухуровневая задача может быть сильно упрощена, если удается найти ее седловое представление [9, 29, 35, 42, 43, 47, 92]. В частности, если функции $\psi(x,u)$, $f(x,u)$ – выпуклы по $x$ и вогнуты по $u$, и $\psi(x,u)$ – простой структуры, то можно (для достаточно большого $\lambda$) заменить исходную задачу на следующую (предложенную Ю.Е. Нестеровым)

$$\min_{x} \max_{u \in U} \left[ f(x,u) + \lambda \psi(x,u) \right].$$

[52] Во всех этих задачах можно считать полные градиенты, и строить на их базе различные методы. То есть для таких задач выбор покомпонентного метода – осмысленный выбор наиболее быстрого способа решения, а не следствие каких-то (в том числе вычислительных) ограничений на задачи.

Полученную седловую задачу стоит решать методами композитной оптимизации (см. замечание 2.1.6 и [92, 259]), чтобы параметр $\lambda$ либо совсем не входил в оценки числа итераций, либо входило очень слабо.

К сожалению, седловое представление возможно далеко не всегда.

**Пример 2.1.6 (огромная скрытая размерность).** Пусть $y(x) \in \mathbb{R}^m$, $x \in \mathbb{R}^n$, $n \ll m$. Требуется решить задачу

$$f\left(x, y(x)\right) \to \min_x .$$

Мы предполагаем, что можем эффективно посчитать с необходимой точностью $y(x)$ и $f\left(x, y(x)\right)$ за $O(m)$. Если для решения этой задачи оптимизации мы будем использовать безградиентный метод, то общее число сделанных арифметических операций пропорционально $mn$ (см. (2.1.20), (2.1.23)). Заметим, что если бы мы могли использовать обычный градиентный метод, то общее число сделанных арифметических операций также было бы пропорционально $mn$, однако вычисление градиента по разным причинам может быть затруднено (см. пример 2.1.5 и метод MCMC [46] для расчета PageRank в подходе [154, 155]). В действительности, часто имеет место следующее полезное наблюдение [154]: если мы можем вычислить значения $y(x)$ и $f\left(x, y(x)\right)$ за $O(m)$, то мы можем с такой же по порядку сложностью (и затратами памяти) вычислить и производные по фиксированному направлению $h$:

$$\frac{dy(x)}{dx}h,\ \left\langle \frac{\partial f(x,y)}{\partial x}, h \right\rangle + \left\langle \frac{\partial f(x,y)}{\partial y}, \frac{dy(x)}{dx}h \right\rangle .$$

Тем не менее, тут требуется много оговорок, в том числе про точность расчетов. Если не вдаваться в детали, то такие рассуждения также приводят к затратам пропорциональным $mn$, где $n$ возникло в виду оценок (2.1.20), (2.1.23) для спусков по направлению. Только в отличие от полно-градиентного метода для покомпонентного метода константа Липшица градиента функционала в оценке числа итераций уже будет рассчитываться не по худшему направлению, а в среднем (это может давать выгоду, по порядку равную корню квадратному из размерности пространства – см. [45, 282], да и, как правило, будет ощутимая выгода в затрачиваемых ресурсах машинной памяти [154]. Оговорки о точности здесь все же необходимы, поскольку для безградиентных методов и спусков по направлению требования к точности могут существенно отличаться (об этом ранее уже было немного написано в данном пункте). Как следствие, в оценку $O(m)$ необходимо явно вводить зависи-мость от точности вычисления $y(x)$ и $f(x, y(x))$. Основные технические детали тут были проработаны П.Е. Двуреченским [154, 155].

В примерах 2.1.5, 2.1.6, в действительности, требуются некоторые оговорки о невозможности или неэффективности использования БАД для полно-градиентного метода (см. подраздел 2.1.2, а также [65, 76]). Нам известны случаи, подпадающие под разобранные примеры, в которых не понятно, как можно было бы воспользоваться БАД [21, 154, 155]. В частности, в работах [154, 155], соответствующей примеру 2.1.6, сложность в том, что БАД хочется использовать для ускорения вычисления матрицы Якоби отображения (вектор функции) $y(x) \in S_m(1)$, неявно заданного уравнением $y = P(x)y$, со стохастической (по столбцам) эргодической матрицей $P(x)$ со спектральной щелью $\alpha$ и числом ненулевых элементов $sm$. Метод простой итерации позволяет с точностью $\varepsilon$ найти $y(x)$ за время $O\left(sm\left(n + \alpha^{-1}\ln\left(\varepsilon^{-1}\right)\right)\right)$ с затратами памяти $O(sm)$. Нам не известно более эффективного способа расчета матрицы Якоби отображения $y(x)$, чем естественное обобщение метода простой итерации (для продифференцированного по $x$ уравнения $y(x) = P(x)y(x)$), требующее затрат времени $O\left(smn\alpha^{-1}\ln\left(\varepsilon^{-1}\right)\right)$ и памяти $O(smn)$. Для реальных приложений [154, 155]: $m \sim 10^9$, $s \sim 10^2$, $n \sim 10^3$, $\alpha \sim 10^{-1}$, $\varepsilon \sim 10^{-12}$ Отсюда ясно, что при использовании полно-градиентного метода просто невозможно будет выделить даже у 64-битной операционной системы, стоящей на самом современном персональном компьютере, необходимой памяти под работающую программу, в основе которой лежит полно-градиентный подход.

Кроме того, в примерах 2.1.5, 2.1.6 важно уметь эффективно пересчитывать значения $u(x)$ или $y(x)$, а не рассчитывать их каждый раз заново (на каждой итерации внешнего цикла). Поясним сказанное. Предположим, что мы уже как-то посчитали, скажем, $u(x)$, решив с какой-то точностью соответствующую задачу оптимизации. Тогда для вычисления $u(x + \Delta x)$ (на следующей итерации внешнего цикла) у нас будет хорошее начальное приближение $u(x)$. А как известно (см. подразделы 2.1.2–2.1.4) расстояние от точки старта до решения (не в сильно выпуклом случае) существенным образом определяет время работы алгоритма оптимизции. Эта конструкция (hot/warm start) напоминает фрагмент обоснования сходимости методов внутренней точки при изучении движения по центральному пути [93, 163, 259, 276]. Тем не менее, известные нам приложения (пример 2.1.4 подраздела 2.1.3 и [42, 43, 154, 155]) пока как раз всецело соответствуют сильно выпуклой ситуации. Связано это с тем, что если расчет $u(x)$ или $y(x)$ с точностью $\varepsilon$ осуществляется за $O\left(C\ln\left(R/\varepsilon\right)\right)$ операций, то для внешней задачи можно выбирать самый быстрый метод (а, стало быть, и самый требовательный к точности), и с точностью до того, что стоит под логарифмом, общая трудоемкость будет просто прямым произведением трудоемкостей решения внутренней и внешних задач по отдельности. Как правило, такое сочетание оказывается недоминируемым.

Также необходимо отметить, что, как правило, итоговые задачи оптимизации (после подстановки зависимости $u(x)$ или $y(x)$ в задачу верхнего уровня) в этих примерах получаются не выпуклыми. В этой связи можно лишь говорить о локальной сходимости к стационарной точке. В отсутствие выпуклости даже если ограничиться локальной сходимостью многое из того, что описано в данном разделе, требует отдельного рассмотрения – см. [270, 282].

Отметим в заключение, что если немного по-другому посмотреть на описанное в этом пункте, то можно заметить следующее. Какой бы большой (но равномерно ограниченный по итерациям) шум ни был, если $\delta(x+\tau s,\xi)$ имеет распределение не зависящее от $s$, то (возможно через очень большое число итераций) мы сможем сколь угодно точно (по функции) решить задачу! Аналогичное можно сказать, если мы изначально исходим из концепции оракула, выдающего зашумленное значение $\left\langle \nabla f(x), s\right\rangle$, причем зашумленность не зависит от $s$. Все это восходит к идеям Р. Фишера, развитым О.Н. Граничиным и Б.Т. Поляком [59].

## 2.2 Универсальный метод для задач стохастической композитной оптимизации

### 2.2.1 Введение

В цикле работ [39, 44, 48, 181, 183, 184, 185, 266, 271, 274] была предложена линейка быстрых градиентных методов для задач композитной оптимизации (выпуклый и сильно выпуклый случаи). Рассматривались приложения этих методов, в том числе к задачам стохастической оптимизации [39, 48, 181]. Однако открытым оставался вопрос (см. предыдущий раздел): возможно ли на базе этих методов предложить универсальный метод (настраивающийся на гладкость задачи) для задач стохастической оптимизации и(или) в сильно выпуклом случае? В данном разделе получены ответы на эти вопросы. В основе предлагаемого подхода лежит оригинальный вариант классического быстрого градиентного метода (см. [271]): метод треугольника (МТ), изложенный в подразделе 2.2.2 (материал, приведенный в подразделе 2.2.2, был предложен в апреле 2016 года Ю.Е. Нестеровым). Особенностью МТ является необходимость выполнения всего одного "проектирования" на каждой итерации. Как следствие, предложенный метод оказался заметно проще многих своих аналогов. Это не только позволило упростить и привести в этом разделе известные результаты для быстрых градиентных методов, но и продвинуться в решении отмеченных выше вопросов. В подразделе 2.2.3 МТ распространяется на задачи сильной выпуклой композитной оптимизации. В подразделе 2.2.4 предлагается универсальный вариант МТ для задач выпуклой и сильно выпуклой композитной оптимизации. В заключительном подразделе 2.2.5 универсальный МТ из подраздела 2.2.4 переносится на задачи стохастической оптимизации.

### 2.2.2 Метод треугольника для задач композитной оптимизации

Рассматривается задача выпуклой композитной оптимизации [266]

$$F(x) = f(x) + h(x) \to \min_{x \in Q}. \qquad (2.2.1)$$

Положим $R^2 = V(x_*, y^0)$, где прокс-расстояние (расстояние Брэгмана) определяется формулой (см., например, [271], главу 2 [181], [183, 259])

$$V(x,z) = d(x) - d(z) - \langle \nabla d(z), x - z \rangle;$$

прокс-функция $d(x) \geq 0$ ($d(y^0) = 0$, $\nabla d(y^0) = 0$) считается сильно выпуклой относительно выбранной нормы $\| \ \|$, с константой сильной выпуклости $\geq 1$; $x_*$ – решение задачи (1) (если решение не единственно, то выбирается то, которое доставляет минимум $V(x_*, y^0)$).

**Предположение 2.2.1.** *Для любых* $x, y \in Q$ *имеет место неравенство*

$$\left\|\nabla f(y) - \nabla f(x)\right\|_* \le L\left\|y - x\right\|.$$

Опишем вариант быстрого градиентного метода для задачи (2.2.1) с одной “проекцией”, который мы далее будем называть “метод треугольника” (МТ).

Положим

$$\varphi_0(x) = V(x, y^0) + \alpha_0\left[f(y^0) + \left\langle \nabla f(y^0), x - y^0\right\rangle + h(x)\right],$$

$$\varphi_{k+1}(x) = \varphi_k(x) + \alpha_{k+1}\left[f(y^{k+1}) + \left\langle \nabla f(y^{k+1}), x - y^{k+1}\right\rangle + h(x)\right], \tag{2.2.2}$$

$$A_k = \sum_{i=0}^{k}\alpha_i, \quad \alpha_0 = L^{-1}, \quad A_k = \alpha_k^2 L, \quad k = 0,1,2,\ldots, \quad x^0 = u^0 = \arg\min_{x\in Q}\varphi_0(x). \tag{2.2.3}$$

**<u>Метод Треугольника</u>**

$$\boxed{\begin{aligned} y^{k+1} &= \frac{\alpha_{k+1}u^k + A_k x^k}{A_{k+1}}, \\ u^{k+1} &= \arg\min_{x\in Q}\varphi_{k+1}(x), \\ x^{k+1} &= \frac{\alpha_{k+1}u^{k+1} + A_k x^k}{A_{k+1}}. \end{aligned}} \tag{2.2.4}$$

**Лемма 2.2.1 (см. [271]).** *Последовательность* $\{\alpha_k\}$*, определяемую формулой (2.2.3), можно задавать рекуррентно*

$$\alpha_{k+1} = \frac{1}{2L} + \sqrt{\frac{1}{4L^2} + \alpha_k^2}\,.$$

*При этом*

$$A_k \ge \frac{(k+1)^2}{4L}.$$

**Лемма 2.2.2.** *Пусть справедливо предположение 2.2.1. Тогда для любого* $k = 0,1,2,\ldots$ *имеет место неравенство*

$$A_k F(x^k) \le \varphi_k^* = \min_{x\in Q}\varphi_k(x) = \varphi_k(u^k). \tag{2.2.5}$$

**Доказательство.** Проведем по индукции. При $k = 0$ формула (2.2.5) следует из того, что для любого $x \in Q$

$$F(x) = f(x) + h(x) \le f(y^0) + \left\langle \nabla f(y^0), x - y^0\right\rangle + \frac{L}{2}\left\|x - y^0\right\|^2 + h(x) \le$$

$$\le LV(x, y^0) + f(y^0) + \left\langle \nabla f(y^0), x - y^0\right\rangle + h(x) = \frac{1}{A_0}\varphi_0(x).$$

Последнее неравенство следует из того, что для любых $x, z \in Q$

$$V(x,z) \geq \frac{1}{2}\|x-z\|^2.$$

Это неравенство, в свою очередь, следует из $\geq 1$-сильной выпуклости $d(x)$ в норме $\| \ \|$.

Итак, пусть формула (2.2.5) установлена для $k$, покажем, что тогда она будет справедлива и для $k+1$. По определению (2.2.2)

$$\varphi_{k+1}^* = \min_{x \in Q} \varphi_{k+1}(x) = \varphi_{k+1}(u^{k+1}) =$$

$$= \varphi_k(u^{k+1}) + \alpha_{k+1}\left[ f(y^{k+1}) + \langle \nabla f(y^{k+1}), u^{k+1} - y^{k+1} \rangle + h(u^{k+1}) \right]. \quad (2.2.6)$$

Поскольку по предположению индукции $A_k F(x^k) \leq \varphi_k(u^k)$, и $\varphi_{k+1}(x)$ – сильно выпуклая в $\| \ \|$-норме функция с константой $\geq 1$ (это следует из аналогичного свойства функции $V(x, y^0)$, что, в свою очередь, следует из аналогичного свойства функции $d(x)$), то

$$\varphi_k(u^{k+1}) \geq \varphi_k(u^k) + \frac{1}{2}\|u^{k+1} - u^k\|^2 \geq A_k \cdot \left( f(x^k) + h(x^k) \right) + \frac{1}{2}\|u^{k+1} - u^k\|^2.$$

Из выпуклости $f(x)$ отсюда имеем

$$\varphi_k(u^{k+1}) \geq A_k f(y^{k+1}) + \langle \nabla f(y^{k+1}), A_k \cdot (x^k - y^{k+1}) \rangle + A_k h(x^k) + \frac{1}{2}\|u^{k+1} - u^k\|^2. \quad (2.2.7)$$

Подставляя (2.2.7) в (2.2.6), получим

$$\varphi_{k+1}^* \geq A_{k+1} \cdot \underbrace{\left( \frac{A_k}{A_{k+1}} h(x^k) + \frac{\alpha_{k+1}}{A_{k+1}} h(u^{k+1}) \right)}_{\geq h(x^{k+1})} + A_{k+1} f(y^k) +$$

$$+ \langle \nabla f(y^{k+1}), \underbrace{\alpha_{k+1} \cdot (u^{k+1} - y^{k+1}) + A_k \cdot (x^k - y^{k+1})}_{= A_{k+1} \cdot (x^{k+1} - y^{k+1})} \rangle + \underbrace{\frac{1}{2}\|u^{k+1} - u^k\|^2}_{= \frac{A_{k+1}^2}{2\alpha_{k+1}^2}\|x^{k+1} - y^{k+1}\|^2}. \quad (2.2.8)$$

Исходя из выпуклости функции $h(x)$ и описания МТ (4), формулу (2.2.8) можно переписать следующим образом

$$\varphi_{k+1}^* \geq A_{k+1}\left[ f(y^{k+1}) + \langle \nabla f(y^{k+1}), x^{k+1} - y^{k+1} \rangle + \frac{A_{k+1}}{2\alpha_{k+1}^2}\|x^{k+1} - y^{k+1}\|^2 + h(x^{k+1}) \right]. \quad (2.2.9)$$

Из предположения 2.2.1 следует, что если $A_{k+1}/\alpha_{k+1}^2 \geq L$, то

$$f(y^{k+1}) + \langle \nabla f(y^{k+1}), x^{k+1} - y^{k+1} \rangle + \frac{A_{k+1}}{2\alpha_{k+1}^2}\|x^{k+1} - y^{k+1}\|^2 \geq f(x^{k+1}). \quad (2.2.10)$$

Согласно (2.2.3) $A_{k+1}/\alpha_{k+1}^2 = L$, поэтому формула (2.2.10) имеет место. С помощью формулы (2.2.10) формулу (2.2.9) можно переписать следующим образом

$$\varphi_{k+1}^* \geq A_{k+1}\left[ f\left(x^{k+1}\right) + h\left(x^{k+1}\right)\right] = A_{k+1}F\left(x^{k+1}\right).$$

Таким образом, шаг индукции установлен. Следовательно, лемма 2.2.2 доказана. ■

Из лемм 2.2.1, 2.2.2 получаем следующий результат, означающий, что МТ сходится как обычный быстрый градиентный метод, см., например, [181, 183, 259, 266, 271], т.е. МТ сходится оптимальным образом для рассматриваемого класса задач [91].

**Теорема 2.2.1.** *Пусть справедливо предположение 2.2.1. Тогда МТ (2) – (4) для задачи (1) сходится согласно оценке*

$$F\left(x^N\right) - \min_{x\in Q} F\left(x\right) \leq \frac{4LR^2}{\left(N+1\right)^2}. \tag{2.2.11}$$

**Доказательство.** Из леммы 2.2.2 следует, что (в третьем неравенстве используется выпуклость функции $f\left(x\right)$)

$$A_N F\left(x^N\right) \leq \min_{x\in Q}\left\{ V\left(x,y^0\right) + \sum_{k=0}^{N}\alpha_k\left[ f\left(y^k\right) + \left\langle \nabla f\left(y^k\right), x - y^k\right\rangle + h\left(x\right)\right]\right\} \leq$$

$$\leq V\left(x_*,y^0\right) + \sum_{k=0}^{N}\alpha_k \underbrace{\left[ f\left(y^k\right) + \left\langle \nabla f\left(y^k\right), x_* - y^k\right\rangle + h\left(x_*\right)\right]}_{\leq f\left(x_*\right)+h\left(x_*\right)} \leq$$

$$\leq V\left(x_*,y^0\right) + \sum_{k=0}^{N}\alpha_k F\left(x_*\right) = R^2 + A_N F\left(x_*\right). \tag{2.2.12}$$

Заметим, что из второго неравенства следует, что если решение задачи (2.2.1) $x_*$ не единственно, то можно выбирать то, которое доставляет минимум $V\left(x_*,y^0\right)$. Именно таким образом возникает $R^2$ в оценке (2.2.12). Для доказательства теоремы осталось подставить нижнюю оценку на $A_N$ из леммы 2.2.1 в формулу (2.2.12). ■

**Замечание 2.2.1.** В действительности в формуле (2.2.12) содержится более сильный результат, чем в формуле (2.2.11). А именно, формула (2.2.12) еще означает, что МТ – прямо-двойственный метод (см., например, [6, 269] и главу 3). Мы не будем здесь подробно на этом останавливаться, отметим лишь, что это свойство позволяет получать эффективные критерии останова для МТ. Критерий останова позволяет не делать предписанного формулой (2.2.11) числа итераций и останавливаться раньше (по достижению желаемой точности). Это замечание можно распространить и на все последующее изложение.

МТ получил такое название, поскольку на базе трех точек $u^k$, $u^{k+1}$, $x^k$ можно построить треугольник, а точки $y^{k+1}$ и $x^{k+1}$ лежат на сторонах этого треугольника ( $y^{k+1}$ лежит на стороне $u^k x^k$, а $x^{k+1}$ на стороне $u^{k+1}x^k$ ), причем прямая, проходящая через точки $y^{k+1}$ и $x^{k+1}$, параллельна прямой, проходящей через точки $u^k$ и $u^{k+1}$.

Из описанной геометрической интерпретации получаются следующие результаты, распространимые и на все последующее изложение этого раздела 2.2 (с некоторыми оговорками в п. 5 [44])

**Следствие 2.2.1.** *Для любого* $k = 0,1,2,...$ *имеют место неравенства*

$$\left\|u^k - x_*\right\|^2 \le 2V\left(x_*, y^0\right),$$

$$\max\left\{\left\|x^k - x_*\right\|^2, \left\|y^k - x_*\right\|^2\right\} \le 4V\left(x_*, y^0\right) + 2\left\|x^0 - y^0\right\|^2 = \tilde{R}^2.$$

**Доказательство.** Второе неравенство следует из первого, описанной выше геометрической интерпретации и неравенства $\left\|a+b\right\|^2 \le 2\left\|a\right\|^2 + 2\left\|b\right\|^2$.

Докажем первое неравенство. Для этого воспользуемся леммой 2.2.1 и тем, что $\varphi_k(x)$ – сильно выпуклая в $\| \ \|$-норме функция с константой $\ge 1$. Для любого $x \in Q$

$$A_k F\left(x^k\right) + \frac{1}{2}\left\|x - u^k\right\|^2 \le \varphi_k^* + \frac{1}{2}\left\|x - u^k\right\|^2 \le \varphi_k(x) \le$$

$$\le \underbrace{\sum_{i=0}^{k} \alpha_i \left[ f\left(y^i\right) + \left\langle \nabla f\left(y^i\right), x - y^i \right\rangle + h(x) \right]}_{\le A_k F(x)} + V\left(x, y^0\right) \le A_k F(x) + V\left(x, y^0\right).$$

Выбирая $x = x_*$ и используя то, что $F\left(x^k\right) \ge F\left(x_*\right)$, получим первое неравенство следствия 1. ■

Ранее такого типа результаты (как следствие 2.2.1) устанавливались только для евклидовой прокс-структуры в некомпозитном случае [6] (см. также раздел 3.2 главы 3). Следствие 2.2.1 играет важную роль в случае неограниченных множеств $Q$, на которых нельзя равномерно ограничить константу $L$, поскольку гарантирует, что как бы “плохо” себя не вела функция $F(x)$ вне шара конечного радиуса $\tilde{R}$ (зависящего от качества начального приближения) с центром в решении $x_*$, это никак не скажется на скорости сходимости метода, поскольку итерационный процесс никогда не выйдет за пределы этого шара.

**Следствие 2.2.2.** *Пусть* $h(x) \equiv 0$ *(т.е.* $F(x) = f(x)$*) и* $\nabla f\left(x_*\right) = 0$. *Тогда*

$$\max\left\{F\left(x^N\right),F\left(y^N\right),F\left(z^N\right)\right\}-\min_{x\in Q}F\left(x\right)\le\frac{L\tilde{R}^2}{N^2}.$$

### 2.2.3 Метод треугольника для сильно выпуклых задач композитной оптимизации

В данном пункте будем считать, что $f\left(x\right)$ в задаче (2.2.1) обладает следующим свойством.

**Предположение 2.2.2.** *$f\left(x\right)$ – $\mu$-сильно выпуклая функция в норме $\|\;\|$, т.е. для любых $x,y\in Q$ имеет место неравенство*

$$f\left(y\right)+\left\langle\nabla f\left(y\right),x-y\right\rangle+\frac{\mu}{2}\|x-y\|^2\le f\left(x\right). \tag{2.2.13}$$

Введем (в евклидовом случае $\tilde{\omega}_n=1$)

$$\tilde{\omega}_n=\sup_{x,y\in Q}\frac{2V\left(x,y\right)}{\|y-x\|^2}.$$

**Замечание 2.2.2.** Из работы [228] следует, что в сильно выпуклом случае (когда сильно выпукла гладкая часть функционала $f\left(x\right)$, как в предположении 2.2.2) естественно выбирать именно евклидову норму и прокс-структуру, т.е. в большинстве случаев можно считать $\tilde{\omega}_n=1$. Поясним это примером из [228]. Число обусловленности (отношение константы Липшица градиента $L_1\left(f\right)$ к константе сильной выпуклости $\mu_1\left(f\right)$), например, для квадратичных функций

$$f\left(x\right)=\frac{1}{2}x^TAx-b^Tx,\; x\in\mathbb{R}^n,$$

посчитанное, скажем, в 1-норме не может быть меньше $n$. В то время как число обусловленности, посчитанное в евклидовой норме, может при этом равняться 1. Действительно, пусть $\xi=\left(\xi_1,....,\xi_n\right)$, где $\xi_k$ – независимые одинаково распределенные случайные величины $P\left(\xi_k=1/n\right)=P\left(\xi_k=-1/n\right)=1/2$. Тогда, учитывая, что $\|\xi\|_1\equiv1$,

$$\mu_1\left(f\right)\le E_\xi\left[\xi^TA\xi\right]=\frac{1}{n^2}\operatorname{tr}\left(A\right)\le\frac{1}{n}\max_{i,j=1,...,n}\left|A_{ij}\right|=\frac{1}{n}L_1\left(f\right).$$

Кроме того, из замечания 2.2.2 следует также и то, что если $\|\;\|=\|\;\|_1$, то $\tilde{\omega}_n\ge n$. Это обстоятельство также говорит в пользу выбора евклидовой прокс-структуры. Тем не ме-нее, для общности далее мы будем допускать, что используется прокс-структура отличная от евклидовой.

Перепишем формулы (2.2.2), (2.2.3) следующим образом ( $\tilde{\mu} = \mu / \tilde{\omega}_n$ )

$$\varphi_0(x) = V(x, y^0) + \alpha_0 \left[ f(y^0) + \langle \nabla f(y^0), x - y^0 \rangle + \tilde{\mu} V(x, y^0) + h(x) \right],$$

$$\varphi_{k+1}(x) = \varphi_k(x) + \alpha_{k+1} \left[ f(y^{k+1}) + \langle \nabla f(y^{k+1}), x - y^{k+1} \rangle + \tilde{\mu} V(x, y^k) + h(x) \right], \quad (2.2.14)$$

$$A_k = \sum_{i=0}^{k} \alpha_i , \quad \alpha_0 = L^{-1}, \; A_{k+1} \cdot (1 + A_k \tilde{\mu}) = \alpha_{k+1}^2 L , \; k = 0,1,2,... \quad (2.2.15)$$

Сам метод по-прежнему будет иметь вид (4) с $x^0 = u^0 = \arg\min_{x \in Q} \varphi_0(x)$.

**Лемма 2.2.3 (см. [184]).** *Последовательность* $\{\alpha_k\}$, *определяемую формулой (2.2.15), можно задавать рекуррентно*

$$\alpha_{k+1} = \frac{1 + A_k \tilde{\mu}}{2L} + \sqrt{\frac{1 + A_k \tilde{\mu}}{4L^2} + \frac{A_k \cdot (1 + A_k \tilde{\mu})}{L}} , \; A_{k+1} = A_k + \alpha_{k+1} . \quad (2.2.16)$$

*При этом*

$$A_k \ge \frac{1}{L} \left( 1 + \frac{1}{2} \sqrt{\frac{\tilde{\mu}}{L}} \right)^{2k} \ge \exp\left( \frac{k}{2} \sqrt{\frac{\tilde{\mu}}{L}} \right).$$

**Лемма 2.2.4.** *Пусть справедливы предположения 2.2.1, 2.2.2. Тогда для любого* $k = 0,1,2,...$ *имеет место неравенство*

$$A_k F(x^k) \le \varphi_k^* = \min_{x \in Q} \varphi_k(x) = \varphi_k(u^k).$$

**Доказательство.** Доказательство аналогично доказательству леммы 2.2.2. В основе лежит неравенство

$$V(x, y^k) \ge \frac{1}{2} \|x - y^k\|^2 ,$$

с помощью которого ключевая формула (2.2.8) перепишется следующим образом

$$\varphi_{k+1}^* \ge A_{k+1} \cdot \left( \frac{A_k}{A_{k+1}} h(x^k) + \frac{\alpha_{k+1}}{A_{k+1}} h(u^{k+1}) \right) + A_{k+1} f(y^k) +$$

$$+ \langle \; \nabla f(y^{k+1}), \alpha_{k+1} \cdot (u^{k+1} - y^{k+1}) + A_k \cdot (x^k - y^{k+1}) \; \rangle + \frac{(1 + A_k \tilde{\mu})}{2} \|u^{k+1} - u^k\|^2 .$$

Отличие от формулы (2.2.8) в следующем

$$\frac{1}{2} \|u^{k+1} - u^k\|^2 \to \frac{(1 + A_k \tilde{\mu})}{2} \|u^{k+1} - u^k\|^2 .$$

Рассуждая дальше точно также как при доказательстве леммы 2.2.2, получим

$$\varphi_{k+1}^{*} \geq A_{k+1}\left[ f\left(y^{k+1}\right) + \left\langle \nabla f\left(y^{k+1}\right), x^{k+1} - y^{k+1} \right\rangle + \frac{A_{k+1}\cdot\left(1 + A_k\tilde{\mu}\right)}{2\alpha_{k+1}^{2}}\left\| x^{k+1} - y^{k+1} \right\|^{2} + h\left(x^{k+1}\right) \right]. \quad (2.2.17)$$

Из предположения 1 следует, что если $A_{k+1}\cdot\left(1 + A_k\tilde{\mu}\right)/\alpha_{k+1}^{2} \geq L$, то

$$f\left(y^{k+1}\right) + \left\langle \nabla f\left(y^{k+1}\right), x^{k+1} - y^{k+1} \right\rangle + \frac{A_{k+1}\cdot\left(1 + A_k\tilde{\mu}\right)}{2\alpha_{k+1}^{2}}\left\| x^{k+1} - y^{k+1} \right\|^{2} \geq f\left(x^{k+1}\right). \quad (2.2.18)$$

Согласно (2.2.15) $A_{k+1}\cdot\left(1 + A_k\tilde{\mu}\right)/\alpha_{k+1}^{2} = L$, поэтому формула (2.2.18) имеет место. С помощью формулы (2.2.18) формулу (2.2.17) можно переписать следующим образом

$$\varphi_{k+1}^{*} \geq A_{k+1}\left[ f\left(x^{k+1}\right) + h\left(x^{k+1}\right) \right] = A_{k+1}F\left(x^{k+1}\right). \ \blacksquare$$

Из лемм 2.2.3, 2.2.4 получаем следующий результат, означающий, что МТ в сильно выпуклом случае сходится как обычный быстрый градиентный метод (с двумя проекциями), т.е. МТ сходится оптимальным образом для рассматриваемого класса задач.

**Теорема 2.2.2.** *Пусть справедливы предположения 1, 2. Тогда МТ (2.2.14), (2.2.15), (2.2.4) для задачи (2.2.1) сходится согласно оценке*

$$F\left(x^{N}\right) - \min_{x\in Q} F\left(x\right) \leq LR^{2}\exp\left( -\frac{N}{2}\sqrt{\frac{\tilde{\mu}}{L}} \right). \quad (2.2.19)$$

**Доказательство.** Из леммы 2.2.4 следует, что (в третьем неравенстве используется то, что $\tilde{\mu} = \mu/\tilde{\omega}_n$ и сильная выпуклость функции $f\left(x\right)$ – см. формулу (2.2.13) предположения 2.2.2)

$$A_N F\left(x^{N}\right) \leq \min_{x\in Q}\left\{ V\left(x, y^{0}\right) + \sum_{k=0}^{N}\alpha_k\left[ f\left(y^{k}\right) + \left\langle \nabla f\left(y^{k}\right), x - y^{k} \right\rangle + \tilde{\mu}V\left(x, y^{k}\right) + h\left(x\right) \right] \right\} \leq$$

$$\leq V\left(x_*, y^{0}\right) + \sum_{k=0}^{N}\alpha_k \underbrace{\left[ f\left(y^{k}\right) + \left\langle \nabla f\left(y^{k}\right), x_* - y^{k} \right\rangle + \frac{\mu}{2}\left\| x_* - y^{k} \right\|^{2} + h\left(x_*\right) \right]}_{\leq f\left(x_*\right) + h\left(x_*\right)} \leq$$

$$\leq V\left(x_*, y^{0}\right) + \sum_{k=0}^{N}\alpha_k F\left(x_*\right) = R^{2} + A_N F\left(x_*\right). \quad (2.2.20)$$

Для того чтобы получить оценку (2.2.19) осталось подставить нижнюю оценку на $A_N$ из леммы 2.2.3 в формулу (2.2.20). ■

В действительности, выше установлено более сильное утверждение.

**Теорема 2.2.3.** *Пусть справедливы предположения 2.2.1, 2.2.2. Тогда МТ (2.2.14), (2.2.16), (2.2.4) для задачи (2.2.1) сходится согласно оценке*

$$F\left(x^N\right)-\min_{x\in Q}F\left(x\right)\le \min\left\{\frac{4LR^2}{\left(N+1\right)^2}, LR^2\exp\left(-\frac{N}{2}\sqrt{\frac{\mu}{L\tilde{\omega}_n}}\right)\right\}. \tag{2.2.21}$$

Теорема 2.2.3 означает, что МТ (2.2.14), (2.2.16), (2.2.4) непрерывен по параметру $\mu$. К сожалению, при этом в (2.2.16) явно входит этот параметр $\mu$. Если значение этого параметра неизвестно, то с помощью рестартов (см., например, [2, 37]) можно получить оценку (2.2.21) увеличив константы не более чем в 4 раза, т.е. число вычислений градиента $\nabla f\left(x\right)$ (обычно именно это является самым затратным в шаге), необходимых для достижения заданной точности, увеличится не более чем в 4 раза.

Между сильно выпуклым и просто выпуклым случаями имеется глубокая связь, позволяющая, например, получить оценку (2.2.19) с помощью оценки (2.2.11) и наоборот. Другими словами, имея эффективные алгоритмы решения выпуклых / сильно выпуклых задач, можно предложить на их базе алгоритмы решения сильно выпуклых / выпуклых задач. Покажем это, следуя, например, [48] (приводимые далее конструкции давно и хорошо известны, мы здесь их приводим для полноты изложения).

Введем семейство $\mu$-сильно выпуклых в норме $\|\ \|$ задач ($\mu > 0$)

$$F^{\mu}\left(x\right)=F\left(x\right)+\mu V\left(x,y^0\right)\to\min_{x\in Q}. \tag{2.2.22}$$

**Теорема 2.2.4.** *Пусть*

$$\mu\le\frac{\varepsilon}{2V\left(x_*,y^0\right)}=\frac{\varepsilon}{2R^2}, \tag{2.2.23}$$

*и удалось найти $\varepsilon/2$-решение задачи (2.2.22), т.е. нашелся такой $x^N\in Q$, что*

$$F^{\mu}\left(x^N\right)-F_*^{\mu}\le\varepsilon/2.$$

*Тогда*

$$F\left(x^N\right)-\min_{x\in Q}F\left(x\right)=F\left(x^N\right)-F_*\le\varepsilon.$$

**Доказательство.** Действительно,

$$F\left(x^N\right)-F_*\le F^{\mu}\left(x^N\right)-F_*\le F^{\mu}\left(x^N\right)-F_*^{\mu}+\varepsilon/2\le\varepsilon.$$

Здесь использовалось определение $F_*^{\mu}$ и формула (2.2.23)

$$F_*^{\mu}=\min_{x\in Q}\left\{F\left(x\right)+\gamma V\left(x,y^0\right)\right\}\le F\left(x_*\right)+\gamma V\left(x_*,y^0\right)\le F_*+\varepsilon/2. \ \blacksquare$$

Приведем в некотором смысле обратную (рестарт-)конструкцию.

**Теорема 2.2.5.** *Пусть функция* $F(x)$ *–* $\mu$ *-сильно выпуклая в норме* $\|\ \|$*. Пусть точка* $x^{\bar{N}}\left(y^0\right)$ *выдается МТ (2.2.2) – (2.2.4), стартующим с точки* $y^0$*, после*

$$\bar{N} = \sqrt{\frac{8L\omega_n}{\mu}} \tag{2.2.24}$$

*итераций, где*

$$\omega_n = \sup_{x\in Q}\frac{2V\left(x,y^0\right)}{\left\|x-y^0\right\|^2}.$$

*Положим*

$$\left[x^{\bar{N}}\left(y^0\right)\right]^1 = x^{\bar{N}}\left(y^0\right)$$

*и определим по индукции*

$$\left[x^{\bar{N}}\left(y^0\right)\right]^{k+1} = x^{\bar{N}}\left(\left[x^{\bar{N}}\left(y^0\right)\right]^k\right),\ k=1,2,\ldots.$$

*При этом на* $k+1$ *перезапуске (рестарте) также корректируется прокс-функция (считаем, что так определенная функция корректно определена на* $Q$ *с сохранением свойства сильной выпуклости)*

$$d^{k+1}(x) = d\left(x-\left[x^{\bar{N}}\left(y^0\right)\right]^k + y^0\right)\ge 0,$$

*чтобы*

$$d^{k+1}\left(\left[x^{\bar{N}}\left(y^0\right)\right]^k\right)=0,\ \nabla d^{k+1}\left(\left[x^{\bar{N}}\left(y^0\right)\right]^k\right)=0$$

*Тогда*

$$F\left(\left[x^{\bar{N}}\left(y^0\right)\right]^k\right)-F_* \le \frac{\mu\left\|y^0-x_*\right\|^2}{2^{k+1}}. \tag{2.2.25}$$

**Доказательство.** МТ (2.2.2) – (2.2.4) согласно теореме 2.2.1 (см. формулу (2.2.11)) после $\bar{N}$ итераций выдает такой $x^{\bar{N}}$, что

$$\frac{\mu}{2}\left\|x^{\bar{N}}-x_*\right\|^2 \le F\left(x^{\bar{N}}\right)-F_* \le \frac{4LV\left(x_*,y^0\right)}{\bar{N}^2}.$$

Отсюда имеем

$$\left\|x^{\bar{N}}-x_*\right\|^2 \le \frac{8LV\left(x_*,y^0\right)}{\mu\bar{N}^2} \le \frac{1}{2}\left\|y^0-x_*\right\|^2\frac{8L\omega_n}{\mu\bar{N}^2}.$$

Поскольку

$$\bar{N} = \sqrt{\frac{8L}{\mu}\omega_n},$$

то

$$\left\|x^{\bar{N}} - x_*\right\|^2 \le \frac{1}{2}\left\|y^0 - x_*\right\|^2.$$

Повторяя эти рассуждения, по индукции получим

$$F\left(\left[x^{\bar{N}}\left(y^0\right)\right]^k\right) - F_* \le \left(\frac{1}{2}\right)^k \left\|y^0 - x_*\right\|^2 \frac{4L\omega_n}{\bar{N}^2} = \frac{\mu\left\|y^0 - x_*\right\|^2}{2^{k+1}}. \blacksquare$$

**Замечание 2.2.3.** Оценка (2.2.24), (2.2.25) в итоге получается похожей на оценку (2.2.19). Однако в подходе, описанном в теореме 2.2.5, всегда можно добиться, чтобы $\omega_n = O(\ln n)$, в том числе и в случае $\| \; \| = \| \; \|_1$ [228] (следует сравнить в этом случае с оценкой $\tilde{\omega}_n \ge n$, приведенной выше). Кроме того, предложенная в теореме 2.2.5 конструкция позволяет рассматривать более общий класс сильно выпуклых задач, в которых $f(x)$ уже не обязательно – сильно выпуклая функция (см. предположение 2.2.2), что позволяет избавиться в оценках числа обусловленности $L/\mu$ от ограничения, описанного в замечании 2.2.2. В ряде приложений эти степени свободы оказываются чрезвычайно важными [48] (см. также следующий раздел). Однако, как показывают численные эксперименты (проводимые в случае евклидовой прокс-структуры и с априорно известной точной оценкой параметра $\mu$), подход, описанный в теореме 2.2.5, может проигрывать в скорости МТ (2.2.14), (2.2.16), (2.2.4) один-два порядка, т.е. для достижения той же точности подход из теоремы 2.2.5 может потребовать до 100 раз больше арифметических операций. Из оценок это никак не следует. Но если метод из теоремы 2.2.5 работает по полученным верхним оценкам, то МТ (2.2.14), (2.2.16), (2.2.4) (с критерием останова связанным с контролем малости нормы градиентного отображения [93, 181]) из-за отсутствия рестартов (т.е. необходимости делать предписанное число итераций) может остановиться раньше, что и происходит на практике. К тому же МТ (2.2.14), (2.2.16), (2.2.4) еще и непрерывен по параметру $\mu$. Заметим также, что недавно появились две работы, которые удачно уточняют описанные здесь конструкции регуляризации и рестартов arXiv:1603.05643, arXiv:1609.07358.

Из написанного выше, кажется, что в процедуре рестартов (теорема 2.2.5) можно использовать вместо предписанного числа итераций $\bar{N}$ на каждом рестарте какой-нибудь критерий останова. В частности, дожидаться, когда норма (или квадрат нормы) градиента (а в общем случае, когда минимум достигается не в точке экстремума, – норма градиент-ного отображения) уменьшиться вдвое (см., например, раздел 4.3 главы 4). С таким критерием останова нет необходимости делать предписанного числа итераций на каждом рестарте. Однако пока не известен способ рассуждений (это же проблема “вылезает” и в первых пунктах раздела 5.1 главы 5), который позволял бы показать, что такая процедура сохраняет при перенесении оптимальность оценок (метод, работающий оптимально не в сильно выпуклом случае, порождает оптимальный метод и в сильно выпуклом случае). Впрочем, имеются различные эффективные на практике способы более раннего выхода с каждого рестарта (см., например, [286]), позволяющие (в случае задач безусловной оптимизации с евклидовой прокс-структурой и отсутствием композитного члена) ускорить описанную выше конструкцию на порядок.

### 2.2.4 Универсальный метод треугольника

В ряде приложений значение константы $L$, необходимой МТ для работы (см. формулу (2.2.16)), не известно. Однако, как следует из формул (2.2.10), (2.2.18), знание константы $L$ необязательно, если разрешается на одной итерации запрашивать значение функции в нескольких точках. Опишем соответствующий адаптивный вариант МТ (2.2.14), (2.2.16), (2.2.4) (АМТ) (см., например, [266]).

Положим $A_0 = \alpha_0 = 1/L_0^0$ – параметр метода (считаем здесь и везде в дальнейшем $L_0^0 \le L$, иначе во всех приводимых далее оценках следует полагать $L := \max\left\{L_0^0, L\right\}$; обычно полагают $L_0^0 = 1$),

$$k = 0,\ j_0 = 0;\ x^0 = u^0 = \arg\min_{x\in Q} \varphi_0(x).$$

До тех пор пока

$$f\left(x^0\right) > f\left(y^0\right) + \left\langle \nabla f\left(y^0\right), x^0 - y^0 \right\rangle + \frac{L_0^{j_0}}{2}\left\|x^0 - y^0\right\|^2,$$

где

$$x^0 := u^0 := \arg\min_{x\in Q} \varphi_0(x),\ (A_0 :=)\,\alpha_0 := \frac{1}{L_0^{j_0}},$$

выполнять

$$j_0 := j_0 + 1;\ L_0^{j_0} := 2^{j_0} L_0^0.$$

<u>**Адаптивный Метод Треугольника**</u>

1. $L_{k+1}^0 = L_k^{j_k}/2,\ j_{k+1} = 0.$

2. $\begin{cases} \alpha_{k+1} := \dfrac{1+A_k\tilde{\mu}}{2L_{k+1}^{j_{k+1}}} + \sqrt{\dfrac{1+A_k\tilde{\mu}}{4\left(L_{k+1}^{j_{k+1}}\right)^2} + \dfrac{A_k\cdot(1+A_k\tilde{\mu})}{L_{k+1}^{j_{k+1}}}},\ A_{k+1} := A_k + \alpha_{k+1}; \\ y^{k+1} := \dfrac{\alpha_{k+1}u^k + A_k x^k}{A_{k+1}},\ u^{k+1} := \arg\min\limits_{x\in Q} \varphi_{k+1}(x),\ x^{k+1} := \dfrac{\alpha_{k+1}u^{k+1} + A_k x^k}{A_{k+1}}. \end{cases}$ (*)

До тех пор пока

$$f\left(y^{k+1}\right) + \left\langle \nabla f\left(y^{k+1}\right), x^{k+1} - y^{k+1}\right\rangle + \frac{L_{k+1}^{j_{k+1}}}{2}\left\|x^{k+1} - y^{k+1}\right\|^2 < f\left(x^{k+1}\right),$$

выполнять

$$j_{k+1} := j_{k+1} + 1;\ L_{k+1}^{j_{k+1}} = 2^{j_{k+1}} L_{k+1}^0;\ (*).$$

3. Если не выполнен критерий останова, то $k := k+1$ и **go to** 1.

В качестве критерия останова, например, можно брать условие

$$\left\| x^{k+1} - \arg\min_{x\in Q}\left\{ \left\langle \nabla f\left(x^{k+1}\right), x - x^{k+1}\right\rangle + \frac{L_{k+1}^{j_{k+1}}}{2}\left\|x - x^{k+1}\right\|^2 \right\} \right\| \le \tilde{\varepsilon}.$$

Здесь и везде в дальнейшем под “До тех пор пока … выполнять …” подразумевается, что после каждого $j_{k+1} := j_{k+1} + 1$ при следующей проверке условия выхода из этого цикла должным образом меняется не только $L_{k+1}^{j_{k+1}}$, но и $x^{k+1}$, $y^{k+1}$, также входящие в это условие.

**Теорема 2.2.6.** *Пусть справедливы предположения 2.2.1, 2.2.2. Тогда АМТ для задачи (2.2.1) сходится согласно оценке*

$$F\left(x^N\right) - \min_{x\in Q} F(x) \le \min\left\{ \frac{8LR^2}{(N+1)^2}, 2LR^2 \exp\left( -\frac{N}{2}\sqrt{\frac{\mu}{2L\tilde{\omega}_n}} \right) \right\}. \qquad (2.2.26)$$

*При этом среднее число вычислений значения функции на одной итерации будет* $\approx 4$*, а градиента функции* $\approx 2$.

**Доказательство.** Нетривиальным в виду оценки (2.2.21) и свойства, что все $L_k^{j_k} \le 2L$, представляется только последняя часть формулировки теоремы. Докажем именно её. Оценим общее число обращений за значениями функции (аналогично получается оценка общего числа обращений за значением градиента функции)

$$\sum_{k=1}^N 2\left(j_k + 1\right) = \sum_{k=1}^N 2\left[\left(j_k - 1\right) + 2\right] = \sum_{k=1}^N 2\left[\log_2\left(\frac{L_k^{j_k}}{L_{k-1}^{j_{k-1}}}\right) + 2\right] = 4N + \log_2\left(\frac{L_N^{j_N}}{L_0^0}\right) \le 4N + \log_2\left(\frac{2L}{L_0^0}\right).$$

Деля обе части на $N$, получим в правой части приблизительно 4. ■

В действительности, оценка (2.2.26) оказывается, как правило, сильно завышенной, поскольку метод адаптивно настраивается на константу Липшица градиента $L$ на данном участке своего пребывания, а константа $L$, входящая в оценку (2.2.26), соответствует (согласно предположению 2.2.1) самому плохому случаю (самому плохому участку).

Предположим теперь, что по каким-то причинам невозможно получить точные значения функции и градиента. Тогда соотношение (аналогичное неравенство выписывается и при $k=0$, см. начало этого пункта)

$$f\left(y^{k+1}\right)+\left\langle\nabla f\left(y^{k+1}\right),x^{k+1}-y^{k+1}\right\rangle+\frac{L_{k+1}^{j_{k+1}}}{2}\left\|x^{k+1}-y^{k+1}\right\|^2\geq f\left(x^{k+1}\right)$$

может не выполниться не при каком $L_{k+1}^{j_{k+1}}$. Допустим, однако, что при этом имеет место

$$f\left(y^{k+1}\right)+\left\langle\nabla f\left(y^{k+1}\right),x^{k+1}-y^{k+1}\right\rangle+\frac{L}{2}\left\|x^{k+1}-y^{k+1}\right\|^2+\frac{\alpha_{k+1}}{A_{k+1}}\varepsilon\geq f\left(x^{k+1}\right).$$

Тогда заменим в АМТ соответствующую часть шага 2 на

$$f\left(y^{k+1}\right)+\left\langle\nabla f\left(y^{k+1}\right),x^{k+1}-y^{k+1}\right\rangle+\frac{L_{k+1}^{j_{k+1}}}{2}\left\|x^{k+1}-y^{k+1}\right\|^2+\frac{\alpha_{k+1}}{A_{k+1}}\varepsilon\geq f\left(x^{k+1}\right). \tag{2.2.27}$$

**Теорема 2.2.7.** *Пусть справедливо предположение 2.2.2 и существует такое число* $L>0$, *что любого* $k=1,\dots,N$ *справедливо неравенство*

$$f\left(y^{k+1}\right)+\left\langle\nabla f\left(y^{k+1}\right),x^{k+1}-y^{k+1}\right\rangle+\frac{L}{2}\left\|x^{k+1}-y^{k+1}\right\|^2+\frac{\alpha_{k+1}}{A_{k+1}}\varepsilon\geq f\left(x^{k+1}\right). \tag{2.2.28}$$

*Тогда АМТ с (2.2.27) для задачи (2.2.1) сходится согласно оценке*

$$F\left(x^N\right)-\min_{x\in Q}F\left(x\right)\leq\min\left\{\frac{8LR^2}{\left(N+1\right)^2},2LR^2\exp\left(-\frac{N}{2}\sqrt{\frac{\mu}{2L\tilde{\omega}_n}}\right)\right\}+\varepsilon. \tag{2.2.29}$$

*При этом среднее число вычислений значения функции на одной итерации будет* $\approx 4$, *а градиента функции* $\approx 2$

**Доказательство.** Ключевым элементом в доказательстве является следующее уточнение леммы 4

$$A_kF\left(x^k\right)\leq\varphi_k^*+A_k\varepsilon, \tag{2.2.30}$$

из которого будет следовать формула (2.2.29). Чтобы доказать (2.2.30) будем рассуждать по индукции. База индукции $k=0$ очевидна. Итак, по предположению индукции

$$A_kF\left(x^k\right)-A_k\varepsilon\leq\varphi_k^*=\varphi_k\left(u^k\right),$$

поэтому

$$\varphi_k\left(u^{k+1}\right)\geq\varphi_k\left(u^k\right)+\frac{1+A_k\tilde{\mu}}{2}\left\|u^{k+1}-u^k\right\|^2\geq A_kF\left(x^k\right)-A_k\varepsilon+\frac{1+A_k\tilde{\mu}}{2}\left\|u^{k+1}-u^k\right\|^2.$$

Отсюда

$$\varphi_{k+1}^{*} \ge A_{k+1}\cdot\left(\frac{A_k}{A_{k+1}}h\left(x^k\right)+\frac{\alpha_{k+1}}{A_{k+1}}h\left(u^{k+1}\right)\right)+A_{k+1}f\left(y^k\right)-A_k\varepsilon$$

$$+\left\langle \nabla f\left(y^{k+1}\right),\alpha_{k+1}\cdot\left(u^{k+1}-y^{k+1}\right)+A_k\cdot\left(x^k-y^{k+1}\right)\right\rangle+\frac{\left(1+A_k\tilde{\mu}\right)}{2}\left\|u^{k+1}-u^k\right\|^2.$$

Следовательно,

$$\varphi_{k+1}^{*}+A_k\varepsilon \ge A_{k+1}\left[f\left(y^{k+1}\right)+\left\langle\nabla f\left(y^{k+1}\right),x^{k+1}-y^{k+1}\right\rangle+\frac{A_{k+1}\cdot\left(1+A_k\tilde{\mu}\right)}{2\alpha_{k+1}^2}\left\|x^{k+1}-y^{k+1}\right\|^2+h\left(x^{k+1}\right)\right].$$

Отсюда и из условия (2.2.27), $A_{k+1}\cdot\left(1+A_k\tilde{\mu}\right)/\alpha_{k+1}^2=L_{k+1}^{j_{k+1}}$ (с учетом (2.2.28)) получаем

$$\varphi_{k+1}^{*}+A_{k+1}\varepsilon=\varphi_{k+1}^{*}+A_k\varepsilon+\alpha_{k+1}\varepsilon \ge A_{k+1}F\left(x^{k+1}\right),\ L_{k+1}^{j_{k+1}}\le 2L. \ ■$$

В действительности, выше установлено более сильное утверждение – в оценке (2.2.29) можно улучшить константу $L$.

**Теорема 2.2.8.** *Пусть справедливо предположение 2. Тогда АМТ с (2.2.27) для задачи (2.2.1) сходится согласно оценке*

$$F\left(x^N\right)-\min_{x\in Q}F\left(x\right)\le\frac{R^2}{A_N}+\varepsilon\le\min\left\{\frac{4LR^2}{\left(N+1\right)^2},LR^2\exp\left(-\frac{N}{2}\sqrt{\frac{\mu}{L\tilde{\omega}_n}}\right)\right\}+\varepsilon, \qquad (2.2.31)$$

*где* $L=\max_{k=0,\dots,N}L_k^{j_k}$. *При этом среднее число вычислений значения функции на одной итерации будет* $\approx 4$, *а градиента функции* $\approx 2$.

Попробуем (подобно [48, 274]) сыграть на условии (2.2.27), искусственно вводя неточность.

**Лемма 2.2.5** (см. [183, 274]). *Пусть*

$$\left\|\nabla f\left(y\right)-\nabla f\left(x\right)\right\|_*\le L_\nu\left\|y-x\right\|^\nu \qquad (2.2.32)$$

*при некотором* $\nu\in\left[0,1\right]$. *Тогда*

$$f\left(y\right)+\left\langle\nabla f\left(y\right),x-y\right\rangle+\frac{L}{2}\left\|x-y\right\|^2+\delta\ge f\left(x\right),\ L=L_\nu\left[\frac{L_\nu}{2\delta}\frac{1-\nu}{1+\nu}\right]^{\frac{1-\nu}{1+\nu}}. \qquad (2.2.33)$$

Основным результатом данного пункта является описание и последующая оценка скорости сходимости нового варианта универсального метода [274] на базе МТ (УМТ).

Положим

$$A_0=\alpha_0=1/L_0^0,\ k=0,\ j_0=0;\ x^0=u^0=\arg\min_{x\in Q}\varphi_0\left(x\right).$$

До тех пор пока

$$f\left(x^0\right) > f\left(y^0\right) + \left\langle \nabla f\left(y^0\right), x^0 - y^0 \right\rangle + \frac{L_0^{j_0}}{2}\left\|x^0 - y^0\right\|^2 + \frac{\alpha_0}{2A_0}\varepsilon,$$

где

$$x^0 := u^0 := \arg\min_{x\in Q} \varphi_0\left(x\right),\ \left(A_0 :=\right)\alpha_0 := \frac{1}{L_0^{j_0}},$$

выполнять

$$j_0 := j_0 + 1;\ L_0^{j_0} := 2^{j_0} L_0^0 .$$

**Универсальный Метод Подобных Треугольников**

1. $L_{k+1}^0 = L_k^{j_k}/2,\ j_{k+1} = 0$.

2. $$\begin{cases} \alpha_{k+1} := \dfrac{1 + A_k\tilde{\mu}}{2L_{k+1}^{j_{k+1}}} + \sqrt{\dfrac{1 + A_k\tilde{\mu}}{4\left(L_{k+1}^{j_{k+1}}\right)^2} + \dfrac{A_k \cdot \left(1 + A_k\tilde{\mu}\right)}{L_{k+1}^{j_{k+1}}}},\ A_{k+1} := A_k + \alpha_{k+1}; \\ y^{k+1} := \dfrac{\alpha_{k+1}u^k + A_k x^k}{A_{k+1}},\ u^{k+1} := \arg\min_{x\in Q}\varphi_{k+1}\left(x\right),\ x^{k+1} := \dfrac{\alpha_{k+1}u^{k+1} + A_k x^k}{A_{k+1}}. \end{cases} \qquad (*)$$

До тех пор пока

$$f\left(y^{k+1}\right) + \left\langle \nabla f\left(y^{k+1}\right), x^{k+1} - y^{k+1} \right\rangle + \frac{L_{k+1}^{j_{k+1}}}{2}\left\|x^{k+1} - y^{k+1}\right\|^2 + \frac{\alpha_{k+1}}{2A_{k+1}}\varepsilon < f\left(x^{k+1}\right),$$

выполнять

$$j_{k+1} := j_{k+1} + 1;\ L_{k+1}^{j_{k+1}} = 2^{j_{k+1}} L_{k+1}^0;\ (*).$$

3. Если не выполнен критерий останова, то $k := k+1$ и **go to** 1.

**Теорема 2.2.9.** *Пусть выполняется условие (2.2.32) хотя бы для* $\nu = 0$, *и справедливо предположение 2.2.2 с* $\mu \ge 0$ *(допускается брать* $\mu = 0$*). Тогда УМТ для задачи (2.2.1) сходится согласно оценке*

$$F\left(x^N\right) - \min_{x\in Q} F\left(x\right) \le \varepsilon,$$

$$N \approx \min\left\{ \inf_{\nu\in[0,1]} \left( \frac{L_\nu \cdot \left(16R\right)^{1+\nu}}{\varepsilon} \right)^{\frac{2}{1+3\nu}},\ \inf_{\nu\in[0,1]} \left\{ \left( \frac{8L_\nu^{\frac{2}{1+\nu}}\tilde{\omega}_n}{\mu\varepsilon^{\frac{1-\nu}{1+\nu}}} \right)^{\frac{1+\nu}{1+3\nu}} \ln^{\frac{2+2\nu}{1+3\nu}} \left( \frac{16L_\nu^{\frac{4+6\nu}{1+\nu}} R^2}{\left(\mu/\tilde{\omega}_n\right)^{\frac{1+\nu}{1+3\nu}} \varepsilon^{\frac{5+7\nu}{2+6\nu}}} \right) \right\} \right\}. \qquad (2.2.34)$$

*При этом среднее число вычислений значения функции на одной итерации будет* $\approx 4$, *а градиента функции* $\approx 2$.

**Доказательство.** Рассмотрим два случая, когда $\mu \geq 0$ – мало: $\mu \ll \varepsilon/(2R^2)$, $\mu$ – велико: $\mu \gg \varepsilon/(2R^2)$, см. формулу (2.2.23).

В первом случае будем считать, что

$$A_{k+1}/\alpha_{k+1}^2 \approx A_{k+1} \cdot (1 + A_k \tilde{\mu})/\alpha_{k+1}^2 = L_{k+1}^{j_{k+1}}, \text{ т.е. } \varepsilon \frac{\alpha_{k+1}}{2A_{k+1}} \approx \frac{\varepsilon}{2}\sqrt{\frac{1}{L_{k+1}^{j_{k+1}} A_{k+1}}}, \tag{2.2.35}$$

а во втором случае

$$A_{k+1}^2 \tilde{\mu}/\alpha_{k+1}^2 \approx A_{k+1} \cdot (1 + A_k \tilde{\mu})/\alpha_{k+1}^2 = L_{k+1}^{j_{k+1}}, \text{ т.е. } \varepsilon \frac{\alpha_{k+1}}{2A_{k+1}} \approx \frac{\varepsilon}{2}\sqrt{\frac{\tilde{\mu}}{L_{k+1}^{j_{k+1}}}}. \tag{2.2.36}$$

Из формулы (2.2.31) (см. теорему 2.2.8) имеем

$$\frac{R^2}{A_N} + \frac{\varepsilon}{2} \approx \varepsilon,$$

т.е. $A_N \approx 2R^2/\varepsilon$, а также ( $L = \max\limits_{k=0,\ldots,N} L_k^{j_k}$ )

$$N^2 \approx \frac{8LR^2}{\varepsilon} \quad \text{(в первом случае)}, \tag{2.2.37}$$

$$N^2 \approx 4\frac{L}{\tilde{\mu}} \ln^2\left(\frac{2LR^2}{\varepsilon}\right) \quad \text{(во втором случае)}. \tag{2.2.38}$$

Из формул (2.2.33), (2.2.35) – (2.2.38) имеем, что в первом случае

$$L \leq 2L_\nu \left[ \frac{L_\nu}{2\frac{\varepsilon}{2}\sqrt{\frac{1}{LA_N}}} \frac{1-\nu}{1+\nu} \right]^{\frac{1-\nu}{1+\nu}} \leq 2L_\nu \left[ \frac{L_\nu \sqrt{LA_N}}{\varepsilon} \frac{1-\nu}{1+\nu} \right]^{\frac{1-\nu}{1+\nu}} \leq 2L_\nu^{\frac{2}{1+\nu}} \left[ \frac{N}{2\varepsilon} \frac{1-\nu}{1+\nu} \right]^{\frac{1-\nu}{1+\nu}}, \tag{2.2.39}$$

а во втором случае

$$L \leq 2L_\nu \left[ \frac{L_\nu}{2\frac{\varepsilon}{2}\sqrt{\frac{\tilde{\mu}}{L}}} \frac{1-\nu}{1+\nu} \right]^{\frac{1-\nu}{1+\nu}} \leq 2L_\nu \left[ \frac{L_\nu \sqrt{L/\tilde{\mu}}}{\varepsilon} \frac{1-\nu}{1+\nu} \right]^{\frac{1-\nu}{1+\nu}} \leq 2L_\nu^{\frac{2}{1+\nu}} \left[ \frac{N}{2\varepsilon} \frac{1-\nu}{1+\nu} \right]^{\frac{1-\nu}{1+\nu}}. \tag{2.2.40}$$

Подставляя (2.2.39) в (2.2.37), а (2.2.40) в (2.2.38) и учитывая, что параметр $\nu \in [0,1]$ можно выбирать произвольно (допускается, что $L_\nu = \infty$ при некоторых $\nu$ – важно, чтобы существовало хотя бы одно значение $\nu$ при котором $L_\nu < \infty$; по условию $L_0 < \infty$), получим соответственно,

$$N^2 \approx \frac{16L_\nu^{\frac{2}{1+\nu}}\left[\frac{N}{2\varepsilon}\frac{1-\nu}{1+\nu}\right]^{\frac{1-\nu}{1+\nu}}R^2}{\varepsilon} \Rightarrow N^{\frac{1+3\nu}{1+\nu}} \approx \frac{16L_\nu^{\frac{2}{1+\nu}}R^2}{\varepsilon^{\frac{2}{1+\nu}}} \Rightarrow N \approx \inf_{\nu\in[0,1]}\left(\frac{L_\nu\cdot(16R)^{1+\nu}}{\varepsilon}\right)^{\frac{2}{1+3\nu}}, \quad (2.2.41)$$

$$N^2 \approx \frac{8L_\nu^{\frac{2}{1+\nu}}\left[\frac{N}{2\varepsilon}\frac{1-\nu}{1+\nu}\right]^{\frac{1-\nu}{1+\nu}}}{\tilde{\mu}}\ln^2\left(\frac{2L_\nu^2R^2N}{\varepsilon^{3/2}}\right) \Rightarrow N^{\frac{1+3\nu}{1+\nu}} \approx \frac{8L_\nu^{\frac{2}{1+\nu}}}{\tilde{\mu}\varepsilon^{\frac{1-\nu}{1+\nu}}}\ln^2\left(\frac{2L_\nu^2R^2N}{\varepsilon^{3/2}}\right) \Rightarrow$$

$$\Rightarrow N \approx \left(\frac{8L_\nu^{\frac{2}{1+\nu}}}{\tilde{\mu}\varepsilon^{\frac{1-\nu}{1+\nu}}}\right)^{\frac{1+\nu}{1+3\nu}}\ln^{\frac{2+2\nu}{1+3\nu}}\left(\frac{2L_\nu^2R^2N}{\varepsilon^{3/2}}\right) \Rightarrow N \approx \inf_{\nu\in[0,1]}\left\{\left(\frac{8L_\nu^{\frac{2}{1+\nu}}}{\tilde{\mu}\varepsilon^{\frac{1-\nu}{1+\nu}}}\right)^{\frac{1+\nu}{1+3\nu}}\ln^{\frac{2+2\nu}{1+3\nu}}\left(\frac{16L_\nu^{\frac{4+6\nu}{1+\nu}}R^2}{\tilde{\mu}^{\frac{1+\nu}{1+3\nu}}\varepsilon^{\frac{5+7\nu}{2+6\nu}}}\right)\right\}. \quad (2.2.42)$$

Из формул (2.2.41), (2.2.42) получаем оценку (2.2.34). ■

Более аккуратные рассуждения позволяют в несколько раз уменьшить константы, входящие в оценку (2.2.34). Оценка (2.2.34) согласуется с нижними оценками для соответствующих классов задач [212].

Практические эксперименты, проведенные А.Ю. Горновым с УМТ (при $\mu = 0$) показали, что этот метод работает подобно другим версиям универсального метода [48, 274]. К сожалению, и в такой (предложенной в статье – УМТ) модификации метод типично проигрывает на практике методу сопряженных градиентов на квадратичных задачах безусловной оптимизации. Впрочем, то что обычные (не универсальные) быстрые градиентные методы проигрывают на таких задачах методу сопряженных градиентов было известно и ранее [285]. Новое здесь то, что универсальные быстрые градиентные методы также проигрывают.

**Замечание 2.2.4.** Все описанные методы могут быть распространены (насколько нам известно, это пока еще не сделано в общем случае) на задачи условной оптимизации. Для этого сначала, следуя п. 2.3 [93], стоит рассмотреть минимаксную задачу (ввести правильную лианеризацию исходного функционала и градиентное отображение). Далее использовать идею метода нагруженных функционалов п. 2.3.4 [93], приводящую к рестартам по неизвестному параметру п. 2.3.5 [93] (оптимальное значение функционала задачи), введение которого, позволяет свести задачу условной оптимизации к минимаксной. Дополнительная плата за такое “введение” (т.е. за рестарты) будет всего лишь логарифмическая, и с точностью до этой “платы” оценки (на число итераций) будут оптимальными. Однако, к сожалению, итерации метода получаются слишком дорогими. В отличие от п. 2.3 [93] сложность итерации растет с ростом её номера. От этого можно избавиться (с точностью до рестартов), используя вместо линейки описанных в этом разделе методов, быстрый градиентный метод из следующего раздела. Другой способ состоит в перезаписи исходного метода в немного другом виде (см. дипломную работу Александра Тюрина, ФКН ВШЭ, 2016). А именно с заменой $u^{k+1} = \arg\min_{x \in Q} \varphi_{k+1}(x)$ на

$$u^{k+1} = \arg\min_{x \in Q} \left\{ \alpha_{k+1} \cdot \left( \left\langle \nabla f\left(y^{k+1}\right), x - u^k \right\rangle + h(x) \right) + V\left(x, u^k\right) \right\}.$$

Описанный план недавно был независимо реализован G. Lan'ом для другого метода https://pdfs.semanticscholar.org/032b/c9cd7068b053f1506907c30938cc80f88033.pdf .

**Замечание 2.2.5 (А.И. Тюрин).** Все описанные методы могут быть распространены на случай, когда шаг $u^{k+1} = \arg\min_{x \in Q} \varphi_{k+1}(x)$ может быть осуществлен с погрешностью. Соответствующие выкладки практически дословно повторяют рассуждения п. 5.5 [259].

**Замечание 2.2.6.** В большинстве приложений "стоимость" (время) получения от оракула (роль которого, как правило, играют нами же написанные подпрограммы вычисления градиента) градиента функционала заметно превышает время, затрачиваемое на то, чтобы сделать шаг итерации, исходя из выданного оракулом вектора. Желание сбалансировать это рассогласование (усложнить итерации, сохранив при этом старый порядок сложности, и выиграть за счет этого в сокращении числа итераций), привело к возникновению композитной оптимизации [266], в которой (аддитивная) часть функционала задачи переносится без лианеризации (запроса градиента) в итерации. Другой способ перенесения части сложности задачи на итерации был описан в замечании 2.2.4. Здесь остается еще много степеней свободы, позволяющих играть на том насколько "дорогой" будет оракул и соответствующая (этому оракулу) "процедура проектирования", и том сколько (внешних) итераций потребуется методу для достижения заданной точности. В частности, если обращение к оракулу за градиентом и последующее проектирование требуют, в свою очередь, решения вспомогательных оптимизационных задач, то можно "сыграть" на том, насколько точно надо решать эти вспомогательные задачи, пытаясь найти "золотую середину" между стоимостью итерации и числом итераций (см. замечание 2.2.5). Также можно сыграть и на том, как выделять эти вспомогательные задачи. Другими словами, что понимать под оракулом и под итерацией метода. Общая идея "разделяй и властвуй", применительно к численным методам выпуклой оптимизации может принимать довольно неожиданные и при этом весьма эффективные формы (ярким примером являются методы внутренней точки [93, 259]). Разные варианты описанной игры в связи с транспортно-сетевыми приложениями уже разбирались нами в других работах [42, 43, 44]. Соответствующие примеры также вошли в раздел 1.6 главы 1, разделы 2.1, 2.3 этой главы 2 и раздел 5.1 главы 5.

### 2.2.5 Универсальный метод треугольника для задач стохастической композитной оптимизации

Предположим теперь, что вместо настоящих градиентов нам доступны только стохастические градиенты $\nabla f(x) \to \nabla f(x,\xi)$ (см., например, [39, 181]).

**Предположение 2.2.3.** *Для всех* $x \in Q$

$$E_{\xi}\left[\nabla f(x,\xi)\right] = \nabla f(x) \text{ и } E_{\xi}\left[\left\|\nabla f(x,\xi) - \nabla f(x)\right\|_*^2\right] \le D. \tag{2.2.43}$$

В ряде приложений бывает полезно рассматривать модификацию условия (2.2.43)

$$L \cdot E_{\xi}\left[\max_{x,y\in Q}\left\{\left\langle \nabla f(y,\xi) - \nabla f(y), x-y\right\rangle - \frac{L}{2}\|x-y\|^2\right\}\right] \le \tilde{D}.$$

Далее приводится стохастический вариант УМТ (СУМТ). По-видимому, это первая попытка перенесения универсального метода на задачи стохастической оптимизации.

Предварительно введём обозначения

$$\tilde{\nabla}^m f(x) = \frac{1}{m}\sum_{k=1}^{m} \nabla f\left(x,\xi^k\right), \tag{2.2.44}$$

где $\xi^k$ – независимые одинаково распределенные (так же как $\xi$) случайные величины. В принципе, можно было бы аналогично ввести

$$\bar{f}^{\tilde{m}}(x) = \frac{1}{\tilde{m}}\sum_{k=1}^{\tilde{m}} f\left(x,\xi^k\right)$$

и распространить приведенные далее результаты на случай, когда и значение функции $f(x)$ необходимо оценивать. Однако для большей наглядности ограничимся далее случаем, когда значение $f(x)$ можно точно посчитать (т.е. нет необходимости его оценивать).

Переопределим последовательность (2.2.14)

$$\varphi_0(x) = V\left(x,y^0\right) + \alpha_0\left[f\left(y^0\right) + \left\langle \tilde{\nabla}^m f\left(y^0\right), x-y^0\right\rangle + \tilde{\mu}V\left(x,y^0\right) + h(x)\right],$$

$$\varphi_{k+1}(x) = \varphi_k(x) + \alpha_{k+1}\left[f\left(y^{k+1}\right) + \left\langle \tilde{\nabla}^m f\left(y^{k+1}\right), x-y^{k+1}\right\rangle + \tilde{\mu}V\left(x,y^k\right) + h(x)\right].$$

Положим

$$A_0 = \alpha_0 = 1/L_0^0,\ m_0 := \frac{2DA_0}{L_0^0\alpha_0\varepsilon},\ k=0,\ j_0 = 0;\ x^0 = u^0 = \arg\min_{x\in Q}\varphi_0(x).$$

До тех пор пока

$$f\left(x^0\right) > f\left(y^0\right) + \left\langle \tilde{\nabla}^{m_0} f\left(y^0\right), x^0 - y^0 \right\rangle + \frac{L_0^{j_0}}{2}\left\|x^0 - y^0\right\|^2 + \frac{3\alpha_{k+1}}{2A_{k+1}}\varepsilon,$$

где

$$x^0 := u^0 := \arg\min_{x\in Q} \varphi_0\left(x\right) =$$

$$= \arg\min_{x\in Q}\left\{V\left(x, y^0\right) + \alpha_0\left[f\left(y^0\right) + \left\langle \tilde{\nabla}^{m_0} f\left(y^0\right), x - y^0 \right\rangle + \tilde{\mu}V\left(x, y^0\right) + h\left(x\right)\right]\right\},$$

$$\left(A_0 :=\right)\alpha_0 := \frac{1}{L_0^{j_0}},\ m_0 := \frac{2DA_0}{L_0^{j_0}\alpha_0\varepsilon},$$

выполнять

$$j_0 := j_0 + 1;\ L_0^{j_0} := 2^{j_0} L_0^0.$$

**<u>Стохастический Универсальный Метод Подобных Треугольников</u>**

1. $L_{k+1}^0 = L_k^{j_k}/2,\ j_{k+1} = 0$.

2. $$\begin{cases} \alpha_{k+1} := \dfrac{1 + A_k\tilde{\mu}}{2L_{k+1}^{j_{k+1}}} + \sqrt{\dfrac{1 + A_k\tilde{\mu}}{4\left(L_{k+1}^{j_{k+1}}\right)^2} + \dfrac{A_k\cdot\left(1 + A_k\tilde{\mu}\right)}{L_{k+1}^{j_{k+1}}}},\ A_{k+1} := A_k + \alpha_{k+1},\ m_{k+1} := \dfrac{2DA_{k+1}}{L_{k+1}^{j_{k+1}}\alpha_{k+1}\varepsilon}; \\ y^{k+1} := \dfrac{\alpha_{k+1}u^k + A_k x^k}{A_{k+1}},\ u^{k+1} := \arg\min_{x\in Q}\varphi_{k+1}\left(x\right),\ x^{k+1} := \dfrac{\alpha_{k+1}u^{k+1} + A_k x^k}{A_{k+1}}. \end{cases} \quad (*)$$

До тех пор пока

$$f\left(y^{k+1}\right) + \left\langle \tilde{\nabla}^{m_{k+1}} f\left(y^{k+1}\right), x^{k+1} - y^{k+1} \right\rangle + \frac{L_{k+1}^{j_{k+1}}}{2}\left\|x^{k+1} - y^{k+1}\right\|^2 + \frac{3\alpha_{k+1}}{2A_{k+1}}\varepsilon < f\left(x^{k+1}\right),$$

выполнять

$$j_{k+1} := j_{k+1} + 1;\ L_{k+1}^{j_{k+1}} = 2^{j_{k+1}} L_{k+1}^0;\ (*).$$

3. Если не выполнен критерий останова, то $k := k + 1$ и **go to** 1.

**Теорема 2.2.10.** *Пусть выполняется условие (2.2.32) хотя бы для* $\nu = 0$*, справедливо предположение 2.2.2 с* $\mu \ge 0$ *(допускается брать* $\mu = 0$*), справедливо предположение 2.2.3. Тогда СУМТ для задачи (2.2.1) сходится согласно оценке*

$$E\left[F\left(x^N\right)\right] - \min_{x\in Q} F\left(x\right) \le 2\varepsilon,$$

$$N \approx \min\left\{ \inf_{\nu\in[0,1]} \left( \frac{L_\nu \cdot (32R)^{1+\nu}}{\varepsilon} \right)^{\frac{2}{1+3\nu}}, \inf_{\nu\in[0,1]} \left\{ \left( \frac{16 L_\nu^{\frac{2}{1+\nu}} \tilde{\omega}_n}{\mu \varepsilon^{\frac{1-\nu}{1+\nu}}} \right)^{\frac{1+\nu}{1+3\nu}} \ln^{\frac{2+2\nu}{1+3\nu}} \left( \frac{32 L_\nu^{\frac{4+6\nu}{1+\nu}} R^2}{(\mu/\tilde{\omega}_n)^{\frac{1+\nu}{1+3\nu}} \varepsilon^{\frac{5+7\nu}{2+6\nu}}} \right) \right\} \right\}. \qquad (2.2.45)$$

*Оценка (2.2.45) – это оценка числа итераций. При этом среднее число вычислений значения функции на одной итерации будет* $\approx 4$. *Однако не менее интересна оценка числа обращений за стохастическим градиентом*

$$\mathrm{Q} \approx 2\min\left\{ \frac{4DR^2}{\varepsilon^2}, \frac{2D\tilde{\omega}_n}{\mu\varepsilon} \ln\left( \frac{2L_0^{j_0} R^2}{\varepsilon} \right) \right\} + 2N. \qquad (2.2.46)$$

**Доказательство.** По неравенству Фенхеля (см. главу 7 [181])

$$\left\langle \tilde{\nabla}^{m_{k+1}} f\left(y^{k+1}\right) - \nabla f\left(y^{k+1}\right), x^{k+1} - y^{k+1} \right\rangle - \frac{L_{k+1}^{j_{k+1}}/2}{2} \left\| x^{k+1} - y^{k+1} \right\|^2 \le \frac{2}{L_{k+1}^{j_{k+1}}} \left\| \tilde{\nabla}^{m_{k+1}} f\left(y^{k+1}\right) - \nabla f\left(y^{k+1}\right) \right\|_*^2 .$$

Поэтому ключевое неравенство (2.2.27)

$$f\left(y^{k+1}\right) + \left\langle \tilde{\nabla}^{m_{k+1}} f\left(y^{k+1}\right), x^{k+1} - y^{k+1} \right\rangle + \frac{L_{k+1}^{j_{k+1}}}{2} \left\| x^{k+1} - y^{k+1} \right\|^2 + \frac{\alpha_{k+1}}{A_{k+1}} \frac{\varepsilon}{2} \ge f\left(x^{k+1}\right)$$

переписывается следующим образом

$$f\left(y^{k+1}\right) + \left\langle \nabla f\left(y^{k+1}\right), x^{k+1} - y^{k+1} \right\rangle + L_{k+1}^{j_{k+1}} \left\| x^{k+1} - y^{k+1} \right\|^2 + \frac{\alpha_{k+1}}{A_{k+1}} \frac{\varepsilon}{2} + \frac{\alpha_{k+1}}{A_{k+1}} \varepsilon \ge f\left(x^{k+1}\right). \qquad (2.2.47)$$

При получении неравенства (2.2.47) для большей наглядности заменяем в рассуждениях правую часть в неравенстве Фенхеля оценкой его математического ожидания, равной согласно условию (2.2.43)

$$\frac{2D}{L_{k+1}^{j_{k+1}} m_{k+1}} = \frac{\alpha_{k+1}}{A_{k+1}} \varepsilon .$$

В действительности, тут требуются более громоздкие рассуждения, см., например, [39, 181] и arXiv:1705.09809, которые приведут к необходимости увеличения в несколько раз (во сколько именно раз, зависит от “тяжести” хвостов распределения $\nabla f(x,\xi)$ и от выбранного доверительного уровня) константы 2 в формуле выбора $m_{k+1}$ в описании СУМТ, и к соответствующему увеличению $N$ и $\mathrm{Q}$.

Оценим число обращений за стохастическим градиентом $\mathrm{Q}$, используя схему доказательства теоремы 2.2.9. Для этого, прежде всего, заметим, что $A_N \approx 2R^2/\varepsilon$. Рассмотрим два случая, когда $\mu \ge 0$ – мало: $\mu \ll \varepsilon/\left(2R^2\right)$, $\mu$ – велико: $\mu \gg \varepsilon/\left(2R^2\right)$.

В первом случае будем считать, что (см. формулу (2.2.35))

$$A_{k+1}/\alpha_{k+1}^2 \approx A_{k+1}\cdot(1+A_k\tilde{\mu})/\alpha_{k+1}^2 = L_{k+1}^{j_{k+1}}.$$

а во втором случае, что (см. формулу (2.2.36))

$$A_{k+1}^2\tilde{\mu}/\alpha_{k+1}^2 \approx A_{k+1}\cdot(1+A_k\tilde{\mu})/\alpha_{k+1}^2 = L_{k+1}^{j_{k+1}}.$$

В первом случае число обращений за стохастическим градиентом оценивается, соответственно, как

$$Q \approx \sum_{k=0}^{N} \frac{2DA_k}{L_k^{j_k}\alpha_k\varepsilon} = \frac{2D}{\varepsilon}\sum_{k=0}^{N}\frac{A_k}{L_k^{j_k}\alpha_k} \approx \frac{2D}{\varepsilon}\sum_{k=0}^{N}\alpha_k = \frac{2D}{\varepsilon}A_N \approx \frac{4DR^2}{\varepsilon^2}, \qquad (2.2.48)$$

$$Q \approx \sum_{k=0}^{N} \frac{2DA_k}{L_k^{j_k}\alpha_k\varepsilon} = \frac{2D}{\varepsilon}\sum_{k=0}^{N}\frac{A_k}{L_k^{j_k}\alpha_k} \approx \frac{2D}{\tilde{\mu}\varepsilon}\sum_{k=0}^{N}\frac{\alpha_k}{A_k} \approx \frac{2D}{\tilde{\mu}\varepsilon}\int_1^{A_N/\alpha_0}\frac{dA}{A} \approx \frac{2D}{\tilde{\mu}\varepsilon}\ln\left(\frac{A_N}{\alpha_0}\right) \approx \frac{2D}{\tilde{\mu}\varepsilon}\ln\left(\frac{2L_0^{j_0}R^2}{\varepsilon}\right). \qquad (2.2.49)$$

Из формул (2.2.48), (2.2.49) получаем оценку (2.2.46). ■

Оценка (2.2.46) с точностью до логарифмических множителей соответствует нижней оценки [91]. В случае $\|\ \| \neq \|\ \|_2$ в приведенных выше оценках может возникать дополнительный множитель $\sim \ln n$, $n = \dim x$ (Proposition 6, arXiv:1601.07592).

**Замечание 2.2.7.** Все приведенные в этом разделе результаты допускают обобщение на случай, когда вместо точного оракула, выдающего значения функции $f(x)$ и ее градиента $\nabla f(x)$ (или их несмещенные реализации), используется (стохастический) $(\delta, L, \mu)$-оракул с $\delta = O(\varepsilon/N)$ и $L = O\left(\max_{k=0,\dots,N} L_k^{j_k}\right)$ [39, 44, 48, 181, 183, 184] (см. также предыдущий раздел). Такое обобщение существенным образом используется, например, в замечании 2.2.6 выше (а также в разделе 1.6 главы 1).

## 2.3 Пример задачи композитной оптимизации (сильно выпуклый случай)

Рассмотрим конкретный пример задачи выпуклой композитной оптимизации [33, 138] (см. также раздел 3.2 главы 3 и раздел 5.1 главы 5)

$$F(x) = \frac{1}{2}\|Ax-b\|_2^2 + \mu\sum_{k=1}^{n} x_k \ln x_k \to \min_{\sum_{k=1}^{n} x_k = 1,\ x\ge 0} . \qquad (2.3.1)$$

Вместо ограничения $\sum_{k=1}^{n} x_k = 1$ можно рассматривать ограничение $\sum_{k=1}^{n} x_k \le 1$.

Разберем два случая а) $0 < \mu \ll \varepsilon/(2\ln n)$ – мало (сильную выпуклость композита в 1-норме можно не учитывать); б) $\mu \gg \varepsilon/(2\ln n)$ – достаточно большое (сильную выпуклость композита в 1-норме необходимо учитывать).

Выберем норму в прямом пространстве $\|\ \| = \|\ \|_1$. Положим

$$f(x) = \frac{1}{2}\|Ax-b\|_2^2,\ h(x) = \mu\sum_{k=1}^{n} x_k \ln x_k,\ Q = S_n(1) = \left\{x \ge 0 : \sum_{k=1}^{n} x_k = 1\right\},\ L = \max_{k=1,\dots,n}\left\|A^{\langle k\rangle}\right\|_2^2,$$

где $A^{\langle k\rangle}$ – $k$-й столбец матрицы $A$.

Введем два оператора

$$\mathrm{Grad}_{f,h}^{L}\left(x^{k+1}\right) = \arg\min_{x\in Q} \bar{F}\left(x; x^{k+1}\right),$$

$$\bar{F}\left(x; x^{k+1}\right) = f\left(x^{k+1}\right) + \left\langle \nabla f\left(x^{k+1}\right), x - x^{k+1}\right\rangle + \frac{L}{2}\left\|x - x^{k+1}\right\|_1^2 + h(x);$$

$$\mathrm{Mirr}_{f,h,z^k}^{\alpha}\left(\nabla f\left(x^{k+1}\right)\right) = \arg\min_{x\in Q}\left\{\left\langle \nabla f\left(x^{k+1}\right), x - z^k\right\rangle + \frac{1}{\alpha}V\left(x, z^k\right) + h(x)\right\},$$

где прокс-расстояние (расстояние Брэгмана) определяется формулой [135, 259, 271]

$$V(x,z) = d(x) - d(z) - \left\langle \nabla d(z), x - z\right\rangle,$$

прокс-функция $d(x) \ge 0$ считается сильно выпуклой относительно выбранной нормы $\|\ \| = \|\ \|_1$, с константой сильной выпуклости $\ge 1$. Для случая а) можно выбирать

$$d(x) = \ln n + \sum_{k=1}^{n} x_k \ln x_k .$$

Тогда

$$V(x,z) = \sum_{k=1}^{n} x_k \ln\left(x_k/z_k\right),\ R^2 \le \ln n .$$

Существуют варианты быстрого градиентного метода (например, Метод Треугольника из предыдущего раздела), в которых шаг $\mathrm{Grad}_{f,h}^{L}\left(x^{k+1}\right)$ заменяется его проксимальным аналогом, т.е., грубо говоря, $\left\|x-x^{k+1}\right\|_1^2$ заменяется в выражении для $\bar{F}\left(x;x^{k+1}\right)$ на $V\left(x,x^{k+1}\right)$. В таких вариантах метода (оценки скорости сходимости аналогичны оценкам, приводимым в утверждении 2.3.1) мы имеем ситуацию, когда композит совпадает по форме с прокс-расстоянием (энтропийного типа), и шаг итерации осуществим по явным формулам (см., например, [125, 259]). Таким образом, стоимость итерации будет $\mathrm{O}\left(nnz\left(A\right)\right)$, где $nnz\left(A\right)$ – число ненулевых элементов в матрице $A$ (считаем, что это число $\geq n$).

Для случая б) планируется использовать рестарт-технику (см. предыдущий раздел). Но для выбранной функции $V\left(x,z\right)$ (расстояние Кульбака–Лейблера [125]) процедура рестартов некорректна. Однако существует другой способ выбора прокс-функции (детали, см., например, в работах [91, 125, 228, 259])

$$d\left(x\right)=\frac{1}{2\left(a-1\right)}\|x\|_a^2,\ a=\frac{2\ln n}{2\ln n-1}. \tag{2.3.2}$$

В этом случае $R^2=\mathrm{O}\left(\ln n\right)$, $\omega_n=\mathrm{O}\left(\ln n\right)$.

Опишем используемый в этом разделе вариант быстрого градиентного метода (БГМ) Ю.Е. Нестерова в форме [135] (на самом деле от выбора конкретного варианта БГМ ключевое утверждение 2.3.2 не зависит). Здесь мы распространяем подход работы [135] на задачи композитной оптимизации. Фактически предложенный далее алгоритм есть сочетание БГМ работы [135] и конструкции композитной оптимизации работы с адаптивным подбором константы Липшица градиента [266].

Определим две числовые последовательности $\{\alpha_k,\tau_k\}$:

$$\alpha_1=\frac{1}{L_1},\ \alpha_{k+1}=\frac{1}{2L_{k+1}}+\sqrt{\frac{1}{4L_{k+1}^2}+\alpha_k^2\frac{L_k}{L_{k+1}}},\ \tau_k=\frac{1}{\alpha_{k+1}L_{k+1}}.$$

Заметим, что при $k\gg 1$

$$\alpha_k\sim\frac{k}{2L},\ \tau_k\sim\frac{2}{k}.$$

БГМ

1. $x^{k+1}=\tau_k z^k+\left(1-\tau_k\right)y^k$.

2. *Положим* $L_{k+1} := L_k/2$. *До тех пор пока*[53]

$$F\left(\mathrm{Grad}_{f,h}^{L_{k+1}}\left(x^{k+1}\right)\right) > \bar{F}_{L_{k+1}}\left(\mathrm{Grad}_{f,h}^{L_{k+1}}\left(x^{k+1}\right); x^{k+1}\right),$$

*выполнять*

$$L_{k+1} := 2L_{k+1}.$$

*Положить*

$$y^{k+1} = \mathrm{Grad}_{f,h}^{L_{k+1}}\left(x^{k+1}\right).$$

3. $z^{k+1} = \mathrm{Mirr}_{f,h,z^k}^{\alpha_{k+1}}\left(\nabla f\left(x^{k+1}\right)\right)$.

4. *Если не выполняется критерий останова (можно по-разному определять [44]), положить*

$$k := k+1$$

*и перейти к п.1. Иначе остановиться и выдать* $y^{k+1}$.

**Утверждение 2.3.1.** *Для задачи (2.3.1) БГМ генерирует такую последовательность точек* $\left\{x^k, y^k, z^k\right\}_{k=0}^N$, *что имеют место следующие неравенства (второе неравенство означает, что описанный вариант БГМ является прямо-двойственным методом)*

$$F\left(y^N\right) - F_* \le \frac{8LR^2}{\left(N+1\right)^2}.$$

$$\alpha_N^2 LF\left(y^N\right) \le \min_{x\in Q}\left\{\sum_{k=0}^{N-1}\alpha_{k+1}\left\{f\left(x^{k+1}\right) + \left\langle \nabla f\left(x^{k+1}\right), x - x^{k+1}\right\rangle + h\left(x\right)\right\} + V\left(x, x^0\right)\right\}.$$

Доказательство этого утверждения (полученное совместно с Антоном Родомановым и Павлом Двуреченским) доступно по ссылке

http://arxiv.org/pdf/1606.08988.pdf .

**Замечание 2.3.1 (см. также главу 3).** Заметим, что если условие $x \in Q$ можно записать, например, как $g\left(x\right) \le 0$, и ввести двойственную функцию

$$G\left(\lambda\right) = \min_x\left\{f\left(x\right) + h\left(x\right) + \left\langle \lambda, g\left(x\right)\right\rangle\right\},\ \lambda \ge 0,$$

то, поскольку

$$\alpha_N^2 L = \sum_{k=0}^{N-1}\alpha_{k+1} \overset{def}{=} S_N \ \text{ и } \ f\left(x^{k+1}\right) + \left\langle \nabla f\left(x^{k+1}\right), x - x^{k+1}\right\rangle \le f\left(x\right),$$

---

[53] Можно так переписать это условие (подобно [266]), чтобы не требовалось рассчитывать значения функции. Для варианта БГМ из предыдущего раздела (Метод Треугольника) нам неизвестна такая перезапись.

получим

$$F\left(y^N\right)-G\left(\bar{\lambda}^N\right)\le\frac{8L\tilde{R}^2}{\left(N+1\right)^2},$$

где $\bar{\lambda}^N=\lambda^N/S_N$, $\lambda^N$ – множитель Лагранжа к ограничению $g\left(x\right)\le 0$ в задаче

$$\sum_{k=0}^{N-1}\alpha_{k+1}\left\{f\left(x^{k+1}\right)+\left\langle\nabla f\left(x^{k+1}\right),x-x^{k+1}\right\rangle+h\left(x\right)\right\}+V\left(x,x^0\right)\to\min_{g(x)\le 0},$$

а $\tilde{R}^2=V\left(\tilde{x},x^0\right)$, где $\tilde{x}$ – решение задачи

$$f\left(x\right)+h\left(x\right)+\left\langle\bar{\lambda}^N,g\left(x\right)\right\rangle\to\min_x.$$

Последнее условие, к сожалению, в общем случае не дает возможности как-то разумно оценивать сверху $\tilde{R}^2$, как следствие, возникает проблема с теоретической оценкой зазора двойственности. Проблема решается, если удается обосновать возможность компактификации. Пример того, как эту компактификацию можно делать (на основе “Слейтеровских соображений”) будет описан далее (см. формулы (2.3.6), (2.3.7)).

Используя описанную в предыдущем разделе технику рестартов можно получить из утверждения 1 его аналог в случае сильно выпуклой постановки задачи (2.3.1) – случай б). Мы опускаем соответствующие рассуждения, и остановимся подробнее на том, как осуществлять шаг итерации описанного БГМ в случае б), т.е. когда прокс-функция выбирается согласно (2.3.2). Сложность выполнения одной итерации (дополнительная к вычислению градиента гладкой части функционала $\mathrm{O}\left(nnz\left(A\right)\right)$) определяется тем, насколько эффективно можно решить задачу следующего вида

$$\tilde{F}\left(x\right)=\left\langle c,x\right\rangle+\left\|x\right\|_a^2+\bar{\mu}\sum_{k=1}^n x_k\ln x_k\to\min_{x\in S_n(1)}. \tag{2.3.3}$$

Задачу (3) удобно переписать следующим почти “сепарабельным” образом

$$\left\langle c,x\right\rangle+t+\bar{\mu}\sum_{k=1}^n x_k\ln x_k\to\min_{\substack{x\in S_n(1),\,\|x\|_a^a\le t^{a/2},\\ 0\le t\le n^{2/a},\,0\le x_k\le 1,\,k=1,\dots,n}}$$

Слово “почти” можно убрать, если с помощью метода множителей Лагранжа переписать задачу следующим образом

$$\tilde{G}\left(\lambda\right)=$$

$$=\min_{\substack{0\le t\le n^{2/a},\\ 0\le x_k\le 1,\,k=1,\dots,n}}\left\{\sum_{k=1}^n c_k x_k+t+\lambda_1\cdot\left(\sum_{k=1}^n x_k-1\right)+\lambda_2\cdot\left(\sum_{k=1}^n x_k^a-t^{a/2}\right)+\bar{\mu}\sum_{k=1}^n x_k\ln x_k\right\}\to\max_{\lambda_1\in\mathbb{R},\,\lambda_2\ge 0}. \tag{2.3.4}$$

Поиск минимума $(x(\lambda), t(\lambda))$, где

$$t(\lambda) = \min\left\{ \left(\frac{\lambda_2 a}{2}\right)^{\frac{2}{2-a}}, n^{\frac{2}{a}} \right\},$$

сводится к решению $n$ одномерных задач сильно выпуклой оптимизации на отрезке $[0,1]$. Таким образом, если задаться некоторой точностью $\sigma > 0$, то за время $\mathrm{O}\big(n \ln(n/\sigma)\big)$ методом деления отрезка пополам (или, скажем, методом золотого сечения [22]) можно найти такой $\tilde{x}^{\sigma}(\lambda)$, что

$$\left\| \tilde{x}^{\sigma}(\lambda) - x(\lambda) \right\|_1 = \mathrm{O}(\sigma). \tag{2.3.5}$$

Далее попробуем (следуя [152]) оценить "запас" в условии Слейтера, чтобы, исходя из этого, оценить сверху размер решения $\lambda = \lambda_*$ двойственной задачи (2.3.4) (в приводимой далее выкладке, приводящей к формуле (2.3.6), для упрощения записи мы опускаем нижний индекс "*" у $\lambda$). Из сильной двойственности имеем

$$-\|c\|_{\infty} - \bar{\mu} \ln n \le \tilde{F}_* = \tilde{G}_* \le \sum_{k=1}^{n} c_k \bar{x}_k + \bar{t} + \lambda_1 \cdot \left( \sum_{k=1}^{n} \bar{x}_k - 1 \right) + \lambda_2 \cdot \left( \sum_{k=1}^{n} \bar{x}_k^{a} - \bar{t}^{a/2} \right) + \bar{\mu} \sum_{k=1}^{n} \bar{x}_k \ln \bar{x}_k .$$

Если $\lambda_1 \ge 0$, то положим $\bar{t} = 1$, $\bar{x}_k = 1/(2n)$, $k = 1,...,n$. Тогда

$$\frac{1}{2}\lambda_1 + \frac{1}{2}\lambda_2 \le 2\|c\|_{\infty} + 2\bar{\mu}\ln(n) + 1,$$

Если $\lambda_1 < 0$, то положим $\bar{t} = 8$, $\bar{x}_k = 2/n$, $k = 1,...,n$. Тогда

$$|\lambda_1| + \lambda_2 \le 3\|c\|_{\infty} + 2\bar{\mu}\ln(2n) + 8.$$

В любом случае, с хорошим запасом можно гарантировать, что

$$\|\lambda_*\|_1 \le 4\|c\|_{\infty} + 4\bar{\mu}\ln(2n) + 8 \overset{def}{=} C. \tag{2.3.6}$$

Таким образом, чтобы решить задачу (2.3.3), мы должны решить двойственную задачу (2.3.4), которую (в виду формулы (2.2.6)) можно переписать следующим образом

$$\breve{G}(\lambda) = -\tilde{G}(\lambda) \to \min_{\substack{\lambda_1 \in \mathbb{R},\, \lambda_2 \ge 0 \\ \|\lambda\|_1 \le C}} . \tag{2.3.7}$$

Поскольку эта задача оптимизации на двумерной плоскости (т.е. в пространстве малой размерности), то ее можно решать, скажем, методом эллипсоидов [91]. При этом для расчета градиента $\breve{G}(\lambda)$ мы должны решить задачу (2.3.3) и воспользоваться формулой Демьянова–Данскина [152]

$$\frac{\partial \breve{G}}{\partial \lambda_1} = 1 - \sum_{k=1}^{n} x_k(\lambda), \ \frac{\partial \breve{G}}{\partial \lambda_2} = t(\lambda)^{a/2} - \sum_{k=1}^{n} x_k(\lambda)^a .$$

К сожалению, точно решить задачу (2.3.4) мы не можем, зато можем найти приближенное значение градиента. Точнее говоря, в виду (2.3.5), (2.3.6), мы можем найти для задачи (2.3.7) $\delta = O(C\sigma)$-градиент $\nabla_\delta G(\lambda)$ (см., например, [99]). Если использовать в методе эллипсоидов в пространстве размерности $r$ (в нашем случае $r = 2$) вместо градиента $\delta$-градиент (чаще говорят $\delta$-субградиент, но в нашем случае можно говорить о градиенте), то имеют место следующие оценки [91]

$$\breve{G}(\lambda^N) - \breve{G}_* \le \varepsilon, \ N = O(r^2 \ln(C/\varepsilon)), \ \delta \le O(\varepsilon), \tag{2.3.8}$$

При этом стоимость одной итерации будет $O(r^2)$. Число итераций можно сократить в $\sim r$ раз, сохранив сложность итерации (см., например, [152]). В нашем случае стоимость одной итерации будет $O(n\ln(nC/\varepsilon))$.

Однако решение задачи (2.3.4) (или (2.3.7)), в смысле (2.3.8) еще не гарантирует возможность точного восстановления решения задачи (2.3.3). Для того чтобы показать, что метод эллипсоидов с той же по порядку точностью $\varepsilon$ позволяет восстанавливать (без каких бы то ни было существенных дополнительных затрат) решение задачи (2.3.3) нужно воспользоваться прямо-двойственностью этого метода [261] (см. также главу 3). В виду компактности множества (единичный симплекс), на котором ведется оптимизация в прямом пространстве и сильной выпуклости функционала прямой задачи (2.3.3) мы не просто восстанавливаем из прямо-двойственной процедуры метода эллипсоидов решение задачи (2.3.3) с точностью по функционалу (прямой задачи) порядка $\varepsilon$, но и делаем это в нужном нам более сильном смысле – см. п. 5.5.1 (следует сравнить с п. 4.6 [181]). Формула 5.5.15 [259] гарантирует при этом справедливость следующего результата.

**Утверждение 2.3.2.** *Для задачи (2.3.1) в случае б) БГМ с рестартами и с прокс-функцией (2.3.2) приводит к необходимости на каждой итерации наряду с расчетом градиента гладкой части функционала (*$O(nnz(A))$ *операций) два раза решать задачу типа (2.3.3) с помощью перехода к двойственной задаче и ее решения с помощью прямо-двойственной версии метода эллипсоидов (*$O(n\ln(C/\varepsilon)\ln(nC/\varepsilon))$ *операций). При этом*

$$F(y^N) - F_* \le \varepsilon,$$

*если общее число итераций (обращений к оракулу за градиентом)*

$$N = \mathrm{O}\left(\sqrt{\frac{L}{\mu}\omega_n}\left\lceil \ln\left(\frac{\mu}{\varepsilon}\right)\right\rceil\right) = \mathrm{O}\left(\sqrt{\frac{\max\limits_{k=1,\ldots,n}\left\|A^{\langle k\rangle}\right\|_2^2 \ln n}{\mu}}\left\lceil \ln\left(\frac{\mu}{\varepsilon}\right)\right\rceil\right). \qquad (2.3.9)$$

Заметим, что в “пороговой” ситуации, отвечающей регуляризации (см. предыдущий раздел), $\mu \simeq \varepsilon/(2\ln n)$. В этом случае формула (2.3.9) примет вид

$$N = \mathrm{O}\left(\sqrt{\frac{\max\limits_{k=1,\ldots,n}\left\|A^{\langle k\rangle}\right\|_2^2 \ln^2 n}{\varepsilon}}\right),$$

что с точностью до $\sim\sqrt{\ln n}$ соответствует оценке в случае а). Отличие случая а) и б) также и в том, что в случае а) существует способ добиться стоимости итерации $\mathrm{O}\big(nnz(A)\big)$, а в случае б) нам не известно более эффективного способа, чем способ (описанный выше) со стоимостью итерации

$$\mathrm{O}\big(nnz(A) + n\ln(C/\varepsilon)\ln(nC/\varepsilon)\big).$$

Поскольку в типичных приложениях первое слагаемое заметно доминирует второе, то можно не сильно задумываться о плате за невозможность выполнения “проектирования” по явным формулам, а также не сильно задумываться с какой точностью нужно решать вспомогательную задачу, делая это с точностью длины мантиссы (описанный подход позволяет решать ее с очень хорошей точностью, и сложность решения вспомогательной задачи практически не чувствительна к этой точности). По-видимому, этот тезис имеет достаточно широкий спектр практических приложений. Настоящий раздел имел одной из своих целей на конкретном примере более подробно, чем это принято на практике, продемонстрировать, что для большого класса задач наличие явных формул для шага итерации не есть сколько-нибудь сдерживающие обстоятельство для использования метода. Используемая при этом техника и способ рассуждений характерным образом (на наш взгляд) демонстрируют современный арсенал средств решения задач выпуклой оптимизации в пространствах больших размеров.

# Глава 3 Прямо-двойственные градиентные методы для задач выпуклой оптимизации

## 3.1 О численных методах решения задач энтропийно-линейного программирования с помощью решения регуляризованной двойственной задачи

### 3.1.1 Введение

Данный раздел посвящен разработке эффективных численных методов решения задач энтропийно-линейного программирования (ЭЛП) [198] и получению оценок их скоростей сходимости. Имеется огромное число приложений, в которых возникают задачи ЭЛП – см. [26, 30, 49, 106, 207, 221, 230], а также работы на ежегодных International workshops on Bayesian inference and maximum entropy methods in science and engineering (AIP Conf. Proceedings). Например, одно из таких приложений нам уже встречалось ранее в связи с расчетом матрицы корреспонденций (разделы 1.1, 1.2 главы 1).

Основная идея (предложенная в конце декабря 2012 года Ю.Е. Нестеровым) предлагаемого в разделе подхода заключается в решении специальным образом регуляризованной (по Тихонову) двойственной задачи к задаче ЭЛП. Зная решение двойственной задачи, по явным формулам можно найти решение исходной задачи. В разделе приводятся точные формулы, определяющие сколько достаточно сделать итераций быстрого градиентного метода для регуляризованной двойственной задачи, чтобы гарантировано (с заданной точностью) восстановить решение исходной задачи. Насколько нам известно, для задач ЭЛП такого рода формулы выписываются впервые. Собственно, по известным нам остальным методам решения задач ЭЛП в основном имеются только теоремы о сходимости [55, 106, 198] (без оценок скорости сходимости).

Отметим также, что полученные в данном разделе оценки при определенных условиях лучше имеющихся нижних оценок. Поясним сказанное.

Известно, что задача поиска такого $x^* \in \mathbb{R}^n$, что $Ax^* = b$ сводится к задаче выпуклой гладкой оптимизации

$$f(x) = \|Ax - b\|_2^2 \to \min_x .$$

Нижние оценки (при $k \le n$) для такой задачи приводят к оценке:

$$f(x_k) \ge \Omega\left(L_x R_x^2 / k^2\right),$$

т.е. существует такое $\chi$, что при $1 \le k \le n$ имеет место неравенство $f(x_k) \ge \chi L_x R_x^2 / k^2$. Здесь $k$ – число умножений матрицы $A$ на столбец/строку, $x_k$ – то, что можно получить на основе результатов этих умножений для самого лучшего метода на классе всевозможных задач такого типа; $L_x = \sigma_{\max}(A) = \lambda_{\max}\left(A^T A\right)$, $R_x = \left\|x^*\right\|_2 = \left\|A^T\left(AA^T\right)^{-1} b\right\|_2$. Откуда следует, что только при $k \geq \Omega\left(\sqrt{L_x} R_x / \varepsilon\right)$ можно гарантировать выполнение неравенства $f(x_k) \leq \varepsilon^2$, т.е. $\left\|Ax_k - b\right\|_2 \leq \varepsilon$. Детали см. в [стр. 264–274, 91]. Далее в разделе будет показано (см., например, теорему 3.1.2), что при достаточно общих условиях можно улучшить приведенную оценку, попутно добившись, чтобы полученный $x_k$ был бы приближенным решением задачи ЭЛП, т.е. не только удовлетворял бы “возмущенным” аффинным ограничениям $\left\|Ax_k - b\right\|_2 \leq \varepsilon$ в задаче ЭЛП.

Структура раздела следующая. В подразделе 3.1.2 приведен простейший пример возникновения задачи ЭЛП (более подробно о важности задач ЭЛП при поиске равновесий макросистем написано в приложении в конце диссертации). В подразделе 3.1.3 рассматривается постановка задачи ЭЛП, приводятся основные определения. В подразделе 3.1.4 приводятся необходимые вспомогательные результаты. В подразделе 3.1.5 формулируются основные результаты данного раздела. В подразделе 3.1.6 эти результаты обсуждаются. В подразделе 3.1.7 делаются заключительные замечания.

### 3.1.2 Парадокс Эренфестов

Следуя учебнику [75] опишем известный в физике парадокс Эренфестов [73] в немного вольной трактовке.

Рядом стоят две собаки с номерами ***1*** и ***2***. На собаках как-то расположились $M \gg 1$ блох. Каждая блоха в промежутке времени $[t, t+h)$ с вероятностью $\lambda h + o(h)$ ($\lambda = 1$) независимо от остальных перескакивает на соседнюю собаку. Пусть в начальный момент все блохи собрались на собаке с номером ***1***. Тогда для всех $t \geq \chi M$ ($\chi \sim 10$)

$$P\left(\frac{\left|n_1(t) - n_2(t)\right|}{M} \leq \frac{5}{\sqrt{M}}\right) \geq 0.99,$$

где $n_1(t)$ – число блох на первой собаке в момент времени $t$, а $n_2(t)$ – на второй (случайные величины). То есть относительная разность числа блох на собаках будет иметь порядок малости $\mathrm{O}\left(1/\sqrt{M}\right)$ на больших временах ($t \geq \chi M$). Далее мы поясним связь этого результата с принципом максимума энтропии, который мы будем записывать в виде принципа минимума минус энтропии.

Описанная выше марковская динамика имеет закон сохранения числа блох: $n_1(t)+n_2(t)\equiv M$ , и это будет единственный закон сохранения. Стационарная мера имеет вид (теорема Санова):

$$\nu(n_1,n_2)=\nu(c_1M,c_2M)=M!\frac{(1/2)^{n_1}}{n_1!}\frac{(1/2)^{n_2}}{n_2!}=C_M^{n_1}2^{-M}\simeq\frac{2^{-M}}{\sqrt{2\pi c_1c_2}}\exp\left(-M\cdot H(c_1,c_2)\right),$$

где $H(c_1,c_2)=\sum_{i=1}^{2}c_i\ln c_i$, $c_i$ – концентрация блох на собаки $i$ в равновесии (т.е. при $t\to\infty$). Кстати, сказать, из такого вида стационарной меры, можно получить, что если в начальный момент все блохи находились на одной собаке, то математическое ожидание времени первого возвращения макросистемы в такое состояние будет порядка $2^M$ . Равновесие данной макросистемы естественно определять как состояние, в малой окрестности которого сконцентрирована стационарная мера (принцип максимума энтропии Больцмана–Джейнса)

$$c^*=\begin{pmatrix}1/2\\1/2\end{pmatrix}=\arg\min_{\substack{c_1+c_2=1\\c\geq 0}}H(c).$$

Поиск равновесия привел к необходимости решения задачи ЭЛП. Этот же результат можно было получить и при другом порядке предельных переходов (обратном к рассмотренному выше порядку: $t\to\infty,\ M\to\infty$). А именно, сначала считаем, что при $t=0$ существует предел $c_i(t)\overset{\text{п.н.}}{=}\lim_{M\to\infty}n_i(t)/M$ . Тогда (теорема Т. Куртца) этот предел существует при любом $t>0$, причем $c_1(t),\ c_2(t)$ – детерминированные (не случайные) функции, удовлетворяющие СОДУ

$$\frac{dc_1}{dt}=\lambda\cdot(c_2-c_1),$$

$$\frac{dc_2}{dt}=\lambda\cdot(c_1-c_2).$$

Глобально устойчивым положением равновесия этой СОДУ будет $c^*$ , а $H(c)$ – функция Ляпунова этой СОДУ (убывает на траекториях СОДУ, и имеет минимум в точке $c^*$ ). Все это можно перенести на общие модели стохастической химической кинетики с (обобщенным) условием детального баланса [30, 24, 32, 86].

Приведенный парадокс Эренфестов является, пожалуй, простейшим примером того как возникают задачи ЭЛП при поиске равновесий макросистем. Подробнее об этом можно прочитать, например, в [30] и в приложении в конце диссертации. В частности, в большинстве реальных приложений возникают задачи ЭЛП, в которых не “две собаки” и не “один за-кон сохранения числа блох”, а на много больше (см., например, [36]), и возникающая задача ЭЛП не решается по явным формулам. Требуется разработка эффективных численных методов.

### 3.1.3 Задача энтропийно-линейного программирования[1]

Рассматривается задача ЭЛП [198]

$$f(x)=\sum_{i=1}^{n} x_i \ln(x_i/\xi_i) \to \min_{x\in S_n(1);\, Ax=b}, \tag{3.1.1}$$

где $S_n(1)=\left\{x\in\mathbb{R}^n : x_i\geq 0,\ i=1,...,n,\ \sum_{i=1}^{n} x_i=1\right\}$ – единичный симплекс[2] в $\mathbb{R}^n$, $\xi\in \operatorname{ri} S_n(\Xi)$.

Мы считаем матрицу $A$ разреженной с не более чем $s \ll n$ элементами в каждой строке. Число строк в матрице $A$ есть $m \ll n$. Построим двойственную задачу

$$\min_{x\in S_n(1);\, Ax=b} \sum_{i=1}^{n} x_i \ln(x_i/\xi_i) = \min_{x\in S_n(1)} \max_{\lambda\in\mathbb{R}^m} \left\{ \sum_{i=1}^{n} x_i \ln(x_i/\xi_i) + \langle \lambda, b-Ax\rangle \right\} =$$

$$= \max_{\lambda\in\mathbb{R}^m} \min_{x\in S_n(1)} \left\{ \sum_{i=1}^{n} x_i \ln(x_i/\xi_i) + \langle \lambda, b-Ax\rangle \right\} = \max_{\lambda\in\mathbb{R}^m} \left\{ \langle \lambda, b\rangle - \ln\left( \sum_{i=1}^{n} \xi_i \exp\left(\left[A^T\lambda\right]_i\right)\right)\right\}.$$

Таким образом, двойственная задача имеет вид

$$\varphi(\lambda)=\langle \lambda, b\rangle - \ln\left( \sum_{i=1}^{n} \xi_i \exp\left(\left[A^T\lambda\right]_i\right)\right) \to \max_{\lambda\in\mathbb{R}^m}. \tag{3.1.2}$$

Решения прямой и двойственной задачи связаны:

$$x_i(\lambda)=\frac{\xi_i \exp\left(\left[A^T\lambda\right]_i\right)}{\sum_{k=1}^{n} \xi_k \exp\left(\left[A^T\lambda\right]_k\right)}.$$

**Определение 3.1.1.** *Под $(\varepsilon_f, \varepsilon)$-решением задачи (3.1.1) будем понимать такой вектор $x$, что*

$$f(x)-f^* \leq \varepsilon_f,\quad \|Ax-b\|_2 \leq \varepsilon,$$

*где* $f^* = \min_{x\in S_n(1);\, Ax=b} f(x)$.

---

[1] Результаты этого и следующего подраздела были установлены совместно с Ю.Е. Нестеровым. Теорема 3.1.1 в подразделе 3.1.5 была установлена Ю.Е. Нестеровым.

[2] Более общая ситуация, когда $x\in S_n(\Lambda)$, может быть сведенная к $x\in S_n(1)$ с помощью замены переменных. При этом $\tilde{x}:=x/\Lambda$ $\tilde{\xi}:=\xi/\Lambda$, $\tilde{b}:=b/\Lambda$, $\tilde{A}:=A$, $\tilde{f}:=f\Lambda$.

**Лемма 3.1.1.** *Пусть* $-\langle \lambda, \nabla\varphi(\lambda)\rangle \le \varepsilon_f$, $\|\nabla\varphi(\lambda)\|_2 \le \varepsilon$, *тогда* $x(\lambda)$ – $(\varepsilon_f, \varepsilon)$*-решением задачи (3.1.1).*

**Доказательство.** Соотношение $\nabla\varphi(\lambda) = b - Ax(\lambda)$ следует из представления двойственной задачи $\varphi(\lambda) = \min_{x \in S_n(1)} \{ f(x) + \langle \lambda, b - Ax\rangle \}$ и теоремы Демьянова–Рубинова–Данскина.

Обозначим через $x^*$ решение задач (3.1.1). Тогда

$$f(x(\lambda)) + \langle \lambda, b - Ax(\lambda)\rangle \le f(x^*) + \langle \lambda, b - Ax^*\rangle = f(x^*) = f^*.$$

Откуда

$$f(x(\lambda)) = f(x^*) - \langle \lambda, b - Ax(\lambda)\rangle = f^* - \langle \lambda, \nabla\varphi(\lambda)\rangle \le f^* + \varepsilon_f. \bullet$$

**Лемма 3.1.2 (см. [271]).** *Имеет место следующее неравенство*

$$\|\nabla\varphi(\lambda_2) - \nabla\varphi(\lambda_1)\|_2 \le L\|\lambda_2 - \lambda_1\|_2,$$

*где* $L = \max_{1\le i\le n} \left\| [A]^{\langle i\rangle} \right\|_2^2$, $[A]^{\langle i\rangle}$ – $i$*-й столбец матрицы* $A$.

К сожалению, несмотря на установленную в лемме 3.1.2 гладкость функционала (3.1.2), мы имеем лишь довольно грубые оценки его сильной вогнутости. Быстрый градиентный метод (БГМ) [93] (стартующий с $\lambda = 0$) для задачи (3.1.2) после $O\left(\sqrt{LR^2/\tilde{\varepsilon}}\right)$ итераций, где $R = \|\lambda^*\|_2$ – размер решения задачи (3.1.2), $\tilde{\varepsilon} = \min\left\{\varepsilon^2/(2L), \varepsilon_f^2/(2LR^2)\right\}$ – точность решения задачи (3.1.2), гарантирует в виду леммы 3.1.1 и следующих неравенств ($\varphi^*$ – максимальное значение функционала в двойственной задаче)

$$\frac{1}{2L}\|\nabla\varphi(\lambda)\|_2^2 \le \varphi^* - \varphi(\lambda), \quad \frac{\left|-\langle \lambda, \nabla\varphi(\lambda)\rangle\right|^2}{2LR^2} \le \frac{1}{2L}\|\nabla\varphi(\lambda)\|_2^2 \le \varphi^* - \varphi(\lambda),$$

что $x(\lambda)$ – $(\varepsilon_f, \varepsilon)$-решение задачи (3.1.1). Итоговая оценка на число итераций примет вид

$$O\left(\max\left\{LR/\varepsilon, LR^2/\varepsilon_f\right\}\right). \qquad (3.1.3)$$

Эту оценку можно улучшить. Этому и будет посвящена последующая часть данного раздела.

### 3.1.4 Вспомогательные результаты для регуляризованной двойственной задачи к задаче энтропийно-линейного программирования

Регуляризуем функционал (3.1.2) по А.Н. Тихонову (см. также раздел 2.2 главы 2):

$$\varphi_\delta(\lambda) = \varphi(\lambda) - \frac{\delta}{2}\|\lambda\|_2^2,$$

и вместо задачи (3.1.2) будем решать регуляризованную задачу

$$\varphi_\delta(\lambda) \to \max_{\lambda \in \mathbb{R}^m}.$$

Параметр $\delta$ будет оптимально подобран позже.

Следующие три леммы решают задачу восстановления $(\varepsilon_f, \varepsilon)$-решением задачи (3.1.1) по решению этой регуляризованной задачи. В отличие от нерегуляризованного случая здесь возникают некоторые технические места, требующие аккуратной проработки.

**Лемма 3.1.3.** *Имеют место следующие неравенства*

$$\|\nabla\varphi(\lambda)\|_2 \le \|\nabla\varphi_\delta(\lambda)\|_2 + \delta\|\lambda\|_2,$$

$$-\langle\lambda, \nabla\varphi(\lambda)\rangle \le \frac{L_\delta}{\delta}\left(\varphi_\delta^* - \varphi_\delta(\lambda)\right), L_\delta = L + \delta,$$

$$\frac{\delta}{2}\|\lambda_\delta^*\|_2^2 \le \ln\left(\Xi \big/ \min_{i=1,\dots,n} \xi_i\right) \overset{def}{=} \Delta_\varphi,$$

*где $\lambda_\delta^*$ – решение регуляризованной двойственной задачи, а $\varphi_\delta^*$ – значение функционала в регуляризованной двойственной задаче на этом решении.*

**Доказательство.** Не очень тривиально лишь второе и третье неравенство. Второе неравенство следует из следующей выкладки

$$\varphi_\delta^* - \varphi_\delta(\lambda) \ge \frac{\|\nabla\varphi_\delta(\lambda)\|_2^2}{2L_\delta} = \frac{\|\nabla\varphi(\lambda) - \delta\lambda\|_2^2}{2L_\delta} \ge -\frac{\delta\langle\nabla\varphi(\lambda), \lambda\rangle}{L_\delta}.$$

Третье неравенство следует из следующей выкладки

$$\frac{\delta}{2}\|\lambda_\delta^*\|_2^2 \le \varphi_\delta^* - \varphi_\delta(0) \le \varphi^* - \varphi(0) =$$

$$= \min_{x \in S_n(1);\, Ax = b} \sum_{i=1}^n x_i \ln(x_i/\xi_i) - \min_{x \in S_n(1)} \sum_{i=1}^n x_i \ln(x_i/\xi_i) \le \ln\left(\Xi \big/ \min_{i=1,\dots,n} \xi_i\right).$$

Первое неравенство в этой цепочке следует из того, что $\nabla\varphi_\delta(\lambda_\delta^*) = 0$ и $\varphi_\delta(\lambda)$ – сильно вогнутая функция, с константой сильной вогнутости $\ge \delta$. ●

Для решения регуляризованной задачи воспользуемся БГМ для сильно вогнутых задач [93] ($\lambda_0 = u_0 = 0$):

$$\begin{cases} \lambda_{k+1} = u_k + \dfrac{1}{L_\delta}\left(b - Ax(\lambda_k) - \delta\lambda_k\right), \\ u_{k+1} = \lambda_k + \dfrac{\sqrt{L_\delta} - \sqrt{\delta}}{\sqrt{L_\delta} + \sqrt{\delta}}(\lambda_{k+1} - \lambda_k). \end{cases} \tag{3.1.4}$$

**Лемма 3.1.4 (см. [93]).** *Имеют место следующие неравенства* ($k = 0,1,\dots$)

$$\varphi_\delta^* - \varphi_\delta(\lambda_k) \le 2\Delta_\varphi \exp\left(-k\sqrt{\frac{\delta}{L_\delta}}\right),$$

$$\left\|\nabla\varphi_\delta(\lambda_k)\right\|_2^2 \le 4L_\delta\Delta_\varphi \exp\left(-k\sqrt{\frac{\delta}{L_\delta}}\right),$$

$$\left\|\lambda_k - \lambda_\delta^*\right\|_2^2 \le \min\left\{\left\|\lambda_\delta^*\right\|_2^2, \frac{4\Delta_\varphi}{\delta}\exp\left(-k\sqrt{\frac{\delta}{L_\delta}}\right)\right\}.$$

Из лемм 3.1.3, 3.1.4 следует

**Лемма 3.1.5.** *Имеет место следующее неравенство*

$$\left\|\lambda_k\right\|_2 \le \left\|\lambda_k - \lambda_\delta^*\right\|_2 + \left\|\lambda_\delta^*\right\|_2 \le 2\left\|\lambda_\delta^*\right\|_2 \le \sqrt{8\Delta_\varphi/\delta}.$$

### 3.1.5 Основные результаты

Следствием лемм 3.1.1, 3.1.3 – 3.1.5 является следующая

**Теорема 1.** *Если выбрать* $\sqrt{\delta} \simeq \varepsilon\big/\sqrt{9\Delta_\varphi}$, *то после*

$$N \simeq \frac{\sqrt{9L\Delta_\varphi}}{\varepsilon}\ln\left(\frac{9L\Delta_\varphi}{\varepsilon_f\varepsilon^2}\right) \tag{3.1.5}$$

*итераций метода (3.1.4) получим такой* $\lambda_N$, *что* $x(\lambda_N)$ – $(\varepsilon_f, \varepsilon)$*-решением задачи (1).*

**Доказательство.** Из лемм 3.1.1, 3.1.3 – 3.1.5 имеем

$$\exp\left(-N\sqrt{\frac{\delta}{L_\delta}}\right) = \min\left\{\frac{\delta\varepsilon_f}{2L_\delta\Delta_\varphi}, \frac{\left(\varepsilon - \sqrt{8\Delta_\varphi\delta}\right)_+^2}{4L_\delta\Delta_\varphi}\right\}.$$

Откуда $\sqrt{2\varepsilon_f}\sqrt{\delta} = \varepsilon - \sqrt{8\Delta_\varphi}\sqrt{\delta}$, следовательно $\sqrt{\delta} = \varepsilon\big/\left(\sqrt{8\Delta_\varphi} + \sqrt{2\varepsilon_f}\right)$. Для упрощения формул будем выбирать $\delta$ немного не оптимально: $\sqrt{\delta} \simeq \varepsilon\big/\sqrt{9\Delta_\varphi}$. Тогда

$$N \simeq \frac{\sqrt{9L\Delta_\varphi}}{\varepsilon}\ln\left(\frac{18L\Delta_\varphi}{\varepsilon_f\varepsilon^2}\right). \bullet$$

В отличие от полученной ранее оценки (3.1.3), в оценку (3.1.5) не входит потенциально большой размер $R = \left\| \lambda^* \right\|_2$. О том, что в приложениях этот размер, действительно, может быть большим, говорят результаты численных экспериментов [56]. Оценка (3.1.5) доминирует оценку (3.1.3) даже для небольших значений $R$. Естественным образом возникает вопрос: можно ли улучшить оценку (3.1.5) если известно, что $R$ – не очень большое число? Разберем эту ситуацию.

Обозначим через $R_\delta = \max_{k=0,1,\dots} \left\| \lambda_k \right\|_2$. Лемма 3.1.5 дает оценку на $R_\delta$, но эта оценка может быть сильно завышена, поэтому завышенной может получиться и оценка (3.1.5). Представим, что нам известно значение $R$, которое мажорирует $R_\delta$ (равномерно по $\delta \ge 0$). Выберем оптимально $\delta$ исходя из известного значения $R$. Для этого в виду лемм 3.1.3, 3.1.5 нужно решить задачу (см. также доказательство теоремы 3.1.1 с $R = \sqrt{8\Delta_\varphi / \delta}$)

$$\max\left\{2\varepsilon_f \delta, \left(\varepsilon - R\delta\right)_+^2\right\} \to \min_{\delta \ge 0}.$$

Эту задачу можно решить явно

$$\sqrt{\delta} = \frac{\varepsilon}{\sqrt{\varepsilon_f/2} + \sqrt{\varepsilon_f/2 + R\varepsilon}} .$$

Однако нам будет удобнее взять

$$\sqrt{\delta} \simeq \frac{\varepsilon}{2\sqrt{\varepsilon_f/2 + R\varepsilon}}. \tag{3.1.6}$$

Такая замена не сильно скажется на оценках необходимого числа итераций, но несколько упростит формулы.

**Теорема 3.1.2.** *Если выбрать $\delta$ согласно формуле (3.1.6), то после*

$$N \simeq \sqrt{\frac{2L\cdot\left(\varepsilon_f + 2R\varepsilon\right)}{\varepsilon^2}} \ln\left(\frac{4L\Delta_\varphi \cdot \left(\varepsilon_f + 2R\varepsilon\right)}{\varepsilon_f \varepsilon^2}\right) \tag{3.1.7}$$

*итераций метода (3.1.4) получим такой $\lambda_N$, что $x(\lambda_N)$ – $(\varepsilon_f, \varepsilon)$-решением задачи (3.1.1).*

*Основной вклад в оценку общей трудоемкости вносит расчет градиента в (3.1.4)*

$$\nabla \varphi_\delta(\lambda) = b - Ax(\lambda) - \delta\lambda,$$

*требующий $O(n + sm)$ арифметических операций. Общие затраты будут*

$$O\left((n + sm)\sqrt{\frac{2L\cdot\left(\varepsilon_f + 2R\varepsilon\right)}{\varepsilon^2}} \ln\left(\frac{4L\Delta_\varphi \cdot \left(\varepsilon_f + 2R\varepsilon\right)}{\varepsilon_f \varepsilon^2}\right)\right)$$

Доказательство теоремы 3.1.2 аналогично доказательству теоремы 3.1.1.

### 3.1.6 Обсуждение основных результатов

Оценка (3.1.7) доминирует оценку (3.1.5) для не очень больших значений $R$. Для очень больших значений $R$ оценки (3.1.5) и (3.1.7) неплохо соответствуют друг другу. Покажем это.

**Лемма 3.1.6 (см. [182]).** *Пусть в задаче выпуклой оптимизации* $f(x) \to \min_{Ax=b,\, x\in Q}$ *функция* $f(x)$ *обладает ограниченной вариацией на множестве* $Q$ : $\max_{x,y\in Q}\left(f(x)-f(y)\right)\le\Delta$ . *Предположим, что* $B_2(0,r)$ *– евклидов шар (в двойственном пространстве* $\lambda$ *– множитель Лагранжа к ограничению* $Ax=b$*) радиуса* $r$ *с центром в точке 0 полностью принадлежит множеству* $\Xi_{b,A}=\{\lambda:\ \lambda=b-Ax,\, x\in Q\}$. *Тогда имеет место следующая оценка на размер решения двойственной задачи:* $\|\lambda^*\|_2\le\Delta/r$.

Для задачи (3.1.1) $\Delta=\Delta_\varphi,\ Q=S_n(1)$. Нужно оценить $r$, чтобы можно было воспользоваться этой леммой. Следующую идею (слейтеровской релаксации) сообщил нам А.С. Немировский. Можно довольно грубо оценить снизу размер вписанного шар $r$. Для этого “подменим” исходную задачу немного другой, в которой вектор $b$ заменяется на такой вектор $b_\varepsilon$, что $\|b_\varepsilon-b\|_2\le\varepsilon$ и $B_2(0,\varepsilon)\subset\Xi_{b_\varepsilon,A}$. Тогда, для задачи с вектором $b_\varepsilon$ мы имеем оценку $r\ge\varepsilon$. Но $(\varepsilon_f,\varepsilon)$-решение этой обновленной задачи гарантировано будет $(\varepsilon_f,2\varepsilon)$-решением исходной задачи. Таким образом, в самом плохом случае можно считать, что $r\simeq\varepsilon$. Используя это наблюдение несложно сопоставить теоремы 3.1.1, 3.1.2 в случае, когда $R$ очень большое и $\varepsilon_f=\mathrm{O}(R\varepsilon)$. Следует сравнить лемму 3.1.6 и последующий текст с Слейтеровскими конструкциями, встречающимися в разделе 2.3 главы 2 и разделе 3.4 этой главы 3.

Остается только одна проблема: неизвестность $R$. Хотя в шаг БГМ (3.1.4) $R$ явно не входит, $R$ входит в $\delta$ (см. формулу (3.1.6)), от которого размер шага уже зависит. Также $R$ входит в критерий останова метода (3.1.7), заключающегося в выполнении точно рассчитанного числа итераций. Далее (следуя [47, 182]) будет описана процедура рестартов, позволяющая справиться с этой сложностью.

Полагаем $R=R_0=100$ (при $n\gg10^4$, $m\gg10^2$ в различных численных экспериментах [56] получалось $R\gg10^3$) делаем предписанное этому $R$ число итераций (3.1.7) и проверяем критерий останова – лемма 3.1.1. Заметим, что в этом критерии останова не надо знать $\varphi^*$.

Если этот критерий не выполняется, то полагаем $R := 4R$ и повторяем процедуру. Через не более чем $\lceil \log_4(4R/R_0) \rceil$ перезапусков мы остановимся.

Поясним, почему была выбрана именно константа 4. Оптимально выбирать такой коэффициент $\beta$ ( $R := \beta R$ ), который доставляет минимум следующему выражению $\beta\sqrt{R/R_0}/\left(\sqrt{\beta}-1\right)$, получающемуся при оценке сверху общего количества сделанных итераций с учетом всех необходимых перезапусков (считаем $R/R_0 \gg 1$, $\varepsilon_f = O(R\varepsilon)$ ).[3] Решением этой задачи будет $\beta = 4$ (детали см. также в разделе 4.3 главы 4).

### 3.1.7 Заключительные замечания[4]

Приведенные в разделе результаты могут быть перенесены на другие сильно выпуклые функционалы (см. следующий раздел). Выбор в качестве функционала минус энтропии лишь один из возможных вариантов. Напомним, что согласно неравенству Пинскера минус энтропия – сильно выпуклая функция в 1-норме с константой сильной выпуклости 1. При этом не так важно, чтобы решение прямой и двойственной задачи были связаны явными формулами, как это имеет место в данном разделе, и, в целом, довольно типично для сепарабельных функционалов. Достаточно, чтобы имела место сильная выпуклость функционала (см. пример 4 [44]) или его сепарабельность вместе с ограничениями. Это, в частности, обеспечивает возможность с геометрической скоростью сходимости находить приближенную зависимость $x(\lambda)$.

Если размеры матрицы в аффинных ограничениях настолько большие, что умножение такой матрицы на столбец/строку не выполнимо за разумное время (напомним, что современный процессор может делать до миллиарда операций с плавающий точкой в секунду), то требуется другая техника решения отмеченных задач. В качестве одного из вариантов альтернативного подхода, укажем на идею седлового представления задачи (как правило, это требует и некоторого пересмотра самой постановки задачи), и последующее решение седло-вой задачи рандомизированным вариантом метода зеркального спуска или проксимальным зеркальным методом А.С. Немировского [227]. При этом рандомизация возникает на этапе расчета стохастического градиента, т.е. на этапе умножения матрицы на столбец/строку (см. также раздел 2.1 главы 2). Выгода от такой рандомизации – существенное удешевление стоимости одной итерации, а плата за это – увеличение числа итераций. К сожалению, такая рандомизация приводит к заметным результатам лишь тогда, когда в прямом и двойственном пространствах переменные “живут” на 1-шарах (или симплексах), что имеет место в лишь в очень небольшом числе приложений (типа задач поиска селектора Данцига [227]). В других случаях рандомизация также используется, но эффект, как правило, заметно меньше (см. концовку работы [227], п. 6.5 обзора [163], новую работу [89]).

---

[3] Действительно, при $\varepsilon_f = O(R\varepsilon)$ теоремы 3.1.1 и 3.1.2 дают одинаковую зависимость $N \sim \sqrt{R}$. Таким образом, общее число сделанных итераций при $k$ рестартах пропорционально $N \sim \sqrt{\beta^{k+1}}/\left(\sqrt{\beta}-1\right)$, где натуральное число $k$ (большое при $R/R_0 \gg 1$) определяется из условия $\beta^{k-1}R_0 \le R \le \beta^k R_0$, т.е. $\sqrt{\beta^{k-1}} \le \sqrt{R/R_0}$. Используя это соотношение, можно оценить сверху $\sqrt{\beta^{k+1}}/\left(\sqrt{\beta}-1\right)$ как $\beta\sqrt{R/R_0}/\left(\sqrt{\beta}-1\right)$.

[4] Этот пункт написан совместно с А.В. Черновым, который получил все приведенные здесь рисунки.

Интересно было бы сравнить предложенные в разделе методы с другими известными численными методами решения задач ЭЛП (см. [55, 106, 198], а также работы M. Cuturi (http://www.iip.ist.i.kyoto-u.ac.jp/member/cuturi/ )). Этому планируется посвятить отдельную публикацию. Для более основательного знакомства с приложениями, в которых возникают задачи ЭЛП можно рекомендовать – см. [26, 30, 49, 207, 221, 230], а также работы на International workshops on Bayesian inference and maximum entropy methods in science and engineering (AIP Conf. Proceedings), проводимых каждый год. Особенно интересно сопоставить предложенные в данном разделе методы с методом балансировки [198] (говорят также методом Шелейховского, Брэгмана–Шелейховского [49], Синхорна [175]) применительно к задаче расчета матрицы корреспонденций по энтропийной модели [26, 36, 49]:

$$f(x)=\sum_{i,j=1}^{n'} x_{ij}\ln x_{ij}+\alpha\sum_{i,j=1}^{n'} c_{ij}x_{ij}\to \min_{\substack{\sum_{j=1}^{n'} x_{ij}=L_i,\ \sum_{i=1}^{n'} x_{ij}=W_j \\ i,j=1,\dots,n;\ x\in S_{n'^2}(1)}} .$$

Если $\alpha\to\infty$ и $f(x):=f(x)/\alpha$, то рассматриваемая задача ЭЛП переходит в классическую транспортную задачу ЛП, которая в общем случае решается за $\mathrm{O}\left(n'^3\ln n'\right)$ арифметических операций, причем данная оценка не улучшаема [293]. Задача же ЭЛП может быть решена намного эффективнее, в частности, методом балансировки.

Для этой задачи, которую несложно привести к виду (3.1.1), двойственная задача будет иметь вид

$$\varphi(\lambda,\mu)=\langle\lambda,L\rangle+\langle\mu,W\rangle-\ln\left(\sum_{i,j=1}^{n'}\exp\left(-\alpha c_{ij}+\lambda_i+\mu_j\right)\right)\to\max_{\lambda,\mu\in\mathbb{R}^{n'}}, \qquad (3.1.8)$$

где $\sum_{i=1}^{n'}L_i=\sum_{j=1}^{n'}W_j=1$. Метод балансировки примет вид ($[\lambda]_0=[\mu]_0=0$):

$$[\lambda_i]_{k+1} = -\ln\left(\frac{1}{L_i}\sum_{j=1}^{n'}\exp\left(-\alpha c_{ij} + [\mu_j]_k\right)\right),$$

$$[\mu_j]_{k+1} = -\ln\left(\frac{1}{W_j}\sum_{i=1}^{n'}\exp\left(-\alpha c_{ij} + [\lambda_i]_k\right)\right) \text{ или } [\mu_j]_{k+1} = -\ln\left(\frac{1}{W_j}\sum_{i=1}^{n'}\exp\left(-\alpha c_{ij} + [\lambda_i]_{k+1}\right)\right).$$

Эти формулы можно получить, если заметить, что задача оптимизации (3.1.8) может быть явно решена по $\lambda$ при фиксированном $\mu$, и наоборот. Отметим при этом, что если $(\lambda,\mu)$ – решение задачи (3.1.8), то $(\lambda + c_1 e_{n'}, \mu + c_2 e_{n'})$, где $e_{n'}$ – вектор из $n'$ единиц, $c_1$, $c_2$ – произвольные числа, также будет решением задачи (3.1.8). Собственно, выписанные формулы – есть не что иное как, метод простой итерации для явно выписываемых условий экстремума (принципа Ферма) для задачи (3.1.8):

$$\lambda = \Lambda(\mu), \ \mu = \mathrm{M}(\lambda).$$

Оператор $(\lambda,\mu) \to (\Lambda(\mu), \mathrm{M}(\lambda))$ является сжимающим в метрике Биркгофа–Гильберта [204]. Причем можно оценить коэффициент сжатия, что в итоге приводит к следующим оценкам скорости сходимости (по числу итераций) метода балансировки

$$N = \mathrm{O}\left(\sqrt{\theta}\ln\left(\frac{\varepsilon_0}{\varepsilon}\right)\right),$$

где точность $\varepsilon$ – расстояние по метрике Биркгофа–Гильберта от того, что выдает метод до решения (неподвижной точки), $\varepsilon_0$ – расстояние от точки старта до неподвижной точки, а

$$\theta = \max_{i,j,p,q=1,\ldots,n'} \exp\left(\alpha\cdot\left(c_{iq} + c_{pj} - \left(c_{ij} + c_{pq}\right)\right)\right).$$

Заметим, что параметр $\theta$ может быть содержательно проинтерпретирован исходя из эволюционного вывода модели расчета матрицы корреспонденций [36, 49].

Был проведен ряд численных экспериментов для решения задачи расчета матрицы корреспонденций (Рисунки 3.1.1 – 3.1.4). Эксперименты проводились на ЭВМ с процессором Intel Core i5, 2.5 ГГц и оперативной памятью 2 Гб в среде Matlab 2012® (8.0) под управлением операционной системой Microsoft Windows 7 (64 разрядная). Значительная часть параметров заполнялась автоматически: компоненты матрицы $\|c_{ij}\|_{i,j=1,1}^{n',n'}$ – случайные числа из интервала $(0,1)$; аналогично выбирались (а потом шкалировались) $L_i$, $W_j$, $\alpha = 100$ (на Рисунках 3.1.1, 3.1.2, 3.1.3). Требования к точности задавались относительной ошибкой в 1% (на Рисунках 3.1.1, 3.1.2, 3.1.4): $\varepsilon_f = 0.01 f(x_0)$, $\varepsilon = 0.01\|Ax_0 - b\|_2$, где $x_0 = x(\lambda = 0, \mu = 0)$ – точка старта.

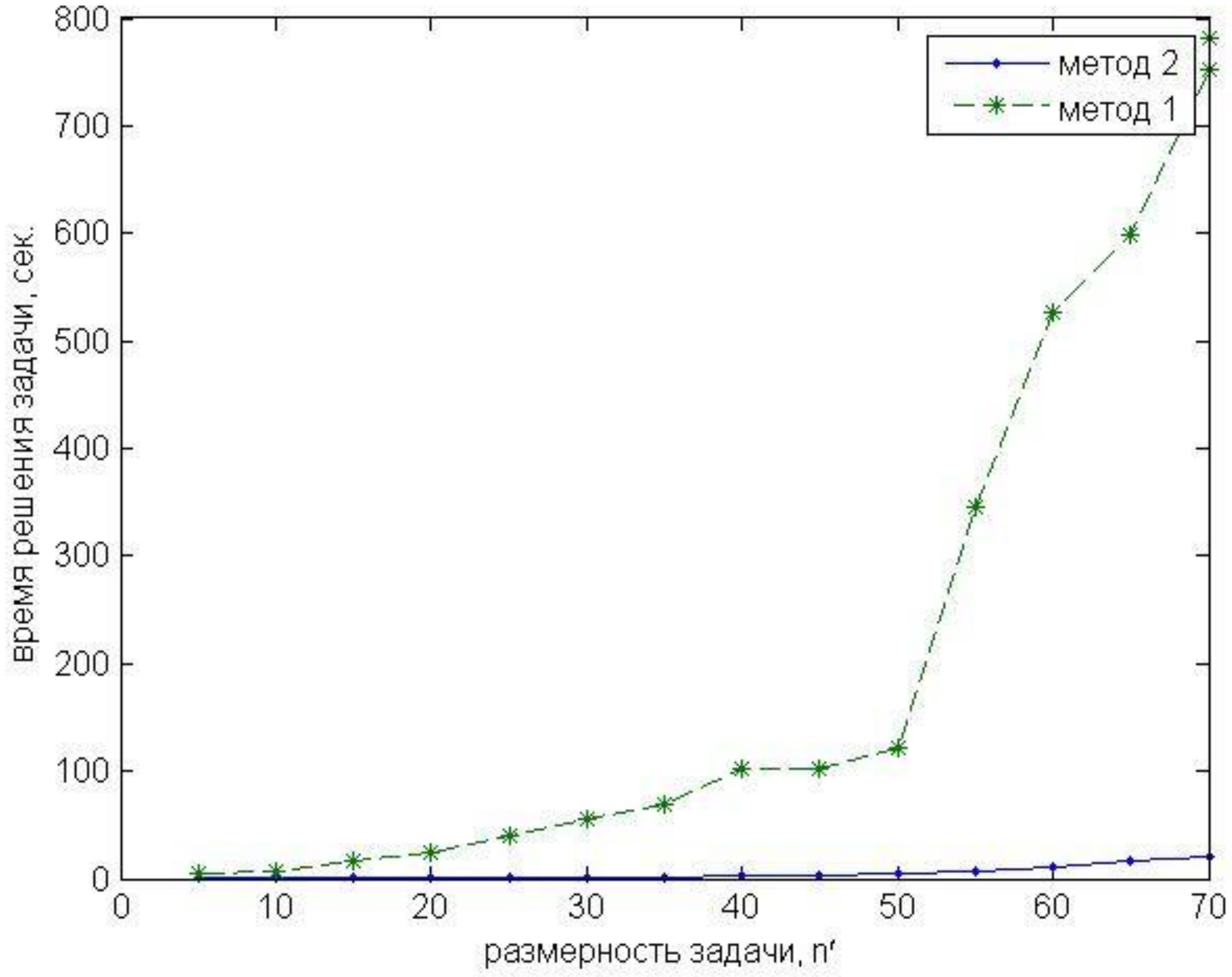

Рисунок 3.1.1. Сравнение зависимостей времен работ методов из теорем 1 и 2 от $n'$ (точность решения 1%)

График зависимости времени работы метода балансировки от точности (на той же задаче, что и на Рисунке 3.1.3) имеет вид горизонтальной прямой на уровне одной секунд, т.е. зависимость времени работы метода балансировки от точности крайне слабая (логарифмическая).

Количество итераций, которые делают методы при $\alpha = 100$, $n' = 30$ на относительной точности 1%, следующие: балансировка (58), метод из теоремы 3.1.2 (10 346), метод из теоремы 3.1.1 (221 481). Видно, что метод 3.1.1 делает много “лишних” итераций. Сложность одной итерации у методов 3.1.1 и 3.1.2 одинакова, но в методе 3.1.1 предписано сделать большее число итераций согласно полученной в теореме 3.1.1 оценке. Если бы в качестве критерия остановки в методе 3.1.1 мы взяли бы проверку условий леммы 3.1.1, то это немного увеличило бы стоимость одной итерации, но при этом число итераций стало бы приблизительно таким же, как и в методе 3.1.2. Численные эксперименты показывают, что в таком случае методы 3.1.1 и 3.1.2 требуют количества итераций, отличающиеся не более чем в 3 – 4 раза.

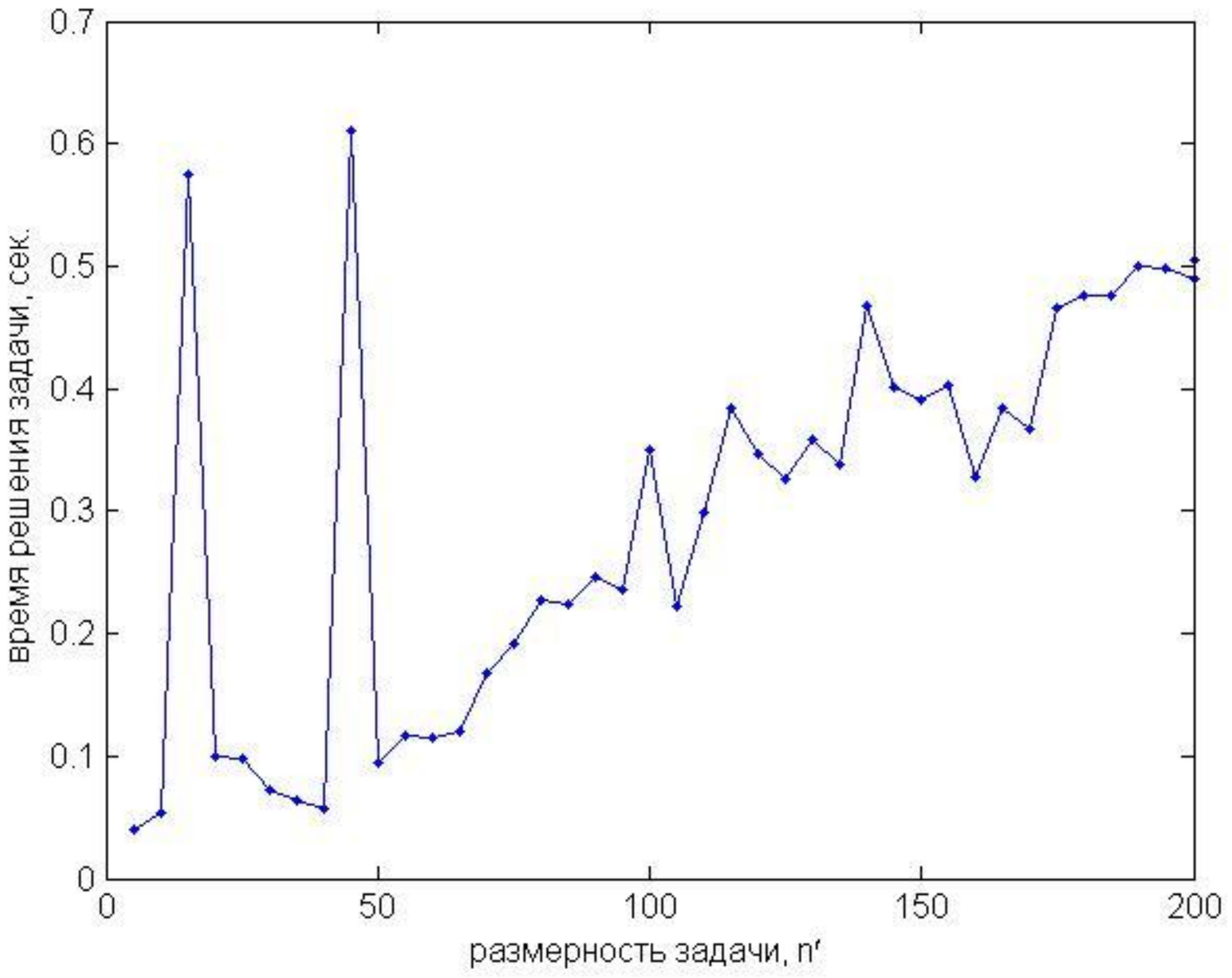

Рисунок 3.1. 2. Зависимость времени работы метода балансировки от $n'$

(точность решения 1%)

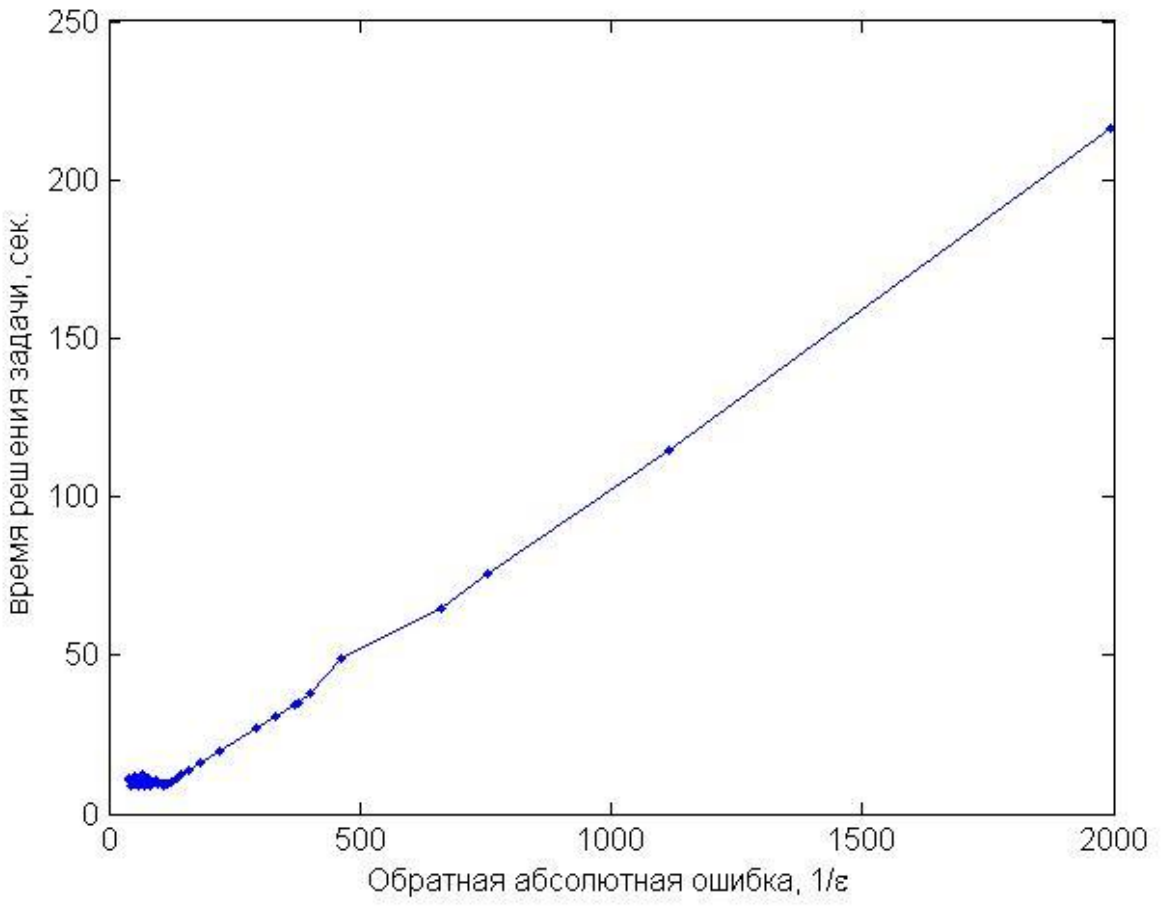

Метод из теоремы 3.1.1

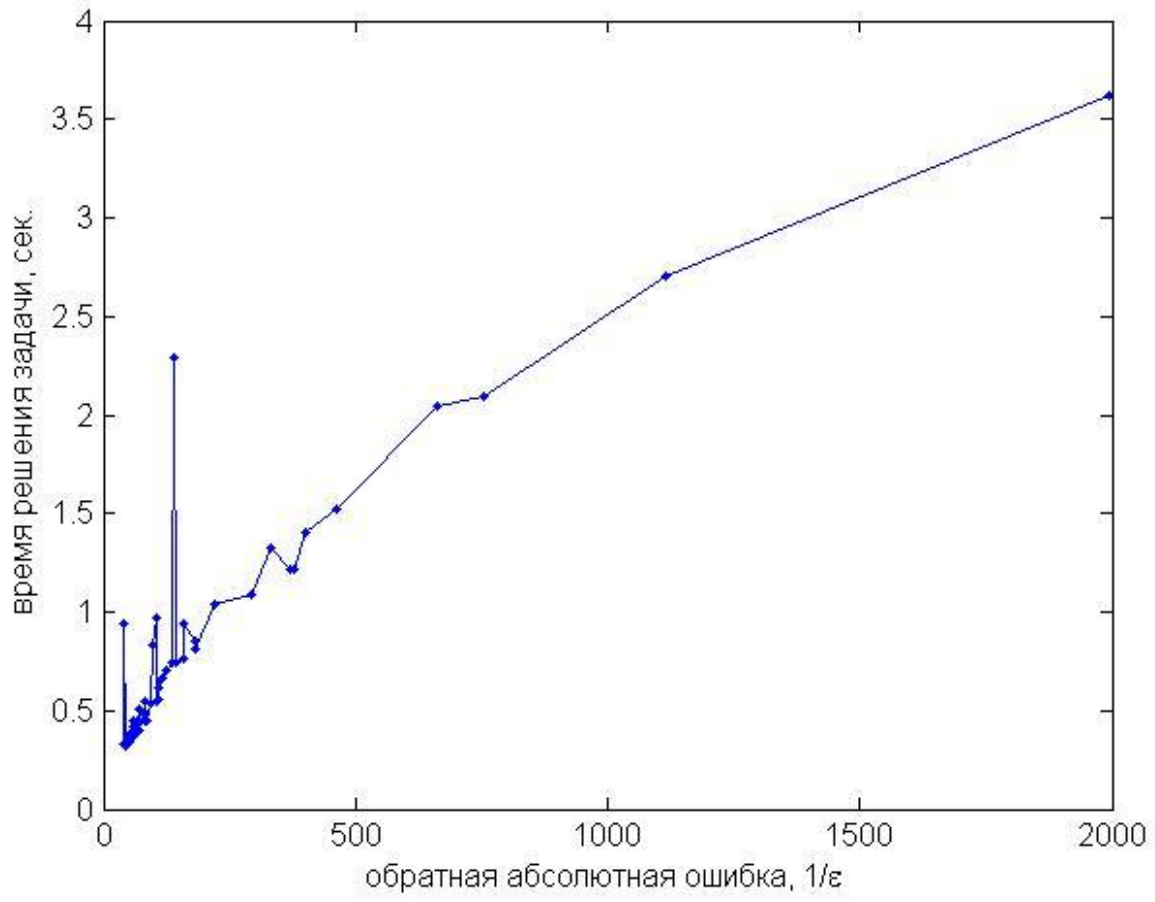

Метод из теоремы 3.1.2

Рисунок 3.1.3 Зависимости времен работ методов из теорем 3.1.1 и 3.1.2 от точности $\varepsilon$

$(n' = 25)$

Интересно отметить, что время работы метода балансировки растет с ростом $\alpha$, в то время как время работы методов 3.1.1 и 3.1.2 с ростом $\alpha$ уменьшается (см. Рисунок 3.1.4) в наблюдаемом диапазоне (далее, при больших значениях $\alpha$, наблюдалась стабилизация).

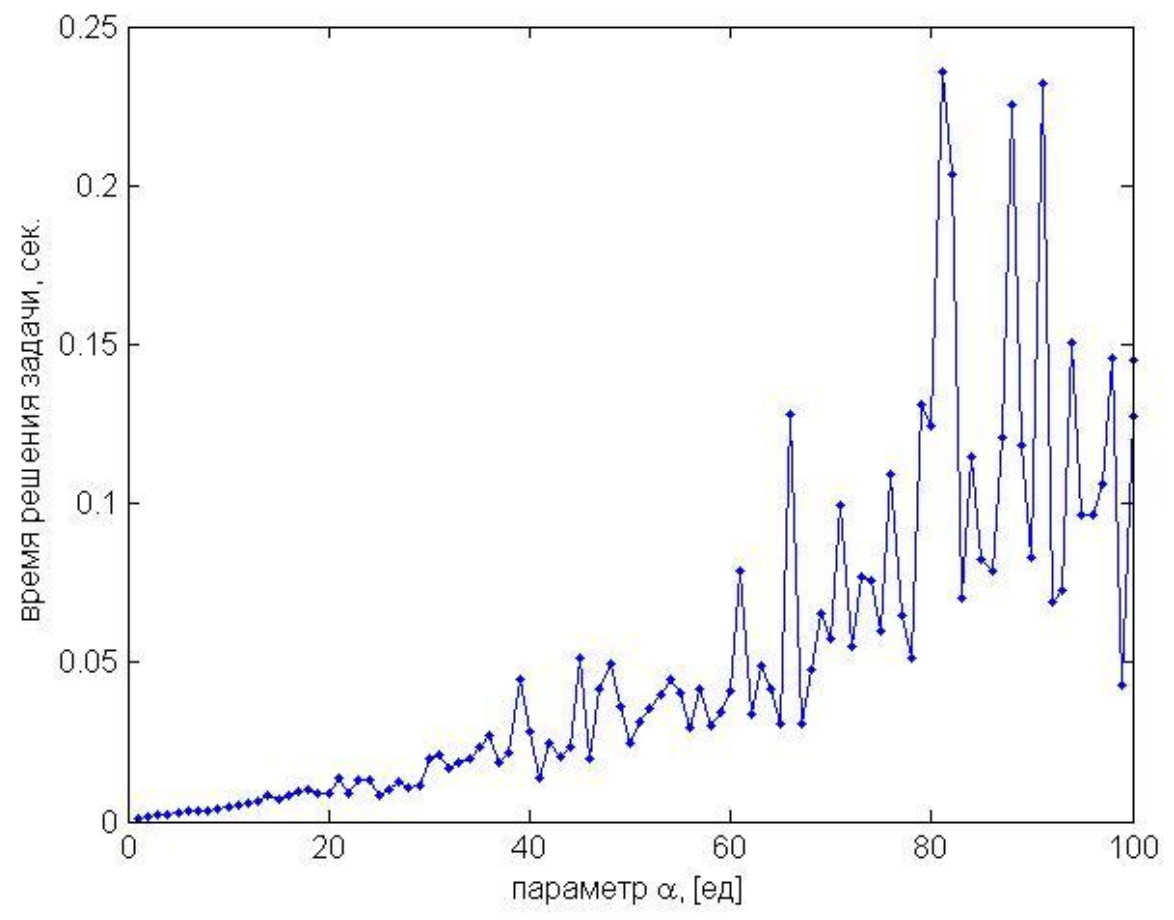

Метод балансировки

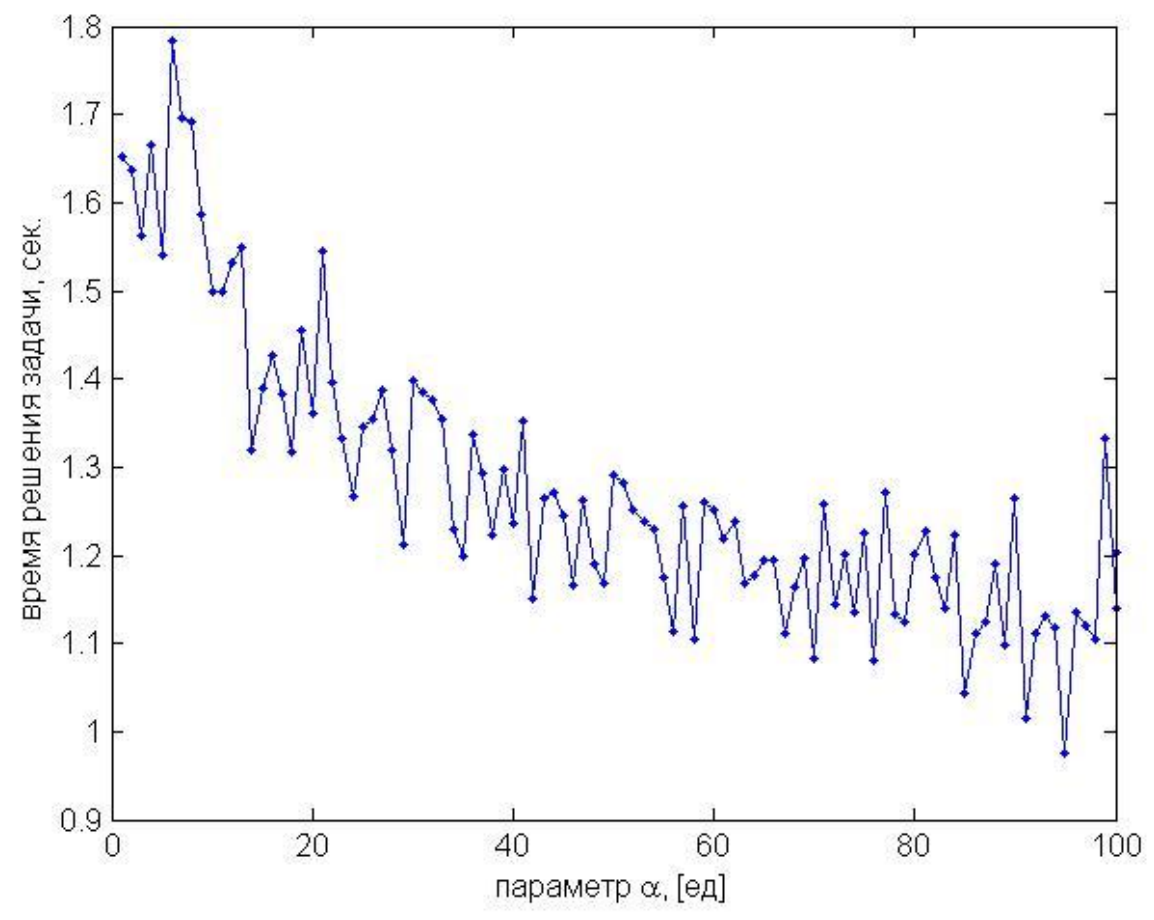

Метод из теоремы 3.1.2

Рисунок 3.1.4 Зависимости времен работ метода балансировки и метода из теорем 2 от $\alpha$ ($n' = 30$)

Проведенные численные эксперименты показали, что предложенные в разделе методы уступают по эффективности методу балансировки. Это не удивительно, поскольку балансировка считается наилучшим методам для специального класса задач ЭЛП (задач расчета матриц корреспонденций по энтропийной модели). Предложенные же в разделе методы применимы для более широкого класса задач ЭЛП, в том числе с ограничениями в виде линейных неравенств (и даже более общих конусных неравенств вида $Ax - b \in K$, с двойственным конусом $K^*$ просто структуры). Также результаты раздела не сложно перенести на случай, когда вместо единичного симплекса рассматривается множество

$$\left\{ x \in \mathbb{R}^n : x_i \geq 0,\ i = 1,\dots,n,\ \sum_{i=1}^{n} x_i \leq 1 \right\}.$$

## 3.2 Двойственные подходы к задачам минимизации сильно выпуклых функционалов простой структуры при аффинных ограничениях

### 3.2.1 Введение

В данном разделе мы существенно обобщаем результаты недавней работы [37] и предыдущего раздела (в том числе улучшаем оценку скорости сходимости метода), в котором был предложен способ решения задачи энтропийно-линейного программирования (ЭЛП) с помощью решения специальным образом регуляризованной двойственной задачи, и восстановлении по явным формулам решения прямой задачи, исходя из полученного приближенного решения двойственной задачи. Настоящий раздел развивает идеи цикла статей А.С. Немировского, Ю.Е. Нестерова с соавторами и их последователей [33, 40, 42, 43, 44, 92, 138, 181, 261, 263, 269, 281], в которых предлагались различные прямо-двойственные методы для широкого класса задач. Название "прямо-двойственные методы" было предложено в диссертации [92] для методов, которые позволяют, исходя из решения прямой (двойственной) задачи, без существенных дополнительных затрат восстанавливать (с той же точностью) решение соответствующей двойственной (прямой) задачи. Для удобства мы постарались привести в разделе все необходимые выкладки, хотя в ряде случаев они и не являются оригинальными.

В подразделе 3.2.2, следуя работе [135], мы приводим быстрый градиент метод (впоследствии, в разделе 5.1 главы 5, именно на базе такого варианта быстрого градиентного метода будет построен его покомпонентный вариант, никакие другие известные нам варианты быстрого градиентного метода, к сожалению, не позволяют также наглядно осуществить отмеченное обобщение в нужной общности). В отличие от [37] мы также исследуем прямо-двойственность [269] быстрого градиентного метода и следующее свойство: генерируемая методом последовательность точек лежит в шаре с центром в решении и радиуса равного расстоянию от точки старта метода до решения задачи (насколько нам известно, это был первый результат такого типа для быстрых градиентных методов, впоследствии, с помощью Ю.Е. Нестерова удалось распространить этот результат и на более общие ситуации – см. раздел 2.2 главы 2). Оба эти свойства необходимы в подразделе 3.2.3, чтобы обосновать способ восстановления решения прямой задачи (минимизации сильно выпуклой функции простой структуры при аффинных ограничениях), исходя из сгенерированной методом последовательности в двойственном пространстве. В конце раздела мы описываем непосредственное обобщение конструкции работы [37], связанной с регуляризацией двойственной задачи. Это обобщение не требует прямо-двойственности метода, но приводит к некоторым потерям в оценках скорости сходимости. В частности, предложенный в данном разделе метод, приме-ненный к задаче ЭЛП, находит на несколько порядков быстрее решения различных задач ЭЛП, по сравнению с методом из предыдущего раздела.

### 3.2.2 Прямо-двойственность быстрого градиентного метода

Рассмотрим задачу выпуклой оптимизации

$$f(x) \to \min_x . \tag{3.2.1}$$

Под решением этой задачи будем понимать такой $\bar{x}^N$, что

$$f(\bar{x}^N) - f_* \le \varepsilon,$$

где $f_* = f(x_*)$ – оптимальное значение функционала в задаче (3.2.1), $x_*$ – решение задачи (3.2.1). Определим множество

$$B_R(x_*) = \left\{ x: \ \|x - x_*\|_2 \le R \right\}.$$

Пусть

$$x^{k+1} = x^k - h\nabla f(x^k), \tag{3.2.2}$$

и при $x \in B_{\sqrt{2}R}(x_*)$, где

$$R = \left\| x^0 - x_* \right\|_2 = \|x_*\|_2,$$

выполняется условие

$$\|\nabla f(x)\|_2 \le M .$$

Тогда из (3.2.2) с учетом этого неравенства и (3.2.4) (см. ниже), имеем

$$\left\| x - x^{k+1} \right\|_2^2 = \left\| x - x^k + h\nabla f(x^k) \right\|_2^2 = \left\| x - x^k \right\|_2^2 + 2h\left\langle \nabla f(x^k), x - x^k \right\rangle + h^2 \left\| \nabla f(x^k) \right\|_2^2 \le$$

$$\le \left\| x - x^k \right\|_2^2 + 2h\left\langle \nabla f(x^k), x - x^k \right\rangle + h^2 M^2 .$$

Отсюда (при $x = x_*$) следует, что

$$f\left( \frac{1}{N} \sum_{k=0}^{N-1} x^k \right) - f_* \le \frac{1}{N} \sum_{k=0}^{N-1} f(x^k) - f(x_*) \le \frac{1}{N} \sum_{k=0}^{N-1} \left\langle \nabla f(x^k), x^k - x_* \right\rangle \le$$

$$\le \frac{1}{2hN} \sum_{k=0}^{N-1} \left\{ \left\| x_* - x^k \right\|_2^2 - \left\| x_* - x^{k+1} \right\|_2^2 \right\} + \frac{hM^2}{2} = \frac{1}{2hN} \left( \left\| x_* - x^0 \right\|_2^2 - \left\| x_* - x^N \right\|_2^2 \right) + \frac{hM^2}{2} .$$

Выбирая

$$h = \frac{R}{M\sqrt{N}}$$

и полагая

$$\bar{x}^N = \frac{1}{N}\sum_{k=0}^{N-1} x^k ,$$

получим

$$f\left(\bar{x}^N\right) - f_* \le \frac{MR}{\sqrt{N}} . \qquad (3.2.3)$$

Заметим, что

$$0 \le \frac{1}{2hk}\left(\left\|x_* - x^0\right\|_2^2 - \left\|x_* - x^k\right\|_2^2\right) + \frac{hM^2}{2} ,$$

Поэтому при $k = 0,...,N$

$$\left\|x_* - x^k\right\|_2^2 \le \left\|x_* - x^0\right\|_2^2 + h^2M^2k \le 2\left\|x_* - x^0\right\|_2^2 ,$$

т.е.

$$\left\|x^k - x_*\right\|_2 \le \sqrt{2}\left\|x^0 - x_*\right\|_2 ,\ k = 0,...,N . \qquad (3.2.4)$$

Для не гладких задач оценка (3.2.3) является неулучшаемой (здесь и далее неулучшаемость оценок подразумевает выполнение предположения, что размерность пространства, в котором происходит оптимизация, достаточно большая: количество итераций, которые сделает метод, не больше размерности пространства) с точностью до мультипликативного множителя [91]. Однако, если дополнительно известно, что градиент $f\left(\bar{x}^N\right)$ липшицев

$$\left\|\nabla f(y) - \nabla f(x)\right\|_2 \le L\left\|y - x\right\|_2 ,$$

где $x, y \in B_R\left(x_*\right)$ (см. (3.2.6)), то

$$\frac{1}{2L}\left\|\nabla f\left(x^k\right)\right\|_2^2 \le f\left(x^k\right) - f_* . \qquad (3.2.5)$$

Это неравенство является формальной записью простого геометрического факта, что если в точке $x^k$ к функции $f(x)$ провести касательную

$$f\left(x^k\right) + \left\langle \nabla f\left(x^k\right), x - x^k \right\rangle$$

и на основе такой касательной построить параболу

$$f\left(x^k\right) + \left\langle \nabla f\left(x^k\right), x - x^k \right\rangle + \frac{L}{2}\left\|x - x^k\right\|_2^2 ,$$

то эта парабола будет мажорировать функцию $f(x)$, т.е.

$$f(x) \le f\left(x^k\right) + \left\langle \nabla f\left(x^k\right), x - x^k \right\rangle + \frac{L}{2}\left\| x - x^k \right\|_2^2.$$

В частности выписанное неравенство имеет место и в точке минимума параболы

$$x^k - \frac{1}{L}\nabla f\left(x^k\right).$$

Поскольку приращение параболы при переходе аргумента от точки $x^k$ к точке минимума параболы составляет

$$\frac{1}{2L}\left\|\nabla f\left(x^k\right)\right\|_2^2,$$

то получаем неравенство (3.2.5). Это неравенство позволяет уточнить проведенные выше рассуждения. Как и раньше запишем

$$\left\| x^{k+1} - x_* \right\|_2^2 = \left\| x^k - x_* \right\|_2^2 - 2h\left\langle \nabla f\left(x^k\right), x^k - x_* \right\rangle + h^2\left\|\nabla f\left(x^k\right)\right\|_2^2 \le$$

$$\le \left\| x^k - x_* \right\|_2^2 - 2h\left( f\left(x^k\right) - f_* \right) + 2Lh^2\left( f\left(x^k\right) - f_* \right) = \left\| x^k - x_* \right\|_2^2 + 2h(Lh-1)\left( f\left(x^k\right) - f_* \right).$$

При $h \le 1/L$ отсюда имеем

$$\left\| x^{k+1} - x_* \right\|_2^2 \le \left\| x^k - x_* \right\|_2^2, \; k = 0,\ldots,N-1,$$

следовательно

$$\left\| x^k - x_* \right\|_2 \le \left\| x^0 - x_* \right\|_2, \; k = 0,\ldots,N. \tag{3.2.6}$$

Полагая

$$h = \frac{1}{2L},$$

получим

$$f\left(\bar{x}^N\right) - f_* \le \frac{1}{N}\sum_{k=0}^{N-1} f\left(x^k\right) - f(x_*) \le \frac{2L}{N}\sum_{k=0}^{N-1}\left\{\left\| x^k - x_* \right\|_2^2 - \left\| x^{k+1} - x_* \right\|_2^2\right\} \le \frac{2LR^2}{N}.$$

**Замечание 3.2.1 (это замечание развивает аналогичный комментарий в разделе 2.1 главы 2).** То что $f\left(\bar{x}^N\right) - f_*$ должно зависеть только от $MR/\sqrt{N}$ и(или) $LR^2/N$ является следствием П-теоремы теории размерностей [69]. Вводя точность

$$f\left(\bar{x}^N\right) - f_* \le \varepsilon,$$

можно показать, что существует всего две (независимые) безразмерные величины, сконструированные из введенных параметров:

$$\frac{M^2R^2}{\varepsilon^2} \text{ и } \frac{LR^2}{\varepsilon}.$$

Согласно П-теореме, любая безразмерная величина должна функционально выражаться через эти две (базисные). В частности,

$$N = G\left(\frac{M^2R^2}{\varepsilon^2}, \frac{LR^2}{\varepsilon}\right).$$

В случае, когда нельзя гарантировать липшицевость градиента ситуация упрощается

$$N = \tilde{G}\left(\frac{M^2R^2}{\varepsilon^2}\right).$$

Полученная методом (3.2.2) оценка (3.2.3) соответствует этой формуле. Более того, как уже отмечалось, это оценка является неулучшаемой на классе негладких выпуклых задач. К сожалению, полученная оценка скорости сходимости метода (3.2.2) с шагом $h = 1/(2L)$ в гладком случае уже не будет оптимальной.

Кстати сказать, из П-теоремы также следует, что шаг метода (3.2.2) $h$ в негладком случае должен вычисляться по формуле

$$h = c\frac{\varepsilon}{M^2}, \; c > 0,$$

которую можно получить из использовавшейся нами ранее формулы

$$h = \frac{R}{M\sqrt{N}},$$

если выразить $N$ через $\varepsilon$ с помощью формулы (3.2.3), полагая

$$\varepsilon = \frac{MR}{\sqrt{N}}.$$

В гладком случае $h$ определяться из соотношения вида

$$W\left(h\frac{M^2}{\varepsilon}, hL\right) = 1.$$

В стохастическом случае [181] (вместо градиента получаем стохастический градиент с дисперсией $\sigma^2$) из

$$\tilde{W}\left(h\frac{M^2}{\varepsilon}, hL, h\frac{\sigma}{R}\right) = 1.$$

При условии липшицевости градиента выписанная оценка скорости сходимости может быть улучшена [91]. Например, при использовании метода сопряженных градиентов [91, 99]. Локальные оценки скорости сходимости метода тяжелого шарика также говорят об этом [99]. Однако среди большого многообразия “ускоренных методов”, которые сходятся по нижним оценкам, мы особо выделим быстрый градиентный метод, предложенный в 1983 г. в кандидатской диссертации Ю.Е. Нестерова (научным руководителем был проф. Б.Т. Поляк). Помимо того, что это был один из первых методов (без использования вспомогательных одномерных и двумерных оптимизаций [91]), для которого было получено строгое доказательство глобальной сходимости согласно нижним оценкам (в гладком случае), метод оказался обладающим хорошими свойствами типа (3.2.6). Но наиболее важное свойство для нас в данном разделе – его прямо-двойственность. Развитие и использование этого метода отражено в диссертации [92].

Перейдем к построению быстрого градиентного метода (БГМ). Далее мы будем во многом следовать недавно предложенному способу понимания БГМ [135]. Тем не менее, нам потребуется из этих рассуждений получить свойство типа (3.2.6) и прямо-двойственность. Наличие у БГМ этих свойств в работе [135] не исследовалось, поэтому далее мы приведем все необходимые рассуждения в нужном объеме.

Предварительно определим две числовые последовательности шагов $\{\alpha_k, \tau_k\}$:

$$\alpha_1 = \frac{1}{L},\ \alpha_k^2 L = \alpha_{k+1}^2 L - \alpha_{k+1},\ \tau_k = \frac{1}{\alpha_{k+1} L}.$$

Можно написать явные формулы. Нам так же пригодится упрощенный вариант этих последовательностей [135], который определяется следующим образом:

$$\alpha_1 = \frac{1}{L},\ \alpha_k^2 L = \alpha_{k+1}^2 L - \alpha_{k+1} + \frac{1}{4L},\ \tau_k = \frac{1}{\alpha_{k+1} L}.$$

В таком случае

$$\boxed{\alpha_{k+1} = \frac{k+2}{2L},\ \tau_k = \frac{1}{\alpha_{k+1} L} = \frac{2}{k+2}.}$$

$$\text{БГМ}\left(x^0 = y^0 = z^0\right)$$

1. $x^{k+1} = \tau_k z^k + (1 - \tau_k) y^k$;
2. $y^{k+1} = x^{k+1} - \frac{1}{L}\nabla f\left(x^{k+1}\right)$;
3. $z^{k+1} = z^k - \alpha_{k+1}\nabla f\left(x^{k+1}\right)$.

Из последней формулы в доказательстве леммы 4.3 [135] имеем (для всех $x$)

$$\alpha_{k+1}^2 L f\left(y^{k+1}\right) - \left(\alpha_{k+1}^2 L - \alpha_{k+1}\right) f\left(y^k\right) \le$$

$$\le \alpha_{k+1}\left\{ f\left(x^{k+1}\right) + \left\langle \nabla f\left(x^{k+1}\right), x - x^{k+1}\right\rangle\right\} + \frac{1}{2}\left\|z^k - x\right\|_2^2 - \frac{1}{2}\left\|z^{k+1} - x\right\|_2^2.$$

Просуммировав то что получается по $k = 0, \ldots, N-1$, получим

$$\alpha_N^2 Lf\left(y^N\right) \le \min_x \left\{ \sum_{k=0}^{N-1} \alpha_{k+1} \left\{ f\left(x^{k+1}\right) + \left\langle \nabla f\left(x^{k+1}\right), x - x^{k+1} \right\rangle \right\} + \frac{1}{2}\left\| z^k - x \right\|_2^2 - \frac{1}{2}\left\| z^{k+1} - x \right\|_2^2 \right\} \le$$

$$\le \left( \sum_{k=0}^{N-1} \alpha_{k+1} \right) f_* + \frac{1}{2}\left\| z^0 - x_* \right\|_2^2 - \frac{1}{2}\left\| z^N - x_* \right\|_2^2. \qquad (3.2.7)$$

Заметим, что

$$\max\left\{ \left\| x^k - x_* \right\|_2, \left\| y^k - x_* \right\|_2, \left\| z^k - x_* \right\|_2 \right\} \le \left\| x^0 - x_* \right\|_2, \ k = 0,...,N. \qquad (3.2.8)$$

Действительно, полагая в формуле (3.2.7) $N := k$, $k := i$, и учитывая, что $\alpha_i$ не зависят от $k$, а $f\left(y_k\right) \ge f_*$ и $\alpha_k^2 L = \sum_{i=0}^{k-1} \alpha_{i+1}$, получим (в точности эту же формулу можно получить и для упрощенной схемы выбора шагов)

$$\left\| z^k - x_* \right\|_2^2 \le \left\| z^0 - x_* \right\|_2^2.$$

Используя это неравенство, неравенство (3.2.6) для $y^{k+1}$ и выпуклость квадрата евклидовой нормы, получим

$$\left\| y^{k+1} - x_* \right\|_2^2 \le \left\| x^{k+1} - x_* \right\|_2^2 = \left\| \tau_k \cdot \left( z^k - x_* \right) + \left(1 - \tau_k\right)\left( y^k - x_* \right) \right\|_2^2 \le$$

$$\le \tau_k \left\| z^k - x_* \right\|_2^2 + \left(1 - \tau_k\right)\left\| y^k - x_* \right\|_2^2 \le \tau_k \left\| z^0 - x_* \right\|_2^2 + \left(1 - \tau_k\right)\left\| y^k - x_* \right\|_2^2 =$$

$$= \tau_k \left\| x^0 - x_* \right\|_2^2 + \left(1 - \tau_k\right)\left\| y^k - x_* \right\|_2^2 = \tau_k \left\| y^0 - x_* \right\|_2^2 + \left(1 - \tau_k\right)\left\| y^k - x_* \right\|_2^2.$$

Отсюда по индукции получаем (3.2.8).

Возвращаясь к формуле (3.2.7), установим прямо-двойственность БГМ. Для этого перепишем формулу (3.2.7) для упрощенного варианта выбора шагов

$$\frac{\left(N+1\right)^2}{4L} f\left(y^N\right) + \sum_{k=1}^{N-1} \frac{1}{4L} f\left(y^k\right) \le \min_x \left\{ \sum_{k=0}^{N-1} \frac{k+2}{2L} \left\{ f\left(x^{k+1}\right) + \left\langle \nabla f\left(x^{k+1}\right), x - x^{k+1} \right\rangle \right\} + \frac{1}{2}\left\| z^0 - x \right\|_2^2 \right\},$$

т.е.

$$f\left(\breve{y}^N\right) \le \frac{4L}{N \cdot \left(N+3\right)} \min_x \left\{ \sum_{k=0}^{N-1} \frac{k+2}{2L} \left\{ f\left(x^{k+1}\right) + \left\langle \nabla f\left(x^{k+1}\right), x - x^{k+1} \right\rangle \right\} + \frac{1}{2}\left\| z^0 - x \right\|_2^2 \right\}, \qquad (3.2.9)$$

где

$$\breve{y}^N = \frac{1}{N \cdot \left(N+3\right)} \left( \sum_{k=0}^{N-1} y^k + \left(N+1\right)^2 y^N \right).$$

Собственно именно неравенство (3.2.9) и даст нам возможность в следующем пункте восстанавливать решение прямой задачи, исходя из решения двойственной методом БГМ.

Сформулируем основной результат данного пункта.

**Теорема 3.2.1.** *Пусть функционал задачи (3.2.1) обладает свойством*

$$\left\|\nabla f(y)-\nabla f(x)\right\|_2 \le L\left\|y-x\right\|_2,\ x,y\in B_R(x_*). \tag{3.2.10}$$

*Тогда БГМ генерирует такую последовательность точек* $\left\{x^k, y^k, z^k\right\}_{k=0}^N$, *что имеют место соотношения (3.2.8) и (3.2.9), причем формулу (3.2.9) можно переписать следующим образом*

$$f\left(\breve{y}^N\right)-f_* \le \frac{2L\left\|z^0-x_*\right\|_2^2}{N\cdot(N+3)} \le \frac{2L\left\|z^0-x_*\right\|_2^2}{(N+1)^2} = \frac{2LR^2}{(N+1)^2}.$$

**Замечание 3.2.2.** Новизна в этой теореме по сравнению со всеми известными ее аналогами заключается в том, что хотя задача оптимизации (3.2.1) решается на неограниченном множестве, параметры, входящие в оценки скорости сходимости, определяются расстоянием от точки старта до решения. То есть весь итерационный процесс будет находиться в (евклидовом) шаре с центром в решении и радиуса равного расстоянию от точки старта до решения. Поскольку мы априорно, как правило, не знаем это расстояние, то может показаться, что ценность в таком замечании небольшая. Однако, как будет продемонстрировано в следующем пункте, сделанное замечание играет важную роль в обосновании предлагаемого прямо-двойственного подхода (впрочем, существует другой способ рассуждений, недавно предложенный П.Е. Двуреченским, при котором условие типа (3.2.8) не используется [171]). Интересно также заметить, что если решение задачи (3.2.1) не единственно (обозначим множество решений через $X$ ), то под $R^2$ можно понимать $\min_{x_*\in X}\left\|z^0-x_*\right\|_2^2$.

**Замечание 3.2.3.** Приведенную выше теорему 3.2.1 можно распространить и на случай, когда множество, на котором происходит оптимизация, не совпадает со всем пространством. Например, является неотрицательным ортантом или симплексом. В общем случае это требует также введения прокс-структуры [92, 259] в задачу, отличной от использованной нами выше – евклидовой. Все это несколько усложняет выкладки. В частности, в описании БГМ нужно использовать шаги типа прямого градиентного метода, а их проксимальные варианты [181]. Насколько нам известно, пока это все сделано (в общности теоремы 1) только для евклидовой прокс-структуры, но с произвольными множествами, на которых ведется оптимизация.

**Замечание 3.2.4.** В действительности, описанную выше конструкцию, с сохранением основного результата – теоремы 3.2.1, можно перенести на задачи композитной оптимизации [259, 266], на задачи, в которых не известна константа $L$, и ее требуется подбирать по ходу процесса [92, 259]. Также на основе БГМ и только что сделанного замечания, можно постро-ить и соответствующий вариант универсального метода Ю.Е. Нестерова [274]. Все описанные обобщения можно сделать и в концепции неточного оракула [39, 42, 44, 48, 181,193].

### 3.2.3 Приложение к задаче минимизации сильно выпуклого функционала простой структуры при аффинных ограничениях

Пусть требуется решать задачу

$$g(x) \to \min_{Ax=b,\, x\in Q}, \qquad (3.2.11)$$

где функция $g(x)$ – 1-сильно выпуклая в $p$-норме $(1 \le p \le 2)$. Построим двойственную задачу

$$F(y) = \max_{x\in Q}\{\langle y, b - Ax\rangle - g(x)\} \to \min_{y}. \qquad (3.2.12)$$

Во многих важных приложениях основной вклад в вычислительную сложность внутренней задачи максимизации дает умножение $Ax$ ( $A^T y$ ). Это так, например, для сепарабельных функционалов

$$g(x) = \sum_{k=1}^{n} g_k(x_k)$$

и параллелепипедных ограничениях $Q$. В частности, это имеет место для задач энтропийно-линейного программирования (ЭЛП) [37], в которых имеется явная формула $x(y)$ (см. также предыдущий раздел).

В общем случае внутренняя задача максимизации не решается точно (по явным формулам). Тем не менее, за счет сильной выпуклости $g(x)$ (аналогичное можно сказать в случае сепарабельности $g(x)$, но отсутствии сильной выпуклости) точность решения этой вспомогательной задачи (на каждой итерации внешнего метода) входит в оценку сложности ее решения логарифмическим образом, как следствие, оговорки о неточности оракула, выдающего градиент для внешней задачи минимизации $F(y)$, можно опустить. Аккуратный учет этого всего приводит лишь к логарифмическим поправкам в итоговых оценках сложности метода (см., например, [44, 48, 138, 259], а также раздел 2.3 главы 2). Поэтому для большей наглядности (следуя совету А.С. Немировского) мы далее в рассуждениях будем просто считать, что есть явная формула $x(y)$.

Применим к задаче (3.2.12) БГМ из п. 2 ( $z^0 = 0$ ), получим по формуле (3.2.9) (были сделаны переобзначения: $x \to y$, $\breve{y} \to \tilde{y}$ )

$$F\left(\tilde{y}^N\right) \le \frac{4L}{N\cdot(N+3)} \min_y \left\{ \sum_{k=0}^{N-1} \frac{k+2}{2L} \left\{ F\left(y^{k+1}\right) + \left\langle \nabla F\left(y^{k+1}\right), y - y^{k+1} \right\rangle \right\} + \frac{1}{2}\left\| z^0 - y \right\|_2^2 \right\}, \qquad (3.2.13)$$

где (см., например, [92])

$$L = \max_{\|x\|_p \le 1} \|Ax\|_2^2 .$$

В частности, для задачи ЭЛП $p = 1$ [37, 44]

$$L = \max_{k=1,\dots,n} \left\| A^{\langle k \rangle} \right\|_2^2 ,$$

где $A^{\langle k \rangle}$ – $k$ -й столбец матрицы $A^{\langle k \rangle}$ . Для задачи PageRank (см. замечание 3.2.10 ниже) $p = 2$

$$L = \lambda_{\max}\left(A^T A\right) = \sigma_{\max}\left(A\right).$$

Из неравенства (3.2.13) имеем ( $R^2 = \left\| z^0 - y_* \right\|_2^2 = \left\| y_* \right\|_2^2$, где $y_*$ – решение задачи (3.2.12))

$$F\left(\tilde{y}^N\right) - \min_{y \in B_{3R}(0)} \left\{ \sum_{k=0}^{N-1} \frac{2(k+2)}{N\cdot(N+3)} \left\{ F\left(y^{k+1}\right) + \left\langle \nabla F\left(y^{k+1}\right), y - y^{k+1} \right\rangle \right\} \right\} \le \frac{18LR^2}{(N+1)^2} \overset{def}{=} \gamma_N . \qquad (3.2.14)$$

Введем

$$\lambda_k = \frac{2(k+2)}{N\cdot(N+3)},$$

$$x^N = \sum_{k=0}^{N-1} \lambda_k x\left(y^{k+1}\right) = \frac{2}{N\cdot(N+3)} \sum_{k=0}^{N-1} (k+2)\, x\left(y^{k+1}\right) = \frac{N^2+N-2}{N\cdot(N+3)} x^{N-1} + 2\frac{N+1}{N\cdot(N+3)} x\left(y^N\right).$$

Перепишем неравенство (3.2.14) исходя из определения $F(y)$ (3.2.12) и $x(y)$ (приводимая далее выкладка аналогична рассуждениям из п. 3 работы [263])

$$F\left(\tilde{y}^N\right) - \sum_{k=0}^{N-1} \lambda_k \left\langle y^{k+1}, b - Ax\left(y^{k+1}\right) \right\rangle + \sum_{k=0}^{N-1} \lambda_k g\left(x\left(y^{k+1}\right)\right) -$$

$$- \min_{y \in B_{3R}(0)} \left\{ \sum_{k=0}^{N-1} \lambda_k \left\langle b - Ax\left(y^{k+1}\right), y - y^{k+1} \right\rangle \right\} \le \gamma_N .$$

Учитывая, что

$$\sum_{k=0}^{N-1} \lambda_k = 1,$$

получим

$$F\left(\tilde{y}^N\right)+g\left(\sum_{k=0}^{N-1}\lambda_k x\left(y^{k+1}\right)\right)+\max_{y\in B_{3R}(0)}\left\{\left\langle A\sum_{k=0}^{N-1}\lambda_k x\left(y^{k+1}\right)-b,y\right\rangle\right\}\le\gamma_N,$$

т.е.

$$F\left(\tilde{y}^N\right)+g\left(x^N\right)+3R\left\|Ax^N-b\right\|_2\le\gamma_N.$$

Отсюда (приводимые далее рассуждения во многом повторяют рассуждениями из п. 6.11 работы [261]), используя то, что

$$Ax_*=b,$$

и слабую двойственность

$$-g\left(x_*\right)\le F\left(y_*\right),$$

получаем

$$g\left(x^N\right)-g\left(x_*\right)\le g\left(x^N\right)+F\left(y_*\right)\le g\left(x^N\right)+F\left(\tilde{y}^N\right)\le g\left(x^N\right)+F\left(\tilde{y}^N\right)+3R\left\|Ax^N-b\right\|_2\le\gamma_N.$$

Исходя из определения $F(y)$ (3.2.12) и свойства (3.2.8), имеем

$$-g\left(x_*\right)=\left\langle y_*,b-Ax_*\right\rangle-g\left(x_*\right)=F\left(y_*\right)\ge\left\langle y_*,b-Ax^N\right\rangle-g\left(x^N\right)\;\Rightarrow\;g\left(x_*\right)-g\left(x^N\right)\le R\left\|Ax^N-b\right\|_2,$$

$$R\left\|Ax^N-b\right\|_2\le-g\left(x^N\right)+\left\langle\tilde{y}^N,b-Ax^N\right\rangle+g\left(x^N\right)+3R\left\|Ax^N-b\right\|_2\le$$

$$\le F\left(\tilde{y}^N\right)+g\left(x^N\right)+3R\left\|Ax^N-b\right\|_2\le\gamma_N.$$

Отсюда следует, что

$$\left|g\left(x^N\right)-g\left(x_*\right)\right|\le\gamma_N,\; R\left\|Ax^N-b\right\|_2\le\gamma_N.$$

Поскольку

$$g\left(x^N\right)-g\left(x_*\right)\le F\left(\tilde{y}^N\right)+g\left(x^N\right)\le\gamma_N,$$

то мы приходим к следующему результату.

**Теорема 3.2.2.** *Пусть нужно решить задачу (3.2.11) посредством перехода к задаче (3.2.12), исходя из выписанных выше формул. Выбираем в качестве критерия останова БГМ выполнение следующих условий*

$$F\left(\tilde{y}^N\right)+g\left(x^N\right)\le\varepsilon,\;\left\|Ax^N-b\right\|_2\le\tilde{\varepsilon}.$$

*Тогда БГМ гарантированно остановится, сделав не более чем*

$$\max\left\{\sqrt{\frac{18LR^2}{\varepsilon}},\sqrt{\frac{18LR}{\tilde{\varepsilon}}}\right\}$$

*итераций.*

**Замечание 3.2.5.** Обратим внимание, что при описании подхода одновременного решения прямой и двойственной задачи мы имели неизвестный параметр $R$. Однако этот параметр не входит в описание алгоритма и в критерий его останова. Он входит только в оценку числа итераций. К сожалению, такого удается добиться далеко не всегда. Нетривиально то, что нам удалось этого добиться в данном контексте. Обычно при решении двойственной задачи искусственно компактифицируют [42, 261, 269] множество, на котором происходит оптимизация (для двойственной задачи это, как правило, либо все пространство, либо прямое произведение пространства на неотрицательный ортант). В результате используются методы, в которых требуется проектирование на шар заранее неизвестного радиуса. Эту проблему (неизвестности размера двойственного решения) обычно решают либо с помощью процедуры рестартов [37, 44, 181], что, как правило, приводит к увеличению числа итераций как минимум на порядок, либо с помощью слейтеровской релаксации, которая часто еще более затратна в смысле требуемого числа итераций [37].

**Замечание 3.2.6.** К сожалению, в ряде приложений имеет место лишь строгая (не сильная) выпуклость $g(x)$. В этом случае (хотя двойственная задача и будет гладкой) ничего нельзя сказать о константе Липшица градиента, которая явно входит в шаг БГМ. Как уже отмечалось ранее, проблема решается адаптивным подбором $L$, а в более общем случае (когда и гладкость $g(x)$ нельзя гарантировать) с использованием универсального метода [42, 48, 274]. Тем не менее, ранее открытым был вопрос, насколько все эти конструкции (описанные выше) работают при решении гладкой двойственной задачи (на неограниченном множестве) с не ограниченными равномерно константами Липшица (на этом множестве). Ведь если предположить, что метод может отдаляться от решения по ходу итерационного процесса больше, чем на расстояние, которое было в момент старта, и расстояние, на которое он может отдаляться, зависит от свойств гладкости функционала, то образуется “порочный круг”. При неаккуратном оценивании так и получается. В данном разделе было показано, что для детерминированных постановок задач при правильном подборе прямо-двойственных методов таких проблем можно избежать. Причем избежать естественным образом, т.е. не за счет искусственной компактификации (как это общепринято), приводящей к дополнительным затратам на рестарты.

**Замечание 3.2.7.** Класс задач, к которым применим описанный выше подход можно существенно расширить, допуская, например, в постановке задачи (3.2.11) ограничения вида неравенств $Cx \le d$ (это недавно было сделано П.Е. Двуреченским и А.В. Черновым [171]), и в более общем случае вида $Cx - d \in K$, где конус $K$ имеет простое двойственное описание [289]. Можно исходить не из задач вида (3.2.11), а сразу из задачи минимизации функционала, имеющего лежандрово представления вида (3.2.12) [44]. При этом минимизация по $y$ теперь может вестись по произвольному выпуклому множеству.

**Замечание 3.2.8.** В случае, когда размерность двойственного пространства небольшая можно использовать вместо БГМ метод эллипсоидов (не требующий гладкости двойственного функционала). Этот метод также является прямо-двойственным [261]. Интересные примеры в связи с этой конструкцией возникают, когда размерность прямого пространства огромна, но при этом есть эффективный линейный минимизационный оракул, позволяющий (несмотря на огромную размерность прямого пространства) эффективно вычислять градиент двойственного функционала [33, 40, 173].

**Замечание 3.2.9.** Полезно заметить, что описанный в этом пункте подход (особенно в связи с замечанием 3.2.8) позволяет во многих важных случаях решать задачу поиска градиентного отображения, возникающую (при проектировании на допустимое множество не очень простой структуры) на каждой итерации большинства итерационных методов [33, 92, 259]. Собственно, общая идея “разделяй и властвуй” применительно к численному решению задач выпуклой оптимизации приобрела форму итерационного процесса: на каждом шаге которого решается более простая задача, чем исходная. Можно играть на том насколько сложная задача решается на каждом шаге и том, как много таких шагов надо сделать. Хороший пример здесь – это композитная оптимизация [259, 266]. Часть сложности постановки задачи (в виде композита) перенесена полностью (без лианеризации) в задачу, которую требуется решать на каждой итерации. Если композит не очень плохой, то это не сильно меняет стоимость итерации, зато может приводить к существенному сокращению необходимого числа итераций (например, для негладкого композита). Другие примеры на эту же тему имеются здесь [42, 43, 44, 48]. Общая линия рассуждений тут приблизительно такая: любое усложнение итерации, как правило, приводит к сокращению их числа, с другой стороны в стоимость одной итерации всегда входит расчет (пересчет) градиента или его (скажем, стохастического) аналога, использующегося в методе. Для детерминированных методов, когда используется градиент, основная стоимость итерации формируется как раз за счет расчета этого градиента (как правило, это умножение матрицы на вектор, т.е. порядка $O(n^2)$ арифметических операций). При рассчитанном градиенте сделать шаг итерационного процесса стоит обычно не более $O(n\ln(n/\varepsilon))$ арифметических операций. Таким образом, остается довольно большой зазор, в который можно занести дополнительные вычисления, перенеся еще какую-то часть сложности постановки задачи на каждый шаг, в надежде сократить число ша-гов. В частности, если рассматривать задачу [48, 138] (см. также раздел 2.3 главы 2 и раздел 5.1 главы 5)

$$\frac{1}{2}\|Ax-b\|_2^2+\mu\sum_{k=1}^{n}x_k\ln x_k\to\min_{x\in S_n(1)},$$

с достаточно большим $\mu>0$, то необходимо учитывать сильную выпуклость энтропийного композита. Обычный БГМ в композитном варианте для этой задачи требует на каждой итерации решения некоторой почти-сепарабельной задачи. При решении такой задачи описанные в этом пункте конструкции оказываются весьма полезными. Более того, мы привели здесь эту задачу также и потому, что при определенном значении параметра $\mu>0$ эта задача "эквивалентна" исходной задаче (3.2.11) с функционалом вида энтропии. В работе [138] обсуждаются (с содержательной стороны и со стороны практических вычислений) эти разные способы понимания одной задачи, а также приводятся некоторые конструкции (близкие к описанным в этом пункте), которые позволяют за небольшую дополнительную плату определять параметр $\mu>0$, чтобы имело место отмеченное соответствие (см. также [48]).

**Замечание 3.2.10 (PageRank и нижние оценки – это замечание повторяет аналогичное замечание из раздела 2.1 главы 2).** Приведенные в теореме 3.2.2 оценки при первом взгляде могут приводить к противоречиям. Поясним это следующим примером [44], принадлежащим Ю.Е. Нестерову. Задача поиска такого $x^*\in\mathbb{R}^n$, что

$$Ax^*=b$$

сводится к задаче выпуклой гладкой оптимизации

$$f(x)=\|Ax-b\|_2^2\to\min_x.$$

Нижняя оценка для скорости решения такой задачи [91] имеет вид

$$f(x^N)\ge\Omega\left(L_xR_x^2/N^2\right),\ L_x=\sigma_{\max}(A),\ R_x=\|x^*\|_2.$$

Откуда следует, что только при

$$N\ge\Omega\left(\sqrt{L_x}R_x/\varepsilon\right)$$

можно гарантировать выполнение неравенства

$$f(x^N)\le\varepsilon^2,\ \text{т.е.}\ \|Ax^N-b\|_2\le\varepsilon.$$

Однако эта нижняя оценка для специальных матриц может быть улучшена. Рассмотрим задачу поиска вектора PageRank [46] ($n\sim10^{10}$), которую мы перепишем как

$$Ax=\begin{pmatrix}(P^T-I)\\ 1\ldots\ldots1\end{pmatrix}x=\begin{pmatrix}0\\1\end{pmatrix}=b\ ,$$

где $I$ – единичная матрица. По теореме Фробениуса–Перрона [95] решение такой системы с неразложимой стохастической матрицей $P$ единственно и положительно $x > 0$. Сведем решение этой системы уравнений к вырожденной задаче выпуклой оптимизации

$$\frac{1}{2}\|x\|_2^2 \to \min_{Ax=b}.$$

Построим двойственную к ней задачу (поскольку система $Ax = b$ совместна, то по теореме Фредгольма не существует такого $y$, что $A^T y = 0$ и $\langle b, y\rangle > 0$, следовательно, двойственная задача имеет конечное решение)

$$\min_{Ax=b}\frac{1}{2}\|x\|_2^2 = \min_x \max_y \left\{\frac{1}{2}\|x\|_2^2 + \langle b - Ax, y\rangle\right\} =$$

$$= \max_y \min_x \left\{\frac{1}{2}\|x\|_2^2 + \langle b - Ax, y\rangle\right\} = \max_y \left\{\langle b, y\rangle - \frac{1}{2}\left\|A^T y\right\|_2^2\right\}.$$

Зная решение $y^*$ двойственной задачи (например, с минимальной евклидовой нормой)

$$\langle b, y\rangle - \frac{1}{2}\left\|A^T y\right\|_2^2 \to \max_y$$

можно восстановить решение прямой задачи

$$x(y) = A^T y.$$

Кроме того, если численно решать двойственную задачу БГМ, то согласно теореме 3.2.2

$$\left\|Ax^N - b\right\|_2 \le \frac{8L_y R_y}{N^2},$$

где

$$L_y = \sigma_{\max}\left(A^T\right) = \sigma_{\max}(A) = L_x,\ R_y = \left\|y^*\right\|_2.$$

Кажется, что это противоречит нижней оценке

$$\left\|Ax^N - b\right\|_2 \ge \Omega\left(\sqrt{L_x} R_x / N\right).$$

Однако важно напомнить [91], что эта нижняя оценка установлена для всех $N \le n$ ($n$ – размерность вектора $x$), и она будет улучшена, в результате описанной процедуры только если дополнительно предположить, что матрица $A$ удовлетворяет следующему условию

$$L_y R_y \ll n\sqrt{L_x} R_x,$$

что сужает класс, на котором была получена нижняя оценка

$$\Omega\left(\sqrt{L_x} R_x / N\right).$$

В типичных ситуациях можно ожидать $R_y \gg R_x$, что мешает выполнению требуемого условия.

**Замечание 3.2.11.** Определенный нюансы возникают при попытке перенесения результатов данного раздела на случай, когда вместо градиента доступен только стохастический градиент. Мы не будем здесь подробно на этом останавливаться (этому случаю планируется посвятить отдельную работу). Тем не менее, отметим, что для рандомизированных покомпонентных методов ответы на многие вопросы уже удалось получить [44]. В частности, задачу из замечания 3.2.10 (PageRank) можно решать прямым ускоренным покомпонентным методом или двойственным. Оценки получаются следующими (см. [44], а также раздел 5.1 главы 5)

$$E\left[\left\|Ax^N - b\right\|_2^2\right] = E\left[f\left(x^N\right)\right] - f_* = \mathrm{O}\left(n^2 \frac{\bar{L}_x R_x^2}{N^2}\right),\ \bar{L}_x^{1/2} = \frac{1}{n}\sum_{k=1}^{n}\left\|A^{\langle k\rangle}\right\|_2 \le 2,\ R_x^2 \le 2$$

(для прямой задачи),

$$E\left[\left\|Ax^N - b\right\|_2\right] = \mathrm{O}\left(n^2 \frac{\bar{L}_y R_y}{N^2}\right),\ \bar{L}_y^{1/2} = \frac{1}{n+1}\sum_{k=1}^{n+1}\left\|A_k\right\|_2$$

(для двойственной задачи).

При этом если матрица $P$ имеет $sn$ ненулевых элементов ($s \ll n$), то одна итерация у обоих методов в среднем требует $\mathrm{O}(s)$ арифметических операций. Без численных экспериментов (исходя из приведенных выше оценок) сложно определить, какой подход будет предпочтительнее (в основном из-за не знания $R_y$). Этот пример является, чуть ли не единственным примером, когда удается так организовать вычисления, что сполна можно использовать разреженность задачи (итерация выполняется за $\mathrm{O}(s)$). К сожалению, из-за необходимости расчета (пересчета) $x^N$, как правило, всегда приходится тратить на одну итерацию как минимум $\mathrm{O}(n)$ арифметических операций [44]. Оказывается, существует другой (более простой) способ восстановления решения прямой задачи, исходя из решения двойственной, который лучше приспособлен к возможности учета разреженности. К описанию этого способа мы сейчас и переходим.

Распространим приведенный в разделе подход на случай, когда двойственный функционал в задаче (3.2.11), в свою очередь, оказывается $\mu$-сильно выпуклым (вогнутым) в 2-норме. Для этого достаточно, чтобы функционал прямой задачи имел равномерно ограниченную константу Липшица градиента, а сама прямая задача решалась бы на всем пространстве [44, 138] (т.е. нет других ограничений, кроме ограничения $Ax = b$). Такой пример нам уже встречался в замечании 3.2.10. Воспользуемся техникой рестартов (см., например, [48, 92, 135], а также раздел 2.2 главы 2), которая в данном случае будет иметь следующий вид (см. теорему 3.2.2; заметим, что в этой оценке мы не использовали прямо-двойственность БГМ)

$$\frac{\mu}{2}\left\|\tilde{y}^{\bar{N}} - y_*\right\|_2^2 \le F\left(\tilde{y}^{\bar{N}}\right) - F\left(y_*\right) \le \frac{2L\left\|y^0 - y_*\right\|_2^2}{\bar{N}^2}.$$

Выбирая

$$\bar{N} = \sqrt{\frac{8L}{\mu}},$$

получим, что

$$\left\|\tilde{y}^{\bar{N}} - y_*\right\|_2^2 \le \frac{1}{2}\left\|y^0 - y_*\right\|_2^2.$$

Выберем в БГМ в качестве точке старта $\tilde{y}^{\bar{N}}$, и снова сделаем $\bar{N}$ итераций, и т.д. Несложно понять, что если мы хотим достичь точности по функции $\varepsilon$, то число таких рестартов (перезапусков) БГМ достаточно взять (здесь используется достаточное стандартное обозначение $\lceil \cdot \rceil$, которое мы поясним примером $\lceil 0.2 \rceil = 1$)

$$\left\lceil \log_2\left(\frac{\mu R^2}{\varepsilon}\right)\right\rceil.$$

Таким образом, общее число итераций, которое сделает БГМ, можно оценить, как

$$\sqrt{\frac{8L}{\mu}}\left\lceil \log_2\left(\frac{\mu R^2}{\varepsilon}\right)\right\rceil.$$

На первый взгляд, кажется, что при таком подходе мы не контролируем

$$\left\|Ax(y) - b\right\|_2 = \left\|\nabla F(y)\right\|_2.$$

В действительности, мы всегда (не только для прямо-двойственных методов) для задач с Липшицевым градиентом можем контролировать $\left\|\nabla F(y)\right\|_2$, используя неравенство (3.2.5) (если $F(y_*) \ne 0$, то вместо градиента $\nabla F(y)$ можно вычислять градиентное отображение, см., например, [92, 181])

$$\frac{1}{2L}\left\|\nabla F(y)\right\|_2^2 \le F(y) - F(y_*).$$

Проблема, однако, в том, что типично это неравенство довольно грубое. Действительно, в $\mu$-сильно выпуклом случае

$$\frac{1}{2L}\left\|\nabla F(y)\right\|_2^2 \le F(y)-F(y_*) \le \frac{1}{2\mu}\left\|\nabla F(y)\right\|_2^2 .$$

Если $L/\mu \gg 1$ (что типично), то использование неравенства (3.2.5) может приводить (и приводит [44], см. также раздел 5.1 главы 5) к сильно завышенным оценкам. Однако, если наряду с использованием оценки (3.2.5), учитывать геометрическую скорость сходимости БГМ в виду наличия сильной выпуклости, то (считаем, что $\mu R^2 \ge \varepsilon$)

$$\left\|Ax\left(y^N\right)-b\right\|_2 = \left\|\nabla F\left(y^N\right)\right\|_2 \le \sqrt{2L\cdot\left(F\left(y^N\right)-F\left(y_*\right)\right)} \le \sqrt{2L\mu R^2}\exp\left(-\frac{N}{2}\sqrt{\frac{\mu}{8L}}\right),$$

где $\left\{y^N\right\}$ обозначает последовательность, которую генерирует описанный выше БГМ с рестартами.

Перейдем к описанию критерия останова метода БГМ с рестартами. По определению $x(y)$ имеем

$$g\left(x(y)\right)+\left\langle y, Ax(y)-b\right\rangle \le g\left(x_*\right).$$

Откуда

$$g\left(x(y)\right)-g\left(x_*\right) \le \left\|y\right\|_2\left\|Ax(y)-b\right\|_2 .$$

Таким образом, критерий останова будет иметь вид

$$\left\|y^N\right\|_2\left\|Ax\left(y^N\right)-b\right\|_2 \le \varepsilon,\ \left\|Ax\left(y^N\right)-b\right\|_2 \le \tilde{\varepsilon} . \qquad (3.2.15)$$

**Теорема 3.2.3.** *Пусть нужно решить задачу (3.2.11) посредством перехода к задаче (3.2.12) с* $\mu$ *-сильно выпуклым в 2-норме функционалом, исходя из выписанных выше формул. Выбираем в качестве критерия останова (3.2.15). Тогда описанный выше БГМ с рестартами гарантированно остановится, сделав не более чем*

$$\max\left\{\sqrt{\frac{8L}{\mu}}\left\lceil \log_2\left(\frac{2L\mu R^4}{\varepsilon^2}\right)\right\rceil, \sqrt{\frac{8L}{\mu}}\left\lceil \log_2\left(\frac{2L\mu R^2}{\tilde{\varepsilon}^2}\right)\right\rceil\right\}$$

*итераций.*

**Замечание 3.2.12.** В действительности, только что мы описали довольно общий способ решения большого числа задач, путем перехода к гладкой двойственной. Все что нужно добавить к описанному выше подходу – это искусственное введение регуляризации $\mu \approx \varepsilon/R^2$ в двойственную задачу, когда она не сильно выпуклая [48]. При этом возникают рестарты по параметру $\mu$, поскольку $R$ априорно не известно [37, 181]. Собственно, в [37] именно таким образом и предлагалось решать задачу ЭЛП (см. также раздел 3.1 этой главы 3). Теоретические оценки говорят, что оба метода должны работать приблизительно в одно время (метод из данного раздела чуть получше), но численные эксперименты (проведенные А.В. Черновым [171]) совершенно однозначно показали, что метод, который изложен в этом разделе, работает на несколько порядков быстрее, чем метод из [37]. Связано это с тем, что из-за рестартов по $\mu$, мы обязаны сделать (на каждом шаге рестарта) предписанное число итераций метода (меньше нельзя), в то время как предложенный в данном разделе метод, во-первых, не требует рестартов (это сразу экономит почти один порядок), а во-вторых, останавливается по более гибкому критерию (см. теорему 3.2.2), допускающему, что может быть сделано меньше итераций, чем предписано полученной оценкой. Численные эксперименты показывают, что именно в этом месте и происходит основная экономия нескольких порядков в оценке общего времени работы нового метода.

**Замечание 3.2.13.** Все, что описано выше для сильно выпуклого случая, переносится на произвольные методы (прямо-двойственность при этом не нужна). Например, на метод сопряженных градиентов или метод Ньютона и их модификации [99] (есть общий тезис принадлежащий А.С. Немировскому, что практически любой разумный численный метод либо уже является прямо-двойственным, либо имеет соответствующую модификацию – для многих важных методов их прямо-двойственность уже установлена в различных работах, однако, насколько нам известно, для метода сопряженных градиентов и метода Ньютона этого пока не было сделано). В частности, критерий останова (3.2.15) является общим способом контроля точности решения, полученным произвольным методом, с помощью которого хотим одновременно решить прямую и двойственную задачу. Как уже отмечалось выше в случае, когда двойственная задача только гладкая (это необходимо для справедливости проведенных рассуждений), но не сильно выпуклая, то имеющиеся сейчас теоретические способы оценивания скорости сходимости с критерием останова (3.2.15) типично приводят к сильно завышенным оценкам. Это не означает, что соответствующие методы будут плохими, просто сейчас, насколько нам известно, не существует точных способов оценивания. Проблема наличия теоретического обоснования решается путем регуляризации двойственной задачи (см. замечание 3.2.12).

**Замечание 3.2.14.** Описанный выше подход получения БГМ для сильно выпуклых задач, базирующийся на рестартах по расстоянию от текущей точки до решения, имеет один существенный недостаток. На каждом рестарте необходимо сделать предписанное число итераций, которое, как правило, оказывается завышенным. Выйти из цикла раньше возможно, если есть критерий останова. Но $F(y_*)$, как правило, не известно. На практике можно попробовать контролировать малость нормы градиента (в общем случае нормы градиентного отображения), и исходить из уменьшения квадрата этой нормы в два раза, однако для такого подхода пока не доказано, что сохраняются аналогичные (неулучшаемые) по порядку оценки итоговой скорости сходимости. Из этой ситуации есть выход – использовать вместо БГМ с рестартами БГМ без рестартов для сильно выпуклых задач [92]. А еще лучше БГМ без рестартов для сильно выпуклых задач с адаптивным подбором константы Липшица градиента. Такой метод (к тому же непрерывный по параметру сильной выпуклости) описан, например, в работе [53] (см. также раздел 2.2 главы 2). Если параметр сильной выпуклости неизвестен (см. замечание 3.2.12), то, к сожалению, все равно не удается уйти от рестартов, но рестарты уже будут только по этому параметру, и выход из последнего рестарта (самого дорого, вносящего основной вклад в оценку общего времени работы метода) может быть осуществлен (в отличие от описанного выше подхода) по контролю малости нормы градиента (градиентного отображения). Другой способ усовершенствования конструкции рестартов изложен в работе [287].

## 3.3 Параллелизуемые двойственные методы поиска равновесий в смешанных моделях распределения потоков в больших транспортных сетях

Данный подраздел представляет собой развитие смешанной модели подраздела 1.1.7 раздела 1.1 главы 1 и прямо-двойственного метода подраздела 1.5.3 раздела 1.5 главы 1.

Напомним вкратце постановку задачи. Пусть транспортная сеть города представлена ориентированным графом $\Gamma = (V, E)$, где $V$ – узлы сети (вершины), $E \subset V \times V$ – дуги сети (рёбра графа), $O \subseteq V$ – источники корреспонденций ($S = |O|$), $D \subseteq V$ – стоки. В современных моделях равновесного распределения потоков в крупном мегаполисе число узлов графа транспортной сети обычно выбирают порядка $n = |V| \sim 10^4$ [40]. Число ребер $|E|$ получается в три четыре раза больше. Пусть $W \subseteq \{w = (i, j) : i \in O, j \in D\}$ – множество корреспонденций, т.е. возможных пар «источник» – «сток»; $p = \{v_1, v_2, ..., v_m\}$ – путь из $v_1$ в $v_m$, если $(v_k, v_{k+1}) \in E$, $k = 1, ..., m-1,\ m > 1$; $P_w$ – множество путей, отвечающих корреспонденции $w \in W$, то есть если $w = (i, j)$, то $P_w$ – множество путей, начинающихся в вершине $i$ и заканчивающихся в $j$; $P = \bigcup_{w \in W} P_w$ – совокупность всех путей в сети $\Gamma$; $x_p$ [автомобилей/час] – величина потока по пути $p$, $x = \{x_p :\ p \in P\}$; $f_e$ [автомобилей/час] – величина потока по дуге $e$:

$$f_e(x) = \sum_{p \in P} \delta_{ep} x_p \ \ (f = \Theta x), \text{ где } \delta_{ep} = \begin{cases} 1, & e \in p \\ 0, & e \notin p \end{cases};$$

$\tau_e(f_e)$ – удельные затраты на проезд по дуге $e$. Как правило, предполагают, что это – (строго) возрастающие, гладкие функции от $f_e$. Точнее говоря, под $\tau_e(f_e)$ правильнее понимать представление пользователей транспортной сети об оценке собственных затрат (обычно временных в случае личного транспорта и комфортности пути (с учетом времени в пути) в случае общественного транспорта) при прохождении дуги $e$, если поток желающих проехать по этой дуге будет $f_e$.

Рассмотрим теперь $G_p(x)$ – затраты временные или финансовые на проезд по пути $p$. Естественно считать, что $G_p(x) = \sum_{e \in E} \tau_e(f_e(x)) \delta_{ep}$.

Пусть также известно, сколько перемещений в единицу времени $d_w$ осуществляется согласно корреспонденции $w \in W$. Тогда вектор $x$, характеризующий распределение потоков, должен лежать в допустимом множестве:

$$X = \left\{ x \geq 0 : \sum_{p \in P_w} x_p = d_w,\ w \in W \right\}.$$

Рассмотрим игру, в которой каждой корреспонденции $w \in W$ соответствует свой, достаточно большой, набор однотипных “игроков”, осуществляющих передвижение согласно корреспонденции $w$ (относительный масштаб характеризуется числами $d_w$) Чистыми стратегиями игрока служат пути, а выигрышем – величина $-G_p(x)$. Игрок “выбирает” путь следования $p \in P_w$, при этом, делая выбор, он пренебрегает тем, что от его выбора также “немного” зависят $|P_w|$ компонент вектора $x$ и, следовательно, сам выигрыш $-G_p(x)$. Можно показать (см., например, [290]), что поиск равновесия Нэша–Вардропа $x^* \in X$ (макро описание равновесия) равносилен решению задачи

$$\Psi(f) = \sum_{e \in E} \sigma_e(f_e) = \sum_{e \in E} \int_0^{f_e} \tau_e(z)\,dz \to \min_{\substack{f = \Theta x \\ x \in X}} . \qquad (3.3.1)$$

В пределе модели стабильной динамики (Нестеров–деПальма) [47, 275] на части ребер $E' \subseteq E$ может быть сделан следующий предельный переход

$$\tau_e^{\mu}(f_e) \xrightarrow[\mu \to 0+]{} \begin{cases} \bar{t}_e, & 0 \leq f_e < \bar{f}_e \\ [\bar{t}_e, \infty), & f_e = \bar{f}_e \end{cases}, \quad d\tau_e^{\mu}(f_e)/df_e \xrightarrow[\mu \to 0+]{} 0,\ 0 \leq f_e < \bar{f}_e.$$

Задача (3.3.1) примет вид

$$\Psi(f) = \sum_{e \in E} \sigma_e(f_e) \to \min_{\substack{f = \Theta x,\ x \in X \\ f_e \leq \bar{f}_e,\ e \in E'}} ,\ \sigma_e(f_e) = f_e \bar{t}_e,\ \ e \in E'.$$

В данном разделе предлагается хорошо параллелизуемый двойственный численный метод поиска равновесия в смешанной модели (3.3.1), т.е. в модели, в которой часть ребер являются Бэкмановскими, а часть – Нестерова–деПальмы. Такие задачи возникают, например, в однопродуктовой модели грузоперевозок РЖД [21]. К таким “смешанным” моделям классический метод условного градиента (Франк–Вульфа), который зашит практически во все современные пакеты транспортного моделирования, к сожалению, уже не применим. Требуется разработка новых подходов.

Для задачи (3.3.1) можно построить следующую двойственную задачу [40, 47] (см. также разделы 1.1, 1.5 главы 1)

$$\Upsilon(t) = \underbrace{-\sum_{w \in W} d_w T_w(t)}_{F(t)} + \sum_{e \in E} \sigma_e^*(t_e) \to \min_{\substack{t_e \geq \bar{t}_e,\ e \in E' \\ t_e \in \operatorname{dom} \sigma_e^*(t_e),\ e \in E \setminus E'}} , \qquad (3.3.2)$$

где $T_w(t)=\min_{p\in P_w}\sum_{e\in E}\delta_{ep}t_e$ – длина кратчайшего пути из $i$ в $j$ ($w=(i,j)\in W$) на графе $\Gamma$, ребра которого взвешены вектором $t=\{t_e\}_{e\in E}$. При этом решение задачи (3.3.1) $f$ можно получить из формул: $f_e=\bar{f}_e-s_e$, $e\in E'$, где $s_e\ge 0$ – множитель Лагранжа к ограничению $t_e\ge \bar{t}_e$; $\tau_e(f_e)=t_e$, $e\in E\setminus E'$. Заметим, что для ребер $e\in E'$ имеем $\sigma_e^*(t_e)=\bar{f}_e\cdot(t_e-\bar{t}_e)$, а для BPR-функций

$$\tau_e(f_e)=\bar{t}_e\cdot\left(1+\gamma\cdot\left(\frac{f_e}{\bar{f}_e}\right)^{\frac{1}{\mu}}\right)\Rightarrow \sigma_e^*(t_e)=\bar{f}_e\cdot\left(\frac{t_e-\bar{t}_e}{\bar{t}_e\cdot\gamma}\right)^{\mu}\frac{(t_e-\bar{t}_e)}{1+\mu}.$$

В приложениях обычно выбирают $\mu=1/4$. В этом случае приводимый ниже шаг итерационного метода (3.3.3) может быть осуществлен по явным формулам, поскольку существуют квадратурные формулы (формулы Кардано–Декарта–Эйлера–Феррари [82]) для уравнений 4-й степени.

Поиск вектора $t$ представляет самостоятельный интерес, поскольку этот вектор описывает затраты на ребрах графа транспортной сети. Решение задачи (3.3.2) дает вектор затрат $t$ в равновесии.

Для решения двойственной задачи (3.3.2) воспользуемся методом зеркального спуска в композитном варианте[5] [191, 259] ($k=0,...,N$, $t^0=\bar{t}$, ограничение $t_e\in\operatorname{dom}\sigma_e^*(t_e)$, $e\in E\setminus E'$ всегда будет не активным, т.е. его можно не учитывать)

$$t^{k+1}=\arg\min_{\substack{t_e\ge\bar{t}_e,\, e\in E'\\ t_e\in\operatorname{dom}\sigma_e^*(t_e),\, e\in E\setminus E'}}\left\{\gamma_k\left\{\left\langle\partial F(t^k),t-t^k\right\rangle+\sum_{e\in E}\sigma_e^*(t_e)\right\}+\frac{1}{2}\left\|t-t^k\right\|_2^2\right\},\qquad(3.3.3)$$

где $\partial F(t^k)$ – произвольный элемент субдифференциала выпуклой функции $F(t^k)$ в точке $t^k$, а $\gamma_k=\varepsilon/M_k^2$, $M_k=\left\|\partial F(t^k)\right\|_2$, где $\varepsilon>0$ – желаемая точность решения задач (3.3.1) и (3.3.2), см. (3.3.6).

Положим (следует сравнить этот подход с подходом из раздела 1.6 главы 1)

$$\bar{t}^N=\frac{1}{S_N}\sum_{k=0}^{N}\gamma_k t^k,\; S_N=\sum_{k=0}^{N}\gamma_k,$$

[5] К аналогичному методу в связи с решением задачи РЖД [21] пришли С.В. Шпирко и А.С. Бондаренко. Описание предложенного Шпирко–Бондаренко метода вошло в дипломную работу А.С. Бондаренко [15]. См. также работу http://www.jmlr.org/papers/volume11/xiao10a/xiao10a.pdf , в которой рассматривается композитный вариант близкого метода двойственных усреднений.

$$f_e^k \in -\partial_e F\left(t^k\right),\quad \bar{f}_e^N = \frac{1}{S_N}\sum_{k=0}^{N}\gamma_k f_e^k,\ e \in E \setminus E', \tag{3.3.4}$$

$$\bar{f}_e^N = \bar{f}_e - s_e^N,\ e \in E', \tag{3.3.4'}$$

где $s_e^N$ – есть множитель Лагранжа к ограничению $t_e \ge \bar{t}_e$ в задаче

$$\frac{1}{S_N}\left\{\sum_{k=0}^{N}\gamma_k \cdot \left\{\sum_{e\in E'}\partial_e F\left(t^k\right)\cdot\left(t_e - t_e^k\right)\right\} + S_N\sum_{e\in E'}\bar{f}_e\cdot\left(t_e - \bar{t}_e\right) + \frac{1}{2}\sum_{e\in E'}\left(t_e - \bar{t}_e\right)^2\right\} \to \min_{t_e\ge\bar{t}_e,\, e\in E'}.$$

Критерий останова метода (правое неравенство)

$$(0 \le)\,\Upsilon\left(\bar{t}^N\right) + \Psi\left(\bar{f}^N\right) \le \varepsilon. \tag{3.3.5}$$

В данном разделе (в основном следуя работе [40], разделу 1.5 главы 1, см. также выступление на семинаре “Математическое моделирование транспортных потоков” в МЦНМО 14 апреля 2012 г. Ю.Е. Нестерова “Модели равновесных транспортных потоков и алгоритмы их нахождения” (видео доступно на сайте http://www.mathnet.ru/) – в этих источниках рассматривался случай $E = E'$), получен следующий результат о сходимости метода (3.3.3).

**Теорема 3.3.1.** *Пусть*

$$\tilde{M}_N^2 = \left(\frac{1}{N+1}\sum_{k=0}^{N}M_k^{-2}\right)^{-1},$$

$$R_N^2 := \frac{1}{2}\sum_{e\in E\setminus E'}\left(\tau_e\left(\bar{f}_e^N\right) - \bar{t}_e\right)^2 + \frac{1}{2}\sum_{e\in E'}\left(\tilde{t}_e^N - \bar{t}_e\right)^2,$$

$$\left\{\tilde{t}_e^N\right\}_{e\in E'} = \arg\min_{\{t_e\}_{e\in E'}\ge 0}\left\{\underbrace{-\sum_{w\in W}d_w T_w\left(\left\{\tau_e\left(\bar{f}_e^N\right)\right\}_{e\in E\setminus E'},\left\{t_e\right\}_{e\in E'}\right)}_{F\left(\left\{\tau_e\left(\bar{f}_e^N\right)\right\}_{e\in E\setminus E'},\left\{t_e\right\}_{e\in E'}\right)} + \sum_{e\in E'}\bar{f}_e^N\cdot\left(t_e - \bar{t}_e\right)\right\}.$$

*Тогда при* $N \ge 2\tilde{M}_N^2 R_N^2/\varepsilon^2$ *имеет место неравенство (3.3.5) и, как следствие,*

$$0 \le \Upsilon\left(\bar{t}^N\right) - \Upsilon_* \le \varepsilon,\ 0 \le \Psi\left(\bar{f}^N\right) - \Psi_* \le \varepsilon. \tag{3.3.6}$$

**Доказательство.** Из работ [191, 269], имеем (также в выкладках используется то, что $d\sigma_e^*\left(t_e\right)/dt_e = f_e \Leftrightarrow t_e = \tau_e\left(f_e\right),\ e \in E \setminus E';\ \sigma_e^*\left(t_e\right) \ge 0,\ \sigma_e^*\left(t_e^0\right) = \sigma_e^*\left(\bar{t}_e\right) = 0,\ e \in E$)

$$\Upsilon\left(\bar{t}^N\right) \le \frac{1}{S_N}\sum_{k=0}^{N}\gamma_k \Upsilon\left(t^k\right) \le \frac{1}{2S_N}\sum_{k=0}^{N}\gamma_k^2 M_k^2 +$$

$$+\frac{1}{S_N}\min_{t_e\ge\bar{t}_e,\, e\in E'}\left\{\sum_{k=0}^{N}\gamma_k\left\{F\left(t^k\right) + \left\langle \partial F\left(t^k\right), t - t^k\right\rangle\right\} + S_N\sum_{e\in E}\sigma_e^*\left(t_e\right) + \frac{1}{2}\left\|t - t^0\right\|_2^2\right\} \le \frac{\varepsilon}{2} +$$

$$+\min_t\left\{\frac{1}{S_N}\left\{\sum_{k=0}^{N}\gamma_k\left\{F\left(t^k\right) + \left\langle \partial F\left(t^k\right), t - t^k\right\rangle\right\} + S_N\sum_{e\in E}\sigma_e^*\left(t_e\right) + \frac{1}{2}\left\|t - \bar{t}\right\|_2^2\right\} + \sum_{e\in E'}s_e^N\cdot\left(\bar{t}_e - t_e\right)\right\} \le$$

$$\le \frac{\varepsilon}{2} + \min_{\{t_e\}_{e\in E'}} \min_{\{t_e\}_{e\in E\setminus E'}} \frac{1}{S_N}\left\{\sum_{k=0}^{N} \gamma_k \left\{F\left(t^k\right) + \left\langle \partial F\left(t^k\right), t - t^k\right\rangle\right\} + \right.$$

$$\left. + \left(S_N \sum_{e\in E\setminus E'} \sigma_e^*\left(t_e\right) + \frac{1}{2} \sum_{e\in E\setminus E'} \left(t_e - \bar t_e\right)^2\right) + \left(S_N \sum_{e\in E'} \bar f_e^N \cdot \left(t_e - \bar t_e\right) + \frac{1}{2} \sum_{e\in E'} \left(t_e - \bar t_e\right)^2\right)\right\} \le$$

$$\le \frac{\varepsilon}{2} + \min_{\{t_e\}_{e\in E'}\ge 0} \frac{1}{S_N}\left\{ \sum_{k=0}^{N} \gamma_k \underbrace{\overbrace{\left\{F\left(t^k\right) + \left\langle \partial F\left(t^k\right), \breve t^N - t^k\right\rangle\right\}}^{\breve t^N = \left(\left\{\tau_e\left(\bar f_e^N\right)\right\}_{e\in E\setminus E'}, \{t_e\}_{e\in E'}\right)}}_{\le F\left(\left\{\tau_e\left(\bar f_e^N\right)\right\}_{e\in E\setminus E'}, \{t_e\}_{e\in E'}\right)} + \right.$$

$$\left. + \left(S_N \sum_{e\in E\setminus E'} \sigma_e^*\left(\tau_e\left(\bar f_e^N\right)\right) + \frac{1}{2} \sum_{e\in E\setminus E'} \left(\tau_e\left(\bar f_e^N\right) - \bar t_e\right)^2\right) + \left(S_N \sum_{e\in E'} \bar f_e^N \cdot \left(t_e - \bar t_e\right) + \frac{1}{2} \sum_{e\in E'} \left(t_e - \bar t_e\right)^2\right)\right\} \le$$

$$\le \frac{\varepsilon}{2} + \min_{\{t_e\}_{e\in E'}\ge 0} \left\{F\left(\left\{\tau_e\left(\bar f_e^N\right)\right\}_{e\in E\setminus E'}, \left\{t_e\right\}_{e\in E'}\right) + \right.$$

$$+ \sum_{e\in E'} \bar f_e^N \cdot \left(t_e - \bar t_e\right) + \frac{1}{2S_N} \sum_{e\in E'} \left(t_e - \bar t_e\right)^2 \Bigg\} +$$

$$+ \sum_{e\in E\setminus E'} \sigma_e^*\left(\tau_e\left(\bar f_e^N\right)\right) + \frac{1}{2S_N} \sum_{e\in E\setminus E'} \left(\tau_e\left(\bar f_e^N\right) - \bar t_e\right)^2 \le$$

$$\le \frac{\varepsilon}{2} + \min_{\{t_e\}_{e\in E'}\ge 0} \left\{F\left(\left\{\tau_e\left(\bar f_e^N\right)\right\}_{e\in E\setminus E'}, \left\{t_e\right\}_{e\in E'}\right) + \sum_{e\in E\setminus E'} \sigma_e^*\left(\tau_e\left(\bar f_e^N\right)\right) + \sum_{e\in E'} \bar f_e^N \cdot \left(t_e - \bar t_e\right)\right\} +$$

$$+ \frac{1}{2S_N} \sum_{e\in E'} \left(\tilde t_e^N - \bar t_e\right)^2 + \frac{1}{2S_N} \sum_{e\in E\setminus E'} \left(\tau_e\left(\bar f_e^N\right) - \bar t_e\right)^2 =$$

$$= \frac{\varepsilon}{2} - \Psi\left(\bar f^N\right) + \frac{R_N^2}{S_N} = \frac{\varepsilon}{2} + \frac{\tilde M_N^2 R_N^2}{\varepsilon \cdot (N+1)} - \Psi\left(\bar f^N\right) \le \varepsilon - \Psi\left(\bar f^N\right). \ \square$$

**Следствие 3.3.1.** *Пусть* $t^*$ *– решение задачи (3.3.2). Положим*

$$R^2 = \frac{1}{2}\left\|t^* - t^0\right\|_2^2 = \frac{1}{2}\left\|t_* - \bar t\right\|_2^2 .$$

*Тогда при* $N = 2\tilde M_N^2 R^2 / \varepsilon^2$ *справедливы неравенства*

$$\frac{1}{2}\left\|t - t^{k+1}\right\|_2^2 \le 2R^2, \ k = 0, \ldots, N \tag{3.3.7}$$

$$0 \le \Upsilon\left(\bar t^N\right) - \Upsilon_* \le \varepsilon . \tag{3.3.8}$$

**Доказательство.** Формула (3.3.8) – стандартный результат, см., например, [259]. Формула (3.3.7) также является достаточно стандартной [269] (см. также раздел 3.1 этой главы 3), однако далее приводится схема ее вывода. Из доказательства теоремы 3.3.1 имеем для любого $k = 0, \ldots, N$

$$0 \le \frac{1}{2}\sum_{l=0}^{k}\gamma_l^2 M_l^2 + \frac{1}{2}\left\|t - t^0\right\|_2^2 - \frac{1}{2}\left\|t - t^{k+1}\right\|_2^2 .$$

Отсюда следует, что

$$\frac{1}{2}\left\|t - t^{k+1}\right\|_2^2 \le R^2 + \frac{1}{2}\sum_{l=0}^{k}\varepsilon^2 M_l^{-2} \le R^2 + \frac{1}{2}\sum_{k=0}^{N}\varepsilon^2 M_k^{-2} = 2R^2 . \ \square$$

**Замечание 3.3.1.** Преимуществом подхода (3.3.3), (3.3.4), (3.3.4’) над подходом п. 3 работы [40] (см. также раздел 1.5 главы 1) является простота описания (отсутствие необходимости делать рестарты по неизвестным параметрам) и наличие эффективно проверяемого критерия останова (3.3.5). К недостаткам стоит отнести вхождение в оценку скорости сходимости плохо контролируемого $R_N^2$, которое может оказаться большим даже в случае, когда $E = E'$. Далее мы опишем другой способ (см. также работы [6, 263], в которых описаны близкие конструкции) восстановления решения прямой задачи (3.3.1), отличный от (3.3.4), (3.3.4’), в части (3.3.4’), который в случае $E = E'$ позволяет использовать $R^2$ вместо $R_N^2$.

Положим (следует сравнить этот подход с подходом из предыдущего раздела)

$$f_e^k \in -\partial_e F\left(t^k\right),\quad \bar{f}_e^N = \frac{1}{S_N}\sum_{k=0}^{N}\gamma_k f_e^k ,\ e \in E . \tag{3.3.9}$$

$$\tilde{R}^2 = \frac{1}{2}\sum_{e\in E'}\left(t_e^* - t_e^0\right)^2 = \frac{1}{2}\sum_{e\in E'}\left(t_e^* - \bar{t}_e\right)^2 .$$

**Теорема 3.3.2.** *Пусть*

$$\tilde{R}_N^2 := \frac{1}{2}\sum_{e\in E\setminus E'}\left(\tau_e\left(\bar{f}_e^N\right) - \bar{t}_e\right)^2 + 5\tilde{R}^2 . \tag{3.3.10}$$

*Тогда при* $N \ge 4\tilde{M}_N^2\tilde{R}_N^2/\varepsilon^2$ *имеют место неравенства*

$$\left|\Upsilon\left(\bar{t}^N\right) - \Upsilon_*\right| \le \varepsilon ,\ \left|\Psi\left(\bar{f}^N\right) - \Psi_*\right| \le \varepsilon .$$

*Более того, также имеют место неравенства*

$$\sqrt{\sum_{e\in E'}\left(\left(\bar{f}_e^N - \bar{f}_e\right)_+\right)^2} \le \tilde{\varepsilon} ,\ \tilde{\varepsilon} = \varepsilon/\tilde{R} , \tag{3.3.11}$$

$$\Psi\left(\bar{f}^N\right) - \Psi_* \le \Upsilon\left(\bar{t}^N\right) + \Psi\left(\bar{f}^N\right) \le \varepsilon ,$$

*которые можно использовать для критерия останова метода (задавшись парой* $(\varepsilon, \tilde{\varepsilon})$*).*

**Доказательство.** Рассуждая аналогично доказательству теоремы 3.3.1 (см. также [6, 263]), получим

$$\Upsilon\left(\bar{t}^N\right) \le \frac{1}{S_N}\sum_{k=0}^{N}\gamma_k \Upsilon\left(t^k\right) \overset{(3.3.12)}{\le} \frac{1}{2S_N}\sum_{k=0}^{N}\gamma_k^2 M_k^2 +$$

$$+\frac{1}{S_N}\min_{\substack{t_e,\,e\in E\setminus E';\,t_e\ge\bar t_e,\,e\in E'\\ \frac{1}{2}\sum_{e\in E'}(t_e-\bar t_e)^2\le 5\tilde R^2}}\left\{\sum_{k=0}^{N}\gamma_k\left\{F\left(t^k\right)+\left\langle\partial F\left(t^k\right),t-t^k\right\rangle\right\}+S_N\sum_{e\in E}\sigma_e^*\left(t_e\right)\right\}+\frac{\tilde R_N^2}{S_N}\overset{(3.3.13)}{\le}$$

$$\overset{(3.3.13)}{\le}\frac{\varepsilon}{2}-\Psi\left(\bar f^N\right)-\max_{\substack{t_e\ge\bar t_e,\,e\in E'\\ \frac{1}{2}\sum_{e\in E'}(t_e-\bar t_e)^2\le 5\tilde R^2}}\left\{\frac{1}{S_N}\sum_{k=0}^{N}\gamma_k\sum_{e\in E'}\left(f_e^k-\bar f_e\right)\left(t_e-\bar t_e\right)\right\}+\frac{\tilde R_N^2}{S_N}\le$$

$$\le\frac{\varepsilon}{2}-\Psi\left(\bar f^N\right)-3\tilde R\sqrt{\sum_{e\in E'}\left(\left(\bar f_e^N-\bar f_e\right)_+\right)^2}+\frac{\tilde R_N^2}{S_N}\le\varepsilon-\Psi\left(\bar f^N\right)-3\tilde R\sqrt{\sum_{e\in E'}\left(\left(\bar f_e^N-\bar f_e\right)_+\right)^2}.$$

Неравенство (3.3.12) было получено на базе следующего соотношения:

$$\left\{\tau_e\left(\bar f_e^N\right)\right\}_{e\in E\setminus E'}=\arg\min_{t_e,\,e\in E\setminus E'}\left\{\frac{1}{S_N}\sum_{k=0}^{N}\gamma_k\left\{F\left(t^k\right)+\left\langle\partial F\left(t^k\right),t-t^k\right\rangle\right\}+\sum_{e\in E}\sigma_e^*\left(t_e\right)\right\}.$$

Неравенство (3.3.13) было получено на базе следующих соотношений:

$$\partial F\left(t^k\right)=-f^k,\ F\left(t^k\right)=-\left\langle f^k,t^k\right\rangle,\ \min_{t_e}\left\{-f_e^k t_e+\sigma_e^*\left(t_e\right)\right\}=-\sigma_e\left(f_e^k\right),\ e\in E\setminus E',$$

$$-f_e^k t_e+\sigma_e^*\left(t_e\right)=-\left(f_e^k-\bar f_e\right)\left(t_e-\bar t_e\right)-f_e^k\bar t_e=-\left(f_e^k-\bar f_e\right)\left(t_e-\bar t_e\right)-\sigma_e\left(f_e^k\right),\ e\in E',$$

$$-\frac{1}{S_N}\sum_{k=0}^{N}\gamma_k\sum_{e\in E}\sigma_e\left(f_e^k\right)\le-\Psi\left(\bar f^N\right).$$

Таким образом,

$$\Upsilon\left(\bar t^N\right)+\Psi\left(\bar f^N\right)+3\tilde R\sqrt{\sum_{e\in E'}\left(\left(\bar f_e^N-\bar f_e\right)_+\right)^2}\le\varepsilon.$$

Повторяя рассуждения п. 6.11 [261] и п. 3 [6], получим искомые неравенства. □

**Замечание 3.3.2.** О преимуществе использования формулы (3.3.9) вместо (3.3.4), (3.3.4') в случае $E=E'$ написано в замечании 3.3.1. Сейчас отметим возникающие при таком подходе недостатки: 1) возможность нарушения ограничения $f_e\le\bar f_e$, $e\in E'$ в прямой задаче, 2) отсутствие левых неравенств в двойных неравенствах (3.3.5), (3.3.6).

**Замечание 3.3.3.** Формулы (3.3.4), (3.3.9) вынужденно (в случае (3.3.4)) или осознано (в случае (3.3.9) при $e\in E'$) восстанавливают решение прямой задачи исходя из "модели" – явной формулы, связывающей прямые и двойственные переменные. Наличие таких переменных неизбежно приводит к возникновению в оценках зазора двойственности трудно контролируемых размеров решений вспомогательных задач. При этом в случае, когда наличие модели сопряжено с каким-то ограничением в прямой задаче (формулы (3.3.9) при $e\in E'$, и ограничение в прямой задаче $f_e\le\bar f_e$, $e\in E'$), допускается нарушение этих ограничений, которые необходимо контролировать (3.3.11). Зато по этим переменным имеется полный контроль соответствующих этим переменным частей оценок зазора двойственности (3.3.10) (см. также предыдущий раздел). Подход (3.3.4') связанный с наличием ограничений в решаемой задаче ($t_e \geq \bar{t}_e$, $e \in E'$) также приводит к возникновению в оценках зазора двойственности трудно контролируемых размеров решений вспомогательных задач, однако уже не приводит к нарушению никаких ограничений в самой задаче и сопряженной к ней (в нашем случае исходная задача – двойственная (3.3.2), а сопряженная к ней – прямая (3.3.1)). Эти два прямо-двойственных подхода дополняются другим прямо-двойственным подходом, в котором шаги осуществляются "по функционалу", если не нарушены или слабо нарушены ограничения в задаче, и "по нарушенному ограничению" в противном случае. Подробнее об этом см., например, в работах [261, 281] и следующем разделе.

**Замечание 3.3.4.** Оба описанных подхода можно распространить, сохраняя вид формул восстановления и структуру рассуждений, на практически произвольные пары прямая / двойственная задача, поскольку выбранный нами пример пары взаимно-сопряженных задач (3.3.1), (3.3.2) и так содержал в себе практически все основные нюансы, которые могут возникать при таких рассуждениях. Более того, вместо МЗС можно было бы использовать любой другой прямо-двойственный метод. Например, композитный универсальный градиентный метод Ю.Е. Нестерова. В частности, варианты этого метода из работы [53].

**Замечание 3.3.5.** Возможность эффективного распараллеливания предложенных методов связана с возможностью эффективного вычисления самой затратной части шага описанного итерационного метода: расчет элемента субдифференциала $\partial F\left(t^k\right)$ (см. формулы (3.3.3), (3.3.4), (3.3.9) в которых этот субдифференциал используется):

$$\partial F(t) = -\sum_{i \in O} \sum_{j \in D:\,(i,j) \in W} d_{ij} \partial T_{ij}(t).$$

Вычисление $\left\{\partial T_{ij}(t)\right\}_{j \in D:\,(i,j) \in W}$ может быть осуществлено алгоритмом Дейкстры [132] (и его более современными аналогами [145]) за $\mathrm{O}(n \ln n)$. При этом под $\partial T_{ij}(t)$ можно понимать описание произвольного (если их несколько) кратчайшего путь из вершины $i$ в вершину $j$ на графе $\Gamma$, ребра которого взвешены вектором $t = \{t_e\}_{e \in E}$. Под "описанием" понимается $\left[\partial T_{ij}(t)\right]_e = 1$, если $e$ попало в кратчайший путь и $\left[\partial T_{ij}(t)\right]_e = 0$ иначе. Таким образом, вычисление $\partial F(t)$ может быть распараллелено на $S$ процессорах.

Студентами 6-го курса ФУПМ и ФАКИ МФТИ были проведены численные эксперименты с описанными алгоритмам https://github.com/vikalijko/transport. Данные для экспериментов брались с ресурса https://github.com/bstabler/TransportationNetworks.

.

## 3.4 Прямо-двойственный метод зеркального спуска для условных задач стохастической оптимизации

Рассмотрим задачу выпуклой условной оптимизации (к такому виду, например, приводится задача из подраздела 1.5.3 раздела 1.5 главы 1 и раздела 3.3 этой главы 3)

$$f(x) \to \min_{g(x)\le 0,\, x\in Q}. \qquad (3.4.1)$$

Под решением этой задачи будем понимать такой $\bar{x}^N \in Q \subseteq \mathbb{R}^n$, что с вероятностью $\ge 1-\sigma$ имеет место неравенство (определение $M_f$, $M_g$ см. ниже)

$$f\left(\bar{x}^N\right) - f_* \le \varepsilon_f = \frac{M_f}{M_g}\varepsilon_g,\ g\left(\bar{x}^N\right) \le \varepsilon_g, \qquad (3.4.2)$$

где $f_* = f(x_*)$ – оптимальное значение функционала в задаче (3.4.1), $x_*$ – решение задачи (3.4.1).

Введем норму $\|\ \|$ в прямом пространстве (сопряженную норму будем обозначать $\|\ \|_*$) и прокс-функцию $d(x)$ сильно выпуклую относительно этой нормы, с константой сильной выпуклости $\ge 1$. Выберем точку старта

$$x^1 = \arg\min_{x\in Q} d(x),$$

считаем, что

$$d\left(x^1\right) = 0,\ \nabla d\left(x^1\right) = 0.$$

Введем брэгмановское “расстояние”

$$V_x(y) = d(y) - d(x) - \langle \nabla d(x), y - x\rangle.$$

Определим “размер” решения

$$d(x_*) = V_{x^1}(x_*) = R^2$$

и “размер” множества $Q$ (для большей наглядности считаем множество ограниченным – в общем случае приводимые далее рассуждения также можно провести подобно [40], см. также раздел 1.5 главы 1)

$$\max_{x,y\in Q} V_x(y) = \bar{R}^2.$$

Будем считать, что имеется такая последовательность независимых случайных величин $\left\{\xi^k\right\}$ и последовательности $\left\{\nabla_x f\left(x,\xi^k\right)\right\}$, $\left\{\nabla_x g\left(x,\xi^k\right)\right\}$, $k = 1,\dots,N$, что имеют место следующие соотношения

$$E_{\xi^k}\left[\nabla_x f\left(x,\xi^k\right)\right]=\nabla f\left(x\right),\ E_{\xi^k}\left[\nabla_x g\left(x,\xi^k\right)\right]=\nabla g\left(x\right); \tag{3.4.3}$$

$$\left\|\nabla_x f\left(x,\xi^k\right)\right\|_*^2\le M_f^2,\ \left\|\nabla_x g\left(x,\xi^k\right)\right\|_*^2\le M_g^2 \tag{3.4.4}$$

или

$$E_{\xi^k}\left[\left\|\nabla_x f\left(x,\xi^k\right)\right\|_*^2\right]\le M_f^2,\ E_{\xi^k}\left[\left\|\nabla_x g\left(x,\xi^k\right)\right\|_*^2\right]\le M_g^2\ . \tag{3.4.4'}$$

На каждой итерации $k=1,...,N$ нам доступен стохастический (суб-)градиент $\nabla_x f\left(x,\xi^k\right)$ или $\nabla_x g\left(x,\xi^k\right)$ в одной, выбранной нами (методом), точке $x^k$.

Опишем стохастический вариант метода зеркального спуска (МЗС) для задач с функциональными ограничениями (этот метод восходит к А.С. Немировскому, 1977 [91]).

Определим оператор “проектирования” согласно этому расстоянию

$$\mathrm{Mirr}_{x^k}\left(\mathrm{v}\right)=\arg\min_{y\in Q}\left\{\left\langle \mathrm{v},y-x^k\right\rangle+V_{x^k}\left(y\right)\right\}.$$

МЗС для задачи (3.4.1) будет иметь вид (см., например, [3])

$$\boxed{\begin{array}{l} x^{k+1}=\mathrm{Mirr}_{x^k}\left(h_f\nabla_x f\left(x^k,\xi^k\right)\right),\ \text{если}\ g\left(x^k\right)\le\varepsilon_g, \\ x^{k+1}=\mathrm{Mirr}_{x^k}\left(h_g\nabla_x g\left(x^k,\xi^k\right)\right),\ \text{если}\ g\left(x^k\right)>\varepsilon_g, \end{array}} \tag{3.4.5}$$

где $h_g=\varepsilon_g/M_g^2$, $h_f=\varepsilon_g/\left(M_f M_g\right)$, $k=1,...,N$. Обозначим через $I$ множество индексов $k$, для которых $g\left(x^k\right)\le\varepsilon_g$. Введем также обозначения

$$\left[N\right]=\left\{1,...,N\right\},\ J=\left[N\right]\setminus I\ ,\ N_I=\left|I\right|,\ N_J=\left|J\right|,\ \bar{x}^N=\frac{1}{N_I}\sum_{k\in I}x^k\ .$$

В сформулированных далее теоремах, предполагается, что последовательность $\left\{x^k\right\}_{k=1}^{N+1}$ генерируется методом (3.4.5).

**Теорема 3.4.1**. *Пусть справедливы условия (3.4.3), (3.4.4’). Тогда при*

$$N\ge\frac{2M_g^2R^2}{\varepsilon_g^2}+1$$

*выполняются неравенства* $N_I\ge1$ *с вероятностью* $\ge1/2$ *и*

$$E\left[f\left(\bar{x}^N\right)\right]-f_*\le\varepsilon_f\ ,\ g\left(\bar{x}^N\right)\le\varepsilon_g\ .$$

*Пусть справедливы условия (3.4.3), (3.4.4). Тогда при*

$$N \geq \frac{81 M_g^2 \bar{R}^2}{\varepsilon_g^2} \ln\left(\frac{1}{\sigma}\right) \tag{3.4.6}$$

*с вероятностью* $\geq 1-\sigma$ *выполняются неравенства* $N_I \geq 1$ *и*

$$f\left(\bar{x}^N\right) - f_* \leq \varepsilon_f , \; g\left(\bar{x}^N\right) \leq \varepsilon_g ,$$

*т.е. выполняются неравенства (3.4.2).*

**Доказательство.** Первая часть теоремы установлена в работе [3] (см. также раздел 4.4 главы 4). Докажем вторую часть. Согласно [191] имеет место неравенство (для любого $x \in Q$, $g(x) \leq 0$)

$$h_f N_I \cdot \left( f\left(\bar{x}^N\right) - f(x) \right) \leq$$

$$\leq h_f \sum_{k \in I} \left\langle E_{\xi^k}\left[ \nabla_x f\left(x^k, \xi^k\right) \right], x^k - x \right\rangle \leq \frac{h_f^2}{2} \sum_{k \in I} \left\| \nabla_x f\left(x^k, \xi^k\right) \right\|_*^2 +$$

$$+ h_f \sum_{k \in I} \left\langle E_{\xi^k}\left[ \nabla_x f\left(x^k, \xi^k\right) \right] - \nabla_x f\left(x^k, \xi^k\right), x^k - x \right\rangle -$$

$$- h_g \sum_{k \in J} \underbrace{\left\langle E_{\xi^k}\left[ \nabla_x g\left(x^k, \xi^k\right) \right], x^k - x \right\rangle}_{\geq g\left(x^k\right) - g(x) > \varepsilon_g} + \frac{h_g^2}{2} \sum_{k \in J} \left\| \nabla_x g\left(x^k, \xi^k\right) \right\|_*^2 +$$

$$+ h_g \sum_{k \in J} \left\langle E_{\xi^k}\left[ \nabla_x g\left(x^k, \xi^k\right) \right] - \nabla_x g\left(x^k, \xi^k\right), x^k - x \right\rangle +$$

$$+ \sum_{k \in [N]} \left( V_{x^k}(x) - V_{x^{k+1}}(x) \right).$$

Положим $x = x_*$, и введем

$$\delta_N = h_f \sum_{k \in I} \left\langle \nabla f\left(x^k\right) - \nabla_x f\left(x^k, \xi^k\right), x^k - x_* \right\rangle +$$

$$+ h_g \sum_{k \in J} \left\langle \nabla g\left(x^k\right) - \nabla_x g\left(x^k, \xi^k\right), x^k - x_* \right\rangle .$$

Тогда

$$h_f N_I \cdot \left( f\left(\bar{x}^N\right) - f\left(x_*\right) \right) \leq$$

$$\leq \frac{1}{2} h_f^2 M_f^2 N_I - \frac{1}{2 M_g^2} \varepsilon_g^2 N_J + V_{x^1}\left(x_*\right) - V_{x^{N+1}}\left(x_*\right) + \delta_N =$$

$$= \frac{1}{2}\left( h_f^2 M_f^2 + \frac{\varepsilon_g^2}{M_g^2} \right) N_I - \frac{1}{2 M_g^2} \varepsilon_g^2 N + R^2 - V_{x^{N+1}}\left(x_*\right) + \delta_N =$$

$$= \varepsilon_f h_f N_I - \frac{1}{2M_g^2}\varepsilon_g^2 N + R^2 - V_{x^{N+1}}(x_*) + \delta_N \le$$

$$\le \varepsilon_f h_f N_I + \left( R^2 + \delta_N - \frac{1}{2M_g^2}\varepsilon_g^2 N \right). \qquad (3.4.7)$$

По неравенству Азума–Хефдинга [159] (см. также раздел 1.5 главы 1 и раздел 5.1 главы 5)

$$P\left(\delta_N \ge 2\sqrt{2}\bar{R}\Lambda\sqrt{h_f^2 M_f^2 N_I + h_g^2 M_g^2 N_J}\right) \le \exp\left(-\Lambda^2/2\right),$$

т.е. с вероятностью $\ge 1-\sigma$

$$P\left( \delta_N \ge \frac{4\bar{R}\varepsilon_g}{M_g}\sqrt{N \ln\left(\frac{1}{\sigma}\right)} \right) \le \sigma .$$

Будем считать, что (константу 81 можно уменьшить до $\left(4+\sqrt{18}\right)^2$)

$$N \ge \frac{81 M_g^2 \bar{R}^2}{\varepsilon_g^2} \ln\left(\frac{1}{\sigma}\right).$$

Тогда с вероятностью $\ge 1-\sigma$ выражение в скобке в формуле (3.4.7) строго меньше нуля, поэтому выполнены неравенства

$$f\left(\bar{x}^N\right) - f_* \le \varepsilon_f , \; g\left(\bar{x}^N\right) \le \varepsilon_g .$$

Последнее неравенство следует из того, что по построению $g\left(x^k\right) \le \varepsilon_g$ , $k \in I$ и из выпуклости функции $g(x)$. □

Пусть $g(x) = \max_{l=1,\dots,m} g_l(x)$. Рассмотрим двойственную задачу

$$\varphi(\lambda) = \min_{x \in Q} \left\{ f(x) + \sum_{l=1}^{m} \lambda_l g_l(x) \right\} \to \max_{\lambda \ge 0} . \qquad (3.4.8)$$

Всегда имеет место следующее неравенство (слабая двойственность)

$$0 \le f(x) - \varphi(\lambda) \overset{def}{=} \Delta(x,\lambda), \; x \in Q, \; g(x) \le 0, \; \lambda \ge 0.$$

Обозначим решение задачи (3.4.8) через $\lambda_*$. Будем считать, что выполняются условия Слейтера, т.е. существует такой $\tilde{x} \in Q$, что $g(\tilde{x}) < 0$. Тогда

$$f_* = f(x_*) = \varphi(\lambda_*) \overset{def}{=} \varphi_* .$$

В этом случае “качество” пары $\left(x^N, \lambda^N\right)$ естественно оценивать величиной зазора двойственности $\Delta\left(x^N, \lambda^N\right)$. Чем он меньше, тем лучше.

Пусть (ограничимся рассмотрением детерминированного случая)

$$g\left(x^k\right) = g_{l(k)}\left(x^k\right),\ \nabla g\left(x^k\right) = \nabla g_{l(k)}\left(x^k\right),\ k \in J .$$

Положим

$$\lambda_l^N = \frac{1}{h_f N_I} \sum_{k \in J} h_g I\left[l(k) = l\right],$$

$$I\left[\text{predicat}\right] = \begin{cases} 1, & \text{predicat} = true \\ 0, & \text{predicat} = false \end{cases} .$$

**Теорема 3.4.2**. *Пусть*

$$\left\|\nabla f(x)\right\|_*^2 \le M_f^2,\ \left\|\nabla g(x)\right\|_*^2 \le M_g^2 .$$

*Тогда при*

$$N \ge \frac{2M_g^2 \bar{R}^2}{\varepsilon_g^2} + 1 .$$

*выполняются неравенства* $N_I \ge 1$ *и*

$$\Delta\left(\bar{x}^N, \lambda^N\right) \le \varepsilon_f,\ g\left(\bar{x}^N\right) \le \varepsilon_g .$$

**Доказательство.** Согласно [191] имеет место неравенство

$$h_f N_I f\left(\bar{x}^N\right) \le$$

$$\le \min_{x \in Q}\left\{ h_f N_I f(x) + h_f \sum_{k \in I}\left\langle \nabla f\left(x^k\right), x^k - x \right\rangle \right\} \le$$

$$\le \min_{x \in Q}\left\{ h_f N_I f(x) + \frac{h_f^2}{2} \sum_{k \in I}\left\|\nabla f\left(x^k\right)\right\|_*^2 + \right.$$

$$-h_g \sum_{k \in J} \underbrace{\left\langle \nabla g\left(x^k\right), x^k - x \right\rangle}_{\ge g_{l(k)}\left(x^k\right) - g_{l(k)}(x)} + \frac{h_g^2}{2} \sum_{k \in J}\left\|\nabla g\left(x^k\right)\right\|_*^2 +$$

$$\left. + \sum_{k \in [N]}\left(V_{x^k}(x) - V_{x^{k+1}}(x)\right) \right\} \le$$

$$\le \frac{1}{2} h_f^2 M_f^2 N_I - \frac{1}{2M_g^2} \varepsilon_g^2 N_J + \bar{R}^2 +$$

$$+ h_f N_I \min_{x \in Q}\left\{ f(x) + \sum_{l=1}^m \lambda_l^N g_l(x) \right\} =$$

$$= \varepsilon_f h_f N_I + \left( \bar{R}^2 - \frac{1}{2M_g^2} \varepsilon_g^2 N \right) + h_f N_I \varphi\left(\lambda^N\right) .$$

Последующие рассуждения повторяют соответствующие рассуждения в доказательстве теоремы 3.4.1 (см. формулу (3.4.7) и следующий за ней текст). □

**Замечание 3.4.1.** Результаты теорем 3.4.1, 3.4.2 содержатся в работах [281] в детерминированном случае. Причем эти результаты были установлены для других методов (близких к (3.4.5), но все же отличных от (3.4.5)). Также как и в работе [281] основным достоинством описанного метода (3.4.5) является отсутствие в оценках скорости его сходимости размера двойственного решения, которое входит в оценки других прямо-двойственных методов / подходов (см., например, [6]).

**Замечание 3.4.2.** В работах [173, 261] предлагается другой способ получения близких результатов в детерминированном случае. В основу подхода этих работ положен метод эллипсоидов вместо МЗС (3.4.5). Отметим, что в п. 5 работы [261] показывается, как можно избавиться от нарушения ограничения $g\left(\bar{x}^N\right) \le \varepsilon_g$. К сожалению, в описанном в разделе подходе (следуя циклу работ [3, 281]) ограничение возмущалось, чтобы обеспечить должные оценки скорости убывания зазора двойственности. В работах [173, 261] использовался метод эллипсоидов, который гарантировал нужные оценки скорости сходимости сертификата точности (мажорирующего зазор двойственности) и без релаксации ограничений.

**Замечания 3.4.3.** С помощью конструкции работы [50] не сложно перенести полученные выше результаты на случай наличия малых шумов не случайной природы.

**Замечания 3.4.4 (Анастасия Баяндина).** Описанный метод (3.4.5) можно распространить на случай произвольных $\varepsilon_g$ и $\varepsilon_f$, не связанных соотношением $\varepsilon_f = M_f \varepsilon_g / M_g$. Подобно, например, J.C. Duchi http://stanford.edu/~jduchi/PCMIConvex/Duchi16.pdf можно предложить адаптивный вариант метода (3.4.5), не требующего априорного знания оценок $M_f$ и $M_g$. Метод можно распространить на задачи композитной оптимизации в случае, когда у функции и функциональных ограничений одинаковый композит. С помощью конструкции рестартов можно распространить метод на сильно выпуклые постановки задач (функционал и ограничения сильно выпуклые). Отчасти этот план был недавно независимо реализован в работе G. Lan & Z. Zhou arXiv:1604.03887. Однако, по-прежнему, открытым остается вопрос о распространении описанного подхода на задачи с неограниченным множеством $Q$.

**Пример 3.4.1 (см. также раздел 4.4 главы 4).** Основные приложения описанного подхода, это задачи вида

$$f(x) \to \min_{\max\limits_{k=1,\dots,m} \sigma_k\left(A_k^T x\right) \le 0},$$

с разреженной матрицей

$$A = \left[ A_1, \ldots, A_m \right]^T,$$

где

$$A_k = A_k^+ - A_k^-,$$

и каждый из векторов $A_k^+$, $A_k^-$ имеет не отрицательные компоненты. В частности, задачи вида

$$f(x) \to \min_{Ax \le b}$$

и приводящиеся к такому виду задачи

$$f(x) \to \min_{\substack{Ax \le b \\ Cx = d}}.$$

В этих задачах, как правило, выбирают $\| \; \| = \| \; \|_2$, $d(x) = \|x\|_2^2 / 2$. При этом в определенных ситуациях выгодно (для сокращения стоимости итерации) привносить рандомизацию [3]. Например, такую

$$\nabla_x g\left(x, \xi^k\right) = \left\| A_{k(x)}^+ \right\|_1 e_{i\left(\xi^k\right)} - \left\| A_{k(x)}^- \right\|_1 e_{j\left(\xi^k\right)},$$

где

$$k(x) \in \operatorname{Arg} \max_{k=1,\ldots,m} \sigma_k \left( A_k^T x \right),$$

причем не важно, какой именно представитель $\operatorname{Arg}\max$ выбирается;

$$e_i = \overbrace{(\underbrace{0,\ldots,0,1}_{i},0,\ldots,0)}^{n};$$

$$i\left(\xi^k\right) = i \text{ с вероятностью } A_{k(x)i}^+ \Big/ \left\| A_{k(x)}^+ \right\|_1,\; i = 1,\ldots,n;$$

$$j\left(\xi^k\right) = j \text{ с вероятностью } A_{k(x)j}^- \Big/ \left\| A_{k(x)}^- \right\|_1,\; j = 1,\ldots,n. \;\blacksquare$$

# Глава 4 Решение задач выпуклой оптимизации в пространствах сверхбольших размеров. Рандомизация и разреженность

## 4.1 Об эффективных рандомизированных алгоритмах поиска вектора PageRank

### 4.1.1 Введение

Известно, что поисковая система Google была создана в качестве учебного проекта студентов Л. Пейжда и С. Брина из Стэнфордского университета, см. [162]. В 1996 году авторы работали над поисковой системой BackRub, а в 1998 году на её основе создали новую поисковую систему Google [141, 241]. В [162] был предложен определенный способ ранжирования web-страниц. Этот способ, также как и довольно большой класс задач ранжирования, возникающих, например, при вычислении индексов цитирования ученых или журналов [202], сводится [38] к нахождению левого собственного вектора (нормированного на единицу: $\sum_{k=1}^{n} p_k = 1$), отвечающего собственному значению 1, некоторой стохастической (по строкам) матрицы $P = \left\| p_{ij} \right\|_{i,j=1}^{n,n}$: $p^T = p^T P$, $n \gg 1$.

**Замечание 4.1.1.** К поиску такого вектора, который иногда называют вектором Фробениуса–Перрона, сводится (например, в модели де Гроота) задача поиска консенсуса. Подробнее об этом, и в целом о моделях консенсуса, написано в обзоре [1].

**Замечание 4.1.2.** Решение $p^T = p^T P$ всегда существует по теореме Брауера (непрерывный (ограниченный) оператор $P$ отображает выпуклый компакт (симплекс) в себя), и единственно в классе распределений вероятностей тогда и только тогда, когда имеется всего “один класс сообщающихся (существенных) состояний”, при возможном наличии “несущественных состояний” [246]. Другими словами, если мы поставим в соответствие матрице $P$ такой ориентированный граф, что вершины $i$ и $j$ соединены ребром тогда и только тогда, когда $p_{ij} > 0$, то в таком графе любая вершина может принадлежать только одному из двух типов: несущественная – стартуя из этой вершины, двигаясь по ребрам с учетом их ориентации, всегда можно забрести в такую вершину, из которой обратно никогда не вернемся; существенная – стартуя из любой такой вершины, мы можем добраться в любую другую существенную вершину (в частности вернуться в исходную).

Приведем, следуя [38] (см. также приложение в конце диссертации), обоснование такому способу. Пусть имеется ориентированный граф $G = \langle V, E \rangle$ сети Интернет (вершины – web-страницы, ребра – ссылки: запись $(i, j) \in E$ означает, что на $i$-й странице имеется ссылка на $j$-ю страницу), $N$ – число пользователей сети (это число не меняется со вре-менем). Пусть $n_i(t)$ – число посетителей web-страницы $i$ в момент времени *t*. За один такт времени каждый посетитель этой web-страницы независимо ни от чего с вероятностью $p_{ij}$ переходит по ссылке на web-страницу *j*. Считаем стохастическую матрицу $P$ неразложимой и апериодической (см. [246]). Ниже приведен результат из [38] (обоснование будет приведено в подразделе 4.1.4), позволяющий по-другому интерпретировать вектор $p$ (PageRank), согласно которому и происходит ранжирование web-страниц:

$$\exists\ \lambda_{0.99} > 0, T_G > 0: \ \forall\ t \ge T_G$$

$$P\left(\left\|\frac{n(t)}{N} - p\right\|_2 \le \frac{\lambda_{0.99}}{\sqrt{N}}\right) \ge 0.99,$$

где $p^T = p^T P$ (решение единственно в классе распределений вероятностей в виду неразложимости).

В ряде случаев считают, что

$$P = (1-\delta) I + \delta \tilde{P},$$

где $\delta \in (0,1]$, $I = \mathrm{diag}\{1,...,1\}$ – единичная матрица,

$$\tilde{p}_{ij} = \begin{cases} \left|\{k:\ (i,k) \in E\}\right|^{-1}, i \ne j \\ 0, \quad \text{иначе} \end{cases}.$$

Такая специфика матрицы $P$ нами будет использоваться в подразделе 4.1.4. Отметим также, что вместо $I$ часто берут стохастическую матрицу с одинаковыми элементами – матрицу телепортации (см. [241]). Это сразу дает оценку снизу $\delta$ на спектральную щель матрицы $P$.

Опишем вкратце структуру раздела. В подразделе 4.1.2 приводится обзор наиболее популярных известных ранее способов численного поиска вектора PageRank. Новизна заключается в том, что этот обзор делается одновременно с декларированием двух новых способов поиска вектора PageRank: на основе Markov chain Monte Carlo и на основе алгоритма Григориадиса–Хачияна. В подразделе 4.1.3, также носящем обзорный характер, кратко описываются основные необходимые в дальнейшем факты о методе Markov chain Monte Carlo (MCMC). Новых результатов в этом пункте нет. Тем не менее, интерес может представлять обзор литературы. В подразделе 4.1.4 метод MCMC применяется для поиска вектора PageRank. Новыми здесь являются оценки скорости сходимости такого метода, а также оценки общего числа затраченных арифметических операций. Подчеркнем, что в отличие от общего случая, мы ограничиваемся в данном разделе изучением эффективности метода MCMC на специальном классе задач. За счет этого удается получить сущест-венно лучшие оценки, чем можно было бы ожидать в общем случае. Новизна метода также заключается в том, как организуются случайные блуждания. Отметим, что предложена хорошо параллелизуемая версия метода MCMC. Описанный в предыдущем пункте метод MCMC будет хорошо работать, если спектральная щель матрицы $P$ достаточно велика. В подразделе 4.1.5 предложен новый способ поиска вектора PageRank, не требующий ограничений на спектральную щель (см. также следующий раздел этой главы и раздел 6.1 главы 6). Этот способ сводит поиск ранжирующего вектора к поиску равновесия Нэша в антагонистической матричной игре. Для поиска равновесия мы используем метод Григориадиса–Хачияна, который позволяет учитывать разреженность матрицы $P$ и хорошо распараллеливается.

### 4.1.2 Обзор и обсуждение известных ранее и новых способов численного поиска вектора PageRank

В работах [44, 88, 102, 120, 225, 264, 273, 276, 317] предложены различные способы численного поиска вектора PageRank $p_*$.

**Замечание 4.1.3.** Строго говоря, в [120, 225] решались другие задачи, но несложно перенести алгоритмы этих работ на задачу поиска вектора PageRank: в случае [225] это делается тривиально [44], а вот в случае [120] потребовалось немного более точное исследование сходимости – поправка $\ln\left(\sigma^{-1}\right)$ появилась у нас из-за того, что мы избавились от математического ожидания в критерии качества (цели).

Приведем краткое резюме сложностных оценок алгоритмов работ [44, 88, 102, 120, 225, 264, 273, 276, 317] и алгоритмов, предложенных в этом разделе. “Сложность” понимается как количество арифметических операций типа умножения двух чисел, которые достаточно осуществить, чтобы с вероятностью не меньше $1-\sigma$ достичь точности решения $\varepsilon$ по “Целевому” функционалу.

Таблица 4.1.1

| Метод | Условие | Сложность | Цель (min) |
|---|---|---|---|
| Назина–Поляка [88] | Нет | $O\left(\frac{n\ln(n/\sigma)}{\varepsilon^2}\right)$ | $\left\|P^T p-p\right\|_2^2$ |
| Ю.Е. Нестерова [44, 264, 273] | $\bar{s}$ | $O\left(\frac{s^2\ln n}{\varepsilon^2}\right)$ | $\left\|P^T p-p\right\|_2$ |
| Юдицкого–Лана–Немировского–Шапиро [44, 225] | Нет | $\mathrm{O}\left(\frac{n\ln(n/\sigma)}{\varepsilon^2}\right)$ | $\left\Vert P^T p - p\right\Vert_\infty$ |
| Григориадиса–Хачияна [120] | $\bar{S}$ | $\mathrm{O}\left(\frac{s\ln n\ln(n/\sigma)}{\varepsilon^2}\right)$ | $\left\Vert P^T p - p\right\Vert_\infty$ |
| Нестерова–Немировского [276] | G, S | $\frac{sn}{\alpha}\ln\left(\frac{2}{\varepsilon}\right)$ | $\left\Vert p - p_*\right\Vert_1$ |
| Поляка–Трембы [102] | S | $\frac{2sn}{\varepsilon}$ | $\left\Vert P^T p - p\right\Vert_1$ |
| Д. Спилмана [317] | G, S | $\mathrm{O}\left(\frac{s^2}{\alpha\varepsilon}\ln\left(\frac{1}{\varepsilon}\right)\right)$ | $\left\Vert p - p_*\right\Vert_\infty$ |
| MCMC | SG | $\mathrm{O}\left(\frac{\ln n\ln(n/\sigma)}{\alpha\varepsilon^2}\right)$ | $\left\Vert p - p_*\right\Vert_2$ |

Поясним основные сокращения, использованные в Таблице 4.1.1:

- *G-условие* – наличие такой web-страницы (например, страницы, отвечающей самой поисковой системе **G**oogle), на которую можно перейти из любой другой web-страницы с вероятностью не меньшей, чем $\alpha \gg n^{-1}$, не ограничивая общности, будем считать, что вершина, отвечающая этой web-странице, имеет номер 1 (см. алгоритм MCMC);
- *S-условие* – из каждой web-страницы (в среднем) выходит не более $s \ll n$ ссылок на другие, то есть имеет место разреженность матрицы $P$ (**S**parsity); если из каждой web-страницы одновременно выходит и входит не более $s \ll n$ ссылок, то будем говорить об $\bar{S}$ условии;

**Замечание 4.1.4.** В ряде случаев (например, для метода из [102]) можно релаксировать *S*-условие: “выходит в среднем”.

- *SG-условие* – спектральная щель $\alpha$ (**S**pectral **G**ap) матрицы $P$ удовлетворяет условию $\alpha \gg n^{-1}$, где $\alpha$ – расстояние между максимальным по величине модуля собственным значением (числом Фробениуса–Перрона) матрицы $P$ (равным 1) и модулем следующего (по величине модуля) собственного значения. Если выполняется G-условие, то выполняется и SG-условие с $\alpha$ не меньшим, чем в G-условии [241, 317].

Отметим, что приведенная таблица немного огрубляет результаты процитированных работ, в частности, работ [264, 273] (один из методов Ю.Е. Нестерова, доставляющий аналогичную оценку, нами уже рассматривался в примере 1.1.2 подраздела 1.1.2 раздела 1.1 главы 1, в следующем разделе мы описываем новый метод с аналогичной оценкой, но ра-ботающий на практике быстрее). Это сделано для большей наглядности. Алгоритм из [88] на практике работает не очень быстро из-за большой константы в $O(\cdot)$. Отчасти похожая ситуация и с алгоритмами из [120, 225], они также работают не так быстро, как можно было бы ожидать. Связано это также и с тем, что в стандартных пакетах (использовался MatLab) довольно не эффективно реализована возможность работы со случайностью в огромных размерностях. Алгоритмы из [102, 276] и MCMC работают в точности по приведенным оценкам. Заметно лучше приведенной оценки на практике работает алгоритм из [317], причем речь идет не о константе в $O(\cdot)$. Тем не менее, условия, при которых этот алгоритм эффективен довольно обременительные, и даже при этих условиях он, как правило, доминируем методом из [276]. Практический анализ (проведенный, на тот момент студентами ФУПМ МФТИ, Денисом Дмитриевым, Алексеем Золотаревым и Айсой Имеевой) всех приведенных в таблице (и не только) алгоритмов показал эффективность использования методов из [102, 120], когда мы не можем гарантировать то, что спектральная щель не мала. Если же мы можем это обеспечить, то неплохо работает метод MCMC при условии, что блок рандомизации пишется самостоятельно, т.е. не используется готовые способы генерирования дискретных случайных величин (последнее замечание было основательно проработано А.С. Аникиным). Отметим также цикл экспериментов Сергея Кима, Даниила Смирнова и Степана Плаунова (школьников 2007 Московской школы, выполнявших в июле 2016 г. в летнем лагере “Сириус” проект, посвященный эффективным методам поиска вектора PageRank), дополняющих алгоритмы из таблицы 1:

https://github.com/KoIIdun/PageRank-gradient .

В ходе совместных экспериментов со школьниками, в частности, было показано, что для задачи с функционалом $\left\|P^T p - p\right\|_2^2$ (матрица $P$ порождена моделью Бакли–Остгуса [108]) адаптивные быстрые градиентные методы из подраздела 2.2.4 и раздела 2.3 главы 2 при $n = 10^4$ за несколько секунд работы (на процессоре Intel Pentium CPU 2.10 ГГц, ОЗУ 2 Гб на 64 разрядной операционной системе Windows 10; код написан на языке Python 3.5) достигают точности $\sim 10^{-6} - 10^{-7}$ (при значении функционала в точке старта $\sim 10^{-1} - 10^{-2}$).

Обращает на себя внимание то, что у разных алгоритмов разные “цели”. В связи с этим полезно отметить, что “типично” (см. [74]):

$$\|\cdot\|_1 \sim \sqrt{n}\|\cdot\|_2,\ \|\cdot\|_2 \sim \sqrt{n}\|\cdot\|_\infty,$$

причем это соответствует векторам с одинаковыми по порядку компонентами. В нашем случае это как раз не типично (имеет место степенной закон убывания компонент (см. [108])), поэтому можно ожидать, что приведенные оценки будут с более слабой, чем $\sqrt{n}$ зависимостью, то есть при переходе от одной нормы к другой фактор $\sqrt{n}$ будет заменен чем-то более близким к $O(1)$. Впрочем, мы не располагаем аккуратным обоснованием этого наблюдения.

По классификации главы 6 учебника [13] методы из [88, 264, 273] – являются вариационными, т.е. в этих методах решение системы $p^T = p^T P$ сводится к решению задач выпуклой оптимизации, которые, в свою очередь, решаются различными вариантами метода градиентного спуска (рандомизированный зеркальный спуск, метод Поляка–Шора, рандомизированный покомпонентный спуск).

Методы из [120, 225] – являются вариационными, но с игровым аспектом (оптимизационная задача понимается как задача поиска седловой точки – равновесия Нэша в матричной игре). Действительно, рассматриваемую нами задачу можно переписать (с помощью теории Фробениуса–Перрона [95]) как задачу:

$$f(p) = \max_{u \in S_n(1)} \langle u, Ap \rangle \to \min_{p \in S_n(1)},$$

где

$$A = P^T - I, \; A = \|a_{ik}\|, \; S_n(1) = \left\{ p \geq 0 : \sum_{k=1}^{n} p_k = 1 \right\}.$$

Отметим, что $0 \leq f(p) \leq \|Ap\|_\infty$, $p \in S_n(1)$, и на векторе PageRank (и только на нем) $f(p) = 0$. Фактически, этот целевой функционал $f(p)$ оказывается очень близким к $\|P^T p - p\|_\infty$. В основе обоих игровых методов лежат рандомизированные варианты метода зеркального спуска поиска равновесия в матричных играх.

Методы из [88, 120, 225] могут быть распараллелены на $\sim \log_2(1/\sigma)$ процессорах. Если параллельно запустить $\sim \log_2(1/\sigma)$ независимых траекторий любого из этих методов, то, используя точные оценки числа итераций, гарантирующих заданную точность и доверительный уровень, который в этой параллельной схеме выбираем $1/2$, останавливаем траектории после выполнения предписанного числа итераций и проверяем выполнение критерия малости целевого функционала. С вероятностью $1-\sigma$ хотя бы одна из траекторий выдаст требуемое по точности решение. К сожалению, проверка выполнения критерия малости целевого функционала требует осуществления умножения $P^T p$, что может быть более затратным, чем получить кандидата на решение. Поэтому, такое распараллеливание не всегда осмыслено. Кроме того, отмеченная проблема проверки выполнения критерия малости целевого функционала приводит к необходимости точного оценивания числа итераций, гарантирующих заданную точность и доверительный уровень (см., например, оценку на число итераций в теореме 4.1.1 подраздела 4.1.5). Отметим, что метод из [120] без учета разреженности допускает распараллеливание на $n/\ln n$ процессорах (см. [120]).

Остальные методы из таблицы можно интерпретировать как вариации метода простых итераций из [13, 110]. При этом сразу же возникает вопрос: почему бы не попытаться решать систему $p^T = p^T P$ каким-то уже известным способом, например, методом из [13]? Или, скажем, таким (см. [170]): $p^T = p^T(0)P^{\infty}$, где $P^{\infty}$ – то, что выдает следующий итерационный процесс, требующий конечного числа итераций:

$$P(0) = I\,,\ P(k) = I - k\frac{P(k-1)(I-P)}{\mathrm{tr}\left(P(k-1)(I-P)\right)},\ k \in \mathbb{N}.$$

Итерации заканчиваются, когда первый раз выполнится условие $\mathrm{tr}\left(P(m)(I-P)\right) = 0$, где $\mathrm{tr}(B)$ – след матрицы $B$, т.е. сумма диагональных элементов. При этом $P^{\infty} = P(m)$.

Постараемся здесь ответить на этот вопрос. Прежде всего, многие итерационные методы требуют, чтобы матрица $P$ была симметричной и положительно определенной, что, как правило, место не имеет. Но даже если все необходимые условия будут выполнены и матрица при этом будет еще и разрежена, то будет наблюдаться такая же сходимость, как и в методе из [276], только вместо константы $\alpha$ будет фигурировать некий её аналог (из SG-условия), который в общем случае намного сложнее оценить. Таким образом, вместо того, чтобы приводить в таблице линейку классических итерационных методов с их оценками, мы ограничились тем, что привели наиболее приспособленные из этой линейки методы (условия применимости которых хорошо интерпретируемы) для решения рассматриваемого (довольно узкого) класса задач. Метод из [118, 170] мы не можем использовать, потому что $n \gg 10^4$ (т.е. $n^3 \gg 10^{12}$ – такое количество арифметических операций на одном персональном компьютере может быть выполнено за час).

В связи с упоминанием методов простой итерации, можно обратить внимание на то, что в описанных методах возможны проблемы накапливания ошибок округления (конечности длины мантиссы), возникающие при (см., например, § 4 главы 6 [13] или § 1 главы 5 [110]) условии $\|P\| > 1$ (при этом спектр матрицы может лежать в единичном круге, и метод простой итерации теоретически (без ошибок округления) должен сходиться со скоростью геометрической прогрессии – см. также раздел 4.3 этой главы 4). Однако в рассмат-риваемом нами случае в естественной норме, подчиненной $l_1^n$, $\|P\|=1$, и таких проблем не возникает.

Важно также обратить внимание на то, что в ряде случаев, рассмотренных в таблице, целью является получение такого вектора распределения вероятностей $\vec{p}$, который давал бы малую невязку “по функции”. Но из того, что $\|P^T p - p\|$ малó не следует, что будет малó $\|p - p_*\|$.

**Замечание 4.1.5.** Напомним, что вектор $p_*$ – решение уравнения: $P^T p = p$. Мы считаем это решение единственным в классе распределений вероятностей.

Более того, вполне естественно ожидать обратное, что $\|p - p_*\|$ окажется в $\alpha^{-1}$ раз больше (см. [13, 80]). Это оценка сверху, но, по-видимому, с большой вероятностью (если выбирать точку старта равновероятно) она, в действительности, превратится практически в точное соотношение (см.[80, 121]). Таким образом, добиться малости $\|P^T p - p\|$ – совсем не значит полностью решить задачу поиска вектора PageRank. Это обстоятельство, а также тот факт, что решение ищется на сильно ограниченном множестве (единичном симплексе) и дополнительно делаются всякие упрощающие предположения, отчасти объясняют, почему сложность описанных алгоритмов оказывается столь небольшой. Ведь, например, из того, что $\|p - p_*\|_\infty \le \varepsilon$ в случае, когда $\varepsilon \gg n^{-1}$ и истинное распределение $p_*$ близко к равномерному, совсем не ясно, что именно выдает алгоритм с точностью $\varepsilon$ по функции, и как это можно использовать для ранжирования web-страниц. К счастью, такие ситуации, когда истинное распределение $p_*$ близко к равномерному – нетипичны на практике (см. [107]), особенно в случае выполнения G-условия (SG-условия). Как правило, выделяются компоненты вектора $p_*$, которые достаточно велики. А поскольку содержательно задача формулируется как ранжирование web-страниц, то по-сути, речь идет о восстановлении нескольких первых по величине компонент вектора $p_*$. Другими словами, точно восстанавливать малые компоненты вектора $p_*$ не требуется, если мы знаем, что они достаточно малы и есть достаточное количество не малых компонент. Скажем вектор, выдаваемый алгоритмом [317], имеет отличными от нуля не более чем $3\varepsilon^{-1}$ компонент. Отмеченные обстоятельства объясняют, почему даже при $\varepsilon \gg n^{-1}$ может быть полезен вектор $p$, доставляющий оценку $\|P^T p - p\|_\infty \le \varepsilon$ – в одной из самых “плохих” норм: $l_\infty^n$.

### 4.1.3 Метод Markov Chain Monte Carlo

Идея решать линейные уравнения с помощью MCMC столь же стара, сколь и обычный метод Монте-Карло из [66]. Однако мы имеем дело не с линейным уравнением общего вида, а с уравнением со стохастической матрицей $P$, причем специальным образом заполненной, эти обстоятельства позволят нам более экономно организовать случайное блуждание по графу, соответствующему матрице $P$.

Прежде чем излагать алгоритм, приведем некоторые вспомогательные факты.

**Алгоритм Кнута–Яо (см. [**66, 78**]).** С помощью бросаний симметричной монетки требуется сгенерировать распределение заданной дискретной случайной величины (с.в.), принимающей конечное число значений. Предположим, что нам нужно сгенерировать распределение с.в., принимающей три значения 1, 2, 3 с равными вероятностями 1/3. Действуем таким образом. Два раза кидаем монетку: если выпало 00, то считаем что выпало значение 1, если 01, то 2, если 11, то 3. Если 10, то еще два раза кидаем монетку и повторяем рассуждения. Можно показать, что описанную выше схему можно так обобщить, чтобы сгенерировать распределение дискретной с.в., принимающей, вообще говоря, с разными вероятностями $n$ различных значений, в среднем с помощью не более чем $\log_2(n-1)+2$ подбрасываний симметричной монетки. Если эти вероятности одинаковы, то процедура “приготовления” такого алгоритма генерирования также имеет логарифмическую сложность по $n$.

**Алгоритм Markov Chain Monte Carlo (MCMC) (см. [133, 186, 193, 197, 222, 224, 246, 247, 316, 312]).** Чтобы построить однородный дискретный марковский процесс с конечным числом состояний, имеющий наперед заданную инвариантную (стационарную) меру $\pi$, переходные вероятности ищутся в следующем виде: $p_{ij} = p_{ij}^0 b_{ij}$, $i \neq j$; $p_{ii} = 1 - \sum_{j:\ j\neq i} p_{ij}$, где $p_{ij}^0$ – некоторая “затравочная” матрица, которую будем далее предполагать симметричной. Легко проверить, что матрица $p_{ij}$ имеет инвариантную (стационарную) меру $\pi$:

$$\frac{b_{ij}}{b_{ji}} = \frac{\pi_j p_{ji}^0}{\pi_i p_{ij}^0} = \frac{\pi_j}{\pi_i},\ p_{ij}^0 > 0.$$

Чтобы найти $b_{ij}$ достаточно найти функцию $F : \mathbb{R}_+ \to [0,1]$ такую, что

$$\frac{F(z)}{F(1/z)} = z$$

и положить

$$b_{ij} = F\left(\frac{\pi_j p_{ji}^0}{\pi_i p_{ij}^0}\right) = F\left(\frac{\pi_j}{\pi_i}\right).$$

Пожалуй, самый известный пример такой функции $\tilde{F}(z) = \min\{z,1\}$ – алгоритм (Хастинг-са–)Метрополиса [246]. Заметим, что для любой такой функции $F(z)$ имеем $F(z) \le \tilde{F}(z)$. Другой пример дает функция $F(z) = z/(1+z)$. Заметим также, что $p_{ij}^0$ обычно выбирается равным $p_{ij}^0 = 1/M_i$, где $M_i$ число “соседних” состояний у $i$, или

$$p_{ij}^0 = 1/(2M),\ i \ne j;\ p_{ii}^0 = 1/2,\ i \ne j.$$

При больших значениях времени $t$, согласно эргодической теореме, имеем, что распределение вероятностей близко к стационарному $\pi$. Действительно, при описанных выше условиях имеет место условие детального баланса (марковские цепи, для которых это условие выполняется, иногда называют обратимыми):

$$\pi_i p_{ij} = \pi_j p_{ji},\ i, j = 1,\dots,n,$$

из которого сразу следует инвариантность меры $\pi$, т.е.

$$\sum_i \pi_i p_{ij} = \pi_j \sum_i p_{ji} = \pi_j,\ j = 1,\dots,n.$$

Основное применение замеченного факта состоит в наблюдении, что время выхода марковского процесса на стационарную меру (mixing time (см. [255])) во многих случаях оказывается удивительно малым.

**Замечание 4.1.6.** Более того, задача поиска такого симметричного случайного блуждания на графе (с равномерной инвариантной мерой в виду симметричности) заданной структуры, которое имеет “наименьшее” mixing time (другими словами, наибольшую спектральную щель), сводится к задаче полуопределенного программирования, которая, как известно, полиномиально (от числа вершин этого графа) разрешима [161].

При том, что выполнение одного шага по времени случайного блуждания по графу, отвечающему рассматриваемой марковской цепи, как следует из алгоритма Кнута–Яо, также может быть быстро сделано. Таким образом довольно часто можно получать эффективный способ генерирования распределения дискретной случайной величины с распределением вероятностей $\pi$ за время, полиномиальное от логарифма числа компонент вектора $\pi$.

Для оценки mixing time нужно оценить спектральную щель стохастической матрицы, переходных вероятностей, задающей исследуемую марковскую динамику, то есть нужно оценить расстояние от максимального собственного значения этой матрицы равного единицы (теорема Фробениуса–Перрона) до следующего по величине модуля. Именно это число определяет основание геометрической прогрессии, мажорирующей исследуемую последовательность норм разностей расстояний (по вариации) между распределением в данный момент времени и стационарным (финальным) распределением. Для оценки спектральной щели разработано довольно много методов, из которых мы упомянем лишь некоторые (см. [133, 186, 193, 197, 222, 224, 247, 316]): неравенство Пуанкаре (канонический путь), изопериметрическое неравенство Чигера (проводимость), с помощью техники каплинга [252] (получаются простые, но, как правило, довольно грубые оценки), с помощью каплинга Мертона [291], с помощью дискретной кривизны Риччи и теорем о концентрации меры (Мильмана–Громова [224, 244]). Приведем некоторые примеры применения MCMC [246]: Тасование $n$ карт, разбиением приблизительно на две равные кучи и перемешиванием этих куч (mixing time $\sim \log_2 n$); Hit and Run (mixing time $\sim n^3$); Модель Изинга – $n$ спинов на отрезке, стационарное распределение = распределение Гиббса, Глауберова динамика (mixing time $\sim n^{2\log_2 e/T}$, $0 < T \ll 1$); Проблема поиска кратчайших гамильтоновых путей; Имитация отжига для решения задач комбинаторной оптимизации, MCMC для решения задач перечислительной комбинаторики.

**Замечание 4.1.7.** На примере тасования карт контраст проявляется, пожалуй, наиболее ярко. Скажем для колоды из 52 карт пространство состояний марковской цепи будет иметь мощность 52! (если сложить времена жизней в наносекундах каждого человека, когда либо жившего на Земле, то это число на много порядков меньше 52!). В то время как такое тасование: взять сверху колоды карту, и случайно поместить ее во внутрь колоды, отвечающее определенному случайному блужданию, с очень хорошей точностью выйдет на равномерную меру, отвечающую перемешанной колоде, через каких-то 200–300 шагов. Если брать тасование разбиением на кучки, то и того меньше – за 8–10 шагов [75].

Продемонстрируем сказанное выше двумя примерами (с помощью которых, например, получаются оценки mixing time работы [32]), которые нам пригодятся в дальнейшем.

**Пример 4.1.1 (подход Чигера [197, 222, 316]).** Пусть (напомним, что под $\pi(\cdot)$ понимается инвариантная мера, а под $P = \left\| p_{ij} \right\|_{i,j=1}^{n,n}$ – матрица переходных вероятностей марковской цепи)

$$h(G) = \min_{S \subseteq V_G:\, \pi(S) \le 1/2} P\left(S \to \bar{S} \middle| S\right) = \min_{S \subseteq V_G:\, \pi(S) \le 1/2} \frac{\sum\limits_{(i,j)\in E_G:\, i\in S,\, j\in \bar{S}} \pi(i)\, p_{ij}}{\sum\limits_{i \in S} \pi(i)}, \quad \text{(константа Чигера)}$$

$$T(i,\varepsilon) = \Theta\left(h(G)^{-2}\left(\ln\left(\pi(i)^{-1}\right) + \ln\left(\varepsilon^{-1}\right)\right)\right). \quad \text{(Mixing time)}$$

Тогда

$$\forall\ i=1,...,n,\ t\ge T(i,\varepsilon)\to\left\|P^t(i,\cdot)-\pi(\cdot)\right\|_1=\sum_{j=1}^n\left|P^t(i,j)-\pi(j)\right|\le\varepsilon,$$

где $P^t(i,j)$ – условная вероятность того, что стартуя из состояния $i$ через $t$ шагов, марковский процесс окажется в состоянии $j$. Отметим, что от вектора PageRank $\pi$ мы вправе ожидать степенной закон убывания компонент, отсортированных по возрастанию (см. [107]), поэтому $\max\limits_{i=1,...,n}\ln\left(\pi(i)^{-1}\right)=\mathrm{O}(\ln n)$.

**Пример 4.1.2 (coarse Ricci curvature [224]).** Введем расстояние Монжа–Канторовича между двумя (дискретными) распределениями вероятностей $\mu$ и $\nu$:

$$W_1(\mu,\nu)=\min_{\substack{\xi\ge0:\ \sum_j\xi(i,j)=\mu(i)\\ \sum_i\xi(i,j)=\nu(j)}}\sum_{i,j}d(i,j)\xi(i,j),$$

где каждой паре вершин поставлено в соответствие неотрицательное число $d(i,j)$ (со свойствами расстояния (метрики)). Говорят, что $\kappa$ – дискретная кривизна Риччи, если

$$\exists\ t_0>0:\ \forall\ i,j=1,...,n$$
$$W_1\left(P^{t_0}(i,\cdot),P^{t_0}(j,\cdot)\right)\le(1-\kappa)d(i,j).$$

Пусть существует такое $\kappa>0$, тогда

$$W_1\left(P^t(i,\cdot),\pi(\cdot)\right)\le\left[\kappa^{-1}\sum_{j=1}^n d(i,j)p_{ij}\right]\cdot(1-\kappa)^{t/t_0},$$

$$E\left[\left\|\frac{1}{T}\sum_{t=T_0}^{T+T_0}\vec{x}_t-\vec{\pi}\right\|_2^2\right]=(\mathrm{Bias})^2+\mathrm{Var}=\left(\mathrm{O}\left(n\frac{(1-\kappa)^{T_0}}{\kappa T}\right)\right)^2+\mathrm{O}\left(\frac{1}{\kappa T}\right),$$

где у случайного вектора $\vec{x}_t$ все компоненты нулевые, кроме единственной компоненты, отвечающей вершине, в которой находился марковский случайный процесс на шаге $t$. Считая, что

$$\mathrm{Var}\gg(\mathrm{Bias})^2$$

имеем (этот результат в явном виде не содержится в [224], но к нему можно прийти, используя некоторые идеи работы [291])

$$\exists\ c>0:\ P\left(\left\|\frac{1}{T}\sum_{t=T_0}^{T+T_0}x_t-\pi\right\|_2>c\sqrt{\frac{\ln n+\Omega}{\kappa T}}\right)\le\exp(-\Omega).$$

Результат вида

$$\frac{1}{T}\sum_{t=T_0}^{T+T_0} x_t \xrightarrow[T\to\infty]{} \pi$$

– есть эргодическая теорема для марковских процессов (при этом выше была приведена довольно тонкая оценка того, какая скорость сходимости), вполне ожидаем (см. [102]) и соответствует классическим вариантам эргодических теорем для динамических системах (Биркгофа–Хинчина, фон Неймана). Правда, в отличие от динамических систем здесь удается оценить скорость сходимости.

Отметим связь описанного в примере 4.1.2 подхода с результатами о сжимаемости (в пространстве всевозможных лучей с центром в начале координат неотрицательного ортанта) в метрике Биркгофа положительных линейных операторов (см. [80]), в частности, заданных стохастической матрицей $P$, некоторая степень которой имеет все элементы положительными (это равносильно неразложимости и непериодичности марковской цепи (см. [75])). Кстати сказать, такое понимание эргодической теоремы для марковских цепей позволяет интерпретировать её как теорему о сжимающих отображениях, что выглядит несколько необычно в контексте сопоставления этой теоремы с эргодическими теоремами для динамических систем.

### 4.1.4 Алгоритм MCMC

Теперь можно приступать к изложению нужной нам версии алгоритма MCMC.

**Шаг 1. Инициализация:** $X=0$, вершина $=1$, $t=0$.

**Шаг 2. Счетчик итераций:** $t:=t+1$.

**Шаг 3. Модификация** $X$ **:** если $t \ge T^0_{\varepsilon,\alpha,n}$, то $X_k := X_k + 1$, где $k$ – номер текущей вершины.

**Шаг 4. Модификация вершины:** случайно "переходим" из текущей вершины в одну из "соседних" согласно матрице $P$.

**Шаг 5. Остановка:** если $t < T_{\varepsilon,\sigma,\alpha,n}$ перейти на шаг 2, иначе на шаг 6.

**Шаг 6. Ответ:** $p = X\Big/\left(T_{\varepsilon,\sigma,\alpha,n} - T^0_{\varepsilon,\alpha,n}\right)$.

Самым затратным шагом из первых пяти шагов является шаг 4, но даже этот шаг при самом неблагоприятном раскладе (из текущей вершины выходит порядка $n$ ребер) можно осуществить за $\mathrm{O}(\ln n)$ операций (см., например, алгоритм Кнута–Яо – на самом деле все это можно сделать разными способами, причем сложность $\mathrm{O}(\ln n)$ можно понимать не в среднем, а обычном образом (см. [66]); алгоритм Кнута–Яо был приведен лишь как иллюстративный пример). Из примеров 4.1.1 и 4.1.2 (мы заменяем оценки спектральной щели в этих примерах на саму спектральную щель $\alpha$, что можно делать в таком контексте для обратимых марковских цепей (см. [291]), и нуждается в некоторых оговорках в общем случае), получаем, что при

$$T^0_{\varepsilon,\alpha,n} = \mathrm{O}\left(\frac{1}{\alpha}\ln\left(\frac{n}{\varepsilon}\right)\right)$$

после

$$T_{\varepsilon,\sigma,\alpha,n} = \mathrm{O}\left(\frac{\ln(n/\sigma)}{\alpha\varepsilon^2}\right)$$

итераций с вероятностью не меньшей $1-\sigma$ алгоритм MCMC выдаст $\varepsilon$-оптимальное (по функции) решение исходной задачи. При этом алгоритм MCMC затрачивает в общей сложности

$$\mathrm{O}\left(n + \frac{\ln n \ln(n/\sigma)}{\alpha\varepsilon^2}\right)$$

элементарных арифметических операций (типа умножения двух чисел с плавающей точкой) – слагаемое $\mathrm{O}(n)$ “отражает” стоимость шага 6.

Несложно предложить способ эффективного распараллеливания такого алгоритма. Для этого выпускается не одна траектория, а

$$N_{\varepsilon,\sigma} = \mathrm{O}\left(\frac{\ln(\sigma^{-1})}{\varepsilon^2}\right),$$

независимо блуждающих, и стартующих с вершин, выбираемых случайно и независимо. За каждой траекторией следим время $T^0_{\varepsilon/2,\alpha,n}$. Потом усредняем (см. ниже) по всем этим траекториям, получаем ответ (с возможностью организации параллельных вычислений), но при этом придется затратить чуть больше операций – вместо множителя $\ln(n/\sigma)$ появится $\ln(n/\varepsilon)\ln(\sigma^{-1})$, что не так уж и плохо.

**Замечание 4.1.8.** Более того, это обстоятельство на практике можно частично нивелировать, например, таким образом. Сначала выпускается одна траектория. Блуждающая частица через случайные моменты времени “рождает потомков”, которые также начинают независимо блуждать, стартуя с “места рождения”, причем рожденные частицы также склонны к “спонтанному делению”. Здесь важно правильно подобрать интенсивность такого деления (размножения).

Действительно, мы получаем набор из $N_{\varepsilon,\sigma}$ независимых одинаково распределенных случайных векторов $x^k_{T^0_{\varepsilon,\alpha,n}}$, $k = 1,\dots,N_{\varepsilon,\sigma}$, где $x^k_{T^0_{\varepsilon,\alpha,n}}$ – вектор, все компоненты которого равны нулю кроме одной, равной 1; эта компонента соответствует состоянию, в котором находится $k$-е блуждание на шаге $T^0_{\varepsilon/2,\alpha,n}$. Векторы одинаково распределены, причем (напомним, что $p_* = \pi$ – инвариантная мера)

$$\left\| E\left[ x^k_{T^0_{\varepsilon/2,\alpha,n}} - \pi \right] \right\|_1 \le \varepsilon/2.$$

Применяем далее (этот способ был предложен совместно с Е.Ю. Клочковым [77], см. также раздел 1.2 главы 1) к

$$\frac{1}{N_{\varepsilon,\sigma}} \sum_{k=1}^{N_{\varepsilon,\sigma}} x^k_{T^0_{\varepsilon/2,\alpha,n}}$$

неравенство типа Хефдинга в гильбертовом пространстве $l_2^n$ (см. [159]), получаем

$$P\left( \left\| \frac{1}{N_{\varepsilon,\sigma}} \sum_{k=1}^{N_{\varepsilon,\sigma}} x^k_{T^0_{\varepsilon/2,\alpha,n}} - \pi \right\|_2 \ge \varepsilon \right) = P\left( \left\| \sum_{k=1}^{N_{\varepsilon,\sigma}} \left( x^k_{T^0_{\varepsilon/2,\alpha,n}} - \pi \right) \right\|_2 \ge \varepsilon N_{\varepsilon,\sigma} \right) \le$$

$$\le \exp\left( -\frac{1}{4N_{\varepsilon,\sigma}} \left( \varepsilon N_{\varepsilon,\sigma} - E\left( \left\| \sum_{k=1}^{N_{\varepsilon,\sigma}} \left( x^k_{T^0_{\varepsilon/2,\alpha,n}} - \pi \right) \right\|_2 \right) \right)^2 \right),$$

где

$$E\left( \left\| \sum_{k=1}^{N_{\varepsilon,\sigma}} \left( x^k_{T^0_{\varepsilon/2,\alpha,n}} - \pi \right) \right\|_2 \right) \le \sqrt{E\left( \left\| \sum_{k=1}^{N_{\varepsilon,\sigma}} \left( x^k_{T^0_{\varepsilon/2,\alpha,n}} - \pi \right) \right\|_2^2 \right)} \le \sqrt{2\sqrt{2} N_{\varepsilon,\sigma} + \left( \varepsilon^2/4 \right) N^2_{\varepsilon,\sigma}} .$$

Подберем $N_{\varepsilon,\sigma}$ так, чтобы выполнялось неравенство

$$P\left( \left\| \frac{1}{N_{\varepsilon,\sigma}} \sum_{k=1}^{N_{\varepsilon,\sigma}} x^k_{T^0_{\varepsilon/2,\alpha,n}} - \pi \right\|_2 \ge \varepsilon \right) \le \sigma.$$

Для этого достаточно, чтобы выполнялось неравенство

$$\varepsilon N_{\varepsilon,\sigma} - \sqrt{2\sqrt{2} N_{\varepsilon,\sigma} + \left( \varepsilon^2/4 \right) N^2_{\varepsilon,\sigma}} \ge \sqrt{4N_{\varepsilon,\sigma} \ln\left( \sigma^{-1} \right)},$$

где

$$N_{\varepsilon,\sigma} = \frac{4 + 6\ln\left( \sigma^{-1} \right)}{\varepsilon^2} = \mathrm{O}\left( \frac{\ln\left( \sigma^{-1} \right)}{\varepsilon^2} \right).$$

Отметим, что с точностью до мультипликативной константы эта оценка не может быть улучшена. Это следует из неравенства Чебышёва

$$P\left( X \ge EX - \varepsilon \right) \ge 1 - \frac{\mathrm{Var}\left( X \right)}{\varepsilon^2},$$

при

$$X = \left\| \frac{1}{N_{\varepsilon,\sigma}} \sum_{k=1}^{N_{\varepsilon,\sigma}} x^k_{T^0_{\varepsilon/2,\alpha,n}} - \pi \right\|_2, \ \pi = \left(n^{-1},...,n^{-1}\right)^T, \ N_{\varepsilon,\sigma} \gg n.$$

**Замечание 4.1.9.** Стоит отметить, что при экспоненциальном убывании компонент отранжированного вектора $\pi$ можно заменить полученные далее оценки в 2-норме на аналогичные оценки (с другими константами и дополнительным логарифмическим штрафом), но в 1-норме. Понять это можно из того, что фактически речь идет о задаче восстановления параметров мультиномиального распределения. Точнее, речь идет о влиянии размерности пространства параметров на оценки скорости сходимости в теореме Фишера о методе наибольшего правдоподобия в неасимптотическом варианте [319] (выборочные средние как раз и будут оценками, полученными согласно этому методу). Для данного примера мультиномиальное распределение порождает (вот в этом месте, к сожалению, приходится использовать формулу Стирлинга, то есть требуются некоторые асимптотические оговорки) в показателе экспоненты расстояние Кульбака–Лейблера между выборочным средним и истинным распределением [111]. В свою очередь это расстояние оценивается согласно неравенству Пинскера квадратом 1-нормы [159].

Немного удивляет крайне слабая зависимость общей сложности от размера матрицы $P$, особенно, если учесть, что никаких предположений о разреженности $P$ не делалось. “Подвох” здесь в том, что, в действительности, если речь не идет о каких-то специальных графах (например, экспандерах [109]), на которых рассматривается случайное блуждание, то условие, что $\alpha$ равномерно по $n$ отделимо от нуля – противоестественно. Так, в работе [32] (см. также приложение в конце диссертации) приводится следующая, по-видимому, не улучшаемая оценка $\alpha \sim n^{-1}$ для класса часто встречающихся на практике марковских процессов, которые возникают при описании различных макросистем, “живущих” в приближении среднего поля (в динамику заложено равноправие взаимодействующих агентов – закон действующих масс Гульдберга–Вааге). Также можно заметить, что выигрыш в оценке сложности в зависимости от $n$ происходит за счет довольно плохой зависимости от $\varepsilon$. Вообще ожидать зависимости сложности от $\varepsilon$ лучшей, чем $\varepsilon^{-2}$ в рандомизированных алгоритмах не приходится (см. [91]), поэтому рандомизация осмыслена, как правило, только при $\varepsilon \gg n^{-1}$. Это становится особенно ясно, если сравнить полученную оценку сложности с оценкой сложности, скажем, алгоритма из [276]. То есть, чтобы MCMC имело смысл использовать нужно, чтобы $\varepsilon \gg n^{-1}$.

Основным недостатком метода MCMC является отсутствие точного знания о $T_{\varepsilon,\sigma,\alpha}$ и $T^0_{\varepsilon,\alpha,n}$, даже если известен размер спектральной щели $\alpha$. Более того, эффективность алго-ритма напрямую завязана на оценку спектральной щели $\alpha$, которая, как правило, априорно не известна. Тем не мене, проблема оценки $\alpha$ может быть решена за $Cn$ арифметических операций (к сожалению, с довольно большой константой $C$), например, с помощью $\delta^2$-процесса практической оценки спектральной щели (см., например, [13]). Другой способ оценки $\alpha$ и $T^0_{\varepsilon,\alpha,n}$ имеется в [292]. Проблема определения $T_{\varepsilon,\sigma,\alpha}$ решается с помощью контроля разности $\|p_{t+\tau} - p_t\|_2$. Относительно $T^0_{\varepsilon,\alpha,n}$ можно действовать так: изначально запускать алгоритм с $T^0_{\varepsilon,\alpha,n} = 0$, а затем скорректировать полученный ответ, полагая, скажем, $T^0_{\varepsilon,\alpha,n} = T_{\varepsilon,\sigma,\alpha}/5$.

Отметим также, что изначально необходимо правильным образом разместить в памяти компьютера матрицу $P$. Если работать с этой матрицей обычным образом, то при случайном выборе алгоритмом MCMC на каждом шаге соседней вершины будет тратиться время порядка $n$, а не $\ln n$, как хотелось бы. Чтобы избежать этого, необходимо перед началом работы алгоритма представить граф в виде списка ссылок. Это займет время (память) порядка числа ребер, но сделать это нужно всего один раз.

### 4.1.5 Алгоритм Григориадиса–Хачияна

Рассматриваемую нами задачу перепишем (с помощью теории Фробениуса–Перрона [95]) как задачу поиска седловой точки (равновесия в антагонистической матричной игре):

$$f(p) = \max_{u \in S_n(1)} \langle u, Ap \rangle \to \min_{p \in S_n(1)}, \tag{4.1.1}$$

где $A = P^T - I$, $S_n(1) = \left\{ p \geq 0 : \sum_{k=1}^{n} p_k = 1 \right\}$. Отметим, что $f(p) \geq 0$, $p \in S_n(1)$, и на векторе PageRank (и только на нем) $f(p) = 0$.

С точностью (по функции) $\varepsilon$ и с вероятностью не меньшей $1-\sigma$ равновесие можно найти с помощью стохастического зеркального спуска из работы [120, 225] за количество операций

$$O\left( \frac{n \ln(n/\sigma)}{\varepsilon^2} \right).$$

Однако эти методы не достаточно полным образом учитывают специфику задачи в случае разреженной матрицы $P$. Для решения задачи (4.1.1) воспользуемся методом поиска равновесий в симметричных антагонистических матричных играх из [120] (отметим, что этот метод по сути является определенным образом рандомизированным методом зеркального спуска из [91]). Для этого, предварительно приведем задачу (4.1.1), следуя

Данцигу [120], к симметричному виду (этого можно и не делать, если использовать рандомизированный онлайн метод зеркального спуска из теоремы 2 работы [54]):

$$\max_{u\in S_{2n+1}(1)}\langle u, \mathrm{A}x\rangle \to \min_{x:=(y,p',u)\in S_{2n+1}(1)}, \quad \mathrm{A} := \begin{bmatrix} 0 & A & -e \\ -A^T & 0 & e \\ e^T & -e^T & 0 \end{bmatrix},$$

где $e=(1,...,1)^T$, $\mathrm{A}=\|a_{ik}\|$. Тогда, $f(p)\le 2\varepsilon$, где $p=p'/(e^T p')$, причем $e^T p'\ge 1/2-\varepsilon$, если $\mathrm{A}x\le\varepsilon e$.

**Шаг 1. Инициализация:** $X=0$, $p^T=\frac{1}{2n+1}\underbrace{(1,...,1)}_{2n+1}^T$, $t=0$.

**Шаг 2. Счетчик итераций:** $t:=t+1$.

**Шаг 3. Датчик случайных чисел:** выбираем $k\in\{1,...,2n+1\}$ с вероятностью $p_k$.

**Шаг 4. Модификация** $X$ **:** $X_k := X_k+1$.

**Шаг 5. Модификация** $p$ **:** $i=1,...,2n+1$ $p_i := p_i\exp\left(\frac{\varepsilon a_{ik}}{2}\right)$.

**Шаг 6. Остановка:** если $t<T_{\varepsilon,\sigma,n}$ перейти на шаг 2, иначе на шаг 7.

**Шаг 7. Ответ:** $x=X/t$.

**Замечание 4.1.10 (к шагу 5).** Далее будет пояснено, что достаточно задавать распределение вероятностей с точностью до нормирующего множителя.

Итак, считая матрицу $P$ разреженной с не более чем $s$ элементами в столбце и строке (см. [273]), получаем, что самыми "дорогими" шагами будет пересчет (перегенерирование) распределения вероятностей (шаги 3 и 5). Используя способ генерирования дискретной случайный величины с помощью сбалансированного двоичного дерева, мы получаем, что изменение веса вероятности одного исхода, по сути, равносильно процедуре изменение весов тех вершин дерева, путь через которые ведет к изменяемому значению "листа" дерева.

**Замечание 4.1.11.** Опишем точнее эту процедуру. У нас есть сбалансированное двоичное дерево высоты $\mathrm{O}(\log_2 n)$ с $2n+1$ листом (дабы не вдаваться в технические детали, считаем что число $2n+1$ – есть степень двойки, понятно, что это не так, но можно взять, скажем, наименьшее натуральное $m$ такое, что $2^m>2n+1$ и рассматривать дерево с $2^m$ листом) и с $\mathrm{O}(n)$ общим числом вершин. Каждая вершина дерева (отличная от листа) имеет неотрицательный вес равный сумме весов двух ее потомков. Первоначальная процедура приготовления дерева, отвечающая одинаковым весам листьев, потребует $\mathrm{O}(n)$ операций. Но сделать это придется всего один раз. Процедура генерирования дискретной случайной величины с распределением, с точностью до нормирующего множителя, соответствующим весам листьев может быть осуществлена с помощью случайного блуждания из корня дерева к одному из листьев. Отметим, что поскольку дерево двоичное, то прохождение каждой его вершины при случайном блуждании, из которой идут два ребра в вершины с весами $a>0$ и $b>0$, осуществляется путем подбрасывания монетки (“приготовленной” так, что вероятность пойти в вершину с весом $a$ – есть $a/(a+b)$). Понятно, что для осуществления этой процедуры нет необходимости в дополнительном условии нормировки: $a+b=1$. Если вес какого-то листа алгоритму необходимо поменять по ходу работы, то придется должным образом поменять дополнительно веса тех и только тех вершин, которые лежат на пути из корня к этому листу. Это необходимо делать, чтобы поддерживать свойство: каждая вершина дерева (отличная от листа) имеет вес равный сумме весов двух ее потомков.

Обратим при этом внимание, что нет необходимости заниматься перенормировкой распределения вероятностей (это бы стоило $O(n)$), то есть изменением весов вершин дерева, отличных от тех, путь через которые ведет к изменяемому значению листа дерева. Все это (сгенерировать с помощью дерева и обновить дерево) можно сделать за $O(\ln n)$ операций типа сравнения двух чисел, а в (“типичном”) случае, когда это нужно делать $\sim s$ раз, то за $O(s\ln n)$ операций (и лишь в случаях, когда $k=2n+1$, придется делать $\sim 2n$ операций). Далее будет показано, что после

$$T_{\varepsilon,\sigma,n}=O\left(\frac{\ln(n/\sigma)}{\varepsilon^2}\right)$$

итераций с вероятностью, не меньшей $1-\sigma$ алгоритм выдаст $\varepsilon$-оптимальное (по функции) решение исходной задачи. При этом алгоритм затрачивает в общей сложности

$$O\left(n+\frac{s\ln n\ln(n/\sigma)}{\varepsilon^2}\right)$$

элементарных арифметических операций (типа умножения двух чисел с плавающей точкой), что в случае $s\ll n$ заметно лучше всех известных сейчас оптимизационных аналогов. Слагаемое $O(n)$ “отражает” те случаи, когда $k=2n+1$, а также “стоимость” заключительного шага 7. Таким образом, имеет место следующее утверждение, доказательство которого вынесено в приложение.

**Теорема 4.1.1.** *Алгоритм Григориадиса–Хачияна после*

$$T_{\varepsilon,\sigma,n} = 12\left(\ln(2n+1) + \ln\left(\sigma^{-1}\right)\right)\varepsilon^{-2}$$

*итераций с вероятностью, не меньшей* $1-\sigma$ *выдает такое* $p$*, что* $0 \le f(p) \le \varepsilon$*. При этом затрачивается в общей сложности*

$$\mathrm{O}\left(n + \frac{s \ln n \ln(n/\sigma)}{\varepsilon^2}\right)$$

*операций (вида умножения двух чисел типа double).*

Ранее мы уже отмечали, что улучшить зависимость сложности от $\varepsilon$ не представляется возможным, существенно не ухудшая зависимость сложности от $n$. Можно даже сказать точнее, что для методов, в которых используется рандомизация с дисперсией стохастического градиента $D$, зависимость числа итераций $\sim D\varepsilon^{-2}\ln\left(\sigma^{-1}\right)$ типична и не улучшаема; тем не менее, стоит оговориться, что если от сильно выпуклого функционала берется математическое ожидание (это не наш случай, мы, напротив, избавлялись от этого, чтобы функционал был более информативен), то сходимость может быть улучшена до $\sim D\varepsilon^{-1}$ (см. [164, 251], см. также раздел 2.1 главы 2 и раздел 4.3 этой главы 4). Например, если использовать (детерминированный) быстрый градиентный метод (см. [93, 271]), то можно получить зависимость сложности от $\varepsilon$ вида $\varepsilon^{-1}$, но при этом число операций увеличится не менее чем в $n$ раз. Поскольку для таких задач вполне естественным является соотношение $\varepsilon \gg n^{-1}$, то выгода от этого представляется сомнительной. Отметим также (см. [120]), что в классе детерминированных алгоритмов зависимость от $n$ не может быть лучше, чем $\sim n^2$ (речь идет не о разреженных матрицах $P$).

Если бы мы применили обычный метод зеркального спуска (см. [91, 269]) для поиска (безусловного) минимума негладкой выпуклой функции $f(p) + \left(\sum_{k=1}^{n} p_k - 1\right)^2$, то используя теорему о субдиференциале максимума из [27, 85], и считая матрицу $A$ слабо заполненной, со средним показателем заполненности (по строкам и столбцам) $\chi \ll 1$ (в частности, число элементов матрицы отличных от нуля будет $\le \chi n^2$), получили бы оценку общей сложности в среднем

$$\mathrm{O}\left(\chi n^2 + \frac{\left(n + \chi^2 n^2\right)\ln n}{\varepsilon^2}\right).$$

В определенных ситуациях (например, при $\chi \sim n^{-1/2}$) эта оценка оказывается вполне конкурентноспособной. Эта идея развивается в разделе 4.4 этой главы 4.

Отметим связь описанного алгоритма с онлайн оптимизацией. Грубо говоря, такого типа алгоритмы асимптотически наиболее эффективны в задачах онлайн обучения на основе опыта экспертов даже при сопротивляющейся “Природе” (см. [27, 54, 91, 125, 164, 226, 250, 251, 269]). Подробнее этот вопрос будет рассмотрен в разделе 6.1 главы 6.

Отметим также связь алгоритма Григориадиса–Хачияна с концепцией ограниченной рациональности в контексте Discrete choice theory (см. [54, 137, 250]).

Резюмируем полученные в этом пункте результаты. Описанный выше алгоритм Григориадиса–Хачияна фактически соответствуют методу зеркального спуска для (4.1.1) с кососимметричной матрицей $A$. Для не кососимметричной матрицы см. [54].

Отличие рассмотренного зеркального спуска от зеркального спуска, например, работы [225] в том, что в нашем подходе (также как и в [120] движение итерационного процесса для поиска седловой точки $\langle u, Ax\rangle$ на произведении двух единичных симплексов осуществляется только по $x$, в то время как обычный метод зеркального спуска (рассчитанный не на кососимметричную матрицу $A$) осуществляет движение также и по $u$. Кроме того, очень важно, как именно осуществлять рандомизацию. В этом разделе, также как и в [120], рандомизация происходит при выборе компоненты вектора $x$, по которой осуществляется покомпонентное движение, во всех других работах рандомизацию предлагается вводить на этапе вычисления $Ax = E\left[A^{\langle x\rangle}\right]$, где $x$ – дискретная с.в. с распределением $x$, а $A^{\langle x\rangle}$ есть $x$-столбец матрицы $A$. Оба описанных способа рандомизации требуют одинакового по порядку числа шагов

$$O\left(\frac{\ln(n/\sigma)}{\varepsilon^2}\right)$$

для достижения точности (по функции) $\varepsilon$, но рандомизация Григориадиса–Хачияна позволяет заметно лучше учитывать разреженную специфику. Причина проста – на каждом шаге рандомизированных зеркальных спусков работ [88, 225] требуется обновлять вектор $x$, прибавляя к нему вектор с $2n+1$ ненулевыми компонентами, в независимости от разреженности $A$, стало быть, тратить не менее $2n+1$ операций на один шаг. В то время как у алгоритма Григориадиса–Хачияна в разреженном случае число операций на шаг оказывается в типичной ситуации, как мы видели выше, $O(s\ln n)$, что может быть заметно меньше $n$ в случае $s \ll n$.

В последующих разделах развивается идея о том, как организовывать случайный покомпонентный градиентный спуск для задач выпуклой оптимизации в пространствах огромной размерности, чтобы как можно сильнее учесть разреженную специфику задачи.

В связи с написанным выше интересно также было бы исследовать вопрос о робастном оценивании вектора PageRank (см. [229, 321]) и развить некоторые идеи, связанные с распределенными вычислениями вектора PageRank, см. работы

R. Tempo (http://staff.polito.it/roberto.tempo/).

### 4.1.6 Доказательство теоремы 4.1.1

Докажем теорему 4.1.1, сформулированную выше. Далее (для простоты будем вместо $2n+1$ писать $n$). Определим $r>0$ из соотношения: $\sigma \simeq n^{-r}$. Покажем, во многом следуя работе [120], что алгоритм Григориадиса–Хачияна после

$$T_{\varepsilon,\sigma,n} = 3\left(\ln n + \ln\left(\sigma^{-1}\right)\right)\varepsilon^{-2} = 3(1+r)\varepsilon^{-2}\ln n$$

итераций выдает такое $x$, что с вероятностью не меньшей чем $1-\sigma$, имеет место неравенство

$$\mathrm{A}x \le \varepsilon e\,.$$

Сначала, следуя работам [27, 120, 250], положим

$$p_i(t) = P_i(t)\left(\sum_{j=1}^{n} P_j(t)\right)^{-1},\quad P_i(t) = \exp\left(\varepsilon U_i(t)/2\right)$$

и

$$\Phi(t) = \sum_{i=1}^{n} P_i(t)\text{, где } U(t) = \mathrm{A}X(t).$$

Далее, аналогично [120], имеем

$$\Phi(t+1) = \sum_{i=1}^{n} P_i(t)\exp\left(\varepsilon a_{ik}/2\right) =$$

$$= \Phi(t)\sum_{i=1}^{n} p_i(t)\exp\left(\varepsilon a_{ik}/2\right),$$

$$E\left[\Phi(t+1)\middle| P(t)\right] = \Phi(t)\sum_{i,k=1}^{n} p_i(t)p_k(t)\exp\left(\varepsilon a_{ik}/2\right)$$

и ($|a_{ik}| \le 1$)

$$\exp\left(\varepsilon a_{ik}/2\right) \le 1 + \varepsilon a_{ik}/2 + \varepsilon^2/6,$$

но поскольку

$$\sum_{i,k=1}^{n} p_i(t)\,p_k(t) = \left(\sum_{i=1}^{n} p_i(t)\right)^2 = 1,\quad \sum_{i,k=1}^{n} p_i(t)\,p_k(t)\frac{\varepsilon^2}{6} = \frac{\varepsilon^2}{6},$$

$$\sum_{i,k=1}^{n} p_i(t)\,p_k(t)\,a_{ik} = \left\langle p(t), \mathrm{A}p(t)\right\rangle = 0,$$

то

$$E\left[\Phi(t+1)\middle|\vec{P}(t)\right]\le\Phi(t)\left(1+\varepsilon^2/6\right),$$

$$E\left[\Phi(t+1)\right]\le E\left[\Phi(t)\right]\left(1+\varepsilon^2/6\right).$$

Используя это неравенство и то, что

$$E\left[\Phi(0)\right]=\Phi(0)=n,$$

имеем

$$E\left[\Phi(t)\right]\le n\left(1+\varepsilon^2/6\right)^t.$$

Следовательно,

$$E\left[\Phi(t)\right]\le n\exp\left(t\varepsilon^2/6\right) \text{ и } E\left[\Phi\left(t^*\right)\right]\le n^{3/2+r/2}.$$

Отсюда по неравенству Маркова

$$\forall \text{ с.в. } \xi\ge 0, t>0\to P(\xi\ge t)\le E\xi/t,$$

имеем

$$P\left(\Phi\left(t^*\right)\le n^{3/2(1+r)}\right)=1-P\left(\Phi\left(t^*\right)\ge n^{3/2(1+r)}\right)\ge 1-E\left[\Phi\left(t^*\right)\right]\Big/n^{3/2(1+r)}\ge 1-\sigma.$$

Тогда, логарифмируя обе части неравенства

$$\exp\left(\varepsilon U_i\left(t^*\right)/2\right)=P_i\left(t^*\right)\le\sum_{i=1}^n P_i\left(t^*\right)=\Phi\left(t^*\right)\le n^{3/2(1+r)},$$

имеющего место с вероятностью не меньшей, чем $1-\sigma$, получим

$$P\left(\varepsilon U_i\left(t^*\right)/2\le 3/2(1+r)\ln n,\ i=1,\dots,n\right)\ge 1-\sigma,$$

$$P\left(U_i\left(t^*\right)/\left(3(1+r)\varepsilon^{-2}\ln n\right)\le\varepsilon,\ i=1,\dots,n\right)\ge 1-\sigma.$$

Откуда уже следует то, что требуется

$$P\left(\mathrm{A}x\left(t^*\right)\le\varepsilon e\right)\ge 1-\sigma.$$

## 4.2 Эффективные численные методы решения задачи PageRank для дважды разреженных матриц

### 4.2.1 Введение

В данном разделе мы также сконцентрируемся на решении задачи PageRank и ее окрестностях. Хорошо известно (Брин–Пейдж, 1998 [162, 241]), что задача ранжирования web-страниц приводит к поиску вектора Фробениуса–Перрона стохастической матрицы $x^T P = x^T$. Размеры матрицы могут быть колоссальными (в современных реалиях это сто миллиардов на сто миллиардов). Выгрузить такую матрицу в оперативную память обычного компьютера не представляется возможным. Задачу PageRank можно переписать, как задачу выпуклой оптимизации (см., например, [46, 88, 264, 273], а также раздел 2.1 главы 2 и раздел 4.1 этой главы 4) в разных вариантах: минимизация квадратичной формы $\|Ax-b\|_2^2$ и минимизация бесконечной нормы $\|Ax-b\|_\infty$ на единичном симплексе. Здесь $A = P^T - I$, $I$ – единичная матрица. К аналогичным задачам приводят задачи анализа данных (Ridge regression, LASSO [215]), задачи восстановления матрицы корреспонденций по замерам трафика на линках в компьютерных сетях [330] и многие другие задачи. Особенностью постановок этих задач, также как и в случае задачи PageRank, являются колоссальные размеры.

В данном разделе планируется сосредоточиться на изучении роли разреженности матрицы $A$, на использовании рандомизированных подходов и на специфики множества, на котором происходит оптимизация. Симплекс является (в некотором смысле) наилучшим возможным множеством, которое может порождать (независимо от разреженности матрицы $A$) разреженность решения (см., например, п. 3.3 [165], это тесно связано с l1-оптимизацией [168]).

Все о чем здесь написано хорошо известно на таком уровне грубости. Поясним, дополнительно погружаясь в детали, в чем заключаются отличия развиваемых в этом разделе подходов от известных подходов. Прежде всего, мы вводим специальный класс *дважды разреженных матриц* (одновременно разреженных и по строкам и по столбцам).[1] Такие матрицы, например, возникают в методе конечных элементов, и в более общем контексте при изучении разностных схем.[2] Если считать, что матрица $A$ имеет размеры $n$ на $n$, а число элементов в каждой строке и столбце не больше чем $s \ll n$, то число ненулевых элементов в матрице может быть $sn$. Кажется, что это произведение точно должно возникать в оценках общей сложности (числа арифметических операций, типа умножения или сложения двух чисел типа double) решения задачи (с определенной фиксированной точностью). Оказывается, что для первой постановки (минимизация квадратичной формы на симплексе) сложность может быть сделана пропорциональна $s^2$ (подраздел 4.2.2), а для второй (минимизация бесконечной нормы) и того меньше $s\ln n$ (подраздел 4.2.3).

---

[1] Заметим, что за счет специального “раздутия” изначальной матрицы можно обобщить приведенные в статье оценки и подходы с класса дважды разреженных матриц, на матрицы, у которых есть небольшое число полных строк и/или столбцов. Такие задачи встречаются в приложениях заметно чаще.

[2] Применительно к ранжированию web-страниц, это свойство означает, что входящие и выходящие степени всех web-страниц равномерно ограничены.

Первая задача может решаться обычным прямым градиентным методом с 1-нормой в прямом пространстве [44, 135] (см. также раздел 2.1 главы 2) – нетривиальным тут является, в том числе, организация работы с памятью (этими вопросами занимался А.С. Аникин). В частности, требуется хранение градиента в массиве, обновление элементов которого и вычисления максимального/минимального элемента должно осуществляться за время, не зависящие от размера массива (то есть от $n$). Тут оказываются полезными фибоначчиевы и бродалевские кучи [79]. Основная идея – не рассчитывать на каждом шаге градиент функции заново $A^T A x_{k+1}$, а пересчитывать его, учитывая, что $x_{k+1} = x_k + e_k$, где вектор $e_k$ состоит в основном из нулей:

$$A^T A x_{k+1} = A^T A x_k + A^T A e_k .$$

Аналогичную оценку общей сложности в этом разделе планируется получить методом условного градиента Франк–Вульфа [203], большой интерес к которому появился в последние несколько лет в основном в связи с задачами Big Data ([214, 220, 263]). Аккуратный анализ работы этого метода также при правильной организации работы с памятью позволяет (аналогично предыдущему подходу) находить такой $x$, что $\|Ax - b\|_2 \le \varepsilon$ за число арифметических операций

$$O\left(n + s^2 \ln\left(2 + n/s^2\right)/\varepsilon^2\right),$$

причем оценка оказывается вполне практической. Этому будет посвящен подраздел 4.2.2.

Вторая задача (минимизации $\|Ax - b\|_\infty$ на единичном симплексе) с помощью техники п. 6.5.2.3 [227] может быть переписана как седловая задача на произведение двух симплексов. Если применять рандомизированный метод зеркального спуска (основная идея которого заключается в вычислении вместо градиента $Ax$ стохастического градиента, рассчитываемого по следующему правилу: с вероятностью $x_1$ он равен первому столбцу матрицы $A$, с вероятностью $x_2$ второму и т.д.), то для того, чтобы обеспечить с вероятностью $\ge 1 - \sigma$ неравенство $\|Ax - b\|_\infty \le \varepsilon$ достаточно выполнить

$$O\left(n\ln\left(n/\sigma\right)/\varepsilon^2\right)$$

арифметических операций (см. раздел 2.1 главы 2). К сожалению, этот подход не позволяет в полном объеме учесть разреженность матрицы $A$. В этом разделе предлагается другой путь, восходящий к работе Григориадиса–Хачияна, 1995 [211] (см. также [54]). В основе этого подхода лежит рандомизация не при вычислении градиента, а при последующем проектировании (согласно $KL$-расстоянию, что отвечает экспоненциальному взвешиванию) градиента на симплекс. Предлагается вместо "честной" проекции на симплекс случайно выбирать одну из вершин симплекса (разреженный объект!) так, чтобы математическое ожидание проекции равнялось бы настоящей проекции [54]. В этом разделе планируется показать, что это можно эффективно делать, используя специальный двоичные деревья пересчета компонент градиента [46]. В результате планируется получить (отличным от раздела 4.1 данной главы 4 образом) следующую оценку для дважды разреженных матриц

$$O\left(n + s\ln n\ln\left(n/\sigma\right)/\varepsilon^2\right).$$

Напомним, что при этом число ненулевых элементов в матрице $A$ может быть $sn$. Возникает много вопросов относительно того насколько все это практично. К сожалению, на данный момент наши теоретические оценки говорят о том, что константа в $O(\ )$ может иметь порядок $10^2$. Подробнее об этом будет написано в подразделе 4.2.3.

В данном разделе при специальных предположениях относительно матрицы PageRank $P$ (дважды разреженная) мы предлагаем методы, которые работают по наилучшим известным сейчас оценкам для задачи PageRank в условиях отсутствия контроля спектральной щели матрицы PageRank [46] (см. также раздел 4.1 этой главы 4). А именно, полученная в разделе оценка (см. подраздел 4.2.2)

$$O\left(n + \ln\left(2 + n/s^2\right)/\varepsilon^2\right)$$

сложности поиска такого $x$, что $\|Ax - b\|_2 \le \varepsilon$, является на данный момент наилучшей при $s \ll \sqrt{n}$ для данного класса задач. Быстрее может работать только метод MCMC (для матрицы $P$ с одинаковыми по строкам внедиагональными элементами)

$$O\left(\ln n\ln\left(n/\sigma\right)/\left(\alpha\varepsilon^2\right)\right),$$

требующий, чтобы спектральная щель $\alpha$ была достаточно большой [46] (см. также раздел 4.1 этой главы 4).

Приведенные в разделе подходы, применимы не только к квадратным матрицам $A$ и не только к задаче PageRank. В частности, можно обобщать приведенные в разделе подхо-ды на композитные постановки [266]. Тем не менее, свойство двойной разреженности матрица $A$ является существенным. Если это свойство не выполняется, и имеет место, скажем, только разреженность в среднем по столбцам, то приведенные в разделе подходы могут быть доминируемы рандомизированными покомпонентными спусками. Если матрица $A$ сильно вытянута по числу строк (такого рода постановки характерны для задач анализа данных), то покомпонентные методы применяются к двойственной задаче [307, 331], если по числу столбцов (такого рода постановки характерны для задач поиска равновесных конфигураций в больших транспортных/компьютерных сетях), то покомпонентные методы применяются к прямой задаче [199, 295] (см. также раздел 5.1 главы 5). Подчеркнем, что при условии двойной разреженности $A$ с $s \ll \sqrt{n}$ нам не известны более эффективные (приведенных в данной разделе) способы решения описанных задач. В частности, это относится и к упомянутым выше покомпонентным спускам.

### 4.2.2 Метод условного градиента Франк–Вульфа

Перепишем задачу поиска вектора PageRank следующим образом

$$f(x) = \frac{1}{2}\|Ax\|_2^2 \to \min_{x \in S_n(1)},$$

где $S_n(1)$ – единичный симплекс в $\mathbb{R}^n$. Для решения этой задачи будем использовать метод условного градиента Франк–Вульфа [203, 214, 220, 263] (ранее мы уже встречались с этим методом в разделе 1.5 главы 1). Напомним, в чем состоит метод.

Выберем одну из вершин симплекса и возьмем точку старта $x_1$ в этой вершине. Далее по индукции, шаг которой имеет следующий вид. Решаем задачу

$$\langle \nabla f(x_k), y \rangle \to \min_{y \in S_n(1)} .$$

Обозначим решение этой задачи через

$$y_k = (0, \dots, 0, 1, 0, \dots, 0),$$

где 1 стоит на позиции

$$i_k = \arg\min_{i=1,\dots,n} \partial f(x_k) / \partial x^i .$$

Положим

$$x_{k+1} = (1 - \gamma_k) x_k + \gamma_k y_k, \; \gamma_k = \frac{2}{k+1}, \; k = 1, 2, \dots,$$

Имеет место следующая оценка [203, 214, 220, 263]

$$f(x_N) = f(x_N) - f_* \le \frac{2 L_p R_p^2}{N+1},$$

где

$$R_p^2 \le \max_{x,y\in S_n(1)} \|y-x\|_p^2,\ L_p \le \max_{\|h\|_p\le 1} \langle h, A^T A h\rangle = \max_{\|h\|_p\le 1} \|Ah\|_2^2,\ 1\le p\le\infty .$$

С учетом того, что оптимизация происходит на симплексе, мы выберем $p=1$. Не сложно показать, что этот выбор оптимален. В результате получим, что $R_1^2 = 4$,

$$L_1 = \max_{i=1,\dots,n} \left\| A^{\langle i\rangle} \right\|_2^2 \le 2 .$$

Таким образом, чтобы $f(x_N) \le \varepsilon^2/2$ достаточно сделать $N = 32\varepsilon^{-2}$ итераций. В действительности, число $32$ можно почти на порядок уменьшить немного более тонкими рассуждениями (детали мы вынуждены здесь опустить).

Сделав дополнительные вычисления стоимостью $O(n)$, можно так организовать процедуру пересчета $\nabla f(x_k)$ и вычисления $\arg\min_{i=1,\dots,n} \partial f(x_k)/\partial x^i$, что каждая итерация будет иметь сложность $O\left(s^2 \ln\left(2+n/s^2\right)\right)$. Для этого вводим (отметим, что А.С. Аникиным была продемонстрирована эффективность этой конструкции на практике)

$$\beta_k = \prod_{r=1}^{k-1}(1-\gamma_r),\ z_k = x_k/\beta_k,\ \ \tilde{\gamma}_k = \gamma_k/\beta_{k+1} .$$

Тогда итерационный процесс можно переписать следующим образом

$$z_{k+1} = z_k + \tilde{\gamma}_k y_k .$$

Пересчитывать $A^T A z_{k+1}$ при известном значении $A^T A z_k$ можно с помощью бинарной кучи за $O\left(s^2 \ln\left(2+n/s^2\right)\right)$. Далее, задачу

$$i_{k+1} = \arg\min_{i=1,\dots,n} \partial f(x_{k+1})/\partial x^i$$

можно переписать как

$$i_{k+1} = \arg\min_{i=1,\dots,n} \left(A^T A z_{k+1}\right)^i .$$

Поиск $i_{k+1}$ можно также осуществить за $O\left(s^2 \ln\left(2+n/s^2\right)\right)$. Определив в конечном итоге $z_N$, мы можем определить $x_N$, затратив дополнительно не более

$$O(n+\ln N) = O(n)$$

арифметических операций. Таким образом, итоговая оценка сложности описанного метода будет

$$O\left(n + s^2 \ln\left(2+n/s^2\right)/\varepsilon^2\right).$$

Стоит отметить, что функционал, выбранный в этом примере, обеспечивает намного лучшую оценку $\|Ax\|_2 \le \varepsilon$ по сравнению с функционалом из подраздела 4.2.3, который обеспечивает лишь $\|Ax\|_\infty \le \varepsilon$. Наилучшая (в разреженном случае без, условий на спектральную щель матрицы $P$ [46]) из известных нам на данный момент оценок

$$\mathrm{O}\left(n + s\ln n \ln\left(n/\sigma\right)/\varepsilon^2\right)$$

(см. подраздел 4.2.3) для $\|Ax\|_\infty$ может быть улучшена приведенной в этом разделе оценкой, поскольку $\|Ax\|_2$ может быть (и так часто бывает) в $\sim\sqrt{n}$ раз больше $\|Ax\|_\infty$, а $s \ll \sqrt{n}$.

### 4.2.3 Седловое представление задачи PageRank и рандомизированный зеркальный спуск в форме Григориадиаса–Хачияна

Перепишем задачу поиска вектора PageRank следующим образом

$$f(x) = \|Ax\|_\infty \to \min_{x \in S_n(1)}.$$

Следуя [227], эту задачу можно переписать седловым образом

$$\min_{x \in S_n(1)} \max_{\|y\|_1 \le 1} \langle Ax, y \rangle.$$

В свою очередь последнюю задачу можно переписать как

$$\min_{x \in S_n(1)} \max_{\omega \in S_{2n}(1)} \langle Ax, J\omega \rangle,$$

где

$$J = \|I_n, -I_n\|.$$

В итоге исходную задачу можно переписать, с сохранением свойства разреженности, как

$$\min_{x \in S_n(1)} \max_{\omega \in S_{2n}(1)} \langle \omega, \tilde{A}x \rangle.$$

Следуя [54] (см. также раздел 6.1 главы 6), опишем эффективный и интерпретируемый способ решения этой задачи. Пусть есть два игрока А и Б. Задана матрица (антагонистической) игры $\tilde{A} = \|\tilde{a}_{ij}\|$, где $|\tilde{a}_{ij}| \le 1$, $\tilde{a}_{ij}$ – выигрыш игрока А (проигрыш игрока Б) в случае когда игрок А выбрал стратегию $i$, а игрок Б стратегию $j$. Отождествим себя с игроком Б. И предположим, что игра повторяется $N \gg 1$ раз (это число может быть заранее неизвестно [54], однако для простоты изложения будем считать это число известным). Введем функцию потерь (игрока Б) на шаге $k$

$$f_k(x) = \langle \omega^k, \tilde{A}x \rangle, \; x \in S_n(1),$$

где $\omega^k \in S_{2n}(1)$ – вектор (вообще говоря, зависящий от всей истории игры до текущего момента включительно, в частности, как-то зависящий и от текущей стратегии (не хода) игрока Б, заданной распределением вероятностей (результат текущего разыгрывания (ход Б) игроку А не известен)) со всеми компонентами равными 0, кроме одной компоненты, соответствующей ходу А на шаге $k$, равной 1. Хотя функция $f_k(x)$ определена на единичном симплексе, по “правилам игры” вектор $x^k$ имеет ровно одну единичную компоненту, соответствующую ходу Б на шаге $k$, остальные компоненты равны нулю. Обозначим цену игры (в нашем случае $C=0$)

$$C = \max_{\omega\in S_{2n}(1)} \min_{x\in S_n(1)} \left\langle \omega, \tilde{A}x \right\rangle = \min_{x\in S_n(1)} \max_{\omega\in S_{2n}(1)} \left\langle \omega, \tilde{A}x \right\rangle.$$ (теорема фон Неймана о минимаксе)

Пару векторов $(\omega, x)$, доставляющих решение этой минимаксной задачи (т.е. седловую точку), назовем равновесием Нэша. По определению (это неравенство восходит к Ханнану [250])

$$\min_{x\in S_n(1)} \frac{1}{N}\sum_{k=1}^{N} f_k(x) \le C.$$

Тогда, если мы (игрок Б) будем придерживаться следующей рандомизированной стратегии (см., например, [46, 54, 211]), выбирая $\{x^k\}$:

1. $p^1 = \left(n^{-1},...,n^{-1}\right)$;
2. Независимо разыгрываем случайную величину $j(k)$ такую, что $P\left(j(k)=j\right) = p_j^k$;
3. Полагаем $x_{j(k)}^k = 1$, $x_j^k = 0$, $j \ne j(k)$;
4. Пересчитываем

$$p_j^{k+1} \sim p_j^k \exp\left(-\sqrt{\frac{2\ln n}{N}}\tilde{a}_{i(k)j}\right),\ j=1,...,n,$$

где $i(k)$ – номер стратегии, которую выбрал игрок А на шаге $k$;[3]

[3] Заметим, что эта стратегия имеет естественную интерпретацию. Мы (игрок Б) описываем на текущем шаге игрока А вектором, компоненты которого – сколько раз игрок А использовал до настоящего момента соответствующую стратегию. Согласно этому вектору частот мы рассчитываем вектор своих ожидаемых потерь (при условии, что игрок А будет действовать согласно этому частотному вектору). Далее, вместо того, чтобы выбирать наилучшую (для данного вектора частот А) стратегию, дающую наименьшие потери, мы используем распределение Гиббса с параметром $\sqrt{2\ln n/N}$ (экспоненциально взвешиваем с этим параметром вектор ожидаемых потерь с учетом знака). С наибольшей вероятностью будет выбрана наилучшая стратегия, но с ненулевыми вероятностями могут быть выбраны и другие стратегии.

то с вероятностью $\geq 1-\sigma$

$$\frac{1}{N}\sum_{k=1}^{N} f_k\left(x^k\right)-\min_{x\in S_n(1)}\frac{1}{N}\sum_{k=1}^{N} f_k\left(x\right)\leq\sqrt{\frac{2}{N}}\left(\sqrt{\ln n}+2\sqrt{2\ln\left(\sigma^{-1}\right)}\right),$$

т.е. с вероятностью $\geq 1-\sigma$ наши (игрока Б) потери ограничены

$$\frac{1}{N}\sum_{k=1}^{N} f_k\left(x^k\right)\leq C+\sqrt{\frac{2}{N}}\left(\sqrt{\ln n}+2\sqrt{2\ln\left(\sigma^{-1}\right)}\right).$$

Самый плохой для нас случай (с точки зрения такой оценки) – это когда игрок А тоже действует (выбирая $\left\{\omega^k\right\}$) согласно аналогичной стратегии (только игрок А решает задачу максимизации).[4] Если и А и Б будут придерживаться таких стратегий, то они сойдутся к равновесию Нэша (седловой точке), причем довольно быстро [54]: с вероятностью $\geq 1-\sigma$

$$0\leq\left\|A\bar{x}^N\right\|_\infty=\max_{\omega\in S_{2n}(1)}\left\langle\omega,\tilde{A}\bar{x}^N\right\rangle-\max_{\omega\in S_{2n}(1)}\min_{x\in S_n(1)}\left\langle\omega,\tilde{A}x\right\rangle\leq\max_{\omega\in S_{2n}(1)}\left\langle\omega,\tilde{A}\bar{x}^N\right\rangle-\min_{x\in S_n(1)}\left\langle\bar{\omega}^N,\tilde{A}x\right\rangle\leq$$

$$\leq\max_{\omega\in S_{2n}(1)}\left\langle\omega,\tilde{A}\bar{x}^N\right\rangle-\frac{1}{N}\sum_{k=1}^{N}\left\langle\omega^k,\tilde{A}x^k\right\rangle+\frac{1}{N}\sum_{k=1}^{N}\left\langle\omega^k,\tilde{A}x^k\right\rangle-\min_{x\in S_n(1)}\left\langle\bar{\omega}^N,\tilde{A}x\right\rangle\leq$$

$$\leq\sqrt{\frac{2}{N}}\left(\sqrt{\ln(2n)}+2\sqrt{2\ln(2/\sigma)}\right)+\sqrt{\frac{2}{N}}\left(\sqrt{\ln n}+2\sqrt{2\ln(2/\sigma)}\right)\leq$$

$$\leq 2\sqrt{\frac{2}{N}}\left(\sqrt{\ln(2n)}+2\sqrt{2\ln(2/\sigma)}\right),$$

где

$$\bar{x}^N=\frac{1}{N}\sum_{k=1}^{N}x^k,\quad \bar{\omega}^N=\frac{1}{N}\sum_{k=1}^{N}\omega^k.$$

Таким образом, чтобы с вероятностью $\geq 1-\sigma$ иметь $\left\|A\bar{x}^N\right\|_\infty\leq\varepsilon$, достаточно сделать

$$N=16\frac{\ln(2n)+8\ln(2/\sigma)}{\varepsilon^2}\text{ – итераций;}$$

$$O\left(n+\frac{s\ln n\left(\ln n+\ln\left(\sigma^{-1}\right)\right)}{\varepsilon^2}\right)\text{ – общее число арифметических операций,}$$

где число $s\ll\sqrt{n}$ – ограничивает сверху, число ненулевых элементов в строках и столбцах матрицы $P$.

---

[4] Игрок А пересчитывает

$$p_i^{k+1}\sim p_i^k\exp\left(\sqrt{\frac{2\ln(2n)}{N}}\tilde{a}_{ij(k)}\right),\ i=1,\ldots,2n,$$

где $j(k)$ – номер стратегии, которую выбрал игрок Б на шаге $k$.

Нетривиальным местом здесь является оценка сложности одной итерации $\text{O}(s\ln n)$. Получается эта оценка исходя из оценки наиболее тяжелых действий шаг 2 и 4. Сложность тут в том, что чтобы сгенерировать распределение дискретной случайной величины, принимающей $n$ различных значений (в общем случае) требуется $\text{O}(n)$ арифметических операций – для первого генерирования (приготовления памяти). Последующие генерирования могут быть сделаны за эффективнее – за $\text{O}(\ln n)$. Однако в нашем случае есть специфика, заключающаяся в том, что при переходе на следующий шаг $s$ вероятностей в распределении могли как-то поменяться. Если не нормировать распределение вероятностей, то можно считать, что остальные вероятности остались неизменными. Оказывается, что такая специфика позволяет вместо оценки $\text{O}(n)$ получить оценку $\text{O}(s\ln n)$. Практически это было реализовано А.С. Аникиным. Эксперименты подтвердили справедливость выписанных оценок.

**Замечание 4.2.1 (см. также аналогичные замечания в разделе 4.1 этой главы 4).** Опишем точнее эту процедуру. У нас есть сбалансированное двоичное дерево высоты $O(\log_2 n)$ с $n$ листом (дабы не вдаваться в технические детали, считаем что число $n$ – есть степень двойки, понятно, что это не так, но можно взять, скажем, наименьшее натуральное $m$ такое, что $2^m > n$ и рассматривать дерево с $2^m$ листом) и с $\text{O}(n)$ общим числом вершин. Каждая вершина дерева (отличная от листа) имеет неотрицательный вес равный сумме весов двух ее потомков. Первоначальная процедура приготовления дерева, отвечающая одинаковым весам листьев, потребует $\text{O}(n)$ операций. Но сделать это придется всего один раз. Процедура генерирования дискретной случайной величины с распределением, с точностью до нормирующего множителя, соответствующим весам листьев может быть осуществлена с помощью случайного блуждания из корня дерева к одному из листьев. Отметим, что поскольку дерево двоичное, то прохождение каждой его вершины при случайном блуждании, из которой идут два ребра в вершины с весами $a > 0$ и $b > 0$, осуществляется путем подбрасывания монетки (“приготовленной” так, что вероятность пойти в вершину с весом $a$ – есть $a/(a+b)$). Понятно, что для осуществления этой процедуры нет необходимости в дополнительном условии нормировки: $a+b=1$. Если вес какого-то листа алгоритму необходимо поменять по ходу работы, то придется должным образом поменять дополнительно веса тех и только тех вершин, которые лежат на пути из корня к этому листу. Это необходимо делать, чтобы поддерживать свойство: каждая вершина дерева (отличная от листа) имеет вес равный сумме весов двух ее потомков.

Итак, выше в этом разделе была описана процедура генерирования последовательности $x^k$ со сложностью

$$O\left(n + s\ln n \ln\left(n/\sigma\right)/\varepsilon^2\right),$$

на основе которой можно построить такой частотный вектор

$$\bar{x}^N = \frac{1}{N}\sum_{k=1}^{N} x^k ,$$

компоненты – частоты, с которыми мы использовали соответствующие стратегии, что $\left\|A\bar{x}^N\right\|_\infty \le \varepsilon$ с вероятностью $\ge 1-\sigma$.

А.С. Аникиным и Д.И. Камзоловым [5] проводились численные эксперименты с описанными в этом и предыдущем разделе 4.1 методами. Наилучшие результаты (в условиях малой спектральной щели) показал метод условного градиента Франк–Вульфа в описанной в подразделе 4.2.2 модификации. К сожалению, метод из подраздела 4.2.3 в целом оказался не столь эффективным на практике, как это можно было ожидать из приведенных оценок.

## 4.3 О неускоренных эффективных методах решения разреженных задач квадратичной оптимизации

### 4.3.1 Введение

В работе [5] (см. также [44], раздел 2.1 главы 2 и раздел 4.2 этой главы 4) были предложены два различных способа решения задачи минимизации (на единичном симплексе) квадратичного функционала специального вида, связанного с задачей PageRank [46]. Первый способ (восходящий к Ю.Е. Нестерову) базировался на обычном (неускоренном) прямом градиентном методе в l1-норме [135]. Второй способ базировался на методе условного градиента [220] (Франк–Вульф). Общей особенностью обоих способов является возможность учитывать разреженность задачи в большей степени, чем известные альтернативные способы решения отмеченной задачи. За счет такого свойства, итоговые оценки времени работы упомянутых методов в ряде случаев получаются лучше, чем оценки для ускоренных методов и их различных вариантов [44, 44] (см. также раздел 5.1 главы 5). Естественно было задаться вопросами: на какие классы задач и на какие простые множества (входящие в ограничения) можно перенести отмеченные алгоритмы? Настоящий раздел имеет целью частично ответить на поставленные вопросы.

В подразделе 4.3.2 мы показываем (во многом следуя изначальным идеям Ю.Е. Нестерова [264, 273]), что прямой градиентный метод в l1-норме из [5] (см. также раздел 2.1 главы 2) можно распространить на общие задачи квадратичной минимизации. При этом можно перейти от задач на симплексе к задачам на всем пространстве (это даже заметно упростит сам метод и доказательство оценок). Но наиболее интересно то, что удается распространить этот метод на общие аффинно-сепарабельные задачи с сепарабельными композитами [44]. Такие задачи часто возникают при моделировании компьютерных и транспортных сетей [33, 40, 138] (см. также главу 1), а также в анализе данных [223]. Основными конкурирующими методами здесь являются различные варианты (прямые, двойственные) ускоренных покомпонентных методов [44, 138, 199, 264, 295], которые чаще интерпретируются как различные варианты методов рандомизации суммы [44, 129, 223, 232, 234, 240, 242, 248, 307,331] (см. также раздел 5.1 главы 5). Таким образом, в подразделе 4.3.2 предлагается новый подход к этим задачам.

В подразделе 4.3.3 мы распространяем метод Франк–Вульфа из [5] (см. также раздел 4.2 этой главы 4) на более общий класс задач квадратичной минимизации. К сожалению, выйти за пределы описанного в подразделе 4.3.3 класса задач нам не удалось. Тем не менее, стоит отметить интересную (саму по себе) конструкцию, позволившую перенести метод с единичного симплекса на неотрицательный ортант. Стоит также отметить, что эта конструкция существенным образом опиралась на ряд новых идей, успешно использовавшихся и в других контекстах.

В подразделе 4.3.4 мы предлагаем новый рандомизированный метод (дополнительно предполагая сильную выпуклость функционала задачи), который работает быстрее неускоренного покомпонентного метода, и позволяет для ряда задач в большей степени воспользоваться их разреженностью, чем позволяют это сделать ускоренные покомпонентные методы.

“За кадром” описываемых далее сюжетов стоит задача разработки эффективных численных методов решения системы линейных уравнений $Ax = b$ в пространстве огромной размерности (от десятков миллионов переменных и выше). При этом качество решения оценивается согласно $\|Ax - b\|_2$. Полученные в этом разделе результаты позволяют надеяться, что предложенные в подразделах 4.3.2–4.3.4 методы будут доминировать остальные в случае, когда (приводим наиболее важные условия):

1) Матрица $A$ такова, что в каждом столбце и каждой строке не более $s \ll \sqrt{n}$ отличных от нуля элементов.
2) Решение $x^*$ системы $Ax = b$ таково, что $\|x^*\|_1 \approx \|x^*\|_2$ ну или точнее: всегда имеет место неравенство
$$\|x^*\|_2 \le \|x^*\|_1 \le \sqrt{n}\|x^*\|_2,$$
которое (в нашем случае) можно уточнить следующим образом
$$\|x^*\|_2 \le \|x^*\|_1 \ll \sqrt{n}\|x^*\|_2.$$

Сравнительному анализу предложенных в данном разделе 4.3 методов с альтернативными методами поиска решения системы $Ax = b$ (в частности с популярными сейчас ускоренными покомпонентными методами [44]) посвящено дополнение.

### 4.3.2 Прямой неускоренный градиентный метод в l1-норме

Рассматривается следующая задача
$$f(x) = \frac{1}{2}\langle Ax, x\rangle - \langle b, x\rangle \to \min_{x\in\mathbb{R}^n}, \tag{4.3.1}$$
где квадратная матрица $A$ предполагается симметричной неотрицательно определенной. Также предполагается, что в каждом столбце и каждой строке матрицы $A$ не более $s \ll n$ элементов отлично от нуля. Мы хотим решить задачу (4.3.1) с точностью $\varepsilon$
$$f(x^N) - f(x^*) \le \varepsilon.$$
Заметим, что тогда

$$\left\|Ax^N - b\right\|_2 \le \sqrt{2\lambda_{\max}(A)\varepsilon}.$$

Эту задачу предлагается решать обычным градиентным методом, но не в евклидовой норме, а в l1-норме [135] (см. также раздел 2.1 главы 2) или (что то же самое для данной задачи) неускоренным покомпонентным методом (см. раздел 5.1 главы 5) с выбором максимальной компоненты [264] (если $i^k$ определяется не единственным образом, то можно выбрать любого представителя):

$$x^{k+1} = x^k + \arg\min_{h \in \mathbb{R}^n} \left\{ f\left(x^k\right) + \left\langle \nabla f\left(x^k\right), h \right\rangle + \frac{L}{2}\|h\|_1^2 \right\} =$$

$$= \begin{cases} x_i^k, & i \ne i^k = \arg\max\limits_{j=1,\dots,n} \left|\partial f\left(x^k\right)/\partial x_j\right| \\ x_i^k - \dfrac{1}{L}\dfrac{\partial f\left(x^k\right)}{\partial x_{i^k}}, & i = i^k \end{cases}, \tag{4.3.2}$$

где $L = \max\limits_{i,j=1,\dots,n} \left|A_{ij}\right|$. Точку старта $x^0$ итерационного процесса выберем таким образом, чтобы только одна из компонент была отлична от нуля. Любая итерация такого метода (за исключением самой первой – первая требует $\mathrm{O}(n)$ арифметических операций) может быть осуществлена за $\mathrm{O}(s\log_2 n)$. Действительно, $x^{k+1}$ и $x^k$ отличаются в одной компоненте. Следовательно,

$$\nabla f\left(x^{k+1}\right) = Ax^{k+1} - b = Ax^k - b + \frac{1}{L}\frac{\partial f\left(x^k\right)}{\partial x_{i^k}} Ae_{i^k} = \nabla f\left(x^k\right) + \frac{1}{L}\frac{\partial f\left(x^k\right)}{\partial x_{i^k}} Ae_{i^k},$$

где

$$e_i = (\overbrace{\underbrace{0,\dots,0,1}_{i},0,\dots,0}^{n}).$$

Таким образом, по условию задачи (матрица $A$ разрежена) $\nabla f\left(x^{k+1}\right)$ отличается от $\nabla f\left(x^k\right)$ не более чем в $s$ компонентах. Следовательно, $i^{k+1}$ можно пересчитать за $\mathrm{O}(s\log_2 n)$, храня массив $\left\{\left|\partial f\left(x^k\right)/\partial x_j\right|\right\}_{j=1}^n$ в виде специальной структуры данных (кучи) для поддержания максимального элемента [5, 273].

Согласно [135] необходимое число итераций

$$N = \frac{2\max\limits_{i,j=1,\dots,n}\left|A_{ij}\right| R_1^2}{\varepsilon},$$

где

$$R_1^2 = \sup\left\{\left\|x - x^*\right\|_1^2 :\ f(x) \le f\left(x^0\right)\right\}.$$

Таким образом, общее число арифметических операций (время работы метода) можно оценить следующим образом

$$\mathrm{O}\left(n + s\log_2 n \frac{\max\limits_{i,j=1,\dots,n}\left|A_{ij}\right| R_1^2}{\varepsilon}\right),$$

Рассмотрим теперь задачу вида

$$f(x) = \sum_{k=1}^{m} f_k\left(A_k^T x\right) + \sum_{k=1}^{n} g_k\left(x_k\right) \to \min_{x\in\mathbb{R}^n},$$

где все функции (скалярного аргумента) $f_k$, $g_k$ имеют равномерно ограниченные числом $L$ первые две производные,[5] а матрица $A = \left[A_1,\dots,A_m\right]^T$ такова, что в каждом ее столбце не больше $s_m \ll m$ ненулевых элементов, а в каждой строке – не больше $s_n \ll n$. Описанный выше подход (формула (4.3.2)) позволяет решить задачу за время

$$\mathrm{O}\left(n + s_n s_m \log_2 n \frac{L\max\limits_{i=1,\dots,n}\left\|A^{\langle i\rangle}\right\|_2^2 R_1^2}{\varepsilon}\right),$$

где $A^{\langle i\rangle}$ – $i$-й столбец матрицы $A$.

Собственно, задача из [5] (см. также пример 2.1.2 раздела 2.1 главы 2) является частным случаем приведенной выше общей постановки (в смысле выбора функционала):

$$f(x) = \frac{1}{2}\left\|Ax\right\|_2^2 + \frac{\gamma}{2}\sum_{i=1}^{n}\left(-x_i\right)_+^2 \to \min_{\langle x,e\rangle=1}.$$

В принципе, можно перенести результаты, полученные в этом пункте на случай, когда оптимизация происходит на симплексе. Это частично (но не полностью [5]) решает основную проблему данного метода: большое значение $R_1^2$, даже когда $\left\|x^*\right\|_1^2$ небольшое. Однако мы не будем здесь этого делать.

### 4.3.3 Метод условного градиента (Франк–Вульфа)

Рассмотрим теперь следующую задачу[6]

---

[5] $g_k''$ можно равномерно ограничивать числом $L\max\limits_{i=1,\dots,n}\left\|A^{\langle i\rangle}\right\|_2^2$.

[6] Отличие от задачи (4.3.1) в том, что рассматривается неотрицательный ортант вместо всего пространства – раздутием исходного пространства в два раза к такому ограничению можно прийти из задачи оптимизации на всем пространстве.

$$f\left(x\right)=\frac{1}{2}\left\langle Ax,x\right\rangle-\left\langle b,x\right\rangle\to\min_{x\in\mathbb{R}^n_+},\qquad(4.3.3)$$

где квадратная матрица $A$ предполагается симметричной неотрицательно определенной. Также предполагаем, что в каждом столбце и каждой строке матрицы $A$ не более $s\ll n$ элементов отлично от нуля, и в векторе $b$ не более $s$ элементов отлично от нуля.

Предположим сначала, что мы знаем такой $R$, что решение задачи (4.3.3) удовлетворяет условию

$$x^*\in\bar{S}_n\left(R\right)=\left\{x\in\mathbb{R}^n_+:\ \sum_{i=1}^n x_i^*\le R\right\}.$$

Выберем одну из вершин $\bar{S}_n\left(R\right)$ и возьмем точку старта $x^1$ в этой вершине. Далее будем действовать по индукции, шаг которой имеет следующий вид.

Решаем задачу

$$\left\langle\nabla f\left(x^k\right),y\right\rangle\to\min_{y\in\bar{S}_n\left(R\right)}\ .\qquad(4.3.4)$$

Введем (если $i_k$ определяется не единственным образом, то можно выбрать любого представителя)

$$i_k=\arg\min_{i=1,\ldots,n}\partial f\left(x^k\right)\big/\partial x_i\,.$$

Обозначим решение задачи (4.3.4) через

$$y^k=\begin{cases}R\cdot e_{i_k},\ \text{если}\ \partial f\left(x^k\right)\big/\partial x_{i_k}<0\\ 0,\ \text{если}\ \partial f\left(x^k\right)\big/\partial x_{i_k}\ge 0\end{cases}.$$

Положим

$$x^{k+1}=\left(1-\gamma_k\right)x^k+\gamma_k y^k,\ \gamma_k=\frac{2}{k+1},\ k=1,2,\ldots,$$

Имеет место следующая оценка [214, 220, 263]

$$f\left(x^N\right)-f\left(x^*\right)\le f\left(x^N\right)-\max_{k=1,\ldots,N}\left\{f\left(x^k\right)+\left\langle\nabla f\left(x^k\right),y^k-x^k\right\rangle\right\}\le\frac{2L_pR_p^2}{N+1},\qquad(4.3.5)$$

где

$$R_p^2=\max_{x,y\in\bar{S}_n\left(R\right)}\left\|y-x\right\|_p^2,\ L_p=\max_{\left\|h\right\|_p\le 1}\left\langle h,Ah\right\rangle,\ 1\le p\le\infty\,.$$

С учетом того, что оптимизация происходит на $\bar{S}_n\left(R\right)$, мы выбираем $p=1$. Несложно показать, что этот выбор оптимален. В результате получим, что

$$R_1^2=4R^2,\ L_1=\max_{i,j=1,\ldots,n}\left|A_{ij}\right|.\qquad(4.3.6)$$

Таким образом, чтобы

$$f\left(x^N\right)-f\left(x^*\right)\le\varepsilon,$$

достаточно сделать

$$N=\frac{8\max\limits_{i,j=1,\ldots,n}\left|A_{ij}\right|R^2}{\varepsilon}$$

итераций.

Сделав один раз дополнительные вычисления стоимостью $\mathrm{O}(n)$, можно так организовать процедуру пересчета $\nabla f\left(x^k\right)$ и вычисление $i_k$, что каждая итерация будет иметь сложность $\mathrm{O}(s\log_2 n)$. Для этого вводим

$$\beta_k=\prod_{r=1}^{k-1}\left(1-\gamma_r\right),\; z^k=x^k/\beta_k,\quad \tilde{\gamma}_k=\gamma_k/\beta_{k+1}.$$

Тогда итерационный процесс можно переписать следующим образом

$$z^{k+1}=z^k+\tilde{\gamma}_k y^k.$$

Пересчитывать $Az^{k+1}-b$ при известном значении $Az^k$ можно за $\mathrm{O}(s)$. Далее, задачу

$$i_{k+1}=\arg\min_{i=1,\ldots,n}\partial f\left(x^{k+1}\right)/\partial x_i$$

можно переписать как

$$i_{k+1}=\arg\min_{i=1,\ldots,n}\left(\left[Az^{k+1}\right]_i-\frac{b_i}{\beta_{k+1}}\right).$$

Поиск $i_{k+1}$ можно осуществить за $\mathrm{O}(s\log_2 n)$ (см. подраздел 4.3.2). Определив в конечном итоге $z^N$, мы можем определить $x^N$, затратив дополнительно не более

$$\mathrm{O}(n+\ln N)=\mathrm{O}(n)$$

арифметических операций. Таким образом, итоговая оценка сложности описанного метода будет

$$\mathrm{O}\left(n+s\log_2 n\frac{\max\limits_{i,j=1,\ldots,n}\left|A_{ij}\right|R_1^2}{\varepsilon}\right),$$

где $R_1^2=\left\|x^*\right\|_1^2$ (следует сравнить с аналогичной формулой из подраздела 4.3.2).

Отметим, что при этом мы можем пересчитывать (это не надо делать, если известно значение $f\left(x^*\right)$, например, для задачи (4.3.3) естественно рассматривать постановки с $f\left(x^*\right)=0$)

$$f\left(x^{k+1}\right)+\left\langle\nabla f\left(x^{k+1}\right),y^{k+1}-x^{k+1}\right\rangle$$

также за $\mathrm{O}(s\log_2 n)$. Следуя [214], можно немного упростить приведенные рассуждения за счет небольшого увеличения числа итераций. А именно, можно использовать оценки (при этом следует полагать $\gamma^k = 2(k+2)^{-1}$)

$$f(x^k) - f(x^*) \le \langle \nabla f(x^k), x^k - y^k \rangle,$$

$$\min_{k=1,\ldots,N} \langle \nabla f(x^k), x^k - y^k \rangle \le \frac{7L_1R_1^2}{N+2}.$$

Вернемся теперь к предположению, что изначально известно $R$. На практике, как правило, даже если мы можем как-то оценить $R$, то оценка получается слишком завышенной. Используя формулы (4.3.5), (4.3.6) мы можем воспользоваться следующей процедурой рестартов по параметру $R$.

Выберем сначала $R = R^0 = 1$. Делаем предписанное этому $R$ число итераций и проверяем (без дополнительных затрат) критерий останова (все необходимые вычисления уже были сделаны по ходу итерационного процесса)

$$f(x^N) - \max_{k=1,\ldots,N}\left\{ f(x^k) + \langle \nabla f(x^k), y^k - x^k \rangle \right\} \le \frac{8\max\limits_{i,j=1,\ldots,n}|A_{ij}|R^2}{N+1} \le \varepsilon.$$

Если он выполняется, то мы угадали и получили решение. Если не выполняется, то полагаем $R := \chi R$ ($\chi > 1$), и повторяем рассуждения. Остановимся поподробнее на вопросе оптимального выбора параметра $\chi$ [37]. Обозначим $R^* = \|x^*\|_1$. Пусть

$$R^0\chi^{r-1} < R^* \le R^0\chi^r.$$

Общее число итераций будет пропорционально

$$1 + \chi^2 + \chi^4 + \chi^{2r} = \frac{\chi^{2r+2} - 1}{\chi^2 - 1} \le \frac{\chi^4}{\chi^2 - 1}\left(\frac{R^*}{R^0}\right)^2.$$

Выберем $\chi = \sqrt{2}$, исходя из минимизации правой части этого неравенства. При этом общее число итераций возрастет не более чем в четыре раза по сравнению со случаем, когда значение $R^*$ известно заранее. Следует сравнить эту конструкцию, с аналогичной конструкцией из раздела 3.1 главы 3.

Описанный выше подход распространяется и на задачи

$$f(x) = \frac{1}{2}\|Ax - b\|_2^2 \to \min_{x \in \mathbb{R}_+^n},$$

Матрица $A$ такова, что в каждом ее столбце не больше $s_m \ll m$ ненулевых элементов, а в каждой строке – не больше $s_n \ll n$. В векторе $b$ не более $s_m$ элементов отлично от нуля.

Описанный выше подход (на базе метода Франк–Вульфа) позволяет решить задачу за время

$$\mathrm{O}\left(n+s_n s_m \log_2 n \frac{\max\limits_{i=1,\dots,n}\left\|A^{\langle i\rangle}\right\|_2^2 R_1^2}{\varepsilon}\right),$$

где $R_1^2=\left\|x^*\right\|_1^2$ (следует сравнить с аналогичной формулой из подраздела 4.3.2, полученной для более общего класса задач).[7]

В связи написанным выше в этом пункте, заметим, что задача может быть не разрежена, но свойство разреженности появляется в решении при использовании метода Франк–Вульфа, что также может заметно сокращать объем вычислений. Интересные примеры имеются в работах п. 3.3 [163], [277].

В последнее время методы условного градиента переживают бурное развитие в связи с многочисленными приложениями к задачам машинного обучения. В связи с этим появились интересные обобщения классического метода Франк–Вульфа. Отметим, например, работы [235, 240]. Интересно было бы понять: возможно ли перенести (а если возможно, то как именно и в какой степени) результаты этого пункта на более общий класс задач (чем класс задач с квадратичным функционалом) и на более общий класс методов?

### 4.3.4 Неускоренный рандомизированный градиентный спуск в сильно выпуклом случае

Рассмотрим для большей наглядности снова задачу (4.3.1). Будем считать, что $f(x)$ $\mu$-сильно выпуклый функционал в 2-норме, т.е. $\lambda_{\min}(A)\ge\mu>0$. Сочетая методы [6] и [236], введем следующий рандомизированный метод

$$x^{k+1}=x^k-\frac{2}{\mu\cdot(k+1)}\left\|\nabla f\left(x^k\right)\right\|_1 \operatorname{sign}\left(\frac{\partial f\left(x^k\right)}{\partial x_{i(x^k)}}\right)e_{i(x^k)},\ x^1=0,$$

где

[7] Здесь также как и в задаче (4.3.3) требуются рестарты по неизвестному параметру $R_1$, который явно используется в методе в качестве размера симплекса, фигурирующего при решении вспомогательной задачи ЛП на каждой итерации. Однако, как было продемонстрировано выше, все это приводит к увеличению числа итераций не более чем в четыре раза.

$$e_i = (\overbrace{\underbrace{0,\dots,0,1}_{i},0,\dots,0}^{n}),$$

$$i(x^k) = i \text{ с вероятностью } \frac{1}{\left\|\nabla f(x^k)\right\|_1}\left|\frac{\partial f(x^k)}{\partial x_i}\right|, \; i = 1,\dots,n.$$

Тогда[8] [236]

$$E\left[f(y^k)\right] - f(x^*) \le \frac{\frac{2}{N}\sum_{k=1}^{N} E\left[\left\|\nabla f(x^k)\right\|_1^2\right]}{\mu\cdot(N+1)} \le \frac{2\max\limits_{k=1,\dots,N} E\left[\left\|\nabla f(x^k)\right\|_1^2\right]}{\mu\cdot(N+1)},$$

$$E\left[\left\|\nabla f(y^k)\right\|_2^2\right] \le \frac{4L\max\limits_{k=1,\dots,N} E\left[\left\|\nabla f(x^k)\right\|_1^2\right]}{\mu\cdot(N+1)},$$

где $L = \lambda_{\max}(A)$, а

$$y^N = \frac{2}{N\cdot(N+1)}\sum_{k=1}^{N} k\cdot x^k.$$

На базе описанного метода построим новый метод, который “следит” за последовательностью $\sum_{k=1}^{N} k\cdot x^k$, пересчитывая (с некоторой частотой)

$$\left\|\nabla f\left(\frac{2}{N\cdot(N+1)}\sum_{k=1}^{N} k\cdot x^k\right)\right\|_1^2.$$

Метод “дожидается” момента[9] $N = O(nL/\mu)$, когда

$$\left\|\nabla f(y^N)\right\|_1^2 \le \frac{1}{2}\left\|\nabla f(x^1)\right\|_1^2.$$

и перезапускается с $x^1 := y^N$. Можно показать, что необходимое число таких перезапусков для достижения точности $\varepsilon$ (по функции) будет $\sim \log_2(\varepsilon^{-1})$. Более того, подобно [44], можно так организовать описанную выше процедуру, чтобы попутно получить и оценки вероятностей больших отклонений (детали мы вынуждены здесь опустить).

Предположим теперь, что $\nabla f(x + he_i)$ отличается от $\nabla f(x)$ (при произвольных $x$, $h$ и $e_i$) не более чем в $s$ компонентах ($s \ll n$).[10] Для задачи (4.3.1) это имеет место (мож-

---

[8] Второе неравенство может быть достаточно грубым [45].

[9] Это довольно грубая оценка на $N$, поскольку использующееся при ее получении правое неравенство

$$1 \le \left\|\nabla f(y^N)\right\|_1^2 \Big/ \left\|\nabla f(y^N)\right\|_2^2 \le n$$

может быть сильно завышенным.

но также говорить и о задаче (4.3.3), с очевидной модификацией описанного метода – все оценки сохраняются). Можно так организовать процедуру выбора момента $N$ (с сохранением свойства $N = O(nL/\mu)$), что амортизационная (средняя) сложность итерация метода (с учетом затрат на проверку условий остановки на каждом перезапуске) будет $O(s\log_2 n)$ [5, 46]. Казалось бы, что мы ничего не выиграли по сравнению с обычными (неускоренным) покомпонентными методами [45]. Количество итераций и стоимость одной итерации одинаковы для обоих методов (в категориях $O(\ )$ с точностью до логарифмических множителей). В действительности, ожидается, что предложенный нами алгоритм будет работать заметно быстрее (неускоренного) покомпонентного метода, поскольку, как уже отмечалось, при получении этих оценок мы пару раз использовали потенциально довольно грубые неравенства. С другой стороны пока мы говорили только о задаче (4.3.1) (в разреженном случае), для которой эффективно работают ускоренные покомпонентные методы [45] (см. также раздел 5.1 главы 5) с такой же оценкой стоимости одной итерации, но заметно лучшей оценкой для числа итераций $N = O\left(n\sqrt{L/\mu}\right)$ (с точностью до логарифмического множителя). В связи с этим может показаться, что предложенный в этом пункте метод теоретически полностью доминируем ускоренными покомпонентными методами. На самом деле, это не совсем так. Хорошо известно, см., например, [45], что существующие сейчас всевозможные модификации ускоренных покомпонетных методов, которые могут сполна учитывать разреженность задачи (что проявляется в оценке стоимости итерации $O(s)$), применимы лишь к специальному классу задач [45] (см. также раздел 5.1 главы 5). Для общих задач стоимость одной итерации ускоренного покомпонентного метода будет $O(n)$ независимо от разреженности задачи. Например, на данный момент не известны такие модификации ускоренных покомпонентных методов, которые бы позволяли сполна учитывать разреженность в задаче [6, 37, 45]

$$f(x) = \ln\left(\sum_{i=1}^{m}\exp\left(A_k^T x\right)\right) - \langle x, b\rangle + \frac{\mu}{2}\|x\|_2^2 \to \min_{x\in\mathbb{R}^n}.$$

---

[10] Это условие можно обобщить с сохранением оценок, например, на случай, когда $f(x) = g(x) + h(x)$, и существует такая (скалярная) функция $\alpha(x, he_i) > 0$, что $\nabla f(x + he_i) = \nabla g(x + he_i) + \nabla h(x + he_i)$ отличается от $\alpha(x, he_i)\cdot\nabla g(x) + \nabla h(x)$ не более чем в $s$ компонентах. Для этого приходится перезаписать исходный метод: отличие в том, что теперь вводится двойная рандомизация (рандомизация согласно вектору $\nabla g(x)$ и независимая рандомизация согласно вектору $\nabla h(x)$).

В то время как описанный выше метод распространим и на эту задачу с оценкой амортизационной стоимости одной итерации $\mathrm{O}\left(s_m s_n \log_2 n\right)$.

### 4.3.5 Дополнение

В этом подразделе мы напомним некоторые полезные факты о поиске решения (псевдорешения) системы линейных уравнений $Ax = b$ [110, 258]. К близким задачам и методам приводит изучение проекции точки на аффинное множество $Ax = b$ [6, 138].

Хорошо известно, что задача поиска такого вектора $x$, что $Ax = b$, может быть эффективно полиномиально решена даже в концепции битовой сложности (например, простейшим методом Гаусса или более современными методами внутренней точки). Однако требуется, учитывая специфику задачи, так подобрать метод, чтобы можно было искать решения систем огромных размеров. Несмотря на более чем двухвековую историю, эта область математики (эффективные численные методы решения систем линейных уравнений) до сих пор бурно развивается. О чем, например, говорит недавний доклад Д. Спильмана на международном математическом конгрессе [315].

Пусть наблюдается вектор

$$b = Ax + \varepsilon,$$

где матрица $A$ – известна, $\varepsilon_k \in N\left(0, \sigma^2\right)$ – независимые одинаково распределенные случайные величины, $k = 1, ..., m$ (ненаблюдаемые). Оптимальная оценка[11] неизвестного вектора параметров $x$ определяется решением задачи

$$f\left(x\right) = \frac{1}{2}\left\|Ax - b\right\|_2^2 \to \min_{x \in \mathbb{R}^n}. \qquad (4.3.7)$$

Введем псевдообратную матрицу (Мура–Пенроуза) $A^{\dagger} = \left(A^T A\right)^{-1} A^T$, если столбцы матрицы $A$ линейно независимы (первый случай) и $A^{\dagger} = A^T\left(AA^T\right)^{-1}$, если строки матрицы $A$ линейно независимы (второй случай). В первом случае $x^* = A^{\dagger} b$ – единственное решение задачи (4.3.7), во втором случае $x^* = A^{\dagger} b$ – решение задачи (4.3.7) с наименьшим значением 2-нормы.[12] Вектор $x^*$ называют псевдорешением задачи $Ax = b$.

---

[11] То есть несмещенная и с равномерно (по $x$) минимальной дисперсией. При такой (статистической) интерпретации также стоит предполагать, что $n \le m$.

[12] В приложениях во втором случае чаще ищут не псевдорешение, а вводят в функционал регуляризатор [138], например, исходя из байесовского подхода (дополнительного априорного вероятностного предположения относительно вектора неизвестных параметров) или исходя из желания получить наиболее разреженное решение. В таком случае многое (но не все) из того, что написано в данном разделе удается сохранить

Поиск решения системы $Ax = b$ был сведен выше к решению задачи (4.3.7). Для задачи выпуклой оптимизации (4.3.7) существуют различные эффективные численные методы. Например, метод сопряженных градиентов, сходящийся со скоростью (стоимость одной итерации этого метода равна $O(sn)$ – числу ненулевых элементов в матрице $A$)

$$\frac{1}{2}\left\|Ax^N - b\right\|_2^2 = f\left(x^N\right) - f\left(x^*\right) \le \varepsilon,$$

$$N = O\left(\sqrt{\frac{\lambda_{\max}\left(A^T A\right)\left\|x^*\right\|_2^2}{\varepsilon}}\right),$$

где $N \le n$. На практике при не очень больших размерах матрицы $A$ ($\max\{m,n\} \le 10^6$) эффективно работают квазиньютоновские методы типа L-BFGS [285]. А при совсем небольших размерах ($\max\{m,n\} \le 10^3 - 10^4$) и методы внутренней точки [258]. Однако нам интересны ситуации, когда $\max\{m,n\} \gg 10^6$. В таких случаях эффективно использовать прямые и двойственные ускоренные покомпонентные методы для задачи (4.3.7) (см. раздел 3.2 главы 3 и раздел 5.1 главы 5). Эти методы дают следующие оценки [44] (одна итерация у обоих методов в среднем требует $O(s)$ арифметических операций)[13]

$$E\left[\frac{1}{2}\left\|Ax^N - b\right\|_2^2\right] \le \varepsilon, \qquad \text{(прямой метод)}$$

$$N = O\left(n\sqrt{\frac{\bar{L}_x R_x^2}{\varepsilon}}\right), \quad \bar{L}_x^{1/2} = \frac{1}{n}\sum_{k=1}^{n}\left\|A^{\langle k\rangle}\right\|_2, \quad R_x = \left\|x^*\right\|_2,$$

$$E\left[\left\|Ax^N - b\right\|_2\right] \le \varepsilon, \qquad \text{(двойственный метод)}$$

$$N = O\left(n\sqrt{\frac{\bar{L}_y R_y}{\varepsilon}}\right), \quad \bar{L}_y^{1/2} = \frac{1}{m}\sum_{k=1}^{m}\left\|A_k\right\|_2, \quad R_y = \left\|y^*\right\|_2,$$

где $y^*$ – решение (если решение не единственно, то можно считать, что $y^*$ – решение с наименьшей 2-нормой) “двойственной” задачи

---

(должным образом модифицировав). При этом появляются и дополнительные новые возможности за счет появления новых степеней свободы в выборе регуляризатора [277].

[13] Далее для наглядности мы дополнительно предполагаем, что $f\left(x^*\right) = 0$. Если $f\left(x^*\right) = f^* > 0$, то из того что $f\left(x^N\right) - f^* \le \varepsilon$ будет следовать $\left\|Ax^N - b\right\|_2 - f^* \le \varepsilon / f^*$ (этот простое, но полезное наблюдение было сделано Дмитрием Островским (INRIA Гренобль)).

$$-\frac{1}{2}\left\|A^T y\right\|_2^2 + \langle b, y\rangle \to \max_{y\in\mathbb{R}^m}.$$

Напомним, что $(A_k)^T$ – $k$-я строка матрицы $A$, а $A^{\langle k\rangle}$ – $k$-й столбец.

Из написанного выше не очевидно, что для получения решения $x^*$ необходимо исходить именно из решения задачи оптимизации (4.3.7) и, как следствие, исходить из определяемого этой задачей критерия точности решения (2-норма невязки). Например, в этом разделе мы исходили в основном из задачи (4.3.1), которая хотя и похожа по свойствам на задачу (4.3.7), но все же приводит к отличному от (4.3.7) критерию точности решения. Вместо задачи (4.3.7) в определенных ситуациях вполне осмысленно бывает рассматривать даже негладкую (но выпуклую) задачу

$$f(x) = \|Ax - b\|_\infty \to \min_{x\in\mathbb{R}^n}.$$

Для такой задачи можно предложить специальным образом рандомизированный вариант метода зеркального спуска [285] (см. также следующий раздел), который гарантируют выполнение (считаем $f(x^*) = 0$)

$$E\left[\left\|Ax^N - b\right\|_\infty\right] \le \varepsilon,$$

за время

$$\mathrm{O}\left(s_m n + s_m \log_2 m \frac{\max_{k=1,\dots,m}\|A_k\|_1^2 R_x^2}{\varepsilon^2}\right),$$

если в каждом столбце матрицы $A$ не больше $s_m \le m$ ненулевых элементов. Отметим [3], что обычный (не рандомизированный) метод зеркального спуска гарантирует выполнение

$$\left\|Ax^N - b\right\|_\infty \le \varepsilon$$

за время

$$\mathrm{O}\left(m + n + s_n s_m \log_2 m \frac{\max_{k=1,\dots,m}\|A_k\|_2^2 R_x^2}{\varepsilon^2}\right),$$

если в каждом столбце матрицы $A$ не больше $s_m \le m$ ненулевых элементов, а в каждой строке – не больше $s_n \le n$.

Выше мы обсуждали, так называемые, вариационные численные методы, которые сводили поиск решения системы $Ax = b$, к решению задачи выпуклой оптимизации. Однако имеется большой класс итерационных численных методов (типа метода простой итерации), которые исходят из перезаписи системы $Ax = b$ в эквивалентном виде. Например, в виде $x = \tilde{A}x + b$, где $\tilde{A} = I - A$, и организации вычислений по формуле[14] $x^{k+1} = \tilde{A}x^k + b$ [97]. Для ряда таких методов также удается получить оценки скорости сходимости, причем не только для нормы невязки в системе, но и непосредственно для оценки невязки в самом решении $\left\|x^N - x^*\right\|_2$. К сожалению, такие оценки получаются весьма пессимистичными, причем эти оценки типично реализуются на практике (закон Мерфи). Поясним сказанное примером, восходящим к работе Красносельского–Крейна [103]. Пусть

$$x^{k+1} = \tilde{A}x^k + b,$$

где $\tilde{A} = I - A$ – положительно определенная матрица $n \times n$ с собственными числами $0 < \lambda_1 \le ... \le \lambda_n < 1$, а $x^0$ – выбирается равновероятно из $B_R(0)$ – 2-шара радиуса $R$ с центром в $0$. Обозначим через

$$\delta^k = x^{k+1} - x^k = \tilde{A}x^k + b - x^k = b - Ax^k$$

контролируемую невязку. В качестве критерия останова итерационного процесса выберем момент $N$ первого попадания $\delta^N \in B_\varepsilon(0)$. Заметим, что

$$\left\|x^k - x^*\right\|_2 = \left\|\left(I - \tilde{A}\right)^{-1}\delta^k\right\|_2 = \left\|A^{-1}\delta^k\right\|_2.$$

Имеет место следующий результат. При $R \to \infty$ вероятность выполнения неравенств

$$0.999 \cdot \frac{\varepsilon\lambda_n}{1-\lambda_n} \le \left\|x^N - x^*\right\|_2 \le \frac{\varepsilon}{1-\lambda_n}$$

стремится к единице. Поскольку типично, что при больших $n$ число $\lambda_n$ близко к 1,то вряд ли можно рассчитывать на то, что $\left\|x^N - x^*\right\|_2$ мало, если установлена только малость $\left\|Ax^N - b\right\|_2$.

---

[14] Интересно заметить, что критерием сходимости такого итерационного процесса является: $\lambda_{\max}\left(\tilde{A}\right) < 1$ [110]. Однако, если при этом $\left\|\tilde{A}\right\|_{2,2} = \sigma_{\max}\left(\tilde{A}\right) = \lambda_{\max}\left(A^T A\right) > 1$ (к счастью, такие ситуации не типичны [3]), то существуют такие $x^0$ (причем на них довольно легко попасть), стартуя с которых, итерационная процедура на первых итерациях приводит к резкому росту нормы невязки $\left\|x^N - x^*\right\|_2$, т.е. к наличию “горба” у такой последовательности. Высота этого горба может быть (в зависимости от конкретной задачи) сколь угодно большой [102, 110]. Поскольку все вычисления проходят с конечной длиной мантиссы, то такой резкий рост может приводить к вычислительной неустойчивости процедуры, что неоднократно наблюдалось в численных экспериментах [110]. Если $\sigma_{\max}\left(\tilde{A}\right) < 1$, то таких проблем нет, и последовательность $\left\|x^N - x^*\right\|_2$ мажорируется геометрической прогрессией с основанием $\sigma_{\max}\left(\tilde{A}\right)$.

Отметим, что оценки скорости сходимости многих итерационных методов (типа простой итерации) исходят из первого метода Ляпунова, т.е. из спектральных свойств матрицы $A$. При плохих спектральных свойствах вводится регуляризация. По-сути, большинство таких методов можно проинтерпретировать как численный метод решения соответствующей вариационной переформулировки задачи. Регуляризация при этом имеет четкий смысл (восходящий к пионерским работам А.Н. Тихонова) – сделать функционал сильно выпуклым и использовать его сильную выпуклость при выборе метода и получении оценки скорости сходимости [300]. Однако вариационный подход позволяет сполна использовать и второй метод Ляпунова, что дает возможность распространить всю мощь современных численных методов выпуклой оптимизации на решение задачи $Ax = b$. В частности, использовать рандомизацию [285] или(и) свойства неравнозначности компонент в решении $x^*$ и разреженность $A$, которые использовались при разработке методов из данного раздела.

Ю.Е. Нестеровым недавно было отмечено, что возможные прорывы в разработке новых эффективных численных методов (в том числе применительно к решению системы $Ax = b$) можно ожидать от рандомизированных квазиньютоновских методов. К сожалению, здесь имеются большие проблемы, унаследованные от обычных (не рандомизированных) квазиньютоновских методов, с получением оценок скоростей сходимости, адекватных реальной (наблюдаемой на практике) скорости сходимости (все известные сейчас оценки никак не объясняют быстроту сходимости этих методов на практике). Тем не менее, недавно появилось несколько интересных работ Гаверса–Ричтарика в этом направлении [208, 209, 210].

## 4.4 Рандомизация и разреженность в задачах huge-scale оптимизации на примере работы метода зеркального спуска

### 4.4.1 Введение

В недавнем цикле работ Ю.Е. Нестерова с соавторами [5, 44, 273, 281] был введен класс задач huge-scale оптимизации (задачи выпуклой оптимизации, для которых размерность прямого и(или) двойственного пространства не меньше десятков миллионов), и исследовалась роль разреженности в таких задачах. Данный раздел посвящен изучению конкретного (но, пожалуй, наиболее важного) метода решения таких задач – метода зеркального спуска (МЗС). Это метод представляет собой обобщение метода проекции градиента. Он был предложен в конце 70-х годов прошлого века А.С. Немировским [259]. С тех пор метод получил повсеместное распространение для решения задач больших размерностей, прежде всего, в связи с идей рандомизации. Метод оказался слабо чувствительным к замене настоящего градиента его несмещенной оценкой. Это обстоятельство активно используется на практике, поскольку построить несмещенную оценку в ряде случаев удается намного дешевле, чем посчитать градиент. Как правило, выгода получается пропорциональной размерности пространства, в котором происходит оптимизация. Однако для разреженных задач тут возникают нюансы, связанные с проработкой правильного сочетания рандомизации и разреженности. Настоящий раздел посвящен изучению такого сочетания для МЗС. Данный раздел является продолжением работы [5] (см. разделы 4.1 и 4.2 этой главы 4), в которой уже был рассмотрен один сюжет на эту тему (метод Григориадиса–Хачияна для задачи PageRank).

Структура раздела 4.4 следующая. В подразделе 4.4.2 мы описываем рандомизированый МЗС. Стоит обратить внимание на случай, когда оптимизация происходит на неограниченном множестве, и при этом выводятся оценки вероятностей больших уклонений. Тут есть некоторые тонкости, проработка которых делает п. 2 не просто вводным материалом для последующего изложения, но и представляющим самостоятельный интерес. В подразделе 4.4.3 результаты подраздела 4.4.2 переносятся на задачи с функциональными ограничениями, на которые мы не умеем эффективно проектироваться (см. также раздел 3.4 главы 3). В детерминированном случае (когда используется обычный градиент) такого типа задачи рассматривались достаточно давно. Разработаны эффективные способы их редуцирования к задаче с простыми ограничениями. Предлагались различные эффективные методы (Поляк–Шор–Немировский–Нестеров). Однако в случае, когда вместо градиента мы используем его несмещенную оценку (для функционала и ограничений) нам не известны оценки, поэтому в подразделе 4.4.3 приводится соответствующее обобщение МЗС и устанавливаются необходимые в дальнейшем оценки. В подразделе 4.4.4 на основе теоретических заготовок подразделов 4.4.2, 4.4.3 мы описываем класс разреженных задач (обобщающих задачу PageRank), для которых удается за счет рандомизации получить дополнительную выгоду.

### 4.4.2 Рандомизированный метод зеркального спуска

Рассмотрим задачу выпуклой оптимизации

$$f(x) \to \min_{x \in Q}. \tag{4.4.1}$$

Под решением этой задачи будем понимать такой $\bar{x}^N \in Q \subseteq \mathbb{R}^n$, что с вероятностью $\geq 1-\sigma$ имеет место неравенство

$$f(\bar{x}^N) - f_* \leq \varepsilon,$$

где $f_* = f(x_*)$ – оптимальное значение функционала в задаче (1), $x_* \in Q$ – решение задачи (4.4.1). На каждой итерации $k = 1,...,N$ нам доступен стохастический (суб-)градиент $\nabla_x f_k(x, \xi^k)$ в одной, выбранной нами (методом), точке $x^k$.

Опишем метод зеркального спуска (МЗС) для решения задачи (4.4.1) (мы в основном будем следовать работам [135, 259]). Введем норму $\| \ \|$ в прямом пространстве (сопряженную норму будем обозначать $\| \ \|_*$) и прокс-функцию $d(x)$ сильно выпуклую относительно этой нормы, с константой сильной выпуклости $\geq 1$. Выберем точку старта

$$x^1 = \arg\min_{x \in Q} d(x),$$

считаем, что

$$d(x^1) = 0, \ \nabla d(x^1) = 0.$$

Введем брэгмановское “расстояние”

$$V_x(y) = d(y) - d(x) - \langle \nabla d(x), y - x \rangle.$$

Определим “размер” решения

$$d(x_*) = V_{x^1}(x_*) = R^2.$$

Определим оператор “проектирования” согласно этому расстоянию

$$\mathrm{Mirr}_{x^k}(\mathrm{v}) = \arg\min_{y \in Q} \left\{ \langle \mathrm{v}, y - x^k \rangle + V_{x^k}(y) \right\}.$$

МЗС для задачи (4.4.1) будет иметь вид, см., например, [135]

$$x^{k+1} = \mathrm{Mirr}_{x^k}\left( \alpha \nabla_x f_k(x^k, \xi^k) \right), \ k = 1,...,N.$$

Будем считать, что $\left\{\xi^k\right\}_{k=1}^N$ представляет собой последовательность независимых случайных величин, и для всех $x \in Q$ имеют место условия ($k = 1,...,N$)

4.4.1 $E_{\xi_k}\left[\nabla_x f_k\left(x,\xi^k\right)\right] = \nabla f\left(x\right)$;

4.4.2 $E_{\xi^k}\left[\left\|\nabla_x f_k\left(x,\xi^k\right)\right\|_*^2\right] \le M^2$.

В ряде случаев нам также понадобится более сильное условие

4.4.3 $\left\|\nabla_x f_k\left(x,\xi^k\right)\right\|_*^2 \le \tilde{M}^2$ почти наверное по $\xi^k$.

При выполнении условия 4.4.1 для любого $u \in Q$, $k = 1,...,N$ имеет место неравенство, см., например, [135]

$$\alpha\left\langle \nabla_x f_k\left(x^k,\xi^k\right), x^k - u\right\rangle \le \frac{\alpha^2}{2}\left\|\nabla_x f_k\left(x^k,\xi^k\right)\right\|_*^2 + V_{x^k}\left(u\right) - V_{x^{k+1}}\left(u\right).$$

Это неравенство несложно получить в случае евклидовой прокс-структуры [269]

$$d\left(x\right) = \left\|x\right\|_2^2/2,\ V_x\left(y\right) = \left\|y - x\right\|_2^2/2.$$

В этом случае МЗС для задачи (4.4.1) есть просто вариант обычного метода проекции градиента (см. примеры 4.4.1, 4.4.2 ниже).

Разделим сначала выписанное неравенство на $\alpha$ и возьмем условное математическое ожидание $E_{\xi^k}\left[\,\cdot\,\middle|\,\Xi^{k-1}\right]$ ($\Xi^{k-1}$ – сигма алгебра, порожденная $\xi^1$, …, $\xi^{k-1}$), затем просуммируем то, что получится по $k = 1,...,N$, используя условие 4.4.1. Затем возьмем от того, что получилось при суммировании, полное математическое ожидание, учитывая условие 4.4.2. В итоге, выбрав $u = x_*$, и определив

$$\bar{x}^N = \frac{1}{N}\sum_{k=1}^N x^k,$$

получим

$$N\cdot\left(E\left[f\left(\bar{x}^N\right)\right] - f_*\right) \le \frac{V_{x^1}\left(x_*\right)}{\alpha} - \frac{E\left[V_{x^{N+1}}\left(x_*\right)\right]}{\alpha} + \frac{1}{2}M^2\alpha N \le$$

$$\le \frac{R^2}{\alpha} + \frac{1}{2}M^2\alpha N.$$

Выбирая[15]

$$\alpha = \frac{R}{M}\sqrt{\frac{2}{N}},$$

[15] Можно получить и адаптивный вариант приводимой далее оценки, для этого потребуется использовать метод двойственных усреднений [135, 269]. Впрочем, в [259] имеется адаптивный вариант МЗС.

получим

$$E\left[ f\left(\bar{x}^N\right)\right] - f_* \le MR\sqrt{\frac{2}{N}} . \qquad (4.4.2)$$

Заметим, что в детерминированном случае вместо $\bar{x}^N$ можно брать

$$\breve{x}^N = \arg\min_{k=1,\dots,N} f\left(x^k\right).$$

Немного более аккуратные рассуждения (использующие неравенство Азума–Хефдинга) с

$$\alpha = \frac{R}{\tilde{M}}\sqrt{\frac{2}{N}}$$

позволяют уточнить оценку (4.4.2) следующим образом (см., например, [40, 227] и раздел 1.5 главы 1):

$$f\left(\bar{x}^N\right) - f_* \le \tilde{M}\sqrt{\frac{2}{N}}\left(R + 2\tilde{R}\sqrt{\ln\left(2/\sigma\right)}\right) \qquad (4.4.3)$$

с вероятностью $\ge 1-\sigma$, где

$$\tilde{R} = \sup_{x\in\tilde{Q}}\left\|x - x_*\right\|,$$

$$\tilde{Q} = \left\{x \in Q :\ \left\|x - x_*\right\|^2 \le 65R^2\ln\left(4N/\sigma\right)\right\}.$$

Собственно, для справедливости оценки (4.4.3) достаточно требовать выполнение условий 4.4.1 – 4.4.3 лишь на множестве $\tilde{Q} \subseteq Q$. Это замечание существенно, когда рассматриваются неограниченные множества $Q$ (см., например, подраздел 4.4.4).

Оценки (4.4.2), (4.4.3) являются неулучшаемыми с точностью до мультипликативного числового множителя по $N$ и $\sigma$. Наряду с этим можно обеспечить и их неулучшаемость по размерности пространства $n$, путем “правильного” выбора прокс-функции [259] (такой выбор всегда возможен, и известен для многих важных в приложениях случаев выбора множества $Q$). Собственно прокс-структура (новая степень свободы по сравнению с классическим методом проекции градиента) и вводилась для того, чтобы была возможность обеспечить последнее свойство.

В виду того, что мы используем рандомизированный метод и всегда

$$f\left(\bar{x}^N\right) - f_* \ge 0,$$

то используя идею амплификации (широко распространенную в Computer Science) можно немного “улучшить” оценку (4.4.3). Для этого сначала перепишем оценку (4.4.2) в виде

$$E\left[ f\left(\bar{x}^N\right)\right] - f_* \le \varepsilon,$$

где (здесь $N$, конечно, должно быть натуральным числом, поэтому эту формулу и последующие формулы такого типа надо понимать с точностью до округления к наименьшему натуральному числу, большему написанного)

$$N = \frac{2M^2R^2}{\varepsilon^2}. \qquad (4.4.4)$$

Отсюда по неравенству Маркова

$$P\left(f\left(\bar{x}^N\right) - f_* \geq 2\varepsilon\right) \leq \frac{E\left[f\left(\bar{x}^N\right)\right] - f_*}{2\varepsilon} \leq \frac{1}{2}.$$

Можно параллельно (независимо) запустить $\log_2\left(\sigma^{-1}\right)$ траекторий метода. Обозначим $\bar{x}^N_{\min}$ тот из $\bar{x}^N$ на этих траекториях, который доставляет минимальное значение $f\left(\bar{x}^N\right)$. Из выписанного неравенства Маркова получаем, что имеет место неравенство

$$P\left(f\left(\bar{x}^N_{\min}\right) - f_* \geq 2\varepsilon\right) \leq \sigma.$$

Таким образом, можно не более чем за

$$N = \frac{8M^2R^2}{\varepsilon^2}\log_2\left(\sigma^{-1}\right)$$

обращений за стохастическим градиентом и не более чем за $\log_2\left(\sigma^{-1}\right)$ обращений за значением функции найти решение задачи (4.4.1) $\bar{x}^N$ с требуемой точностью $\varepsilon$ и доверительным уровнем $\sigma$.

Как уже отмечалось, во многих приложениях множество $Q$ неограниченно (см., например, подраздел 4.4.4). Поскольку $x_*$ априорно не известно, то это создает проблемы для определения $R$, которое входит в формулу для расчета шага метода

$$\alpha = \frac{R}{M}\sqrt{\frac{2}{N}}.$$

Однако если мы заранее выбрали желаемую точность $\varepsilon$, то с помощью формулы (4.4.4) можно выразить шаг следующим образом

$$\alpha = \frac{\varepsilon}{M^2}.$$

Рассмотрим три конкретных примера множества $Q$, которые нам понадобятся в дальнейшем (см. подраздел 4.4.4). В примерах 4.4.1, 4.4.2 мы не приводим оценки скорости сходимости, поскольку они будут иметь вид (4.4.2), (4.4.3), т.е. никакой уточняющей информации тут не появляется, в отличие от примера 4.4.3.

**Пример 4.4.1 (все пространство).** Предположим, что $Q = \mathbb{R}^n$. Выберем

$$\|\ \| = \|\ \|_2,\ d(x) = \frac{1}{2}\|x\|_2^2.$$

Тогда МЗС примет следующий вид ($\alpha = \varepsilon/M^2$, $x^1 = 0$):

$$x^{k+1} = x^k - \alpha \nabla_x f_k\left(x^k, \xi^k\right),\ k = 1, \dots, N. \ \square$$

**Пример 4.4.2 (неотрицательный ортант).** Предположим, что

$$Q = \mathbb{R}_+^n = \left\{x \in \mathbb{R}^n: \quad x \geq 0\right\}.$$

Выберем

$$\|\ \| = \|\ \|_2,\ d(x) = \frac{1}{2}\|x - \bar{x}\|_2^2,\ \bar{x} \in \operatorname{int} Q.$$

Тогда МЗС примет следующий вид ($\alpha = \varepsilon/M^2$, $x^1 = \bar{x}$):

$$x^{k+1} = \left[x^k - \alpha \nabla_x f_k\left(x^k, \xi^k\right)\right]_+ = \max\left\{x^k - \alpha \nabla_x f_k\left(x^k, \xi^k\right), 0\right\},\ k = 1, \dots, N,$$

где $\max\{\ \}$ берется покомпонентно. □

**Пример 4.4.3 (симплекс).** Предположим, что

$$Q = S_n(1) = \left\{x \geq 0: \quad \sum_{i=1}^n x_i = 1\right\}.$$

Выберем

$$\|\ \| = \|\ \|_1,\ d(x) = \ln n + \sum_{i=1}^n x_i \ln x_i.$$

Тогда МЗС примет следующий вид:

$$x_i^1 = 1/n,\ i = 1, \dots, n,$$

при $k = 1, \dots, N$, $i = 1, \dots, n$

$$x_i^{k+1} = \frac{\exp\left(-\sum_{r=1}^k \alpha \frac{\partial f_r\left(x^r, \xi^r\right)}{\partial x_i}\right)}{\sum_{l=1}^n \exp\left(-\sum_{r=1}^k \alpha \frac{\partial f_r\left(x^r, \xi^r\right)}{\partial x_l}\right)} = \frac{x_i^k \exp\left(-\alpha \frac{\partial f_k\left(x^k, \xi^k\right)}{\partial x_i}\right)}{\sum_{l=1}^n x_l^k \exp\left(-\alpha \frac{\partial f_k\left(x^k, \xi^k\right)}{\partial x_l}\right)}.$$

Оценки скорости сходимости будут иметь вид:

$$E\left[f\left(\bar{x}^N\right)\right] - f_* \leq M\sqrt{\frac{2\ln n}{N}} \ \text{(при } \alpha = M^{-1}\sqrt{2\ln n/N}\text{)};$$

$$f\left(\bar{x}^N\right) - f_* \leq \tilde{M}\sqrt{\frac{2}{N}}\left(\sqrt{\ln n} + 4\sqrt{\ln\left(\sigma^{-1}\right)}\right) \ \text{(при } \alpha = \tilde{M}^{-1}\sqrt{2\ln n/N}\text{)}$$

с вероятностью $\geq 1 - \sigma$. □

Представим себе, что задача (4.4.1) видоизменилась следующим образом

$$f(x) \to \min_{\substack{g(x)\le 0 \\ x\in Q}}.$$

Можно ли к ней применить изложенный выше подход, если “проектироваться” на множество (в отличие от $Q$)

$$\{x \in Q:\ g(x) \le 0\}$$

мы эффективно не умеем? В случае, когда мы знаем оптимальное значение $f_*$, то мы можем свести новую задачу к задаче (4.4.1)

$$\min\{f(x) - f_*, g(x)\} \to \min_{x\in Q}.$$

Несложно записать рандомизированный МЗС для такой задачи и применить к ней все, что изложено выше. Однако мы не будем здесь этого делать, поскольку в следующем пункте мы приведем более общий вариант МЗС, который не предполагает, что известно $f_*$.

### 4.4.3 Рандомизированный метод зеркального спуска с функциональными ограничениями

Рассмотрим задачу

$$f(x) \to \min_{\substack{g(x)\le 0 \\ x\in Q}}. \tag{4.4.5}$$

Под решением этой задачи будем понимать такой $\bar{x}^N \in Q$, что имеют место неравенства

$$E\left[f\left(\bar{x}^N\right)\right] - f_* \le \varepsilon_f = \frac{M_f}{M_g}\varepsilon_g,\ g\left(\bar{x}^N\right) \le \varepsilon_g,$$

где $f_* = f(x_*)$ – оптимальное значение функционала в задаче (5), $x_* \in Q$ – решение задачи (4.4.5). Будем считать, что имеется такая последовательность независимых случайных величин $\{\xi^k\}$ и последовательности $\{\nabla_x f_k(x,\xi^k)\}$, $\{\nabla_x g_k(x,\xi^k)\}$, $k = 1,\dots,N$, что для всех $x \in Q$ имеют место следующие соотношения (можно считать их выполненными при $x \in \tilde{Q}$, см. подраздел 4.4.2 и последующий текст в этом подразделе)

$$E_{\xi^k}\left[\nabla_x f_k(x,\xi^k)\right] = \nabla f(x),\ E_{\xi^k}\left[\nabla_x g_k(x,\xi^k)\right] = \nabla g(x);$$

$$E_{\xi^k}\left[\left\|\nabla_x f_k(x,\xi^k)\right\|_*^2\right] \le M_f^2,\ E_{\xi^k}\left[\left\|\nabla_x g_k(x,\xi^k)\right\|_*^2\right] \le M_g^2.$$

МЗС для задачи (4.4.5) будет иметь вид (см., например, в детерминированном случае двойственный градиентный метод из [281])

$$x^{k+1} = \mathrm{Mirr}_{x^k}\left(h_f \nabla_x f_k\left(x^k,\xi^k\right)\right),\ \text{если}\ g\left(x^k\right) \le \varepsilon_g,$$

$$x^{k+1}=\mathrm{Mirr}_{x^k}\left(h_g\nabla_x g_k\left(x^k,\xi^k\right)\right),\ \text{если}\ g\left(x^k\right)>\varepsilon_g,$$

где $h_g=\varepsilon_g/M_g^2$, $h_f=\varepsilon_g/\left(M_fM_g\right)$, $k=1,...,N$. Обозначим через $I$ множество индексов $k$, для которых $g\left(x^k\right)\le\varepsilon$. Введем также обозначения

$$[N]=\{1,...,N\},\ J=[N]\setminus I,\ N_I=|I|,\ N_J=|J|,\ \bar{x}^N=\frac{1}{N_I}\sum_{k\in I}x^k.$$

Тогда имеет место неравенство

$$h_fN_I\cdot\left(E\left[f\left(\bar{x}^N\right)\right]-f_*\right)\le$$

$$\le h_fE\left[\sum_{k\in I}\left\langle E_{\xi^k}\left[\nabla_xf_k\left(x^k,\xi^k\right)\right],x^k-x_*\right\rangle\right]\le\frac{h_f^2}{2}\sum_{k\in I}E\left[\left\|\nabla_xf_k\left(x^k,\xi^k\right)\right\|_*^2\right]+$$

$$-h_gE\Big[\sum_{k\in J}\underbrace{\left\langle E_{\xi^k}\left[\nabla_xg_k\left(x^k,\xi^k\right)\right],x^k-x_*\right\rangle}_{\ge g\left(x^k\right)-g\left(x_*\right)>\varepsilon_g}\Big]+\frac{h_g^2}{2}\sum_{k\in J}E\left[\left\|\nabla_xg_k\left(x^k,\xi^k\right)\right\|_*^2\right]+$$

$$+\sum_{k\in[N]}\left(E\left[V_{x^k}\left(x_*\right)\right]-E\left[V_{x^{k+1}}\left(x_*\right)\right]\right)\le$$

$$\le\frac{1}{2}h_f^2M_f^2N_I-\frac{1}{2M_g^2}\varepsilon_g^2N_J+E\left[V_{x^1}\left(x_*\right)\right]-E\left[V_{x^{N+1}}\left(x_*\right)\right]=$$

$$=\frac{1}{2}\left(h_f^2M_f^2+\frac{\varepsilon_g^2}{M_g^2}\right)N_I-\frac{1}{2M_g^2}\varepsilon_g^2N+R^2-E\left[V_{x^{N+1}}\left(x_*\right)\right].$$

Будем считать, что (следует сравнить с формулой (4.4.4))

$$N=N\left(\varepsilon_g\right)=\frac{2M_g^2R^2}{\varepsilon_g^2}+1. \qquad (4.4.6)$$

Тогда $N_I\ge1$ с вероятностью $\ge1/2$ и

$$E\left[f\left(\bar{x}^N\right)\right]-f_*\le\frac{1}{2}\left(h_fM_f^2+\frac{\varepsilon_g^2}{M_g^2h_f}\right)=\frac{M_f}{M_g}\varepsilon_g=\varepsilon_f.$$

Соотношение

$$g\left(\bar{x}^N\right)\le\varepsilon_g$$

следует из того, что по построению $g\left(x^k\right)\le\varepsilon_g$, $k\in I$ и из выпуклости функции $g(x)$.

Заметим, что в детерминированном случае вместо $\bar{x}^N$ можно брать

$$\breve{x}^N=\arg\min_{k\in I}f\left(x^k\right).$$

Если известно, что для всех $x\in\tilde{Q}$ и почти наверное по $\xi^k$

$$\left\|\nabla_x f_k\left(x,\xi^k\right)\right\|_*^2 \le \tilde{M}^2,\ \left\|\nabla_x g_k\left(x,\xi^k\right)\right\|_*^2 \le \tilde{M}^2,\ k=1,\dots,N,$$

то для описанного в этом пункте метода (с $\varepsilon = \varepsilon_f = \varepsilon_g$ и $h_f = h_g = \varepsilon/\tilde{M}^2$ ) вид оценки вероятностей больших уклонений (4.4.3) из подраздела 4.4.2 сохранится (оценка получается чуть лучше, чем нижняя оценка из работы [91], когда ограничений-неравенств больше одного, поскольку мы имеем доступ к точному значению $g(x)$)

$$f\left(\bar{x}^N\right) - f_* \le 2\tilde{M}\sqrt{\frac{2}{N}}\left(R + 2\tilde{R}\sqrt{\ln(2/\sigma)}\right),$$

где $N = 2N(\varepsilon)$ (см. формулу (4.4.6), а также раздел 3.4 главы 3). К сожалению, трюк с амплификацией (см. подраздел 4.4.2) здесь уже не проходит в том же виде, как и раньше, поскольку теперь уже нельзя гарантировать

$$f\left(\bar{x}^N\right) - f_* \ge 0.$$

Однако если ввести обозначение

$$\varepsilon_\Delta = f_* - \min_{\substack{g(x)\le\varepsilon_g \\ x\in Q}} f(x) = \min_{\substack{g(x)\le 0 \\ x\in Q}} f(x) - \min_{\substack{g(x)\le\varepsilon_g \\ x\in Q}} f(x),$$

то

$$P\left(f\left(\bar{x}^N\right) - f_* + \varepsilon_\Delta \ge 2\left(\varepsilon_f + \varepsilon_\Delta\right)\right) \le \frac{E\left[f\left(\bar{x}^N\right)\right] - f_* + \varepsilon_\Delta}{2\left(\varepsilon_f + \varepsilon_\Delta\right)} \le \frac{1}{2}.$$

Можно параллельно (независимо) запустить $\log_2\left(\sigma^{-1}\right)$ траекторий метода. Обозначим $\bar{x}_{\min}^N$ тот из $\bar{x}^N$ на этих траекториях, который доставляет минимальное значение $f\left(\bar{x}^N\right)$. Из выписанного неравенства Маркова получаем, что имеет место неравенство

$$P\left(f\left(\bar{x}_{\min}^N\right) - f_* \ge 2\varepsilon_f + \varepsilon_\Delta\right) \le \sigma.$$

К сожалению, этот подход требует малости $\varepsilon_\Delta$, что, вообще говоря, нельзя гарантировать из условий задачи.

Немного более аккуратные рассуждения (без новых идей) позволяют развязать во всех приведенных выше в подразделе 4.4.3 рассуждениях $\varepsilon_f$ и $\varepsilon_g$, допуская, что они могут выбираться независимо друг от друга. Детали мы вынуждены здесь опустить.

Основные приложения описанного подхода, это задачи вида

$$f(x) \to \min_{\max\limits_{k=1,\dots,m} \sigma_k\left(A_k^T x\right) \le 0},$$

с разреженной матрицей

$$A = \left[A_1,\dots,A_m\right]^T.$$

В частности, задачи вида

$$f(x) \to \min_{Ax \le b}$$

и приводящиеся к такому виду задачи

$$f(x) \to \min_{\substack{Ax \le b \\ Cx = d}}.$$

В этих задачах, как правило, выбирают $\| \; \| = \| \; \|_2$, $d(x) = \|x\|_2^2 / 2$. Подобно [281] можно попутно восстанавливать (без особых дополнительных затрат) и двойственные множители к этим ограничениям. Причем эта процедура позволяет сохранить дешевизну итерации даже в разреженном случае.

### 4.4.4 Примеры решения разреженных задач с использованием рандомизированного метода зеркального спуска

Начнем с известного примера [88], демонстрирующего практически все основные способы рандомизации, которые сейчас активно используются в самых разных приложениях. Рассмотрим задачу поиска левого собственного вектора $x$, отвечающего собственному значению 1, стохастической по строкам матрицы $P = \|p_{ij}\|_{i,j=1,1}^{n,n}$ (такой вектор называют вектором Фробениуса–Перрона, а с учетом контекста PageRank вектором – см. раздел 4.1 этой главы 4). Изложение рандомизации, связанной с ускоренными покомпонентными методами мы опускаем, поскольку она не завязана на МЗС. Тем не менее, приведем ссылки на работы, в которых такой подход к поиску PageRank описан: замечание 5 [44] и пример 4 [44] (см. также замечания 10, 11 [6]).

Перепишем задачу поиска вектора PageRank следующим образом [5] (см. разделы 4.1, 4.2 этой главы 4)

$$f(x) = \frac{1}{2}\|Ax\|_2^2 \to \min_{x \in S_n(1)},$$

где $S_n(1)$ – единичный симплекс в $\mathbb{R}^n$, $A = P^T - I$, $I$ – единичная матрица. Далее будем использовать обозначения $A^{\langle k \rangle}$ – $k$-й столбце матрицы $A$, $A_k$ – транспонированная $k$-я строка (то есть $A_k$ – это вектор) матрицы $A$. Следуя [88], воспользуемся для решения этой задачи рандомизированным МЗС со стохастическим градиентом[16]

---

[16] Сначала согласно вектору $x \in S_n(1)$ случайно разыгрываем один из столбцов матрицы $A = P^T - I$ (пусть это будет $\xi^k$-й столбец). Пользуясь тем, что столбцы матрицы $P^T$ сами представляют собой распределения вероятностей, независимо разыгрываем (согласно выбранному столбцу $\xi^k$) случайную величину, и выбира-

$$\nabla_x f_k\left(x,\xi^k\right)=\left(P-I\right)^{\left\langle j\left(\xi^k\right)\right\rangle}-\left(P-I\right)^{\left\langle \xi^k\right\rangle},$$

где

$$\xi^k=i \text{ с вероятностью } x_i,\ i=1,...,n;$$

$$j\left(\xi^k\right)=j \text{ с вероятностью } p_{\xi^k j},\ j=1,...,n.$$

Несложно проверить выполнение условия 4.4.1 подраздела 4.4.2, если генерирование использующихся вспомогательных случайных величин осуществляется независимо. В виду симплексных ограничений, естественно следовать при выборе прокс-структуры примеру 4.4.3 подраздела 4.4.2. Таким образом, можно оценить

$$\tilde{M}^2=\max_{x\in S_n(1),\xi^k}\left\|\nabla_x f_k\left(x,\xi^k\right)\right\|_\infty^2\le 4.$$

Даже в случае, когда матрица $P$ полностью заполнена амортизационная (средняя) стоимость одной итерации будет $\mathrm{O}(n)$ (вместо $\mathrm{O}\left(n^2\right)$, в случае честного расчета градиента). Таким образом, общее число арифметических операций будет $\mathrm{O}\left(n\ln n/\varepsilon^2\right)$.

К худшей оценке приводит другой способ рандомизации (рандомизации суммы [44]). Чтобы его продемонстрировать, перепишем исходную задачу следующим образом

$$f(x)=n\sum_{k=1}^n\frac{1}{n}\frac{1}{2}\left(A_k^Tx\right)^2\to\min_{x\in S_n(1)}.$$

Из такого представления следует, что можно определить стохастический градиент следующим образом

$$\nabla_x f_k\left(x,\xi^k\right)=nA_{\xi^k}A_{\xi^k j(x)},$$

где

$$\xi^k=i \text{ с вероятностью } 1/n,\ i=1,...,n;$$

$$j(x)=j \text{ с вероятностью } x_j,\ j=1,...,n.$$

Амортизационная (средняя) стоимость одной итерации будет по-прежнему $\mathrm{O}(n)$, но вот оценка $\tilde{M}^2$ получается похуже. Здесь мы имеем пример, когда $M^2$ и $\tilde{M}^2$ существенно отличаются – в действительности, можно вводить промежуточные условия, не такие жесткие, как условие 4.4.3, и получать более оптимистичные оценки вероятностей больших уклонений [44].

---

ем соответствующий столбец матрицы $A^T=P-I$, из которого вычитаем $A^{T\left\langle\xi^k\right\rangle}$ (из-за наличия матрицы $I$), таким образом, конструируется несмещенная оценка градиента $\nabla f(x)=A^TAx$.

К сожалению, эти методы не позволяют полноценно воспользоваться разреженностью матрицы $P$, которая, как правило, имеет место. Собственно, этот пункт отчасти и будет посвящен тому, как можно сочетать рандомизацию и разреженность. В частности, если переписать задачу PageRank следующим образом

$$\|Ax\|_\infty \to \min_{x \in S_n(1)},$$

что равносильно (факт из теории неотрицательных матриц [95])

$$\max_{k=1,\dots,n} A_k^T x \to \min_{x \in S_n(1)},$$

то исходя из примера 4.4.3 (в детерминированном случае), можно получить следующую оценку [46] на общее число арифметических операций $O\left(n \ln n / \varepsilon^2\right)$, при условии, что число элементов в каждой строке и столбце матрицы $P$ не больше $O\left(\sqrt{n/\ln n}\right)$. Здесь не использовалась рандомизация, а использовалась только разреженность матрицы $P$ (следовательно и $A$). По-сути, способ получения этой оценки всецело базируется на возможности организации эффективного пересчета субградиента функционала, подобно [5, 273, 281]. Далее мы распространим этот пример на более общий класс задач, и постараемся привнести в подход рандомизацию.

Итак, рассмотрим сначала класс задач с $Q$ из примера 4.4.1 или 4.4.2 подраздела 4.4.2

$$\max_{k=1,\dots,m} \sigma_k\left(A_k^T x\right) \to \min_{x \in Q}, \qquad (4.4.7)$$

где $\sigma_k(\ )$ – выпуклые функции с константой Липшица равномерно ограниченной известным числом $\mathrm{M}$, (суб-)градиент каждой такой функции (скалярного аргумента) можно рассчитать за $O(1)$. Введем матрицу

$$A = \left[A_1, \dots, A_m\right]^T$$

и будем считать, что в каждом столбце матрицы $A$ не больше $s_m \le m$ ненулевых элементов, а в каждой строке – не больше $s_n \le n$. Заметим, что некоторую обременительность этим условиям создает требование, что в “каждом” столбце/строке. Это требования можно ослаблять, приближаясь к некоторым средним показателям разреженности (численные эксперименты в этой связи также проводились [5]), однако в данном разделе для большей наглядности и строгости рассуждений мы ограничимся случаем, когда именно в каждом столбце/строке имеет место такая (или еще большая) разреженность.

Из работ [5, 273, 281] следует, что МЗС из примеров 4.4.1, 4.4.2 (в детерминированном случае) для задачи (4.4.7) будет требовать

$$\mathrm{O}\left(\frac{\mathrm{M}^2 \max\limits_{k=1,\ldots,m} \|A_k\|_2^2 R_2^2}{\varepsilon^2}\right)$$

итераций, где $R_2^2$ – квадрат евклидова расстояния от точки старта до решения, а одна итерация (кроме первой) будет стоить

$$\mathrm{O}\left(s_n s_m \log_2 m\right).$$

И все это требует препроцессинг (предварительных вычислений, связанных с “правильным” приготовлением памяти) объема $\mathrm{O}(m+n)$. Таким образом, в интересных для нас случаях общее число арифметических операций для МЗС из примеров 4.4.1, 4.4.2 будет

$$\mathrm{O}\left(s_n s_m \log_2 m \frac{\mathrm{M}^2 \max\limits_{k=1,\ldots,m} \|A_k\|_2^2 R_2^2}{\varepsilon^2}\right). \qquad (4.4.8)$$

Постараемся ввести рандомизацию в описанный подход. Для этого осуществим дополнительный препроцессинг, заключающийся в приготовлении из векторов $A_k$ вектора распределения вероятностей. Представим

$$A_k = A_k^+ - A_k^-,$$

где каждый из векторов $A_k^+$, $A_k^-$ имеет не отрицательные компоненты. Согласно этим векторам приготовим память таким образом, чтобы генерирование случайных величин из распределений $A_k^+ / \|A_k^+\|_1$ и $A_k^- / \|A_k^-\|_1$ занимало бы $\mathrm{O}(\log_2 n)$. Это всегда можно сделать [5]. Однако это требует хранение в “быстрой памяти” довольно большого количества соответствующих “деревьев”. Весь этот препроцессинг и затраченная память будут пропорциональны числу ненулевых элементов матрицы $A$, что в случае huge-scale задач сложно осуществить из-за ресурсных ограничений. Тем не менее, далее мы будем считать, что такой препроцессинг можно осуществить, и (самое главное) такую память можно получить. Введем стохастический (суб-)градиент

$$\nabla_x f_k\left(x,\xi^k\right) = \left\|A_{k(x)}^+\right\|_1 e_{i(\xi^k)} - \left\|A_{k(x)}^-\right\|_1 e_{j(\xi^k)},$$

где

$$k(x) \in \operatorname{Arg}\max_{k=1,\ldots,m} \sigma_k\left(A_k^T x\right),$$

причем не важно, какой именно представитель $\operatorname{Arg}\max$ выбирается;

$$e_i = \overbrace{(\underbrace{0,\ldots,0,1}_{i},0,\ldots,0)}^{n};$$

$$i\left(\xi^k\right)=i \text{ с вероятностью } A^{+}_{k(x)i}\Big/\left\|A^{+}_{k(x)}\right\|_1,\ i=1,\dots,n;$$

$$j\left(\xi^k\right)=j \text{ с вероятностью } A^{-}_{k(x)j}\Big/\left\|A^{-}_{k(x)}\right\|_1,\ j=1,\dots,n.$$

Легко проверить выполнение условия 4.4.1 подраздела 4.4.2 (заметим, что $\nabla f\left(x\right)=A_{k(x)}$). Также легко оценить

$$\tilde{M}^2 \le \mathrm{M}^2 \max_{k=1,\dots,m}\left\|A_k\right\|_1^2.$$

И получить из примеров 4.4.1, 4.4.2 следующую оценку числа итераций (ограничимся для большей наглядности сходимостью по математическому ожиданию, т.е. без оценок вероятностей больших уклонений)

$$\mathrm{O}\left(\frac{\mathrm{M}^2 \max_{k=1,\dots,m}\left\|A_k\right\|_1^2 R_2^2}{\varepsilon^2}\right)$$

Основная трудоемкость тут в вычислении $k\left(x\right)$. Однако, за исключением самой первой итерации можно эффективно организовать перерешивание этой задачи. Действительно, предположим, что уже посчитано $k\left(x^l\right)$, а мы хотим посчитать $k\left(x^{l+1}\right)$. Поскольку согласно примерам 4.4.1, 4.4.2 $x^{l+1}$ может отличаться $x^l$ только в двух компонентах, то пересчитать $\max_{k=1,\dots,m}\sigma_k\left(A_k^T x^{l+1}\right)$, исходя из известного $\max_{k=1,\dots,m}\sigma_k\left(A_k^T x^l\right)$, можно за (см., например, [5, 273], раздел 2.1 главы 2 и раздел 4.2 этой главы 4) $\mathrm{O}\left(2s_m\log_2 m\right)$. Таким образом, общее ожидаемое число арифметических операций нового рандомизированного варианта МЗС из примеров 4.4.1, 4.4.2 для задачи (4.4.7) будет

$$\mathrm{O}\left(s_m\log_2 m\frac{\mathrm{M}^2 \max_{k=1,\dots,m}\left\|A_k\right\|_1^2 R_2^2}{\varepsilon^2}\right). \tag{4.4.9}$$

Для матриц $A$, все отличные от нуля элементы которых одного порядка, скажем $\mathrm{O}\left(1\right)$, имеем

$$\max_{k=1,\dots,m}\left\|A_k\right\|_2^2=s_n,\ \max_{k=1,\dots,m}\left\|A_k\right\|_1^2=s_n^2.$$

В таком случае не стоит ожидать выгоды (формулы (4.4.8) и (4.4.9) будут выглядеть одинаково). Но если это условие (ненулевые элементы $A$ одного порядка) выполняется не очень точно, то можно рассчитывать на некоторую выгоду.

Рассмотрим теперь более общий класс задач, возникающих, например, при поиске равновесий в транспортных сетях [40] (см. главу 1 диссертации)

$$\frac{1}{r}\sum_{k=1}^{r}\max_{l=a_k+1,\dots,b_k}\sigma_l\left(A_l^T x\right)\to\min_{x\in Q}, \tag{4.4.10}$$

$$0=a_1<b_1=a_2<b_2=a_3<\dots<b_{r-1}=a_r<b_r=m.$$

Матрица $A$ и числа $s_n$, $s_m$ определяются аналогично. Привнося (при расчете стохастического градиента) к описанным выше двум подходам для задачи (4.4.7) сначала равновероятный (и независимый от других рандомизаций) выбор одного из слагаемых в этой сумме, получим соответствующие обобщения (для задачи (4.4.10)) оценок (4.4.8), (4.4.9), которые будут иметь точно такой же вид. Только матрица $A$ собирается теперь из всех слагаемых суммы (4.4.10).

Возвращаясь к примеру 4.4.3, заметим, что все описанные выше конструкции (в том числе, связанные с задачей (4.4.10)) можно перенести на этот пример, в случае, когда

$$\sigma_k\left(A_k^T x\right)=A_k^T x-b_k.$$

При этом

$$R_2^2\to\ln n,\ \mathrm{M}=1\ \text{(для (4.4.8), (4.4.9))}$$

$$\max_{k=1,\dots,m}\left\|A_k\right\|_2^2\to\max_{\substack{i=1,\dots,m\\ j=1,\dots,n}}\left|A_{ij}\right|,\ \min\left\{s_n s_m\log_2 m,n\right\}\to\max\left\{s_n s_m\log_2 m,n\right\}\ \text{(для (4.4.8))}$$

$$\max_{k=1,\dots,m}\left\|A_k\right\|_1^2\to\max_{k=1,\dots,m}\left\|A_k\right\|_1^2,\ s_m\log_2 m\to\max\left\{s_m\log_2 m,n\right\}\ \text{(для (4.4.9))}.$$

Собственно, примера PageRank, изложенный в начале этого пункта, как раз подходил под применение оценки (4.4.8).

С помощью подраздела 4.4.3 все написанное выше переносится и на задачи вида

$$f(x)\to\min_{\substack{Ax\le b\\ Cx=d}},$$

с разреженными матрицами. Такие задачи играют важную роль, например, при проектировании механических конструкций [281] (Truss topology design). Мы не будем здесь приводить соответствующие рассуждения, поскольку они достаточно очевидны.

# Глава 5 Покомпонентные и безградиентные методы для задач выпуклой оптимизации

## 5.1 О нетривиальности быстрых (ускоренных) рандомизированных методов

### 5.1.1 Введение

Для первоначального погружения в описываемую далее проблематику можно рекомендовать пп. 6.2 – 6.5 обзора [163].

В данном разделе мы хотим подчеркнуть нетривиальность таких методов как, например, ускоренный (быстрый) покомпонентный метод Ю.Е. Нестерова [264], APPROX или ALPHA [199, 295], которые, в частности, являются покомпонентными вариантами быстрого градиентного метода (БГМ) [163] – эту ссылку можно также рекомендовать с точки зрения интересной подборки ссылок на работы, в которых объясняется, что такое БГМ. Нетривиальность в том, что эти методы являются рандомизированными и при этом ускоренными. Число необходимых итераций (как функция от желаемой точности) для таких ускоренных покомпонентных методов увеличивается в число раз $n$, равное размерности пространства, по сравнению с классическим БГМ, что и не удивительно, поскольку вместо всех $n \gg 1$ компонент градиента на каждой итерации используется только одна. Также нетривиальность в том, что если полный расчет градиента, скажем, требовал полного умножения разреженной матрицы на вектор – $sn$ операций, то пересчет (важно, что именно пересчет, а не расчет) компоненты градиента в определенных ситуациях требует всего $s$ операций (см. стр. 16–17 [199] и подразделы 5.1.4, 5.1.5). Таким образом, увеличение числа итераций в $n$ раз компенсируется уменьшением стоимости одной итерации в $n$ раз (в не разреженном случае оговорка об “определенных ситуациях” существенно ослабляется, см. подразделы 5.1.4, 5.1.5). Но выгода от использования покомпонентных методов, как правило, есть из-за того, что в таких методах вместо константы Липшица градиента по худшему направлению (максимального собственного значения матрицы Гессе функционала) в оценки числа итераций входит “средняя” константа Липшица (оценивающаяся сверху средним арифметическим суммы диагональных элементов (следа) матрицы Гессе, т.е. средним арифметическим всех собственных чисел матрицы Гессе). Разница в этих константах для матриц Гессе, состоящих из элементов одного порядка, может равняться по порядку $n$ (см. пример 5.1.2 подраздела 5.1.5). На данный момент известно довольного много примеров применения покомпонентных методов для решения задач ог-ромных размеров, в частности, приложений для задач моделирования сетей больших размеров и анализе данных [138].

В подразделах 5.1.2, 5.1.3 мы демонстрируем те сложности, которые возникают при попытках получить ускоренные покомпонентные методы из оптимальных методов для задач стохастической оптимизации без погружения в то, как устроены эти оптимальные методы. Мы не ставили себе в этих пунктах цель: получить и подробно исследовать какие-то новые эффективные методы, поэтому изложение в этих пунктах ведется на “физическом” уровне строгости. В подразделе 5.1.4 мы приводим новое доказательство оценки скорости сходимости ускоренного покомпонентного метода, базирующееся на конструкции линейного каплинга [135]: БГМ = “выпуклая комбинация” прямого градиентного метода (ПГМ) и метода зеркального спуска (МЗС). Основная идея получения ускоренного покомпонентного метода: заменить в таком представлении БГМ в методах ПГМ и МЗС градиенты на соответствующие несмещенные оценки градиентов, полученные на основе покомпонентной рандомизации. Несмотря на то, что основной результат подраздела 5.1.4 (теорема 5.1.2 о сходимости предложенного метода и замечания к ней) – не есть полностью новый результат, подобные оценки (в различных частных случаях) ранее уже встречались в литературе, тем не менее, способ их получения (и его универсальность) представляется новым и весьма перспективным с точки зрения возможных последующих обобщений и приложений (некоторые примеры таких приложений и обобщений приведены в подразделах 5.1.4, 5.1.5). Описанный способ также позволяет устанавливать различные новые факты об ускоренных покомпонентных методах. Наброски приведены в цикле замечаний к теореме 5.1.2 в подразделе 5.1.4 и в примерах подраздела 5.1.5. В подразделе 5.1.6 кратко резюмируются результаты этого раздела, приводятся заключительные замечания.

### 5.1.2 Нетривиальность ускоренных покомпонентных методов

Рассматривается задача гладкой выпуклой оптимизации

$$f(x) \to \min_{x \in Q}.$$

Мы постараемся сначала пояснить, как получить в неускоренном случае для этой задачи оценки для покомпонентных спусков из оценок рандомизированных методов решения этой задачи. Оказывается это можно довольно изящно сделать. К сожалению, при этом даже из оптимальных рандомизированных методов не удается “вытащить” оценки для ускоренных покомпонентных методов. В этом-то и заключается нетривиальность ускоренных покомпонентных методов.

Для простоты считаем, что везде в дальнейшем (в подразделах 5.1.2, 5.1.3) мы говорим о 2-норме и евклидовой прокс-структуре (интересно было бы понять, как все, что далее будет написано, распространяется на более общие нормы/прокс-структуры [44]). Считаем, что функция $f(x)$ имеет Липшицев градиент с константой $L$, является $\mu$-сильно выпуклой, а множество $Q$ имеет диаметр $R$, при этом в точке минимума $\nabla f(x_*)=0$. Последнее предположение – обременительное для задач условной оптимизации. К сожалению, мы пока не знаем как от него отказаться.

Будем считать, что на каждой итерации оракул выдает нам несмещенную оценку градиента с дисперсией $D$. Определим зависимость $N(\varepsilon)$ для изучаемого итерационного процесса: $N(\varepsilon)$ – наименьшее $N$ такое, что ( $f_* = f(x_*)$, где $x_*$ – решение задачи)

$$E\left[f(x_N)\right]-f_*\le\varepsilon.$$

**Теорема 5.1.1.** *Существуют такие неускоренные методы (см., например, п. 6.2 [163]), которые работают по оценкам (* $\Delta f^0 = f(x_0)-f_*$ *)*

$$\varepsilon=\min\left\{\mathrm{O}\left(\frac{LR^2}{N}+\sqrt{\frac{DR^2}{N}}\right),\mathrm{O}\left(\Delta f^0\exp\left(-N\frac{\mu}{L}\right)+\frac{D}{\mu N}\right)\right\}.$$

*Существуют такие ускоренные методы (например, линейки SIGMA, см., например, [44, 181, 193] и раздел 2.1 главы 2), которые работают по следующим не улучшаемым оценкам*

$$\varepsilon=\min\left\{\mathrm{O}\left(\frac{LR^2}{N^2}+\sqrt{\frac{DR^2}{N}}\right),\mathrm{O}\left(\Delta f^0\exp\left(-N\sqrt{\frac{\mu}{L}}\right)+\frac{D}{\mu N}\right)\right\}.$$

Будем говорить о сильно выпуклом случае, если минимум в этих формулах достигается на втором аргументе.

Далее заметим, что если мы вместо обычного градиента используем его аппроксимации, возникающие в безградиентных и покомпонентных подходах [44] (когда оракул может на каждой итерации выдавать только значение функции в двух точках или производную по указанному нами направлению)

$$g_\tau(x,s)=\frac{n}{\tau}\big(f(x+\tau s)-f(x)\big)s \text{ или } g(x,s)=n\langle\nabla f(x),s\rangle s,$$

где $s$ – случайный вектор, равномерно распределенный на $S_2^n(1)$ – единичной сфере в 2-норме в пространстве $\mathbb{R}^n$, то имеет место следующий простой факт (см. также следующий раздел), являющийся следствием явления концентрации равномерной меры на сфере вокруг экватора [52, 155, 244] (северный полюс задается градиентом).

**Утверждение 5.1.1.** *Имеют место следующие формулы (см. [44, 52, 155])*

$$E_s\left[\left\|g_\tau(x,s)\right\|_2^2\right] \le 4n\left\|\nabla f(x)\right\|_2^2 + L^2\tau^2 n^2,$$

$$E_s\left[\left\|g(x,s)\right\|_2^2\right] = n\left\|\nabla f(x)\right\|_2^2.$$

Далее мы ограничимся оценкой для покомпонентного метода. “Почувствовать” эту оценку можно на еще более простом примере, когда несмещенная оценка градиента

$$g(x,s) = n\left\langle \nabla f(x), s\right\rangle s,$$

получается за счет другого выбора случайного вектора $s$. Мы считаем, что $s$ принимает равновероятно одно из $n$ направлений, соответствующих единичных ортов. В таком случае также имеет место соотношение

$$E_s\left[\left\|g(x,s)\right\|_2^2\right] = n\left\|\nabla f(x)\right\|_2^2.$$

Попробуем теперь исходя из полученных оценок и теоремы 5.1.1 в сильно выпуклом случае с $\nabla f(x_*) = 0$ получить оценку

$$N(\varepsilon) = \mathrm{O}\left(n\frac{L}{\mu}\ln\left(\frac{\Delta f^0}{\varepsilon}\right)\right)$$

для спусков по направлению. Итак, мы считаем, что вместо градиента оракул на каждой итерации (в точке $x_k$) может нам выдавать только производную по указанному нами направлению. Тогда мы имеем несмещенную оценку градиента с дисперсией (приводимая ниже оценка является не улучшаемой с точностью до мультипликативной константы)

$$D = \mathrm{O}\left(n\left\|\nabla f(x_k)\right\|_2^2\right).$$

Используя тот факт, что для любой гладкой выпуклой функции (в предположении, что $\nabla f(x_*) = 0$; для последнего неравенства еще нужно потребовать, чтобы $k$ не был слишком маленьким)

$$\left\|\nabla f(x_k)\right\|_2^2 \le 2L\cdot\left(f(x_k) - f_*\right) \le 2L\Delta f^0,$$

получаем из теоремы 5.1.1, что после

$$N = \mathrm{O}\left(5nL/\mu\right)$$

итераций

$$f\left(x_N\right) - f_* \le \mathrm{O}\left(\Delta f^0 \exp\left(-N\frac{\mu}{L}\right) + \frac{2nL\Delta f^0}{\mu N}\right) \le \mathrm{O}\left(\frac{1}{2}\Delta f^0\right)$$

и, тем более,

$$f\left(x_N\right) - f_* \le \mathrm{O}\left(\Delta f^0 \exp\left(-N\sqrt{\frac{\mu}{L}}\right) + \frac{2nL\Delta f^0}{\mu N}\right) = \mathrm{O}\left(\frac{1}{2}\Delta f^0\right).$$

Тут также можно пользоваться методами, которые работают по оценкам [218, 228, 298] (собственно, именно так и делается в разделе 4.3 главы 4)

$$\varepsilon = \mathrm{O}\left(M^2/\left(\mu N\right)\right),$$

где

$$M^2 = E_s\left[\left\|g_{\tau,\delta}\left(x_k, s\right)\right\|_2^2\right] = \mathrm{O}\left(n\left\|\nabla f\left(x_k\right)\right\|_2^2\right),$$

что не удивительно, поскольку мы, фактически, при данном подходе и работаем с $M^2$, а не с $D \le M^2$. В любом случае при таком способе рассуждений возникает неаккуратность, связанная с тем, что мы лишь обеспечили

$$f\left(x_N\right) - f_* \le \mathrm{O}\left(\frac{1}{2}\Delta f^0\right).$$

В действительности, тут нужно аккуратно выписывать константы, которые в итоге увеличат константу “5” в ожидаемой сейчас формуле $N = 5nL/\mu$ в несколько раз. Наконец, необходимо проводить рассуждения с оценками вероятностей больших уклонений (здесь помогают самые грубые неравенства типа Буля, поскольку, имеются субгауссовские хвосты, а точнее вообще финитный носитель у стохастического градиента). Далее мы уже не будем делать такие оговорки, поскольку в этом и следующем пункте мы преследуем цель – продемонстрировать нетривиальность ускоренных рандомизированных методов, а не точного выписывания методов, которые получаются по ходу рассуждений. Эти методы не очень интересны, поскольку заведомо не являются оптимальными.

Делая $\log_2\left(\Delta f^0/\varepsilon\right)$ таких перезапусков (стартуем в новом цикле с той точки, на которой остановились на прошлом цикле) с $N = \mathrm{O}\left(5nL/\mu\right)$ итерациями на каждом перезапуске (цикле), в итоге получим оценку общего числа итераций

$$N(\varepsilon) = \mathrm{O}\left( n\frac{L}{\mu}\ln\left(\frac{\Delta f^0}{\varepsilon}\right)\right).$$

Здесь можно "поиграться" на, так называемом, mini-batch'инге (см., например, п. 6.2 [163]), для этого нужна уже формула с дисперсией, то есть $\varepsilon = \mathrm{O}\left(M^2/(\mu N)\right)$ не подходит.

Тем не менее, даже при использовании не улучшаемых (с точностью до мультипликативной константы) рандомизированных методов из теоремы 5.1.1 мы не смогли получить оценки работ [199, 264, 295], для ускоренных покомпонентных методов

$$N(\varepsilon) = \mathrm{O}\left( n\sqrt{\frac{L}{\mu}}\ln\left(\frac{\Delta f^0}{\varepsilon}\right)\right).$$

Строго говоря, в работах [199, 295] таких оценок для покомпонентных методов (в сильно выпуклом случае) мы и не видели (отметим при этом, что такая оценка есть для метода из работы [264] и для безградиентного метода из работы [270]), однако в [199, 295] есть аналогичные "ускоренные" оценки в несильно выпуклом случае. Кроме того, Питер Ричтарик (см. работы на личной странице http://www.maths.ed.ac.uk/~richtarik/ ) сообщил нам, что он умеет устанавливать эти оценки (для ускоренных покомпонентных методов и в сильно выпуклом случае), и сейчас готовит статью на эту тему (свой способ получения таких оценок мы изложим в подразделе 5.1.4). Проблема тут в том, что мы пользовались правым неравенством (считаем $\nabla f(x_*) = 0$)

$$2\mu\cdot\left(f(x_k) - f_*\right) \le \left\|\nabla f(x_k)\right\|_2^2 \le 2L\cdot\left(f(x_k) - f_*\right).$$

Мы специально здесь написали и левое неравенство. Отсюда видно, что в принципе при использовании правого неравенства мы можем потерять $L/\mu$. Не удивительно, что в итоге мы, действительно, теряем $\sqrt{L/\mu} \ll L/\mu$. К сожалению, такого рода рассуждения не позволяют никак "вытащить" оценки оптимального (ускоренного) покомпонентного метода из соответствующих оптимальных полных градиентных методов, не погружаясь в детальный анализ доказательства их сходимости. Нетривиально то, что это один из тех редких примеров (другой см. ниже в подразделе 5.1.3), когда такая философия переноса не сработала. Обычно все удается перенести без особых погружений в детали доказательства. То есть работает принцип: *оптимальный метод порождает оптимальный.*[1]

---

[1] Хорошо известные примеры тут: 1) регуляризация [264], позволяющая переносить оптимальные методы, работающие в сильно выпуклом случае, на просто выпуклый случай; 2) техника рестартов [135, 228] (см.

### 5.1.3 Нетривиальность ускоренных методов рандомизации суммы

Рассмотрим теперь в задаче подраздела 5.1.2 случай когда (этот случай разбирается, например, в п. 6.3 [163] и работах [223, 232, 234, 242])

$$f(x)=\frac{1}{m}\sum_{k=1}^{m} f_k(x),$$

где все функции гладкие с константой Липшица градиента $L$. Также как и раньше считаем $f(x)$ $\mu$-сильно выпуклой. В качестве несмещенной оценки градиента будем брать вектор (по поводу определения SIGMA см. [193])

$$\nabla f\left(x_t^s,\xi\right)=\nabla f_\xi\left(x_t^s\right)-\nabla f_\xi\left(y^s\right)+\nabla f\left(y^s\right),\ y^s=x_N^{s-1},$$

$$x_{t+1}^s=\mathrm{SIGMA}\left(x_t^s,\nabla f\left(x_t^s,\xi\right)\right),\ t=0,\dots,N-1,$$

где случайная величина $\xi$ принимает равновероятно одно из значений $1,\dots,m$; параметр $N$ будет выбран позже как $N=\mathrm{O}\left(4L/\mu\right)$. Здесь по $t$ идет внутренний цикл, а по $s$ внешний.

Приведенный метод (как и метод из п. 6.3 [163]) можно обобщить (с сохранением всех последующих оценок и способов их получения) на стохастический случай, когда $f_k(x):=E_\eta\left[f_k\left(x;\eta_k\right)\right]$, где $f_k\left(x;\eta_k\right)$ – выпуклые по $x$ функции с равномерно (по $\eta_k$) ограниченными (числом $L$) константами Липшица градиентов. При этом (случайная величина $\eta_\xi^{s,t}$ имеет такое же распределение, как и $\eta_\xi$; также считаем, что $\eta_\xi^{s,t}$ ни от чего не зависит, в частности, от других $\left\{\eta_\xi^{s,t}\right\}$ и от $\xi$)

$$\nabla f\left(x_t^s,\xi;\eta_\xi^{s,t}\right)=\nabla f_\xi\left(x_t^s;\eta_\xi^{s,t}\right)-\nabla f_\xi\left(y^s;\eta_\xi^{s,t}\right)+\nabla f\left(y^s\right).$$

Это очевидно, для случая, когда

$$f_k(x):=E_\eta\left[f_k\left(x;\eta_k\right)\right]=\frac{1}{l}\sum_{i=1}^{l} f_k(x;i),$$

поскольку все сводится к исходной постановке с $m:=ml$. Далее мы ограничимся рассмотрением только детерминированного случая.

**Утверждение 5.1.2.** *Имеет место следующая оценка (*$\Delta f^s=f\left(y^s\right)-f_*$*)*

---

также подраздел 5.1.4 далее), позволяющая из оптимального метода для выпуклой задачи, получить оптимальный метод для сильно выпуклой задачи.

$$D^s = E_\xi\left[\left\|\nabla f\left(x_t^s,\xi\right)-E_\xi\left[\nabla f\left(x_t^s,\xi\right)\right]\right\|_*^2\right]=\mathrm{O}\left(L\cdot\left(f\left(y^s\right)-f_*\right)+L\cdot\left(f\left(x_t^s\right)-f_*\right)\right)=\mathrm{O}\left(L\Delta f^s\right).$$

Для доказательства утверждения 5.1.2 в случае $\nabla f\left(x_*\right)=0$ см., например, формулу (6.3) и лемму 6.1 [163] и цитированную в п. 6.3 [163] литературу (в частности, [223, 232, 242]). Причем это утверждение можно формулировать с точными константами вместо $\mathrm{O}(\ )$, чтобы ей можно было далее практически воспользоваться (однако, мы здесь не будем этого делать). В случае $\nabla f\left(x_*\right)\neq 0$ доказательство утверждения 5.1.2 нам не известно (не известно даже останется ли оно верным).

Возьмем теперь в теореме 5.1.1 $N=\mathrm{O}\left(4L/\mu\right)$ и воспользуемся утверждением 5.1.2

$$\Delta f^{s+1}=\mathrm{O}\left(\Delta f^s\exp\left(-N\sqrt{\frac{\mu}{L}}\right)+\frac{D^s}{\mu N}\right),$$

$$\frac{D^s}{\mu N}=\mathrm{O}\left(\frac{L\Delta f^s}{\mu N}\right)=\mathrm{O}\left(\frac{1}{4}\Delta f^s\right),\ \Delta f^s\exp\left(-N\sqrt{\frac{\mu}{L}}\right)=\mathrm{O}\left(\frac{1}{4}\Delta f^s\right).$$

Получим

$$\Delta f^{s+1}\le\mathrm{O}\left(\frac{1}{2}\Delta f^s\right).$$

Здесь можно “поиграться” на mini-batch’инге [234] (см. также раздел 2.2 главы 2), вычисляя вместо $\nabla f\left(x_t^s,\xi\right)$ агрегат

$$\frac{1}{r}\sum_{i=1}^{r}\nabla f\left(x_t^s,\xi_i\right).$$

Таким образом, если у нас есть возможность параллельно на одной итерации вычислять градиенты $\nabla f_{\xi_i}\left(x_t^s\right)$, то (поскольку дисперсия этого агрегата будет $\mathrm{O}\left(L\Delta f^s/r\right)$) можно выбирать $r=\mathrm{O}\left(2\sqrt{L/\mu}\right)$ (выбрать $r$ большим нельзя) и, соответственно, сократить число итераций на цикле до $N=\mathrm{O}\left(2\sqrt{L/\mu}\right)$. Далее, также как и раньше, делаем $\log_2\left(\Delta f^0/\varepsilon\right)$ перезапусков (циклов), на каждом из которых вначале надо посчитать полный градиент (это стоит $m$ вычислений градиентов слагаемых), а потом еще сделать $N$ итераций, на каждой из которых дополнительно требуется рассчитывать в новой точке градиент одного (или нескольких, если используется mini-batch’инг) слагаемого. Таким образом, общая сложность (измеряемая на этот раз не в итерациях, а числе вычислений градиентов слагаемых в сумме в представлении $f(x)$, при этом мы считаем, что сложность вычисления разных слагаемых одинакова в категориях $O(\ )$) будет

$$O\left(\left(m+\frac{L}{\mu}\right)\cdot\ln\left(\frac{\Delta f^0}{\varepsilon}\right)\right),$$

что соответствует части нижней оценки в классе детерминированных методов (см. [129])

$$O\left(\left(m+\min\left\{\sqrt{m\frac{L}{\mu}},\frac{L}{\mu}\right\}\right)\cdot\ln\left(\frac{\Delta f^0}{\varepsilon}\right)\right).$$

S. Shalev-Shwartz-ем и T. Zhang-ом был поставлен вопрос (см., например, п. 6.3 [163]): возможно ли достичь такой (нижней) оценки в целом каким-нибудь методом? То есть речь опять (как и раньше) на самом деле о том, можно ли сохранить ускоренность метода? Только пока, наверное, не очень понятно, причем здесь (ускоренные) рандомизированные покомпонентные методы. Ситуация проясняется, если мы выделим строго выпуклое слагаемое из $f(x)$ в отдельный композит $g(x)$ (см., например, п. 5.1 [163]):

$$f(x)=\frac{1}{m}\sum_{k=1}^{m}f_k(x)+g(x),$$

и построим специальную двойственную задачу (см. замечание 5 [138], а также пример 5.1.4 подраздела 5.1.5 ниже). Тогда (при некоторых дополнительных предположениях [138, 307, 330]) удается показать, что в некотором смысле выписанная нижняя оценка, действительно, достигается на ускоренном покомпонентном методе для двойственной задачи (другой подход описан в http://arxiv.org/pdf/1603.05953.pdf ). Только для сопоставления потребуется перейти от анализа числа вычислений слагаемых к общему числу арифметических операций (см. пример 5.1.4 подраздела 5.1.5 ниже). Собственно, S. Shalev-Shwartz и T. Zhang в работе [307] сами таким образом (правда, не подчеркивая, что, по сути, используют для двойственной задачи ускоренный покомпонентный спуск) и привели пример достижимости (опять оговоримся, что в некотором смысле) нижней оценки.

В связи со сказанным выше отметим, что приведенные в этом разделе оценки (подобно оценкам для покомпонентных методов, см., например, подраздел 5.1.4) можно рассматривать в случае не равномерной рандомизации (выбора слагаемых), а также в случае разных свойств гладкости у разных функций. Сейчас в оценки методов, описанных в этом пункте, входит худшая (по всем слагаемым) константа Липшица градиентов (мы выбрали все константы Липшица градиентов слагаемых для наглядности одинаковыми). В действительности, можно перейти к некоторым их средним вариантам. Однако мы не будем здесь этого делать, поскольку, как уже отмечалось, природа описанных рандомизированных методов вскрывается применением покомпонентных методов к двойственной задачи, для которых все эти нюансы хорошо проработаны.

В этом примере мы опять видим, что попытка из оптимального метода для задач стохастической оптимизации SIGMA (см. теорему 5.1.1) вытащить оптимальные оценки для задачи со специальной (скрытой, через двойственную задачу, покомпонентной) структурой не привели к успеху (на mini-batch'инг не стоит обращать внимание, он просто позволяет параллелить вычисления, не более того). Таким образом, это лишний раз подчеркивает некоторую самостоятельность и важность отдельного изучения ускоренных покомпонентных методов. Именно такого типа методы (например, APPROX, ALPHA [199, 295]) позволяют получать наилучшие оценки. И оценки скорости сходимости этих методов представляют отдельный интерес (в смысле их получения). Насколько нам известно (лучше всего следить за этой областью по работам P. Richtarik-a и T. Zhang-a (см. работы на личной странице http://www.stat.rutgers.edu/home/tzhang/ ), и цитированной ими литературы), сейчас для таких методов используется только евклидова прокс-структура, используются только простые ограничения (сепарабельные), которые обычно зашивают в композитный член [199]. Интересно также было бы понять (охарактеризовать) класс задач, в которых возможно эффективно организовать пересчет компоненты градиента для ускоренных покомпонентных методов. Кое-что на эту тему имеется вот здесь стр. 16–17 [199], [245, 282]. Подробнее все это будет рассмотрено далее.

### 5.1.4 Получение ускоренных покомпонентных методов с помощью каплинга неускоренных прямых покомпонентных методов и покомпонентного метода зеркального спуска

Как уже отмечалось на текущий момент не до конца ясно, насколько все, что сейчас известно для методов, в которых доступен полный градиент, имеет свои аналоги и в (блочно-)покомпонентных методах. Скажем, не все понятно с тем, как можно играть на выборе прокс-функции в ускоренных (блочно-)покомпонентных методах, не до конца ясно: можно ли (если можно, то каким образом) использовать ускоренные покомпонентные методы, если рассматривается задача условной минимизации с множеством специальной простой структуры (в смысле прокс-проектирования, а точнее (блочно-)покомпонентной версии этой операции), но, вообще говоря, не сепарабельной структуры – в частности, в таких задачах $\nabla f\left(x_*\right) \neq 0$ [199, 264]; есть ли аналог универсального метода Ю.Е. Нестеро-ва [274] в покомпонентном варианте; имеют ли покомпонентные методы прямо-двойственную структуру [138]; как перенести на покомпонентные методы концепцию неточного оракула Деволдера–Глинера–Нестерова (см., например, [44, 181, 193]); верно ли, что для ускоренного покомпонентного метода расстояние от любо точки итерационного процесса до решения всегда ограничено некоторой универсальной небольшой константой (меньшей 10), умноженной на расстояние от точки старта до решения, как это имеет место для БГМ (замечание 4 [44])? Список можно продолжить, однако мы здесь остановимся, и сформулируем общий тезис: *все, что сейчас известно для методов первого порядка, в которых доступен полный градиент, имеет (с оговорками, о возможности перенесения результатов на не сепарабельные множества) свои аналоги и в (блочно-)покомпонентных методах; более того, константы (Липшица градиента), фигурирующие в обычных градиентных методах, рассчитанные на худший случай (худшее направление), в покомпонентных методах заменяются "средними" значениями, что в определенных ситуациях может давать ускорение в корень из размерности пространства раз (не говоря о том, что покомпонентные методы при этом могут еще и хорошо параллелиться) [138, 199, 282, 327]*.

К сожалению, исходя из всех известных нам опубликованных на данный момент способов вывода (доказательства сходимости) покомпонентных методов (наиболее, конечно, интересны тут ускоренные варианты) сформулированный выше тезис (гипотеза) совсем не кажется очевидным. Однако совсем недавно, в работе [135] был предложен изящный и перспективный во многих отношениях[2] способ получения БГМ с помощью выпуклого каплинга (комбинации) обычного (неускоренного) прямого градиентного метода[3] (ПГМ) и метода зеркального спуска (МЗС). Ряд "хороших" свойств (например, прямо-двойственность) "наследуются" при таком представлении от зеркального спуска. Естественно, возникает идея попробовать использовать соответствующие легко исследуемые в отдельности покомпонентные аналоги этих двух структурных блоков, чтобы получить ус-коренный покомпонентный метод. Оказывается, что это, действительно, можно сделать. Далее, основываясь на результатах работы [135], мы приведем соответствующие выкладки.

[2] В том числе в отношении более простого обоснования отмеченной выше возможности перенесения свойств с полноградиентных методов на покомпонентные.

[3] Использование в БГМ в качестве одного из структурных блоков именно ПГМ (это заметно упрощает рассуждения в случае $Q = \mathbb{R}^n$ по сравнению с другими возможными вариантами) не является обязательным атрибутом. По-видимому, можно построить (в схожем ключе) аналог БГМ (с аналогичными оценками скорости сходимости) на базе МЗС (или его "сходящегося" варианта [280]) и прямого проксимального градиентного метода (ППГМ) и(или) двойственного градиентного метода [181], который будет лишен отмеченных недостатков.

Исходя из написанного в предыдущих пунктах, можно сказать, что для получения ускоренных покомпонентных методов требуется более тонкая игра (на каждой итерации) на правильном сочетании базовых методов со специальным выбором параметров. Оптимальный метод порождается выпуклой комбинацией неоптимальных методов для класса гладких задач, и именно из этого стоит исходить (распространяя конструкцию на покомпонентные методы), чтобы получить ускоренный покомпонентный метод. Этот тезис нам также представляется полезным, поскольку он подтверждает, что оптимальные методы порождают оптимальные, просто в ряде случаев требуется заглядывание в структуру (базис) метода, чтобы иметь возможность из него породить что-то новое оптимальное.

Сначала мы постараемся в максимально упрощенной ситуации пояснить, как можно получить ускоренный покомпонентный метод, исходя из конструкции п. 3 работы [135]. Все, что далее будет написано, допускает серьезные обобщения, о которых мы упомянем ближе к концу этого пункта.

Итак, рассмотрим задачу

$$f(x) \to \min_{x \in Q}.$$

Введем необходимые в дальнейшем обозначения/определения:

$$e_i = \overbrace{(\underbrace{0,\ldots,0,1}_{i},0,\ldots,0)}^{n};$$

$$\left|\partial f(x+he_i)/\partial x_i - \partial f(x)/\partial x_i\right| \le L_i h \text{ для всех } x \in \mathbb{R}^n \text{ и } h \in \mathbb{R};$$

$$\|x\|^2 = \sum_{i=1}^{n} L_i x_i^2, \quad \|\nabla f(x)\|_*^2 = \sum_{i=1}^{n} L_i^{-1}\left(\partial f(x)/\partial x_i\right)^2;$$

$$d(x) = \frac{1}{2}\|x\|^2, \quad V_x(y) = d(y) - \langle \nabla d(x), y - x\rangle - d(x) = \frac{1}{2}\|y - x\|^2;$$

$$\nabla_i f(x) = (\underbrace{0,\ldots,0,\partial f(x)/\partial x_i}_{i},0,\ldots,0);$$

$$i \in [1,\ldots,n] \text{ – означает, что } P(i = j) = n^{-1},\ j = 1,\ldots,n;$$

$$E_i\left[G(i)\right] \text{ – математическое ожидание по } i \in [1,\ldots,n];$$

$$E_{i_{k+1}}\left[G(i_1,\ldots,i_{k+1}) \middle| i_1,\ldots,i_k\right] = g(i_1,\ldots,i_k) \text{ – условное математическое ожидание по } i_{k+1} \in [1,\ldots,n];$$

$E_{i_1,...,i_k}\left[E_{i_{k+1}}\left[G\left(i_1,...,i_{k+1}\right)\middle|i_1,...,i_k\right]\right]=E\left[G\left(i_1,...,i_{k+1}\right)\right]$ – полное математическое ожидание по всему набору $i_1,...,i_{k+1}\in\left[1,...,n\right]$;

$$\text{Grad}_i\left(x\right)=\arg\min_{\tilde{x}\in Q}\left\{\left\langle\nabla_i f\left(x\right),\tilde{x}-x\right\rangle+\frac{1}{2}\left\|\tilde{x}-x\right\|^2\right\}=x-\frac{1}{L_i}\nabla_i f\left(x\right);$$

$$\text{Mirr}_z\left(\xi\right)=\arg\min_{y\in Q}\left\{\left\langle\xi,y-z\right\rangle+V_z\left(y\right)\right\}=\left(\left\{z_i-\frac{1}{L_i}\xi_i\right\}_{i=1}^n\right).$$

Приведенные формулы специально были записаны таким образом, чтобы их легко было перенести на случай, когда выбирается не одна компонента $i$, а целый блок компонент и $Q\neq\mathbb{R}^n$.

Опишем костяк покомпонентного ускоренного метода (Accelerated by Coupling Randomized Coordinate Descent – ACRCD) на базе специального каплинга покомпонентных вариантов ПГМ (Grad) и МЗС (Mirr) ($x_0=y_0=z_0$)

$$\text{ACRCD}\left(\alpha,\tau;\Theta,x_0,f\left(x_0\right)-f_*\leq d\right)$$

1. $x_{k+1}=\tau z_k+\left(1-\tau\right)y_k$, $\tau\in\left[0,1\right]$ – будет выбрано позже;
2. $i_{k+1}\in\left[1,...,n\right]$ – независимо от предыдущих розыгрышей;
3. $y_{k+1}=\text{Grad}_{i_{k+1}}\left(x_{k+1}\right)$;
4. $z_{k+1}=\text{Mirr}_{z_k}\left(\alpha n\nabla_{i_{k+1}}f\left(x_{k+1}\right)\right)$, $\alpha>0$ – будет выбрано позже.

Поскольку

$$E_i\left[n\nabla_i f\left(x\right)\right]=\nabla f\left(x\right),$$

то шаг 4 (согласно формуле (3.1) [135]; в евклидовом случае можно ограничиться более простыми рассуждениями – см., например, стр. 223 [269]; другой способ получить оценки для МЗС, показывающий дополнительную связь МЗС и проксимального ПГМ (ППГМ), – воспользоваться оценками для метода ППГМ с неточным оракулом из [181, 274]) влечет[4]

[4] Отметим, что первое неравенство специально записано таким образом (в достаточно общем виде) чтобы была видна возможность рассмотрение прокс-структур отличных от евклидовой. По-видимому, для покомпонентных методов в подавляющем большинстве приложений можно ограничиться рассмотрением только евклидовой прокс-структуры. Нюансы могут возникать, когда вместо одной компоненты разрешается сразу случайно выбирать целый блок компонент (необязательно постоянного размера) [199, 264, 295], что для евклидовой прокс-структуры можно также понимать как mini-batch’инг [163143]. Иногда в

$$\alpha n\left\langle \nabla_{i_{k+1}} f\left(x_{k+1}\right), z_k - u\right\rangle \le \frac{\alpha^2 n^2}{2}\left\|\nabla_{i_{k+1}} f\left(x_{k+1}\right)\right\|_*^2 + V_{z_k}\left(u\right) - V_{z_{k+1}}\left(u\right) =$$

$$= \frac{\alpha^2 n^2}{2L_{i_{k+1}}}\left|\partial f\left(x_{k+1}\right)/\partial x_{i_{k+1}}\right|^2 + V_{z_k}\left(u\right) - V_{z_{k+1}}\left(u\right) \overset{\text{шаг 3}}{\le} \alpha^2 n^2\left(f\left(x_{k+1}\right) - f\left(y_{k+1}\right)\right) + V_{z_k}\left(u\right) - V_{z_{k+1}}\left(u\right).$$

Возьмем от этого неравенства условное математическое ожидание $E_{i_{k+1}}\left[\,\cdot\,\middle|i_1,...,i_k\right]$:

$$\alpha\left\langle \nabla f\left(x_{k+1}\right), z_k - u\right\rangle \le \alpha^2 n^2\left(f\left(x_{k+1}\right) - E_{i_{k+1}}\left[f\left(y_{k+1}\right)\middle|i_1,...,i_k\right]\right) + V_{z_k}\left(u\right) - E_{i_{k+1}}\left[V_{z_{k+1}}\left(u\right)\middle|i_1,...,i_k\right].$$

Согласно формуле (3.2) [135], которая используется в совершенно таком же виде, как и в [135], из последнего неравенства при

$$\frac{1-\tau}{\tau} = \alpha n^2,$$

получаем

$$\alpha\left\langle \nabla f\left(x_{k+1}\right), x_{k+1} - u\right\rangle \le \alpha^2 n^2\left(f\left(y_k\right) - E_{i_{k+1}}\left[f\left(y_{k+1}\right)\middle|i_1,...,i_k\right]\right) + V_{z_k}\left(u\right) - E_{i_{k+1}}\left[V_{z_{k+1}}\left(u\right)\middle|i_1,...,i_k\right].$$

Положим $u = x_*$, и возьмем математическое ожидание $E_{i_1,...,i_k}\left[\,\cdot\,\right]$ (если $k \ge 1$) от каждого такого неравенства, и просуммируем то что получается по $k = 0,...,K-1$[5]

$$\alpha K\left(E\left[f\left(\bar{x}_K\right)\right] - f_*\right) \le \alpha\sum_{k=1}^{K} E\left[\left\langle \nabla f\left(x_k\right), x_k - x_*\right\rangle\right] \le$$

$$\le \alpha^2 n^2\left(f\left(x_0\right) - E\left[f\left(y_K\right)\right]\right) + V_{x_0}\left(x_*\right) - V_{z_K}\left(x_*\right) \le \alpha^2 n^2\left(f\left(x_0\right) - f_*\right) + V_{x_0}\left(x_*\right),$$

$$\boxed{\bar{x}_K = \bar{x}_K\left(x_0\right) = \frac{1}{K}\sum_{k=1}^{K} x_k\,.}$$

Пусть

$$V_{x_0}\left(x_*\right) \le \Theta,$$

$$\boxed{f\left(x_0\right) - f_* \le d\,.}$$

Выбирая

---

приложениях это (использовать сразу блок случайно выбранных компонент, причем в понятие “случайно” тут можно много что вкладывать) бывает полезно [40].

[5] Из этого неравенства устанавливается прямо-двойственная природа метода ACRCD [261]. Прямо-двойственность ускоренных покомпонентных методов требуется в ряде приложений, см., например, [138] (и подраздел 5.1.5 этого раздела 5.1). Впрочем, в сильно выпуклом случае (к которому все можно сводить за дополнительную логарифмическую плату, см. далее) прямо-двойственность оказывается уже не нужна (см., например, главу 3 [181], раздел 3.2 главы 3 и подраздел 5.1.5 этого раздела 5.1).

$$\boxed{\tau = \frac{1}{\alpha n^2 + 1},\ \alpha = \frac{1}{n}\sqrt{\frac{\Theta}{d}},\ K = K(d) = 8n\sqrt{\frac{\Theta}{d}}}$$

получим

$$E\left[f\left(\bar{x}_K\right)\right] - f_* \le \frac{2n\sqrt{\Theta d}}{K} \le \frac{d}{4}.$$

Для получения сходимости по вероятности, воспользуемся следующим приемом [44], о котором мы узнали от А.С. Немировского (этот прием также можно встретить в разделе 2.1 главы 2 и в разделе 4.4 главы 4). Из

$$E\left[f\left(\bar{x}_K\right)\right] - f_* \le d/4$$

по неравенству Маркова

$$X = f\left(\bar{x}_K\right) - f_* \ge 0,\ P\left(X \ge t\right) \le E\left[X\right]/t,\ t = d/2,$$

имеем

$$P\left(f\left(\bar{x}_K\right) - f_* \ge d/2\right) \le 1/2.$$

Отсюда следует, что если мы независимо (можно параллельно) запустим $\left\lceil \log_2\left(\sigma^{-1}\right)\right\rceil$ траекторий $\mathrm{ACRCD}\left(\alpha,\tau;\Theta,x_0,d\right)$ (определив тот $\bar{x}_K$, для которого значение $f\left(\bar{x}_K\right)$ будет наименьшим), то с вероятностью $\ge 1-\sigma$ хотя бы на одной траектории будем иметь

$$f\left(\bar{x}_K\right) - f_* \le d/2.$$

К сожалению, это требует расчета в $\left\lceil \log_2\left(\sigma^{-1}\right)\right\rceil$ точках значения функции $f\left(x\right)$. Расчет функции в точке $f\left(x\right)$ может быть заметно дороже стоимости одной итерации метода ACRCD. Однако это все равно не изменит по порядку оценку общего числа арифметических операций.

Итак, пусть $\mathrm{ACRCD}\left(\alpha,\tau;\Theta,x_0,d\right)$ выдал такое $\bar{x}_{K(d)}\left(x_0\right)$, что с вероятностью $\ge 1-\sigma$ имеет место неравенство

$$f\left(\bar{x}_K\left(x_0\right)\right) - f_* \le d/2.$$

Важно заметить, что при этом с вероятностью $\ge 1-\sigma$

$$V_{\bar{x}_{K(d)}(x_0)}\left(x_*\right) \le \max_{x\in Q}\left\{V_x\left(x_*\right):\ f\left(x\right) - f_* \le d/2\right\} \le \max_{x\in Q}\left\{V_x\left(x_*\right):\ f\left(x\right) - f_* \le d\right\}.$$

При этом выписанная оценка не запрещает, например, что $V_{\bar{x}_{K(d)}(x_0)}\left(x_*\right) \gg \Theta$. В этой связи для правильной работы описываемой далее процедуры перезапусков, к сожалению, необходимо переопределить $\Theta$ следующим образом (все это может существенно ухудшить итоговую оценку,[6] однако впоследствии с помощью регуляризации исходной постановки задачи мы покажем, как можно практически полностью нивелировать эту проблему)

$$\max_{x\in Q}\left\{V_x\left(x_*\right):\ f\left(x\right)-f_*\le d\right\}\le\Theta.$$

Запустим далее $\mathrm{ACRCD}\left(\alpha,\tau;\Theta,\bar{x}_K\left(x_0\right),d/2\right)$, получим такой $\bar{x}_{K(d/2)}\left(\bar{x}_K\left(x_0\right)\right)$, что с вероятностью $\ge 1-2\sigma$ (воспользовались неравенством Буля) имеет место неравенство

$$f\left(\bar{x}_{K(d/2)}\left(\bar{x}_K\left(x_0\right)\right)\right)-f_*\le d/4.$$

Процесс можно продолжать ... Для достижения по функции точности $\varepsilon$ с вероятностью

$$\ge 1-\left\lceil \log_2\left(d/\varepsilon\right)\right\rceil\sigma$$

будет достаточно $\left\lceil \log_2\left(d/\varepsilon\right)\right\rceil$ таких итераций (перезапусков). Требуемое при этом общее число итераций (число обращений к компонентам вектора градиента) оценивается сверху следующим образом:[7]

$$N\le 8n\sqrt{\frac{\Theta}{\varepsilon}}\left(1+2^{-1/2}+2^{-1}+2^{-3/2}+\ldots\right)\log_2\left(\sigma^{-1}\right)<30n\sqrt{\frac{\Theta}{\varepsilon}}\log_2\left(\sigma^{-1}\right).$$

**Теорема 5.1.2.** *После* $\left\lceil \log_2\left(d/\varepsilon\right)\right\rceil$ *описанных выше рестартов, метод ACRCD выдает такой* $x^N$, *что с вероятностью* $\ge 1-\sigma$ *имеет место неравенство*

$$f\left(x^N\right)-f_*\le\varepsilon.$$

*При этом методу требуется для этого сделать*

$$N=N\left(\varepsilon\right)=27n\sqrt{\frac{\Theta}{\varepsilon}}\log_2\left(\frac{\log_2\left(d/\varepsilon\right)}{\sigma}\right)$$

*итераций.*

---

[6] Например, для выпуклой квадратичной функции (безобидной, с первого взгляда, в виду равномерной ограниченности всех коэффициентов), о которой нам сообщил Ю.Е. Нестеров: $f\left(x\right)=x_1^2+\sum_{k=1}^{n-1}\left(x_{k+1}-2x_k\right)^2$ множества Лебега оказываются сильно сплющенными (плохо обусловленными). Также эту функцию интересно прооптимизировать с помощью ПГМ методов с не евклидовой нормой [135]. Эти методы релаксационные (то есть значение функции монотонно убывает на итерациях), но при этом точки, генерируемые методами, по ходу итерационного процесса могут уходить намного дальше от решения, чем точка старта.

[7] Выбирая $K=9n\sqrt{\Theta/d}$, можно уменьшить константу 30 до 27, последняя константа уже не улучшаемая при таком способе рассуждений.

**Замечание 5.1.1. (достоинства и недостатки ACRCD)** Из описанной конструкции ACRCD, как уже отмечалось ранее, следует его прямо-двойственность [261]. Заметив, что ПГМ и МЗС (Grad и Mirr) легко могут быть обобщены на композитные постановки задач [266] (к этому случаю, в частности, можно свести и минимизацию на параллелепипеде), с помощью описанной выше конструкции можно получить соответствующий композитный вариант ACRCD для задач с сепарабельным композитом (см. также [199]). К сожалению, ряд других свойств ACRCD уже не так просто "вытащить" (и не всегда понятно даже, возможно ли это в принципе, и имеют ли вообще нужные свойства место здесь). В частности, например, не понятно, как можно адаптивно (по ходу итерационного процесса) подбирать константы Липшица по разным направлениям, подобно п. 6.1 работы [264]. Не понятно, как можно "бороться" с проблемой неизвестности одновременно двух параметров $\Theta$ и $d$, нужных методу для работы.[8] Описанный нами вариант метод ACRCD работает в предположении $Q = \mathbb{R}^n$, и не гарантирует свойство равномерной ограниченности (в вероятностях категориях [40, 44]) последовательности расстояний от решения до точек генерируемых методом, значением этого расстояния в начальный момент, умноженным не небольшую универсальную константу (в частности, не зависящую от свойств функционала задачи). Эти плохие свойства ACRCD[9] "унаследовал" от БГМ, описанного в п. 3 [135]. В п. 4 [135] описан вариант немного другой вариант БГМ, который лишен этих недостатков. Оказывается можно распространить и его на покомпонентный случай, что далее будет сделано (см. замечание 5.1.2).

Уже отмеченные свойства (прямо-двойственность, обобщение на композитные задачи) и все далее изложенные свойства (обобщения) ACRCD допускают всевозможные сочетания друг с другом. Детали мы вынуждены опустить (планируется посвятить этому отдельную работу), но, в большинстве случаев, все это является довольно простыми фактами (впрочем, как правило, требующие для аккуратного доказательства довольно громоздких, но вполне стандартных рассуждений). Можно сказать по-другому: далее приводится "базис" для всевозможных последующих обобщений.

---

[8] Стандартные приемы рестартартов по неизвестному параметру разработаны сейчас только для случая одного неизвестного параметра [44], формальная попытка перенесения на случай двух и более неизвестных параметров (без дополнительных предположений [40]) приводит к резкому увеличению сложности процедуры.

[9] Вместе с уже отмеченной ранее проблемой вхождения $\Theta$ в итоговую оценку (числа итераций $N(\varepsilon)$) вместо $V_{x_0}(x_*)$, как это можно было ожидать [199].

**Замечание 5.1.2. (ACRCD*)** Используемая при построении ACRCD техника рестартов позволила довольно просто получить оценки вероятностей больших уклонений. Однако эта же техника создала ряд проблем (см. замечание 5.1.1), многие из которых, в первую очередь, связаны с некоторым запаздыванием в обновлении параметров $\tau$ и $\alpha$. Они обновляются только на новом рестарте. Основная идея (см. п. 4 [135]) – сделать эти параметры зависящими от шага. Тогда удастся избавиться от рестартов и приобрести ряд хороших свойств. Далее описывается соответствующая модификация метода ACRCD. Предварительно определим две числовые последовательности

$$\alpha_1 = \frac{1}{n^2},\ \alpha_k^2 n^2 = \alpha_{k+1} n^2 - \alpha_{k+1},\ \tau_k = \frac{1}{\alpha_{k+1} n^2}.$$

Можно написать явные формулы. Также можно, следуя [135], брать близкие последовательности (теоретически немного проще исследовать первый вариант, но второй вариант более нагляден, а итоговые оценки скорости сходимости практически идентичны)

$$\boxed{\alpha_{k+1} = \frac{k+2}{n^2},\ \tau_k = \frac{1}{\alpha_{k+1} n^2} = \frac{2}{k+2}.}$$

$$\text{ACRCD*}\left(x_0 = y_0 = z_0\right)$$

1. $x_{k+1} = \tau_k z_k + \left(1-\tau_k\right) y_k$;
2. $i_{k+1} \in \left[1,...,n\right]$ – независимо от предыдущих розыгрышей;
3. $y_{k+1} = \mathrm{Grad}_{i_{k+1}}\left(x_{k+1}\right)$;
4. $z_{k+1} = \mathrm{Mirr}_{z_k}\left(\alpha_{k+1} n \nabla_{i_{k+1}} f\left(x_{k+1}\right)\right)$.

Оценка скорости сходимости такого метода

$$\boxed{N\left(\varepsilon\right) = \mathrm{O}\left(n\sqrt{\frac{\Theta}{\varepsilon}}\ln\left(\frac{1}{\sigma}\right)\right),}$$

где

$$\boxed{\Theta = V_{x_0}\left(x_*\right).}$$

Получается эта оценка из следующей формулы (см. последнюю формулу в доказательстве леммы 4.3 [135], см. также раздел 3.2 главы 3)

$$\alpha_{k+1}^2 n^2 E_{i_{k+1}}\left[f\left(y_{k+1}\right)\middle| i_1,...,i_k\right] - \left(\alpha_{k+1} n^2 - \alpha_{k+1}\right) f\left(y_k\right) \le$$

$$\le \alpha_{k+1}\left\{f\left(x_{k+1}\right) + \left\langle \nabla f\left(x_{k+1}\right), u - x_{k+1}\right\rangle\right\} + V_{z_k}\left(u\right) - E_{i_{k+1}}\left[V_{z_{k+1}}\left(u\right)\middle| i_1,...,i_k\right].$$

Взяв математическое ожидание $E_{i_1,\dots,i_k}\left[\cdot\right]$ (если $k \ge 1$) от каждого такого неравенства, и просуммировав то что получается по $k = 0,\dots,N-1$, получим

$$\alpha_N^2 n^2 E_{y_N}\left[f\left(y_N\right)\right] \le \min_{u\in Q}\left\{\sum_{k=0}^{N-1}\alpha_{k+1}E_{x_{k+1}}\left[f\left(x_{k+1}\right)+\left\langle\nabla f\left(x_{k+1}\right),u-x_{k+1}\right\rangle\right]+V_{z_0}\left(u\right)-E_{z_N}\left[V_{z_N}\left(u\right)\right]\right\} \le$$

$$\le\left(\sum_{k=0}^{N-1}\alpha_{k+1}\right)f_* + V_{z_0}\left(x_*\right)-E_{z_N}\left[V_{z_N}\left(x_*\right)\right].$$

Из последнего неравенства получается нужная оценка $N\left(\varepsilon\right)$ (только для сходимости в среднем, для получения оценки вероятностей больших уклонений необходимо использовать[10] неравенство концентрации Азума–Хефдинга для последовательностей мартингал-разностей, см., например, главу 7 [181], раздел 1.5 главы 1 и раздел 6.1 главы 6). Приведенное неравенство сразу показывает прямо-двойственность метода [261, 269], что это означает (и какая от этого польза), хорошо можно продемонстрировать конкретными примерами [33, 40, 42, 138] (см. также пример 5.1.3 ниже). Также из приведенной оценки следует, что

$$E_{z_k}\left[\frac{1}{2}\left\|z_k - x_*\right\|^2\right] = E_{z_k}\left[V_{z_k}\left(x_*\right)\right] \le E_{z_{k-1}}\left[V_{z_{k-1}}\left(x_*\right)\right] \le \dots \le E_{z_{k-1}}\left[V_{z_1}\left(x_*\right)\right] \le V_{z_0}\left(x_*\right).$$

Можно привести и более точные вероятностные оценки на субмартингал

$$\left\|z_k - x_*\right\|^2.$$

Аналогичными субмартингальными свойствами обладают и последовательности[11]

$$\left\|y_k - x_*\right\|^2,\ \left\|x_k - x_*\right\|^2,$$

что доказывается по индукции, исходя из подразделов 5.1.1, 5.1.3 в определении ACRCD*, выпуклости квадрата нормы и неравенства Йенсена. Детали мы вынуждены здесь опустить. В качестве “сухого остатка” можно сформулировать следующий результат (см. также [40, 44], раздел 1.5 главы 1, раздел 2.1 главы 2, раздел 3.2 главы 3 и раздел 4.4 главы 4): с вероятностью $\ge 1-\sigma$

$$\boxed{\max_{k=1,\dots,N}\left\{\left\|y_k - x_*\right\|^2, \left\|z_k - x_*\right\|^2, \left\|x_k - x_*\right\|^2\right\} \le R^2,}$$

где

[10] Впрочем, можно получить неравенства на вероятности больших уклонений с помощью неравенства Маркова подобно тому, как это было описано выше в подразделе 5.1.4 для ACRCD.

[11] При доказательстве этого факта существенно используется евклидовость нормы, к сожалению, для не евклидовых норм похоже, что результат перестает быть верным (см., например, приложение В.1 в [135]) если по-прежнему исходить из ПГМ в представлении БГМ (не пытаясь заменить его, например, на ППГМ).

$$\boxed{R^2 = CV_{z_0}\left(x_*\right)\ln\left(N/\sigma\right),}$$

а $C < 100$ – некоторая универсальная константа. По сути это означает, что если заранее знать $V_{z_0}\left(x_*\right)$, то, например, константы Липшица можно определять не на всем $Q$ (если $Q$ не ограничено, то и константы могут быть не ограничены), а на пересечении $Q$ с “шаром” в $\|\;\|$-норме с центром в точке $x_0$ и радиуса $R$. Вместе с прямо-двойственностью, это оказывается полезным инструментом для использования покомпонентных методов при решении двойственных задач [138] (см. также пример 5.1.3 ниже). Из описания ACRCD* также следует, что метод позволяет адаптивно подбирать константы Липшица по разным направлениям, подобно п. 6.1 работы [264].[12] Теперь уже нет проблемы с завышенной оценкой параметра $\Theta$, входящего в оценку $N\left(\varepsilon\right)$, в виду $\Theta = V_{x_0}\left(x_*\right)$. И из двух потенциально неизвестных априорно параметров $\Theta$ и $d$ теперь остается только один $\Theta$.

Можно показать, что в выписанной формуле для $N\left(\varepsilon\right)$ константа в $O(\;)$ не больше, чем в теореме 5.1.2. Используя эту явную формул для $N\left(\varepsilon\right)$, подобно п. 5 работы [135], с помощью техники рестартов (по расстоянию от текущей точки до решения) можно перенести полученные результаты на случай $\mu$-сильно выпуклой в норме $\|\;\|$ функции (заметим, что при таком перенесении можно сохранить возможность метода адаптивно настраиваться на константы Липшица). Соответствующая оценка числа итераций будет иметь следующий вид (для евклидовой нормы $\|\;\|$, для неевклидовой под корнем может возникнуть дополнительный логарифмический по $n$ множитель [44])

$$\boxed{N\left(\varepsilon\right) = O\left(n\sqrt{\frac{1}{\mu}\ln\left(\frac{\ln\left(\mu\Theta/\varepsilon\right)}{\sigma}\right)}\ln\left(\frac{\mu\Theta}{\varepsilon}\right)\right).}$$

**Замечание 5.1.3. (обобщение на блочно-компонентные методы и на более общие прокс-структуры и множества)** Описанный метод допускает следующее обобщение.[13] Пусть (см. также теорему[14] 5 [264])

---

[12] Чтобы сохранить дешевизну итерации (эффективность метода) для композитных постановок или в случае когда множество $Q$ не параллелепипедного типа здесь требуются некоторые оговорки (подобные сделанным в замечаниях 5.1.7, 5.1.8 ниже) о возможности эффективно пересчитывать значения функции (см. также главу 4).

[13] Мы не будем подробно пояснять все используемые далее обозначения – они стандартны и должны быть понятны из контекста, детали см., например, в [200, 264].

$$x=(x_1,...,x_n),\ Q=\prod_{i=1}^{n} Q_i .$$

Каждый $x_i \in Q_i$, в свою очередь, является вектором (размерности у этих векторов могут быть разными). Пусть в соответствующих подпространствах (отвечающих различным блокам) введены нормы $\left\{\sqrt{L_i}\left\|x_i\right\|_i\right\}_{i=1}^{n}$ и соответствующие этим нормам "расстояния" Брегмана $\left\{V_{x_i}^i(y_i)\right\}_{i=1}^{n}$ (см., например, [44, 181]). Положим

$$\|x\|^2=\sum_{i=1}^{n} L_i\|x_i\|_i^2,\ \|\nabla f(x)\|_*^2=\sum_{i=1}^{n} L_i^{-1}\left\|\mathrm{grad}_{x_i} f(x)\right\|_{i,*}^2,\ V_x(y)=\sum_{i=1}^{n} V_{x_i}^i(y_i).$$

Будем считать, что для всех $x, x+h\tilde{e}_i \in Q$

$$\left\|\mathrm{grad}_{x_i} f(x+h\tilde{e}_i)-\mathrm{grad}_{x_i} f(x)\right\|_{*,i} \le L_i h\left\|[\tilde{e}_i]_i\right\|_i,$$

где вектор $\tilde{e}_i$ имеет все нули в компонентах, не соответствующих $i$-му блоку. Введенные обозначения позволяют переписать сам метод (см. также [135]). При этом оценки будут иметь точно такой же вид, меняется только интерпретация параметров, норм (расстояний) в этих оценках. Заметим, что в приложениях к поиску равновесий в популяционных играх загрузок с большим числом популяций (в частности, задачах поиска равновесного распределения потоков по путям в графе транспортной сети [33, 40, 42], см. также главу 1), часто возникают множества имеющие вид прямого произведения симплексов. До настоящего момента было не понятно, можно ли (а если можно, то как) применять к таким задачам покомпонентные методы.

**Замечание 5.1.4. (обобщение на задачи стохастической оптимизации)** Предположим, что исходная задача имеет вид

$$E_{\xi}\left[f(x;\xi)\right]\to \min_{x\in Q}.$$

Если $f(x;\xi)$ – выпуклая по $x$ функция (при всех $\xi$) с константами Липшица (равномерно не только по $x, x+he_i \in Q$, но и по $\xi$)

$$\left\|\mathrm{grad}_{x_i} f(x+h\tilde{e}_i;\xi)-\mathrm{grad}_{x_i} f(x;\xi)\right\|_{*,i} \le L_i h\left\|[\tilde{e}_i]_i\right\|_i,$$

где вектор $e_i$ имеет все нули в компонентах, не соответствующих $i$-му блоку. Введем

---

[14] Отметим также, что везде в этой теореме можно вместо $R_1^2(x_0)$ писать $2\|x_0-x_*\|_1^2$, что немного улучшает оценку теоремы.

$$D=\max_{x\in Q}E_{\xi}\left[\left\|\nabla f\left(x;\xi\right)-E_{\xi}\left[\nabla f\left(x;\xi\right)\right]\right\|_{*}^{2}\right].$$

Тогда если вместо $\nabla_i f\left(x;\xi\right)$ можно рассчитывать только на $\nabla_i f\left(x;\xi\right)$:

$$E_{\xi}\left[\nabla_i f\left(x;\xi\right)\right]\equiv\nabla_i E_{\xi}\left[f\left(x;\xi\right)\right],$$

то оценка в теореме 5.1.2 изменится следующим образом (см. также [44] и разделы 2.1, 2.2 главы 2)

$$\boxed{N\left(\varepsilon\right):=\mathrm{O}\left(\max\left\{N\left(\varepsilon\right),n\frac{D\Theta}{\varepsilon^{2}}\ln\left(\frac{1}{\sigma}\right)\right\}\right),}$$

в случае $\mu$-сильно выпуклой в норме $\|\ \|$ функции

$$\boxed{N\left(\varepsilon\right):=\mathrm{O}\left(\max\left\{N\left(\varepsilon\right),n\frac{D}{\mu\varepsilon}\ln\left(\frac{\ln\left(N\left(\varepsilon\right)\right)}{\sigma}\right)\right\}\right).}$$

**Замечание 5.1.5. (учет ошибок в вычислении компонент градиента)** Немного специфицируя концепцию неточного оракула ($\delta\geq 0$ – уровень шума) из главы 4 [181], введем следующее предположение (векторы $\tilde{e}_i$, $\nabla_i f_{\delta,L_i}\left(x\right)$ имеет все нули в компонентах, не соответствующих $i$-му блоку): для любого $x\in Q$ существуют такие $f_{\delta,L_i}\left(x\right)$ и $\nabla_i f_{\delta,L_i}\left(x\right)$, что для всех $y=x+h\tilde{e}_i\in Q$ выполняется

$$0\leq f\left(y\right)-f_{\delta,L_i}\left(x\right)+\left\langle\nabla_i f_{\delta,L_i}\left(x\right),y-x\right\rangle\leq\frac{L_i}{2}\left\|\left[y-x\right]_i\right\|_i^{2}+\delta,$$

т.е.

$$0\leq f\left(x+h\tilde{e}_i\right)-f_{\delta,L_i}\left(x\right)+h\left\langle\nabla_i f_{\delta,L_i}\left(x\right),\tilde{e}_i\right\rangle\leq\frac{L_i h^{2}}{2}\left\|\left[\tilde{e}_i\right]_i\right\|_i^{2}+\delta.$$

При $\delta=0$ отсюда получаем определение констант Липшица $L_i$ в блочно-покомпонентном методе из замечания 5.1.3. Полезными следствиями введенного определения являются следующие неравенства

$$0\leq f\left(x\right)-f_{\delta,L_i}\left(x\right)\leq\delta,$$

$$\frac{1}{2}\left\|\nabla_i f_{\delta,L_i}\left(x\right)\right\|_{*}^{2}\leq f\left(x\right)-f\left(\mathrm{Grad}_i^{\delta}\left(x\right)\right)+\delta,$$

где

$$\mathrm{Grad}_i^{\delta}\left(x\right)=\arg\min_{\tilde{x}\in Q}\left\{\left\langle\nabla_i f_{\delta,L_i}\left(x\right),\tilde{x}-x\right\rangle+\frac{1}{2}\left\|\tilde{x}-x\right\|^{2}\right\}=x-\frac{1}{L_i}\nabla_i f_{\delta,L_i}\left(x\right),$$

здесь для наглядности мы ограничились случаем евклидовой нормы и $Q = \mathbb{R}^n$. С помощью этих неравенств можно скорректировать (см. также [44] и раздел 2.1 главы 2) оценку теоремы 5.1.2 (и различные ее обобщения) на случай, когда мы можем вычислять вместо “честных” компонент градиента $\nabla_i f(x)$ только их приближенные (в указанном выше смысле) аналоги $\nabla_i f_{\delta,L_i}(x)$:

$$\boxed{f\left(x^{N(\varepsilon)}\right) - f_* \le \varepsilon + \mathrm{O}\left(N(\varepsilon)\delta\right).}$$

Вообще эта формула типична для всех известных нам ускоренных методов (полноградиентных, покомпонентных, прямых). И это соответствует самому худшему (быстрому) варианту накопления ошибки. Для неускоренных методов $\mathrm{O}(N\delta) \to \mathrm{O}(\delta)$, что соответствует самому лучшему варианту накопления ошибки (то есть когда такое накопление отсутствует). Подобно полноградиентным методам для покомпонентных методов (и прямых) можно предложить так называемы промежуточные методы (см., например, [44, 181, 193]) с накоплением ошибки $\mathrm{O}\left(N^p\delta\right)$, $p \in [0,1]$. Из этого замечания (а также замечания 5.1.2) возникает гипотеза о возможности создания универсального покомпонентного метода [274] (см. также разделы 2.1, 2.2 главы 2).

Интересно было бы объединить замечания 5.1.4, 5.1.5 с целью получения покомпонентной версии результатов главы 7 [181] и [193].

**Замечание 5.1.6. (обобщение на взвешенную рандомизацию)** Предположим, что вместо одной компоненты можно выбирать блок компонент (всего $n$ блоков), причем, вообще говоря, с разными вероятностями: выбираем блок компонент $i$ с вероятностью[15]

$$\boxed{p_i = \frac{L_i^\beta}{\sum_{j=1}^n L_j^\beta},\ i = 1,\ldots,n,}$$

где параметр степени $\beta \in [0,1]$. При этом необходимо будет переопределить норму

$$\boxed{\|x\|^2 = \sum_{i=1}^n L_i x_i^2 \to \sum_{i=1}^n L_i^{1-2\beta}\|x_i\|_i^2,}$$

---

[15] Приготовление памяти для генерирования из описанного распределения стоит $\mathrm{O}(n)$. Это делается один раз (строится соответствующее двоичное дерево Л.В. Канторовича [264] – см. разделы 4.1, 4.2 главы 4). Случайные разыгрывания $i$ при наличии правильно подготовленной памяти будут стоить $\mathrm{O}(\log_2 n)$ – каждое.

а, соответственно, также прокс-функцию и параметр $\Theta$. При этом во всех приведенных выше формулах, которые определяют метод, в частности, для ACRCD это

$$\tau = \frac{1}{\alpha n^2 + 1},\ \alpha = \frac{1}{n}\sqrt{\frac{\Theta}{d}},\ K = 8n\sqrt{\frac{\Theta}{d}},\ N = 27n\sqrt{\frac{\Theta}{\varepsilon}}\log_2\left(\frac{\log_2(d/\varepsilon)}{\sigma}\right),$$

необходимо будет сделать замену

$$\boxed{n \to \sum_{i=1}^{n} L_i^{\beta}.}$$

Выше в этом пункте мы рассматривали случай $\beta = 0$ (в другом ключе этот случай также рассматривался в [264]). Можно ожидать, что на практике этот вариант предпочтительнее. Заметим, что ранее уже рассматривались отдельно случаи $\beta = 1/2$ [282] и $\beta = 1$ [245].

В связи с замечанием 5.1.6 возникает вопрос: существуют ли еще более общие способы (с большим числом степеней свободы) сочетания выбора рандомизации и нормы? Положительный ответ, более менее, очевиден (см. также [295]), но интересно было бы предложить такие способы, которые в определенных ситуациях позволяли бы еще более ускориться по сравнению с методом, порожденным замечаниями 5.1.2, 5.1.6.

**Замечание 5.1.7. (стоимость итерации / неразреженный случай)** Рассмотрим, следуя Ю.Е. Нестерову [282]), следующий случай

$$f(x) = F(Ax, x),\ x \in \mathbb{R}^n,\ y = Ax \in \mathbb{R}^m.$$

Будем считать, что значение $F(y,x)$ (а, следовательно, и градиент [76]) можно посчитать за $O(m+n)$. Пусть верно хотя бы одно из следующих условий: 1) $n = O(m)$; 2) расчет $\mathrm{grad}_y F(y,x)$ стоит $O(m)$, а $\partial F(y,x)/\partial x_j$ – $O(m)$.[16] Тогда амортизационная (средняя) сложность одной итерации будет $O(m)$. Обоснование этого факта можно получить как простое следствие более общих рассуждений, проводимых в следующем замечании.

**Замечание 5.1.8. (стоимость итерации / разреженный случай)** Из первого пункта описания алгоритма ACRCD кажется, что всегда один шаг этого алгоритма будет требо-

---

[16] В ряде приложений посчитать одну компоненту градиента оказывается в $n$ раз дешевле, чем сам градиент. Например, это так для $f(x) = x^T Ax$. Но все же верно это далеко не всегда. Например, для функции (см. пример 5.1.3 подраздела 5.1.5 ниже)

$$f(x) = \ln\left(\sum_{k=1}^{n}\exp(x_k)\right)$$

стоимость расчета самой функции, ее градиента и любой компоненты градиента одинаковы по порядку.

вать, как минимум $\geq n$ арифметических операций. Однако замечание 5.1.7 показывает (при $m \ll n$), что это совсем не обязательно. Естественно возникает вопрос: а можно ли получить еще больше (еще более дешевую итерацию)? В определенных (разреженных задачах специальной структуры) ответ оказывается положительным. Пояснению этого тезиса и будет посвящена оставшаяся часть данного замечания. Оказывается, что при наличии у задачи определенной структуры (например, в случае $f(x) = x^T Ax$ или $f(x) = \|Ax - b\|_2^2$), нет необходимости выполнять первый пункт честно (в полном объеме). Далее мы описываем идею, заимствованную из работ [199, 245, 282], близкие идеи мы использовали в разделах 4.2, 4.3 главы 4 для перезаписи метода условного градиента. Предварительно перепишем алгоритм ACRCD в рассматриваемом нами случае (рассматривается задача безусловной оптимизации) следующим (эквивалентным) образом ($x_0 = u_0 = v_0$)

$$\mathrm{ACRCD}'\left(\alpha, \tau; \Theta, x_0, f(x_0) - f_* \leq d\right)$$

1. $x_{k+1} = (1-\tau)^{k+1} v_k + u_k$;
2. $i_{k+1} \in [1, \ldots, n]$ – независимо от предыдущих розыгрышей;
3. $v_{k+1} = v_k + \dfrac{1}{L_{i_{k+1}}} \dfrac{\alpha n - 1}{(1-\tau)^k} \nabla_{i_{k+1}} f(x_{k+1})$;
4. $u_{k+1} = u_k - \dfrac{\alpha n}{L_{i_{k+1}}} \nabla_{i_{k+1}} f(x_{k+1})$.
5. если $k = N - 1$, то выдаем
$y_N = y_{k+1} = \mathrm{Grad}_{i_{k+1}}(x_{k+1}) = \mathrm{Grad}_{i_N}(x_N)$.

Предположим, что

$$f(x) = \sum_{r=1}^{m} \varphi_r\left(a_r^T x\right),$$

где функции $\varphi_r$ – простой структуры, т.е. дифференцируемые за $O(1)$ каждая, причем $A = \|a_1 \ldots a_m\|^T$ – разреженная матрица (число ненулевых элементов $sn$, т.е. в каждом столбце в среднем $s \ll m$ ненулевых элементов). Тогда подобно [199] один шаг метода $\mathrm{ACRCD}'$ (кроме самого первого шага, который может стоить $O(n)$) может быть осуществлен в среднем за $O(s)$ (амортизационная сложность). Действительно, пункт 2 $\mathrm{ACRCD}'$ можно осуществить за $O(\ln n)$ (считаем $\ln n = O(s)$). Пункты 3, 4 за $O(s)$. Та-кая стоимость этих пунктов обусловлена необходимостью пересчета $\partial f(x_{k+1})/\partial x_i$. Осуществлять этот пересчет необходимо, используя пункт 1. Покажем, как можно это эффективно делать. Прежде всего, заметим, что если мы уже посчитали $a_r^T v_k$ и $a_r^T u_k$, то посчитать дополнительно $a_r^T x_{k+1}$ будет стоить $\mathrm{O}(1)$. Также заметим, что если мы уже посчитали $Av_k$ и $Au_k$, то посчитать дополнительно $Av_{k+1}$ и $Au_{k+1}$ будет стоить $\mathrm{O}(s)$. Чтобы посчитать $\partial f(x_{k+1})/\partial x_i$ нужно вычислить частные производные (по $x_i$) в среднем у $s$ слагаемых в сумме. Каждая такая частная производная (по предположению о простоте структуры функций $\varphi_r$) рассчитывается за $\mathrm{O}(1)$, в предположении известности всех аргументов этих функций. Учитывая, что на пересчет всех аргументов уйдет $\mathrm{O}(s)$ (см. выше), то в среднем, общие трудозатраты будут $\mathrm{O}(s)+\mathrm{O}(s)=\mathrm{O}(s)$. Таким образом, один шаг метода будет стоить $\mathrm{O}(s)$. Чтобы посчитать выход алгоритма:

$$\bar{x}_K = \frac{1}{K}\sum_{k=1}^{K} x_k = \frac{1}{K}\sum_{k=1}^{K}\sum_{\tilde{K}=k}^{K}\left(\frac{\alpha n-1}{L_{i_k}}\frac{(1-\tau)^{\tilde{K}}}{(1-\tau)^{k}} - \frac{\alpha n}{L_{i_k}}\right)\nabla_{i_k} f(x_k) =$$

$$= \sum_{k=1}^{K}\frac{1}{K}\left(\frac{\alpha n-1}{L_{i_k}}\frac{1-\tau}{\tau}\left(1-(1-\tau)^{K-k}\right) - \frac{\alpha n}{L_{i_k}}(K-k+1)\right)\nabla_{i_k} f(x_k),$$

достаточно $\mathrm{O}(K)$ арифметических операций. Все написанное выше переносится и на ACRCD*.

$$\text{ACRCD*}'(x_0 = u_0 = v_0 = w_0)$$

1. $x_{k+1} = \left\{\prod_{j=1}^{k}(1-\tau_j)\right\}\cdot\left[u_k + v_k + \left\{\sum_{j=1}^{k}\frac{\tau_j}{\prod_{l=1}^{j}(1-\tau_l)}\right\}w_k\right]$;
2. $i_{k+1} \in [1,\ldots,n]$ – независимо от предыдущих розыгрышей;
3. $u_{k+1} = u_k + \frac{\alpha_{k+1} n}{L_{i_{k+1}}}\sum_{j=1}^{k}\frac{\tau_j}{\prod_{l=1}^{j}(1-\tau_l)}\nabla_{i_{k+1}} f(x_{k+1})$;
4. $v_{k+1} = v_k - \frac{1}{L_{i_{k+1}}}\frac{1}{\prod_{j=1}^{k}(1-\tau_j)}\nabla_{i_{k+1}} f(x_{k+1})$;

| |
|---|
| 5. $w_{k+1} = w_k - \frac{\alpha_{k+1} n}{L_{i_{k+1}}} \nabla_{i_{k+1}} f\left(x_{k+1}\right)$; |
| 6. если $k = N-1$, то выдаем |
| $y_N = y_{k+1} = \mathrm{Grad}_{i_{k+1}}\left(x_{k+1}\right) = \mathrm{Grad}_{i_N}\left(x_N\right)$. |

### 5.1.5 Примеры применения ускоренных покомпонентных методов

Начнем этот пункт с очень простого примера, демонстрирующего, что покомпонентные методы (в том числе неускоренные), вообще говоря, не применимы к произвольной задаче выпуклой оптимизации.

**Пример 5.1.1.** Рассмотрим выпуклую задачу

$$f(x) = (x_1 - 2)^2 + (x_2 - 1)^2 \to \min_{x \in Q},$$

$$Q = \left\{x = (x_1, x_2) \geq 0: \quad x_1 + x_2 \leq 2\right\}.$$

Предположим, что покомпонентный метод стартует с точки

$$x^0 = (1, 1) \in Q.$$

Тогда метод (если в методе жестко прописано оптимизировать по выбранному направлению – у нас это не так) не сдвинется с места по обоим направлениям, задаваемым ортами и проходящим через $x^0$ (поскольку $f(x)$ внутри $Q$ имеет минимум в $x^0$ по этим направлениям), в то время как в точке

$$x^* = (0.5, 1.5) \in Q$$

$f(x)$ достигает минимума на $Q$:

$$f\left(x^*\right) = 0.5 < 1 = f\left(x^0\right).$$

Если в методе жестко прописано оптимизировать по выбранному направлению (в подавляющем большинстве существующих вариантов покомпонентных методов именно так и сделано), то для возможности рассмотрения случая $Q \neq \mathbb{R}^n$ нужно дополнительно потребовать, чтобы для любого $i = 1, ..., n$ имело место условие: для любого $x = (x_1, ..., x_n) \in Q$ выполняется

$$(x_1, ..., x_{i-1}, x_{*i}, x_{i+1}, ..., x_n) \in Q.$$

Если множество $Q = \prod_{k=1}^{n} [a_k, b_k]$ – сепарабельно, то это условие, очевидно, выполняется. □

Перейдем теперь непосредственно к приложениям покомпонентных методов. За более подробной информацией о покомпонентных методах и примерах их приложений можно рекомендовать обратиться к [289], а также работам T. Zhang-а и P. Richtarik-а.

**Пример 5.1.2.** Возьмем функцию (в изложении этого примера мы во многом следуем Ю.Е. Нестерову, см. [282])

$$f(x) = \frac{1}{2}\langle x, Sx\rangle - \langle b, x\rangle,$$

где $S$ – симметричная матрица, все элементы которой числа от 1 до 2. Возьмем метод ACRCD* в варианте замечания 5.1.6 с $\beta = 1/2$. Выберем евклидову норму. Константа Липшица этой функции по определению есть

$$L = \lambda_{\max}(S) \geq \lambda_{\max}\left(1_n 1_n^T\right) = n.$$

При этом покомпонентный метод дает константы Липшица $L_i = S_{ii} \leq 2$. Таким образом, получаем ускорение в $\sim \sqrt{n}$ раз. Действительно (см. замечание 5.1.6),

$$n \to \sum_{i=1}^{n} \sqrt{L_i} \leq \sqrt{2}n,$$

поэтому оценка числа итераций соответствующего покомпонентного быстрого градиентного метода (ПБГМ) будет

$$N_{ПБГМ}(\varepsilon) = \mathrm{O}\left(n\sqrt{\frac{\Theta}{\varepsilon}}\ln\left(\frac{1}{\sigma}\right)\right),$$

а стоимость одной итерации $\mathrm{O}(n)$. Итого

$$T_{ПБГМ} = \mathrm{O}\left(n^2\sqrt{\frac{\Theta}{\varepsilon}}\ln\left(\frac{1}{\sigma}\right)\right).$$

Для обычного (не покомпонентного) БГМ соответствующая оценка числа итераций имеет вид (отметим, что при выборе $\beta = 1/2$ можно считать $\Theta$ в обеих формулах одинаковым)

$$N_{БГМ}(\varepsilon) = \mathrm{O}\left(\sqrt{\frac{L\Theta}{\varepsilon}}\right) = \mathrm{O}\left(\sqrt{\frac{n\Theta}{\varepsilon}}\right).$$

Зато одна итерация стоит $\mathrm{O}\left(n^2\right)$. Итого

$$T_{БГМ} = \mathrm{O}\left(n^{5/2}\sqrt{\frac{\Theta}{\varepsilon}}\right).$$

В более общем случае полезно иметь в виду следующие неравенства

$$\frac{1}{n}\operatorname{tr}(S) \le \lambda_{\max}(S) \le \operatorname{tr}(S),\ \frac{1}{n}\sum_{i=1}^{n}\sqrt{L_i} \le \sqrt{\frac{1}{n}\sum_{i=1}^{n}L_i} = \sqrt{\frac{1}{n}\operatorname{tr}(S)}.$$

Таким образом (здесь мы опустили логарифмический множитель в оценке $T_{БГМ}$, поэтому вместо $O(\ )$ ввели $\tilde{O}(\ )$),

$$T_{ПБГМ} = \tilde{O}\left(n^2\sqrt{\frac{(\operatorname{tr}(S)/n)\Theta}{\varepsilon}}\right) \le O\left(n^2\sqrt{\frac{\lambda_{\max}(S)\Theta}{\varepsilon}}\right) = T_{БГМ}.$$

В разреженном случае, согласно замечанию 5.1.8, пропорции сохраняются

$$T_{ПБГМ} = \tilde{O}\left(sn\sqrt{\frac{(\operatorname{tr}(S)/n)\Theta}{\varepsilon}}\right) \le O\left(sn\sqrt{\frac{\lambda_{\max}(S)\Theta}{\varepsilon}}\right) = T_{БГМ}.$$

Обратим внимание, что выгода в $\sim\sqrt{n}$ раз является максимально возможной. Достигается она в ситуациях, когда $\lambda_{\max}(S)$ и $\operatorname{tr}(S)$ одного порядка. Скажем, если собственные значения матрицы $S$: $\{1,\dots,n\}$, то $\lambda_{\max}(S) = n$, а $\operatorname{tr}(S) \sim n^2$, т.е. нужна большая (более резкая) асимметрия. Если под матрицей $S$ понимать гессиан функционала задачи в “худшей” (с точки зрения рассматриваемых оценок) точке, то выписанные формулы не изменятся. Однако в разреженном случае потребуются большие оговорки, чтобы можно было сполна учесть разреженность в стоимости итерации. Поскольку в приложениях довольно типично выполнение неравенства

$$\frac{1}{n}\operatorname{tr}(S) \ll \lambda_{\max}(S),$$

то из приведенных оценок следует, что во многих случаях, получается ускорить вычисления за счет использования ПБГМ вместо БГМ (не говоря уже о возможности распараллеливания [199, 295]). Как уже отмечалось, в ряде случаев это выгода может достигать $\sim\sqrt{n}$ раз. Другие примеры, когда похожие пропорции имеют место можно посмотреть в работах [40, 138]. □

**Пример 5.1.3.** Рассмотрим следующую задачу энтропийно-линейного программирования (см., например, [37, 282], а также разделы 1.1, 1.2 главы 1 и разделы 3.1, 3.2 главы 3)

$$f(x) = \sum_{i=1}^{n} x_i \ln(x_i) \to \min_{x\in S_n(1);\ Ax=b},$$

$$S_n(1) = \left\{x \in \mathbb{R}^n : x_i \ge 0,\ i = 1,\dots,n,\ \sum_{i=1}^{n} x_i = 1\right\},$$

причем будем считать (в связи с различными транспортными приложениями это представляется довольно естественным [33, 37, 40, 42]), что условие $\sum_{i=1}^{n} x_i = 1$ является следствием системы $Ax = b$. Построим двойственную задачу

$$\min_{x\in S_n(1);\, Ax=b} \sum_{i=1}^{n} x_i \ln(x_i) = \min_{x\in S_n(1)} \max_{y\in\mathbb{R}^m} \left\{ \sum_{i=1}^{n} x_i \ln(x_i) + \langle y, b - Ax\rangle \right\} =$$

$$= \max_{y\in\mathbb{R}^m} \min_{x\in S_n(1)} \left\{ \sum_{i=1}^{n} x_i \ln(x_i) + \langle y, b - Ax\rangle \right\} = \max_{y\in\mathbb{R}^m} \left\{ \langle y, b\rangle - \ln\left( \sum_{i=1}^{n} \exp\left( \left[A^T y\right]_i \right) \right) \right\}.$$

Но с учетом написанного выше, двойственную задачу можно строить и по-другому

$$\min_{x\in S_n(1);\, Ax=b} \sum_{i=1}^{n} x_i \ln(x_i) = \min_{Ax=b} \sum_{i=1}^{n} x_i \ln(x_i) = \min_{x\in\mathbb{R}_+^n} \max_{y\in\mathbb{R}^m} \left\{ \sum_{i=1}^{n} x_i \ln(x_i) + \langle y, b - Ax\rangle \right\} =$$

$$= \max_{y\in\mathbb{R}^m} \min_{x\in\mathbb{R}_+^n} \left\{ \sum_{i=1}^{n} x_i \ln(x_i) + \langle y, b - Ax\rangle \right\} = \max_{y\in\mathbb{R}^m} \left\{ \langle y, b\rangle - \sum_{i=1}^{n} \exp\left( \left[A^T y\right]_i - 1 \right) \right\}.$$

В обоих случаях решение прямой задачи можно восстановить по решению двойственной:

1) $x_i(y) = \dfrac{\exp\left( \left[A^T y\right]_i \right)}{\sum_{k=1}^{n} \exp\left( \left[A^T y\right]_k \right)}$, $i = 1,\dots,n$,

2) $x_i(y) = \exp\left( \left[A^T y\right]_i - 1 \right)$, $i = 1,\dots,n$.

Таким образом, можно работать с двумя различными двойственными задачами:

1) $\varphi_1(y) = \ln\left( \sum_{i=1}^{n} \exp\left( \left[A^T y\right]_i \right) \right) - \langle y, b\rangle \to \min_{y\in\mathbb{R}^m}$,

2) $\varphi_2(y) = \sum_{i=1}^{n} \exp\left( \left[A^T y\right]_i - 1 \right) - \langle y, b\rangle \to \min_{y\in\mathbb{R}^m}$.

Исходя из явного двойственного (Лежандрова) представления этих функционалов, можно вычислить константу Липшица градиента и соответствующие константы Липшица по направлениям. Во втором случае, к сожалению, константы получаются не ограниченными, а вот в случае 1 они ограничены,[17] и их можно оценить, соответственно, как:

[17] За счет сильной выпуклости энтропии в 1-норме с константой 1 на единичном симплексе (на положительном ортанте энтропия строго выпукла, но не сильно выпукла), из теоремы 1 [271] имеем:

$$\left\| \nabla\varphi_1(y_2) - \nabla\varphi_1(y_1) \right\|_q = \left\| Ax(y_2) - Ax(y_1) \right\|_q \le \left\| A^T \right\|_{p,1}^2 \left\| y_2 - y_1 \right\|_p, \quad \left\| A^T \right\|_{p,1}^2 = \max_{\|y\|_p \le 1,\, \|x\|_1 \le 1} \left\langle A^T y, x \right\rangle^2 = \max_{k=1,\dots,n} \left\| A^{\langle k\rangle} \right\|_q^2,$$

$$L_{БГМ} = \max_{k=1,...,n} \left\| A^{\langle k \rangle} \right\|_2^2 \text{ и } L_{ПБГМ}^k \le L_{ПБГМ} = \max_{\substack{i=1,...,m \\ j=1,...,n}} \left| A_{ij} \right|^2 \le L_{БГМ},\ k = 1,...,n .$$

Будем считать, что все элементы матрицы $A_{ij}$ удовлетворяют условию: $1 \le A_{ij} \le 2$. Тогда $L_{БГМ} \ge m$, а $L_{ПБГМ} \le 4$. Выберем в двойственном пространстве евклидову прокс-структуру.[18] Решая первую двойственную задачу БГМ (стартуем в точке 0), получим следующую оценку времени работы метода

$$T_{БГМ}^1 = \mathrm{O}\left( mn\sqrt{\frac{L_{БГМ}\Theta}{\varepsilon}} \right) = \mathrm{O}\left( mn\sqrt{\frac{m\Theta}{\varepsilon}} \right).$$

Если же применить к первой двойственной задаче ACRCD* (стартуем в точке 0) с $\beta = 1/2$ или $\beta = 0$ (см. замечания 5.1.2, 5.1.6, 5.1.7), то получим, что с вероятностью $\ge 1-\sigma$

$$T_{ПБГМ}^1 = \tilde{\mathrm{O}}\left( nm\sqrt{\frac{L_{ПБГМ}\Theta}{\varepsilon}} \right) = \tilde{\mathrm{O}}\left( mn\sqrt{\frac{\Theta}{\varepsilon}} \right).$$

В обеих формулах $\Theta$ – квадрат евклидового размера двойственного решения. Таким образом, за счет использования ПБГМ удается ускориться в $\sim \sqrt{m}$ раз. Все изложенное в этом примере распространяется и на случай когда вместо ограничений в виде равенств (или наряду с ними) мы имеем ограничения в виде неравенств $Cx \le d$. Если по-прежнему обозначать общее число ограничений через $m$, то выписанные формулы останутся справедливыми. Пока мы решили только двойственную задачу с точностью по функции $\varepsilon$. То есть в двойственном пространстве используемый нами метод сгенерировал последовательности $\{y_k\}_{k=1}^N$, $\{z_k\}_{k=1}^N$, $\{\alpha_k\}_{k=1}^N$ (см. обозначения замечания 5.1.2 с заменой $x_k \to z_k$, сделанной во избежание путаницы) такие, что (для БГМ математическое ожидание можно не писать)

$$E_{y_N}\left[ \varphi_1(y_N) \right] - \varphi_* \le \varepsilon,\ \varphi_* = f_* .$$

---

где $1/p + 1/q = 1$. Беря $p = 2$, получим константу Липшица градиента $\varphi_1(y)$: $\max_{k=1,...,n} \left\| A^{\langle k \rangle} \right\|_2^2$. Беря $p = 1$ можно получить, что константа Липшица производной $\varphi_1(y)$ по каждому направлению не больше $\max_{\substack{i=1,...,m \\ j=1,...,n}} \left| A_{ij} \right|^2$.

[18] Вообще при решении двойственных задач это совершенно естественно [37], поскольку оптимизация происходит либо на всем пространстве, либо на прямом произведении пространства и неотрицательно ортанта.

Оказывается, что не использующаяся в этой формуле и накопленная методом информация $\{z_k\}_{k=1}^N$, $\{\alpha_k\}_{k=1}^N$ позволяет восстанавливать с такой же точностью решение прямой задачи (детали см., например, в работе [138] и разделе 3.2 главы 3):

$$\bar{x}^N = \frac{1}{S_N}\sum_{k=1}^N \alpha_k x(z_k),\; S_N = \sum_{k=1}^N \alpha_k ,$$

$$\sqrt{\Theta}\left\|A\bar{x}^N - b\right\|_2 \le \varepsilon,\; \left|f\left(\bar{x}^N\right) - f_*\right| \le \varepsilon .$$

Несложно показать, что можно так организовать вычисление $\bar{x}^N$, что выписанные ранее оценки трудоемкости методов БГМ и ПБГМ (в категориях $\mathrm{O}(\;)$) не изменятся. Если теперь рассмотреть разреженный случай ($s \ll m$ – среднее число ненулевых элементов в столбце матрицы $A$, $\tilde{s} = sn/m$ – в строчке), то оценка БГМ улучшится

$$T_{БГМ}^1 = \mathrm{O}\left(sn\sqrt{\frac{L_{БГМ}\Theta}{\varepsilon}}\right),$$

в то время как оценка ПБГМ останется неизменной (см. замечание 5.1.8)

$$T_{ПБГМ}^1 = \tilde{\mathrm{O}}\left(mn\sqrt{\frac{L_{ПБГМ}\Theta}{\varepsilon}}\right).$$

Получилось это из-за наличия (при построении двойственной задачи) связывающего переменные (симплексного) ограничения – не имеющего полностью сепарабельную структуру (то есть не распадающегося в прямое произведение ограничений на отдельные компоненты). Другое дело, если мы рассмотрим вторую двойственную задачу. Она полностью сепарабельная (подходит под замечание 5.1.7 в этом смысле). То есть для второй двойственной задачи можно найти ее решение с точностью по функции $\varepsilon$:

$$E_{y_N}\left[\varphi_2\left(y_N\right)\right] - \varphi_* \le \varepsilon$$

за время (обратим внимание, что $\bar{L}_{ПБГМ}$, $\bar{\Theta}$ соответствуют функционалу и решению второй двойственной задачи, и, вообще говоря, отличаются от введенных ранее $L_{ПБГМ}$, $\Theta$)

$$T_{ПБГМ}^2 = \tilde{\mathrm{O}}\left(\tilde{s}m\sqrt{\frac{\bar{L}_{ПБГМ}\bar{\Theta}}{\varepsilon}}\right) = \tilde{\mathrm{O}}\left(sn\sqrt{\frac{\bar{L}_{ПБГМ}\bar{\Theta}}{\varepsilon}}\right).$$

При этом следует использовать метод ACRCD* из замечания 5.1.2 (если смотреть через замечание 5.1.6, то следует полагать $\beta = 0$) с адаптивным подбором констант Липшица, поскольку ограничить их не представляется возможным в виду априорного отсутствия информации о локализации решения двойственной задачи. Тем не менее, согласно заме-чанию 2 можно быть уверенным, что несмотря на неограниченность множества, на котором происходит оптимизация, и неограниченности констант Липшица на этом множестве, существуют “эффективные” константы, на которые мы адаптивно настраиваемся по ходу работы метода. К сожалению, восстановить за то же по порядку время решение прямой задачи $\bar{x}^N$ в данном случае не получается, даже если пытаться использовать соответствующие наработки замечания 5.1.8 по расчету $\bar{x}_K$. Мы снова возвращаемся к оценке типа

$$T_{ПБГМ}^2 = \tilde{O}\left( nm\sqrt{\frac{\bar{L}_{ПБГМ}\bar{\Theta}}{\varepsilon}} \right).$$

Если отказаться от ускоренности метода, то обычный неускоренный покомпонентный метод (ПМ) [163, 264] позволяет сохранить (с учетом необходимости восстановления решения прямой задачи) дешевую итерацию[19] $O(\tilde{s})$

$$T_{ПМ}^2 = \tilde{O}\left( m\tilde{s}\frac{\bar{L}_{ПБГМ}\bar{\Theta}}{\varepsilon} \right) = \tilde{O}\left( sn\frac{\bar{L}_{ПБГМ}\bar{\Theta}}{\varepsilon} \right).$$

В данном случае (и это довольно типично), этому методу скорее стоит предпочесть БГМ с итоговой оценкой

$$T_{БГМ}^2 = \tilde{O}\left( sn\sqrt{\frac{\bar{L}_{БГМ}\bar{\Theta}}{\varepsilon}} \right),$$

но в таком случае уж лучше применять БГМ к первой двойственной задачи, обладающей лучшими свойствами. Получается, что необходимость восстановления решения прямой задачи для ПБГМ накладывает дополнительные ограничения на структуру задачи, чтобы можно было полноценно воспользоваться разреженностью. К сожалению, похоже, что эти дополнительные ограничения, фактически, оставляют возможность только для задач вида (и небольших “аффинных” релаксаций, например, добавление разреженных аффинных неравенств)

$$f(x) = \frac{1}{2}\left\| x - x_g \right\|_2^2 \to \min_{Ax=b},$$

полноценно использовать в разреженном случае описанный выше подход. При этом решение прямой и двойственной задачи можно восстанавливать по описанному выше меха-низму (с учетом линейности зависимости $x(y)$ удается воспользоваться техникой пересчета $\bar{x}_K$ из замечания 5.1.8 для восстановления $\bar{x}^N$). При этом оценки, соответствующих методов БГМ и ПБГМ, применённых к двойственной задаче (и стартующих с точки 0), примут следующий вид[20]

$$\tilde{T}_{БГМ} = \tilde{O}\left( sn\sqrt{\frac{\tilde{L}_{БГМ}\tilde{\Theta}}{\varepsilon}} \right) \text{ и } \tilde{T}_{ПБГМ} = \tilde{O}\left( sn\sqrt{\frac{\tilde{L}_{ПГМ}\tilde{\Theta}}{\varepsilon}} \right).$$

---

[19] Это типично для неускоренных покомпонентных методов (в том числе прямо-двойственных), т.е. в отличие от ускоренных методов, тут требуется намного более слабые предположения, чтобы обеспечить выполнение условия: стоимость итерации покомпонентного метода дешевле стоимости итерации соответствующего полноградиентного метода в число раз по порядку равному размерности пространства.

Написанное выше может навести на мысль, что ускоренные покомпонентные методы для двойственной задачи, как правило, не позволяют учитывать разреженность задачи. На самом деле это не так. За небольшую дополнительную плату (логарифмический множитель) можно специальным образом регуляризовать двойственную задачу (с помощью техники рестартов подобрать правильный параметр регуляризации, см. главу 3 [181], [37]), и использовать ПБГМ для регуляризованного функционала двойственной задачи, т.е. использовать ПБГМ в сильно выпуклом случае (см. текст сразу после замечания 5.1.2). При таком подходе достаточно просто решить (с желаемой точностью) двойственную задачу, а решение прямой задачи (в том же смысле, что и выше – с той же точностью) получается просто при подстановке найденного решения двойственной задачи в формулу $x(y)$. Тем не менее, здесь необходимо оговориться, что хотя описанный только что прием и “спасает положение”, все же получается это за упомянутую дополнительную плату. Хотя по теоретическим оценкам это плата, действительно, небольшая, численные эксперименты показывают, что реальные потери при использовании такой регуляризации вместо прямо-двойственности могут быть существенны. □

Результаты, изложенные в примере 5.1.3, допускают серьезные обобщения. В частности, можно переносить (частично), изложенное в примере 5.1.3, на сепарабельные задачи типа проектирования на аффинное многообразие

$$f(x) = \sum_{i=1}^{n} f_i(x_i) \to \min_{\substack{Ax=b,\, Cx\le d \\ x\in Q_x}}$$

и более общий класс сепарабельных задач лежандровского типа (включающий проекционный класс)

---

[20] Используем следующие обозначения: $\tilde{\Theta}$ – квадрат евклидового размера решения двойственной задачи,

$$\tilde{L}_{БГМ} = \max_{\|y\|_2\le 1,\, \|x\|_2\le 1} \left\langle A^T y, x \right\rangle^2 = \max_{\|x\|_2\le 1} \|Ax\|_2^2 = \lambda_{\max}\left(A^T A\right) \text{ и } \tilde{L}_{ПБГМ} \le \max_{\|y\|_1\le 1,\, \|x\|_2\le 1} \left\langle A^T y, x \right\rangle^2 = \max_{\|x\|_2\le 1} \|Ax\|_\infty^2 = \max_{k=1,..,m} \|A_k\|_2^2 .$$

$$\max_{x\in Q_x}\left\{\langle y, Ax\rangle - \sum_{i=1}^{n} f_i(x_i)\right\} + \tilde{g}(y) \to \min_{y\in Q_y}.$$

Здесь $Q_x$ – множество простой структуры (в смысле проектирования), $Q_y$ – множество, подходящее для эффективного использования (блочно-)покомпонентных методов (см. подраздел 5.1.4), $\tilde{g}(y)$ – "хорошая" функция для покомпонентных методов (см. подраздел 5.1.4). Можно также не предполагать явной формулы, связывающей $x(y)$ (тогда потребуется еще воспользоваться замечанием 5.1.5). В определенных ситуациях можно даже пытаться отказаться от сепарабельности $f(x)$ (к сожалению, вот тут пока мало что удалось получить). Все это порождает довольно много разных сочетаний (вариантов) и требует большого числа оговорок. Этому планируется посвятить отдельную работу. Далее мы ограничимся одним специальном классом задач, играющих важную роль в моделировании компьютерных и транспортных сетей (см., например, [33, 138]).

**Пример 5.1.4.** Рассмотрим задачу ($Q$ – множество просто структуры, скажем, неотрицательный ортант)

$$\sum_{k=1}^{m} f_k\left(A_k^T x\right) + g(x) \to \min_{x\in Q},$$

где $g(x) = \sum_{i=1}^{n} g_i(x_i)$ (впрочем, часть изложенных далее конструкций не требует выполнения этого условия). Градиенты функции $f_k$ вычислимы за $O(1)$ и имеют равномерно ограниченные константы (числом $L$) Липшица производной в 2-норме. Функция $g(x)$ предполагается сильно выпуклой в $p$-норме с константой $\mu_p$. Вводя матрицу $A = [A_1, ..., A_m]^T$ и вспомогательный вектор $z = Ax$ мы можем переписать эту задачу в "раздутом" пространстве $x := (x, z)$, как задачу типа проектирования на аффинное многообразие [33], рассмотренную ранее.[21] Однако для полноты картины[22] нам представляется

[21] В связи с этим, можно добавить к ограничениям, например, такого типа неравенства $A_k^T x \ge c_k$. Сложность задачи это не изменит. Этот факт можно использовать при численном поиске стохастических равновесий в модели стабильной динамики [33] (см. главу 1 диссертации).

[22] Приводимая далее конструкция позволяет (с помощью перехода к двойственной задаче и ее последующего изучения покомпонентными методами) в некотором смысле перейти от игры на разной гладкости по разным направлениям для исходной задачи к и игре на разной сильной выпуклости функционала исходной задачи по разным направлениям (при переходе к двойственной задаче эта игра

полезнее провести для этой задачи рассуждения немного в другом ключе (следуя, например, [33, 138]). Прежде всего, заметим, что в эту схему погружаются следующие задачи [138]:

1) $\frac{L}{2}\|Ax-b\|_2^2+\frac{\mu}{2}\|x-x_g\|_2^2 \to \min_{x\in\mathbb{R}^n}$,

2) $\frac{L}{2}\|Ax-b\|_2^2+\mu\sum_{k=1}^{n} x_k \ln x_k \to \min_{x\in S_n(1)}$.

Константа Липшица производных одинаковы $L^1 = L^2 = L$, константы сильной выпуклости (считаются в разных нормах) также одинаковы $\mu_2^1 = \mu_1^2 = \mu$. Опишем далее довольно общий способ построения двойственной задачи:[23]

$$\min_{x\in Q}\left\{\sum_{k=1}^{m} f_k\left(A_k^T x\right)+g(x)\right\}=\min_{\substack{x\in Q\\ z=Ax}}\left\{\sum_{k=1}^{m} f_k\left(z_k\right)+g(x)\right\}=$$

$$=\min_{\substack{x\in Q\\ z=Ax,z'}}\max_{y}\left\{\langle z-z',y\rangle+\sum_{k=1}^{m} f_k\left(z'_k\right)+g(x)\right\}=$$

$$=\max_{y\in\mathbb{R}^m}\left\{-\max_{\substack{x\in Q\\ z=Ax}}\left\{\langle -z,y\rangle-g(x)\right\}-\max_{z'}\left\{\langle z',y\rangle-\sum_{k=1}^{m} f_k\left(z'_k\right)\right\}\right\}=$$

$$=\max_{y\in\mathbb{R}^m}\left\{-\max_{x\in Q}\left(\langle -A^T y,x\rangle-g(x)\right)-\sum_{k=1}^{m}\max_{z'_k}\left(z'_k y_k - f_k\left(z'_k\right)\right)\right\}=$$

$$=\max_{y\in\mathbb{R}^m}\left\{-g^*\left(-A^T y\right)-\sum_{k=1}^{m} f_k^*\left(y_k\right)\right\}=-\min_{y\in\mathbb{R}^m}\left\{g^*\left(-A^T y\right)+\sum_{k=1}^{m} f_k^*\left(y_k\right)\right\}.$$

Для упомянутых задач, получим:

---

переходит в игру на гладкости двойственного функционала, которая уже неплохо проработана покомпонентными методами). К сожалению, все это возможно не в общем случае.

[23] Здесь важна оговорка о возможности "эффективно" решать задачу вида

$$\langle c,x\rangle+g(x)\to\max_{x\in Q}.$$

Вообще говоря, оговорка нетривиальная. В общем случае эта задача по сложности может соответствовать исходной. Стоит, однако, оговориться, что $g(x)$ в таких постановках в типичных приложениях является, как правило, "регуляризатором" исходной задачи (введенном нами с целью получения сильно выпуклой постановки или, например, возникшем при байесовском оценивании в качестве прайера (prior) или просто как пенализация (penalization) за сложность модели и т.д.). В любом случае, мы, как правило, имеем достаточно степеней свободы, чтобы добиться нужной простоты этой вспомогательной задачи. Сепарабельность $g(x)$ здесь довольно часто является "ключом к успеху".

1) $\frac{1}{2\mu}\left(\left\|x_g - A^T y\right\|_2^2 - \left\|x_g\right\|_2^2\right) + \frac{1}{2L}\left(\left\|y+b\right\|_2^2 - \left\|b\right\|_2^2\right) \to \min_{y\in\mathbb{R}^m}$,

2) $\frac{1}{\mu}\ln\left(\sum_{i=1}^{n}\exp\left(\frac{\left[-A^T y\right]_i}{\mu}\right)\right) + \frac{1}{2L}\left(\left\|y+b\right\|_2^2 - \left\|b\right\|_2^2\right) \to \min_{y\in\mathbb{R}^m}$.

В общем случае можно утверждать, что $\sum_{k=1}^{m} f_k^*\left(y_k\right)$ (композитный член в двойственной задаче) является сильно выпуклым в стандартной евклидовой норме (2-норме) с константой сильной выпуклости равной $L^{-1}$. Легко понять, что изучение свойств гладкости $g^*\left(-A^T y\right)$ (с точностью до множителя $\mu^{-1}$) совершенно аналогично тому, что мы уже делали в примере 5.1.3. То есть можно утверждать, что для двойственной задачи в стандартной евклидовой норме (2-норме)

$$\breve{L}_{БГМ} = \frac{1}{\mu}\max_{\|y\|_2\le 1,\,\|x\|_p\le 1}\left\langle A^T y, x\right\rangle^2 = \frac{1}{\mu}\max_{\|x\|_p\le 1}\left\|Ax\right\|_2^2 = \frac{1}{\mu}\begin{cases} 1)\ \lambda_{\max}\left(A^T A\right) \\ 2)\ \max\limits_{k=1,\dots,n}\left\|A^{\langle k\rangle}\right\|_2^2 \end{cases},$$

$$\breve{L}_{ПБГМ} \le \frac{1}{\mu}\max_{\|y\|_1\le 1,\,\|x\|_p\le 1}\left\langle A^T y, x\right\rangle^2 = \frac{1}{\mu}\max_{\|x\|_p\le 1}\left\|Ax\right\|_\infty^2 = \frac{1}{\mu}\begin{cases} 1)\ \max\limits_{k=1,..,m}\left\|A_k\right\|_2^2 \\ 2)\ \max\limits_{\substack{i=1,\dots,m \\ j=1,\dots,n}}\left|A_{ij}\right|^2 \end{cases}.$$

Для метода ACRCD* из замечания 5.1.2 (если смотреть через замечание 5.1.6, то следует полагать $\beta = 0$), примененного к двойственной задаче в не разреженном случае имеем следующие оценки времени работы:

1) $T_1 = \tilde{O}\left(nm\sqrt{\frac{L\max_{k=1,..,m}\left\|A_k\right\|_2^2}{\mu}}\right)$,

2) $T_2 = \tilde{O}\left(nm\sqrt{\frac{L\max_{i,j}\left|A_{ij}\right|^2}{\mu}}\right)$.

Если теперь посмотреть на исходную прямую задачу (с $L := L/m$)

$$\frac{1}{m}\sum_{k=1}^{m} f_k\left(A_k^T x\right) + g\left(x\right) \to \min_{x\in Q},$$

и оценить время работы ускоренного метода рандомизации суммы из подраздела 5.1.3 данного раздела 5.1, то получим анонсированное в подразделе 5.1.3 соответствие с приве-денными только что оценками (с $L := L/m$). Действительно, с учетом того, что константа Липшица градиентов $f_k\left(A_k^T x\right)$, посчитанные в соответствующих нормах (соответствующей норме, в которой сильно выпукл композит прямой задачи), равномерно оцениваются следующим образом:

1) $L \max_{k=1,..,m} \|A_k\|_2^2$,

2) $L \max_{i,j} |A_{ij}|^2$,

а сложность вычисления $\nabla f_k\left(A_k^T x\right)$ равна $\mathrm{O}(n)$, то согласно подразделу 5.1.3 имеем следующие оценки времени работы (далее считаем, что $m$ меньше $\min\{\ \}$):[24]

$$1)\quad \tilde{T}_1 = \tilde{\mathrm{O}}\left(n \cdot \left(m + \min\left\{\sqrt{m\frac{L \max_{k=1,..,m} \|A_k\|_2^2}{\mu}}, \frac{L \max_{k=1,..,m} \|A_k\|_2^2}{\mu}\right\}\right)\right) = \tilde{\mathrm{O}}\left(nm\sqrt{\frac{(L/m) \max_{k=1,..,m} \|A_k\|_2^2}{\mu}}\right),$$

$$2)\quad \tilde{T}_2 = \tilde{\mathrm{O}}\left(n \cdot \left(m + \min\left\{\sqrt{m\frac{L \max_{i,j} |A_{ij}|^2}{\mu}}, \frac{L \max_{i,j} |A_{ij}|^2}{\mu}\right\}\right)\right) = \tilde{\mathrm{O}}\left(nm\sqrt{\frac{(L/m) \max_{i,j} |A_{ij}|^2}{\mu}}\right).$$

Таким образом, имеет место полное соответствие (с точностью до опущенных в рассуждениях логарифмических множителей).[25] Интересно заметить, что для первой задачи здесь также как и в примере 5.1.3 можно сполна использовать разреженность матрицы $A$. Более того, эту задачу (также с полным учетом разреженности) можно решать и прямым ПБГМ (см. замечание 5.1.8). Соответствующие оценки имеют вид [138] (мы снова возвращаемся к исходному пониманию параметра $L$):

---

[24] Можно показать, развивая конструкцию п. 6.3 [163] с помощью идей, изложенных в замечании 5.1.8, что неускоренная составляющая этой оценки допускает для определенного класса задач в разреженном случае стоимости итерации $\tilde{\mathrm{O}}(s)$.

[25] Соответствие имеет место и для неускоренной составляющей выписанных оценок (собственно, именно этот случай рассматривался в подразделе 5.1.3 данного раздела 5.1). Чтобы это понять, нужно в оценки для неускоренных покомпонентных методов $T_{ПМ}$ (см. текст сразу после замечания 5.1.2 и пример 5.1.3 с $s = m$, $\tilde{s} = n$)

$$T_{ПМ}^1 = \tilde{\mathrm{O}}\left(nm\frac{L \max_{k=1,..,m} \|A_k\|_2^2}{\mu}\right),\quad T_{ПМ}^2 = \tilde{\mathrm{O}}\left(nm\frac{L \max_{i,j} |A_{ij}|^2}{\mu}\right)$$

подставить $L := L/m$.

$$T_1^{прям} = \tilde{O}\left( sn\sqrt{\frac{L\max\limits_{k=1,\dots,n}\left\|A^{\langle k\rangle}\right\|_2^2}{\mu}} \right),$$

$$T_1^{двойств} = \tilde{O}\left( \tilde{s}m\sqrt{\frac{L\max\limits_{k=1,\dots,m}\left\|A_k\right\|_2^2}{\mu}} \right) = \tilde{O}\left( sn\sqrt{\frac{L\max\limits_{k=1,\dots,m}\left\|A_k\right\|_2^2}{\mu}} \right).$$

В действительности, обе эти оценки оказываются завышенными.[26] Более аккуратные рассуждения, позволяют обобщить результаты примеров 5.1.3, 5.1.4 на случай неравноправия слагаемых. Все это приведет к замене максимума на некоторое (в зависимости от выбора $\beta$) среднее. Скажем, в упомянутых уже ранее транспортных приложениях [33, 37, 40, 42, 138] матрица $A$ не просто разреженная, но еще и битовая (состоит из нулей и единиц). В таком случае приведенные оценки переписываются следующим образом:

$$T_1^{прям} = \tilde{O}\left( sn\sqrt{\frac{Ls}{\mu}} \right),\ T_1^{двойств} = \tilde{O}\left( sn\sqrt{\frac{L\tilde{s}}{\mu}} \right).$$

Отсюда можно сделать довольно неожиданный вывод [138]: при $m \ll n$ стоит использовать прямой ПБГМ, а в случае $m \gg n$ – двойственный. Первый случай соответствует приложениям к изучению больших сетей (компьютерных, транспортных). Второй случай соответствует задачам, приходящим из анализа данных. □

### 5.1.6 Заключительные замечания

Если сравнить общие трудозатраты быстрого градиентного метода и его покомпонентного варианта, то довольно часто удается ускориться в $\sim\sqrt{n}$ раз (где $n$ – размерность пространства, в котором происходит оптимизация). Собственно, значительная часть данного раздела (подразделы 5.1.4, 5.1.5) была посвящена обсуждению того, в каких ситуациях можно рассчитывать на такое ускорение.

Отмеченное ускорение происходит за счет “обмана” потенциального сопротивляющегося оракула, корректирующего исходя из оставшихся у него свобод по ходу итерационного процесса оптимизируемую функцию таким образом, чтобы нам доставались наиболее плохие направления градиента (с большой константой Липшица – это соответствует пилообразному движению по дну растянутого оврага, с медленным приближением к середине оврага, в которой расположился минимум).

---

[26] Это легко усмотреть из способа рассуждений, в котором мы заменяем константы Липшица по разным направлениям, на худшую из них, это также позволяет эффективно использовать оценку ACRCD* с $\beta = 0$, см. замечания 5.1.2, 5.1.6.

Введение рандомизации в метод – это универсальный рецепт гарантированно обезопасить себя от худшего случая. Причем важно отметить, что это не теоретический трюк, который позволяет просто гарантировать лучшую теоретическую оценку. К сожалению, овражность функций – это довольно типичное свойство задач больших размеров, поэтому даже если мы возьмем сложность в среднем (по множеству типичных входов) для быстрого градиентного метода, то оценка получится все хуже (поскольку типичные входы не столь хороши), чем для рандомизированного метода. В обоих случаях мы считаем средние затраты (математическое ожидание), только в разных пространствах и по разным вероятностным мерам. В случае рандомизированного метода мы частично диверсифицируем себя от всего того, что может быть на входе, но это происходит не бесплатно, а с помощью препроцессинга, требующего расчет констант Липшица градиента по всем направлениям. К счастью, такой препроцессинг можно делать адаптивно и эффективно.[27]

В начале 2016 года появился электронный препринт [136] (см. также другие работы Z. Allen-Zhu этого периода http://arxiv.org/find/math/1/au:+Allen_Zhu_Z/0/1/0/all/0/1), в котором получены довольно близкие результаты к части результатов, приведенных в данном разделе (замечание 5.1.6 при $\beta = 1/2$ и пример 5.1.2). Однако для автора первоисточником идеи (для замечания 5.1.6 при $\beta = 1/2$ и примера 5.1.2) и отправной точкой в развитие этого направления послужило выступление Ю.Е. Нестерова в мае 2015 г (на конференции в Москве, приуроченной к 80-и летию Б.Т. Поляка). Заметим также, что в марте 2016 появился препринт [282], в основу которого положено это выступление.

В феврале 2016 года появился электронный препринт [192], в котором обсуждаются общие прямо-двойственные подходы для постановки задачи из примера 5.1.4.

Численные эксперименты с методами типа ACRCD, подтвердившие основные полученные в этом разделе теоретические оценки, проводила И.Н. Усманова (доклад на 58-й научной конференции МФТИ, 28 ноября 2015 г.). На практике оказалось, что метод ACRCD без рестартов сходится почти точно также как метод ACRCD с рестартами, а метод ACRCD* (на задачах больших размеров) сходится подобно ACRCD.

[27] Это можно делать за небольшую дополнительную плату – мультипликативный фактор порядка 4 в числе обращений за $\partial f(x)/\partial x_i$, если рассматривается композитная постановка или $Q$ не параллелепипедного типа, то требуется еще уметь рассчитывать (пересчитывать) вместе с $\partial f(x)/\partial x_i$ и значение функции.

## 5.2 Безградиентные прокс-методы с неточным оракулом для негладких задач выпуклой стохастической оптимизации на симплексе

### 5.2.1 Введение

Представим, следуя Ю.Е. Нестерову (см. выступление "Алгоритмические модели человеческого поведения" на математическом кружке, Москва, МФТИ и МЦНМО, 14 сентября 2012 г., видео доступно на сайте http://www.mathnet.ru/ ), что некоторый человек может характеризовать свое состояние вектором

$$x \in S_n(1) = \left\{ x \geq 0 : \sum_{i=1}^{n} x_i = 1 \right\}.$$

Насколько это состояние хорошее он может оценить, посчитав значение своей функции потерь $f(x)$ на этом векторе. К сожалению, рассматриваемый человек существенно ограничен в своих возможностях, поэтому посчитать (суб-)градиент этой функции он не может. Более того, значение функции он может посчитать лишь с неконтролируемым им шумом уровня $\delta$. Функция потерь (заданная в некоторой окрестности симплекса) предполагается выпуклой, но необязательно гладкой, с равномерно ограниченной нормой субградиента $\|\nabla f(x)\|_\infty \leq M$. Человек стремится оказаться в состоянии с наименьшими потерями $f_* = \min_{x \in S_n(1)} f(x)$, действуя итерационно по следующему простому правилу:

- выбрать случайно направление;
- сдвинуться с некоторым (небольшим) шагом из текущего состояния по этому направлению;
- посчитать значение функции в новом состоянии;
- исходя из значения функции в этих двух состояниях (состояния, в котором находились, и состояния, в котором оказались), определить новое состояние (из симплекса), в которое следует перейти.

В связи с описанным естественным способом действий возникает ряд вопросов. Например, как выбирать новое состояние, с целью минимизации числа обращений к «оракулу» (см. подраздел 5.2.3) за значением функции? Как именно "случайно" стоит выбирать направление? Как скажется зашумленность выдаваемых значений функции на это число обращений? Можно ли приблизиться к нижним оценкам требуемого числа обращений для достижения $f_*$ с точностью (по функции) $\varepsilon$?

В разделе мы рассматриваем еще более общую постановку, когда оракул может выдавать не значение функции, а лишь несмещенную (или не сильно смещенную, смещение контролируется уровнем шума $\delta$) оценку этого значения $f(x;\eta)$:

$$E_\eta\left[f(x;\eta)\right]=f(x).$$

В такой общности мы постараемся ответить на сформулированные вопросы (см. также раздел 6.2 главы 6). В частности, будет предложена процедура, требующая в случае гладкой функции $f(x)$ (точнее липшицевости градиента $f(x;\eta)$ по $x$)

$$O\left(\frac{M^2 n\ln n}{\varepsilon^2}\right)$$

обращений к оракулу за реализацией функции $f(x;\eta)$, что с точностью до логарифмического множителя соответствует нижней оценке [190]. Предложенный в разделе подход, позволяет также при некоторых дополнительных предположениях заметно улучшить приведенную оценку.

Отметим, что если бы вместо значения функции оракул выдавал стохастический градиент или хотя бы стохастическую производную функции по направлению, то ответы были бы, соответственно:

$$O\left(\frac{M^2\ln n}{\varepsilon^2}\right),\ O\left(\frac{M^2 n\ln n}{\varepsilon^2}\right),$$

что также с точностью до логарифмических множителей соответствует нижним оценкам (см., например, [128, 190]).

В целом проблематика раздела восходит к статье [231] (см. также [59, 91, 99, 233]). В подразделе 5.2.2 мы описываем известные результаты о сходимости метода зеркального спуска (МЗС) для задач стохастической оптимизации [91], которые нам понадобятся в дальнейшем. В подразделе 5.2.3 мы вводим неточный оракул, выдающий зашумленное значения реализации функции. Исходя из такой (частичной) информации в подразделе 5.2.3 предлагаются различные рандомизированные (безградиентные) обобщения МЗС. Рандомизация заключается в выборе случайного направления и вычислении (с помощью оракула) вместо стохастического субградиента стохастической дискретной производной функции по этому направлению [231]. Основные степени свободы, на которых можно играть: способ выбора случайного направления (в разделе обсуждаются равномерное распределение на евклидовом шаре и равномерные распределения на шарах в $l_1$ и $l_\infty$ нормах) и выбор шага дискретизации. В отсутствии шума выгоднее всего этот шаг стремить к ну-лю, т.е. просто вычислять стохастическую производную по направлению. Однако мы допускаем шум, и хотим понять, при каком максимально допустимом уровне шума оценки сохранят свой вид, скажем, в таких категориях: число итераций возрастет не более чем в два раза. В подразделах 5.2.2, 5.2.3 гладкость не предполагается. В подразделе 5.2.4 на примере изучения стохастических спусков по случайным направлением демонстрируется увеличение скорости сходимости, связанное с наличием гладкости. В подразделе 5.2.5 результаты подраздела 5.2.4 переносятся на стохастические безградиентные методы, т.е. по сути, на гладкий вариант постановки задачи из подраздела 5.2.3. Наличие гладкости дает ускорение в подразделах 5.2.4, 5.2.5 приблизительно в $n$ раз.

### 5.2.2 Метод зеркального спуска для задач стохастической оптимизации с неточным оракулом

Рассмотрим задачу стохастической оптимизации

$$f(x) = E_\eta\left[ f(x;\eta)\right] \to \min_{x\in S_n(1)}. \tag{5.2.1}$$

Здесь $\eta$ – случайная величина, $E_\eta\left[ f(x;\eta)\right]$ – математическое ожидание “взятое по $\eta$”, то есть при фиксированном $x$, при этом далее допускается, что в такой записи $x$ может быть случайным вектором. В таком случае математическое ожидание берётся только по $\eta$ (случайность в $x$ “фиксируется”). Если математическое ожидание берётся по $x$ (первое неравенство в теореме 5.2.1), то нижний индекс $\eta$ опускаем.

Обозначим

$$f_* = \min_{x\in S_n(1)} f(x) = \min_{x\in S_n(1)} E_\eta\left[ f(x;\eta)\right].$$

**Замечание 5.2.1.** Везде далее мы будем использовать обозначения обычного градиента для субградиента. Запись $\nabla_x f(x;\eta)$ в вычислительном контексте (например, в итерационной процедуре (5.2.2) ниже) означает какой-либо измеримый селектор стохастического субдифференциала [308], а если в контексте проверки условий (например, в условии 5.2.2 или условии 5.2.3 ниже), то $\nabla_x f(x;\eta)$ пробегает все элементы стохастического субдифференциала.

Для формулировки основной теоремы этого пункта нам понадобятся следующие **условия**:

5.2.1 $f(x;\eta)$ – выпуклая функция по $x$ (в действительности, с некоторыми оговорками [308], достаточно только выпуклости $f(x)$);

5.2.2 Стохастический субградиент $\nabla_x f(x;\eta)$ [309] удовлетворяет условию (тождественно по $x$):

$$E_\eta\left[\nabla_x f(x;\eta)\right] \equiv \nabla_x E_\eta\left[f(x;\eta)\right];$$

5.2.3 $\left\|\nabla_x f(x;\eta)\right\|_\infty \le M$ – равномерно, с вероятностью 1.

Для справедливости части утверждений достаточно требовать вместо условия 5.2.3 одно из следующих (более слабых) условий:

$$\text{а) } E_\eta\left[\left\|\nabla_x f(x;\eta)\right\|_\infty^2\right] \le M^2; \quad \text{б) } E_\eta\left[\exp\left(\frac{\left\|\nabla_x f(x;\eta)\right\|_\infty^2}{M^2}\right)\right] \le \exp(1).$$

Для решения задачи (5.2.1) воспользуемся методом зеркального спуска (точнее двойственных усреднений) в форме [125, 269]. Положим $x_i^1 = 1/n$, $i = 1,\ldots,n$. Пусть $t = 1,\ldots,N-1$.

$$x_i^{t+1} = \frac{\exp\left(-\frac{1}{\beta_{t+1}}\sum_{k=1}^{t}\frac{\partial f\left(x^k;\eta^k\right)}{\partial x_i}\right)}{\sum_{l=1}^{n}\exp\left(-\frac{1}{\beta_{t+1}}\sum_{k=1}^{t}\frac{\partial f\left(x^k;\eta^k\right)}{\partial x_l}\right)}, \quad i = 1,\ldots,n, \quad \beta_t = \frac{M\sqrt{t}}{\sqrt{\ln n}}. \tag{5.2.2}$$

Здесь $\left\{\eta^k\right\}$ – независимые одинаково распределенные (также как $\eta$) случайные величины.

Приводимая ниже теорема фактически установлена в работах [125, 225, 269]. Однако здесь мы непосредственно воспользовались формулировкой из работы [54] (см. также раздел 6.1 главы 6).

**Теорема 1.** *Пусть справедливы условия 5.2.1, 5.2.2, 5.2.3.а, тогда*

$$E\left[f\left(\frac{1}{N}\sum_{k=1}^{N}x^k\right)\right] - f_* \le \frac{1}{N}\sum_{k=1}^{N}E\left[f\left(x^k\right)\right] - f_* \le 2M\sqrt{\frac{\ln n}{N}}.$$

*Пусть справедливы условия 5.2.1, 5.2.2, 5.2.3, тогда при* $\Omega \ge 0$

$$P_{x^1,\ldots,x^N}\left\{\frac{1}{N}\sum_{k=1}^{N}f\left(x^k\right) - f_* \ge \frac{2M}{\sqrt{N}}\left(\sqrt{\ln n} + \sqrt{8\Omega}\right)\right\} \le$$

$$\le P_{x^1,\ldots,x^N}\left\{f\left(\frac{1}{N}\sum_{k=1}^{N}x^k\right) - f_* \ge \frac{2M}{\sqrt{N}}\left(\sqrt{\ln n} + \sqrt{8\Omega}\right)\right\} \le \exp(-\Omega).$$

**Замечание 5.2.2.** Если вместо условия 5.2.3 имеет место более слабое условие 5.2.3.б, то последняя формула останется верной, при небольшой корректировке (см. раздел 6.1 главы 6):

$$\frac{2M}{\sqrt{N}}\left(\sqrt{\ln n}+\sqrt{8\Omega}\right)\to C\frac{M}{\sqrt{N}}\left(\sqrt{\ln n}+\Omega\right),$$

где константа $C\sim 10$. Приведенный результат можно обобщить и на более тяжелые хвосты [44].

### 5.2.3 Безградиентная модификация метода зеркального спуска для задач стохастической оптимизации с неточным оракулом

Введем понятие оракула, выдающего зашумленное значение функции $f(x)$, определенной[28] в $\mu_0$-окрестности $S_n(1)$.[29] При этом везде в дальнейшем под $f(x)$ и $f(x,\eta)$ мы будем понимать, соответственно, не зашумленные значение функции и ее несмещенной (тождественно по $x$) реализации. Наличие шума мы будем явно указывать, вводя его аддитивным образом.

**Предположение 5.2.1.** *Оракул выдает (на запрос, в котором указывается только одна точка $x$) $f(x,\eta)+\tilde{\delta}(\eta)$, где с.в. $\eta$ независимо разыгрывается из одного и того же распределения, фигурирующего в постановке (1); случайная величина $\tilde{\delta}(\eta)$ (случайность может быть обусловлена не только зависимостью от $\eta$) не зависит от $x$ и ограничена по модулю известным нам числом $\delta$ – допустимым уровнем шума.*

Приведем одну из возможных мотивировок такого оракула (см. также раздел 2.1 главы 2). Предположим, что оракул может считать абсолютно точно значение (или реализацию) функции, но вынужден нам выдавать лишь конечное (предписанное) число первых бит (конечная мантисса). Таким образом, в последнем полученном бите есть некоторая неточность (причем мы не знаем по какому правилу оракул формирует этот последний выдаваемый значащий бит). Однако мы всегда можем прибавить (по mod 1) к этому биту случайно приготовленный (независимый) бит. В результате, не ограничивая общности, можно считать, что оракул последний бит выбирает просто случайно в независимости от отброшенного остатка.

---

[28] Везде далее в разделе мы будем предполагать, что $f(x)$ не просто определена в достаточно большой $\mu_0$-окрестности исходного множества, но и сохраняет все свои свойства в этой окрестности, в частности, выпуклость и константы Липшица.

[29] Все, что будет написано далее, можно перенести (без изменений итоговых формул с точностью до константного множителя) на случай более общего оракула, описанного в п. 4 работы [44] (см. также раздел 2.1 главы 2 и раздел 6.2 главы 6). К сожалению, в [44], все равно, относительно оракула делаются обременительные предположения. Впрочем, в этой же работе схематично показано, как можно распространить (с ужесточением условий на допустимый уровень шума) все, что далее будет написано на случай самого общего оракула, выдающего зашумленное значение функции (реализации функции). Об этом также написано в работе [155] и разделе 6.2 главы 6.

**Предположение 5.2.2.** *В случае задач стохастической оптимизации принципиально важно, что разрешается на каждом шаге (итерации) обратиться к оракулу за значениями функции на одной реализации ($\eta$ одно и то же), но в двух разных точках. В не стохастическом случае достаточно иметь возможность одного обращения на каждом шаге.*

Число итераций (с точностью до множителя 2 в стохастическом случае) – это число обращений к такому оракулу. Наша цель, обращаясь к оракулу на одном шаге (итерации) не более двух раз, так организовать итерационную процедуру, чтобы сгенерированная на основе опроса оракула последовательность $\{x^k\}$ с вероятностью $\geq 1-\sigma$ удовлетворяла неравенству

$$f\left(\frac{1}{N}\sum_{k=1}^{N}x^k\right)-f_*\leq\frac{1}{N}\sum_{k=1}^{N}f\left(x^k\right)-f_*\leq\varepsilon$$

с как можно меньшим значением $N$.

**Замечание 5.2.3.** На самом деле, не очень важно, сколько раз разрешено обращаться к оракулу, важно только, что не менее двух раз [190]. Приведенные в подразделах 5.2.4, 5.2.5 результаты легко переписываются, если вместо двух точек (на одной реализации) разрешается использовать $k\leq n+1$ точек (на одной реализаций): грубо говоря, оценки числа итераций от желаемой точности улучшатся в $k$ раз $N(\varepsilon)\to N(\varepsilon)/k$ [41, 190]. Если же разрешается обращаться только один раз, то картина принципиально меняется [190]. В этом случае на данный момент имеется достаточно большой зазор (для детерминированных постановок задач) между нижними оценками и тем, что сейчас дают лучшие методы [91, 148, 190].

Изложим далее общую схему, позволяющую свести описанную выше постановку к постановке подразделе 5.2.2. Тогда можно будет воспользоваться теоремой 5.2.1.

Пусть $e\in RS_p^n(1)$ ($\tilde{e}\in RB_p^n(1)$) – случайный вектор, равномерно распределенный на сфере (шаре) единичного радиуса в $l_p$ норме в $\mathbb{R}^n$ (далее мы ограничимся рассмотрением случаев: $p=1$, $p=2$, $p=\infty$). Сгладим (следуя [91]) исходную функцию с помощью локального усреднения по шару радиуса $\mu>0$ ($\mu\leq\mu_0$), который будет выбран позже,

$$f^{\mu}(x;\eta)=E_{\tilde{e}}\left[f(x+\mu\tilde{e};\eta)\right],$$

$$f^{\mu}(x)=E_{\tilde{e},\eta}\left[f(x+\mu\tilde{e};\eta)\right].$$

Заменим исходную задачу (5.2.1) следующей задачей

$$f^{\mu}(x)\to\min_{x\in S_n(1)}. \qquad (5.2.3)$$

Легко проверить (см., например, [91, 166, 270] для $p'=2$, в общем случае рассуждения в точности такие же), что если выполняется **условие** (это условие обобщает условие 5.2.3 подраздела 5.2.2, в частности, $M_1$ в условии 5.2.4 соответствует $M$ в условии 5.2.3.а подраздела 5.2.2)

5.2.4 $\left|f(x;\eta)-f(y;\eta)\right|\le M_{p'}(\eta)\|x-y\|_{p'}$, $M_{p'}=\sqrt{E_{\eta}\left[M_{p'}(\eta)^2\right]}<\infty$,

то

$$0\le f^{\mu}(x;\eta)-f(x;\eta)\le M_{p'}(\eta)\mu,$$

$$0\le f^{\mu}(x)-f(x)\le M_{p'}\mu.$$

Если выполняется **условие**

5.2.5 $\left\|\nabla_x f(x;\eta)-\nabla_x f(y;\eta)\right\|_{q'}\le L_{p'}(\eta)\|x-y\|_{p'}$, $L_{p'}=\sqrt{E_{\eta}\left[L_{p'}(\eta)^2\right]}<\infty$,

то

$$0\le f^{\mu}(x;\eta)-f(x;\eta)\le L_{p'}(\eta)\mu^2/2,$$

$$0\le f^{\mu}(x)-f(x)\le L_{p'}\mu^2/2,$$

где $1/p'+1/q'=1$. Предположим, что (в этом пункте мы не предполагаем гладкости, поэтому можно просто положить $L_{p'}=\infty$)

$$\min\left\{M_{p'}\mu,L_{p'}\mu^2/2\right\}\le\varepsilon/2, \qquad (5.2.4)$$

и с вероятностью $\ge 1-\sigma$ удалось получить следующее неравенство (например, воспользовавшись каким-то образом для задачи (5.2.3) теоремой 5.2.1):

$$\frac{1}{N}\sum_{k=1}^{N}f^{\mu}\left(x^k\right)-\min_{x\in S_n(1)}f^{\mu}(x)\le\frac{\varepsilon}{2}.$$

Тогда с вероятностью $\ge 1-\sigma$:

$$f\left(\frac{1}{N}\sum_{k=1}^{N}x^k\right)-f_*\le\frac{1}{N}\sum_{k=1}^{N}f\left(x^k\right)-\min_{x\in S_n(1)}f(x)\le\frac{1}{N}\sum_{k=1}^{N}f^{\mu}\left(x^k\right)-\min_{x\in S_n(1)}f^{\mu}(x)+\frac{\varepsilon}{2}\le\varepsilon.$$

Таким образом, при условии (5.2.4), решение задачи (5.2.3) с точностью $\varepsilon/2$ является решением задачи (5.2.1) с точностью $\varepsilon$.

Сглаживание было введено для того, чтобы для сглаженной задачи с помощью описанного оракула можно было получить несмещенную оценку субградиента. К сожалению, без сглаживания не понятно, как это можно было бы сделать. Итак, введем (при $p=2$, см., например, [314]) аналог стохастического субградиента

$$g^{\mu}(x;e,\eta)=\frac{\mathrm{Vol}\left(S_p^n(\mu)\right)}{\mathrm{Vol}\left(B_p^n(\mu)\right)}\left(f(x+\mu e;\eta)-f(x;\eta)\right)\breve{e},$$

где $e$ – случайный вектор, равномерно распределенный на сфере радиуса 1 в $l_p$ норме (обозначим такую сферу $S_p^n(1)$); $\mathrm{Vol}\left(B_p(\mu)\right)$ – объем шара радиуса $\mu$ в $l_p$ норме, аналогично определяется $\mathrm{Vol}\left(S_p(\mu)\right)$; $\breve{e}=\breve{e}(e)$ – вектор, с $l_2$ нормой равной 1, ортогональный поверхности $S_p^n(1)$ в точке $e$. Например,

Таблица 5.2.1

| $p$ | Аналог стохастического субградиента | Выбор направления |
|---|---|---|
| 1 | $\frac{n}{\mu}\left(f(x+\mu e;\eta)-f(x;\eta)\right)\begin{pmatrix}\mathrm{sign}\, e_1\\ \ldots\ldots\ldots\\ \mathrm{sign}\, e_n\end{pmatrix}$ | $e\in RS_1^n(1)$ |
| 2 | $\frac{n}{\mu}\left(f(x+\mu e;\eta)-f(x;\eta)\right)e$ | $e\in RS_2^n(1)$ |
| $\infty$ | $\frac{n}{\mu}\left(f(x+\mu e;\eta)-f(x;\eta)\right)\breve{e}_{i(e)}$ (п.н.)<br>$\breve{e}_{i(e)}=(\underbrace{0,\ldots,0,1}_{i(e)},0,\ldots,0)$, $i(e)=\arg\max\limits_{i=1,\ldots,n}\lvert e_i\rvert$ | $e\in RS_\infty^n(1)$ |

Основное свойство $g^{\mu}(x;e,\eta)$ заключается в том, что (воспользовались векторным вариантом теоремы Стокса, подобно [166])

$$E_{e,\eta}\left[g^{\mu}(x;e,\eta)\right]\equiv\nabla f^{\mu}(x).$$

Причем, это свойство сохраняется и в случае, когда вместо “идеального” значения реализаций $f(x+\mu e;\eta)$ и $f(x;\eta)$ оракул выдает зашумленные

$$E_{e,\eta}\left[g_{\delta}^{\mu}(x;e,\eta)\right]\equiv\nabla f^{\mu}(x).$$

Чтобы можно было воспользоваться теоремой 5.2.1 для сглаженной задачи (5.2.3) необходимо оценить $\left\|g_\delta^\mu(x;e,\eta)\right\|_\infty$, где

$$g_\delta^\mu(x;e,\eta)=\frac{\operatorname{Vol}\left(S_p^n(\mu)\right)}{\operatorname{Vol}\left(B_p^n(\mu)\right)}\left(f(x+\mu e;\eta)+\tilde{\delta}_{x+\mu e}(\eta)-\left(f(x;\eta)+\tilde{\delta}_x(\eta)\right)\right)\bar{e}.$$

Из определения оракула следует, что при $p=1$ и условии 5.2.3 подраздела 5.2.2

$$\left\|g_\delta^\mu(x;e,\eta)\right\|_\infty\le\left(M+\frac{2\delta}{\mu}\right)n. \tag{5.2.5}$$

При $p=2$ и $p=\infty$ оценка (5.2.5) получается хуже (см. Таблицу 5.2.2 ниже).

Выберем согласно условию (5.2.4) $\mu=\varepsilon/(2M)$ и будем считать, что (условие на допустимый уровень шума)

$$\delta\le\varepsilon/4.$$

Тогда условие (5.2.5) перепишется следующим образом

$$\left\|g_\delta^\mu(x;e,\eta)\right\|_\infty\le 2Mn.$$

Подобно алгоритму (5.2.2) опишем алгоритм решения задачи (5.2.3) для $p=1$. Положим $x_i^1=1/n$, $i=1,...,n$. Пусть $t=1,...,N-1$.

$$x_i^{t+1}=\frac{\exp\left(-\frac{1}{\beta_{t+1}}\sum_{k=1}^{t}\left[g_\delta^\mu\left(x^k;e^k,\eta^k\right)\right]_i\right)}{\sum_{l=1}^{n}\exp\left(-\frac{1}{\beta_{t+1}}\sum_{k=1}^{t}\left[g_\delta^\mu\left(x^k;e^k,\eta^k\right)\right]_l\right)},\quad i=1,...,n,\quad \beta_t=\frac{2Mn\sqrt{t}}{\sqrt{\ln n}},$$

где $[z]_i$ – $i$-я координата вектора $z$.

**Теорема 5.2.2.** *Пусть мы располагаем оракулом из предположений 5.2.1, 5.2.2 с* $\delta\le\varepsilon/4$ *и справедливы условия 5.2.1, 5.2.2, 5.2.3.а подразлдела 5.2.2,* $p=1$, *тогда для*

$$N=\left\lceil\frac{64M^2n^2\ln n}{\varepsilon^2}\right\rceil$$

*имеет место оценка*

$$E\left[f\left(\frac{1}{N}\sum_{k=1}^{N}x^k\right)\right]-f_*\le\varepsilon.$$

*Если (дополнительно) справедливо условие 5.2.3 подраздела 5.2.2, тогда для*

$$N = \left\lceil \frac{128M^2n^2}{\varepsilon^2}\left(\ln n + 8\ln\left(\sigma^{-1}\right)\right)\right\rceil$$

*с вероятностью* $\geq 1-\sigma$ *имеет место оценка*

$$f\left(\frac{1}{N}\sum_{k=1}^{N} x^k\right) - f_* \leq \varepsilon.$$

**Доказательство.** Применим теорему 5.2.1 к функции $f^{\mu}(x)$ с

$$N = \left\lceil \frac{4(2Mn)^2}{(\varepsilon/2)^2}\ln n\right\rceil$$

для оценки скорости сходимости по математическому ожиданию, и с

$$N = \left\lceil \frac{8(2Mn)^2}{(\varepsilon/2)^2}\left(\ln n + 8\ln\left(\sigma^{-1}\right)\right)\right\rceil = \left\lceil \frac{128M^2n^2}{\varepsilon^2}\left(\ln n + 8\ln\left(\sigma^{-1}\right)\right)\right\rceil$$

для оценки скорости сходимости с учетом вероятностей больших уклонений. В последнем случае мы еще воспользовались неравенством $\left(\sqrt{a}+\sqrt{b}\right)^2 \leq 2a+2b$ . □

Резюмируем полученные результаты в виде Таблицы 5.2.2. При этом считаем выполненными условия 5.2.1, 5.2.2 подраздела 5.2.2 и условие 5.2.4 (константы $M_1$, $M_2$, $M_\infty$ определяются в условии 5.2.4). Во второй строчке Таблицы 5.2.2 приведены математические ожидания числа итераций. Заметим при этом, что

$$M_1^2 \leq M_2^2 \leq nM_1^2,\; M_2^2 \leq M_\infty^2 \leq nM_2^2.$$

Таблица 5.2.2

| $p=1$ | $p=2$ | $p=\infty$ |
|---|---|---|
| $O\left(\frac{M_1^2n^2\ln n}{\varepsilon^2}\right)$ | $O\left(\frac{M_2^2n^2\ln n}{\varepsilon^2}\right)$ | $O\left(\frac{M_\infty^2n^2\ln n}{\varepsilon^2}\right)$ |

Нам не известно, оптимальна ли выписанная в Таблице 5.2.2 оценка для $p=1$ при наложенных условиях на уровень шума $\delta \leq \varepsilon/4$. Однако, имеется гипотеза, что полученная оценка оптимальна с точностью до мультипликативной константы при заданном уровне шума. В условиях отсутствия шума ($\delta = 0$) приведенная оценка (и, тем более, остальные: $p=2$, $p=\infty$), вообще говоря, не является оптимальной для негладких задач стохастической оптимизации. В работе [190] получена оценка (с помощью техники двойного сглаживания), которая позволяет сократить число итераций в оценке из второго столбца Таблицы 5.2.2 ( $p = 2$ ) в $\sim n/\ln n$ раз (причем, по-видимому, логарифмический множитель тут можно убрать). Однако предложенный в [190] метод не практичный (в отличие от предложенного в данном пункте метода), поскольку чрезвычайно чувствителен даже к очень небольшим шумам.

### 5.2.4 Модификация метода зеркального спуска для гладких задач стохастической оптимизации при спусках по случайному направлению

К сожалению, описанный в подразделе 5.2.3 подход дает оценку в $n$ раз большую нижней оценки в гладком случае [190]. Поскольку нам интересны ситуации, в которых $n \gg 1$, то необходимо этот зазор как-то устранить. Естественно попытаться найти в рассуждениях подраздела 5.2.3 наиболее грубое место и попробовать провести более точные рассуждения. К счастью, такое место всего одно – неравенство (5.2.5).

Считаем далее выполненными условия 5.2.1, 5.2.2 подраздела 5.2.2 и условие 5.2.4 подраздела 5.2.3 (с $p' = 2$ ).

Чтобы пояснить, в чем заключается грубость, рассмотрим для большей наглядности случай с $\delta = 0$. Тогда можно устремить $\mu \to 0+$ и получить

$$g^{\mu}(x;e,\eta) \to g(x;e,\eta) = \frac{\mathrm{Vol}\left(S_p^n(1)\right)}{\mathrm{Vol}\left(B_p^n(1)\right)} \langle \nabla_x f(x;\eta), e \rangle \bar{e} .$$

Аналогично подразделу 5.2.3, имеем

$$E_{e,\eta}\left[g(x;e,\eta)\right] \equiv \nabla f(x) .$$

Оценим сначала $E_{e,\eta}\left[\left\|g(x;e,\eta)\right\|_{\bar{q}}^2\right]$ ( $2 \le \bar{q} \le \infty$ выбирается исходя из структуры множества, на котором происходит оптимизация,[30] см. [128]) при $p = 1$ (этот параметр отве-

---

[30] В разбираемом в разделе 5.2 случае, когда ограничение в виде симплекса, выбирают $\bar{q} = \infty$, см. подраздел 5.2.2 (обоснование такому выбору имеется, например, в [259]). Как уже отмечалось, выбор $\bar{q}$ осуществляется исходя из структуры множества, на котором происходит оптимизация. Вместо используемого нами варианта метода зеркального спуска (МЗС) из подраздела 5.2.2, “настроенного” на то, что оптимизация происходит на симплексе, можно использовать вариант, подходящий для любого другого выпуклого множества $Q$, в котором в прямом пространстве выбрана норма $l_{\bar{p}}$ ( $1/\bar{p} + 1/\bar{q} = 1$ ) и определена неотрицательная сильно выпуклая (с константой $\ge 1$ ) относительно этой нормы функция $d(x)$, задающая “расстояние” Брэгмана

$$V(x,y) = d(x) - d(y) - \langle \nabla d(y), x - y \rangle .$$

чает за выбор способа рандомизации в Таблице 5.2.1 подраздела 5.2.3) в категориях $\mathrm{O}(\ )$. Для этого заметим,[31] что случайный вектор $e \in RS_1^n(1)$ можно представить как $e = a/\|a\|_1$, где компоненты вектора – независимые лапласовские случайные величины, т.е. с плотностью $e^{-|y|}/2$. Согласно Таблице 5.2.1, имеем

$$E_{e,\eta}\left[\left\|g(x;e,\eta)\right\|_{\bar q}^2\right] = n^{2+2/\bar q} E_{e,\eta}\left[\left\langle \nabla_x f(x;\eta), e\right\rangle^2\right] = n^{2+2/\bar q} E_{e,\eta}\left[\frac{\left\langle \nabla_x f(x;\eta), a\right\rangle^2}{\|a\|_1^2}\right].$$

Далее воспользуемся тем, что $n \gg 1$. Тогда исходя из явления концентрации меры[32] [159, 244], имеем: $\|a\|_1^2$ – сконцентрирован (с хвостами вида $e^{-\sqrt{y}}$) около своего математического ожидания $c_1 n$, $\left\langle \nabla_x f(x;\eta), a\right\rangle^2$ – $e^{-\sqrt{y}}$-сконцентрирован (при зафиксированном $\eta$) около своего математического ожидания $c_2\left\|\nabla_x f(x;\eta)\right\|_2^2$.[33] В результате получается следующая оценка (отметим, что здесь и далее $M_2$ определяется условием 5.2.4 подраздела 5.2.3)

$$E_{e,\eta}\left[\left\|g(x;e,\eta)\right\|_{\bar q}^2\right] = \mathrm{O}\left(n^{1+2/\bar q} M_2^2\right) \text{ (при } p = 1\text{)}. \qquad (5.2.6)$$

В действительности, можно показать, что и $\left\|g(x;e,\eta)\right\|_{\bar q}^2$ имеет $e^{-\sqrt{y}}$-концентрацию около своего математического ожидания, если в условии 5.2.4 подраздела 5.2.3 $M_2(\eta) \equiv M_2$.

---

Итоговая оценка ожидаемого числа итераций для соответствующего МЗС (в случае, когда на каждой итерации доступен несмещенный стохастический субградиент, математическое ожидание квадрата $l_{\bar q}$ нормы которого равномерно по $x$ ограничено числом $M_{\bar p}^2$) будет [128, 259]: $\mathrm{O}\left(M_{\bar p}^2 V\left(x^*, x^1\right)/\varepsilon^2\right)$. Причем если $Q = B_{\bar p}(R)$, то оптимально (с точностью до мультипликативного множителя) выбирать:

$$d(x) = \frac{1}{2(a-1)}\|x\|_a^2,\ 1 \le \bar p \le a;\ d(x) = \frac{1}{2(\bar p-1)}\|x\|_{\bar p}^2,\ a \le \bar p \le 2;\ d(x) = \frac{1}{2}\|x\|_2^2,\ 2 \le \bar p \le \infty.$$

В подразделе 5.2.2 использовался вариант МЗС с прокс-функцией $d(x) = \ln n + \sum_{k=1}^n x_k \ln x_k$, что также приводит к неулучшаемой (с точностью до числового множителя) оценке числа итераций [259].

[31] Приводимая далее в этом абзаце схема рассуждений была нам сообщена Александром Содиным.

[32] Впрочем, можно приведенные ниже результаты получить и без тонких оценок плотности концентрации, исходя из классических вариантов закона больших чисел, центральной предельной теоремы, и их идемпотентных аналогов [84].

[33] Точные значения положительных констант $c_1$ и $c_2$ (аналогично $c_3$, $c_4$, $c_5$) нас не интересуют, для нас сейчас важно только то, что они не зависят от $n$. Здесь и далее для большей наглядности мы предполагаем выполненным максимально сильное условие 5.2.3 подраздела 5.2.2.

Еще более геометрически наглядные рассуждения, восходящие к Пуанкаре–Леви [244], связанные с концентрацией равномерной меры на евклидовой сфере, позволяют получить следующую оценку

$$E_{e,\eta}\left[\left\|g(x;e,\eta)\right\|_{\bar{q}}^{2}\right]=\mathrm{O}\left(n^{2/\bar{q}}\ln n\,M_2^2\right)\ \text{(при } p=2\text{)}. \qquad (5.2.7)$$

Отличие в рассуждениях в том, что $e\in RS_2^n(1)$ стоит представлять как $e=a/\|a\|_2$, где $a\in \mathrm{N}(0,I_n)$, где $I_n$ – единичная матрица (на диагонали 1, остальные элементы нули) размера $n\times n$. Тогда

$$E_{e,\eta}\left[\left\|g(x;e,\eta)\right\|_{\bar{q}}^{2}\right]=n^2E_{e,\eta}\left[\left\langle\nabla_x f(x;\eta),e\right\rangle^2\|e\|_{\bar{q}}^{2}\right]=n^2E_{e,\eta}\left[\frac{\left\langle\nabla_x f(x;\eta),a\right\rangle^2\|a\|_{\bar{q}}^{2}}{\|a\|_2^4}\right],$$

где $\|a\|_2^4$ – $e^{-\sqrt{y}}$-сконцентрирован около своего математического ожидания $c_3n^2$, $\|a\|_{\bar{q}}^2$ – $e^{-y}$-сконцентрирован (экспоненциально сконцентрирован) около своего математического ожидания, которое оценивается сверху[34] $c_4n^{2/\bar{q}}\ln n$ [84], $\left\langle\nabla_x f(x;\eta),a\right\rangle^2$ – экспоненциально сконцентрирован (при зафиксированном $\eta$) около своего математического ожидания $c_5\left\|\nabla_x f(x;\eta)\right\|_2^2$, которое оценивается сверху $c_5M_2^2$, если в условии 5.2.4 подраздела 5.2.3 $M_2(\eta)\equiv M_2$.

Наиболее же просто исследуется случай $p=\infty$. Основным здесь является следующее наблюдение: практически весь объем многомерного куба сосредоточен на его границе [244].[35] Таким образом, в предположении $n\gg 1$ с хорошей точностью мы можем заменить условие $e\in RS_\infty^n(1)$ условием $e\in RB_\infty^n(1)$. Последнее распределение тривиально исследуется [244]. Аналогично вышеизложенному

$$E_{e,\eta}\left[\left\|g(x;e,\eta)\right\|_{\bar{q}}^{2}\right]=n^2E_{e,\eta}\left[\left\langle\nabla_x f(x;\eta),e\right\rangle^2\|e\|_{\bar{q}}^{2}\right]=n^2E_{e,\eta}\left[\left\langle\nabla_x f(x;\eta),e\right\rangle^2\right].$$

Таким образом,

---

[34] Эту оценку можно уточнить. В частности (бакалаврский диплом И.Н. Усмановой [117]),

$$E_e\left[\|e\|_{\bar{q}}^2\right]\le(\bar{q}-1)n^{2/\bar{q}-1},\ E_e\left[\|e\|_{\infty}^2\right]\le(4\ln n)/n,\ e\in RS_2^n(1).$$

Мы используем это далее, см. Таблицу 5.2.4 и выкладки в подразделе 5.2.5.

[35] Действительно, объем $n$-мерного куба со стороной 1 равен 1, а со стороной $1-\delta$ равен $(1-\delta)^n\ll 1$ – при достаточно больших $n$.

$$E_{e,\eta}\left[\left\|g(x;e,\eta)\right\|_{\bar{q}}^2\right]=\mathrm{O}\left(n^2M_2^2\right).\ (\text{при } p=\infty) \qquad (5.2.8)$$

В действительности, можно показать, что и $\left\|g(x;e,\eta)\right\|_{\bar{q}}^2$ имеет экспоненциальную концентрацию (при зафиксированном $\eta$) около своего математического ожидания, если в условии 5.2.4 подраздела 5.2.3 $M_2(\eta)\equiv M_2$.

Исходя из несмещенной оценки субградиента $g(x;e,\eta)$, можно построить алгоритм, аналогичный (5.2.2): в (5.2.2) заменяем $\nabla_x f(x;\eta)$ на $g(x;e,\eta)$. Подставляя в оценки (5.2.6) – (5.2.8) $\bar{q}=\infty$ (что соответствует рассматриваемой в данном разделе оптимизации на симплексе [259]), получим итоговые оценки среднего числа итераций такого алгоритма:

Таблица 5.2.3

| $p=1$ | $p=2$ | $p=\infty$ |
|---|---|---|
| $\mathrm{O}\left(\frac{nM_2^2\ln n}{\varepsilon^2}\right)$ | $\mathrm{O}\left(\frac{M_2^2\ln^2 n}{\varepsilon^2}\right)$ | $\mathrm{O}\left(\frac{n^2M_2^2\ln n}{\varepsilon^2}\right)$ |

Из Таблицы 5.2.3 хорошо видно, какая рандомизация предпочтительнее – на евклидовой сфере ($p=2$). Отсюда, с учетом того, что $M_1^2\le M_2^2\le nM_1^2$, получаем оценку

$$\mathrm{O}\left(\frac{M_2^2\ln^2 n}{\varepsilon^2}\right)\le\mathrm{O}\left(\frac{M_1^2 n\ln^2 n}{\varepsilon^2}\right),$$

которая с точностью до логарифмического множителя соответствует нижней оценке [190]. Однако если предположить, что $M_2^2\ll nM_1^2$, то получается, что можно превзойти нижнюю оценку, т.е. быстрее достичь желаемой точности, чем предписано нижней оракульной оценкой [91, 190]. Но никакого противоречия здесь, конечно, нет, поскольку нижняя оценка была получена без всяких дополнительных предположений. Делая такое предположение ($M_2^2\ll nM_1^2$), мы уже не в праве говорить об оценке [190], как о нижней оценке для этого нового класса.

Мы рассмотрели только три значения $p$ и только симплекс в качестве множества, на котором происходит оптимизация. Можно показать, что общий вывод сохранится при рассмотрении всевозможных $1\le p\le\infty$ и всевозможных выпуклых множеств $Q$, на которых происходит оптимизация: *наиболее предпочтительная рандомизация* $e\in RS_2^n(1)$.

В действительности, мы получили намного более общий результат (см. также [44]). Пусть рассматривается задача

$$f(x) = E_\eta\left[f(x;\eta)\right] \to \min_{x\in Q},$$

где $Q$ – выпуклое множество (необязательно ограниченное). Пусть в прямом пространстве выбрана $l_{\bar{p}}$ норма,[36] $1/\bar{p} + 1/\bar{q} = 1$. Введена соответствующая этой норме прокс-функция [44, 128, 259]. Пусть $R_{\bar{p}}^2$ – "расстояние" Брэгмана от точки старта до решения, посчитанное согласно этой прокс-функции [44, 128, 259]. Приводимая ниже Таблица 5.2.4 была ранее известна при $\bar{p} = \bar{q} = 2$ (при $\bar{q} = 2$ в сильно выпуклом случае можно убрать множитель $\ln n$).

Таблица 5.2.4

| Выполнены условия 5.2.1, 5.2.2 подраздела 5.2.2 и условие[37] 5.2.4 подраздела 5.2.3 | $f(x)$ – выпуклая функция | $f(x)$ – $\gamma_{\bar{p}}$-сильно выпуклая функция в $l_{\bar{p}}$ норме |
|---|---|---|
| $2 \le \bar{q} \le \ln n$ | $O\left(\dfrac{\bar{q}M_2^2R_{\bar{p}}^2n^{2/\bar{q}}}{\varepsilon^2}\right)$ | $O\left(\dfrac{\bar{q}M_2^2n^{2/\bar{q}}\ln n}{\gamma_{\bar{p}}\varepsilon}\right)$ |
| $\ln n \le \bar{q} \le \infty$ | $O\left(\dfrac{M_2^2R_{\bar{p}}^2\ln n}{\varepsilon^2}\right)$ | $O\left(\dfrac{M_2^2\ln^2 n}{\gamma_{\bar{p}}\varepsilon}\right)$ |

Все сказанное выше, относилось не к безградиентным методам, а к методам спуска по случайному направлению, и притом в гладком случае. Однако нижние оценки тут с точностью до логарифмических множителей одинаковы. Выше было показано, как можно для спусков по случайному направлению в гладком случае приблизиться, а в определенных ситуациях, и превзойти нижнюю оценку. Естественно, возникает желание перенести предложенный здесь оптимальный метод и на безградиентные методы так, чтобы сохранить полученную оценку. При этом необходимо определить уровень допустимого шума, при котором это возможно. Собственно, этому и посвящен следующий пункт.

---

[36] $\bar{p} \in [1,2]$ – другие значения, как правило, не интересны [44, 128, 259].

[37] С выполнением условия 5.2.4, есть нюанс, когда $Q$ не ограничено [44]. Однако, можно искусственно компактифицировать $Q$, исходя из того, что по ходу итерационного процесса "расстояние" от текущей точки до решения может быть оценено сверху "расстоянием" от точки старта до решения, умноженным на некоторую константу (см. доказательство теоремы 4 в работе [40], а также раздел 5.1 этой главы 5).

Сейчас же мы остановимся на одном интересном обстоятельстве, выявленном в подразделах 5.2.3, 5.2.4. Получается довольно неожиданная ситуация: оказывается, имеет место сильная зависимость скорости сходимости метода от того, какой способ рандомизации (а по сути сглаживания) выбирать. Причем, как это видно из Таблицы 5.2.3, разница очень существенная. К сожалению, в своем желании сохранить несмещенность оценки субградиента мы “перегнули палку” в случае $p=1$, и, особенно, $p=\infty$. Несмещенность в этих случаях досталась нам дорогой ценой – большой оценкой дисперсии соответствующих оценок. Собственно, при предельном переходе $\mu \to 0+$ мы унаследовали большую дисперсию, что и наблюдали в Таблице 5.2.3. Естественно, в этой связи задаться вопросом: а может быть рандомизация $e \in RS_2^n(1)$ оптимальна только в классе несмещенных оценок? А если допускать смещение (bias), то, возможно, можно будет добиться лучшего, как, скажем, в случае оптимальных оценок в математической статистике [70] (см. также [166])? Оказывается, что если допускать смещение, рандомизация $e \in RS_2^n(1)$ по-прежнему будет оптимальной (с точностью до логарифмического множителя). Чтобы это пояснить, мы продолжим рассмотрение гладкого случая, с возможностью получения на каждом шаге (итерации) от оракула незашумленной производной по указанному нами направлению. Рассмотрим более общую схему (см., например, [190]). Пусть $Z$ – случайный вектор с корреляционной матрицей $E_Z\left[ZZ^T\right]=I_n$. Тогда

$$g(x;Z,\eta)=\left\langle \nabla_x f(x;\eta),Z\right\rangle Z = ZZ^T\nabla_x f(x;\eta).$$

Очевидно, что

$$E_{Z,\eta}\left[g(x;Z,\eta)\right]=\nabla f(x).$$

Оказывается, можно улучшить оценку, соответствующую $p=1$, выбирая в этом подходе случайный вектор $Z$ так, чтобы каждая компонента принимала независимо и равновероятно одно из двух значений $1$, $-1$ (равномерное распределение на Хэмминговском кубе) – см. Таблицу 5.2.1. Тогда[38] [190] (см. Таблицу 5.2.3)

---

[38] Впрочем, основная выкладка, поясняющая формулу, достаточно тривиальна

$$E_{Z,\eta}\left[\left\|g(x;Z,\eta)\right\|_\infty^2\right]=E_{Z,\eta}\left[\left\|\left\langle \nabla_x f(x;\eta),Z\right\rangle Z\right\|_\infty^2\right]=$$

$$=E_{Z,\eta}\left[\left\|\left\langle \nabla_x f(x;\eta),Z\right\rangle\right\|_\infty^2\right]=E_{Z,\eta}\left[\left\langle \nabla_x f(x;\eta),Z\right\rangle^2\right]=E_\eta\left[\nabla_x f(x;\eta)^T \underbrace{E_Z\left[ZZ^T\right]}_{I_n}\nabla_x f(x;\eta)\right]\le M_2^2.$$

$$\mathrm{O}\left(\frac{nM_2^2\ln n}{\varepsilon^2}\right)\to\mathrm{O}\left(\frac{M_2^2\ln n}{\varepsilon^2}\right),$$

что улучшает приведенную ранее оценку с рандомизацией на евклидовой сфере на логарифмический множитель.

Если выбрать $Z\in\mathrm{N}(0,I_n)$, то получим в точности те же самые оценки, что получали ранее с рандомизацией на евклидовой сфере.

Если выбрать[39] (см. Таблицу 5.2.1)

$$Z=\sqrt{n}(\underbrace{0,...,0,1}_{i},0,...,0),$$

где случайная величина $i$ независимо и равновероятно принимает значения $1,...,n$, т.е.

$$P(i=k)=\frac{1}{n},\ k=1,...,n,$$

и считать, что в прямом пространстве выбрана (в связи со свойствами множества $Q$) норма $l_2$, то такая покомпонентная рандомизация (см. также раздел 5.1 этой главы 5) приводит к аналогичным оценкам, даваемым рандомизацией на евклидовой сфере при $\overline{p}=\overline{q}=2$, что отражено в Таблице 5.2.5 [190].

Таблица 5.2.5

| $\overline{p}=\overline{q}=2$ | $f(x)$ – выпуклая функция | $f(x)$ – $\gamma_2$-сильно выпуклая функция в $l_2$ норме |
|---|---|---|
| Рандомизация на евклидовой сфере | $\mathrm{O}\left(\frac{M_2^2R_2^2n}{\varepsilon^2}\right)$ | $\mathrm{O}\left(\frac{M_2^2n}{\gamma_2\varepsilon}\right)$ |
| Покомпонентная рандомизация | $\mathrm{O}\left(\frac{M_2^2R_2^2n}{\varepsilon^2}\right)$ | $\mathrm{O}\left(\frac{M_2^2n}{\gamma_2\varepsilon}\right)$ |

Таким образом, мы видим, что в случае рассмотрения методов спуска по случайному направлению (покомпонентных методов) вполне можно рассчитывать на альтернативный способ получения оптимальных методов (оценок). Причем в последнем случае (покомпонентной рандомизации), на самом деле, мы можем существенно выиграть в стоимости одной итерации. Ранее мы себе такую задачу в этом разделе не ставили. Мы хотели минимизировать число обращений к оракулу (за значением функции, за производной по направлению), гарантирующих достижения заданной точности по функции. Если же минимизировать общую вычислительную сложность (число арифметических операций), то покомпонентные методы для большого класса важных в приложениях задач позволяют эффективно организовывать пересчет компонент градиента, т.е. не рассчитывать их каждый раз заново, что позволяет серьезно сэкономить в общих трудозатратах по сравнению с рандомизацией на евклидовой сфере (см., например, [44, 44, 138, 327], а также раздел 1 этой главы). Отметим также, что в последней строчке Таблицы 5.2.5 константу $M_2^2$ можно оценивать как среднее значение по направлениям координатных осей, в то время как в предпоследней строчке Таблицы 5.2.5 $M_2^2$ оценивается по худшему направлению [44].

[39] Собственно, именно с дискретных аналогов такого подхода и начиналось изучение безградиентных методов [99, 231, 314].

Сказанное выше относилось к методам спуска по случайному направлению. Оказывается [190], что эти результаты можно перенести и на безградиентные методы. Для этого вводится аналог $g(x;Z,\eta)$:

$$g^{\tau}(x;Z,\eta)=\frac{f(x+\tau Z;\eta)-f(x;\eta)}{\tau}Z,$$

аналогично подразделу 5.2.3 можно ввести и шумы $g_{\delta}^{\tau}(x;Z,\eta)$. К сожалению, даже при $\delta=0$ мы не получаем несмещенность, т.е. не выполняется условие 5.2.2 подраздела 5.2.2, необходимое для справедливости теоремы 5.2.1 подраздела 5.2.2, которой мы пользуемся. К счастью, у теоремы 5.2.1 есть обобщение (см., например, [227]) не только на произвольные выпуклые множества $Q$ (что мы ранее уже неявно использовали при заполнении Таблиц 5.2.4, 5.2.5), но и на случай, когда условие 5.2.2 выполняется неточно (это как раз сейчас наш случай). Именно, исходя из такого обобщения [227], можно перенести (без изменения) выписанные оценки (при условии достаточной малости $\tau$ и $\delta=0$) на безградиентные методы [190], причем рассуждения [190] можно обобщить и на случай $\delta>0$, контролируя уровень шума (мы не будем в этом разделе приводить соответствующие выкладки). Более того, отмеченное обобщение (из работы [227]) теоремы 5.2.1 позволяет не делать никаких ограничений (типа предположения 5.2.1 подраздела 5.2.3) на шум, кроме должной малости уровня шума. Сами оценки (числа итераций) при этом удается сохранить, но за счет ужесточения требований к уровню шума. Схематично детали такого обобщения описаны в [44, 155] (см. также раздел 2.1 главы 2 и раздел 6.2 главы 6).

### 5.2.5 Перенесение результатов подраздела 5.2.4 на безградиентные методы

Итак, рассмотрим

$$g_\delta^\mu(x;e,\eta)=\frac{n}{\mu}\Big(f(x+\mu e;\eta)+\tilde{\delta}_{x+\mu e}(\eta)-\big(f(x;\eta)+\tilde{\delta}_x(\eta)\big)\Big)e,$$

где $e\in RS_2^n(1)$. Поскольку (см. подраздел 5.2.3)

$$E_{e,\eta}\left[g_\delta^\mu(x;e,\eta)\right]\equiv\nabla f^\mu(x),$$

то для возможности использования теоремы 5.2.1 и схемы подраздела 5.2.3 нужно аккуратно ограничить сверху (в случае $Q=S_n(1)$ имеем $\bar{q}=\infty$) $E_{e,\eta}\left[\left\|g_\delta^\mu(x;e,\eta)\right\|_{\bar{q}}^2\right]$. Далее мы сконцентрируемся именно на этой задаче. Здесь мы не будем бороться за то, чтобы получить оценки вероятностей больших уклонений.

Рассмотрим гладкий случай, в данном случае это подразумевает, что дополнительно к условиям 5.2.1, 5.2.2 подраздела 5.2.2 и условию 5.2.4 подраздела 5.2.3 имеет место условие 5.2.5 подраздела 5.2.3 (в обоих условиях предполагается, что $p'=2$).

Из определения $g_\delta^\mu(x;e,\eta)$ и предположения 5.2.1 подраздела 5.2.2 имеем

$$E_{e,\eta}\left[\left\|g_\delta^\mu(x;e,\eta)\right\|_{\bar{q}}^2\right]=\frac{n^2}{\mu^2}E_{e,\eta}\left[\Big(f(x+\mu e;\eta)-f(x;\eta)+\big(\tilde{\delta}_{x+\mu e}(\eta)-\tilde{\delta}_x(\eta)\big)\Big)^2\|e\|_{\bar{q}}^2\right]=$$

$$=\frac{n^2}{\mu^2}E_{e,\eta}\Big[\Big(\big(f(x+\mu e;\eta)-f(x;\eta)-\mu\langle\nabla_x f(x;\eta),e\rangle\big)+\mu\langle\nabla_x f(x;\eta),e\rangle+$$

$$+\big(\tilde{\delta}_{x+\mu e}(\eta)-\tilde{\delta}_x(\eta)\big)\Big)^2\|e\|_{\bar{q}}^2\Big]. \qquad (5.2.9)$$

Поскольку

$$\left|f(x+\mu e;\eta)-f(x;\eta)-\mu\langle\nabla_x f(x;\eta),e\rangle\right|\le L_2(\eta)\mu^2/2, \qquad (5.2.10)$$

$$\left|\tilde{\delta}_{x+\mu e}(\eta)-\tilde{\delta}_x(\eta)\right|\le 2\delta,$$

$$(a+b+c)^2\le 3a^2+3b^2+3c^2,$$

то

$$E_{e,\eta}\left[\left\|g_\delta^\mu(x;e,\eta)\right\|_{\bar{q}}^2\right]\le\frac{3}{4}n^2L_2^2\mu^2E_e\left[\|e\|_{\bar{q}}^2\right]+3n^2E_{e,\eta}\left[\langle\nabla_x f(x;\eta),e\rangle^2\|e\|_{\bar{q}}^2\right]+12\frac{\delta^2n^2}{\mu^2}E_e\left[\|e\|_{\bar{q}}^2\right].$$

Наиболее интересные ситуации это $\bar{q}=2$ и $\bar{q}=\infty$:

$$E_{e,\eta}\left[\left\|g_\delta^\mu(x;e,\eta)\right\|_{\bar{q}}^2\right]\le 3nM_2^2+\frac{3}{4}n^2L_2^2\mu^2+12\frac{\delta^2n^2}{\mu^2}\ \text{(при } \bar{q}=2\text{)},$$

$$E_{e,\eta}\left[\left\|g_\delta^\mu(x;e,\eta)\right\|_{\bar q}^2\right]\le 4\ln n M_2^2+3n\ln n L_2^2\mu^2+48\frac{\delta^2 n\ln n}{\mu^2}\ \text{(при } \bar q=\infty\text{)}.$$

Выберем $\mu$ согласно условию (5.2.4) подраздела 5.2.3 $L_2\mu^2/2\le\varepsilon/2$, т.е. $\mu\le\sqrt{\varepsilon/L_2}$. Следующее условие на $\mu$ и на допустимый уровень шума $\delta$ получим, исходя из желания обеспечить выполнение неравенства (константа 5 здесь выбрана для определенности)

$$E_{e,\eta}\left[\left\|g_\delta^\mu(x;e,\eta)\right\|_{\bar q}^2\right]\le 5nM_2^2\ \text{(при } \bar q=2\text{)},$$

$$E_{e,\eta}\left[\left\|g_\delta^\mu(x;e,\eta)\right\|_{\bar q}^2\right]\le 5\ln n M_2^2\ \text{(при } \bar q=\infty\text{)}.$$

Отсюда можно получить

$$\mu=\min\left\{\max\left\{\frac{\varepsilon}{2M_2},\sqrt{\frac{\varepsilon}{L_2}}\right\},\frac{M_2}{L_2}\sqrt{\frac{4}{3n}}\right\},\quad \delta\le\frac{M_2\mu}{\sqrt{12n}}\ \text{(при } \bar q=2\text{)},$$

$$\mu=\min\left\{\max\left\{\frac{\varepsilon}{2M_2},\sqrt{\frac{\varepsilon}{L_2}}\right\},\frac{M_2}{L_2}\sqrt{\frac{1}{6n}}\right\},\quad \delta\le\frac{M_2\mu}{\sqrt{96n}}\ \text{(при } \bar q=\infty\text{)},$$

Подобно алгоритму (5.2.2) подраздела 5.2.2 опишем оптимальный алгоритм (см. также подраздел 5.2.3, только в подразделе 5.2.3 используется другая рандомизация) для задачи (5.2.1) и оракула из предположения 5.2.1 подраздела 5.2.2. Положим $x_i^1=1/n$, $i=1,...,n$. Пусть $t=1,...,N-1$.

$$x_i^{t+1}=\frac{\exp\left(-\frac{1}{\beta_{t+1}}\sum_{k=1}^{t}\left[g_\delta^\mu\left(x^k;e^k,\eta^k\right)\right]_i\right)}{\sum_{l=1}^{n}\exp\left(-\frac{1}{\beta_{t+1}}\sum_{k=1}^{t}\left[g_\delta^\mu\left(x^k;e^k,\eta^k\right)\right]_l\right)},\quad i=1,...,n,\quad \beta_t=M_2\sqrt{5t},$$

где $[z]_i$ – $i$-я координата вектора $z$.

**Теорема 5.2.3.** *Пусть мы располагаем оракулом из предположений 5.2.1, 5.2.2 с* $\delta\le M_2\mu/\sqrt{96n}$ *и справедливы условия 5.2.1, 5.2.2 подраздела 5.2.2 и условия 5.2.4, 5.2.5 подраздела 5.2.3 (в которых* $p'=2$*). Тогда для задачи (5.2.1) и описанного выше алгоритма при*

$$N=\left\lceil\frac{80M_2^2\ln^2 n}{\varepsilon^2}\right\rceil$$

*имеет место оценка*

$$E\left[f\left(\frac{1}{N}\sum_{k=1}^{N}x^k\right)\right]-f_*\le\varepsilon.$$

**Доказательство.** Применим теорему 5.2.1 (с учетом выкладок подраздела 5.2.3) с

$$N = \left\lceil \frac{4\cdot 5M_2^2 \ln n}{\left(\varepsilon/2\right)^2} \ln n \right\rceil$$

к функции $f^{\mu}\left(x\right)$. □

Согласно [190] эта оценка оптимальна с точностью до мультипликативного множителя. Более того, подобно подразделу 5.2.4 можно заметить, что в определенных случаях полученная оценка будет лучше нижней (оптимальной) оценки [190].

При дополнительных условиях (уточняющих условия 5.2.4, 5.2.5) здесь также как и в подразделе 5.2.3 можно получить оценки вероятностей больших уклонений, однако мы не будем здесь приводить соответствующие оценки.

### 5.2.6 Заключение

В разделе 5.2 предложены эффективные методы нулевого порядка (также говорят прямые методы или безградиентные методы) для задач стохастической выпуклой оптимизации на симплексе и более общих выпуклых множествах с хорошей проксимальной структурой. Методы строились на базе обычного зеркального спуска для задач стохастической оптимизации. Вместо стохастического градиента в алгоритм зеркального спуска подставлялась специальная конечная разность, аппроксимирующая стохастический градиент. При правильном пересчете размера шага, получался эффективный метод, работающий по известным нижним оценкам, и даже их немного улучшающий при определенных условиях.

Все полученные результаты, кроме третьего столбца Таблицы 5.2.4 (здесь нам известен только результат для $\overline{p} = \overline{q} = 2$) переносятся на онлайн постановки [54]. Детали будут изложены в разделе 6.2 главы 6.

Оригинальность результатов раздела обеспечивается за счет рассмотрения неточного оракула, выдающего на каждой итерации в двух разных точках зашумленные значение оптимизируемой функции на одной и той же реализации. Наличие указанных (дополнительных) шумов “моделирует” практическую неустойчивость конечного дифференцирования, положенного в основу практически всех безградиентных методов. В частности, наличие таких шумов “в первом приближении” может моделировать конечность длины мантиссы.

В отличие от большинства других работ, в данном разделе также прорабатывался вопрос оптимального сочетания способа рандомизации при конструировании “дискретного стохастического градиента”, использующегося в методах вместо недоступного настоящего стохастического градиента, с выбором прокс-структуры, определяемой геометрией выпук-лого множества, на котором происходит оптимизация. В частности, подробно рассматривался, пожалуй, наиболее интересный пример такого множества (после евклидова шара) – симплекс. Ответ оказался достаточно универсальным: для гладких задач оптимальная рандомизация (в независимости от структуры множества) – рандомизация на евклидовой сфере. Для негладких задач в случае достаточно больших шумов, этот ответ уже может быть не верен. Соответствующий пример разбирался в подразделе 5.2.3.

# Глава 6 Онлайн оптимизация с точки зрения выпуклой оптимизации

## 6.1 Об эффективности одного метода рандомизации зеркального спуска в задачах онлайн оптимизации

### 6.1.1 Введение

В конце 70-х годов А.С. Немировским и Д.Б. Юдиным был предложен итерационный метод решения негладких выпуклых задач оптимизации (см., например, [91]), который можно интерпретировать как разновидность метода проекции субградиента, когда проектирование понимается, например, в смысле расстояния Брэгмана (Кульбака–Лейблера) [146], или как прямо-двойственный метод [269] (см. также главу 3). Этот метод, получивший название *метода зеркального спуска* (МЗС), позволяет хорошо учитывать структуру множества, на котором происходит оптимизация (например, симплекса) – ранее этот метод нам уже неоднократно встречался в диссертации. Как и многие другие методы решения негладких выпуклых оптимизационных задач, этот метод требует $O\left(M^2R^2/\varepsilon^2\right)$ итераций, где $\varepsilon$ – точность найденного решения по функции, что соответствует нижним оценкам по $\varepsilon$ [91]. Однако константа $M$, равномерно ограничивающая норму субградиента оптимизируемой функции, и размер множества $R$ зависят от выбора нормы в пространстве, в котором ведется оптимизация. Так, если мы выбрали норму в нашем пространстве $l_p$ ($1 \le p \le \infty$), то $M$ – есть сопряженная $l_q$-норма субградиента ($1/p + 1/q = 1$), а $R^2$ – есть “размер” множества в “метрике” сильно выпуклой относительно $l_p$, с константой сильной выпуклости $\alpha \ge 1$.

**Замечание 6.1.1.** Слово “размер” взято в кавычки, потому что в действительности то, что задается, мы интерпретируем в данном контексте как квадрат размера (физически правильнее квадратом размера называть $R^2/\alpha$), поскольку “метрика” сильно выпуклая относительно нашей “рабочей” нормы в этом пространстве. Слово “метрика” взято в кавычки, потому что может быть не выполнено одно свойств метрики – нет симметричности, например, для расстояния (дивергенции) Брэгмана.

При таких предположениях говорят, что выбранная “метрика” порождает прокс-структуру на множестве. Например, когда множество, на котором происходит оптимизация является симплексом, то, как правило, выбирают $p = 1$, а “метрику” задают расстоянием Брэгмана. При этом “проекция” на симплекс согласно такому расстоянию считается по явным формулам (экспоненциальное взвешивание). В работе [125] была выдвинута гипотеза о том, что применительно к задачам стохастической оптимизации на единичном симплексе (не онлайн) такой выбор нормы и расстояния являются наилучшими с точки зрения зависимости $M^2R^2$ от размера пространства $n$ (в типичных приложениях эта зависимость $\sim \ln n$). Однако в определенных ситуациях (в задачах о многоруких бандитах, когда $M^2R^2 \sim n\ln n$) удается выиграть логарифмический по $n$ фактор (см. пример 6.1.1 раздела 6.1.4), более подходящим образом выбирая расстояние [164]. При этом теряется возможность явного вычисления проекции.

В работах [88, 125, 164, 211, 225, 269] исследовались стохастические (рандомизированные) версии МЗС. В том числе и онлайн [164]. При этом анализировалась ситуация, когда именно градиент функции выдается оракулом со случайными шумами, но несмещенным образом. Такая релаксация детерминированного МЗС оказалась весьма полезной применительно к задачам адаптивного агрегирования оценок [125], оптимизации в пространствах огромных размеров [88, 225], задачах о многоруких бандитах и т.п. [164, 226, 250, 251].

В работе [91] была также отмечена возможность онлайн интерпретации МЗС. Впоследствии, у разных авторов можно найти заметки на эту тему [164, 166, 216, 217, 250, 251, 269, 297, 306]. Наблюдение состоит в том, что ничего не изменится с точки зрения изучения сходимости метода (и его стохастической версии), если на каждом шаге допускать, что функция меняется, причем, возможно, враждебным образом (при этом оставаясь в классе выпуклых функций с ограниченной нормой субградиента).

В данном разделе приводятся две версии стохастического онлайн МЗС. Первая версия - более менее классическая. Приблизительно в таком же виде её уже можно было встретить в литературе у разных авторов [88, 125, 164, 166, 211, 216, 217, 225, 250, 251, 269, 297, 306]. Точнее говоря, предложенная версия аккумулирует в себе в виде частных случаев многие известные ранее версии МЗС.  Вторая неявно была предложена в работе [211] применительно к поиску равновесия в антагонистической матричной игре (онлайн модификация в [211] не была затронута, равно как и связь предложенного метода с МЗС). Согласно работе [211] мы рандомизируем не на этапе вычисления субградиента функции, как это общепринято [88, 225], а на этапе проектирования на симплекс. В результате получается покомпонентный субградиентный спуск со случайным выбором компоненты, который, как будет ниже показано, допускает онлайн интерпретацию. Получив такой метод, мы расширяем множество тех задач онлайн оптимизации, к которым можно применять МЗС.

### 6.1.2 Онлайн метод зеркального спуска со стохастическим субградиентом

Рассмотрим задачу стохастической онлайн оптимизации (запись $E_{\xi^k}\left[f_k\left(x;\xi^k\right)\right]$ означает, что математическое ожидание берется по $\xi^k$, то есть $x$ и $f_k$ понимаются в такой записи не случайными)

$$\frac{1}{N}\sum_{k=1}^{N}E_{\xi^k}\left[f_k\left(x;\xi^k\right)\right]\to\min_{x\in S_n(1)},\ S_n(1)=\left\{x\geq 0:\ \sum_{i=1}^{n}x_i=1\right\},\qquad(6.1.1)$$

при следующих **условиях**:

6.1.1 $E_{\xi^k}\left[f_k\left(x;\xi^k\right)\right]$ – выпуклые функции (по $x$), для этого достаточно выпуклости по $x$ функций $f_k\left(x;\xi^k\right)$ и независимости распределения $\xi^k$ от $x$;

6.1.2 Существует такой вектор $\nabla_x f_k\left(x;\xi^k\right)$, который для компактности будем называть субградиентом, хотя последнее верно не всегда (см. пример 6.1.1 подраздела 6.1.4), что

$$E_{\xi^k}\left(\nabla_x f_k\left(x;\xi^k\right)-\nabla_x E_{\xi^k}\left[f_k\left(x;\xi^k\right)\right]\middle|\Xi^{k-1}\right)\equiv 0,$$

где $\Xi^{k-1}$ – $\sigma$-алгебра, порожденная случайными величинами $\xi^1,\ldots,\xi^{k-1}$. Далее везде в разделе мы будем использовать обозначения обычного градиента для векторов, которые мы назвали здесь субградиентами. В частности, если мы имеем дело с обычным субградиентом, то запись $\nabla_x f_k\left(x;\xi^k\right)$ в вычислительном контексте (например, в итерационной процедуре МЗС, описанной ниже) означает какой-то его элемент (не важно какой именно), а если в контексте проверки условий (например, в условии 6.1.3 ниже), то $\nabla_x f_k\left(x;\xi^k\right)$ пробегает все элементы субградиента (говорят также, субдифференциала);

6.1.3 $\left\|\nabla_x f_k\left(x;\xi^k\right)\right\|_\infty\leq M$ – (равномерно, с вероятностью 1) ограниченный субградиент. Для справедливости части утверждений достаточно требовать одно из следующих (более слабых) условий:

$$\text{а) } E_{\xi^k}\left[\left\|\nabla_x f_k\left(x;\xi^k\right)\right\|_\infty^2\right]\leq M^2;\ \text{ б) } E_{\xi^k}\left[\exp\left(\frac{\left\|\nabla_x f_k\left(x;\xi^k\right)\right\|_\infty^2}{M^2}\right)\middle|\Xi^{k-1}\right]\leq\exp(1).$$

Задача (6.1.1) является лишь компактной (и далеко не полной) записью настоящей постановки задачи стохастической онлайн оптимизации. В действительности, требуется подобрать последовательность $\left\{x^k\right\}$ (в подразделе 6.1.2 $\left\{x^k\right\}\in S_n(1)$, а в подразделе 6.1.3 и ряде примеров подраздела 6.1.4 при дополнительном ограничении, что $\{x^k\}$ выбираются с возможными повторениями среди вершин симплекса $S_n(1)$) исходя из доступной исторической информации ($x^k$ может зависеть только от $\{x^1,\xi^1,f_1(\cdot);...;x^{k-1},\xi^{k-1},f_{k-1}(\cdot)\}$) так, чтобы минимизировать псевдо регрет:

$$\frac{1}{N}\sum_{k=1}^{N}E_{\xi^k}\left[f_k\left(x^k;\xi^k\right)\right]-\min_{x\in S_n(1)}\frac{1}{N}\sum_{k=1}^{N}E_{\xi^k}\left[f_k\left(x;\xi^k\right)\right]$$

или регрет

$$E_{\xi^1,...,\xi^N}\left[\min_{x\in S_n(1)}\frac{1}{N}\sum_{k=1}^{N}\left(f_k\left(x^k;\xi^k\right)-f_k\left(x;\xi^k\right)\right)\right].$$

В данном пункте мы сосредоточимся на минимизации псевдо регрета на основе информации $\{\nabla f_1(x^1;\xi^1);...;\nabla f_{k-1}(x^{k-1};\xi^{k-1})\}$ при расчете $x^k$. Онлайновость постановки задачи допускает, что на каждом шаге $k$ функция $f_k$ может подбираться из рассматриваемого класса функций враждебно по отношению к используемому нами методу генерации последовательности $\{x^k\}$. В частности, в этом пункте $f_k$ может зависеть от $\{x^1,\xi^1,f_1(\cdot);...;x^{k-1},\xi^{k-1},f_{k-1}(\cdot);x^k\}$.

В целом в данном разделе мы ограничиваемся рассмотрением задач минимизации псевдо регрета, когда на каждом шаге мы можем получить независимую реализацию стохастического субградиента в одной указанной нами (допустимой) точке. Приведенные в разделе результаты можно распространить на случай, когда градиент выдается не точно (с не случайной ошибкой), выдается не полностью (скажем, вместо градиента выдается производная по выбранному направлению) или вместо градиента выдается только значение функции [166]. Впрочем, немного об этом написано далее (пример 6.1.1 подраздела 6.1.4). Другим способом релаксации исходной постановки является возможность несколько раз обращаться на одном шаге за значением градиента функции и(или) значением самой функции и взвешивать $f_k(x;\xi^k)$ в (6.1.1) разными весами (см. [166], [311]). Результаты этого раздела также можно распространить и на случай, когда функции $E_{\xi^k}\left[f_k\left(x;\xi^k\right)\right]$ равномерно сильно выпуклые по $x$ с константой $\mu$ [164, 216, 217, 251]. При этом выбирается евклидова прокс-структура, поскольку в сильно выпуклом случае (в стохастическом и не стохастическом) игра на выборе прокс-структуры не может дать выгоды. Неулучшаемая и достижимая оценка в этом случае будет иметь следующий вид $\varepsilon=\mathrm{O}\left(M^2\ln N/(\mu N)\right)$. В не онлайн стохастическом случае эта оценка (с точностью до фактора $\ln N$) также будет неулучшаемой. Отметим при этом, что (для рассмотренных выпуклых и сильно выпуклых задачах) в отличие от не онлайн случая, в онлайн случае игра на гладкости функций $E_{\xi^k}\left[f_k\left(x;\xi^k\right)\right]$ и(или) отсутствии стохастичности ($f_k\left(x;\xi^k\right)\equiv f_k\left(x\right)$) не дает никаких дивидендов (выписанные нижние оценки сохранятся).

Можно обобщить приведенную постановку и последующие результаты на задачи композитной оптимизации [266]. Если функция, которая добавляется ко всем $f_k$, – линейная, то мы даже можем ее не знать (просто знать, что она есть и одна и та же). Тогда опустив первое слагаемое в сумме (6.1.1) и переписав условие 6.1.3 в виде

$$\left\|\nabla_x f_k\left(x;\xi^k\right)-\nabla_x f_m\left(y;\xi^m\right)\right\|_\infty \le M \text{, для всех } x,y\in S_n(1),\ k,m\in\mathbb{N},$$

можно перенести результаты статьи [268] на такой онлайн контекст.

Выше мы исходили из того, что оптимизация ведется на единичном симплексе. Возникает вопрос: насколько все, что приведено в этом разделе, обобщается на более общий случай? Собственно говоря, ответ на этот вопрос частично известен уже давно [91]. Приведенные в подразделе 6.1.2 рассуждения универсальны, то есть если исходить из оптимизации на каком-нибудь другом выпуклом компакте (от условия компактности (ограниченности) множества можно отказаться [269], поскольку, в действительности, в оценку числа итераций входит не размер множества, а “расстояние” от точки старта до решения), то задав норму в прямом пространстве и расстояние (сильно выпуклое относительно этой нормы), согласно которому будет осуществляться проектирование субградиента на этот компакт, можно повторить аналогичные рассуждения [227].

Для решения задачи (6.1.1) воспользуемся адаптивным методом зеркального спуска (точнее двойственных усреднений) в форме [125, 269]. Положим $x_i^1=1/n$, $i=1,...,n$. Пусть $t=1,...,N-1$.

**Алгоритм МЗС1-адаптивный / Метод двойственных усреднений**

$$x_i^{t+1}=\frac{\exp\left(-\frac{1}{\beta_{t+1}}\sum_{k=1}^{t}\frac{\partial f_k\left(x^k;\xi^k\right)}{\partial x_i}\right)}{\sum_{l=1}^{n}\exp\left(-\frac{1}{\beta_{t+1}}\sum_{k=1}^{t}\frac{\partial f_k\left(x^k;\xi^k\right)}{\partial x_l}\right)},\quad i=1,...,n,\quad \beta_t=\frac{M\sqrt{t}}{\sqrt{\ln n}}.$$

Не сложно показать, что этот метод представим также в виде:

$$x^{t+1}=\arg\min_{x\in S_n(1)}\left\{\sum_{k=1}^{t}\left\{f_k\left(x^k;\xi^k\right)+\left\langle\nabla f_k\left(x^k;\xi^k\right),x-x^k\right\rangle\right\}+\beta_{t+1}V(x)\right\}$$

или

$$\begin{cases} y^k = y^{k-1} - \gamma_k \nabla_x f_k\left(x^k;\xi^k\right) \\ x^{k+1} = \nabla W_{\beta_{k+1}}\left(y^k\right) \end{cases}, \; y^0 = 0, \; \gamma_k \equiv 1, \; \beta_k = \frac{M}{\sqrt{\ln n}}\sqrt{k}, \; k = 1,\dots,N, \qquad (6.1.2)$$

где

$$W_\beta(y) = \sup_{x\in S_n(1)}\left\{\langle y,x\rangle - \beta V(x)\right\} = \beta\ln\left(\frac{1}{n}\sum_{i=1}^n \exp\left(y_i/\beta\right)\right),$$

$$V(x) = \ln n + \sum_{i=1}^n x_i \ln x_i .$$

Рассуждая подобно [125, 225, 269], можно получить следующий результат.

**Теорема 6.1.1.** *Пусть справедливы условия 6.1.1, 6.1.2, 6.1.3.а, тогда*

$$\frac{1}{N}\sum_{k=1}^N E\left[f_k\left(x^k;\xi^k\right)\right] - \min_{x\in S_n(1)}\frac{1}{N}\sum_{k=1}^N E_{\xi^k}\left[f_k\left(x;\xi^k\right)\right] \le 2M\sqrt{\frac{\ln n}{N}}.$$

*Если* $f_k \equiv f$ *, а* $\left\{\xi^k\right\}$ *– независимы и одинаково распределены, как* $\xi$ *, то*

$$E\left[f\left(\frac{1}{N}\sum_{k=1}^N x^k;\xi\right)\right] - \min_{x\in S_n(1)} E_\xi\left[f\left(x;\xi\right)\right] \le 2M\sqrt{\frac{\ln n}{N}}.$$

*Пусть справедливы условия 1, 2, 3, тогда при* $\Omega \ge 0$

$$P_{x^1,\dots,x^N}\left\{\frac{1}{N}\sum_{k=1}^N E_{\xi^k}\left[f_k\left(x^k;\xi^k\right)\right] - \min_{x\in S_n(1)}\frac{1}{N}\sum_{k=1}^N E_{\xi^k}\left[f_k\left(x;\xi^k\right)\right] \ge \frac{2M}{\sqrt{N}}\left(\sqrt{\ln n}+\sqrt{8\Omega}\right)\right\} \le \exp(-\Omega).$$

*Если* $f_k \equiv f$ *, а* $\left\{\xi^k\right\}$ *– независимы и одинаково распределены, как* $\xi$ *, то*

$$P_{x^1,\dots,x^N}\left\{E_\xi\left[f\left(\frac{1}{N}\sum_{k=1}^N x^k;\xi\right)\right] - \min_{x\in S_n(1)} E_\xi\left[f\left(x;\xi\right)\right] \ge \frac{2M}{\sqrt{N}}\left(\sqrt{\ln n}+\sqrt{8\Omega}\right)\right\} \le \exp(-\Omega).$$

**Замечание 6.1.2.** Запись

$$\text{“} P_{x^1,\dots,x^N}\left\{\frac{1}{N}\sum_{k=1}^N E_{\xi^k}\left[f_k\left(x^k;\xi^k\right)\right] - \dots\text{”}$$

означает, что под вероятностью мы считаем математическое ожидание по $\xi^k$, которое, вообще говоря, зависит и от $\xi^1,\dots,\xi^{k-1}$ (мы не предполагаем независимости $\left\{\xi^k\right\}$), как бы “замораживая” (считая не случайными) $x^k$, то есть забывая про то, что $x^k$ тоже зависит от $\xi^1,\dots,\xi^{k-1}$. А вероятность берется как раз по $\left\{x^k\right\}$, с учетом того, что такая зависимость есть (см. определение алгоритма МЗС1).

**Доказательство теоремы 6.1.1.** Возьмем за основу обозначения МЗС1, приведенные в (6.1.2). Рассуждая далее аналогично [125, 225, 269], получим (заметим, что $\beta_{k+1} \ge \beta_k > 0$):

$$W_{\beta_{k+1}}\left(y^k\right) \le W_{\beta_k}\left(y^k\right) = W_{\beta_k}\left(y^{k-1}\right) + \int\limits_0^1 \left(y^k - y^{k-1}\right)^T \nabla W_{\beta_k}\left(\tau y^k + (1-\tau) y^{k-1}\right) d\tau =$$

$$= W_{\beta_k}\left(y^{k-1}\right) - \gamma_k \nabla_x f_k\left(x^k;\xi^k\right)^T \nabla W_{\beta_k}\left(y^{k-1}\right) -$$

$$-\gamma_k \nabla_x f_k\left(x^k;\xi^k\right)^T \int\limits_0^1 \left(\nabla W_{\beta_k}\left(\tau y^k + (1-\tau) y^{k-1}\right) - \nabla W_{\beta_k}\left(y^{k-1}\right)\right) d\tau \overset{(*)}{\le}$$

$$\overset{(*)}{\le} W_{\beta_k}\left(y^{k-1}\right) - \gamma_k \nabla_x f_k\left(x^k;\xi^k\right) \nabla W_{\beta_k}\left(y^{k-1}\right) +$$

$$+\gamma_k \left\|\nabla_x f_k\left(x^k;\xi^k\right)\right\|_\infty \int\limits_0^1 \left\|\nabla W_{\beta_k}\left(\tau y^k + (1-\tau) y^{k-1}\right) - \nabla W_{\beta_k}\left(y^{k-1}\right)\right\|_1 d\tau \le$$

$$\le W_{\beta_k}\left(y^{k-1}\right) - \gamma_k \nabla_x f_k\left(x^k;\xi^k\right) \nabla W_{\beta_k}\left(y^{k-1}\right) + \frac{\gamma_k^2 \left\|\nabla_x f_k\left(x^k;\xi^k\right)\right\|_\infty^2}{2\alpha\beta_k},$$

последнее неравенство следует из того, что [125, 225, 269]:

$$\left\|\nabla W_\beta\left(\tilde{y}\right) - \nabla W_\beta\left(y\right)\right\|_1 \le \frac{1}{\alpha\beta}\left\|\tilde{y} - y\right\|_\infty,$$

где $\alpha = 1$ – константа сильной выпуклости $V(x)$ в 1-норме. Мы специально выделили (*) неравенство, которое иногда (например, в задачах о многоруких бандитах) бывает довольно грубым. В работах [164, 166] указан способ уточнения этого неравенства.

Суммируя эти неравенства, учитывая, что $\nabla W_{\beta_k}\left(y^{k-1}\right) = x^k$ и формулу (6.1.2), получим

$$\sum_{k=1}^N \gamma_k \left(x^k - x\right)^T \nabla_x f_k\left(x^k;\xi^k\right) \le W_{\beta_1}\left(y^0\right) + x^T y^N - W_{\beta_{N+1}}\left(y^N\right) + \sum_{k=1}^N \frac{\gamma_k^2}{2\alpha\beta_k}\left\|\nabla_x f_k\left(x^k;\xi^k\right)\right\|_\infty^2,$$

где $x^T$ – означает транспонирование вектора $x$, который мы выбираем так, чтобы он доставлял решение задачи (6.1.1). Поскольку [125]

$$W_{\beta_1}\left(y^0\right) = W_{\beta_1}(0) = 0 \text{ и } \beta_{N+1} V(x) \ge x^T y^N - W_{\beta_{N+1}}\left(y^N\right),$$

то

$$\sum_{k=1}^N \gamma_k \left(x^k - x\right)^T \nabla_x f_k\left(x^k;\xi^k\right) \le \beta_{N+1} V(x) + \sum_{k=1}^N \frac{\gamma_k^2}{2\alpha\beta_k}\left\|\nabla_x f_k\left(x^k;\xi^k\right)\right\|_\infty^2.$$

Тогда, из выпуклости функции $E_{\xi^k}\left[f_k\left(x;\xi^k\right)\right]$ по $x$ (в виду условия 6.1.1) следует, что

$$\sum_{k=1}^N \gamma_k \left\{E_{\xi^k}\left[f_k\left(x^k;\xi^k\right)\right] - E_{\xi^k}\left[f_k\left(x;\xi^k\right)\right]\right\} \le$$

$$\le \sum_{k=1}^N \gamma_k \left(x^k - x\right)^T \nabla_x E_{\xi^k}\left[f_k\left(x^k;\xi^k\right)\right] \le$$

$$\le \beta_{N+1}V(x)-\sum_{k=1}^{N}\gamma_k\left(x^k-x\right)^T\left(\nabla_x f_k\left(x^k;\xi^k\right)-\nabla_x E_{\xi^k}\left[f_k\left(x^k;\xi^k\right)\right]\right)+$$

$$+\sum_{k=1}^{N}\frac{\gamma_k^2}{2\alpha\beta_k}\left\|\nabla_x f_k\left(x^k;\xi^k\right)\right\|_\infty^2. \qquad (6.1.3)$$

Возьмем полное (т.е., в отличие от замечания 6.1.2, с учетом зависимости $x^k$ от $\left\{\xi^1,...,\xi^{k-1}\right\}$) математическое ожидание (в два шага $E[\,\cdot\,]=E\left[E_{\xi^k}\left[\,\cdot\,\middle|\Xi^{k-1}\right]\right]$ – для каждого слагаемого свое $k$) от обеих частей неравенства, учитывая условие 6.1.3.а и то, что

$$E\left(\left(x^k-x\right)^T\left(\nabla_x f_k\left(x^k;\xi^k\right)-\nabla_x E_{\xi^k}\left[f_k\left(x^k;\xi^k\right)\right]\right)\right)=$$

$$=E\left[\left(x^k-x\right)^T E_{\xi^k}\left(\nabla_x f_k\left(x^k;\xi^k\right)-\nabla_x E_{\xi^k}\left[f_k\left(x^k;\xi^k\right)\right]\middle|\Xi^{k-1}\right)\right]=0,$$

поскольку $x^k$ – $\Xi^{k-1}$-измеримый вектор и внутреннее условное математическое ожидание в силу условия 6.1.2 равно 0, получим

$$\sum_{k=1}^{N}\gamma_k\left\{E\left[f_k\left(x^k;\xi^k\right)\right]-E_{\xi^k}\left[f_k\left(x;\xi^k\right)\right]\right\}\le\beta_{N+1}V(x)+\sum_{k=1}^{N}\frac{\gamma_k^2}{2\alpha\beta_k}E\left[\left\|\nabla_x f_k\left(x^k;\xi^k\right)\right\|_\infty^2\right]\le$$

$$\le\beta_{N+1}R^2+M^2\sum_{k=1}^{N}\frac{\gamma_k^2}{2\alpha\beta_k},$$

где $R^2=\max_{x\in S_n(1)}V(x)=\ln n$. Подставляя $\gamma_k\equiv1$, и минимизируя правую часть неравенства по неубывающим последовательностям с положительными элементами $\left\{\beta_k\right\}_{k=1}^{N+1}$, не допуская при этом зависимость $\left\{\beta_k\right\}_{k=1}^{N+1}$ от потенциально неизвестного $N$, получим

$$\beta_k=\sqrt{\frac{M^2}{\alpha R^2}}\sqrt{k},$$

$$\frac{1}{N}\sum_{k=1}^{N}E\left[f_k\left(x^k;\xi^k\right)\right]-\min_{x\in S_n(1)}\frac{1}{N}\sum_{k=1}^{N}E_{\xi^k}\left[f_k\left(x;\xi^k\right)\right]\le2\sqrt{\frac{M^2R^2}{\alpha N}}.$$

Для того чтобы доказать первую часть теоремы, осталось подставить $\alpha=1$ и $R^2=\ln n$.

**Замечание 6.1.3.** Если разрешать $\left\{\beta_k\right\}_{k=1}^{N+1}$ зависеть от $N$ (см., например, алгоритм МЗС2-неадаптивный в следующем пункте), то в последней формуле “2”-у можно занести под знак корня [125]. Все это переносится и на последующие рассуждения с вероятностями больших отклонений.

**Замечание 6.1.4.** Строго говоря, в получено оценке в знаменателе вместо $N$ нужно писать $N^2/(N+1)$. Считая, что $N\gg1$, мы пренебрегли этим для компактности записи.

Для доказательства второй части теоремы вернемся к формуле (6.1.3). Из условия 6.1.3 имеем

$$P\left(\sum_{k=1}^{N}\frac{\gamma_k^2}{2\alpha\beta_k}\left\|\nabla_x f_k\left(x^k;\xi^k\right)\right\|_\infty^2 > M^2\sum_{k=1}^{N}\frac{\gamma_k^2}{2\alpha\beta_k}\right)=0. \qquad (6.1.4)$$

Из неравенства Азума–Хефдинга [159], подобно [225, 238], получаем для ограниченной мартингал–разности

$$\left|\left(x-x^k\right)^T\left(\nabla_x f_k\left(x^k;\xi^k\right)-\nabla_x E_{\xi^k}\left[f_k\left(x^k;\xi^k\right)\right]\right)\right|\le 4M$$

следующее неравенство:

$$P\left(\sum_{k=1}^{N}\gamma_k\left(x-x^k\right)^T\left(\nabla_x f_k\left(x^k;\xi^k\right)-\nabla_x E_{\xi^k}\left[f_k\left(x^k;\xi^k\right)\right]\right)\ge 4M\Lambda\sqrt{\sum_{k=1}^{N}\gamma_k^2}\right)\le\exp\left(-\Lambda^2/2\right).$$

Подставляя $\gamma_k\equiv 1$, $\Lambda=\sqrt{2\Omega}$, получим вторую часть теоремы.

**Замечание 6.1.5.** В статье [225] для задачи стохастической выпуклой оптимизации ($\left\{\xi^k\right\}$ – независимые случайные величины) приводится оценка вероятностей больших уклонений с точностью до констант аналогичная оценке, приведенной в теореме. Аналогичная оценка приводится в [225] и для случая, когда вместо условия 6.1.3 предполагается условие 6.1.3.б. При этом в [225] использовалось условие независимости $\left\{\xi^k\right\}$ (при установлении неравенства типа (4) в общем случае и в неравенстве Азума–Хефдинга). Для онлайн оптимизации, как правило, независимость $\left\{\xi^k\right\}$ место не имеет. Тем не менее, условия 6.1.2, 6.1.3.б обеспечивают выполнения этих неравенств в том же виде, как если бы независимость $\left\{\xi^k\right\}$ имела место. Приведем теперь оценки в случае тяжелых хвостов. Если $\left\|\nabla_x f\left(x,\xi\right)\right\|_\infty^2$ имеет степенной хвост ($\alpha>2$)

$$P\left(\frac{\left\|\nabla_x f\left(x,\xi\right)\right\|_\infty^2}{M^2}\ge t\right)=O\left(\frac{1}{t^\alpha}\right),$$

то существует такая константа $C_\alpha>0$, что с вероятностью $\ge 1-\sigma$

$$\frac{1}{N}\sum_{k=1}^{N}E_{\xi^k}\left[f_k\left(x^k;\xi^k\right)\right]-\min_{x\in S_n(1)}\frac{1}{N}\sum_{k=1}^{N}E_{\xi^k}\left[f_k\left(x;\xi^k\right)\right]\le C_\alpha M\frac{\sqrt{\ln n\ln\left(\sigma^{-1}\right)}+\frac{\left(N/\sigma\right)^{1/\alpha}}{N}}{\sqrt{N}}.$$

Если мы не делаем никаких предположений относительно распределений случайных величин $\left\{\xi^k\right\}$, кроме 6.1.1, 6.1.2 и существования первых двух равномерно ограниченных моментов у $\nabla_x f_k\left(x^k;\xi^k\right)$, то из неравенства Маркова и первого неравенства в теореме 6.1.1 (на математические ожидания) имеем: существует такая константа $C>0$, что с вероятностью $\ge 1-\sigma$

$$\frac{1}{N}\sum_{k=1}^{N}E_{\xi^k}\left[f_k\left(x^k;\xi^k\right)\right]-\min_{x\in S_n(1)}\frac{1}{N}\sum_{k=1}^{N}E_{\xi^k}\left[f_k\left(x;\xi^k\right)\right]\le\frac{CM}{\sigma}\sqrt{\frac{\ln n}{N}}.$$

Труднее обстоит дело, если мы хотим оценить регрет или

$$\frac{1}{N}\sum_{k=1}^{N}f_k\left(x^k;\xi^k\right)-\min_{x\in S_n(1)}\frac{1}{N}\sum_{k=1}^{N}f_k\left(x;\xi^k\right).$$

Тем не менее, при дополнительных оговорках и такие выражения можно вероятностно оценивать [166, 238].

Аналогичное замечание имеет место и для теоремы 6.1.2 ниже.

### 6.1.3 Онлайн метод зеркального спуска со стохастической проекцией

Снова рассмотрим постановку задачи стохастической онлайн оптимизации (6.1.1) из подраздела 6.1.2. Но на этот раз будем допускать, что метод генерирования последовательности $\left\{x^k\right\}$ может допускать (внешнюю, дополнительную) рандомизацию. Это допущение позволит частично перенести результаты подраздела 6.1.2 на не выпуклые функции $E_{\xi^k}\left[f_k\left(x;\xi^k\right)\right]$ (см. пример 6.1.4 подраздела 6.1.4), на ситуации, когда по условию задачи $\left\{x^k\right\}$ должны выбираться среди вершин единичного симплекса (примеры 6.1.1 и 6.1.5 подраздела 6.1.4). Также как и раньше онлайновость постановки задачи допускает, что на каждом шаге $k$ функция $f_k$ может подбираться из рассматриваемого класса функций враждебно по отношению к используемому нами методу генерации последовательности $\left\{x^k\right\}$. В частности, $f_k$ может зависеть от $\left\{x^1,\xi^1,f_1(\cdot);...;x^{k-1},\xi^{k-1},f_{k-1}(\cdot)\right\}$, и даже от распределения вероятностей $p^k$ (многорукие бандиты), согласно которому осуществляется выбор $x^k$. Чтобы можно было работать с таким классом задач, нам придется наложить дополнительное **условие**:

6.1.4 На каждом шаге генерирование случайной величины $x^k$ согласно распределению вероятностей $p^k$ осуществляется независимо ни от чего. Выбор $f_k$ осуществляется без знания реализации $x^k$.

Положим $p_i^1=x_i^1=1/n$, $i=1,...,n$. Пусть $t=1,...,N-1$.

**<u>Алгоритм МЗС2-адаптивный / Метод двойственных усреднений</u>**

*Согласно распределению вероятностей*

$$p_i^{t+1} = \frac{\exp\left( -\frac{1}{\beta_{t+1}} \sum_{k=1}^{t} \frac{\partial f_k\left(x^k;\xi^k\right)}{\partial x_i} \right)}{\sum_{l=1}^{n} \exp\left( -\frac{1}{\beta_{t+1}} \sum_{k=1}^{t} \frac{\partial f_k\left(x^k;\xi^k\right)}{\partial x_l} \right)}, \quad i = 1,\dots,n\,,\, \beta_t = \frac{M\sqrt{t}}{\sqrt{\ln n}},$$

*получаем случайную величину* $i(t+1)$, $x_{i(t+1)}^{t+1} = 1,\ x_j^{t+1} = 0,\ j \ne i(t+1)$.

**<u>Алгоритм МЗС2-неадаптивный</u> (заранее известно $N$ )**

*Согласно распределению вероятностей*

$$p_i^{t+1} = \frac{\exp\left( -\frac{1}{\beta_{t+1}} \sum_{k=1}^{t} \gamma_k \frac{\partial f_k\left(x^k;\xi^k\right)}{\partial x_i} \right)}{\sum_{l=1}^{n} \exp\left( -\frac{1}{\beta_{t+1}} \sum_{k=1}^{t} \gamma_k \frac{\partial f_k\left(x^k;\xi^k\right)}{\partial x_l} \right)}, \quad i = 1,\dots,n,$$

$$\gamma_k \equiv M^{-1}\sqrt{2\ln n/N}\,,\ \beta_t \equiv 1,$$

*получаем случайную величину* $i(t+1)$, $x_{i(t+1)}^{t+1} = 1,\ x_j^{t+1} = 0,\ j \ne i(t+1)$.

**Мотивация (ограничимся детерминированным случаем с $\gamma_k \equiv 1$).** Аппроксимируя

$$\min_{x\in S_n(1)} \frac{1}{N}\sum_{k=1}^{t} f_k(x) \approx \min_{x\in S_n(1)} \frac{1}{N}\sum_{k=1}^{t}\left\{ f_k\left(x^k\right) + \left\langle \nabla f_k\left(x^k\right), x - x^k \right\rangle \right\},$$

получим $x_j^{t+1} = 1;\, x_i^{t+1} = 0, i \ne j$ (для простоты считаем, что имеет место не вырожденный случай)

$$j = \arg\max_{i=1,\dots,n}\left\{ \sum_{k=1}^{t}\left[ -\nabla f_k\left(x^k\right) \right]_i \right\}.$$

Поскольку мы работаем с аппроксимацией (нижними аффинными минорантами) исходной задачи, то предлагается немного видоизменить это правило

$$P\left(x_j^{t+1} = 1;\, x_i^{t+1} = 0, i \ne j\right) \overset{def}{=} P_\varsigma\left( j = \arg\max_{i=1,\dots,n}\left\{ \left( \sum_{k=1}^{t}\left[ -\nabla f_k\left(x^k\right) \right]_i \right) + \varsigma_{t,i} \right\} \right),$$

где $\varsigma_{t,i}$ – независимые одинаково распределенные случайные величины по закону Гумбеля [84, 137] с параметром $\beta_{t+1}$, характеризующим среднеквадратичное отклонение $\varsigma_{t,i}$:

$$P\left(\varsigma_{t,i} < \tau\right) = \exp\left\{ -e^{-\tau/\beta_{t+1}} \right\}.$$

Можно показать [137, 250, 269], что

$$E_\varsigma\left[ x^{t+1} \right] = \nabla W_{\beta_{t+1}}\left( -\sum_{k=1}^{t} \nabla f_k\left(x^k\right) \right).$$

**Замечание 6.1.6.** Естественно задаться вопросом: а какое распределение “наиболее подходит” для $\varsigma_{t,i}$, чтобы в случае “враждебной Природы” (то есть в минимаксном смысле) иметь наилучшие оценки [250]? Ответом будет [250]: показательное распределение (точнее $\beta_{t+1}\mathrm{Exp}(1)$, где $\mathrm{Exp}(1)$ – случайная величина, имеющая показательное распределение с параметром 1), которое ведет себя в интересном для анализа диапазоне подобно распределению Гумбеля, но в случае Гумбеля мы явно можем посчитать интересующие нас вероятности. Как правило, такого рода задачи явно не решаются, и распределение Гумбеля является приятным исключением, для которого есть явные формулы.

Приведенная выше мотивация имеет одно интересное приложение в содержательной интерпретации равновесного распределения транспортных потоков. Не много об этом написано [47] (см. также главу 1).

К сожалению, не делая относительно функций $f_k\left(x;\xi^k\right)$ дополнительно никаких предположений, не удается доказать для МЗС2 аналог теоремы 6.1.1. Чтобы можно было сформулировать такой аналог, мы вынуждены будем предполагать, что $f_k\left(x;\xi^k\right)$ – линейные функции по $x$ (можно обобщить и на сублинейные). С одной стороны это существенно сужает класс задач, к которым применим МЗС2. С другой стороны, как будет продемонстрировано в следующем пункте, даже такой узкий класс функций за счет “онлайновости” позволяет применять МЗС2 к довольно широкому кругу задач. Для того чтобы лучше чувствовалась преемственность методов и доказательств их сходимости, далее мы по-прежнему будем использовать общие обозначения $f_k\left(x;\xi^k\right)$, не подчеркивая в формулах линейность.

**Теорема 6.1.2.** *Пусть справедливы условия 6.1.1, 6.1.2, 6.1.3.a, 6.1.4 и* $f_k\left(x;\xi^k\right)$ *– линейные функции по* $x$*, тогда*

$$\frac{1}{N}\sum_{k=1}^{N}E\left[f_k\left(x^k;\xi^k\right)\right]-\min_{x\in S_n(1)}\frac{1}{N}\sum_{k=1}^{N}E_{\xi^k}\left[f_k\left(x;\xi^k\right)\right]\le 2M\sqrt{\frac{\ln n}{N}}.$$

*Для неадаптивного метода “2”-у перед* $M$ *можно занести под знак корня.*

*Кроме того, если справедливы условия 1, 2, 3, 4, то при* $\Omega\ge 0$

$$P_{x^1,\dots,x^N}\left\{\frac{1}{N}\sum_{k=1}^{N}E_{\xi^k}\left[f_k\left(x^k;\xi^k\right)\right]-\min_{x\in S_n(1)}\frac{1}{N}\sum_{k=1}^{N}E_{\xi^k}\left[f_k\left(x;\xi^k\right)\right]\ge\frac{2M}{\sqrt{N}}\left(\sqrt{\ln n}+\sqrt{18\Omega}\right)\right\}\le\exp(-\Omega).$$

*Если* $f_k\left(x;\xi^k\right)\equiv f_k(x)$*, то это неравенство можно уточнить*

$$\frac{2M}{\sqrt{N}}\left(\sqrt{\ln n}+\sqrt{18\Omega}\right)\to\frac{2M}{\sqrt{N}}\left(\sqrt{\ln n}+\sqrt{2\Omega}\right),$$

*при этом же условии для неадаптивного метода можно еще больше уточнить*

$$\frac{2M}{\sqrt{N}}\left(\sqrt{\ln n}+\sqrt{2\Omega}\right)\to\frac{\sqrt{2}M}{\sqrt{N}}\left(\sqrt{\ln n}+2\sqrt{\Omega}\right).$$

**Схема доказательства теоремы 6.1.2.** Доказательство фактически дословно повторяет доказательство теоремы 6.1.1. Небольшая разница лишь в том, что основная формула (6.1.3) перепишется следующим образом:

$$\sum_{k=1}^{N}\gamma_k\left\{E_{\xi^k}\left[f_k\left(x^k;\xi^k\right)\right]-E_{\xi^k}\left[f_k\left(x;\xi^k\right)\right]\right\}\le\sum_{k=1}^{N}\gamma_k\left(x^k-x\right)^T\nabla_x E_{\xi^k}\left[f_k\left(x^k;\xi^k\right)\right]\le$$
$$\le\beta_{N+1}V(x)+\sum_{k=1}^{N}\gamma_k\left(x^k-E_{x^k}\left[x^k\right]\right)^T\nabla_x f_k\left(x^k;\xi^k\right)-$$
$$-\sum_{k=1}^{N}\gamma_k\left(x^k-x\right)^T\left(\nabla_x f_k\left(x^k;\xi^k\right)-\nabla_x E_{\xi^k}\left[f_k\left(x^k;\xi^k\right)\right]\right)+\sum_{k=1}^{N}\frac{\gamma_k^2}{2\alpha\beta_k}\left\|\nabla_x f_k\left(x^k;\xi^k\right)\right\|_\infty^2.$$

Здесь мы просто не много по-другому (по сравнению с доказательством теоремы 6.1.1) переписали неравенство

$$\sum_{k=1}^{N}\gamma_k\left(x^k-x\right)^T f_k\left(x^k;\xi^k\right)\le\beta_{N+1}V(x)+\sum_{k=1}^{N}\frac{\gamma_k^2}{2\alpha\beta_k}\left\|\nabla_x f_k\left(x^k;\xi^k\right)\right\|_\infty^2,$$

используя выпуклость функции $E_{\xi^k}\left[f_k\left(x;\xi^k\right)\right]$ по $x$ в виду условия 6.1.1.

Введем случайные величины

$$Z_k=\left(x^k-E_{x^k}\left[x^k\right]\right)^T\nabla_x f_k\left(x^k;\xi^k\right),\ \tilde{Z}_k=\left(x^k-x\right)^T\left(\nabla_x f_k\left(x^k;\xi^k\right)-\nabla_x E_{\xi^k}\left[f_k\left(x^k;\xi^k\right)\right]\right)$$

Поскольку $f_k\left(x^k;\xi^k\right)=l_k\left(\xi^k\right)^T x^k$ – линейные функции, то по условиям 6.1.4, 6.1.2

$$E_{x^k}\left[Z_k\middle|\Xi^{k-1}\right]\equiv 0,\ E_{\xi^k}\left[\tilde{Z}_k\middle|\Xi^{k-1}\right]\equiv 0.$$

Именно в этом месте и только в нем используется линейность $f_k\left(x;\xi^k\right)$. К сожалению, предложенный здесь способ рассуждения не позволяет хоть сколько-нибудь ослабить это условие.

Рассуждая далее также как в доказательстве теоремы 6.1.1, получим теорему 6.1.2.

**Замечание 6.1.7.** Как уже отмечалось во введении, теорема 6.1.2 во многом мотивирована работой [211]. Собственно, форма, в который мы представили алгоритм МЗС2-неадаптивный, выбрана именно такой (альтернативным вариантом было положить $\gamma_k\equiv 1$, $\beta_k\equiv M\sqrt{N/(2\ln n)}$ ), чтобы была максимальная близость к алгоритму работы [211].

**Замечание 6.1.8.** Идея рандомизации (искусственного введения случайности), положенная в основу описанных алгоритмов, чрезвычайно продуктивна: против нас играет, возможно, враждебная “Природа”, которая, зная историю игры, и наши текущие намерения старается нам “предложить вариант похуже”. С этим можно “бороться” за счет случайного независимого осуществления своих намерений на каждом шаге, с реализацией неизвестной “Природе”. За счет этой случайности мы переходим от анализа по худшему случаю (роль которого в онлайн оптимизации играет враждебная “Природа”) к анализу “в среднем”. Такая рандомизация, как будет отмечено в следующем пункте, дает возможность получать оценки, которые в детерминированном случае получить невозможно. Причем, если в онлайн постановке такая рандомизация прописывается в “правилах игры”, то применительно к задачам обычной оптимизации все это возникает совершенно естественным образом, как желание с большой вероятностью обезопасить себя от “самых худших случаев” детерминированной версии метода. Отметим, что речь идет о, так называемых, массовых задачах, т.е., исследуя тот или иной метод, мы точно не знаем какой конкретно объект поступит на вход, поэтому, чтобы гарантировано что-то иметь, мы исходим из худшего (наименее благоприятного для данного метода) случая входных данных. Описанный МЗС2 естественно также понимать как покомпонентный метод (стохастического) субградиентного спуска со случайным выбором компоненты. Происходит рандомизация при проектировании (в смысле расстояния Брэгмана) на единичный симплекс. А именно, если проектироваться в указанном выше смысле на единичный симплекс, то получится вектор, который можно проинтерпретировать как распределение вероятностей некоторой дискретной случайной величины, принимающей значения $1,\ldots,n$. Если выбрать вершину симплекса согласно этой дискретной случайной величине, и заменить проекцию этой вершиной, то получим случайную проекцию, математическое ожидание которой равно честной проекции. Как будет отмечено в следующем пункте, такой метод не только оптимален с точки зрения числа итераций, но и в некотором смысле с точки зрения затрат на выполнение одной итерации (см. также [46]).

### 6.1.4 Приложения метода зеркального спуска

В заключительном пункте мы постараемся продемонстрировать некоторые возможности и ограничения описанных в подразделах 6.1.2 и 6.1.3 методов. Мы не будем стремиться здесь к максимальной общности или рассмотрению всех основных приложений МЗС.

**Пример 6.1.1 (многорукие бандиты [164, 166, 250]).** Имеется $n$ различных ручек. Игра повторяется $N \gg 1$ раз (это число может быть заранее неизвестно). На каждом шаге $k$ мы должны выбрать ручку $i(k)$, которую “дергаем”. Дергание ручки приносит нам не-которые, вообще говоря, случайные потери $r_{i(k)}^k$ (считаем, для определенности, что всегда $r_{i(k)}^k \in [0,1]$), зависящие от номера шага, номера ручки и от того, какой стратегии мы придерживались до шага $k$ включительно. Наша стратегия на шаге $k$ описывается вектором распределения вероятностей $x^k \in S_n(1)$, согласно которому мы независимо ни от чего выбираем ручку, которую будем дергать. Все, чем мы располагаем на шаге $k$, это вектором

$$\left(\left(x^1, i(1), r_{i(1)}^1\right);...;\left(x^{k-1}, i(k-1), r_{i(k-1)}^{k-1}\right)\right).$$

Мы считаем, что потери на $k$-м шаге $r^k$ зависят от $x^k$ (но не от результата разыгрывания из распределения $x^k$), зависят от $\left(x^1,...,x^{k-1}\right)$ и результатов соответствующих разыгрываний, а также зависят от $\left(r^1,...,r^{k-1}\right)$. Целью является таким образом организовать процедуру дергания ручек, чтобы ожидаемые суммарные потери были бы минимальны. Введем функцию ($r^k$ и результат разыгрывания, согласно распределению вероятностей, заданному вектором $x$, – независимы; обе эти “случайности” мы обозначаем $\xi^k$)

$$f_k\left(x;\xi^k\right) = r_i^k \text{ с вероятностью } x_i,\ i = 1,...,n,$$

и её обобщенный (в смысле удовлетворения условию 6.1.2) стохастический градиент

$$\nabla_x f_k\left(x;\xi^k\right) = (\underbrace{0,...,r_i^k}_{i} / x_i,...,0)^T \text{ с вероятностью } x_i,\ i = 1,...,n.$$

Тогда выполнены условия 6.1.1 и 6.1.2. Однако имеется проблема: константа $M$ в условии 6.1.3 получается слишком большой (например, в 6.1.3.а $M = \sup\limits_{x\in S_n(1)} \sqrt{\sum\limits_{i=1}^n x_i^{-1}} = \infty$), то есть теорема 6.1.1 ничего дать не может в том виде, в котором она была нами приведена. Возникает желание “что-то подкрутить” в доказательстве теоремы, чтобы можно было ей воспользоваться. Уже по ходу самого доказательства мы отмечали неравенство (*), которое может оказаться довольно грубым в определенных ситуациях. Многорукие бандиты дают пример как раз такой ситуации. Более аккуратный анализ [164, 166, 225], использующий специфику данной задачи, позволяет оценить (см. доказательство теоремы 6.1.1)

$$\gamma_k \nabla_x f_k\left(x^k;\xi^k\right)^T \int\limits_0^1 \left(\nabla W_{\beta_k}\left(\tau y^k + (1-\tau) y^{k-1}\right) - \nabla W_{\beta_k}\left(y^{k-1}\right)\right) d\tau$$

точнее, что приводит в основной формуле (6.1.3) к замене слагаемых вида

$$\frac{\gamma_k^2}{2\alpha\beta_k}\left\|\nabla_x f_k\left(x^k;\xi^k\right)\right\|_\infty^2$$

на (в этом месте, для наглядности, мы намеренно несколько упрощаем и огрубляем)

$$\gamma_k^2 \frac{x_j^k\left(1-x_j^k\right)}{\alpha\beta_k}\left(\frac{r_j^k}{x_j^k}\right)^2,$$

где $j$ – номер ручки, выбранной алгоритмом на $k$-м шаге. В результате мы получаем, что теорема 6.1.1 остается верной с эффективной константой $M = \sqrt{2n}$. Таким образом, действуя согласно МЗС1, наши потери (псевдо регрет [166]) будут

$$\mathrm{O}\left(\sqrt{\frac{n\ln n}{N}}\right) \text{ – в среднем; } \mathrm{O}\left(\sqrt{\frac{n\ln(n/\sigma)}{N}}\right) \text{ – с вероятностью } \geq 1-\sigma,$$

что с точностью до логарифмического фактора соответствует нижним оценкам [164, 166, 250].

**Замечание 6.1.9.** Стоит обратить внимание, что если использовать более специальную прокс-структуру [164, 166], то для псевдо регрета можно получить оценки без логарифмического фактора $\ln n$ под корнем, что уже соответствует нижним оценкам. В частности, это обстоятельство означает, что выбирать “расстояние” Брэгмана для симплекса не всегда оптимально (но близко к оптимуму). Тем не менее, в последующих нескольких примерах мы убедимся, что для ряда других постановок, оценки, полученные с помощью прокс-структуры, порожденной расстоянием Брэгмана, – оптимальные. Кроме того, есть еще плата за избавление от фактора $\ln n$ под корнем – удорожание процедуры вычисления проекции на симплекс в смысле этой прокс-структуры. Другими словами, это ускорение оправдано только для онлайн постановок, в которых, как правило, стремятся минимизировать (псевдо) регрет, не сильно учитывая общую вычислительную трудоемкость.

Труднее обстоит дело с оценкой регрета (в среднем и вероятностей больших уклонений). Пример из лекции 6 [251] показывает $(n=2)$, что МЗС1 может давать регрет $\sim cN^{-1/4}$, что значительно хуже оценки псевдо регрета $\sim cN^{-1/2}$. По сути, речь идет о том, что написано в конце замечания 6.1.5. Здесь уже требуется некая игра “bias–variance trade off”: отказаться от несмещенности оценки градиента для уменьшения дисперсии этой оценки. Этот популярный трюк в математической статистике и машинном обучении позволяет с некоторыми оговорками распространить приведенные выше оценки псевдо регрета $\mathrm{O}\left(\sqrt{n\ln(n/\sigma)/N}\right)$ и на случай оценок регрета. Кое-что на эту тему применительно к многоруким бандитам можно найти в обзоре [166].

Интересно заметить, что результаты, описанные в примере 6.1.1 можно получить (с аналогичными оговорками) с помощью МЗС2 и теоремы 6.1.2. Для этого нужно взять

$$f_k\left(x;\xi^k\right) = \left\langle r^k, x\right\rangle$$

и её обобщенный (в смысле удовлетворения условию 6.1.2) стохастический градиент

$$\nabla_x f_k\left(x;\xi^k\right)=(\underbrace{0,...,r_i^k}_{i}/\,p_i,...,0)^T\text{, если } x=(\underbrace{0,...,1}_{i},...,0)^T,$$

$$\text{где } x=(\underbrace{0,...,1}_{i},...,0)^T \text{ с вероятностью } p_i,\ i=1,...,n,$$

здесь $\xi^k$ отражает только случайность, сидящую в $r^k$. Это определение стохастического градиента (в отличие от рассмотренного выше) явно учитывает распределение вероятностей $p$, из которого генерируется номер единственной ненулевой (единичной) компоненты вектора $x$. Также как и раньше для выполнения условия 2 необходимо предполагать независимость $\xi^k$ и процедуры разыгрывания единичный компоненты вектора $x$ согласно распределению $p$.

**Пример 6.1.2 (взвешивание экспертных решений, линейные потери [27, 250]).** Рассмотрим задачу взвешивание экспертных решений, следуя [27, 250]. Имеется $n$ различных Экспертов. Каждый Эксперт играет на рынке. Игра повторяется $N \gg 1$ раз (это число может быть заранее неизвестно). Пусть $l_i^k$ – проигрыш Эксперта $i$ на шаге $k$ ($\left|l_i^k\right| \le M$ ). На каждом шаге $k$ мы распределяем один доллар между Экспертами, согласно вектору $x^k \in S_n\left(1\right)$. Потери, которые мы при этом несем, рассчитываются по потерям экспертов $\left\langle l^k, x^k\right\rangle$. Целью является таким образом организовать процедуру распределения доллара на каждом шаге, чтобы наши суммарные потери были бы минимальны. Допускается, что потери экспертов $l^k$ могут зависеть еще и от текущего хода $x^k$. Легко проверить, что для данной постановки применима теорема 6.1.1 в детерминированном варианте с функциями

$$f_k\left(x;\xi^k\right) \equiv f_k\left(x\right)=\left\langle l^k, x\right\rangle.$$

При этом оценка, даваемая теоремой 1,

$$\mathrm{O}\left(M\sqrt{\frac{\ln n}{N}}\right)$$

– оптимальна для данного класса задач [27, 250].

**Пример 3 (взвешивание экспертных решений, выпуклые потери [27, 250]).** В условиях предыдущего примера предположим, что на $k$-м шаге $i$-й эксперт использует стратегию $\zeta_i^k \in \Delta$ (множество $\Delta$ – выпуклое), дающую потери $\lambda\left(\omega^k, \zeta_i^k\right)$, где $\omega^k$ – “ход”, возможно, враждебной “Природы”, знающей, в том числе, и нашу текущую стратегию. Функция $\lambda\left(\,\cdot\,\right)$ – выпуклая по второму аргументу и $\left|\lambda\left(\,\cdot\,\right)\right| \le M$. На каждом шаге мы должны выбирать свою стратегию

$$x \overset{def}{=} \sum_{i=1}^{n} x_i \cdot \zeta_i^k \in \Delta,$$

дающую потери $\lambda\left(\omega^k, x\right)$ так, чтобы наши суммарные потери были минимальны. Для данной постановки также применима теорема в детерминированном варианте с

$$f_k\left(x;\xi^k\right) \equiv f_k\left(x\right) = \sum_{i=1}^{n} x_i \lambda\left(\omega^k, \zeta_i^k\right) \geq \lambda\left(\omega^k, x\right).$$

Чтобы применить теорему осталось заметить, что функция $\lambda\left(\omega^k, \zeta\right)$ – выпуклая по $\zeta$ для любого $\omega^k$, поэтому

$$\sum_{k=1}^{N} \lambda\left(\omega^k, x^k\right) - \min_{i=1,\ldots,n} \sum_{k=1}^{N} \lambda\left(\omega^k, \zeta_i^k\right) \leq \sum_{k=1}^{N} f_k\left(x^k\right) - \min_{x \in S_n(1)} \sum_{k=1}^{N} f_k\left(x\right).$$

При этом оценка, даваемая теоремой 6.1.1,

$$\mathrm{O}\left(M\sqrt{\frac{\ln n}{N}}\right)$$

– также оптимальна для данного класса задач [27, 250].

Полезно, на наш взгляд, будет здесь привести другой способ (более типичный для данного класса приложений) получения аналогичного результата, не связанный на прямую с МЗС схемой вывода, но фактически, приводящий к точно такому же алгоритму. Этот способ также весьма популярен в машинном обучении, теории игр, теории алгоритмов [27, 139, 250].[1]

Введем обозначение $L_i^N = \sum_{k=1}^{N} \lambda\left(\omega^k, \zeta_i^k\right)$, $\tilde{L}^N = \sum_{k=1}^{N} \lambda\left(\omega^k, x^k\right)$, по определению считаем $L_i^0 \equiv 0$. Рассмотрим

$$W_\beta\left(\left\{-L_i^N\right\}_{i=1}^n\right) = \beta \ln\left(\frac{1}{n}\sum_{i=1}^{n} \exp\left(-L_i^N/\beta\right)\right) \geq -\min_{i=1,\ldots,n} L_i^N - \beta \ln n.$$

С другой стороны, вводя дискретную случайную величину $z^k$, имеющую (независящее ни от чего) распределение $x^k$ (рассчитанное также как и раньше исходя из МЗС1, примененного к набору функций $\left\{f_k\left(x\right)\right\}_{k=1}^N$, определенных выше), можно заметить, что

$$W_\beta\left(\left\{-L_i^N\right\}_{i=1}^n\right) = \sum_{k=1}^{N}\left(W_\beta\left(\left\{-L_i^k\right\}_{i=1}^n\right) - W_\beta\left(\left\{-L_i^{k-1}\right\}_{i=1}^n\right)\right) = \beta \sum_{k=1}^{N} \ln\left(E_z\left(e^{-\lambda\left(\omega^k, z^k\right)/\beta}\right)\right).$$

---

[1] Максимум из независимых случайных величин, который сложно исследовать, заменяется (с хорошей точностью, контролируемой малостью параметра $\beta$) логарифмом от суммы экспонент от этих независимых случайных величин. А сумму независимых случайных величин (их экспонент) исследовать уже на много проще. В оптимизации эту процедуру называют сглаживанием [271].

Используя далее неравенство Хефдинга (для с.в. $X \in [-M, M]$), см. [159]

$$\ln\left(E_X\left(e^{sX}\right)\right) \le sE_X\left(X\right) + s^2\frac{M^2}{2},$$

Получим

$$W_\beta\left(\left\{-L_i^N\right\}_{i=1}^n\right) \le -\tilde{L}^N + \left(2\beta\right)^{-1}M^2N.$$

Таким образом,

$$\tilde{L}^N \le \min_{i=1,\ldots,n} L_i^N + \beta\ln n + \left(2\beta\right)^{-1}M^2N.$$

Минимизация правой части по $\beta > 0$ приводит нас к уже известному ответу. Аналогичные, но чуть более тонкие рассуждения, позволяют избавиться от зависимости $\beta$ от $N$, то есть сделать алгоритм адаптивным.

**Пример 6.1.4 (взвешивание экспертных решений, невыпуклые потери [27, 250]).** Предположим, что в условиях примера 6.1.3 мы не можем гарантировать выпуклость $\lambda(\cdot)$ – по второму аргументу. Тогда мы выбираем стратегию – распределение вероятностей на множестве стратегий Экспертов, и разыгрываем случайную величину согласно этому распределению вероятностей. Другими словами мы просто пользуемся МЗС2 с $f_k\left(x;\xi^k\right) \equiv f_k\left(x\right) = \sum_{i=1}^n x_i\lambda\left(\omega^k, \zeta_i^k\right)$, применимость которого обосновывается теоремой 6.1.2, с оценками

$$\mathrm{O}\left(M\sqrt{\frac{\ln n}{N}}\right) \text{ – в среднем; } \mathrm{O}\left(M\sqrt{\frac{\ln(n/\sigma)}{N}}\right) \text{ – с вероятностью } \ge 1-\sigma,$$

– оптимальными для данного класса задач [27, 250]. Ключевая разница в примерах 6.1.1 и 6.1.4, “сто́ящая” $\sim\sqrt{n}$ в оценке $M$ для многоруких бандитов (пример 6.1.1), заключается в том, что в многоруких бандитах мы имеем только свою историю дергания ручек (нам не известно, какие бы потери нам принесли другие ручки, кабы мы их выбрали), а в постановке взвешивания экспертных решений это все известно, и называется потерями экспертов.

Как будет видно из следующего примера описанный только что подход вполне успешно работает (дает не улучшаемые результаты) и в случае выпуклой по второму аргументу функции $\lambda(\cdot)$.

**Пример 6.1.5 (антагонистические матричные игры [27, 46, 211, 250], см. также разделы 4.1 и 4.2 главы 4).** Пусть есть два игрока А и Б. Задана матрица игры $A = \left\|a_{ij}\right\|$, где $|a_{ij}| \le M$, $a_{ij}$ – выигрыш игрока А (проигрыш игрока Б) в случае когда игрок А выбрал стратегию $i$, а игрок Б стратегию $j$. Отождествим себя с игроком Б. И предположим, что игра повторяется $N \gg 1$ раз (это число может быть заранее неизвестно). Мы находимся в условиях примера 6.1.4 с $\lambda\left(\omega^k, \zeta_j^k\right) = \sum_{i=1}^{n} \omega_i^k a_{ij}$, то есть

$$f_k(x) = \left\langle \omega^k, Ax \right\rangle,\ x \in S_n(1),$$

где $\omega^k$ – вектор (вообще говоря, зависящий от всей истории игры до текущего момента включительно, в частности, как-то зависящий и от текущей стратегии (не хода) игрока Б, заданной распределением вероятностей (результат текущего разыгрывания (ход Б) игроку А не известен)) со всеми компонентами равными 0, кроме одной компоненты, соответствующей ходу А на шаге $k$, равной 1. Хотя функция $f_k(x)$ определена на единичном симплексе, по “правилам игры” вектор $x^k$ имеет ровно одну единичную компоненту, соответствующую ходу Б на шаге $k$, остальные компоненты равны нулю. Обозначим цену игры

$$C = \max_{\omega \in S_n(1)} \min_{x \in S_n(1)} \langle \omega, Ax \rangle = \min_{x \in S_n(1)} \max_{\omega \in S_n(1)} \langle \omega, Ax \rangle. \text{ (теорема фон Неймана о минимаксе)}$$

**Замечание 6.1.10.** Отметим, что с помощью онлайн оптимизации и экспоненциального взвешивания можно похожим образом проинтерпретировать и вариант теоремы о минимаксе для векторнозначной функции выигрыша – теорему Блэкуэлла о достижимости [27, 250], которая используется, например, при построении калибруемых предсказаний.

Пару векторов $(\omega, x)$, доставляющих решение этой минимаксной задачи (т.е. седловую точку), назовем равновесием Нэша. По определению (это неравенство восходит к Ханнану [27, 250])

$$\min_{x \in S_n(1)} \frac{1}{N} \sum_{k=1}^{N} f_k(x) \le C.$$

Тогда, если мы (игрок Б) будем придерживаться рандомизированной стратегии МЗС2, выбирая $\left\{x^k\right\}$, то по теореме 6.1.2 с вероятностью $\ge 1-\sigma$ (в случае когда $N$ заранее известно оценку можно уточнить)

$$\frac{1}{N} \sum_{k=1}^{N} f_k\left(x^k\right) - \min_{x \in S_n(1)} \frac{1}{N} \sum_{k=1}^{N} f_k(x) \le \frac{2M}{\sqrt{N}} \left(\sqrt{\ln n} + \sqrt{2\ln\left(\sigma^{-1}\right)}\right),$$

т.е. с вероятностью $\ge 1-\sigma$ наши потери ограничены

$$\frac{1}{N} \sum_{k=1}^{N} f_k\left(x^k\right) \le C + \frac{2M}{\sqrt{N}} \left(\sqrt{\ln n} + \sqrt{2\ln\left(\sigma^{-1}\right)}\right).$$

Самый плохой для нас случай (с точки зрения такой оценки) – это когда игрок А тоже “знает” теорему 6.1.2, и действует (выбирая $\{\omega^k\}$) согласно МЗС2 (точнее версии МЗС2 для максимизации вогнутых функций на симплексе). Очевидно, что если и А и Б будут придерживаться МЗС2, то они сойдутся к равновесию Нэша (седловой точке), причем чрезвычайно быстро [46]:

$$\mathrm{O}\left(\frac{M\left(\ln n+2\ln\left(\sigma^{-1}\right)\right)}{\varepsilon^2}\right) \text{ – итераций;}$$

$$\mathrm{O}\left(n+M\frac{s\ln n\left(\ln n+\ln\left(\sigma^{-1}\right)\right)}{\varepsilon^2}\right) \text{ – общее число арифметических операций,}$$

где $s \le n$ – среднее число элементов в строках и столбцах матрицы $A$. Отсюда видно, что если $\varepsilon$ (зазор двойственности [225, 269]) – не очень малое, то может случиться, что общее число арифметических операций будет много меньше числа элементов матрицы $A$ отличных от нуля, в то время как любой детерминированный способ поиска равновесия Нэша потребовал бы прочтения как минимум половины элементов матрицы $A$ [46, 211]. Другой способ содержательной интерпретации описанного метода базируется на прямо-двойственности МЗС [269].

**Замечание 6.1.11.** Если использовать методы работ [260, 271] то можно получить зависимость сложности от $\varepsilon$ вида $\varepsilon^{-1}$, но при этом число операций увеличится не менее чем в $n$ раз. Поскольку для таких задач вполне естественным является соотношение $\varepsilon \gg n^{-1}$, то выгоднее использовать описанный в этом примере алгоритм.

Есть основания полагать, что описанный здесь метод оптимален не только с точки зрения числа итераций, но и с точки зрения “стоимости” шага. Особенно ярко это проявляется, когда $s \ll n$ [46]. Нам кажется, что описанная в этом примере методология может оказаться полезной в huge-scale optimization [264, 273] (см. также главу 4 диссертации). Отметим в связи с вышесказанным другой пример задачи выпуклой оптимизации, когда удается получить ответ (с требуемой точностью), с большим запасом не просматривая весь объем имеющихся данных (общая проблема здесь: как понять, что просматривать, а что нет): метод наименьших квадратов с разреженной структурой [165]. Нам представляется, что именно эта ветвь более общего и бурно развивающегося в последнее время направления huge-scale оптимизации является сейчас одной из наиболее интересных, как с практической, так и с теоретической точки зрения, и основные результаты здесь еще впереди. Уточним, что речь идет о задачах, приходящих из: машинного обучения (в частности, распознавания изображений), моделирования различных сетей огромных размеров типа сети Интернет или транспортных сетей (глава 1), биоинформатики, численных методов (проектирование конструкций методом конечных элементов) и ряда других приложений. Все эти задачи отличают колоссальные размеры. Скажем, для задачи ранжирования web-страниц необходимо решать вспомогательную оптимизационную задачу в пространстве, размерность которого больше миллиарда [46] (глава 4). Но помимо размеров их отличает некоторые релаксированные требования к решению. Например, нет никакой необходимости ранжировать абсолютно все web-страницы по заданному запросу и делать это очень точно. Достаточно, чтобы качественно выдавались только первые сто наиболее значимых (больших) компонент ранжирующего вектора, причем важны не столько сами значения этих компонент, сколько их порядок. Также эти задачи довольно специальные, то есть использовать концепцию черного ящика [91] для оценки числа требуемых итераций, как правило, не представляется возможным. Более того, оценки имеются, в основном, только на число шагов (итераций). Поскольку общее время работы алгоритма определяется произведением числа итераций на стоимость одной итерации, то возникает игра между числом итераций и стоимостью итерации. Описанный в этом примере метод как раз играет в эту игру. А именно, он увеличивает за счет рандомизации число итераций в $\ln\left(\sigma^{-1}\right)$ раз, при этом итерация становится в $n/\ln n$ раз дешевле. Наконец, важной составляющей многих задач является разреженная структура. Тут имеется два варианта: разреженная структура решения и данных. В первом случае (частично уже затронутым выше на примере ранжирования web-страниц) часто речь идет о подмене исходной задачи, вычислительно более привлекательной задачей, по решению которой можно получить приближенное решение исходной. Пожалуй, наиболее ярким примером здесь является сжатие измерений (см. [115], и цитированную там литературу). В случае разреженных данных представляется перспективным использование специальных покомпонентных спусков и исследование вычислительных особенностей пересчета различных классов функций многих переменных в случае изменения лишь не большого числа их аргументов [264, 273] (см. также главу 4 и раздел 5.1 главы 5). В заключение отметим, что важными составляющими анализа эффективности методов (в виду специфики описываемых задач) является исследование возможности распараллеливания [199] (в рассматриваемом нами примере 6.1.5 это возможно сделать [46, 211]) и вероятностный анализ в среднем [115] (или для почти всех входов). В отличие от Computer Science [189, 256], в численных методах выпуклой оптимизации такой анализ можно встретить не часто (все же кое-что есть, например, вероятностный анализ симплекс-метода Спилманом и Тэнгом [318]). Хотя уже сейчас (в связи с бурным развити-ем идей концентрации меры [189, 244]) становится все более и более ясно, что в пространствах огромных размеров такой анализ необходим, и может многое дать [115].

Отметим в заключение, что к рассмотренной в примере 6.1.5 антагонистической матричной игре (с $M = 1$) сводится Google problem: задача поиска вектора Фробениуса–Перрона стохастической матрицы $P$ огромных размеров [46] (см. также разделы 1, 2 главы 4 и приложение в конце диссертации): $A = P^T - I$. Разработанный в [46] и описанный в примере 5 подход позволяет учитывать разреженную структуру матрицы $P$, заменяя в алгоритме из [88] $n$, которое с точностью до логарифмических факторов входит линейно в оценку общей трудоемкости метода (то есть с учетом затрат на каждой итерации) на $s \ll n$.

## 6.2 Стохастическая онлайн оптимизация. Одноточечные и двухточечные нелинейные многорукие бандиты. Выпуклый и сильно выпуклый случаи

### 6.2.1 Введение

Данный раздел представляет собой попытку перенесения результатов статьи [52] (см. также раздел 5.2 главы 5) на онлайн контекст [54, 130, 164, 166, 216, 250, 299, 306, 320]. А именно, следуя работе [52] рассматривается постановка задачи выпуклой стохастической онлайн оптимизации, в которой на каждом шаге (итерации) вместо градиента можно получать только реализацию значения соответствующей этому шагу функции. При этом допускается, что эта реализация доступна с шумом уровня $\delta$, вообще говоря, не случайной природы. Рассматривается две возможности: на одном шаге (при одной реализации) получать зашумленное значение в одной точке и в двух точках. В первом случае говорят, что рассматривается задача о нелинейных многоруких бандитах (иногда добавляя, одноточечных) [166]. Во втором случае говорят о нелинейных многоруких двухточечных бандитах [166]. Принципиальная разница есть именно при таком переходе [130, 166]. Последующее увеличение числа точек не меняет принципиально картину, соответствующую двум точкам [190].

Основная идея заключается в специальном сглаживании исходной постановки задачи, и использовании метода зеркального спуска [52, 54, 91, 166, 201]. Оригинальной составляющей здесь, в частности, является предложенное в данном разделе обобщение этой конструкции на случай наличия шумов. Обратим внимание на условие 6.2.1 в подразделе 6.2.2 (следует сравнить, например, с [44, 52, 320]. Это условие позволило с одной стороны изящно распространить известные оценки на случай, когда есть шумы, см. формулы (6.2.2), (6.2.3) подраздела 6.2.2, а с другой стороны это условие хорошо подходит под специфику рассматриваемой в разделе 6.2 постановки (можем получать только зашумленные реализации значений функций), что демонстрируется в подразделе 6.2.3. Основным результатом этого раздела 6.2 является теорема 6.2.1 раздела 6.2.3, в которой результаты статьи [52] переносятся на онлайн контекст.

Во избежание большого количества громоздких выражений, мы опустили часть (наиболее очевидных, но громоздких) выкладок, подробно описав, как они могут быть сделаны. Также в изложении мы не стремились к общности. В частности, для большинства оценок данного раздела можно не только выписать точные константы в оценках сходимости в среднем (для этого вполне достаточно написанного в данном разделе), но и получить оценки вероятностей больших уклонений. Также можно накладывать более общие требования на классы изучаемых семейств функций, делая константы, характеризующие семейство, не универсальными (одинаковыми для всех шагов), а зависящими от номера шага [130, 259].

Полученные оценки, с учетом известных нижних оценок [128,166, 167, 216, 250], позволяют говорить о том, что в данном разделе предложены достаточно эффективные методы, доминирующие в ряде случаев существующие сейчас алгоритмы.

### 6.2.2 Метод зеркального спуска для задач стохастической онлайн оптимизации с неточным оракулом

Сформулируем основную задачу стохастической онлайн оптимизации с неточным оракулом. Требуется подобрать последовательность $\{x^k\} \in Q$ ($Q$ – выпуклое множество) так, чтобы минимизировать псевдо регрет [54, 130, 164, 166, 216, 250, 299, 306, 320] (см. также раздел 1 этой главы):

$$\mathrm{Regret}_N\left(\{f_k(\cdot)\},\{x^k\}\right) = \frac{1}{N}\sum_{k=1}^{N} f_k\left(x^k\right) - \min_{x\in Q}\frac{1}{N}\sum_{k=1}^{N} f_k(x). \qquad (6.2.1)$$

на основе доступной информации

$$\left\{\nabla_x \tilde{f}_1\left(x^1,\xi^1\right);\ldots;\nabla_x \tilde{f}_{k-1}\left(x^{k-1},\xi^{k-1}\right)\right\}$$

при расчете $x^k$. Причем выполнено **условие**[2]

6.2.1 для любых $N \in \mathbb{N}$ ($\Xi^{k-1}$ – сигма алгебра, порожденная $\xi^1$, …, $\xi^{k-1}$)

$$\sup_{\left\{x^k = x^k\left(\xi^1,\ldots,\xi^{k-1}\right)\right\}_{k=1}^{N}} E\left[\frac{1}{N}\sum_{k=1}^{N}\left\langle E_{\xi^k}\left[\nabla_x f_k\left(x^k,\xi^k\right) - \nabla_x \tilde{f}_k\left(x^k,\xi^k\right)\middle|\Xi^{k-1}\right], x^k - x_*\right\rangle\right] \le \sigma,$$

где $x_*$ – решение задачи

$$\frac{1}{N}\sum_{k=1}^{N} f_k(x) \to \min_{x\in Q},$$

$$E_{\xi^k}\left[\nabla_x f_k\left(x,\xi^k\right)\right] = \nabla f_k(x).$$

Здесь случайные величины $\{\xi^k\}$ могут считаться независимыми одинаково распределенными. Онлайновость постановки задачи допускает, что на каждом шаге $k$ функция $f_k(\cdot)$ может выбираться из рассматриваемого класса функций враждебно по отношению к ис-

---

[2] В частности, если

$$\left\|\nabla_x \tilde{f}_k\left(x^k,\xi^k\right) - \nabla_x f_k\left(x^k,\xi^k\right)\right\|_* \le \delta, \quad \max_{x,y\in Q}\|x-y\| \le R,$$

то $\sigma \le \delta R$.

пользуемому нами методу генерации последовательности $\{x^k\}$. В частности, $f_k(\cdot)$ может зависеть от

$$\left\{x^1,\xi^1, f_1(\cdot);...;x^{k-1},\xi^{k-1}, f_{k-1}(\cdot);x^k\right\}.$$

Относительно класса функций, из которого выбираются $\{f_k(\cdot)\}$, в данном разделе в зависимости от контекста будем предполагать выполненными следующие **условия**:

6.2.2 $\{f_k(\cdot)\}$ – выпуклые функции (считаем, что это условие имеет место всегда);

6.2.3 $\{f_k(\cdot)\}$ – $\gamma_2$-сильно выпуклые функции в $l_2$;

6.2.4 для любых $k=1,...,N$, $x\in Q$

$$E_\xi\left[\left\|\nabla_x \tilde{f}_k(x,\xi)\right\|_*^2\right]\le M^2.$$

Выше (и далее в разделе) используется стандартная терминология онлайн оптимизации (см., например, [166, 216, 250, 306]). Однако в отечественной литературе на данный момент имеется определенный дефицит работ по этой (достаточно популярной на западе) тематике. В связи с этим было решено “разбавить” данный раздел несколькими простыми примерами (предложенными А.А. Лагуновской), которые позволят лучше прочувствовать смысл используемых понятий.

**Пример 6.2.1 (взвешивание экспертных решений, линейные потери, см. раздел 6.1 этой главы).** Рассмотрим задачу взвешивания экспертных решений, следуя [250]. Имеется $n$ различных Экспертов. Каждый Эксперт “играет” на рынке. Игра повторяется $N\gg 1$ раз. Пусть $l_i^k$ – проигрыш (выигрыш со знаком минус) Эксперта $i$ на шаге $k$ ($\left|l_i^k\right|\le M$). На каждом шаге $k$ распределяется один доллар между Экспертами, согласно вектору

$$x^k\in Q=S_n(1)=\left\{x\ge 0:\ \sum_{i=1}^n x_i=1\right\}.$$

Потери, которые при этом несем, рассчитываются по потерям экспертов $f_k(x)=\left\langle l^k,x\right\rangle$. Целью является таким образом организовать процедуру распределения доллара на каждом шаге, чтобы суммарные потери (за $N$ шагов) были бы минимальны. Допускается, что потери экспертов $l^k$ могут зависеть еще и от текущего хода $x^k$. Установленные далее в этом разделе результаты (формула (6.2.2) с $R^2=\ln n$) позволяют утверждать, что если на каждом шаге можно наблюдать лишь зашумленные проигрыши Экпертов

$$\frac{\partial \tilde{f}_k\left(x^k,\xi^k\right)}{\partial x_i}=l_i^k+\xi_i^k+\delta_i^k,$$

где $\left\{\xi_i^k\right\}$ – независимые одинаково распределенные случайные величины $\xi_i^k \in N\left(0,1\right)$, $\left|\delta_i^k\right| \le \delta$, то существует такой способ действий $x^k\left(\nabla_x \tilde{f}_1\left(x^1,\xi^1\right),\dots,\nabla_x \tilde{f}_{k-1}\left(x^{k-1},\xi^{k-1}\right)\right)$ (метод зеркального спуска с $\| \ \| = \| \ \|_1$, $d\left(x\right)=\ln n+\sum_{i=1}^{n} x_i \ln x_i$, см. ниже), который позволяет с вероятностью не менее 0.999 после $N$ шагов проиграть лучшему (на этом периоде $k=1,\dots,N$) Эксперту не более $\mathrm{O}\left(\left(M+1\right)\sqrt{N\ln n}+\delta N\right)$ долларов, что означает (см. формулу (6.2.1))

$$\mathrm{Regret}_N\left(\left\{f_k\left(\cdot\right)\right\},\left\{x^k\right\}\right)=\mathrm{O}\left(\left(M+1\right)\sqrt{\frac{\ln n}{N}}+\delta\right).$$

При $\delta=0$ эта оценка оптимальна для данного класса задач [250]. □

Опишем метод зеркального спуска для решения задачи (6.2.1) (здесь можно следовать огромному числу литературных источников, мы в основном будем следовать работам [135, 259]). Введем норму $\| \ \|$ в прямом пространстве (сопряженную норму будем обозначать $\| \ \|_*$) и прокс-функцию $d\left(x\right)$ сильно выпуклую относительно этой нормы, с константой сильной выпуклости $\ge 1$. Выберем точку старта

$$x^1=\arg\min_{x\in Q} d\left(x\right),$$

считаем, что $d\left(x^1\right)=0$, $\nabla d\left(x^1\right)=0$.

Введем брэгмановское “расстояние”

$$V_x\left(y\right)=d\left(y\right)-d\left(x\right)-\left\langle \nabla d\left(x\right),y-x\right\rangle.$$

Везде в дальнейшем будем считать, что

$$d\left(x\right)=V_{x^1}\left(x\right)\le R^2 \text{ для всех } x\in Q.$$

Определим оператор “проектирования” согласно этому расстоянию

$$\mathrm{Mirr}_{x^k}\left(g\right)=\arg\min_{y\in Q}\left\{\left\langle g,y-x^k\right\rangle+V_{x^k}\left(y\right)\right\}.$$

Метод зеркального спуска (МЗС) для задачи (1) будет иметь вид, см., например, [135]

$$x^{k+1}=\mathrm{Mirr}_{x^k}\left(\alpha_k\nabla_x \tilde{f}_k\left(x^k,\xi^k\right)\right),\ k=1,\dots,N.$$

Тогда при выполнении условии (6.2.2) для любого $u\in Q$, $k=1,\dots,N$ имеет место неравенство, см., например, [135]

$$\alpha_k \left\langle \nabla_x \tilde{f}_k\left(x^k,\xi^k\right), x^k - u\right\rangle \le \frac{\alpha_k^2}{2}\left\|\nabla_x \tilde{f}_k\left(x^k,\xi^k\right)\right\|_*^2 + V_{x^k}(u) - V_{x^{k+1}}(u).$$

Это неравенство несложно получить в случае евклидовой прокс-структуры $d(x) = \|x\|_2^2/2$ [269] (в этом случае МЗС для задачи (6.2.1) есть просто вариант обычного метода проекции градиента). Разделим сначала выписанное неравенство на $\alpha_k$ и возьмем условное математическое ожидание $E_{\xi^{k+1}}\left[\,\cdot\,|\Xi^k\right]$, затем просуммируем то, что получится по $k = 1,\dots,N$, используя условие 6.2.1. Затем возьмем от того, что получилось при суммировании, полное математическое ожидание, учитывая условие 6.2.4. В итоге, выбирая $u = x_*$, получим при условиях 6.2.1, 6.2.2, 6.2.4, $\alpha_k \equiv \alpha$ [54]

$$N\cdot E\left[\mathrm{Regret}_N\left(\{f_k(\,\cdot\,)\},\{x^k\}\right)\right] \le \frac{V_{x^1}(x_*)}{\alpha} - \frac{E\left[V_{x^{N+1}}(x_*)\right]}{\alpha} + \left(\frac{1}{2}M^2\alpha + \sigma\right)N \le$$

$$\le \frac{R^2}{\alpha} + \left(\frac{1}{2}M^2\alpha + \sigma\right)N,$$

выбирая[3]

$$\alpha = \frac{R}{M}\sqrt{\frac{2}{N}},$$

получим

$$E\left[\mathrm{Regret}_N\left(\{f_k(\,\cdot\,)\},\{x^k\}\right)\right] \le MR\sqrt{\frac{2}{N}} + \sigma; \qquad (6.2.2)$$

при условиях[4] 6.2.1, 6.2.3, 6.2.4, $\alpha_k \equiv (\gamma_2 k)^{-1}$, $\|\ \| = \|\ \|_2$ [216]

$$E\left[\mathrm{Regret}_N\left(\{f_k(\,\cdot\,)\},\{x^k\}\right)\right] \le \frac{M^2}{2\gamma_2 N}(1 + \ln N) + \sigma. \qquad (6.2.3)$$

Оценки (6.2.2), (6.2.3) являются неулучшаемыми с точностью до мультипликативного числового множителя. Причем верно это и для детерминированных (не стохастических)

---

[3] Можно получить и адаптивный вариант приводимой далее оценки, для этого потребуется использовать метод двойственных усреднений [54], [135], [269].

[4] Отметим, что при условии 6.2.2, мы еще используем неравенство

$$f_k\left(x^k\right) - f(x_*) \le \left\langle \nabla f_k\left(x^k\right), x^k - x_*\right\rangle$$

при преобразовании левой части неравенства в псевдо регрет, а при условии 6.2.3 более точное неравенство

$$2\left(f_k\left(x^k\right) - f(x_*)\right) \le 2\left\langle \nabla f_k\left(x^k\right), x^k - x_*\right\rangle - \gamma_2\left\|x^k - x_*\right\|_2^2.$$

постановок, в которых нет шумов ($\sigma = 0$), в случае оценки (6.2.2) при этом можно ограничиться классом линейных функций [250].

**Пример 6.2.2.** Пусть $Q = B_p^n(1)$ – единичный шар в $l_p$ норме. Относительно оптимального выбора нормы и прокс-структуры можно заметить следующее: если $p \geq 2$, то в качестве нормы $\| \ \|$ оптимально выбирать $l_2$ норму и евклидову прокс-структуру. Определим $q$ из $1/p + 1/q = 1$. Пусть $1 \leq p \leq 2$, тогда $q \geq 2$. Если при этом $q = o(\ln n)$, то оптимально выбирать $l_p$ норму, а прокс-структуру задавать прокс-функцией

$$d(x) = \frac{1}{2(p-1)} \|x\|_p^2 .$$

Во всех этих случаях

$$R^2 = \max_{x \in Q} d(x) = \mathrm{O}(1).$$

Для $q \geq \Omega(\ln n)$, выберем $l_a$ норму, где

$$a = \frac{2 \ln n}{2 \ln n - 1},$$

а прокс-структуру будем задавать прокс-функцией

$$d(x) = \frac{1}{2(a-1)} \|x\|_a^2 .$$

В этом случае $R^2 = \mathrm{O}(\ln n)$. Детали см., например, в работах [128, 259]. □

### 6.2.3 Одноточечные и многоточечные нелинейные многорукие бандиты

Везде в этом подразделе будем считать, что все функции $f_k(x)$ и реализации $f_k(x, \xi)$ определены в $Q_{\mu_0}$ – $\mu_0$-окрестности множества $Q$, и удовлетворяют соответствующим условиям из подраздела 6.2.2 именно в $Q_{\mu_0}$.

Пусть требуется подобрать последовательность $\{x^k\} \in Q$ так, чтобы минимизировать псевдо регрет (6.2.1) на основе доступной информации ($m = 1, 2$)

$$\left\{ \left\{ \tilde{f}_1\left(x_i^1, \xi^1\right) \right\}_{i=1}^m ; \ldots ; \left\{ \tilde{f}_{k-1}\left(x_i^{k-1}, \xi^{k-1}\right) \right\}_{i=1}^m \right\}$$

при расчете $x^k$. Будем предполагать, что имеет место следующее **условие**

6.2.5 для любых $k \in \mathbb{N}$, $i = 1, \ldots, m$, $x_i^k \in Q_{\mu_0}$

$$\left| \tilde{f}_k\left(x_i^k, \xi^k\right) - f_k\left(x_i^k, \xi^k\right) \right| \leq \delta ,$$

$$E_{\xi^k}\left[ f_k\left(x_i^k,\xi^k\right)\right] = f_k\left(x_i^k\right),$$

$$E_{\xi^k}\left[ \tilde{f}_k\left(x_i^k,\xi^k\right)^2\right] \le B^2;$$

и, в зависимости от контекста, **условия**

6.2.6 для любых $k = 1,...,N$, $x, y \in Q_{\mu_0}$ (далее, как правило, это условие будет использоваться при $r = 2$, исключение сделано в Таблице 6.2.2)

$$\left| f_k(x,\xi) - f_k(y,\xi)\right| \le M_r(\xi)\|x-y\|_r,\ M_r = \sqrt{E_\xi\left[ M_r(\xi)^2\right]} < \infty;$$

6.2.7 для любых $k = 1,...,N$, $x, y \in Q_{\mu_0}$

$$\left\|\nabla_x f_k(x,\xi) - \nabla_x f_k(y,\xi)\right\|_2 \le L_2(\xi)\|x-y\|_2,\ L_2 = \sqrt{E_\xi\left[ L_2(\xi)^2\right]} < \infty.$$

Введем аналоги $\nabla_x \tilde{f}_k(x,\xi)$ из подраздела 6.2. 2 ( $\mu \le \mu_0$ )

$$\nabla_x \tilde{f}_k(x;e,\xi) := \frac{n}{\mu}\tilde{f}_k(x+\mu e,\xi)e \quad \text{(при } m=1\text{)},$$

$$\nabla_x \tilde{f}_k(x;e,\xi) := \frac{n}{\mu}\left(\tilde{f}_k(x+\mu e,\xi) - \tilde{f}_k(x,\xi)\right)e \quad \text{(при } m=2\text{)},$$

где $e \in RS_2^n(1)$, т.е. случайный вектор $e$ равномерно распределен на сфере радиуса 1 в $l_2$. Считаем, что разыгрывание $e$ происходит независимо ни от чего. Аналогично можно определить незашумленную оценку стохастического градиента $\nabla_x f_k(x;e,\xi)$, убрав в правой части тильды (волны).

Онлайновость постановки задачи допускает, что на каждом шаге $k$ функция $f_k(\cdot)$ может выбираться из рассматриваемого класса функций враждебно по отношению к используемому нами методу генерации последовательности $\left\{x^k\right\}$. В частности, $f_k(\cdot)$ может зависеть от

$$\left\{x^1,\xi^1, f_1(\cdot);...;x^{k-1},\xi^{k-1}, f_{k-1}(\cdot)\right\}.$$

Более того, при выборе $f_k(\cdot)$ считается полностью известным наша стратегия. Подчеркнем, что поскольку стратегия рандомизированная, то речь идет об описании этой стратегии, а не о реализации. Это означает, что тому, кто подбирает $f_k(\cdot)$, известно, что $e \in RS_2^n(1)$, но не известно как именно мы его разыграем. Это важная оговорка, если допускать, как и в подразделе 6.2.2, что на каждом шаге $k$ реализация $e_k \in RS_2^n(1)$ становится известной тому, кто враждебно подбирает $f_k(\cdot)$, то нельзя получить оценку псевдо регрета лучше чем $\Omega(N)$ [130]. Причины этого, связаны с введением рандомизации, и на более простой задаче (линейные одноточечные многорукие бандиты) поясняются, например, в работе [54].

Сгладим исходную постановку с помощью локального усреднения по евклидову шару радиуса $\mu > 0$, который будет выбран позже,

$$f_k^\mu(x,\xi) = E_{\tilde{e}}\left[f_k(x+\mu\tilde{e},\xi)\right],$$

$$f_k^\mu(x) = E_{\tilde{e},\xi}\left[f_k(x+\mu\tilde{e},\xi)\right],$$

где $\tilde{e} \in RB_2^n(1)$, т.е. случайный вектор $e$ равномерно распределен на шаре радиуса 1 в $l_2$. Заменим исходную задачу (6.2.1) следующей задачей минимизации

$$\text{Regret}_N^\mu\left(\{f_k(\cdot)\},\{x^k\}\right) = \frac{1}{N}\sum_{k=1}^N f_k^\mu(x^k) - \min_{x\in Q}\frac{1}{N}\sum_{k=1}^N f_k^\mu(x). \tag{6.2.4}$$

Это делается для того, чтобы обеспечить выполнение условия 6.2.1 подраздела 6.2.2, см. ниже. Будем считать, что имеют место условия 6.2.6, 6.2.7 (если условие 6.2.7 не выполнено, просто полагаем $L_2 = \infty$). Предположим также, что

$$\min\left\{M_2\mu, L_2\mu^2/2\right\} \le \varepsilon/2,$$

т.е.

$$\mu \le \max\left\{\frac{\varepsilon}{2M_2}, \sqrt{\frac{\varepsilon}{L_2}}\right\}, \tag{6.2.5}$$

где $\varepsilon = \varepsilon(N)$ определятся из условия (можно также сказать, что из этого условия определяется $N = N(\varepsilon)$)

$$E\left[\text{Regret}_N\left(\{f_k(\cdot)\},\{x^k\}\right)\right] \le \varepsilon.$$

Из [52] следует, что при условии (6.2.5), из

$$\text{Regret}_N^\mu\left(\{f_k(\cdot)\},\{x^k\}\right) \le \varepsilon/2$$

для тех же самых последовательностей $\{f_k(\cdot)\},\{x^k\}$, следует

$$\text{Regret}_N\left(\{f_k(\cdot)\},\{x^k\}\right) \le \varepsilon.$$

Далее сконцентрируемся на минимизации сглаженной версии псевдо регрета (6.2.4), контролируя при этом выполнение условия (6.2.5).

Введенные выше $\nabla_x \tilde{f}_k(x;e,\xi)$ для задачи (6.2.4) удовлетворяют условию 6.2.1 с $\sigma$ равным, соответственно,

$$\sigma \le E\left[\frac{\delta n}{N\mu}\sum_{k=1}^{N}\left|\langle e_k, r_k\rangle\right|\right] \le \frac{2\delta R\sqrt{n}}{\mu} \quad (\text{при } m=1), \tag{6.2.6}$$

$$\sigma \le E\left[\frac{2\delta n}{N\mu}\sum_{k=1}^{N}\left|\langle e_k, r_k\rangle\right|\right] \le \frac{4\delta R\sqrt{n}}{\mu} \quad (\text{при } m=2), \tag{6.2.7}$$

где $E\left[r_k^2\right] \le 2R^2$, $e_k \in RS_2^n(1)$ – не зависит от $r_k = x^k - x_*$. Оценки (6.2.6), (6.2.7) следуют из того, что [52, 166, 201]

$$E_e\left[\nabla_x f_k(x; e, \xi)\right] = \nabla f_k^{\mu}(x, \xi),$$

и из явления концентрации равномерной меры на сфере вокруг экватора (при северном полюсе, заданном вектором $r_k$) [244].

Чтобы можно было воспользоваться оценками (6.2.2), (6.2.3) подраздела 6.2.2 осталось в условии 6.2.4 подраздела 6.2.2 оценить константу $M$. Выберем в прямом пространстве норму $l_p$, $1 \le p \le 2$ (см. пример 6.2.2 подраздела 6.2.2). Положим $1/p + 1/q = 1$. При $m = 1$ и условии 6.2.5 имеем оценки [52]

$$M^2 \le \frac{(q-1)n^{1+2/q}B^2}{\mu^2} \quad (\text{при } 2 \le q \le 2\ln n),$$

$$M^2 \le \frac{4n\ln n B^2}{\mu^2} \quad (\text{при } 2\ln n < q \le \infty).$$

Наиболее интересны случаи, когда $q = 2$, $q = \infty$

$$M^2 \le \frac{n^2 B^2}{\mu^2} \quad (\text{при } q = 2), \tag{6.2.8}$$

$$M^2 \le \frac{4n\ln n B^2}{\mu^2} \quad (\text{при } q = \infty). \tag{6.2.9}$$

При $m = 2$ и выполнении условий условия 6.2.5, 6.2.6 имеем оценки [52] (случай $2 < q < \infty$ рассматривается совершенно аналогично)

$$M^2 \le 3nM_2^2 + \frac{3}{4}n^2 L_2^2 \mu^2 + 12\frac{\delta^2 n^2}{\mu^2} \quad (\text{при } q = 2),$$

$$M^2 \le 4\ln n M_2^2 + 3n\ln n L_2^2 \mu^2 + 48\frac{\delta^2 n \ln n}{\mu^2} \quad (\text{при } q = \infty).$$

В частности, если

$$\mu \le \min\left\{\max\left\{\frac{\varepsilon}{2M_2}, \sqrt{\frac{\varepsilon}{L_2}}\right\}, \frac{M_2}{L_2}\sqrt{\frac{4}{3n}}\right\}, \quad \delta \le \frac{M_2 \mu}{\sqrt{12n}} \quad (\text{при } q = 2), \tag{6.2.10}$$

$$\mu \le \min\left\{\max\left\{\frac{\varepsilon}{2M_2},\sqrt{\frac{\varepsilon}{L_2}}\right\},\frac{M_2}{L_2}\sqrt{\frac{1}{6n}}\right\},\ \ \delta \le \frac{M_2\mu}{\sqrt{96n}}\ (\text{при } q=\infty), \tag{6.2.11}$$

то

$$M^2 \le 5nM_2^2\ \ (\text{при } q=2), \tag{6.2.12}$$

$$M^2 \le 5\ln n M_2^2\ \ (\text{при } q=\infty). \tag{6.2.13}$$

Далее, полагая в (6.2.6), (6.2.7), что $\sigma \le \varepsilon/4$, получим дополнительно к (6.2.5) (и (6.2.10), (6.2.11) при $m=2$) условия на $\delta$, $\mu$, $\varepsilon$

$$\frac{2\delta R\sqrt{n}}{\mu} \le \frac{\varepsilon}{4}\ \ (\text{при } m=1),$$

$$\frac{4\delta R\sqrt{n}}{\mu} \le \frac{\varepsilon}{4}\ \ (\text{при } m=2),$$

т.е.

$$\delta \le \frac{\varepsilon\mu}{8R\sqrt{n}}\ \ (\text{при } m=1), \tag{6.2.14}$$

$$\delta \le \frac{\varepsilon\mu}{16R\sqrt{n}}\ \ (\text{при } m=2). \tag{6.2.15}$$

Далее надо воспользоваться оценками (6.2.2), (6.2.3), добиваясь, соответственно,

$$MR\sqrt{\frac{2}{N}} \le \frac{\varepsilon}{4}, \tag{6.2.16}$$

$$\frac{M^2}{2\gamma_2 N}\left(1+\ln N\right) \le \frac{\varepsilon}{4}. \tag{6.2.17}$$

Таким образом, при $m=1$ получаем оценки на $\mu(\varepsilon)$ из (6.2.5), на $\delta(\varepsilon)$ из (6.2.14) и оценки $\mu(\varepsilon)$, на $N(\varepsilon)$ из (6.2.8), (6.2.9), (6.2.16), (6.2.17) и оценки $\mu(\varepsilon)$; при $m=2$ получаем оценки на $\mu(\varepsilon)$ из (6.2.10), (6.2.11), на $\delta(\varepsilon)$ из (6.2.10), (6.2.11), (6.2.15) и оценки $\mu(\varepsilon)$, на $N(\varepsilon)$ из (6.2.12), (6.2.13), (6.2.16), (6.2.17).

Не будем здесь выписывать то, что получается – это довольно тривиально, но достаточно громоздко. Вместо этого, резюмируем полученные в разделе результаты в более наглядной форме. Для этого введем $\tilde{O}(\ )$. Будем считать, что $\tilde{O}(\ )$ – с точностью до логарифмического множителя (от $n$ и(или) $N$) совпадает с $O(\ )$.

Напомним (обозначения см. в подразделе 6.2.2 и примере 6.2.2), что

$$x^{k+1} = \mathrm{Mirr}_{x^k}\left(\alpha_k \nabla_x \tilde{f}_k\left(x^k; e^k, \xi^k\right)\right),\ k=1,\dots,N,$$

где $\{e^k\}_{k=1}^N$ – независимые одинаково распределенные случайные векторы $e^k \in RS_2^n(1)$,

$$\nabla_x \tilde{f}_k\left(x^k; e^k, \xi^k\right) := \frac{n}{\mu} \tilde{f}_k\left(x^k + \mu e^k, \xi^k\right) e^k \quad \text{(при } m = 1\text{)},$$

$$\nabla_x \tilde{f}_k\left(x^k; e^k, \xi^k\right) := \frac{n}{\mu}\left( \tilde{f}_k\left(x^k + \mu e^k, \xi^k\right) - \tilde{f}_k\left(x^k, \xi^k\right)\right) e^k \quad \text{(при } m = 2\text{)},$$

$$\alpha_k \equiv \alpha = \frac{R}{M}\sqrt{\frac{2}{N}}$$

– в общем случае и

$$\alpha_k \equiv \left(\gamma_2 k\right)^{-1},$$

если $f_k(x)$ – $\gamma_2$-сильно выпуклые функции в $l_2$ (в этом случае выбирают $p = 2$).

**Теорема 6.2.1.** *Пусть рассматривается задача стохастической онлайн оптимизации (6.2.1), в постановке, описанной в этом разделе (в безградиентном варианте). Пусть выбрана $l_p$-норма, $1 \le p \le 2$, (см. подраздел 6.2.2). Согласно этой норме задана прокс-функция и расстояние Брэгмана $V_x(y)$. Пусть $R^2 = V_{x^1}(x_*)$, где $x^1$ и $x_*$ определены в подразделе 6.2.2. Тогда*

$$E\left[\operatorname{Regret}_{N(\varepsilon)}\left(\{f_k(\cdot)\}, \{x^k\}\right)\right] \le \varepsilon,$$

*где $N(\varepsilon)$ определяется в Таблицах 6.2.1, 6.2.2.*

Таблица 6.2.1

| $m = 1$ | $f_k(x)$ – выпуклые функции | $f_k(x)$ – $\gamma_2$-сильно выпуклые функции в $l_2$ норме и $p = 2$ |
|---|---|---|
| Выполнены условия 6.2.5, 6.2.6 | $\tilde{O}\left(\frac{B^2 M_2^2 R^2 n^{1+2/q}}{\varepsilon^4}\right)$ | $\tilde{O}\left(\frac{B^2 M_2^2 n^2}{\gamma_2 \varepsilon^3}\right)$ |
| Выполнены условия 6.2.5, 6.2.7 | $\tilde{O}\left(\frac{B^2 L_2 R^2 n^{1+2/q}}{\varepsilon^3}\right)$ | $\tilde{O}\left(\frac{B^2 L_2 n^2}{\gamma_2 \varepsilon^2}\right)$ |

Таблица 6.2.2

| $m = 2$ | $f_k(x)$ – выпуклые функции | $f_k(x)$ – $\gamma_2$-сильно выпуклые функции в $l_2$ норме и $p = 2$ |
|---|---|---|
| Выполнены условия 6.2.5, 6.2.6 | $\tilde{O}\left(\frac{M_p^2 R^2 n^2}{\varepsilon^2}\right)$ | $\tilde{O}\left(\frac{M_2^2 n^2}{\gamma_2 \varepsilon}\right)$ |
| Выполнены условия 6.2.5, 6.2.6, 6.2.7 | $\tilde{O}\left(\frac{M_2^2R^2n^{2/q}}{\varepsilon^2}\right)$ | $\tilde{O}\left(\frac{M_2^2n}{\gamma_2\varepsilon}\right)$ |

Обе таблицы заполняются исходя из описанной выше техники. Исключением является вторая строчка Таблицы 6.2.2, ее мы взяли из [52]. Несложно выписать точные формулы вместо $\tilde{O}(\ )$ во всех полях обеих таблиц. Также несложно выписать условие на допустимый уровень шума $\delta$, при котором мультипликативная константа в точной формуле увеличится, скажем, не более чем в два раза.

Оценки в третьей строчке Таблицы 6.2.2 неулучшаемы [190] (соответствуют нижним оценкам). Оценки во второй строчке Таблицы 6.2.2 неулучшаемы по $\varepsilon$ [91, 128]. Все сказанное выше касается и стохастических, но не онлайн постановок [91, 128].

Относительно Таблицы 6.2.1 имеется гипотеза, что приведенные оценки – неулучшаемы по $n$. По $\varepsilon$ оценки могут быть улучшены за счет ухудшения того, как входит $n$ [167], см. также http://arxiv.org/pdf/1607.03084v1.pdf .

В заключение рассмотрим пример, демонстрирующий, что полученные в теореме 6.2.1 результаты представляются интересными не только в онлайн контексте.

**Пример 6.2.3.** Предположим, что “успешность” некоторого человека зависит от того, как он распоряжается своим временем. Имеется $n$ различных родов деятельности. В $k$-день человек распоряжается своим временем согласно вектору $x^k \in S_n(1)$. Этот вектор отражает доли времени, уделенного соответствующим делам. В конце каждого дня человек получает “обратную связь” от “внешнего мира” вида

$$\tilde{f}\left(x^k+\mu e^k,\xi^k\right)=f\left(x^k+\mu e^k\right)\cdot\left(1+\xi^k\right),$$

где $\left\{e^k\right\}$ – независимые одинаково распределенные случайные векторы $e^k \in RS_2^n(1)$, $\mu$ определяется согласно формуле (6.2.5), $\left\{\xi^k\right\}$ – независимые (между собой и от $\left\{e^k\right\}\in RS_2^n(1)$) одинаково распределенные случайные величины $\xi^k \in N(0,1)$, а выпуклая функция $f(x)$, со свойствами

$$\left|f(x)\right|\le B,\ \left|f(x)-f(y)\right|\le M_2\left\|x-y\right\|_2,$$

правильно отражает реальное “положение дел”, т.е. минимум этой функции соответствует оптимальной для данного человека конфигурации. Задача человека, заключается в том, чтобы (получая каждый день описанную выше обратную связь), так организовать “про-цесс своего обучения” (на основе получаемой информации), чтобы как можно быстрее достичь такого состояния[5]

$$\bar{x}^N = \frac{1}{N}\sum_{k=1}^{N} x^k \in S_n(1),$$

что с вероятностью 0.999 имеет место неравенство

$$f\left(\bar{x}^N\right) - \min_{x\in S_n(1)} f(x) \le \text{Regret}_N\left(\left\{f(\cdot)\right\},\left\{x^k\right\}\right) \le \varepsilon .$$

Согласно теореме 6.2.1 человек может этого достичь за

$$N = \mathrm{O}\left(\frac{B^2 M_2^2 n \ln^2 n}{\varepsilon^4}\right)$$

дней. Если есть возможность получать каждый день информацию $\tilde{f}\left(x^k + \mu e^k, \xi^k\right)$ и $\tilde{f}\left(x^k, \xi^k\right)$, где $\mu$ определяется согласно формуле (6.2.11), то тогда можно улучшить оценку до

$$N = \mathrm{O}\left(\frac{M_2^2 \ln^2 n}{\varepsilon^2}\right)$$

дней – здесь прдполагается также гладкость $f(x)$. □

[5] Заметим, что с достаточной точностью и доверительным уровнем (при $N \gg 1$) можно считать, что

$$\bar{x}^N \simeq \frac{1}{N}\sum_{k=1}^{N}\left(x^k + \mu e^k\right).$$

# Заключение

В диссертации разработаны новые подходы к построению многостадийных моделей транспортных потоков и эффективные численные методы поиска равновесий в таких моделях.

Предложен общий способ построения многостадийных моделей транспортных потоков, с помощью которого можно вывести как уже существующие многостадийные модели, так и предложить новые. В частности, предложить многостадийные модели с заменой устаревшего блока Бэкмана более современным блоком Стабильной Динамики. В диссертации продемонстрировано, как сводить поиск равновесий в многостадийных транспортных моделях к решению задач выпуклых оптимизации со специально иерархической (многоуровневой) сетевой структурой. В основе численных подходов к поиску равновесий лежат два базисных (прямо-двойственных) алгоритма: метод зеркального спуска и быстрый градиентный метод. С помощью небольшого набора операций удалось построить наиболее эффективные методы. Эффективность методов оценивалась исходя из известных нижних оценок и структуры задачи.

Развивая численные методы поиска равновесий в многостадийных транспортных моделях, в диссертации единообразно были рассмотрены основные концепции современных численных методов выпуклой оптимизации на множествах простой структуры. Это позволило "овыпуклить" данную область знаний. Были получены новые результаты для методов с неточным оракулом. Новые результаты были получены при изучении прямо-двойственных свойств методов и при исследовании аспектов, связанных с оптимизацией на неограниченных областях. Для класса задач huge-scale оптимизации изучена "игра" на рандомизации и разреженности. Для спусков по направлению и методов нулевого порядка изучена "игра" на структуре множества, на котором происходит оптимизация, и способе выбора случайного направления.

Из-за ограничений на объем (не более 500 стр.) в диссертацию не были включены работы [14, 23, 33], дополняющие главу 1 (последняя статья хорошо связывает главу 1 с главами 2, 3), не вошла работа [51], которая связывает главу 1 с главой 6. Удачно дополнить вторую главу могла бы работа [193], содержащая пропущенные доказательства используемых в этой главе фактов, связанных с неточным оракулом. Работы [171] уточняет результаты главы 3. Цикл работ [154, 155] связывал между собой материалы глав 4 и 5.

Перспективным представляется последующее развитие изученной в диссертации (на конкретных примерах) конструкции сборки итогового "оптимального" метода для задач со сложной структурой из базисных методов, оптимальных для отдельных частей.

За развитием полученных в диссертации результатов, можно следить, например, по следующей ссылке https://arxiv.org/find/all/1/all:+gasnikov/0/1/0/all/0/1.

# СПИСОК ЛИТЕРАТУРЫ

# Приложение А. Теория макросистем с точки зрения стохастической химической кинетики

## А.1 Введение

В данном приложении мы постараемся пояснить важность формализма стохастической химической кинетики в изучении равновесий макросистем. Под макросистемой мы понимаем систему большого количества случайно взаимодействующих агентов. Мы будем интересоваться поведением такой системы на больших временах, в частности, необходимыми и достаточными условиями, при которых на больших временах такая система сходится к равновесию. При этом под равновесием будем понимать такое макросостояние, в малой окрестности которого концентрируется стационарная мера (динамику считаем марковской). Имеется большое количество конкретных примеров макросистем [20, 25, 26, 28, 46, 124, 188, 246, 253, 284, 324], встречающихся в самых разнообразных предметных областях (физика, биология, экономика, транспортное моделирование и т.д.). Однако в подавляющем большинстве случаев, при описании макросистемы не описывается ее эволюция (считается, что система уже находится в стационарном состоянии), не исследуется концентрация в окрестности равновесия. В данном приложении будет изучаться эволюция макросистем и концентрация, возникающих при этом стационарных мер. В отличие от подавляющего большинства работ, мы будем допускать, что макросистема характеризуется вектором чисел заполнений, размерность которого может зависеть от числа агентов.

Структура приложения следующая. В разделе А.2 мы приведем три игрушечных примера макросистем, изучение равновесий в которых может быть единообразно осуществлено. Пункт г) теоремы из А.3 является оригинальным, и представляющим наибольший интерес в этом приложении.

## А.2 Примеры макросистем

Далее мы рассматриваем три хорошо известных примера макросистем. Однако способ преподнесения материала представляется оригинальным. Также отметим, что приводимые ниже количественные оценки мы ранее не встречали в литературе. Эти оценки получаются исходя из формализма, описанного в разделе А3, применимого к большому количеству макросистем.

**Пример А.1 (Кинетика социального неравенства и предельные формы).** В некотором городе живет $N \gg 1$ жителей (четное число). В начальный момент у каждого жителя имеется по $\bar{s}$ монеток. Каждый день жители случайно разбиваются на пары. В каждой паре жители скидываются по монетке (если один или оба участника банкроты, то банкрот не скидывается, в то время, как не банкрот, в любом случае, обязан скинуть монет-ку). Далее в каждой паре случайно разыгрывается победитель, который и забирает “призовой фонд”. Обозначим через $c_s(t)$ – долю жителей города, у которых ровно $s$, $s = 0,...,\overline{s}N$, монеток на $t$-й день. Имеет место следующие оценки

$$\exists\ a>0:\ \forall\ \sigma>0, t \ge aN\ln N \to P\left(\left\|c(t)-c^*\right\|_2 \ge \frac{2\sqrt{2}+4\sqrt{\ln\left(\sigma^{-1}\right)}}{\sqrt{N}}\right) \le \sigma,$$

$$\exists\ b,D>0:\ \forall\ \sigma>0, t \ge bN\ln N \to P\left(\left\|c(t)-c^*\right\|_1 \ge D\sqrt{\frac{\ln^2 N+\ln\left(\sigma^{-1}\right)}{N}}\right) \le \sigma,$$

где $c_s^* \simeq C\exp(-s/\overline{s})$, а $C \simeq 1/\overline{s}$ находится из условия нормировки $\sum_{s=0}^{\overline{s}N} C\exp(-s/\overline{s}) = 1$. Таким образом, кривая (предельная форма [25], равновесие макросистемы [188]) $C\exp(-s/\overline{s})$ характеризует распределение населения по богатству на больших временах. Этот результат восходит к работам итальянского экономиста Вильфредо Парето (см., например, [188]), пытавшегося объяснить расслоение населения по богатству. Насколько нам известно, ранее в таком контексте не приводились оценки скорости сходимости и плотности концентрации.

Для объяснения этого результата полезно рассмотреть схожий процесс, в котором каждой паре жителей приписан свой (независимый) “пуассоновский будильник” (звонки происходят в случайные моменты времени, соответствующие скачкам пуассоновского процесса; параметр интенсивность этого пуассоновского процесса называют интенсивностью/параметром будильника). Все будильники “приготовлены” одинаково: у всех у них одна и та же интенсивность $\lambda N^{-1}$. Далее следует погрузить задачу в модель стохастической химической кинетики с бинарными реакциями и воспользоваться результатом, приведенным в А.3. Наиболее технически сложными моментами является оценка с помощью неравенства Чигера mixing time $\sim N\ln N$ [246] и получение поправки под корнем $\ln^2 N$.

К такой же предельной форме можно было бы прийти и по-другому (подход Булгакова–Маслова “разбрасывание червонцев в варьете”). В некотором городе живет $N \gg 1$ жителей (изначально банкротов). Каждый день одному из жителей (случайно выбранному в этот день) дается одна монетка. Тогда $\left\{c_s\left(\overline{s}N\right)\right\}_{s=0}^{\overline{s}N} \xrightarrow[N\to\infty]{} \left\{C\exp(-s/\overline{s})\right\}_{s=0}^{\overline{s}N}$.

**Пример А.2 (обезьянка и печатная машинка; закон Ципфа–Мандельброта).** На печатной машинке $n+1$ символ, один из символов пробел. Обезьянка на каждом шаге случайно (независимо и равновероятно) нажимает один из символов. Прожив долгую жизнь, обезьянка сгенерировала текст огромной длинны. По этому тексту составили сло-варь. Этот словарь упорядочили по частоте встречаемости слова (слова – это всевозможные набор букв без пробелов, которые хотя бы раз встречались в тексте обезьянки между какими-то двумя пробелами). Так на первом месте в словаре поставили самое часто встречаемое слово, на второе поставили второе по частоте встречаемости и т.д. Номер слова в таком словаре называется рангом и обозначается буквой $r$. Предельная форма кривой, описывающей распределение частот встречаемости слов от рангов, имеет вид

$$\text{Частота}(r) \simeq \frac{C}{(r+B)^{\alpha}},$$

где

$$\alpha = \frac{\log(n+1)}{\log(n)},\ B = \frac{n}{n-1},\ C = \frac{n^{\alpha-1}}{(n-1)^{\alpha}}.$$

Такой вывод закона Ципфа (с поправкой) был одним из двух, предложенных Бенуа Мандельбротом [324]. В работе [124] была предложена довольно общая схема, приводящая к закону Ципфа, в которую можно погрузить и обезьянку с печатной машинкой. А именно, предположим, что динамика (порождения слов в большом тексте) такова, что вероятность того, что в тексте из $N \gg 1$ слов $x_1$ (первое по порядку слово в ранговом словаре) встречалось $N_1$ раз, $x_2$ (второе по порядку слово в ранговом словаре) встречалось $N_2$ раза и т.д. есть

$$\sim \frac{N!}{N_1!N_2!\ldots}\exp\left(-\eta\sum_{k\in\mathbb{N}} N_k E_k\right),$$

где $E_k$ – число букв в слове с рангом $k$. Часто считают, что $\eta = 0$, но зато динамика такова, что число слов и число букв становятся асимптотически (по размеру текста) связанными (закон больших чисел). Таким образом, к закону сохранения $\sum_{k\in\mathbb{N}} N_k = N$ добавляется приближенный закон сохранения $\sum_{k\in\mathbb{N}} N_k E_k \simeq \bar{E}N$ ($\bar{E}$ – среднее число букв в слове). Поиск предельной формы приводит к задаче (воспользовались формулой Стирлинга и методом множителей Лагранжа)

$$\sum_{k\in\mathbb{N}}\{N_k \ln N_k + \lambda E_k N_k\} \to \min_{\substack{N_k \ge 0 \\ \sum_{k\in\mathbb{N}} N_k = N}},$$

где $\lambda$ – либо равняется $\eta$, либо является множителем Лагранжа к ограничению $\sum_{k\in\mathbb{N}} N_k E_k \simeq \bar{E}N$. Стоит отметить, что к аналогичной задаче (с $E_k = k$) приводит поиск предельной формы в модели примера А.1. Решение нашей задачи дает $N_k = \exp(-\mu - \lambda E_k)$, где $\mu$ – множитель Лагранжа к ограничению $\sum_{k\in\mathbb{N}} N_k = N$. Далее, считают, что $r(E)$ – число различных используемых слов с числом букв не большим $E$, приближенно представимо в виде $r(E) \simeq a^E$. Тогда $N_k \sim k^{-\gamma}$, где $\gamma = \lambda/\ln a$.

**Пример А.3 (теорема Гордона–Ньюэлла и PageRank [46]).** Имеется $N \gg 1$ пользователей, которые случайно (независимо) блуждают в непрерывном времени по ориентированному графу (на $m$ вершинах) с эргодической инфинитезимальной матрицей $\Lambda$. Назовем вектор $p$ (из единичного симплекса) PageRank, если $\Lambda p = 0$. Обозначим через $n_i(t)$ – число пользователей на $i$-й странице в момент времени $t \ge 0$. Не сложно показать (теорема Гордона–Ньюэлла [305]), что $n(t)$ асимптотически имеет мультиномиальное распределение с вектором параметров PageRank $p$, т.е.

$$\lim_{t\to\infty} P(n(t) = n) = \frac{N!}{n_1!\cdot\ldots\cdot n_m!} p_1^{n_1} \cdot\ldots\cdot p_m^{n_m}.$$

Следовательно (неравенство Хефдинга в гильбертовом пространстве в форме [36]),

$$\lim_{t\to\infty} P\left( \left\| \frac{n(t)}{N} - p \right\|_2 \ge \frac{2\sqrt{2} + 4\sqrt{\ln(\sigma^{-1})}}{\sqrt{N}} \right) \le \sigma.$$

Этот же результат можно получить, рассмотрев соответствующую систему унарных химических реакций. Переход одного из пользователей из вершины $i$ в вершину $j$ – означает превращение одной молекулы вещества $i$ в одну молекулу вещества $j$, $n_i(t)$ – число молекул $i$-го типа в момент времени $n_i(t)$. Каждое ребро графа соответствует определенной реакции (превращению). Интенсивность реакций определяется матрицей $\Lambda$ и числом молекул, вступающих в реакцию (закон действующих масс). Условие $\Lambda p = 0$ – в точности соответствует условию унитарности в стохастической химической кинетике (см. А.3).

Отметим, что если воспользоваться теоремой Санова о больших уклонениях для мультиномиального распределения [111], то получим

$$\frac{N!}{n_1!\cdot\ldots\cdot n_m!} p_1^{n_1} \cdot\ldots\cdot p_m^{n_m} = \exp\left( -N \sum_{i=1}^{m} \nu_i \ln(\nu_i/p_i) + R \right),$$

где $\nu_i = n_i/N$, $|R| \le \frac{m}{2}(\ln N + 1)$. Однако последующее применение неравенства Пинскера не дает нам равномерной по $m$ оценки в 1-норме. Как и ожидалось, выписанная в 2-норме оценка (правая часть неравенства под вероятностью) и так полученная оценка в 1-норме будут отличаться по порядку приблизительно в $\sqrt{m}$ раз, что соответствует типичному (по

Б.С. Кашину) соотношению между 1 и 2 нормами. Слово “типично” здесь отвечает, грубо говоря, за ситуацию, когда компоненты вектора одного порядка. При этом важно отметить, что для многих приложений, где возникают предельные конфигурации (кривые), описывающиеся вектором с огромным числом компонент, имеет место быстрый закон убывания этих компонент (см., примеры 1, 2 и [253, 284]), т.е. такая ситуация “не типична” и можно ожидать лучшие оценки концентрации в 1-норме (см. пример А.1).

## А.3 Изучение динамики макросистемы с точки зрения стохастической химической кинетики

Предположим, что некоторая макросистема может находиться в различных состояниях, характеризуемых вектором $n$ с неотрицательными целочисленными компонентами. Будем считать, что в системе происходят случайные превращения (химические реакции).

Пусть $n \to n-\alpha+\beta$, $(\alpha,\beta)\in J$ – все возможные типы реакций, где $\alpha$ и $\beta$ – вектора с неотрицательными целочисленными компонентами. Введем интенсивность реакции:

$$\lambda_{(\alpha,\beta)}(n)=\lambda_{(\alpha,\beta)}(n\to n-\alpha+\beta)=N^{1-\sum_i \alpha_i} K_\beta^\alpha \prod_{i:\alpha_i>0} n_i\cdot\ldots\cdot(n_i-\alpha_i+1),$$

где $K_\beta^\alpha \ge 0$ – константа реакции; при этом $\sum_{i=1}^m n_i(t)\equiv N \gg 1$. Другими словами, $\lambda_{(\alpha,\beta)}(n)$ – вероятность осуществления в единицу времени перехода $n\to n-\alpha+\beta$. Здесь не предполагается, что число состояний $m=\dim n$ и число реакций $|J|$ не зависят от числа агентов $N$. Тем не менее, если ничего не известно про равновесную конфигурацию $c^*$ (типа быстрого убывания компонент этого вектора), то дополнительно предполагается, что $m \ll N$ – это нужно для обоснования возможности применения формулы Стирлинга при получении вариационного принципа (максимума энтропии) для описания равновесия макросистемы $c^*$ (в концентрационной форме). Однако часто априорно можно предполагать (апостериорно проверив), что компоненты вектора $c^*$ убывают быстро, тогда это условие можно отбросить. Так, например, обстоит дело в примере А.1 из раздела А.2. Описанный выше марковский процесс считается неразложимым. Имеет место

**Теорема. а)** $\langle \mu, n(t)\rangle \equiv \langle \mu, n(0)\rangle$ (inv) $\Leftrightarrow$ *вектор* $\mu$ *ортогонален каждому вектору семейства* $\{\alpha-\beta\}_{(\alpha,\beta)\in J}$. *Здесь* $\langle\,\cdot,\cdot\,\rangle$ *– обычно евклидово скалярное произведение.*

**б)** *Если существует* $\lim\limits_{N\to\infty} n(0)/N = c(0),\, K_\beta^\alpha := K_\beta^\alpha(n/N)$, *m и* $|J|$ *не зависит от* $N$, *то для любого* $t>0$ *с вероятностью 1 существует* $\lim\limits_{N\to\infty} n(t)/N = c(t)$, *где* $c(t)$ *– не случайная вектор-функция, удовлетворяющая СОДУ Гульдберга–Вааге:*

$$\frac{dc_i}{dt} = \sum_{(\alpha,\beta)\in J} (\beta_i - \alpha_i) K_\beta^\alpha(c) \prod_j c_j^{\alpha_j}. \qquad \text{(ГВ)}$$

*Гиперплоскость* (inv) *(с очевидной заменой* $n \Rightarrow c$*) инвариантна относительно этой динамики.*[1] *Более того, случайный процесс* $n(t)/N$ *слабо сходится при* $N\to\infty$ *к* $c(t)$ *на любом конечном отрезке времени.*

**в)** *Пусть выполняется условие унитарности*[2] *(очевидно, что* $\xi$, *удовлетворяющий условию (U), – неподвижная точка в (ГВ))*

$$\exists\ \xi>0:\ \forall\ \beta \to \sum_{\alpha:(\alpha,\beta)\in J} K_\beta^\alpha \prod_j (\xi_j)^{\alpha_j} = \sum_{\alpha:(\alpha,\beta)\in J} K_\alpha^\beta \prod_j (\xi_j)^{\beta_j}. \qquad \text{(U)}$$

*Тогда неотрицательный ортант* $\mathbb{R}_+^m$ *расслаивается гиперплоскостями* (inv)*, так что в каждой гиперплоскости* (inv) *уравнение (U) (положительно) разрешимо притом единственным образом. Следовательно, существует, притом единственная, неподвижная точка* $c^* \in (\text{inv})$ *у системы (ГВ), являющаяся при этом глобальным аттрактором. Система (ГВ) имеет функцию Ляпунова* $KL(c,\xi) = \sum_{i=1}^m c_i \ln(c_i/\xi_i)$.

*Стационарное распределение описанного марковского процесса имеет носителем множество* (inv) *и (с точностью до нормирующего множителя) имеет вид*

$$\frac{N!}{n_1!\cdot\ldots\cdot n_m!}(\xi_1)^{n_1}\cdot\ldots\cdot(\xi_m)^{n_m} \sim \exp(-N\cdot KL(c,\xi)),$$

*где* $\xi$ *– произвольное решение (U), не важно какое именно (от этого, конечно, будет зависеть нормирующий множитель, но это ни на чем не сказывается). При этом условие унитарности (U) является обобщением условия детального равновесия*[3] *(баланса)*[4]

---

[1] Можно сказать, что “жизнь” нелинейной динамической системы определяется линейными законами сохранения, унаследованными ею при скейлинге (каноническом). Этот тезис имеет, по-видимому, более широкое применение [24].

[2] Иногда это условие называют “условием Штюкельберга–Батищевой–Пироговой”. В частности, такое название, предлагает, В.В. Веденяпин, имеющий к этому условию самое непосредственное отношение. На наш взгляд, если все-таки делать это условие именным, то точнее было бы его тогда назвать “условие Штюкельберга–Батищевой–Веденяпина–Малышева–Пирогова–Рыбко”. Из-за громоздкого вида последнего, было отдано предпочтение варианту, предложенному С.А. Пироговым: “условие унитарности”.

$$\exists\ \xi > 0:\ \forall\ (\alpha,\beta)\in J \to K_{\beta}^{\alpha}\prod_j (\xi_j)^{\alpha_j} = K_{\alpha}^{\beta}\prod_j (\xi_j)^{\beta_j},$$

*принимающего такой вид для мультиномиальной стационарной меры.*

*Существует такая зависимость* $a(m,c(0))$ *(во многих приложениях можно убрать второй аргумент в этой зависимости), что*

$$\exists\ a = a(m,c(0)):\ \forall\ \sigma > 0, t \ge aN\ln N \to P\left(\left\|\frac{n(t)}{N} - c^*\right\|_2 \ge \frac{2\sqrt{2}+4\sqrt{\ln(\sigma^{-1})}}{\sqrt{N}}\right) \le \sigma,$$

*где (принцип максимума энтропии, Больцман–Джейнс)*

$$c^* = \arg\max_{c\in(\text{inv})}\left(-\sum_i c_i \ln(c_i/\xi_i)\right) = \arg\min_{c\in(\text{inv})} KL(c,\xi),$$

*а* $\xi$ *– произвольное решение (U), причем* $c^*$ *определяется единственным образом, т.е. не зависит от выбора* $\xi$*. Геометрически* $c^*$ *– это KL-проекция произвольного* $\xi$*, удовлетворяющего (U), на гиперплоскость (inv), соответствующую начальным данным* $c(0)$*. Независимость этой проекции от выбора* $\xi$ *из (U) просто означает, что кривая (U) проходит KL-перпендикулярно через множество (inv).*

**г)** *Верно и обратное утверждение, то есть условие (U) не только достаточное для того чтобы равновесие находилось из приведенной выше задачи энтропийно-линейного программирования, но и, с небольшой оговоркой (для почти всех* $c(0)$*), необходимое. Также верно и более общее утверждение, связывающее понимание энтропии в смысле Больцмана (функция Ляпунова прошкалированной кинетической динамики) и Санова (функционал действия в неравенствах больших уклонений для стационарной меры): если стационарная мера асимптотически представима в виде* $\sim \exp(-N\cdot V(c))$*, то* $V(c)$ *– функция Ляпунова (ГВ).*

Результаты п. а) взяты из [12], п. б) из [86, 195], п. в) из [12, 86], п. г) из [32]. Выписанная оценка скорости сходимости и плотности концентрации в п. в) ранее нам не встре-чалась. Также отметим, что ранее результаты п. г) были получены при дополнительном предположении $K_{\beta}^{\alpha}(n/N) \equiv K_{\beta}^{\alpha}$.

---

[3] В терминах примера А.3 условие унитарности просто означает, что в равновесии для любой вершины имеет место баланс числа пользователей входящих в единицу времени в эту вершину с числом пользователей, выходящих в единицу времени из этой вершины. В то время как условие детального равновесия означает, что в равновесии для любой пары вершин число пользователей, переходящих в единицу времени из одной вершину в другую равно числу пользователей, переходящих в обратном направлении. Понятно, что второе условие является частным случаем первого.

[4] Много интересных примеров макросистем, в которых $K_{\beta}^{\alpha} := K_{\beta}^{\alpha}(n/N)$, и имеет место детальный баланс, собрано в книге [20].

К п. г) можно сделать следующее оригинальное пояснение. Обозначим через $h(c)$ вектор-функцию, стоящую в правой части СОДУ (ГВ). Тогда (см., например, [28, 195]) при $N \gg 1$ по теореме Т. Куртца $n(t)/N$ будет $O\left(\log N/\sqrt{N}\right)$-близко к $x_t = x(t)$ – решению стохастической системы дифференциальных уравнений (с начальным условием $x_0 = c(0)$)

$$dx_t = h(x_t)dt + \sqrt{\frac{g(x_t)}{N}}dW_t,$$

где функция $g(x_t) > 0$ рассчитывается по набору реакций и константам реакций (которые могут быть не постоянны и зависеть от концентраций), $W_t$ – стандартный винеровский процесс. Стационарная мера $m(x) = \lim\limits_{t\to\infty} p(t,x)$ этого однородного марковского процесса удовлетворяет уравнению

$$\frac{1}{2N}\nabla^2\left(g(x)m(x)\right) - \operatorname{div}\left(h(x)m(x)\right) = 0,$$

поскольку плотность распределения $p(t,x)$ процесса $x_t$ подчиняется уравнению Колмогорова–Фоккера–Планка

$$\frac{\partial p(t,x)}{\partial t} = -\operatorname{div}\left(h(x)p(t,x)\right) + \frac{1}{2N}\nabla^2\left(g(x)p(t,x)\right).$$

Здесь использовалось обозначение: $\nabla^2 f(x) = \sum_{i,j}\partial^2 f(x)/\partial x_i \partial x_j$. Если известно, что $m(x) \simeq \mathrm{const}\cdot\exp\left(-N\cdot V(x)\right)$, то из уравнения на $m(x)$ имеем

$$N\langle h,\nabla V\rangle - \operatorname{div} h - \frac{1}{2}\langle\nabla g,\nabla V\rangle + \frac{1}{2N}V\nabla^2 g - \frac{1}{2}g\nabla^2 V + \frac{N}{2}g\cdot(\nabla V)^2 \simeq 0,$$

следовательно

$$\langle h,\nabla V\rangle \simeq -\frac{1}{2}g\cdot(\nabla V)^2 + O\left(\frac{1}{N}\right) \underset{N\to\infty}{=} -\frac{1}{2}g\cdot(\nabla V)^2 \le 0.$$

Эта выкладка поясняет, почему функция $V(c)$ может быть функцией Ляпунова системы (ГВ) $dc/dt = h(c)$. Более того, модель стохастической химической кинетики здесь может быть заменена и более общими шкалирующимися марковскими моделями [195].

В заключение отметим, что много различных примеров макросистем собрано в главе 6 книги [19].